\documentclass[12pt]{article}
\pdfoutput=1
\usepackage{fullpage}

\usepackage{amsmath,amsfonts,amssymb,graphics,amsthm}
\usepackage{graphicx}
\usepackage{subfigure}
\usepackage{nicefrac}
\usepackage[margin=10pt,font=small,labelfont=bf,labelsep=period]{caption}
\usepackage{multicol}
\allowdisplaybreaks

\usepackage{enumitem}
\usepackage{xcolor}

\usepackage{titlesec}

\usepackage{hyperref}
\hypersetup{
    colorlinks=true,
    linkcolor=blue,
    citecolor=red,
    urlcolor=blue,
    pdfborder={0 0 0}
}

\usepackage{cleveref}
  \crefname{theorem}{Theorem}{Theorems}
  \crefname{thm}{Theorem}{Theorems}
  \crefname{lemma}{Lemma}{Lemmas}
  \crefname{lem}{Lemma}{Lemmas}
  \crefname{remark}{Remark}{Remarks}
    \crefname{rmk}{Remark}{Remarks}
  \crefname{prop}{Proposition}{Propositions}
\crefname{notation}{Notation}{Notations}
\crefname{claim}{Claim}{Claims}
  \crefname{defn}{Definition}{Definitions}
   \crefname{definition}{Definition}{Definitions}
  \crefname{cor}{Corollary}{Corollaries}
  \crefname{section}{Chapter}{Chapters}
  \crefname{subsection}{Section}{sections}
  \crefname{figure}{Figure}{Figures}
  \crefname{example}{Example}{Examples}
    \crefname{assumption}{Assumption}{Assumptions}

\newtheorem{theorem}{Theorem}[section]

\newtheorem{claim}[theorem]{Claim}
\newtheorem{lemma}[theorem]{Lemma}
\newtheorem{cor}[theorem]{Corollary}
\newtheorem{prop}[theorem]{Proposition}
\newtheorem{definition}[theorem]{Definition}
\newtheorem{example}[theorem]{Example}

\numberwithin{equation}{section}

\theoremstyle{definition}
\newtheorem{rmk}[theorem]{Remark}

\newcommand{\dis}[1]{\begin{center} \rule[.2em]{\textwidth}{.1em} {#1} \rule{\textwidth}{.1em} \end{center}}
\usepackage[T1]{fontenc}
\usepackage[utf8]{inputenc}

\def\cW{\mathcal{W}}
\def\cV{\mathcal{V}}

\def\cT{\mathcal{T}}
\def\cS{\mathcal{S}}
\def\cR{\mathcal{R}}
\def\cQ{\mathcal{Q}}
\def\cP{\mathcal{P}}

\def\cN{\mathcal{N}}
\def\cM{\mathcal{M}}
\def\cL{\mathcal{L}}
\def\cK{\mathcal{K}}

\def\cI{\mathcal{I}}
\def\cH{\mathcal{H}}
\def\cG{\mathcal{G}}
\def\cF{\mathcal{F}}
\def\cE{\mathcal{E}}
\def\cD{\mathcal{D}}
\def\cC{\mathcal{C}}

\def\cA{\mathcal{A}}

\def\k{\kappa}
\def\h{\mathbf{h}}

\def\e{\mathrm{e}}
\def\P{\mathbb{P}}
\def\Q{\mathbb{Q}}

\def\H{\mathbb{H}}
\def\D{\mathbb{D}}

\def\E{\mathbb{E}}
\def\C{\mathbb{C}}
\def\R{\mathbb{R}}
\def\Z{\mathbb{Z}}
\def\N{\mathbb{N}}
\def\Sp{\mathbb{S}}

\def\ve{\vec{e}}

\def\la{\langle}
\def\ra{\rangle}
\def\eps{\varepsilon}

\def\pd{\partial}
\def\ph{\varphi}
\def\az{{\alpha,\mathbf{z}}}
\def\haz{h^L_{\az}}

\def\kp{{\kappa'}}

\def\fc{\mathfrak{c}}
\DeclareMathOperator{\var}{Var}

\DeclareMathOperator{\dist}{dist}
\DeclareMathOperator{\Leb}{Leb}
\DeclareMathOperator{\dd}{d\!}

\DeclareMathOperator{\tv}{TV}
\DeclareMathOperator{\rr}{rad}
\DeclareMathOperator{\cc}{circ}
\DeclareMathOperator{\cov}{Cov}
\DeclareMathOperator{\GFF}{GFF}
\DeclareMathOperator{\wg}{wedge}
\DeclareMathOperator{\DN}{DN}

\newcommand{\Supp}{H_0^1}
\DeclareMathOperator{\Harm}{Harm}
\DeclareMathOperator{\Conv}{Conv}

\newcommand{\indic}[1]{\mathbf{1}_{\{#1\}}}

\newcommand{\old}[1]{{}}

\newcommand{\ellen}[1]{{\textcolor{magenta}{#1}}}
\newcommand{\nb}[1]{{\textcolor{orange}{#1}}}

\def\mn{\medskip \noindent}

\usepackage{comment}

\usepackage[quiet]{imakeidx}
\indexsetup{level=\subsection*,noclearpage}
\makeindex[title=Index]
\makeindex[name=notation,title=Notation and Symbols,columns=2]
\usepackage{hyperref}
\hypersetup{
	colorlinks=true,
	linkcolor=blue,
	citecolor=red,
	urlcolor=blue,
	pdfborder={0 0 0}
}

\newcommand{\ind}[1]{\index{#1|BH}}
\newcommand{\indN}[1]{\index[notation]{#1|BH}}

\begin{document}

\title{
Gaussian free field\\
%Gaussian multiplicative chaos\\
and Liouville quantum gravity\\}
%and Gaussian multiplicative chaos}

\author{Nathana\"el Berestycki \and Ellen Powell}

\date{[Draft: \today]}

\maketitle
\vspace{1cm}
\begin{center}
\includegraphics[scale=.3]{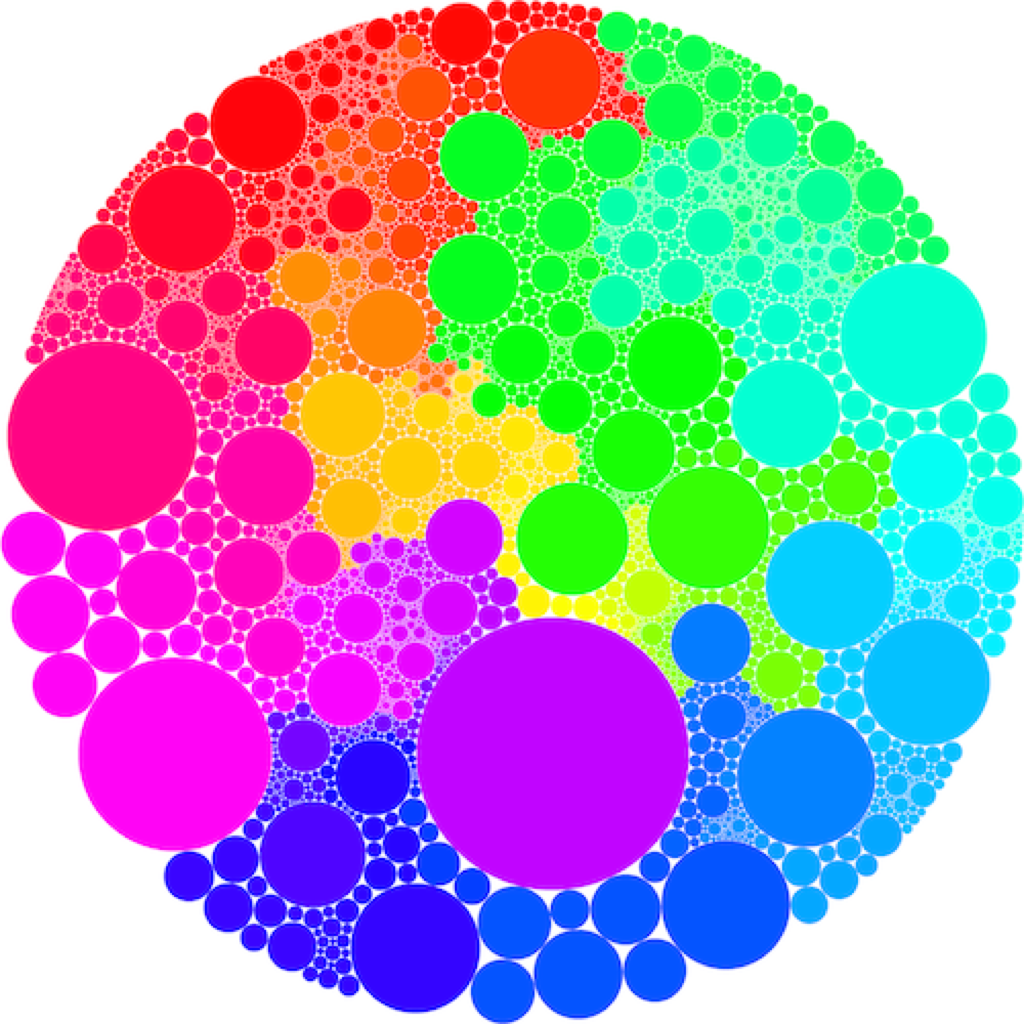}
\end{center}

\vspace{1cm}

{\noindent \footnotesize A book based on this material will be published by Cambridge University Press. This version is free to view and download for personal use only. Not for re-distribution, re-sale or use in derivative works.} 

{\noindent \footnotesize \copyright Nathanaël Berestycki and Ellen Powell, 2024.} 

%\vfill \centerline{[Draft Lecture Notes]}

\newpage

\

\vspace{10cm}

\centerline{\emph{For our families}}

\newpage
\tableofcontents

\newpage

\section*{Preface}

Over fourty years ago, the physicist Polyakov \cite{polyakovstrings} proposed a bold framework for string theory, in which the problem was reduced to the study of certain ``random surfaces''. He further made the tantalising suggestion that this theory could be explicitly solved. 
%Although the mathematical meaning of his proposal was not immediately clear, subsequent attempts by him and other physicists to compute the correlation functions and related quantities led to some fascinating predictions.

Recent breakthroughs from the last fifteen years such as (among many other works) \cite{DKRV, DOZZ,DuplantierSheffield, DuplantierMillerSheffield, Mliouville1, HS19, LeGall, Miermont, GM, DDDF, FKGM, GwynneMillerSAW, BGKRV} 
have not only given a concrete mathematical basis for this theory but also verified some of its most striking predictions -- as well as Polyakov's original vision. This theory, now known in the mathematics literature either as
 \textbf{Liouville quantum gravity (LQG)} or \textbf{Liouville conformal field theory (CFT)}, is based on a remarkable combination of ideas coming from different fields, above all probability and geometry. At its heart is the planar {Gaussian free field} (GFF) $h$, a random distribution on a given reference surface or domain of $\mathbb{R}^2$ whose covariance involves the Green function. A key role is played by the family of measures $\cM^\gamma$ (sometimes referred to as Liouville measures, although this should not be confused with the notion of Liouville measure arising for instance in Hamiltonian dynamics) defined formally as  $\cM^\gamma(\dd x) = \exp ( \gamma h (x)) \dd x$, for a parameter $\gamma$ known as the coupling constant.

\medskip
 This book is intended to be an introduction to these developments assuming as few prerequisites as possible.
 Our starting point is a self-contained and thorough introduction to the two-dimensional continuum \textbf{Gaussian free field} (GFF).
Although surveys and overviews of this object have been written before (notably \cite{Sheffield, WWnotes}), which give plenty of context, both historical and in relation to other topics, the presentation here gives a comprehensive and systematic treatment of some of the analytic subtleties that arise. Many of the details given here for the construction and basic properties of the GFF have perhaps surprisingly not appeared anywhere else before, to the best of our knowledge. 

\medskip  The second basic ingredient and main building block for subsequent chapters is the theory of \textbf{Gaussian multiplicative chaos}. Historically, this theory was first proposed by H{\o}egh-Krohn in \cite{hoeghkrohn} with motivations from constructive quantum field theory not too dissimilar from the ones of this book. 
 In the mathematical literature however it was Kahane, in his seminal contribution \cite{Kahane85}, who introduced it, independently of (and going considerably beyond) \cite{hoeghkrohn}. 
 Kahane was for his part initially motivated by the description of turbulence. In addition to these two distinct motivations, the theory has since found numerous applications in seemingly unrelated areas, such as random matrices, number theory, mathematical finance and planar Brownian motion. 
A useful and early survey of this theory was written in \cite{RhodesVargas_survey} which sketched some of the arguments of the best results available at the time, and also outlined some of these applications. However the state of the art has evolved considerably since then; as a result ours is probably the first unified, systematic and self-contained presentation of this theory.

\medskip  With these tools in hand, the second part of our book is devoted to an exposition of some aspects of Liouville quantum gravity as well as Liouville conformal field theory. These two topics are closely related to one another and they describe, roughly speaking, the same physical theory but with somewhat different perspectives. Essentially, we use the label ``Liouville quantum gravity'' for a random geometric approach highlighting connections with Schramm--Loewner Evolution (SLE). By contrast, we use the label ``Liouville conformal field theory'' for an approach based on the path integral formulation.  
We cover topics such as correlation functions and the so-called Seiberg bounds, Weyl anomaly formula, quantum cones and wedges, quantum zipper, and mating of trees, as well as discrete counterparts to this theory in the form of random planar maps decorated by a model of statistical mechanics (namely, the self-dual Fortuin--Kasteleyn percolation model). 

\medskip These developments require us to work with variants of the (Dirichlet) GFF, respectively the GFF on a Riemannian surface and GFF with Neumann boundary conditions, to which we also provide a systematic introduction. In fact, to the best of our knowledge this is the first place where the analytic details of their construction are given in full.

 %We have chosen to use these two labels when such a distinction appeared useful to us. 

\begin{figure}
\begin{center}
\includegraphics[width=.8\textwidth]{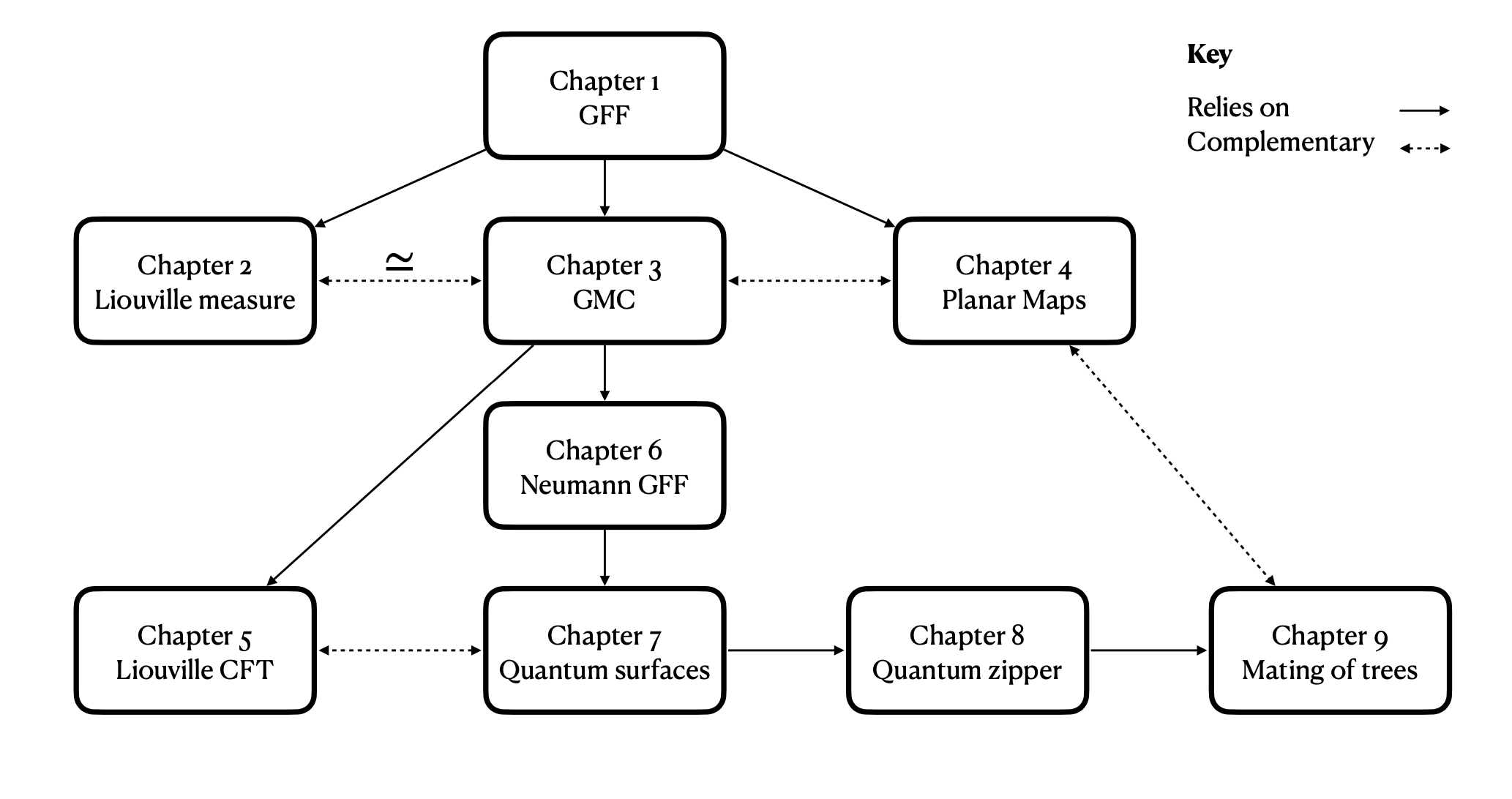}
\caption{Guide to reading. A solid arrow from Chapter $m$ to Chapter $n$ indicates that $m$ is a preqrequisite for $n$. A dashed arrow indicates a complementary perspective on similar/related topics.  The $\simeq$ symbol indicates that Chapter 2 and Chapter 3 are somewhat parallel, with Chapter 2 being focused solely on the construction of the Gaussian Multiplicative Chaos (GMC) measure associated to the (Dirichlet) Gaussian free field, while Chapter 3 gives an exposition of the general theory of GMC.}
\label{F:guide}
\end{center}
\end{figure}

\medskip More specifically, the topics covered include:

{\begin{itemize} \item \textbf{Chapter 1:} the definition and main properties of the GFF with Dirichlet (or zero) boundary conditions;

		\item \textbf{Chapter 2:} the construction of the Liouville measure (in the GFF case), its non degeneracy and change of coordinate formula;
		
		\item \textbf{Chapter 3:} a comprehensive exposition of the construction and properties of general Gaussian multiplicative chaos measures;
		
		\item \textbf{Chapter 4}: an introduction to statistical mechanics on random planar maps -- the discrete counterparts of Liouville quantum gravity -- and Sheffield's bijection with pairs of {trees} \cite{Sheffield_burger};
		
		\item \textbf{Chapter 5}: an introduction to Liouville conformal field theory, as developed in a series of papers starting with \cite{DKRV} by David, Kupiainen, Rhodes and Vargas; 
		%this paper (as well as some related results) forms the main focus of the chapter; 
		
		\item \textbf{Chapter 6}: the definition, construction and main properties of the GFF with Neumann boundary conditions;
		
				\item \textbf{Chapter 7:} an account of the notion of quantum surfaces, and the theory of special quantum surfaces  enjoying scale invariance, such as quantum spheres, discs, wedges and cones; and a proof of equivalence with aspects of the theory developed in Chapter 5;
				
		\item \textbf{Chapter 8:} an exposition of Sheffield's quantum zipper theorem (with novel additional details) and its relation with conformal welding \cite{zipper};
		
		\item \textbf{Chapter 9:} an introduction to the powerful mating-of-trees theory of Duplantier, Miller and Sheffield 
		\cite{DuplantierMillerSheffield}. This includes an extensive and partly novel treatment of space-filling and whole-plane SLE.
	\end{itemize}}
The final three topics above are rather technical, and readers are advised that it will be of most use to people who are actively working in this area. See also Figure \ref{F:guide} for a reading guide.  
%(While the arguments in Chapter 8 follow the main ideas of \cite{zipper}, some are new or simplified, and the overall structure of the proof has been rethought. The result is, {we hope}, that some of the key ideas are more transparent and will be helpful to others.)

\medskip The theory is in full blossom and attempting to make a complete survey of the field would be hopeless, so quickly is it developing. Nevertheless, as the theory grows in complexity and applicability, it has appeared useful to summarise some of its basic and foundational aspects in one place, especially since complete proofs of some facts can be spread over a multitude of papers.

\medskip Clearly, the main drawback of this approach is that many of the important subsequent developments and alternative points of view are not included. {For instance: the expansive body of work on random planar maps and their rigorous connections with Liouville quantum gravity, the Brownian map, Liouville Brownian motion, imaginary geometry, imaginary chaos, and the Liouville quantum gravity metric, do not feature in this book.} For all this we apologise in advance. %{Having said that, a future version of these notes will also touch upon the critical regime for Gaussian multiplicative chaos, the rigorous construction of Liouville conformal field theory, the peanosphere or "mating of trees" description of Liouville quantum gravity surfaces, and Liouville Brownian motion.}
\subsection*{Acknowledgements}  An {initial draft} was written {by the first-named author}, in preparation for the LMS / Clay institute research school on \emph{Modern Developments in Probability} taking place in Oxford, July 2015. The draft was subsequently revised on the occasion of several courses given on this material: at the Spring School on \emph{Geometric Models in Probability} in Darmstadt,  then in July 2016 for the Probability Summer School at Northwestern (for which the chapter on statistical physics on random planar maps was added), in December 2017 for the Lectures on Probability and Statistics (LPS) at ISI Kolkata, and in Berlin at the Stochastic Analysis in Interaction summer school, in 2023. The second-named author lectured on parts of this material in Helsinki in 2022, Santiago de Chile in 2023, and Guanajuato (CIMAT) in 2023.

In all cases {we} thank the organisers (Christina Goldschmidt and Dmitry Beliaev; Volker Betz and Matthias Meiners; Antonio Auffinger and Elton Hsu; Arijit Chakrabarty, Manjunath Krishnapur, Parthanil Roy and especially Rajat Subhra Hazra; Peter Friz and Peter Bank;
Eero Saksman and Eveliina Peltola; Avelio Sep\'ulveda and Daniel Remenik; and Daniel Kious, Andreas Kyprianou, Sandra Palau and Juan Carlos Pardo)
 for their invitations and superb organisation.  Thanks also to Benoit Laslier for agreeing to run the exercise sessions accompanying the lectures at the initial school.

\medskip The Isaac Newton institute's semester on \emph{Random Geometry} in 2015 {was another important influence and motivation for this book}, and we would like to thank the INI for its hospitality. {In fact, this semester served as the second author's initiation into the world of the Gaussian free field and Liouville quantum gravity. The resulting years of discussions between us has led to the present expanded and revised version.}
{We would like to thank many of the} INI programme participants for enlightening discussions related to aspects of the book; especially, Omer Angel, Juhan Aru, St\'ephane Beno\^it, Bertrand Duplantier,  Ewain Gwynne, Nina Holden, Henry Jackson, Benoit Laslier, Jason Miller, James Norris, Gourab Ray, Scott Sheffield, Xin Sun, Wendelin Werner and Ofer Zeitouni. Special thanks in particular to Juhan Aru, Ewain Gwynne, Nina Holden and Xin Sun for many inspiring discussions over the years over a broad range of topics.
%\medskip I am hugely indebted to Juhan Aru and Xin Sun for numerous discussions throughout the semester at INI and while I was writing these notes, be it about the ``big picture'' or about technical details.

\medskip {We} would also like to thank the participants of two reading groups at the University of Bonn and ETH Z\"{u}rich/University of Z\"{u}rich respectively (particularly the organisers Nina Holden and Eveliina Peltola) which followed these notes, and led to many helpful comments. We also received important feedback following graduate courses based on this material which took place at MIT, University of Washington and University of Vienna. Participants and organisers (Scott Sheffield and Zhen-Qing Chen) are warmly acknowledged. We would particularly like to thank Scott Sheffield for his constant encouragement throughout. 

\medskip A substantial part of the writing took place while the authors were invited participants to the semester on \emph{The Analysis and Geometry of Random Spaces} which took place at MSRI in 2022 (now Simons--Laufer Mathematical Sciences Institute) in Berkeley, California. We are very grateful to the organisers of the programme (Mario Bonk, Steffen Rohde, Joan Lind, Eero Saksman, Fredrik Viklund, Jang-Mei Wu) for this amazing opportunity. We also thank the many other participants of this programme for the pleasant atmosphere and the many comments we received while working on Chapter 5 of this book in connection with the reading group on Liouville CFT.

\medskip Finally, comments on versions of this draft have been received at various stages from Morris Ang, Juhan Aru, {Jacopo Borga},  Zhen-Qing Chen, William Da Silva, {Nina Holden}, Henry Jackson, Jakob Klein, Aleksandra Korzhenkova, Benoit Laslier, Joona Oikarinen, L\'{e}onie Papon, {Eveliina Peltola}, Gourab Ray, Mark Sellke, Huy Tran, Joonas Turunen, Fredrik Viklund, Mo-Dick Wong, Henrik Weber and Dapeng Zhan. {We are} grateful for their input which helped to correct minor problems, as well as to emphasise some of the subtle aspects of the arguments; as a result this text is much better than it would have been otherwise. Of course we hasten to add that we retain the entire responsibility for any remaining typo, error, omission, or lack of clarity.

\medskip   This book has been greatly enhanced by a number of beautiful simulations reproduced here with the permission of their authors. These are, respectively:  Jason Miller (Figure \ref{F:CP}, which can also be seen on the cover, and Figure \ref{F:matingtrees});  J\'er\'emie Bettinelli and Benoit Laslier (Figure \ref{F:maps}), Henry Jackson (Figure \ref{F:RSLE4}) and Oskar-Laurin Koiner (Figure \ref{F:dGFF} and Figure \ref{F:Markov}). We thank them wholeheartedly.

\medskip The work of the first author was supported during various stages of the writing by EPSRC (via grants EP/I03372X/1 and EP/L018896/1) and the FWF (via grants 10.55776/P33083 and 10.55776/F1002), while the second author has been supported by the SNF grant 175505 and the UKRI Future Leader’s Fellowship MR/W008513/1. This support is gratefully acknowledged.

\vspace{1cm}
\hfill Nathana\"el Berestycki

\hfill Ellen Powell

\hfill \emph{Vienna and Durham, February 2024}

\newpage

\section{Definition and properties of the GFF}\label{S:GFF}

% !TEX root = master.tex

\subsection{Discrete case}\label{SS:DGFF}

%\subsubsection{Green function; definition}

\dis{The discrete case is included here only for the purpose of guiding intuition when we come to work in the continuum.}

Consider a finite, weighted, undirected graph $\cG = (V,E)$ (with weights $(w_e)_{e \in E}$ on the edges). For instance, $\cG$ could be a finite portion of the Euclidean lattice $\mathbb{Z}^d$ with weights $w_e \equiv 1$. Let $\partial$ be a distinguished set of vertices, called the boundary of the graph, and set $\hat V = V \setminus \partial$. Let $(X_t)_{t\ge 0}$ be the random walk on $\cG$ in continuous time, meaning that it jumps from $x$ to $y$ at rate $w_{x,y}$, and let $\tau $ be the first time that $X$ hits $\partial$ (which we assume to be finite almost surely for every starting point).%,  if the chain is irreducible).

Write $Q=(q_{x,y})_{x,y \in V}$ for the $Q$-matrix  of $X$. That is, its infinitesimal generator, so that for each $x\in V$, $q_{x,y} = w_{x,y}$ for $y \neq x$ and $q_{x,x} = - \sum_{y \sim x} w_{x,y} < \infty$ where $y\sim x$ means that $x$ and $y$ are connected by an edge in $E$.
%$d(x) = \deg(x)$ be the degree, which is a reversible invariant measure for $X$
 Note that the uniform measure $\pi(x) \equiv 1$ is reversible for $X$. We write $\mathbb{P}_x$ for the law of the random walk started and $x\in V$ and $\mathbb{E}_x$ for the corresponding expectation.

\begin{definition} [Green function] \ind{Green function! Discrete} \label{D:dGreen}
The Green function $G(x,y)$ is defined for any $x,y \in V$ by setting
$$G(x,y) =\E_x \left(\int_{0}^\infty \indic{X_t = y;\tau >t} \dd t\right).$$
\end{definition}
In other words $G(x,y)$ is the expected time that $X$ spends at $y$, when started from $x$, before hitting the boundary. Note that with this definition we have $G(x,y) = G(y,x)$ for all $x, y \in \hat V$, since $\P_x(X_t = y ;\, \tau>t) = \P_y(X_t = x; \,\tau >t)$ by reversibility of $X$ with respect to $\pi$.

An equivalent expression for the Green function when working with the  random walk in discrete time $Y=(Y_n)_{n\ge 0}$ (which jumps from $x$ to $y$ with probability proportional to $w_{x,y}$) is
\begin{equation}\label{eq:discreteGreen}
G(x,y) = \frac1{q_y}\E_x \left( \sum_{n=0}^\infty\indic{Y_n = y;\tau(Y) >n} \right),
\end{equation}
where $q_y = \sum_{y \sim x} w_{x,y} = - q_{y,y}$ and $\tau(Y)$ is the first time that $Y$ hits $\partial D$.
Indeed, $X$ and $Y$ can be coupled in such a way that for each $y\in \hat{V}$ and each visit of $Y$ to $y$, $X$ stays at $y$ for an exponentially distributed time with mean $1/q_y$.

The Green function is a basic ingredient in the definition of the Gaussian free field, so the following elementary properties will be important to us.
%\begin{ppl}

\begin{prop}\label{P:basicDGFF}
Let $\hat Q$ denote the restriction of $Q$ to $\hat V \times \hat V$. Then
\begin{enumerate}

\item  $(-\hat Q)^{-1} (x,y)= G(x,y)  $ for all $x,y \in \hat V$.

\ind{Nonnegative definite}
\ind{Definite|see{Nonnegative definite}}
\item  $G$ is a symmetric and non-negative definite function. That is, one has $$G(x,y) = G(y,x)$$ for all $x,y\in V$, and if  $(\lambda_x)_{x\in V}$ is any vector of length $|V|$, then $$\sum_{x,y \in V} \lambda_x \lambda_y G(x, y) \ge 0.$$ Equivalently, $G$ is symmetric and therefore diagonalisable (when viewed as a matrix), and all of the eigenvalues of $G$ are non-negative. Furthermore, restricted to $\hat V$, $G$ is a positive definite function (that is, its eigenvalues are strictly positive).

\item $G(x, \cdot)$ is discrete harmonic in $\hat V \setminus\{x\}$; more precisely $G$ is the unique function of $x,y \in V$ such that $\hat Q G(x, \cdot) = - \delta_x( \cdot)$ for all $x\in \hat{V}$, and satisfies the ``boundary condition'' $G(x,\cdot) = 0 $ on $\partial$ for all $x \in V$.
\end{enumerate}
\end{prop}
%\end{ppl}
Here,  $\delta_x(\cdot)$ denotes the Dirac function at $x$, namely $\delta_x(\cdot) = 1_{\{ \cdot = x \}}$. We also use the natural notation $Qf(x) = \sum_{y \sim x} q_{xy} (f(y) - f(x))$ for the action of the generator $Q$ on functions. Viewed as an operator in this way, $Q$ is often referred to as the discrete Laplacian in continuous time. (Note that by definition, $Q f(x)$ measures the infinitesimal expected change in $f(X_t)$ if the chain starts at $x$).

\begin{rmk}
The proof below is written in the formalism of continuous time Markov chains, which is a little more natural. However, it can equivalently be written using discrete time Markov chains and the definition of the Green function in \eqref{eq:discreteGreen}. 
\end{rmk}

\begin{proof}
Note that since $\hat{Q}$ is symmetric it is diagonalisable, and that all its eigenvalues are negative (this is true of the infinitesimal generator of any Markov chain in continuous time, and here $\hat Q$ is nothing else but the infinitesimal generator of the Markov chain absorbed at $\partial$). Since the chain is absorbed at $\partial$, 0 is not an eigenvalue and all the eigenvalues of $\hat Q$ are therefore strictly negative.

Furthermore, if
$\hat P^t(x,y) = \P_x(X_t = y, \tau>t)$  then  $\hat P_t$ satisfies the backward Kolmogorov equation, namely
$$(d /dt)\hat P^t(x,y) =  \hat Q \hat P^t(x,y),$$
so that $\hat P^t(x,y) = e^{\hat Q t}(x,y)=1+\sum_{j\ge 1}\tfrac1{j!}(\hat Q)^j(x,y)$. It then follows, by Fubini, that
 \begin{align}
G(x,y) &= \E_x (\int_{0}^\infty \indic{X_t = y;\tau >t} \dd t) \nonumber \\
 & = \int_0^\infty \hat P^t(x,y) \dd t  \nonumber \\ &= \int_0^\infty e^{\hat Q t}(x,y)\dd t \nonumber \\
  & = (-\hat Q)^{-1}(x,y). \label{Greeninverse}
 \end{align}
The justification for the last equality comes from thinking about the action of the operator $\int_0^\infty e^{\hat{Q}t}\dd t$ on a single eigenfunction of $\hat{Q}$ (recalling that the corresponding eigenvalue is negative). Since there is a basis of eigenfunctions of $\hat Q$ by symmetry, this suffices to prove the last equality.

\medskip  For the second point, we have already mentioned that $G(x,y) = G(y,x)$. Since $G(x,y) =0$ whenever $y \in \partial$ it suffices to show that the restriction of $G$ to $\hat V$ is positive definite. For this, we can use again that all the eigenvalues of $-\hat Q$, and hence of $(-\hat{Q})^{-1}$ are positive.
  This gives that $G$ is positive definite when restricted to $\hat V$, by \eqref{Greeninverse}. %\ellen{Do we need to add symmetry, and also consider note that $G(x,y)=0$ for $x,y\in \partial$?}

\medskip  Let us finally check the third point. This can be seen as a straightforward consequence of the first point, but we prefer to also include a probabilistic proof which based on the Markov property; effectively, we decompose according to the first jump of the chain.\footnote{The analogous derivation of the same fact using the discrete time chain $Y$ instead of the continuous time chain $X$ is in fact slightly simpler -- we recommend this as an exercise for the reader!} Let $L(x) = \int_0^\infty \indic{X_t = x; \tau > t}\dd t$. Suppose that $y \neq x$ and $t \ge 0$. If $X_0=y$ and $J$ is the first time that $X$ jumps away from $y$ (so $J$ is an exponential random variable with rate $q_y=-q_{y,y}$), we can decompose $\mathbb{E}_y(L(x))$ according to whether $\{ J > t\}$ or $\{ J = s \text{ for some } 0 \le s \le t\}$). Also applying the Markov property at time $J$, we obtain that
  \begin{align*}
    G(x,y)  & = G(y,x) = \E_y (L(x))\\
    & = \E_y ( L(x) | J > t) \P_y(J>t) + \int_0^t\sum_{z \neq y} \P_y( J \in \dd s, X_J =z) \E_y ( L(x) | J =s, X_J = z)\\
    & = G(y,x) e^{-q_y t} + \int_0^t q_y e^{ - q_y s}\dd s \sum_{z \neq y} \frac{q_{y,z}}{q_y} \E_z( L(x)).
  \end{align*}
  Taking the time derivative on both sides at $t=0$ and again invoking symmetry, we arrive at the equality
  $$
  0 = -q_y G(y,x) +  \sum_{z \neq y} q_{y,z} G(z,x) = - q_y G(x,y)+ \sum_{z \neq y} q_{y,z} G(x,z).
  $$
   This means (for fixed $x$, viewing $G(x,y)$ as a function $g(y)$ of $y$ only) that $Q g(y) = \sum_z q_{y,z} g(z) = 0$. Hence $G(x, \cdot)$ is harmonic in $\hat V \setminus\{x\}$.

  When $y = x$, a similar argument can be made, but now the event $\{J>t\}$ contributes to $L(x)$, namely:
  \begin{align*}
    G(x,x) & = \P( J>t) ( t + G(x,x)) + \int_0^t q_x e^{- q_x s} \dd s \sum_{z \neq x} \frac{q_{x,z}}{q_x} \E_z( L(x))\\
    & = e^{ - q_x t} ( t + G(x,x)) + \int_0^t e^{- q_x s} \sum_{z\neq x} q_{x,z} G(x,z) \dd s.
  \end{align*}
  Taking the derivative of both sides at $t =0$ gives
  $$
   0 = -q_x G(x,x) + 1 + \sum_{z\neq x} q_{x,z} G(x,z),
  $$
  and hence
  $$
  \sum_{z}q_{xz} G(x,z) = - 1.
  $$
  The uniqueness comes from the invertibility of $-\hat Q$.
 %
 % by symmetry $G(x,y) = \E_y( \text{time spent at $x$ before hitting $\partial$})$.
  %Decomposing on the position of the walk at time $h$ and letting $h \to 0$, we get
  %$G(x,y) = \sum_{y' \sim y} q(y, y')h  G(x, y') + G(x,y) ( 1 - q_x h) + o(h) $ where $q_x = |q_{x,x}|$.
   %After cancelling $G(x,y)$ on both sides, dividing by $h$ and letting $h \to 0$ this translates into $Q_y G(x, y) = 0$ (here we write $Q_y$ to denote $Q$ acting on the space variable $y$). So $G(x, \cdot)$ is harmonic in $\hat V \setminus \{x\}$ as desired. When $y = x$ we find $Q_y G(x,y) = 1$ to account for the time spent at $x$ initially.
\end{proof}

\begin{rmk}
  An alternative proof of the first point (that is, of \eqref{Greeninverse}) uses the transition matrix $\hat R^n(x,y) =\P_x(Y_n = y , \tau(Y)>n)$ of the jump chain. Indeed, we have already noted that
  \begin{align*}
  G(x,y) &= \frac1{q_y} \E_x( \sum_{n=0}^\infty \indic{Y_n = y , \tau(Y)>n}) \\
  & = \frac1{q_y} \sum_{n=0}^\infty \hat R^n(x,y) \\
  & = \frac1{q_y} (I - \hat R)^{-1}(x,y)\\
  & = ( - \hat Q)^{-1}(x,y)
  \end{align*}
where in jumping from the second to the third line we used the fact that $\sum_{n=0}^\infty \hat R^n   = (I - \hat R)^{-1}$, an identity valid for any matrix of spectral radius (that is, largest eigenvalue modulus) strictly smaller than one, which is the case here. An alternative proof that $G$ is non-negative definite can be obtained using same argument in the proof of \cref{L:posdef} (this is stated in the continuous case but can also easily be adapted to this discrete setting).
  \end{rmk}

\begin{definition} [Discrete GFF] The (zero boundary) discrete Gaussian free field on $\cG =(V,E)$ is the centred Gaussian vector $(h(x))_{x\in V}$ with covariance given by the Green function $G$.
\end{definition}

\begin{rmk}
This definition is justified. Indeed, suppose that $(C(x,y))_{x,y\in V}$ is a given function. Then there exists a centred Gaussian vector $X$ having covariance matrix $C$ if and only if $C$ is symmetric and non-negative definite  (in the sense of property 2 above).
\end{rmk}

Note that if $x \in \partial$, then $G(x,y) = 0$ for all $y\in V$ and hence $h(x) = 0$ almost surely.

\begin{figure}\begin{center}
\includegraphics[width=.5\textwidth]{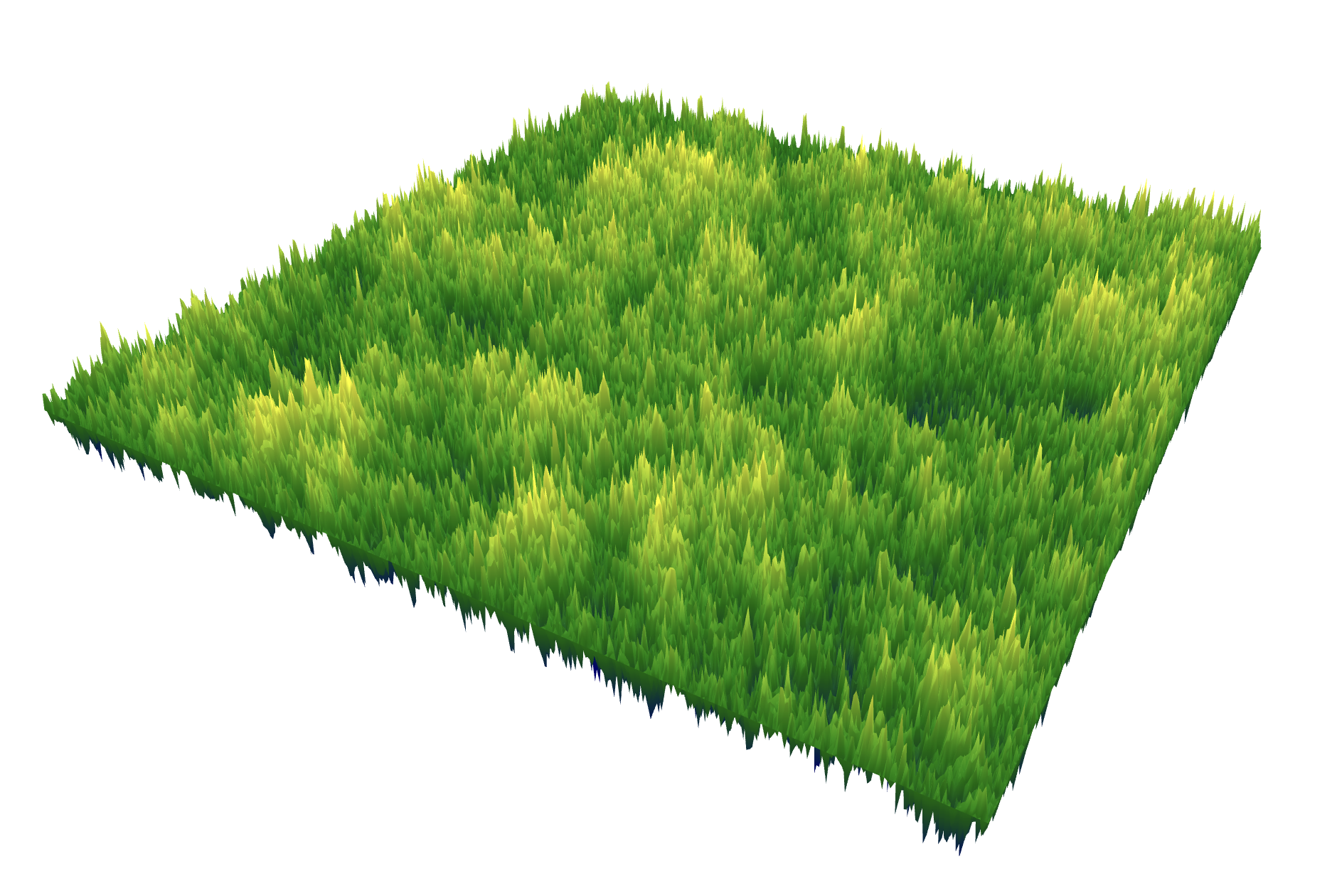}
\end{center}
\caption{A discrete Gaussian free field. Simulation by Oskar-Laurin Koiner.}
\label{F:dGFF}
\end{figure}

\medskip In fact it is possible to provide a concrete construction of the discrete Gaussian free field, in terms of i.i.d. standard Gaussian random variables. This construction has the advantage that it is very easy to implement on a computer to produce simulations, such as the one in Figure \ref{F:dGFF}. We first introduce some notations. Set $N = |\hat V|$ and consider the space of functions $f: \hat V \to \R$, equipped with the inner product
\begin{equation}\label{l2discrete}
(f,g) = \sum_{x\in \hat V} f(x) g(x).
\end{equation}
For this reason (and even though $\hat V$ is finite) we denote this space of functions by $\ell^2(\hat V)$. Any function in $\ell^2(\hat V)$ can canonically be extended to a function on $V$ by setting it to zero on $\partial$. Recall that a function $f\in \ell^2(\hat V)$ is an eigenfunction of $-\hat Q$ with eigenvalue $\lambda$ (necessarily positive) if for all $x \in\hat V$, $-\hat Q f(x) = \lambda f(x)$, that is, 
$$
-\sum_{y \in \hat V} q_{x,y} f(y) = \lambda f(x). 
$$ 
As already mentioned in the proof of Proposition \ref{P:basicDGFF}, since $- \hat Q$ is symmetric, it is diagonalisable in an orthonormal basis of $\ell^2(\hat V)$. Let $e_1, \ldots, e_N$ denote the orthonormal eigenfunctions and let $0< \lambda_1 \le \ldots \le \lambda_N$ denote the corresponding eigenvalues (with multiplicities).

\begin{theorem}
	\label{T:discretefourier}
	Let $(e_m)_{m=1}^{N}$ and $(\lambda_m)_{m=1}^N$ be as above. Then for $x,y \in \hat V$ we have the expansion
	\begin{equation}\label{eq:Gdiscrete}
	G(x,y) = \sum_{m=1}^N \frac1{\lambda_m} e_m(x) e_m(y).
	\end{equation}
This also extends to $\partial$ if we extend $e_m$ by zero on $\partial$. 
	
	Furthermore, let $(X_m)_{m=1}^{N}$ be a sequence of i.i.d.\, standard Gaussians. Set
	%Then one can check (we leave this as Exercise \ref{ex:discretefourier} at the end of this chapter) that $(h(x))_{x\in V}$ defined by 
	\begin{equation}\label{eq:GFFexpansiondiscrete}
	h(x):=\sum_{m=1}^{N} \frac{X_m}{\sqrt{\lambda_m}} e_m(x)  ; \quad x \in  V
	\end{equation}
	Then $h$ is a discrete GFF on $\cG$.
\end{theorem}

\begin{proof}
Since $- \hat Q$ is invertible it suffices to check that for $x\in \hat V$, if we define $g(y) = g_x(y) $ as on the right hand side of \eqref{eq:Gdiscrete}, namely $g(y) = \sum_{m=1}^N (1/\lambda_m) e_m(x) e_m(y)$ viewed as a function of $y \in \hat V$, 
 one has 
\begin{equation}\label{E:Qg}
- \hat Q g = \delta_x.
\end{equation}
By linearity and since $e_m$ is an eigenfunction of $-\hat Q$, we see that (recall that $x \in \hat V$ is fixed)
\begin{align*}
- \hat Q g &= \sum_{m=1}^N \frac1{\lambda_m} e_m (x)  (- \hat Q e_m) \\
 &= \sum_{m=1}^N  e_m (x)  e_m
\end{align*}
On the other hand, expanding $\delta_x$ in the basis $(e_m)_{m=1}^N$ we also find that
$$
\delta_x = \sum_{m=1}^N ( \delta_x, e_m) e_m = \sum_{m=1}^N e_m(x) e_m, 
$$
which is indeed the same as the right hand side of the previous equation. This proves \eqref{E:Qg} and thus \eqref{eq:Gdiscrete}.

Turning to \eqref{eq:GFFexpansiondiscrete}, we simply note that $(h(x))_{x\in \hat V}$ is clearly a centred Gaussian vector, whose covariance is given by 
\begin{align*}
\E[ h(x) h(y) ] & =\E \left(  \sum_{m,m'=1}^N \frac{X_m X_{m'}}{\sqrt{\lambda_m \lambda_{m'}} } e_m(x) e_{m'} (y) \right)\\
& = \sum_{m=1}^N \frac1{\lambda_m} e_m(x) e_m(y)  = G(x,y) 
\end{align*}
by \eqref{eq:Gdiscrete}, as desired. 
\end{proof}

%\subsubsection{Dirichlet energy; Markov property}

\medskip Usually for Gaussian fields, looking at the covariance structure is the most useful way of gaining intuition. However in this case, the joint probability density function of the $|V|$ components of $h$ is perhaps more illuminating.

\begin{theorem}[Law of the GFF and Dirichlet energy]  \label{T:dGFF_Dirichlet}\ind{Dirichlet energy}The law of the discrete GFF is absolutely continuous with respect to Lebesgue measure on $\R^{\hat V}$, with joint density 
	%with respect to Lebesgue measure on $\R^{\hat V}$ is 
	proportional to $$\exp\left( -\frac14 \sum_{x,y \in V }q_{x,y}(h(x) - h(y))^2 \right)$$
	at any point $(h(x))_{x\in V}$ {with $h(x)=0$ for $x\in \partial$, viewed as a fixed element of $\R^{\hat V}$.} (Note that the sum includes the vertices $v \in \partial$.)
%$d\mathbf{x} = \prod_{1\le i \le n} dx_i$. More specifically, we have %the joint pdf is proportional to
%$$
%\P( h(\mathbf{x}) ) = \frac1{Z} \exp\left( -\frac14 \sum_{x,y \in V }q_{x,y}(h(x) - h(y))^2 \right) dh(\mathbf{x})
%$$
%where $Z$ is a normalising constant (called the partition function).
\label{P:gfflaw}
\end{theorem}

\begin{rmk}
The previous formula might seem a little confusing at first, since we are using $(h(x))_{x\in \hat V}$ both to denote the random vector consisting of the values of the discrete Gaussian free field, and for a fixed (deterministic) element in $\R^{\hat V}$ at which we evaluate the density of this random vector. 
To avoid any confusion, the formula above means the following: if we write
%$|V| = n$,
%= \{v_1, \ldots, v_n\}$
%and
$Y_v $ for the random variable $Y_v: = h(v)$, then
$$\P( (Y_v)_{v \in V} \in A) = \int_A  \frac1{Z} \exp( - \frac14 \sum_{v,w \in V} q_{v,w}(x_v - x_w)^2 ) \prod_{v \in \hat V} \dd x_v $$
where $Z=\int_{ \R^N} \exp(-\frac{1}{4}\sum_{v,w\in V}q_{v,w}(x_v-x_w)^2)\prod_{v \in \hat V}  \dd x_v$, where $N = |\hat V|$. This holds for any Borel set $A$ contained in the hyperplane $\{(x_v)_{v \in V}: x_v=0 \text{ for all } v \in \partial\}$ of $\R^{V}$.%such that whenever $\mathbf{y} \in A$ then $y_i = 0$ for indices $i$ corresponding to $x_i \in \partial$.
\end{rmk}

For a given function $f:V \to \R$, the quantity 
\begin{equation}\label{E:Dirichletdiscrete}
\cE(f,f) := \frac12\sum_{x, y\in V} q_{x,y}(f(x) - f(y))^2
\end{equation}
is known as the \textbf{Dirichlet energy} of $f$, and is a discrete analogue of $(1/2)\int_D |\nabla f|^2$.

\begin{proof}[Proof of Theorem \ref{T:dGFF_Dirichlet}]
The result follows from the fact that for a centred Gaussian vector $(Y_1, \ldots, Y_N)$ with invertible covariance matrix $\Sigma$, the joint probability density function on $\R^N$ is proportional to
$$f(x_1, \ldots, x_N) =  \exp( -  \frac12 x^T \Sigma^{-1} x) .$$

For us, the vertices $v \in \hat V$ play the roles of the indices $1\le i \le N$ above with $N = |\hat V|$, and the values $h(v)$ for $v\in V$ play the roles of the $x_i$ (to get a non-degenerate covariance matrix we restrict ourselves to vertices in $ \hat V$, in which case $G$ is invertible by \cref{P:basicDGFF}). Note that since we are only considering $h$ with $h(v)=0$ for $v\in \partial$, it suffices to show that
\begin{equation*}-\frac{1}{2} h(\mathbf{\hat v})^T G^{-1}h(\mathbf{\hat v})=-\frac{1}{4}\sum_{x,y\in V}q_{x,y}(h(x)-h(y))^2, \;\text{ for } \; h(\mathbf{ \hat v})=(h(v))_{v\in \hat{V}}. \end{equation*}
Recall that $(-\hat{Q})^{-1}(x,y) = G(x,y) $ for $x,y \in \hat V$, so that $G^{-1} (x,y) = -q_{xy}$.
%Therefore, inverting this relation and multiplying by $D$, $G^{-1} = D (I- \hat P)$, or $G^{-1}(x,y) = d(x) (I - \hat P)(x,y)$ for $x,y \in \hat V$.
%Recall that $(I-\hat P)^{-1}(x,y) = G(x,y) d(y)$ for $x,y \in \hat V$. In matrix notations, $(I- \hat P)^{-1} = G D$ where $D$ is the diagonal matrix (of dimension $|\hat V|$) with entries $d(x)$. Therefore, inverting this relation and multiplying by $D$, $G^{-1} = D (I- \hat P)$, or $G^{-1}(x,y) = d(x) (I - \hat P)(x,y)$ for $x,y \in \hat V$.
%\red{Factor 2?}
Hence
\begin{equation*}
h(\mathbf{\hat v})^T G^{-1} h(\mathbf{\hat v})  = \sum_{x,y\in \hat V} G^{-1}(x,y) h(x) h(y)=\sum_{x,y\in \hat V} -q_{x,y} h(x) h(y).
\end{equation*}
Moreover, as we only consider $h$ with $h(x)=0$ for $x\in \partial$, this can be rewritten as
\begin{equation*}  -\sum_{x,y \in V} q_{x,y} h(x) h(y)
=\frac{1}{2} \sum_{x,y\in V} q_{x,y}
 (h(x) - h(y))^2 - \frac{1}{2}\sum_{x,y \in V} h(x)^2 q_{x,y} -  \frac{1}{2} \sum_{x,y\in V} h(y)^2 q_{x,y}, \end{equation*}
where since $\sum_{y \in V} q_{x,y} = 0$ and $q_{x,y} = q_{y,x}$ for all $x,y$, the terms
$$
\sum_{x,y\in V} h(x)^2 q_{x,y}  \text{ and } \sum_{x,y\in V} h(y)^2 q_{x,y}$$
are both equal to $0$.
Note that in this final line of reasoning  it is important to sum over all of $V$ and not just $\hat V$. Thus
\begin{equation*}
-\frac{1}{2}h(\mathbf{\hat v})^T G^{-1}h(\mathbf{\hat v})=-\frac{1}{2} \times \frac{1}{2} \sum_{x,y\in V} q_{x,y}(h(x)-h(y))^2,
\end{equation*}  as required.
\end{proof}

%\begin{rmk}
%The trick of exploiting symmetry and changing variables allowed us to write $$\langle Qf, f\rangle = - \cE(f,f)$$ and so can be thought of as a discrete version of the integration by parts formula or Gauss--Green formula\index{Gauss--Green formula} (see Lemma \ref{L:ipp}).
%\end{rmk}

Notice that the Dirichlet energy of functions is minimised by harmonic functions. This means that the Gaussian free field can be viewed as a ``Gaussian perturbation of a harmonic function'': as much as possible, it ``tries'' to be harmonic. In fact, this is a little ironic, given that in the continuum it is not even a function (see the next section).
%\begin{ex}Show this.
%\end{ex}

This heuristic is at the heart of the Markov property, which is without a doubt the most useful property of the GFF. We state it here without proof, as we will soon prove its (very similar) continuum counterpart. \ind{Markov property}

\begin{theorem}\label{T:markovdiscrete}
[Markov property of the discrete GFF] Fix $U \subset V$. The discrete GFF $h=(h(x))_{x\in V}$ can be decomposed as
\[h = h_0 + \varphi,\]
where $h_0$ is Gaussian free field on $U$ and $\varphi$ is harmonic in $U$. Moreover, $h_0$ and $\varphi$ are independent.
\end{theorem}
By a Gaussian free field in $U$ we mean the GFF on the graph $(V,E)$ but now with $\partial = V \setminus U$, in particular $h_0 = 0$ outside of $U$.

In other words, this theorem says that conditionally on the values of $h$ outside of $U$, the field can be written as the sum of two independent terms. One of these is a zero boundary GFF in $U$, and the other is just the harmonic extension into $U$ of the values of $h$ outside $U$. To see this, note that the information about the values of $h$ outside of $U$ is completely contained in $\varphi$, since $h_0$ is zero outside of $U$. Thus conditioning on the values of $h$ outside of $U$ is the same as conditioning on $\varphi$. Since $h_0$ is independent of $\varphi$, the conditional law of $h$ given $\varphi$ is as described.

\subsection{Continuous Green function}
\label{SS:Greencontinuous}
%\begin{ppl}
We will follow a route that is similar to the previous discrete case. First we need to recall the definition of the Green function. %A \textbf{domain} of $\R^d$ is an open, connected set. 
We will only cover the basics here, and readers who want to know more are advised to consult, for instance, Lawler's book \cite{Lawlerbook} which reviews important facts in a very accessible way. The presentation here will be somewhat different.

Let $d\ge 1$. Let $p_t(x,y)$\indN{{\bf Brownian motion}! $p_t(x,y)$; transition probability for speed two Brownian motion on $\R^d$} denote the transition probability of a Brownian motion ${B}$ in $\R^d$ with ``speed'' two (that is, ${B}_t = (X^1_{2t}, \ldots, X^d_{2t})$ for $t \ge 0$, where $X^1, \ldots, X^d$ are independent standard Brownian motions\footnote{This choice ensures that the infinitesimal generator of ${B}$ is the Laplace operator $\Delta$ instead of $\Delta/2$.} in $\R$). Then \ind{Domain}
  \begin{equation}\label{E:transition_wholespace}
 p_t(x,y) = (4\pi t)^{- d/2} \exp ( - | x- y|^2/(4t)),
 \end{equation}
 which by the Markov property is also the density, with respect to Lebesgue measure on $\R^d$, of the law of $B_t$ (when started from $x$).
   For $D\subset \R^d$ an open set, we define $p^D_t(x,y)$\indN{{\bf Brownian motion}! $p^D_t(x,y)$;  transition probability for speed two Brownian motion killed when leaving $D$} to be the transition probability of Brownian motion with speed two, 
 killed when leaving $D$, which is  defined as the density, with respect to Lebesgue measure on $\R^d$, of the law of ${B}_t$, but restricted to the event $\{\tau_D > t\}$, where 
 $$
 \tau_D = \inf\{ t>0: B_t \notin D\}.
 $$ \indN{{\bf Brownian motion}! $\tau_D$; hitting time of $\partial D$}
 In other words, for any Borel set $A$ in $\R^d$, it satisfies
 \begin{equation}\label{eq:densitykilled}
 \P_x( {B}_t \in A, \tau_D > t) = \int_{\R^d} 1_A (y) p_t^D (x,y) \dd y.
 \end{equation}
The (almost everywhere, for a fixed $t\ge 0$) existence of a function satisfying \eqref{eq:densitykilled} follows directly from the Radon--Nikodym derivative theorem, since it is clear that if $A$ has zero Lebesgue measure, then $\P_x({B}_t \in A, \tau_D > t) \le \P_x( {B}_t \in A) = 0$. 

By conditioning on the position at time $t$ of ${B}_t$, it is not hard to check that $p_t^D (x,y)$ can be expressed rather simply in terms of the whole space transition probabilities in \eqref{E:transition_wholespace} and the so-called (speed two) Brownian bridge\footnote{A reader who is unfamiliar with the notion of Brownian bridge may without danger skip to the conclusion immediately following \eqref{E:ptD}.} $({b}_s)_{0\le s \le t}$ of duration $t$ from $x$ to $y$, which describes the law of ${B}$, conditionally given ${B}_0 = x$ and ${B}_t = y$. Namely, if we denote by $\P_{x\to y; t}$ this law, then we see that 
\begin{align*}
\P_x ({B}_t \in A; \tau_D > t) &= \int_{\R^d} \P_{x\to y; t} ( {b}_t \in A ; \tau_D > t ) p_t(x,y) \dd y\\
& = \int_{\R^d} 1_{A} (y) \P_{x\to y ; t} ( \tau_D > t)  p_t(x,y) \dd y. 
\end{align*}
Comparing with \eqref{eq:densitykilled}, we deduce that, for every fixed $t\ge 0$ and almost every $y$, 
\begin{equation}\label{E:ptD} 
p_t^D(x,y) =  \pi_t^D(x,y) p_t(x,y); \text{ where } \pi_t^D(x,y) = \P_{x\to y; t} (\tau_D > t ). 
\end{equation}
 The right hand side is easily seen to be a jointly continuous function in $t>0$ and $x,y \in \bar D$, as this is clearly satisfied by both $\pi_t^D(x,y)$ and $p_t(x,y)$ separately. This defines the transition probability function $p_t^D(x,y)$ of Brownian motion killed when leaving $D$ uniquely. 
%
%p_t^D(x,y) = p_t(x,y) \pi_t^D(x,y) \end{equation} 

%  
%  
 %is the Gaussian transition probability {with covariance matrix $2t I$}, and $\pi_t^D(x,y)$ is the probability that a Brownian bridge ({with speed 2}) from $x$ to $y$ of duration $t$ remains in $D$.

Clearly, by the Markov property of Brownian motion, the transition probabilities satisfy the Chapman--Kolmogorov equation:
\begin{equation}
p_{t+s}^D (x,y) = \int_{\R^d} p_t^D(x,z) p_s^D (z, y) \dd z \quad \text{ for } s,t\ge 0, \, x,y\in D.
\end{equation}

Note also immediately for future reference that, by definition of $p_t^D(x,y)$ and the monotone class theorem, if $\phi$ is any non-negative Borel function and $t\ge 0$, then 
$$
\E_x( \phi( {B}_t) 1_{\{\tau_D > t\} } ) = \int_{\R^d} \phi(y) p_t^D(x,y) \dd y.
$$
Consequently, by Fubini's theorem, 
\begin{align}
\E_x ( \int_0^{\tau_D} \phi({B}_s) \dd s) & = \E_x ( \int_0^\infty \phi(B_s) 1_{\{\tau_D > s\}} \dd s ) \nonumber \\ 
& = \int_0^\infty \int_{\R^d} \phi(y) p_s^D (x,y) \dd y \dd s \nonumber\\
& = \int_{\R^d} \phi(y) (\int_0^\infty p_s^D (x,y) \dd s ) \dd y. \label{eq:preGformula}
\end{align}

The time integral in brackets in \eqref{eq:preGformula} plays a crucial role in this book, and is called the (continuous) Green function. Note the parallel with Definition \ref{D:dGreen}; intuitively, as in the discrete case, the Green function measures the expected amount of time spent ``at'' a point $y$ (that is, near $y$) before exiting $D$. 

\begin{definition}[Continuous Green function] \ind{Green function! Dirichlet boundary conditions} Let $D\subset \R^d$ be an open set. The Green function $G_0(x,y) = G_0^D(x,y)$ \indN{{\bf Green functions}! $G_0^D(\cdot, \cdot)$; for Laplacian with zero boundary conditions in $D$}
	 is defined by
\begin{equation}\label{eq:g0}
G_0(x,y) = \int_0^\infty p_t^D(x,y) \dd t
\end{equation}
for $x\ne y$ in $D$.
\end{definition}

Note in particular that, combining our definition of the Green function with \eqref{eq:preGformula}, we obtain:
\begin{equation}\label{eq:Gformula}
\E_x( \int_0^{\tau_D} \phi(B_s) \dd s) = \int_{\R^d} G_0^D (x,y) \phi(y) \dd y. 
\end{equation}
This agrees with our intuition that the Green function measures the expected amount of time spent by a Brownian motion near a point $y$ before leaving $D$.

%\begin{rmk}\label{R:Ddisconnected}
%	Suppose that $D$ is a disjoint union of connected components $\cup_i D_i$. Then since Brownian motion cannot get from $x\in D_i$ to $y\in D_j$ before time $\tau_D$ when $i\ne j$, we have $G_0(x,y)=0$ in this case. Furthermore, in each $D_i$ we have $G_0^D(x,y)=G_0^{D_i}(x,y)$ whenever $x,y\in D_i$. In other words, the Green function in such a set $D$ is just equal, when restricted to any connected component, to the Green function of that component. Otherwise, it is zero.
%\end{rmk}

\begin{rmk}[Normalisation] We call the attention of the reader to the fact that the normalisation of the Green function is a little arbitrary. We have chosen to normalise it so that $G$, as we will soon see, is the inverse of (minus) the Laplacian, with no multiplicative constant in front. This choice  is consistent with say \cite{WWnotes}. In particular, in two dimensions, our normalisation is chosen so that for $D\subset \C$ simply connected, say, we will have $$G^D_0(x,y)\sim \frac1{2\pi}\log(|x-y|^{-1})$$ as $y\to x$ (see \cref{P:basicGFFc}).
This is however \emph{not} the standard set up for Gaussian multiplicative chaos (see \cref{S:Liouvillemeasure,S:GMC}) or in papers on Liouville quantum gravity, where the Green function is often normalised so that it blows up like $\log (|x-y|^{-1})$ (that is, it differs from our choice by a factor of $2\pi$). This means that the Gaussian free field we are about to define will differ by a factor of $\sqrt{2\pi}$ from the field usually considered in the Gaussian multiplicative chaos literature, and which we will we also switch to in \cref{S:Liouvillemeasure,S:GMC}. While from the point of view of Gaussian multiplicative chaos it is more natural to define the Gaussian free field as a log-correlated field rather than a $(2\pi)^{-1}-$log-correlated field, we have chosen the above normalisation of the Green function for this chapter, since it is more natural from an analytic perspective. In particular, it saves us many tedious powers of $2\pi$ in our subsequent considerations involving Sobolev spaces.

Another commonly used normalisation of the Green function corresponds to the integral of the transition density for Brownian motion with speed 1 rather than speed 2. This Green function differs from ours by a factor of 2, and the resulting Gaussian free field by a factor of $\sqrt{2}$.
\end{rmk}

 We will sometimes drop the notational dependence of $G_0^D$ on $D$ when it is clear from the context.
The subscript 0 refers to the fact that $G$ has \textbf{zero boundary conditions}; equivalently, that $G$ is defined from a Brownian motion killed when leaving $D$.\\

When $d\ge 2$, it is easy to see that $G^D_0(x,x)$ is typically ill defined ($ = \infty$) for all $x\in D$. This is because $\pi_t^D(x,x) \to 1$ as $t \to 0$ and so $({4}\pi t)^{-d/2}\pi_t^D(x,x)$ cannot be integrable. When $x \neq y$ the bound $p_t^D(x,y) \le p_t(x,y)$ suffices to show there can be no problem of integrability near $t=0$. Furthermore, at least if $d \ge 3$, $p_t(x,y)$ is also integrable near $t =\infty$, so that $G_D(x,y) < \infty$. 

The case $d =2$ is more delicate. In fact, it can be shown that in that case, then $G^D_0(x,y) < \infty$ if $x\neq y \in D$  and $D$ has at least one \textbf{regular} boundary point, that is, a point $b\in \partial D$ such that $\P_b ( \tau_D = 0) = 1$ (in other words, starting from a boundary point, a Brownian motion leaves $D$ instantaneously -- equivalently, $b\in \partial D$ is regular if for all $\eps>0$, $\P_z(\tau_D>\eps) \to 0$ as $z \in D$ converges to $b$). This follows for instance from Lemma 2.32 in \cite{Lawlerbook}, which gives a uniform bound of the form 
\begin{equation}\label{eq:dominationpt}
p_t^D(x,y) \le ct^{-1} (\log t)^{-2}
\end{equation}
 uniformly over $x,y\in D$ and $t>1$ for a regular domain $D$. 
% By regular, we mean that Brownian motion starting from a point $x \in \partial D$ will hit $D^c$ instantaneously. 
 \ind{Domain!Regular} 
 Any proper simply connected open set in two dimensions is easily seen to be regular, i.e., every point on the boundary of $D$ is regular. (When the domain is simply connected, an even stronger bound  than \eqref{eq:dominationpt} holds, and is in fact easier to prove, see \cite{SLEnotes}). We also note for further reference that \eqref{eq:dominationpt} implies (through dominated convergence and the continuity of both $\pi_t^D(x,y)$ and that of $p_t(x,y)$) that $(x,y) \in D^2 \mapsto G^D_0(x,y)$ is continuous away from the diagonal 
 $\{(x,y) \in D^2: y=x\}$. (This continuity could also be proved differently from Proposition \ref{P:basicGFFc} and using elliptic regularity arguments.) Although we will not use this, we note that \cite{SLEnotes} proves that the continuity of $y \in D\setminus\{x\} \mapsto G^D_0(x, y)$ extends to regular boundary points: i.e., as $y \to b\in \partial D$ regular, $G^D_0(x,y) \to 0$.  
   
Finally, in dimension $d=1$, we will see that $G^D_0(x,y)$ is actually finite even when $x=y$. In this case $p_t^D(x,y)$ is zero as soon as $x$ or $y$ are in $D^c$ (including on the boundary of $D$), for any $t>0$.

\begin{example}
 Suppose $D= \H\subset \C$ is the upper half plane. Then
it is not hard to see that $p_t^{\H}(x,y) = p_t(x,y) - p_t(x, \bar y)$ by a reflection argument (in fact, by the reflection principle of ordinary one dimensional Brownian motion). Hence one can deduce that
%using the formula $e^{-a/t} - e^{- b /t} = t^{-1} \int_a^b e^{- x  /t} dx $:
\begin{equation}\label{G_H}
G_0^{\H}(x,y) = {\frac1{2\pi}}
\log \left| \frac{x - \bar y}{x-y}\right|
\end{equation}
for $x\ne y$ (see Exercise \ref{Ex:green_half} for a hint on the proof).
\end{example}

In the special case $d=2$, a fundamental property of Brownian motion (also with speed 2) is that it is \textbf{conformally invariant}. That is, suppose that $(B_s)_{s\ge 0}$ is a Brownian motion in the plane with speed 2, and $T$ is an analytic map defined on a simply connected open set $D$ with $T'\ne 0$ on $D$ (at this stage, we do not require $T$ to be one to one). Then $T(B_s)$, considered up until the exit time $\tau_D$ from $D$ by $B$, is a Brownian motion in the image domain $D' = T(D)$, up to a time change. More precisely, if we define 
$$
F(t) = \int_0^{\tau_D \wedge t}  |T'(B_s)|^2 \dd s,
$$
then we can talk about its right continuous inverse (which is also simply its inverse here)
$$
F^{-1} (s) = \inf\{ t>0: F(t) > s\}.
$$
The conformal invariance of Brownian motion states that if we set 
\begin{equation}
B'_s: = T ( B_{F^{-1} (s)} ); \quad \text{ for } 0 \le s < \tau' = F(\tau_D) 
\end{equation}
then $(B'_s)_{0 \le s \le \tau'}$ is another Brownian motion with speed 2, stopped at the time $\tau'$ when it first leaves $D' = T(D)$. This fundamental property, predicted by L\'evy in the 1940s, can be proved relatively easily using the Cauchy--Riemann equations satisfied by $T$ and an application of It\^o calculus (both It\^o's formula and the Dubins--Schwarz theorem).

\medskip Although the conformal invariance of Brownian motion is only up to a time change, and the Green function $G_0$ measures the expected time spent by Brownian motion close to a location before leaving the domain, a remarkable property of the Green function is that it is \emph{completely}  invariant under conformal isomorphisms, in the following sense. \\

We say that $D\subset \R^d$ or $D\subset \C$ is a \textbf{domain} if it is \emph{open} and \emph{connected}.

\begin{prop}[Conformal invariance of the Green function]\label{P:GreenCI}
Let $D,D'\subset \C$ be regular domains. Suppose that $T: D\to D'$ is a conformal isomorphism (that is, analytic with non-vanishing derivative and one to one). Then 
$$G_0^{T(D)} (T(x), T(y)) = G_0^D(x,y).$$
\end{prop}

Note that together with \eqref{G_H} and the Riemann mapping theorem, this allows us to determine $G_0^D$ in any simply connected proper domain $D\subset \C$. 
\ind{Conformal invariance! of Green function}

\begin{proof}
  The proof is a simple application of the change of variable formula. 
  %The Jacobian term $|T'(y)|^2$ arising from the change of variables just cancels the term $|T'(B_s)|^2$ arising from It\^o's formula. %in the conformal invariance of Brownian motion.
  Let $\phi$ be a test function and let $x' = T(x)$. Then, by \eqref{eq:Gformula}, 
  $$
  \int_{D'} G_0^{D'}(x',y') \phi(y') \dd y' = \E_{x'} ( \int_0^{\tau'} \phi(B'_{s}) \dd s )
  $$
  where $B'$ is a Brownian motion and $\tau'$ is its exit time from $D'$. On the other hand, the change of variable formula applied to the left hand side gives us, letting $y' = T(y)$ (a change of variable whose Jacobian derivative evaluates to $\dd y' = |T'(y)|^2 \dd y$):
  \begin{equation}\label{eq:GreenCI1}
  \int_{D'} G_0^{D'}(x',y') \phi(y') \dd y' = \int_{D} G_0^{D'}(T(x), T(y)) \phi(T(y)) |T'(y)|^2 \dd y.
  \end{equation}
  Now let us compute the right hand side of the initial equation in a different way, using the conformal invariance of Brownian motion discussed above. This allows us to write $B'_{s} = T(B_{F^{-1}(s)})$;
  % where $F(t) = \int_0^t |T'(B_s)|^2 \dd s$ for $s \le \tau $, $\tau$ is the first exit time of $D$ by $B$, and $F^{-1}$ is the cadlag inverse of $F$. 
  moreover, in this correspondence one has $\tau' = F^{-1}(\tau_D).$ We apply the change of variable formula, but now to the time parameter $t = F^{-1}(s)$, or (since $F^{-1} $ is the inverse of $F$), $s = F(t)$. 
The Jacobian derivative is thus 
$$\dd s = F'(t) \dd t = | T'(B_t)|^2 \dd t, $$ 
by definition of $F$ and the fundamental theorem of calculus. Thus, 
  \begin{align}
  \E_{x'} ( \int_0^{\tau'} \phi(B'_{s}) \dd s )& = \E_x (\int_0^{F^{-1}(\tau)} \phi(T(B_{F^{-1}(s)})) \dd s) \nonumber \\
%  & =  \E_x ( \int_0^{\tau} \phi(T(B_t)) F'(t) \dd t )\\
  & = \E_x ( \int_0^{\tau} \phi(T(B_t)) |T'(B_t)|^2 \dd t )  \nonumber \\
  & = \int_D G_0^D(x, y) \phi(T(y)) |T'(y)|^2 \dd y. \label{eq:GreenCI2}
  \end{align}
  Identifying the right hand sides of \eqref{eq:GreenCI1} and \eqref{eq:GreenCI2}, 
  since the test function $\phi$ is arbitrary, we conclude that 
  $$
  G_0^{D'}(T(x), T(y)) |T'(y)|^2 
  = G_0^D(x, y) |T'(y)|^2 
  $$ 
first as distributions, and thus by continuity as functions defined for $x\ne y$. The result follows by cancelling the factors of $|T'(y)|^2$ on both sides.
  \end{proof}

\begin{rmk}
\label{R:GreenCI} We have already mentioned that the conformal invariance of the Green function is at first a little surprising, since conformal invariance of Brownian motion holds only up to a time change, whereas the Green function, which measures time spent in a neighbourhood of a point, is a priori very sensitive to the time parametrisation. Having done the proof, we can now \emph{a posteriori} explain this remarkable fact. When we apply the change of variables spatially, we pick up a term $|T'(y)|^2$ because we are in dimension $d =2$. When we apply it temporally, we pick up another term $|T'(y)|^2$ from It\^o's formula. The fact that these two factors match exactly is what gives the conformal invariance of the Green function. 

From this perspective, conformal invariance of the Green function is a miraculous property, unique to the case $d = 2$. In higher dimensions it is not simply a problem of defining conformal maps: if we consider scalings $z \mapsto rz$ (note that this leaves Brownian motion invariant up to time change in any dimension), it is only in dimension $d=2$ that such scalings leave the Green function invariant. 
\end{rmk}

%\hrule 
\begin{rmk}
We will make use of conformal invariance to analyse the Green function in dimension $d=2$, as it often suffices to prove some desired property in a concrete domain (where we have explicit formulae, such as the upper half plane), and use conformal invariance to deduce the desired property in an arbitrary simply connected domain. We believe this to be an elegant approach, appropriate for many potential readers of this book. However, it has the drawback that it does not apply in other dimensions, where instead one must usually rely on hands on estimates. The latter approach also works in dimension $d=2$ of course, and might be more appropriate for readers who do not have a background in complex analysis -- after all, the theory which will be developed throughout the first four chapters of this book depend only very tangentially on complex analysis arguments, and can mostly be read without any such background.
For more on this hands on approach to properties of the Green function, we refer potentially interested readers to Chapter 2.4 in \cite{Lawlerbook}, where none of the arguments appeal to conformal invariance. 
\end{rmk}

%\hrule 

\begin{example}
Having identified the Green function in one simply connected domain (the upper half plane $\H$), the conformal invariance of the Green function can be used in conjunction with the Riemann mapping theorem to evaluate it on an arbitrary simply connected domain $D$. Here is an example in which the Green function becomes very simple. Let $D = \D$ be the unit disc. 
We can find a M\"obius transformation
$$
T(z)= \frac{ i-z}{i+z}
$$
which maps $\H$ to $\D$. (To check this, note that any M\"obius map, that is, any function of the form $z\mapsto (az + b)/ (cz + d)$ with $ad-bc\ne 0$, is always a homeomorphism of the extended plane $\R^2 \cup\{\infty\}$ onto itself, and maps circles to circles -- where by circles we also allow for infinite lines. Here it is easy to check that if $z \in \R$ then $|T(z) | = 1$, so $T$ maps the real line to the unit circle; since $T(i) = 0$ the image of $\H$ must be the unit disc.) From the explicit form of $G_0^\H$ obtained in \eqref{G_H}, we deduce
\begin{equation}\label{Eq:G_D}
G_0^{\D}( 0,z) =- \frac1{2\pi}\log |z|.
\end{equation}
%This makes \eqref{E:cr} obvious for $D=\D$ and $x = 0$, and so \eqref{E:cr} follows immediately in the general case by conformal invariance and definition of the conformal radius.
\end{example}

We state below some basic and fundamental properties of the Green function in two dimensions, which will be used throughout.

\begin{prop} \label{P:basicGFFc} For any regular, simply connected domain $D\subset \R^2$, and any $x\in D$:
	\begin{enumerate}
		\item $G^D_0(x,y)\to 0$ as $y\to y_0\in \partial D$;
		
		\item
$G^D_0(x,y)  =
 - \frac1{2\pi}\log ( | x- y |) + O(1)
$
as $y\to x$.

\item $G^D_0(x,\cdot)$ is harmonic in $D\setminus \{x\}$; and as a distribution 
		\begin{equation}\label{eq:GreenLaplacian}
		\Delta G_0^D(x, \cdot) =-   \delta_{x} ( \cdot) ;
		\end{equation}
	\end{enumerate}
\end{prop}

\begin{proof}
For the first point, observe that on the unit disc $\D$ with $x=0$, $|G^{\D}_0(0,y)| \le C \dist (y, \partial \D)$ for all $y$ with $|y|\ge 1/2$ (say),
so converges to zero, uniformly as $y$ approaches $\partial \D$. Now suppose $D$ is an arbitrary regular simply connected domain and $x \in D$. Fix a conformal isomorphism $f$ from $\D$ to $D$ with $f(0) = x$. Let $y_n\in D$ be a sequence such that $y_n \to y\in \partial D$. Then we claim that $w_n := f^{-1} (y_n) \in \D$ is a sequence approaching the boundary of $\D$, in the sense that $\dist (w_n , \partial \D) \to 0$ (note however that there is no guarantee, 
without additional assumptions on $D$, that $w_n$ will converge to a point on $\partial \D$). Indeed, since $w_n\in \D$ is a bounded sequence, it suffices to check that no subsequence can converge to a point $w\in \D$. But if that were the case, then $f(w_n) $ would converge to $f(w)$ along that subsequence, which contradicts the fact that $y_n $ converges to $y$.
Hence, by conformal invariance of the Green function and the uniformity of the convergence to zero in $\D$, we deduce that $G_0^D(x, y_n) = G_0^{\D} (0, w_n) \to 0$, as required. 

 The second point also follows from the explicit form of the Green function on the unit disc and conformal invariance. In particular, integrals of the form $\int_D G^D_0(x, y) f(y) \dd y$ are well defined for any test function $f \in \cD_0(D)$, so $G_0^D (x, \cdot)$ may be viewed as a distribution.

%, and the Green function in $\H$ is clearly harmonic off the diagonal (as the difference of two harmonic functions). %, since $z \mapsto \log |z| = \Re ( \log z)$ is harmonic.
 For the final point, we can again use the explicit form of $G_0^\D$ on $\D$, which shows that $G^\D_0(0, \cdot)$ is a harmonic function away from 0 (as the real part of a holomorphic function). Furthermore, harmonicity is preserved under conformal isomorphisms. This shows that $G^D_0(x, \cdot)$ is harmonic away from $x$. To prove \eqref{eq:GreenLaplacian} requires a little more. For instance, one can reduce \eqref{eq:GreenLaplacian} by conformal invariance to showing  that $\Delta \log |z| = 2\pi \delta_{0}$ in the sense of distributions. This follows from explicit computations of $\Delta f_\eps(z)$, where $f_\eps(z) = \log (|z| \vee \eps)$, and the fact that $f_\eps(z)$ converges to $\log |z|$ in the sense of distributions, hence $\Delta f_\eps(z)$ converges to $\Delta \log |z|$ in the sense of distributions.  
 
 However perhaps the simplest argument for \eqref{eq:GreenLaplacian} is as follows: since we already know by the second point that $G^D_0(x, \cdot)$ is a distribution, it suffices to show that for each test function $f \in \cD_0(D)$, 
 \begin{equation}\label{eq:GreenLaplacianIto}
 \int_D G^D_0(x,y) \Delta f(y) \dd y = - f(x). 
 \end{equation}
Using \eqref{eq:Gformula} (and using the consequence of the second point above that the integral in \eqref{eq:GreenLaplacianIto} is well defined) the left hand side can be rewritten as 
 $$
 \E_x ( \int_0^{\tau_D} \Delta f ( B_s) \dd s ). 
 $$
 On the other hand, by It\^o's formula, $M^f_t = f(B_{t\wedge \tau_D}) - \int_0^{t\wedge \tau_D} \Delta f(B_s) \dd s$ is a martingale with initial value $M_0^f = f(x)$,
 hence 
 $$
 \E_x ( \int_0^{t \wedge \tau_D}  \Delta f(B_s) \dd s ) = \E_x( f( X_{t \wedge \tau_D} ))  - f(x) .
 $$   
 The result thus follows by letting $t\to \infty$: in the right hand side, the first term tends to zero since $f$ has compact support and $\tau_D < \infty$ almost surely. In the left hand side, we apply the dominated convergence theorem, together with the fact that the expected occupation measure up to time $t$ is dominated by the total expected occupation measure $G_0^D (x, y) \dd y$, which integrates $\Delta f$ by the second point as $\Delta f$ is itself a smooth compactly supported function. This proves \eqref{eq:GreenLaplacianIto} and thus \eqref{eq:GreenLaplacian}.
  
 Instead of such computations, one could also argue that $p^D_t$ solves the heat equation 
 $$
 \frac{\partial}{\partial t} p^D_t(x,y) = \Delta p^D_t(x,y);
 $$
 Integrating this identity over time gives, at least informally, 
 $$
 \Delta G^D_0(x, \cdot) = p^D_\infty (x, \cdot) - p^D_0(x, \cdot) = - \delta_x( \cdot)
 $$
 as desired, where $p^D_\infty$ denotes the limit as $t \to \infty$ of $p^D_t(x,y)$, which is zero because Brownian eventually leaves $D$ in finite time. Of course, justifying this requires some careful arguments too, so the exact computations using the form of $G^{\D}_0$ on $\D$ are more direct.
  %This can be deduced from the form of \eqref{G_H} and conformal invariance.
% but can also be seen directly. For instance, consider point 2. Fix $x \in D$, and let $G_x(\cdot)$ be the harmonic extension, inside $D$, of the boundary values $(-1/\pi) \log (|x- \cdot|)$.
\end{proof}

\begin{rmk}
	\label{GF:generald}
In fact, the above result also holds in other dimensions with appropriate changes. One can check that in any dimension $d \ge 1$, for any regular domain $D$ such that the Green function $\int_0^\infty p_t^D(x,y) \, \dd t$ in $D$ is finite for all $x\ne y$ (note that $D$ does not need to be simply connected), and for any fixed $x \in D$:

	\begin{enumerate}
		\item $G^D_0(x,y)\to 0$ as $ y \to y_0 \in \partial D$;
		\item $G^D_0(x,\cdot)$ is harmonic in $D\setminus \{x\}$ with  $\Delta G_0^D(x, \cdot) =-   \delta_{x} ( \cdot) $ as distributions;
		\item  $$G^D_0(x,y)  = \begin{cases} G_0^D(x,x) + o(1) & d = 1 \\
-(2\pi)^{-1}\log ( | x- y |) + O(1) & d=2 \\\frac{1}{A_d} |x-y|^{2-d} + O(1) & d\ge 3
		\end{cases}$$
		as $y\to x$, where $A_d$ is the $(d-1)$ dimensional surface area of the unit ball in $d$ dimensions.

	\end{enumerate}
\end{rmk}

\begin{rmk}
One can in fact say more than what is contained in Proposition \ref{P:basicGFFc} or Remark \ref{GF:generald}. Firstly, in Theorem \ref{T:GreenCR}, the on diagonal behaviour of the Green function will be estimated more sharply. 
Secondly, the properties contained in Proposition \ref{P:basicGFFc} are in fact sufficient to characterise the Green function, and thus may be used to identify it explicitly. See Exercise \ref{Green_charact} where such a characterisation will be proved, in fact under even weaker assumptions: if $\phi: D \setminus\{ z_0\}\to \R$ is harmonic, converges to 0 near the boundary, and blows up logarithmically near $z_0$, in the sense that $\phi (z) = (1+ o(1)) / (2\pi) \log (|z-z_0|^{-1}) $, then $\phi$ coincides with the Green function. 
 See also \cite[Lemma 3.7]{WWnotes} for more on this, and below for two examples (in dimensions $d=2$ and $d=1$).
\end{rmk}

\begin{example}
Using this characterisation we obtain another, more conceptual, proof of \eqref{Eq:G_D}.
Indeed, it is clear that $z \in \D\setminus\{0\} \mapsto - \log |z|$ is a harmonic function, as the real part of the holomorphic function $\log (z)$ (defined locally, say, which is sufficient for harmonicity). The logarithmic blow up near $z = 0$ is of course obvious in this case.
\end{example}

As another example, this time in dimension $ d= 1$, one can show:
\begin{example} If $d=1$ and $D = (0,1)$, then
\begin{equation}\label{Green1d}
G_0^D (s,t) = s(1-t)
\end{equation}
for $0 <s \le t < 1$ (a more symmetric expression, independent of the relative position of $s$ and $t$, is $(s, t) \in (0,1)^2 \mapsto s\wedge t - st$). Note that this does not blow up on the diagonal.
\end{example}

\begin{comment}In particular we deduce the following.

\begin{prop} \label{P:basicGFFc} We have that for any simply connected domain $D\subset \C$, and any $x\in D$:
\begin{enumerate}
\item $G^D_0(x,\cdot)$ is harmonic in $D\setminus \{x\}$; and as a distribution $\Delta G_0^D(x, \cdot) =- 2\pi  \delta_{x} ( \cdot) $
\item $G^D_0(x,y)  = - \log ( | x- y |) + O(1)$ as $y \to x$.
\end{enumerate}
\end{prop}
\begin{proof}
For the first point, observe that harmonic functions stay harmonic under composition with a conformal map, and the Green function in $\H$ is clearly harmonic off the diagonal (as the difference of two harmonic functions). %, since $z \mapsto \log |z| = \Re ( \log z)$ is harmonic.
 For the second point, we can again use the explicit form of $G_0^\H$ on $\H$ and conformal invariance.
%This can be deduced from the form of \eqref{G_H} and conformal invariance.
% but can also be seen directly. For instance, consider point 2. Fix $x \in D$, and let $G_x(\cdot)$ be the harmonic extension, inside $D$, of the boundary values $(-1/\pi) \log (|x- \cdot|)$.
\end{proof}

The converse to this proposition also holds, that is, these properties characterise the Green function (see Exercise \ref{Green_charact}). In fact, this can be strengthened: it suffices to assume that $G_0(x, \cdot)$ is harmonic and behaves like $- \log |x - \cdot|$ up to a multiplicative error term $1 + o(1)$ to deduce that $G_0(x, \cdot)$ is the Green function.
\end{comment}

\medskip

Actually, one can be slightly more precise about the behaviour of the Green function near the diagonal; that is, one can find a sharper estimate for the error term $O(1)$ in  \cref{P:basicGFFc}:

\begin{theorem}\label{T:GreenCR}
  \ind{Conformal radius}
\begin{equation}\label{E:cr}
G^D_0(x,y)  = - \frac1{2\pi}\log ( | x- y |) + \frac1{2\pi}\log R(x;D) + o(1)
\end{equation}
as $y \to x$, where $R(x;D)$\indN{{\bf Geometry}! $R(x;D)$; conformal radius of $x$ in $D$} is the \textbf{conformal radius} of $x$ in $D$. That is, $R(x;D) = |f'(0)|$ for $f$ any conformal isomorphism taking $\D$ to $D$ and satisfying $f(0) = x$.

Furthermore, we may write 
\begin{equation}\label{E:Green_harmonic_ext}
G_0^D(x,\cdot) = - \frac1{2\pi} \log |x - \cdot| + \xi_x(\cdot),
\end{equation}
 where $\xi_x(\cdot)$ is a harmonic function over all of $D$, which equals the harmonic extension to $D$ of the function $1/(2\pi) \log (|x - \cdot|)$ on $\partial D$.  (Combining with \eqref{E:cr}, we must have $\xi_x(x) = 1/(2\pi) \log R(x;D)$). 
\end{theorem}

\begin{proof} Recall from \eqref{Eq:G_D} that if $D = \D$ is the unit disc, we have
$$
G_0^{\D}( 0,z) =- \frac1{2\pi}\log |z|.
$$
This makes \eqref{E:cr} obvious for $D=\D$ and $x = 0$, and so \eqref{E:cr} follows immediately in the general case by conformal invariance, Taylor approximation, and definition of the conformal radius. To prove \eqref{E:Green_harmonic_ext}, we note that the difference $G_0^D (x, \cdot) + 1/(2\pi) \log |x - \cdot|$ is bounded and harmonic in the sense of distributions in all of $D$. Elliptic regularity arguments, or direct argumentation with planar Brownian motion, imply that this is a smooth function which is harmonic in the usual sense. 
\end{proof}

\begin{rmk}
Note that the conformal radius is unambiguously defined: the value $|f'(0)|$ does not depend on the choice of $f$ ($f$ is unique up to rotation, which does not affect the value of the modulus derivative). Although we will not use this, we note also that by the classical \textbf{K\"obe quarter theorem}, we have
$$
\dist (x, \partial D) \le R(x; D) \le 4 \dist(x, \partial D)
$$
so the conformal radius is essentially a measure of the Euclidean distance to the boundary.
\end{rmk}

The conformal radius appears in Liouville quantum gravity in various formulae which will be discussed later on in the course. The reason it shows up in these formulae is usually because of \eqref{E:cr}.

The last property of $G^D_0$ that we will need, as in the discrete case, is that it is a non-negative definite function. We will see this in the next section.

\subsection{GFF as a stochastic process}\label{DGFF:sp}

From now on, we will always assume that $D\subset \R^d$ is a regular domain; that is, an open connected set with regular boundary.

Essentially, as in the discrete case, we would like to define the GFF as a Gaussian ``random function'' with mean zero and covariance given by the Green function. However (when $d\ge 2$) the divergence of the Green function on the diagonal means that the GFF cannot be defined pointwise, as the variance at any point would have to be infinite. So instead, we define it as a random distribution, or generalised function in the sense of Schwartz\footnote{This conflicts
with the usage of distribution to mean the law of a random variable but is standard and should not cause confusion.}. More precisely, we will take the point of view that it assigns values to certain measures with finite Green energy. In doing so we follow the approach in the two sets of lecture notes \cite{SLEnotes} and \cite{WWnotes}. The latter in particular contains a great deal more about the relationship between the GFF, SLE, Brownian loop soups and conformally invariant random processes in the plane, which will not be discussed in this book. The foundational paper by Dub\'edat \cite{Dubedat} is also an excellent source of information regarding basic properties of the Gaussian free field. We should point out that the rest of this text is particularly focused on the case $d=2$ but we will include results relevant to other dimensions when there is no cost in doing so.

\medskip Recall that if $I$ is an index set, a \textbf{stochastic process indexed by $I$}  is just a collection of random variables $(X_i)_{i\in I}$, defined on some given probability space. The \textbf{law} of the process is a measure on $\R^I$, endowed with the product topology. It is uniquely characterised by its finite dimensional marginals, that is, the law of $(X_{i_1},\dots, X_{i_n})$ for arbitrary $i_1,\dots, i_n$ in $I$, via Kolmogorov's extension theorem.

\medskip Given $n \ge 1$, a random vector $X = (X_i)_{1\le i \le n}$ is called \textbf{Gaussian} if any linear combination of its entries is real Gaussian; that is, if $\langle \lambda, X \rangle $ is a real Gaussian random variable for any $\lambda = (\lambda_1, \ldots, \lambda_n) \in \R^n$, $n\in \N$. The law of $X$ is uniquely specified by its mean vector $\mu = \E(X) \in \R^n$, that is, $\mu_i = \E(X_i)$ for each $1\le i \le n$, and its covariance matrix $\Sigma \in \cM(\R^n)$ given by $\Sigma_{i,j} = \cov (X_i, X_j)$, $1\le i, j \le n$. Conversely, given a vector $\mu \in \R^n$ and a symmetric, non-negative\footnote{Here non-negative is in the sense of matrices, that is, $\sum_{i,j} \lambda_i \lambda_j \Sigma_{i,j} \ge 0$ for each $\lambda \in \R^n$, or equivalently, the eigenvalues of $\Sigma$ are all non-negative.} matrix $\Sigma \in \cM(\R^n)$, there exists a (unique) law on $\R^n$ which is that of a Gaussian vector with mean $\mu$ and covariance matrix $\Sigma$.

Fix a set $I$ and suppose we are given a function $C: I \times I \to \R$, symmetric and non-negative in the sense that 
\begin{equation} \label{E:nonneg} 
	\sum_{i,j=1}^n \lambda_i \lambda_j C(t_i, t_j) \ge 0 \quad \forall n\ge 1, t_1, \ldots, t_n \in I \text{ and } \lambda_1, \ldots, \lambda_n\in \R.
 \end{equation} 
Then associated to this function $C$, for each $t_1,\dots, t_n \in I$ we can define a centred Gaussian vector $(X_{t_1}, \ldots, X_{t_n})$ with covariance matrix $\Sigma_{i,j} = C(t_i, t_j)$, $1\le i,j\le n$. The resulting laws are automatically consistent, in the sense of Kolmogorov as the parameters $t_1, \ldots, t_n \in I$ and $n\ge 1$ are varied. Therefore, by Kolmogorov's extension theorem, the function $C$ defines a unique law on $\R^I$. This is the law of a stochastic process $(X_t)_{t\in I}$ indexed by $I$ such that the restriction of $(X_t)_{t\in I}$ to any $n$ tuple of indices $t_1, \ldots, t_n \in I$ gives us a centred Gaussian vector $(X_{t_1}, \ldots, X_{t_n})$ with the above covariance matrix. The process $(X_t)_{t\in I}$ is called the (centred) \textbf{Gaussian stochastic process on $I$ with covariance function $C$}. Given a real valued function $(\mu(t), t \in I)$ we can also define a Gaussian stochastic process $Y$ on $I$ with mean function $\mu$ and covariance function $C$, simply by shifting $X$ by $\mu(t)$ at each $t\in I$, that is, setting $(Y_t)_{t\in I}:=(X_t+\mu(t))_{t\in I}$. 
\\

 Now, let $D \subset \R^d$ be an open set with \textbf{regular} boundary, and recall from Section \ref{SS:Greencontinuous} that for such $D$ the Green function $G_0^D$ is finite away from the diagonal: $G_0^D(x,y)<\infty$ for $x\ne y$. \ind{Domain!Regular}
% As we have already mentioned, this is the case as soon as $D$ is a proper domain of $\R^d$ (if $d \le 2$) and if $D$ is regular. In particular, $D$ does not need to be simply connected or even connected.
We will define the Gaussian free field in $D$ (with zero boundary conditions) as a centred Gaussian stochastic process indexed by the set $\mathfrak{M}_0$ (defined below) of signed Borel measures with finite logarithmic energy.

\begin{definition}[Index set for the GFF]Let $\mathfrak{M}_{0}^+$ 
	\indN{{\bf Function spaces}! $\mathfrak{M}_0^+$; non-negative measures $\rho$ supported in $D$ with finite integral tested against $G_0^D$} 
	denote the set of (non-negative) Radon measures supported in $D$, such that
$\int \rho(\dd x) \rho(\dd y)  G_0^D(x,y) < \infty$. 
Denote by $\mathfrak{M}_0$ \indN{{\bf Function spaces}! $\mathfrak{M}_0$; signed measures of the form $\rho=\rho^+-\rho^-$ with $\rho_\pm\in \mathfrak{M}_0^+$} 
the set of signed measures of the form $\rho = \rho^+ - \rho^-$ with $\rho^{\pm}\in \mathfrak{M}_0^+$. 
\end{definition}

Note that when $d=2$, due to the logarithmic divergence of the Green function  on the diagonal, $\mathfrak{M}_0^+$  includes the case where $\rho(\dd x) = f(x) \dd x$ and $f$ is continuous, but does not include Dirac point masses.

For test functions $\rho_1,\rho_2 \in \mathfrak{M}_0$, we set
\begin{equation}\label{E:Gamma_0} \Gamma_0(\rho_1,\rho_2) := \int_{D^2}  G_0^D(x,y) \rho_1(\dd x) \rho_2(\dd y)\end{equation} 
\indN{{\bf Green functions}! $\Gamma_0$; bilinear form, doubly integrating against $G_0^D$}
and also define $\Gamma_0 (\rho) = \Gamma_0( \rho, \rho)$. We will see below why these quantities are in fact well defined, but note for now that this is not immediately obvious.
\\

Essentially, our definition will be that the Gaussian free field on $D$ with zero boundary conditions is the centred Gaussian stochastic process $(\Gamma_\rho)_{\rho \in \mathfrak{M}_0}$ indexed by $\mathfrak{M}_0$ such that for $\rho_1, \rho_2 \in \mathfrak{M}_0$ we have 
$$\cov (\Gamma_{\rho_1}, \Gamma_{\rho_2}) = \Gamma_0( \rho_1, \rho_2)  = \int_{D^2} G_0^D(x,y) \rho_1( \dd x) \rho_2( \dd y).$$ 
However in order to do so a few things need to be checked. Namely:
\begin{itemize}
\item $\Gamma_0( \rho_1, \rho_2) $ is well defined whenever $\rho_1, \rho_2 \in \mathfrak{M}_0$. In fact this is not obvious,\footnote{The necessity of such an argument (in the absence of any form of Cauchy--Schwarz inequality at this stage) seems to not have been noticed before; correspondingly Lemma \ref{L:linearitycov}, although not difficult, is new.} even if we assume $\rho_1, \rho_2 \in \mathfrak{M}_0^+$.

\item The function $\Gamma_0 ( \cdot,\cdot)$ is symmetric and non-negative on $\mathfrak{M}_0\times \mathfrak{M}_0$,  in the sense of \eqref{E:nonneg} with $I=\mathfrak{M}_0$, so is a valid covariance function. 
\end{itemize}
\ind{Nonnegative definite}

As we will see, these properties will follow rather easily from the following lemma. 

\begin{lemma}\label{L:linearitycov}
If $\rho_1, \rho_2 \in \mathfrak{M}_0^+$ then $\Gamma_0( \rho_1, \rho_2) < \infty$. Furthermore $\rho_1 + \rho_2 \in \mathfrak{M}_0^+$. 
\end{lemma}

\begin{proof}
%At this stage it is neither clear that $\mathfrak{M}_0$ is a vector space, nor that if $\rho_1, \rho_2 \in \cM_0$ then $\Gamma_0(\rho_1, \rho_2)$ (that is, the right hand side of \eqref{E:Gamma_0}) is well defined, even if $\rho_1,\rho_2 \in \mathfrak{M}_0^+$, or even if $\rho_1 = \rho_2 \in \mathfrak{M}_0 \setminus \mathfrak{M}_0^+$. To see this (and to prepare for the upcoming definition), we make the following observation. 
By the Markov property, we have
  $$
  p_t^D(x,y) = \int_D p^D_{t/2} (x,z) p^D_{t/2}(z,y) \dd z,
  $$
  and hence by symmetry (that is, $p_t^D(x,y) = p^D_t(y,x)$, which follows from the same symmetry in the full plane, and the fact that a Brownian bridge from $x$ to $y$ has as much chance to stay in $D$ as one from $y$ to $x$, as one is the time reversal of the other), we can deduce that
   $$
  G_0^D(x,y) = 2 \int_D   \dd z \int_0^\infty  p^D_u(x,z) p^D_u(y,z) \dd u.
  $$
  Consequently, if $\rho_1, \rho_2 \in \mathfrak{M}_0^+$ are arbitrary, 
  \begin{align}
\Gamma_0 (\rho_1, \rho_2) & = \iint G^D_0(x, y) \rho_1(\dd x) \rho_2(\dd y) \nonumber \\
 & = \int_D 2 \dd z \int_0^\infty \iint \rho_1(\dd x) \rho_2(\dd y) p^D_u(x, z) p^D_u(y, z) \dd u \nonumber \\
&  = \int_D 2 \dd z \int_0^\infty \left( \int \rho_1(\dd x) p^D_u(x, z)\right) \times \left( \int \rho_2(\dd x) p^D_u(x, z)\right) \dd u . \label{eq:offdiagGamma} 
%    & 
  \end{align}
In particular, if $\rho_1 = \rho_2\in \mathfrak{M}_0^+$ then 
\begin{equation}
\Gamma_0(\rho_1, \rho_1) = \int_D 2 \dd z \int_0^\infty \left( \int \rho_1(\dd x) p^D_u(x, z)\right)^2 \dd u <\infty.%\ge 0.%
 \label{eq:diagGamma}
\end{equation}
%for any $\rho_1, \rho_2 \in \mathfrak{M}_0^+$. 
Hence using the inequality $2ab \le a^2 + b^2$, valid for any real numbers $a$ and $b$, we deduce that $\Gamma_0( \rho_1, \rho_2) < \infty$ whenever $\rho_1, \rho_2 \in \mathfrak{M}_0^+$. This proves the first point. 

For the second point, observe that \emph{a priori} $\Gamma_0( \rho_1+ \rho_2) = \Gamma_0 ( \rho_1) + 2 \Gamma_0( \rho_1, \rho_2) + \Gamma_0( \rho_2)$. This is an equality between terms which are non-negative but might be infinite. Nevertheless, from what we have just seen, if $\rho_1, \rho_2 \in \mathfrak{M}_0^+$, all three terms on the right hand side are finite. Thus the left hand side is finite too, which concludes the proof of Lemma \ref{L:linearitycov}.
\end{proof}

Lemma \ref{L:linearitycov} allows us to extend the notion of energy $\Gamma_0(\rho_1, \rho_2)$ onto $\mathfrak{M}_0\times \mathfrak{M}_0$ and not just $\mathfrak{M}_0^+\times \mathfrak{M}_0^+$, justifying the definition in \eqref{E:Gamma_0}. Indeed writing $\rho_i = \rho_i^+ - \rho_i^-$ for $i = 1, 2$, we have 
$$
\Gamma_0(\rho_1, \rho_2) = \Gamma_0( \rho_1^+, \rho_2^+) + \Gamma_0( \rho_1^-, \rho_2^-) - \Gamma_0( \rho_1^+, \rho_2^-) - \Gamma_0( \rho_1^-, \rho_2^+);
$$
where the finiteness of all four terms on the right hand side is guaranteed by Lemma \ref{L:linearitycov}.
Note also that $\mathfrak{M}_0$ is a vector space (again by Lemma \ref{L:linearitycov}), with $\Gamma_0$ a \textbf{bilinear form} on $\mathfrak{M}_0$.

%This allows us to define the right hand side of \eqref{E:Gamma_0} by taking the appropriate linear combination of terms of the form $\Gamma (\rho, \rho')$ with $\rho, \rho' \in \mathfrak{M}_0^+$. The same reasoning also defines $\Gamma (\rho, \rho)$ for arbitrary $\rho \in \mathfrak{M}_0$ (which was \emph{a priori} only well defined for $\rho \in \mathfrak{M}_0^+$). Furthermore, if we write $\rho = \rho^+ - \rho^-$ with $\rho^\pm \in \mathfrak{M}_0^+$, since $\Gamma_0(\rho^-, \rho^+) = \Gamma_0 (\rho^+, \rho^-)$ (which follows from the symmetry of of $G_0^D(\cdot, \cdot)$, itself a consequence once again of the symmetry of $p_t^D ( \cdot, \cdot)$), we have 
%\begin{equation}
%\Gamma_0 ( \rho, \rho)  = \Gamma( \rho^+, \rho^+)  + \Gamma ( \rho^-, \rho^-) \ge 0.
%\label{Covariance_definite}
%\end{equation}
%The nonnegativity of both terms above follows either from \eqref{eq:diagGamma} or from the definition of $\Gamma_0$, since $G^D_0(x,y) \ge 0$ for every $x,y \in D$ by definition of the Green function. 

% Note also that if $\rho, \rho' \in \mathfrak{M}^+_0$ then 
% $$\Gamma_0( \rho + \rho', \rho+ \rho') = \Gamma_0 ( \rho, \rho') + 2\Gamma_0 (\rho, \rho') + \Gamma_0 ( \rho', \rho').$$
%  In particular the left hand side is finite, hence $\mathfrak{M}_0$ is a vector space.  

  \begin{rmk}
  In fact, we will soon see as a consequence of Lemma \ref{L:M0Sobnorm} that $\mathfrak{M}_0$ is the intersection of the Sobolev space $H^{-1}_0(D)$ with the set of signed measures on $D$. 
  \end{rmk}

\begin{lemma}\label{L:posdef}
The bilinear form $\Gamma_0$ is symmetric and non-negative (in the sense of covariance functions) on $\mathfrak{M}_0\times \mathfrak{M}_0$. That is, for every $n\ge 1$ and every $\rho_1, \ldots, \rho_n \in \mathfrak{M}_0$, for every $\lambda_1, \ldots, \lambda_n \in \R$, 
$$
\sum_{i,j=1}^n\lambda_i \lambda_j \Gamma_0( \rho_i, \rho_j) \ge 0. 
$$
In particular, $\Gamma_0$ is a valid covariance function for a Gaussian stochastic process on $\mathfrak{M}_0$.
\end{lemma}
\ind{Nonnegative definite}

\begin{proof}
Since $\Gamma_0$ is a bilinear form, we have:
$$
\sum_{i,j=1}^n \lambda_i \lambda_j \Gamma_0( \rho_i, \rho_j) = \Gamma_0( \rho)
$$
where 
$$
\rho = \sum_{i=1}^n \lambda_i \rho_i \in \mathfrak{M}_0.
$$
The desired non-negativity therefore follows directly from \eqref{eq:diagGamma}. 
\end{proof}

\medskip As a consequence of Lemma \ref{L:posdef}. we can now finally give the definition of a Gaussian free field (with zero boundary conditions) as a stochastic process. 
%The following is both an easy theorem and the definition of the GFF with zero boundary conditions on a domain.

\begin{theorem}[Zero boundary or Dirichlet GFF]\label{D:GFFc}
There exists a unique stochastic process $(\mathbf{h}_\rho)_{\rho \in \mathfrak{M}_0}$, indexed by $\mathfrak{M}_0$,
such that for every choice of $\rho_1, \ldots, \rho_n$, $(\mathbf{h}_{\rho_1}, \ldots, \mathbf{h}_{\rho_n})$ is a
centred Gaussian vector with covariance structure $ \cov (\mathbf{h}_{\rho_i}, \mathbf{h}_{\rho_j}) = \Gamma_0 (\rho_i, \rho_j)$.
\end{theorem}
\ind{GFF! Dirichlet boundary conditions}

\begin{comment}
\begin{proof}
We need to check several things:
\begin{itemize}
\item the finite-dimensional distributions exist and are uniquely specified;
\item they are consistent.
\end{itemize}
\index{Nonnegative semidefinite}
The consistency is an immediate consequence of the Gaussianity of the finite dimensional marginals: indeed, the restriction of a Gaussian vector to a subset of coordinates is still a Gaussian vector with the same covariance structure.

For the first point, symmetry is a consequence of the reversibility of the heat kernel $p_t^D(x,y) = p_t^D(y,x)$.
To check the non--negative semi-definite character of $\Gamma$, we need to check that for every $\rho_1, \ldots, \rho_n \in \mathfrak{M}_0$ and every $\lambda_1, \ldots, \lambda_n \in \R$, we have that
$$
\sum_{i,j} \lambda_i \lambda_j \Gamma_0 (\rho_i, \rho_j) \ge 0.
$$
However, by linearity of integration, this sum is nothing but $\Gamma_0 (\rho) $ for $\rho = \sum_i \lambda_i \rho_i$. Hence it suffices to prove that
\begin{equation}\label{ND}
\Gamma_0(\rho, \rho) \ge 0 \text{ for } \rho \in \mathfrak{M}_0.
\end{equation}
But we already proved this in \eqref{Covariance_definite},
%We give two different proofs; one based on Proposition \ref{P:Green_sd} and another based on integration by parts.
and so finishes the proof of the theorem.
\end{proof}
\end{comment}

Let us emphasise that for a stochastic process, the index set $I$ does not a priori \emph{need} to be a vector space, although this is the case when $I=\mathfrak{M}_0$. Similarly, the covariance function of a Gaussian stochastic process indexed by $I$ does not \emph{need} to be a bilinear non-negative form on $I$, although again this is true for $\Gamma_0$ on $\mathfrak{M}_0$, and this helped us to prove its validity as a covariance function.

\begin{definition}
  The process $(\mathbf{h}_{\rho})_{\rho \in \mathfrak{M}_0}$ is called the Gaussian free field in $D$ (with Dirichlet or zero boundary conditions). We write GFF as shorthand for Gaussian free field.
  %We will write $(h,\rho)$ instead of $h_\rho$ when $\rho \in \mathfrak{M}$, and think of $(h, \rho)$ integrated against $\rho$.
\end{definition}

Note that in such a setting, it might not be possible to ``simultaneously observe'' more than a countable number of random variables, because our $\sigma$-algebra for the stochastic process $(\mathbf{h}_\rho)_{\rho \in \mathfrak{M}_0}$ is the product $\sigma$-algebra, which is generated by the random variables of the form $(\mathbf{h}_{\rho_1}, \ldots, \mathbf{h}_{\rho_n})$, $n\ge 1$, $\rho_1, \ldots, \rho_n \in \mathfrak{M}_0$.
 A good analogy is with the construction of one dimensional Brownian motion $(B_t, t \ge 0)$: so long as it is constructed as a Gaussian stochastic process indexed by time, numerical quantities such as $\sup_{s \in [t_1, t_2]} B_s$ are not measurable with respect to the product $\sigma$ algebra and so are not random variables. In the case of Brownian motion, it is not until a continuous modification is constructed that such quantities can be seen as (measurable) random variable.  
Likewise, in the case of the GFF, we will have to rely on the existence of suitable modifications with nice continuity properties. More precisely, this modification will be a random distribution living in a certain Sobolev space of negative index, see Section \ref{SS:randomdistributions}, whose law as a stochastic process indexed by $\mathfrak{M}_0$ has the same finite dimensional marginals as the GFF $(\mathbf{h}_\rho)_{\rho \in \mathfrak{M}_0}$. In other words, this random distribution defines a \emph{version} of the GFF.

\begin{rmk} (Terminology). We will use the terminology ``Dirichlet GFF'', ``zero boundary GFF'' and ``GFF with zero/Dirichlet boundary conditions'' interchangeably throughout. With a slight abuse of vocabulary, some authors use the term ``Dirichlet boundary condition'' to indicate that the field has \emph{some} specified (deterministic) boundary conditions, which however may not be identically zero. It will be made clear in the sequel if we wish to talk about anything other than the zero boundary condition case.
\end{rmk}

\begin{rmk}
  In Liouville quantum gravity and in Gaussian multiplicative chaos, it is more convenient (as mentioned previously) to work directly with a field which is logarithmically correlated (as opposed to $(2\pi)^{-1}-$logarithmically correlated), that is, with
  \begin{equation}\label{Dnormalisation}
  h = \sqrt{2\pi} \mathbf{h}.
  \end{equation}
We will use the notations $h$ and $\mathbf{h}$ throughout to make the distinction between these two different conventions.
\end{rmk}

The following property of ``almost sure'' linearity is a consequence of the fact that the covariance function $\Gamma_0$ is a bilinear form on $\mathfrak{M}_0$; its proof is left as an exercise. 

\begin{prop}(Linearity).\label{P:lin} If $\lambda, \lambda' \in \R$ and $\rho, \rho' \in \mathfrak{M}_0$ then $\mathbf{h}_{\lambda \rho + \lambda' \rho' } = \lambda \mathbf{h}_\rho + \lambda' \mathbf{h}_{\rho'}$ almost surely. 
\end{prop}

%(Indeed the two random variables are jointly Gaussian with a trivial covariance matrix).
 In the rest of this text, we will abuse notation slightly and write $(\mathbf{h},\rho)$ for $\mathbf{h}_\rho$ when $\rho \in \mathfrak{M}_0$. We will think of $(\mathbf{h}, \rho)$ as ``$\mathbf{h}$ integrated against $\rho$'', as if $\mathbf{h}$ were an actual distribution, and $\rho$ was a test function. In Section \ref{SS:randomdistributions}, we will see that a version of $\mathbf{h}$ can be defined as a random variable taking values in the space of distributions.   

At this stage, simply note that if $\mathbf{h}$ is a GFF, it cannot be evaluated pointwise (because $\rho= \delta_{x}$ does not lie in $\mathfrak{M}_0$). However it may be tested against smooth, compactly supported test functions $\rho \in \cD_0(D)$. In fact, $h$ may be tested against relatively more singular measures: for instance, the ``integral'' (in the above sense) of $\mathbf{h}$ along a one dimensional segment or a circular arc is always well defined, since the Lebesgue measure on such a one dimensional smooth curve is an element of $\mathfrak{M}_0$. Indeed, one can deduce this from the fact that the divergence of the Green function is only logarithmic, and  that in one dimension, $\int_0^1 \log (r^{-1}) dr  < \infty$.

By Proposition \ref{P:lin}, an alternative definition of the GFF is simply as the unique stochastic process $(\mathbf{h}, \rho)$ which: 
\begin{itemize}
	\item is almost surely linear in $\rho$ (in the sense that for every $\lambda_1, \lambda_2 \in \R$ and $\rho_1, \rho_2 \in \mathfrak{M}_0$,  $(\mathbf{h}, \lambda_1\rho_1+\lambda_2 \rho_2  ) = \lambda_1 (\mathbf{h}, \rho_1) + \lambda_2(\mathbf{h}, \rho_2)$ almost surely); and
	\item is such that $(\mathbf{h}, \rho)$ is a centred Gaussian random variable with variance $\Gamma_0(\rho)$ for every $\rho\in \mathfrak{M}_0$.
\end{itemize}

\medskip \noindent \textbf{Example.} Suppose that $d=1$ and $D = (0,1)$. Then by \eqref{Green1d} we know that $G_0^D(x,y) = x(1- y)$ for $0 < x \le y <1$, and this turns out to be the covariance of a (speed one) \textbf{Brownian bridge} $(b_s, 0 \le s \le 1)$ (see Chapter 1.3 of \cite{RevuzYor}). So, a zero boundary Gaussian free field in one dimension is simply a (speed one) Brownian bridge,
at least in the sense of stochastic processes indexed by, say, test functions.

\medskip Other boundary conditions than zero will also be relevant in practice. For this, we make the following definition (in the case $d=2$ for simplicity). Suppose that $f$ is a (possibly random) continuous function on the \emph{conformal boundary} of a simply connected domain $D\subset \C$ (equivalent to the Martin boundary of the domain for Brownian motion). Then the GFF with boundary data given by $f$ is the random variable $\mathbf{h} = \mathbf{h}_0 + \varphi$, where $\mathbf{h}_0$ is an independent Dirichlet GFF, and $\varphi$ is the harmonic extension of $f$ to $D$. 

The reason for this definition will become clear in light of the Markov property discussed in \cref{SS:Markov}. Alternatively, 
it can be justified by the fact that if one defines a discrete GFF with prescribed boundary condition $f$ by modifying \cref{T:dGFF_Dirichlet} in the natural way (that is, taking the same definition but setting $h(y)=f(y)$ for $y$ on the boundary), then for an appropriate sequence of approximating graphs, the discrete GFF with boundary condition $f$ converges to $\mathbf{h}_0+\varphi$ as defined above. See \cref{SS:scalinglimit} for the proof of such a statement in the case $f\equiv 0$.

%\begin{rmk}
%  Note that $(\mathbf{h}, \rho)$ is linear in $\rho$: if $\rho, \rho' \in \mathfrak{M}_0$, and if $\alpha, \beta \in \R$ $(h, \alpha \rho + \beta \rho') = \alpha (h, \rho) + \beta (h, \rho')$, almost surely (which can be seen by checking that the variance and mean of the difference are both zero). Hence an alternative definition of the GFF is simply as the unique stochastic process $(h, \rho)$ which is linear in $\rho$ and such that $(h, \rho)$ is a centred Gaussian random variable with variance $\Gamma_0(\rho)$.
%\end{rmk}

\medskip If we do not specify the boundary conditions, we always mean a Gaussian free field with zero (or Dirichlet) boundary conditions.

\subsection{Random variables and convergence in the space of distributions}\label{SS:randomdistributions}

As we will soon see, the Gaussian free field can be understood as a random distribution. 
However, since the space of distributions is not metrisable, we first need to address a few foundational issues related to measurability and convergence.

Let $\cD_0(D)$\indN{{\bf Function spaces}!$\cD_0(D)$; compactly supported smooth functions in $D$, or test functions} denote the set of compactly supported, $C^\infty$ functions in $D$, also known as \textbf{test functions}.
The set $\cD_0(D)$ is equipped with a topology in which convergence is characterised as follows.
A sequence $(f_n)_{n\ge 0}$ converges to $0$ in $\cD_0(D)$ if and only if there is a compact set $K\subset D$ such that
$\text{supp} f_n\subset K$ for all $n$ and $f_n$ and all its derivatives converge to $0$ uniformly on $K$.
A continuous linear map $u:\cD_0(D)\to\R$ is called a {\bf distribution} on $D$.
Thus, the set of distributions on $D$ is the dual space of $\cD_0(D)$.
It is denoted by $\cD_0'(D)$\indN{{\bf Function spaces}! $\cD_0'(D)$; distributions on $D$, that is, dual space of $\cD_0(D)$} and is equipped with the weak-$*$ topology.
In particular, $u_n\to u$ in $\cD_0'(D)$ if and only if $u_n(\rho)\to u(\rho)$ for all $\rho\in\cD_0(D)$.

\medskip Let $(\Omega, \cF, \P)$ be a probability space. A random variable $X$ in the space of distributions is, as always, a function $X: \Omega \to \cD'_0(D)$ which is measurable with respect to the Borel $\sigma$-field on $\cD'_0(D)$ induced by the weak-$*$ topology. 

Let $(X_n)_{n\ge 1}$ be a sequence of random variables in $\cD'_0(D)$. We will often ask ourselves whether this sequence converges in $\cD'_0(D)$. However, since the topology of convergence on $\cD'_0(D)$ is not metrisable, it is not clear \emph{a priori} if the event (or rather the subset of $\Omega$)
$$
E= \{ \omega \in \Omega: X_n(\omega) \text{ is weak-}*\text{convergent}\}
$$ 
is measurable. We show here that it is.

\begin{lemma}\label{L:measurability}
Let $D$ be a domain of $\R^d$. Let $\Conv$ denote the set of sequences in $\cD'_0(D)$ which are weak$-*$ convergent. Then $\Conv$\indN{{\bf Function spaces}! $\Conv$; the set of sequences in $\cD'_0(D)$ which are weak$-*$ convergent} is a Borel set in $\cD'_0(D)^{\N}$ equipped with the product Borel $\sigma$-algebra. 
\end{lemma}

See Appendix \ref{App:distributions} for the proof.

\subsection{Integration by parts and Dirichlet energy}

 In order to do view the Gaussian free field as a random variable in the space of distributions, our first step is to relate the covariance of the GFF to the Dirichlet energy of a function (as in the discrete case). The following \textbf{Gauss--Green formula}, which is really just an integration by parts formula, will allow us to do so.

\begin{lemma}[Gauss--Green formula]\label{L:ipp}
Suppose that $D$ is a $C^1$ smooth domain. If $f,g$ are smooth functions on $\bar D$, then
 \begin{equation}\label{ipp}
 \int_D \nabla f\cdot  \nabla g = -\int_D f \Delta g + \int_{\partial D} f \frac{\partial g}{\partial n},
 \end{equation}
 where $\tfrac{\partial g}{\partial n}$ denotes the (exterior) normal derivative.
 \end{lemma}
 \ind{Gauss--Green formula}
 
  \begin{rmk} 
  	For general $D$, the formula holds whenever $g\in \cD_0(D)$ and $f$ is continuously differentiable on $D$ (with the boundary term on the right equal to zero). When $f\in \cD_0'(D)$ is a distribution, the distributional derivative $\nabla f$ is defined to be the distribution such that \eqref{ipp} holds (again with zero boundary term) for all $g\in \cD_0(D)$. 
\end{rmk}

 With  \cref{L:ipp} in hand, we can now rewrite the variance $\Gamma_0(\rho,\rho)$ of $(\mathbf{h}, \rho)$ in terms of the Dirichlet energy of an appropriate function $f$. This Dirichlet energy is of course the continuous analogue of the discrete Dirichlet energy which we encountered in Theorem \ref{T:dGFF_Dirichlet} for instance. 

\begin{lemma}\label{L:definite}
  Suppose that $D$ is a regular domain, $f \in \cD_0(D)$ and that $\rho$ is a smooth function such that $- \Delta f = \rho$. Then $\rho \in \mathfrak{M}_0$ and
  \begin{equation}\label{definite}
  \Gamma_0 (\rho, \rho) = \int_D |\nabla f|^2.
  \end{equation}
\end{lemma}

\begin{proof}
By the Gauss--Green formula (\cref{L:ipp}), noting that there are no boundary terms arising in each application, we have that
\begin{equation*}
\Gamma_0 ( \rho) = - \int_x\rho(x) \int_y G^D_0(x,y) \Delta_y f(y) \dd y \dd x 
 = - \int_x  \rho(x) \int_y \Delta_y G^D_0(x,y) f(y) \dd y \dd x. \end{equation*}
Then using that $\Delta G^D_0(x, \cdot) =- \delta_x(\cdot)$ (in the distributional sense, see \cref{P:basicGFFc}), we conclude that this is equal to 
\begin{equation*}   \int_x \rho(x) f(x) \dd x 
=  - \int_D (\Delta f(x)) f(x) \dd x  = \int_D |\nabla f(x)|^2 \dd x
\end{equation*}
as required.
\end{proof}

Note that this gives another proof that $\Gamma_0(\rho, \rho) \ge 0$, and therefore that the GFF is well defined as a Gaussian stochastic process (at least when indexed by smooth functions $\rho$). Indeed, when $\rho$ is smooth one can always find a smooth function $f$ such that $- \Delta f = \rho$: simply define
\begin{equation}
\label{eqn:frho}
f(x) = \int G_0^D(x, y ) \rho(y) \dd y.
\end{equation}

The following lemma will also be useful. 
\begin{lemma}\label{L:M0dual}
	Suppose that $\rho\in \mathfrak{M}_0$ and $g\in \cD_0(D)$. Then 
	$$\left|\int_D g(x) \rho(\dd x)\right|^2\le \Gamma_0(\rho) \int_D |\nabla g(x)|^2 \dd x.$$
\end{lemma}

\begin{proof}
	It is a simple exercise to check, using dominated convergence, that if $\rho_\eps\in \cD_0(D)$ is defined by $\rho_\eps(x)=\int_{D} \eps^{-d}\ph(\eps^{-1}(x-z))\indic{d(z,\partial D)>2\eps} \rho(\dd z)$ for some smooth positive function $\ph$ supported in the unit ball of $\R^d$ with $\int \ph(y) \dd y =1$, then $\Gamma_0(\rho_\eps) \to \Gamma_0(\rho)$ and also $\int_D g(x) \rho_\eps(\dd x)\to \int_D g(x) \rho(\dd x)$ as $\eps\to 0$. Hence, it suffices to prove the inequality for $\rho(\dd x)=\rho(x)  \dd x$ with $\rho\in \cD_0(D)$. In this case, we have that
	\begin{align*}
		\int_D g(x) \rho(x) \dd x & = \int_D \int_D G_0^D(x,y)(-\Delta g(y)) \dd y \rho(x) \dd x \\
		& = \int_D (-\Delta g(y)) f(y) \dd y
	\end{align*}
where $f$ is defined by \eqref{eqn:frho} and satisfies $\Delta f = -\rho$. Applying Gauss--Green, we see that this is equal to $\int_D \nabla g(y) \nabla f(y) \dd y$, whose modulus is bounded above by the square root of $\int_D|\nabla g(y)|^2 \dd y \int_D |\nabla f(y)|^2 \dd y$ using Cauchy--Schwarz. Since $\int_D |\nabla f(y)|^2 \dd y = \Gamma_0(\rho)$ by \cref{L:definite}, this concludes the proof.
\end{proof}

\subsection{Reminders about function spaces}
\label{S:functions}

As we have already mentioned, one drawback of defining the GFF as a stochastic process is that we cannot realise $(h,\rho)$ for all $\rho\in \mathfrak{M}_0$ \emph{simultaneously}. For example, it will not always be possible to define $(h,\rho)$ when $\rho\in \mathfrak{M}_0$ is random.

With this in mind, it is often useful to work with {versions} of the GFF that almost surely live in some ``function'' space. For example, it turns out to be possible to define a version of the GFF that is a random variable taking values in the space of distributions, or generalized functions. In fact, versions of the GFF taking values in much nicer Sobolev spaces (with negative index) can also be defined.

For completeness we include some brief reminders on function spaces here. We continue to assume that $D$ is a regular domain, unless stated otherwise.

\begin{definition}[Dirichlet inner product]\ind{Sobolev space}
\ind{Dirichlet inner product}
\label{D:DIP}
We define the Dirichlet inner product
\begin{equation}\label{E:DE}
(f,g)_\nabla : = \int_D \nabla f(x) \cdot \nabla g(x) \dd x
\end{equation}\indN{{\bf Inner products}!$(\cdot,\cdot)_\nabla$; Dirichlet energy} 
for $f, g \in \cD_0(D)$. It is straightforward to see that $(\cdot, \cdot)_\nabla$ is a valid inner product. 
\end{definition}

\begin{definition}[The space $H_0^1$]\label{D:Sobolev1} \ind{Sobolev space}
We define the space $H_0^1(D)$ to be the completion of $\cD_0(D)$ with respect to the Dirichlet inner product.\indN{{\bf Function spaces}! $H_0^1(D)$; Sobolev space, completion of $\cD_0(D)$ with respect to the Dirichlet inner product}
 \end{definition}

By definition $H_0^1(D)$ is a separable Hilbert space with inner product $(\cdot, \cdot)_\nabla$. 

\begin{rmk}\label{R:DH0} Observe that since any element of $H_0^1(D)$ corresponds, by definition, to (the limit of) a Cauchy sequence of functions $f_n \in \cD_0(D)$ with respect to the Dirichlet inner product, it can be identified with a distribution $f\in \cD_0'(D)$ via $f(\rho):=\lim_{n\to \infty} f_n(\rho):=\int_D f_n(x) \rho(\dd x)$ for each $\rho\in \cD_0(D)$. In fact, due to Lemma \ref{L:M0dual}, this limit exists whenever $\rho\in \mathfrak{M}_0$. In this case we also have $|f(\rho)|\le (f,f)_\nabla \Gamma_0(\rho)$. \end{rmk}

\begin{rmk}	
	For general $D$, the standard definition of the Sobolev space $H_0^1(D)$ (see for example \cite{adams}) is the completion of $\cD_0(D)$ with respect to the inner product $(f,g):=(f,g)_{L^2(D)}+(f,g)_\nabla$. When $D$ is bounded, this coincides with \cref{D:Sobolev1}; indeed, by the Poincar\'{e} inequality, the norms $\|u\|:=(u,u)$ and $\|u\|_{\nabla}=(u,u)_{\nabla}$ are equivalent in this case.
\end{rmk}

\paragraph{Eigenbasis of $H_0^1(D)$.} When $D$ is bounded, it is easy to find a suitable orthonormal eigenbasis for $H_0^1(D)$. Indeed in this case, $H_0^1(D)$ is compactly embedded in $L^2(D)$ \indN{{\bf Function spaces}! $L^2(D)$; square integrable functions in $D$} by Rellich's embedding theorem, which implies that the resolvent of minus the Laplacian with Dirichlet boundary conditions is a compact operator. Note that this does not require any assumption of smoothness on the boundary of $D$. Consequently, there exists an orthonormal basis $(f_n)_{n\ge 1}$ of eigenfunctions of $-  \Delta$ on $D$, with zero (Dirichlet) boundary conditions, having eigenvalues $(\lambda_n)_{n\ge 1}$. That is, $f_n,\lambda_n$ satisfy
$$ \begin{cases}
- \Delta f_n  = \lambda_n f_n & \text{ in $D$}\\
f_n  = 0 & \text{ on $\partial D$}
\end{cases}
$$
for each $n\ge 1$. The $(\lambda_n)_{n\ge 1}$ are positive, ordered in non-decreasing order and $\lambda_n \to \infty$ as $n\to \infty$. Moreover the Gauss--Green formula \eqref{ipp} implies that for $\lambda_n \neq \lambda_m$,
$$
(f_n, f_m)_\nabla =  \lambda_m \int_D f_n f_m = \lambda_n \int_D f_n f_m.
$$
Hence $(f_n, f_m)_{\nabla} = 0$, and the eigenfunctions corresponding to different eigenvalues are orthogonal with respect to $(\cdot, \cdot)_\nabla$. %Indeed, elementary spectral theory of operators (going slightly beyond the scope of these notes) tells us that $(f_n)_{n\ge 1}$, properly normalised, form an orthonormal basis of $H_0^1(D)$.
%Moreover, $f_j = \lambda_j^{-1/2} e_j$, where $e_j$ are normalised so that $\|e_j \|_2 = 1$. Hence, an alternative description of $H_0^1$ is as the set of infinite series of the form $\sum_{n} \alpha_n f_n$, where $\sum_{n\ge 1} \alpha_n^2 < \infty$.

Often, the eigenfunctions of $-\Delta$ are normalised to have unit $L^2$ norm, %(rather than unit norm in $H_0^1(D)$)
since they also form an orthogonal basis of $L^2(D)$ for the standard $L^2$ inner product (again by the Gauss--Green formula). If $(e_j)_{j\ge1}$ are normalised in this way, then the above considerations imply that setting
\begin{equation}\label{fjej}
f_j =\frac{e_j}{\sqrt{\lambda_j}}
\end{equation}
for each $j$, we get an orthonormal basis $(f_j)_j$ of $H_0^1(D)$.

In particular, $f\in L^2(D)$ is an element of $H_0^1(D)$ if and only if \begin{equation}\label{E:nabla1} (f,f)_{\nabla} = \sum_{j\ge 1} (f,f_j)_\nabla^2=\sum_{j\ge 1} \lambda_j (f,e_j)_{L^2(D)}^2 <\infty.\end{equation}

\paragraph{Sobolev spaces of general index $H_0^s(D), s \in \R$.}

The above leads us to define $H_{0}^s$\indN{{\bf Function spaces}! $H_0^s(D)$; Sobolev space of index $s$ in $D$}  for general $s\in \R$, and bounded $D$, to be the Hilbert space completion of $\cD_0(D)$ with respect to the inner product
\begin{equation}\label{E:Hsinnerprod}
	(f,g)_s = \sum_{j\ge 1} \lambda_j^{s}(f,e_j)_{L^2(D)}(g,e_j)_{L^2(D)}.
	\end{equation}\indN{{\bf Inner products}! $(\cdot, \cdot)_s$; $H_0^s$ inner product}
Note that the above series does converge for $f,g\in \cD_0(D)$: this can be seen by applying Cauchy--Schwarz, using that $\cD_0(D)\subset L^2(D)$, and that all derivatives of functions in $\cD_0(D)$ are again elements $\cD_0(D)$, with $(\Delta f,e_n)=-\lambda_n(f,e_n)$ for $f\in \cD_0(D)$ and $n\ge 1$. We have also seen in \eqref{E:nabla1} that it agrees with the previous definition of $H_0^1(D)$ when $s=1$. \\

Let us make a few more straightforward observations.
\begin{itemize}
	\item When $s=0$ the above space is equivalent, by definition, to $L^2(D)$.
\item In general, when $s\ge 0$, it is simple to check that $L^2(D)\supset H^s_0(D)$, and that $f\in L^2(D)$ is an element of $H^s_0(D)$ if and only if $\sum_{j\ge 1} \lambda_j^s (f,e_j)_{L^2(D)}^2 <\infty$.

\item If $s\le 0$, then an element of $H_0^s(D)$ is by definition the limit of a sequence $\{f_n\}_n\in \cD_0(D)$ for which $\sum_{j\ge 1} \lambda_j^s (f_n,e_j)^2_{L^2(D)}$ has a limit as $n\to \infty$. In particular, for any $\phi \in H_0^{-s}(D)$, $\lim_{n\to \infty} (f_n,\phi)_{L^2}=:f(\phi)$ exists by Cauchy--Schwarz, and we can identify our element of $H_0^s(D)$ with the distribution $f\in \cD_0(D)'$, $\phi\mapsto f(\phi)$. Moreover, this distribution $f$ extends to a continuous linear functional on $H^{-s}(D)$. 

  In summary: for $s\le 0$, $H_0^s(D)$ can be identified with a subspace of $\cD_0'(D)$, and is the dual space\footnote{The space $H_0^{s}(D)$ for $s<0$ is usually referred to in the literature as simply $H^{s}(D)$, but we use the notation $H_0^{s}$ to emphasise that it is the dual of $H_0^{-s}(D)$ rather than $H^{-s}(D)$. When $-s\in \mathbb{Z}_{\ge 0}$, the latter is the space of $L^2$ functions with $|s|$ derivatives in $L^2(D)$ and is a strict superspace of $H_0^{-s}(D)$.} of $H_0^{-s}(D)$.
\item It is also clear from the above that convergence in any negative index Sobolev space implies convergence in the space of distributions $\cD_0'(D)$.
\end{itemize}

It will be useful in what follows to rephrase the expression for $\var (\mathbf{h},\rho)$ (when $\mathbf{h}$ is a GFF) in terms of Sobolev norms. Recall that by \eqref{definite}, if $-\Delta f =  \rho$ for $f, \rho \in \cD_0(D)$, then
  \begin{equation}
  \Gamma_0 (\rho) = (f,f)_{\nabla} = \sum_{j\ge 1} \lambda_j (f,e_j)^2_{L^2(D)}.
  \end{equation}
On the other hand, by Gauss--Green we have that $(\rho,e_j)_{L^2(D)}=-\lambda_j (f,e_j)_{L^2(D)}$ for every $j$, so that $(\rho,\rho)_{-1}=\sum_{j\ge 1} \lambda_j^{-1} (\lambda_j (f,e_j)_{L^2(D)} )^2 = (f,f)_{\nabla}$. 
In other words:
\begin{lemma}\label{L:M0Sobnorm} Suppose that $D\subset \R^n$ is bounded and $\rho\in \cD_0(D)$. Then 
\begin{equation}\label{variancenorm}
\var (\mathbf{h}, \rho) = \Gamma_0(\rho)=(\rho,\rho)_{-1}
\end{equation}
\end{lemma}
%From this it is easy to check that if $\rho $ is a measure in $\mathfrak{M}_0$, then in fact $\rho \in H_0^{-1} (D)$. In particular, if $\rho \in \mathfrak{M}_0$ we can find a sequence of smooth functions $\rho_n\in \cD_0(D)$ such that $\rho_n \to \rho$ in $H_0^{-1}(D)$. For this sequence, we then have $\var (\mathbf{h}, \rho_n - \rho) \to 0$, so that $(\mathbf{h}, \rho_n) \to (\mathbf{h}, \rho)$ in $L^2(\P)$.

\subsection{GFF as a random distribution}\label{ZGFFd}

 At this stage we do not yet know that the GFF may be viewed as a random distribution (that is, as a random variable in $\cD_0'(D)$). The goal of this section will be to prove that such a representation exists. Guided by \eqref{definite} (and by Theorem \ref{T:discretefourier}) we will find an expression for the GFF as a random series, which we will show converges in the distribution space $\cD_0'(D)$. In fact, we will show that it converges in a Sobolev space of appropriate index.

The property \eqref{variancenorm} suggests that $\mathbf{h}$ is formally the canonical Gaussian random variable ``in'' the dual space to $H_0^{-1}(D)$, that is, in $H_0^1(D)$ (the quotation marks are added since in fact $\mathbf{h}$ does not live in $H_0^1(D)$). It should thus have the expansion
\begin{equation}\label{GFFseries}
\mathbf{h}= \sum_{n=1}^\infty X_n g_n = \lim_{N\to \infty} \sum_{n=1}^N X_n g_n,
\end{equation}
where $X_n$ are i.i.d. standard Gaussian random variables and $(g_n)_{n\ge 1}$ is an {arbitrary} orthonormal basis of $H_0^1(D)$. (See for example \cite{Janson} for more about the general theory of Gaussian Hilbert spaces, and associated series such as the one above).

 It is not clear at this point in what sense (if any) this series converges. We will see in  \cref{T:GFFseriesSob} below that when $D$ is bounded, it converges in an appropriate Sobolev space and hence in the space of distributions.
 Note however that the series does \textbf{not} converge almost surely in $H_0^1(D)$, since the $H_0^1$ norms of the partial sums tend to infinity almost surely as $N\to \infty$ (by the law of large numbers).

We start with the following observation, where now $D$ can be any open set with regular boundary. Set $\mathbf{h}_N := \sum_{n=1}^N X_n g_n$, and let $f\in \cD_0(D)$ or more generally let $f \in H_0^1(D)$. Then
\begin{equation}\label{GFFconv}
(\mathbf{h}_N,f)_\nabla = \sum_{n=1}^N X_n (g_n, f)_\nabla
\end{equation}
does converge almost surely and in $L^2(\P)$,  by the martingale convergence theorem.  Its limit is a Gaussian random variable with variance $\sum_{n\ge 1} (g_n,f)_\nabla^2 = \|f\|_\nabla^2$ by Parseval's identity. This defines a random variable which we call $(\mathbf{h},f)_\nabla$, which has the law of a mean zero Gaussian random variable with variance $\|f\|_\nabla^2$. Hence while the series \eqref{GFFseries} does \emph{not} converge in $H_0^1$, when we take the inner product with a given $f \in H_0^1$ then this does converge almost surely.

\begin{comment}
Let $\rho \in \cD_0(D)$, and let $f \in \cD_0(D)$ be such that $-\Delta f = \rho$. Then $(\mathbf{h}_N, \rho) = (\mathbf{h}_N, f)_\nabla$ converges almost surely. This limit then defines a stochastic process indexed by smooth test functions, which has the same law as a GFF restricted to $\cD_0(D)$.
%An eigenvalue calculation can be used to show that the series \eqref{GFFseries} converges in the space of distributions (and in fact, in a nice Sobolev space, known as $H^{-\eps}_0(D)$, for any $\eps>0$. See Exercise \ref{Ex:Sobolev}).  \index{Sobolev space} Therefore, for any given smooth test function $\rho \in \cD_0(D)$, $\sum_{n=1}^\infty X_n (f_n , \rho)$ converges almost surely and defines a stochastic process (indexed by $\cD_0(D)$) whose law coincides with the restriction of the GFF to $\cD_0(D)$. Hence the restriction of the stochastic process $h$ to $\cD_0(D)$ has a modification which is a random distribution.
Using a density argument we will now show the following:
%Furthermore, a density argument (see below for details) can be used to show that the process $\tilde h$ defined in this way can be extended to $\mathfrak{M}_0$, where for a given $\rho \in \mathfrak{M}_0$, $(\tilde h,\rho)$ is defined as the limit in probability (and in $L^2$) of $(\tilde h, \rho_\eps)$ for smooth approximations $\rho_\eps$ of $\rho$.
% In fact, the series \eqref{GFFseries} can be used directly to concretely evaluate the Gaussian free field $(h, \rho)$, as follows.
\end{comment}

By a density argument, we can extend this to the following theorem. We use the notation $(f,\ph)$ for the action of a distribution $f$ on a smooth function $\ph$. 

\begin{theorem}[GFF as a random Fourier series]\label{GF}
Let $D$ be a regular domain and let $\mathbf{h}_N = \sum_{n=1}^N X_n g_n$ be the truncated series in \eqref{GFFseries}. Then for any $\rho \in \mathfrak{M}_0$,
$$\lim_{N\to \infty}
(\mathbf{h}_N, \rho)=:(\mathbf{h},\rho) $$
exists in $L^2(\P)$ (and hence in probability as well). The limit $(\mathbf{h},\rho)$ is a Gaussian random variable with variance $\Gamma_0(\rho,\rho)$.
\end{theorem}
Observe that since $\mathbf{h}_N$ is an element of $H_0^1(D)$, $(\mathbf{h}_N,\rho)$ is well defined for every $N$ by \cref{R:DH0}.

\begin{proof} We will first show that for any $\nu\in \mathfrak{M}_0$, we have the upper bound
	\begin{equation}\label{ineqvargen} \var(\mathbf{h}_N,\nu)\le \Gamma_0(\nu)
		\end{equation}
	for all $N\ge 1$. 
To see this,  recall from (the argument of) Lemma \ref{L:M0dual} that for any $\nu\in \mathfrak{M}_0$, there exists a sequence $\nu_k\in \cD_0(D)$ with $\Gamma(\nu_k)\to \Gamma(\nu)$ as $k\to \infty$. Furthermore, for this sequence it holds by \cref{R:DH0} that for each fixed $N$, $\var(\mathbf{h}_N,\nu_k)=\sum_{n=1}^N (g_n,\nu_k)^2 \to \sum_{n=1}^N (g_n,\nu)^2 = \var(\mathbf{h}_N,\nu)$ as $k\to \infty$. Finally, if we define $f_k$ such that $-\Delta f_k=\nu_k$ for each $k$, then the discussion just above implies that $\var(\mathbf{h}_N,\nu_k)=\var(\mathbf{h}_N,f_k)_\nabla \le \var(\mathbf{h},f_k)_\nabla = (f,f)_\nabla = \Gamma_0(\nu_k)$ for each $k$. Combining these observations gives the upper bound for any $\nu\in \mathfrak{M}_0$ and $N\ge 1$
$$ \var(\mathbf{h}_N,\nu)=\lim_{k\to \infty} \var(\mathbf{h}_N,\nu_k) \le \lim_{k\to \infty} \Gamma_0(\nu_k)=\Gamma_0(\nu),$$
as desired.

We will now use this to prove the result. Take $\rho \in \mathfrak{M}_0$ and choose a sequence $\rho_\eps\in \cD_0(D)$ approximating $\rho$ in the sense that $\Gamma_0(\rho_\eps)\to \Gamma_0(\rho)$ (again using  \cref{L:M0dual}). Set $\nu_\eps = \rho - \rho_\eps$ for each $\eps$. Then $\Gamma_0(\nu_\eps)\to 0$, so that applying \eqref{ineqvargen} to $\nu = \nu_\eps$ we deduce that $(\mathbf{h}_N, \nu_\eps)$ converges to 0 in $L^2(\P)$ and in probability as $\eps \to 0$, uniformly in $N$. The result then follows since for smooth $\rho_\eps$, we already know (as a consequence of the martingale convergence argument before the statement of the theorem) that $(\mathbf{h}_N, \rho_\eps) $ converges to a limit in $L^2$, and that this limit has the same law as $(\mathbf{h}, \rho_\eps)$.
\end{proof}

\begin{comment}One immediate consequence of this is the following\footnote{This makes the definition taken in \cite{Sheffield} precise.}:
\begin{cor}\label{D:GFFalternate}
The restriction of a Gaussian free field $h$ with zero boundary conditions to $\cD_0(D)$ is the unique stochastic process such that for all $f \in \cD_0(D)$, $(\mathbf{h},f)_\nabla$ is a Gaussian centred variable with variance $(f,f)_\nabla$, where we define $(\mathbf{h}, f)_\nabla$ by $- (\mathbf{h}, \Delta f)$ for each $f$.
\end{cor}
\end{comment}

We finally address convergence of the series \eqref{GFFseries}:

\begin{theorem}[GFF as a random variable in a Sobolev space]
  \label{T:GFFseriesSob}
  Suppose $D$ is a regular, bounded domain. If $(X_n)_{n\ge 1}$ are i.i.d. standard Gaussian random variables and $(g_n)_{n\ge 1}$ is \textbf{any} orthonormal basis of $H_0^1(D)$, then the series $\sum_{n\ge 1} X_n g_n$ converges almost surely in $H_0^{s}(D)$, where $$s = 1-\tfrac{d}{2} - \eps,$$ for any $\eps>0$. In particular, for $d=2$, the series converges in $H_0^{-\eps}(D)$ for any $\eps>0$.
\end{theorem}

Observe that by \cref{GF}, the law of the limit is uniquely defined, and coincides with the Gaussian free field $\mathbf{h}$ when its index set is restricted to $H_0^{-s}(D)$.
\begin{proof}

Let us take $(e_m)_{m\ge 1}$ to be an orthonormal basis of $L^2(D)$ which are eigenfunctions for $-\Delta$, as in \cref{S:functions}. This is possible since $D$ is bounded. As usual we write $\lambda_m$ for the eigenvalue corresponding to $e_m$; so that $(\lambda_m^{-s/2}e_m)_{m\ge 1}$ is an orthonormal basis of $H^{s}_0(D)$ and $(\lambda_m^{-1/2}e_m)_{m\ge 1}$ is an orthonormal basis of $H_0^1(D)$.
	In some cases $\lambda_m$ can be computed explicitly: for example, when $D$ is a rectangle or the unit disc. In general, we will make use of the following fundamental estimate due to Weyl (see for example  \cite[VI.4, page 155]{Chavel} for a proof):\ind{Weyl law} 

\begin{lemma}\label{lem:weyl}
We have $$\lambda_m \sim c m^{2/d}$$
as $m \to \infty$, in the sense that the ratio of the two sides tends to 1 as $m \to \infty$, where $c = {(2\pi)^2 / ( a_d\Leb(D))^{2/d}}$, where $a_d$ is the volume of the unit ball in $\R^d$.
\end{lemma}
\ind{Weyl law}

Now let $(g_n)_{n\ge 1}$ be \emph{any} orthonormal basis of $H_0^1(D)$. We start by observing that 
\begin{equation}
\label{eqn:exp_fs}
\E(\|\sum_{n=1}^N X_n g_n\|^2_{H_0^{s}})=\E(\sum_{n=1}^N \sum_{m=1}^N X_nX_m(g_n,g_m)_{H_0^{s}})=\sum_{n=1}^N \|g_n\|_{H_0^{s}}^2 .
\end{equation}
In fact, by the same argument, if $n_0\ge 1$ is fixed and $N\ge n_0$ then, setting $S_N = \|\sum_{n={n_0}}^N X_n g_n\|^2_{H_0^{s}}$, we see that $(S_N)_{N\ge n_0}$ is a submartingale (with respect to the natural filtration generated by $(X_n)_{n\ge 1}$.)
Furthermore, by applying Parseval's identity, we have that
\begin{eqnarray*}
 \sum_{n\ge 1} \|g_n\|^2_{H_0^{s}} & =  &\sum_{n\ge 1} \sum_{m\ge1} (g_n,{\lambda_m}^{-s/2}e_m)^2_{H_0^{s}} \\ &=&\sum_{m\ge 1} {\lambda_m}^{-1+s}\sum_{n\ge 1} (g_n,{\lambda_m}^{-1/2}e_m)^2_{H_0^1}\\ &=&\sum_{m\ge 1} {\lambda_m}^{-1+s}<\infty
\end{eqnarray*}
where we have used positivity and Fubini to interchange the order of summation. The finiteness of the last sum follows by \cref{lem:weyl}, since $(2/d)(-1+s)=-1-\eps(2/d)<-1$. Thus $(S_N)_{N\ge n_0}$ is bounded in $L^1(\P)$. In fact, 
\begin{align*}
\E(S_N^2) &= \sum_{n=n_0}^N \E(X_n^4) \| g_n\|_{H_0^s}^4 + 3\sum_{n_0\le m \neq n \le N} \E(X_n^2 X_m^2) \|g_n\|_{H_0^s}^2 \|g_m\|^2_{H_0^s}\\
& = 3 \sum_{n=n_0}^N  \| g_n\|_{H_0^s}^4  + 3\sum_{n_0\le m \neq n \le N}  \|g_n\|_{H_0^s}^2 \|g_m\|^2_{H_0^s}\\
& \le 3 \left( \sum_{n_0\le n \le N}  \| g_n\|_{H_0^s}^2\right)^2,
\end{align*}
so $(S_N)_{N\ge n_0}$ is bounded in $L^2(\P)$. Furthermore, by Doob's maximal inequality, and monotone convergence, for every $t>0$,
\begin{equation}\label{eq:Doobmax}
\P( \sup_{N\ge n_0} S_N \ge t) \le \frac4{t^2} \sup_{N\ge n_0} \E(S_N^2)   \le \frac{12}{t^2} \left( \sum_{n=n_0 }^\infty  \| g_n\|_{H_0^s}^2\right)^2.
\end{equation}
We now explain how this implies almost sure convergence of the series $\sum_{n\ge1} X_n g_n$ in $H_0^s(D)$.  Set  $U(n_0) = \sup_{N\ge n_0} \| \sum_{n=m}^N X_n g_n \|_{H^s}^2$, which is the random variable in the left hand side of \eqref{eq:Doobmax}. 

Fix $k \ge 1$. Since $\sum_{n \ge 1}  \|g_n \|^2_{H^s} < \infty$, there exists $n_k\ge 1$ such that 
$\sum_{n \ge n_k}  \|g_n \|^2_{H^s} \le 4^{-k}.$
 By \eqref{eq:Doobmax},
$$
\P [ U(n_k) \ge 2^{-k}  ] \le 12 \cdot 4^{-k}.
$$
As the right hand side above is summable, we deduce that $U(n_k) \le 2^{-k}$ for $k$ sufficiently large and hence $U(n_k) \to 0$ almost surely. On the other hand if $n \ge n_k$, since 
$
\| a+ b \|^2 \le 2 ( \| a\|^2 + \|b \|^2)
$
on any Hilbert space, we deduce that if $n \ge n_k$ then  $U(n) \le 4 U(n_k)$. Thus $U(n) \to 0$ almost surely. Hence $\sum_{n=1}^N X_n g_n$ is Cauchy in $H_0^s$ and therefore converges. This completes the proof of the theorem. 
\end{proof}

\begin{rmk}\label{R:conv_GFF_nonbounded}
	The above theorem implies that the series $\sum_{n = 1}^N X_n g_n$ converges almost surely in the space of distributions $\cD_0'(D)$ whenever $D$ is bounded. 
	Recall that the measurability of the event 
	$$
	E= \{ \omega \in \Omega: \sum_{n = 1}^N X_n(\omega) g_n \text{ converges in the space } \cD'_0(D) \}
	$$
	is provided by Lemma \ref{L:measurability}. However, even if we do not appeal to this lemma, the statement `` the series $\sum_{n = 1}^N X_n g_n$ converges almost surely in the space of distributions $\cD_0'(D)$ whenever $D$ is bounded'' would still be meaningful. 
	 Indeed, in Theorem \ref{T:GFFseriesSob} we have checked that this series converges almost surely in the space $H^s_0(D)$ for some $s<0$ (an event which is clearly measurable since $H^s_0(D)$ is a metric and indeed Hilbert space). On that (measurable) event of probability one, say $E_s$, it is clear that convergence in the space of distribution holds. Thus $E \supseteq E_s$ where $E_s$ has probability one. Another way to state this is that, given  Theorem \ref{T:GFFseriesSob} (and independently of Lemma \ref{L:measurability}), the event $E$ is measurable on the completed $\sigma$-field $\cF^*$ of the probability space (the completed $\sigma$-field $\cF^*$ is the $\sigma$-field generated by $\cF$ and the null sets). 
	
	Furthermore, by \cref{GF}, this means that the GFF as a stochastic process, when its index set is restricted to smooth test functions, has a version that is almost surely a random element of $\cD_0'(D)$. 
	Moreover in two dimensions, the boundedness assumption can be removed using conformal invariance, see Theorem \ref{thm:CIDGFF}.
%
%	
%	allows us to drop the boundedness assumption. 
	%Indeed, take $D$ any simply connected domain of $\C$, and let $\phi$ be a conformal map from $\D\to D$. Then by conformal invariance of the Dirichlet inner product, if $(g_n)_{n\ge 1}$ is an orthonormal basis of $H_0^1(D)$, the images $(g_n\circ \phi)_{n\ge 1}$ form an orthonormal basis of $H_0^1(\D)$. So if $(X_n)_{n\ge 1}$ are i.i.d.\ standard Gaussians, we have that $\sum_{n\ge 1} X_n (g_n\circ \phi)$ converges a.s.\ in $\cD_0'(\D)$, and consequently that $\sum_{n\ge 1} X_n g_n$ converges in $\cD_0'(D)$. Conformal invariance will be discussed again in more detail in Theorem \ref{thm:CIDGFF}.
\end{rmk}

Let us reiterate one of the important conclusions from \cref{T:GFFseriesSob,GF}.
\begin{cor}
\label{C:ZGFF_sob}
Suppose that $D\subset \R^d$ is a bounded domain, and $s=1-\tfrac{d}{2}-\eps$ for some $\eps>0$. Then there exists a version of the Dirichlet GFF as a stochastic process  $(\mathbf{h},\rho)_{\rho\in H_0^{-s}}$ (with restricted index set) that is almost surely an element of $H^s_0(D)$.
\end{cor}

\subsection{It\=o's isometry for the GFF}\label{S:ito}

\dis{This section will not be used in the rest of the text and the reader may wish to skip it on a first reading.}

In this section we describe an observation which emerged from joint discussions with James Norris.
It is closely linked to \cref{L:M0Sobnorm}, which implies that for a zero boundary GFF $\mathbf{h}$ in a bounded domain $D$, and for any $f\in H_0^{-1}(D)$, the quantity $(\mathbf{h},f)$ makes sense almost surely. That is, as the almost sure (and $L^2(\P)$) limit of $(\mathbf{h},f_n)$ for any sequence $f_n$ converging to $f$ in $H_0^{-1}(D)$.

In other words, even though $h$ is only almost surely defined as a continuous linear functional on $H^{d/2-1+\eps}(D)$ for $\eps>0$ (\cref{T:GFFseriesSob}), we can actually test it against fixed functions that are much less regular. Namely, we can test it against any fixed function  $H_0^{-1}(D)$. Note that this agrees with (in fact slightly extends) our previous definition of $h$ as a stochastic process, since we have seen that $\mathfrak{M}_0$ is precisely the set of signed measures that are elements of $H_0^{-1}$, a consequence of Lemma \ref{L:M0Sobnorm}.

In this section we will essentially formulate the above discussion in terms of an isometry. To motivate this, it is useful to recall the following well known analogy within It\=o's theory of stochastic integration. Let $B$ be a standard Brownian motion. Even though $dB$ does not have the regularity of a function in $L^2$ (in fact, it is essentially an element of $H^{-1/2 - \eps}$ for any $\eps>0$), it makes perfect sense to integrate it against a test function in $L^2$. This is thanks to the fact that the map
$$
f  \mapsto \int f_s \dd B_s
$$
defines an isometry of suitable Hilbert spaces. Thus much flexibility has been gained: \emph{a priori} we don't even have the right to integrate against functions in $H^{1/2}$, and yet, taking advantage of some almost sure properties of Brownian motion -- namely, quadratic variation -- it is possible to integrate against functions in $L^2$ (and actually much more).

 A similar gain  can be seen in the context of the GFF: \emph{a priori}, as an element of $H_0^{1-d/2-\eps}$ ($\eps>0$), it would seem that integrating against an arbitrary test function $f \in L^2$ is not even allowed when $d\ge 2$. %(we would need functions $f \in H^\eps$ for some $\eps>0$ in view of Lemma \ref{L:cont}).
 Yet, as discussed above, we can almost surely integrate against much rougher objects, namely distributions in $H_0^{-1}$:

 \begin{theorem}[It\=o isometry]
   \label{T:ito}
%To a function $f \in \cD_0(D)$ associate the random variable $X_f = (h,f)$.
The map $X$ sending $f\in \cD_0(D)$ to the random variable $X_f=(\mathbf{h},f)$ can be viewed as a linear map between $\cD_0(D)$ and the set of random variables $L^2(\Omega, \cF, \P)$ viewed as a function space. If we endow $\cD_0(D)$ with the $H_0^{-1}(D)$ norm and $L^2(\Omega, \cF, \P)$ with its $L^2$ norm then $X$ is an isometry:
$$
\| f\|_{H_0^{-1}(D)} = \|X_f\|_2 = \E( (\mathbf{h}, f)^2)^{1/2}.
$$
In particular, since $\cD_0(D)$ is dense in $H_0^{-1}(D)$, $X_f$ extends uniquely as an isometry from $H_0^{-1}(D)$ into $L^2(\Omega, \cF, \P)$. Hence if $f \in H_0^{-1}(D)$, then we can set $(\mathbf{h},f)$ to be the unique limit in $L^2(\P)$ of $(\mathbf{h}, f_n)$ where $f_n$ is any sequence of test functions that converge in $H_0^{-1}(D)$ to $f$.
    \end{theorem}

\begin{proof}
%  We first check that $\| f \|_{\gff}$ defines a norm on $\cD_0(D)$ as claimed: the triangle inequality and the scaling are trivial, and for the positivity observe that if $\| f\|_{\gff}^2 = 0$ then $f = 0$
%    almost everywhere (and then everywhere by continuity) by \eqref{Covariance_definite}.
  This is a direct consequence of \cref{L:M0Sobnorm}.
  %?}  For $f \in \cD_0(D) \subset L^2(D)$, let us write $f  = \sum \alpha_n e_n $ where $e_n$ is an orthonormal basis (for $L^2(D)$) of eigenfunctions, and $\alpha_n = (f, e_n)$. Then
  %  $$
  %  \|X_f\|_{2}^2 =\E((h,f)^2) = \E ( ( \sum_n X_n (f_n, f))^2 ) = \sum_n (f_n, f)^2 = \sum_n (e_n, f)^2 \frac{2\pi}{\lambda_n} = \|f\|^2_{H_0^{-1}(D)}
  %  $$
  %  by definition of the inner product on $H_0^{-1}(D)$ (see \eqref{Sobolev}). This finishes the proof of the theorem.
\end{proof}

\begin{rmk}
	Note that although $(\mathbf{h},f)$ makes sense as an almost sure limit for any fixed $f\in H_0^{-1}(D)$, or indeed for any countable collection of such $f$, this does not mean that $\mathbf{h}$ is an element of $H^{1}_0$ or that we can test $\mathbf{h}$ against every element of $H_0^{-1}$ simultaneously.  For example, writing \begin{equation*}\mathbf{h}=\lim_{N \to \infty} \mathbf{h}_N:=\lim_{N \to \infty}\sum_{n=1}^N \frac{X_n}{ \sqrt{\lambda_n}} e_n\end{equation*} with $(X_n)_n\sim \cN(0,1)$ i.i.d.\ and $(e_n)_n$ an orthonormal basis of Laplacian eigenfunctions for $L^2(D)$, we have $\mathbf{h}_N\to \mathbf{h}$ almost surely in $H^{-1}_0(D)$ but $\var(\mathbf{h},\mathbf{h}_N)=\sum_{n=1}^N \lambda_n^{-1} \to \infty $ (at least when $d\ge 2$). So there do exist random elements of $H_0^{-1}(D)$ that cannot be tested against $\mathbf{h}$.
\end{rmk}

\begin{comment}We end this discussion by remarking that defining the Gaussian free field through this It\=o isometry would be equivalent to \cref{D:GFFc} since
$$
\|f\|^2_{H^{-1}_0 (D)} = 2 \pi ( - \Delta^{-1/2} f, - \Delta^{-1/2} f) =
-2 \pi (f, \Delta^{-1} f) = (f, G_0f) = \Gamma_0(f,f)
$$
so in particular any measure in $\mathfrak{M}_0$ is in $H_0^{-1}(D)$. (In fact, the above computation immediately tells us that $H_0^{-1}(D)$ is just the completion of $\mathfrak{M}_0$ with respect to the norm $\iint \rho(dx) \rho(dy) G_0(x,y) $).
\end{comment}

\subsection{Cameron--Martin space of the Dirichlet GFF}
\label{sec:ac}

\dis{This section will not be used until \cref{S:SIsurfaces} and the reader may wish to skip it on a first reading.}

In this section, we will address the following question:
\begin{itemize} \item for $\mathbf{h}$ a Dirichlet (zero) boundary condition GFF in $D$ and $F$ a (deterministic) function on $D$, when are $\mathbf{h}+F$ and $\mathbf{h}$ mutually absolutely continuous?\footnote{as stochastic processes indexed by $\mathfrak{M}_0$.}
\end{itemize}

The answer is that this holds whenever $F\in H_0^1(D)$. This question can be phrased for general Gaussian processes, and the space of $``F"$ for which absolute continuity holds is known as the \textbf{Cameron--Martin space} of the process. Thus, the lemma below says that $H_0^1(D)$ is the Cameron--Martin space of the (Dirichlet boundary condition) GFF.

\begin{prop}\label{lem:CMGFF} Let $\mathbf{h}$ be a GFF in a bounded domain $D$ with Dirichlet (zero) boundary conditions. Then $\mathbf{h}$ and $\mathbf{h}+F$ are mutually absolutely continuous, as stochastic processes indexed by $\mathfrak{M}_0$, if and only if $F\in H_0^1(D)$. When this holds, the Radon--Nikodym derivative of $(\mathbf{h}+F)$ with respect to $\mathbf{h}$ is given by
	\[ \frac{\exp((\mathbf{h},F)_\nabla)}{\exp((F,F)_\nabla/2)}.\]
\end{prop}

\begin{proof}
	%By conformal invariance, we may assume that $D=\D$.
	Let $(e_i)_i$ be an orthonormal basis of $L^2(D)$ consisting of eigenfunctions of the Laplacian, with associated eigenvalues $(\lambda_i)_i$. We write $g_i:=(\sqrt{\lambda_i})e_i$, so that the $(g_i)_i$ form an orthonormal basis of $H_0^{-1}\supset \mathfrak{M}_0$. Recall that $$F\in H_0^{1}(D) \Leftrightarrow \sum (F,g_i)^2  <\infty. $$
	For $n\in \N$, we consider the finite vector $((\mathbf{h},g_i))_{1\le i \le n}$, which by definition of the GFF is just a vector of independent $\cN(0,1)$ random variables.
	
	This is convenient to work with because of the following elementary fact: if $(X_i)_{1\le i \le n}$ are i.i.d. standard normals and $(a_i)_{1\le i \le n}$ are real numbers, then  $(X_1,X_2,...,X_n)$ and $(X_1+a_1, X_2+a_2,... X_n+a_n)$ are mutually absolutely continuous. Moreover, the RN derivative of the latter with respect to the former is given by $e^{\sum a_i X_i}/e^{\sum a_i^2/2}$.
	
	In our context, this means that the law of  $((\mathbf{h},g_i))_{1\le i \le n}$ is mutually absolutely continuous with that of $((\mathbf{h}+F, g_i))_{1\le i \le n}$, if and only if $|(F,g_i)|< \infty$ for $1\le i \le n$. Furthermore when this does hold, the Radon--Nikodym derivative of $((\mathbf{h}+F, g_i))_{1\le i \le n}$ with respect to $((\mathbf{h},g_i))_{1\le i \le n}$ is equal to
	\begin{equation}\label{eq:approx_RN} \frac{\exp(\sum_{i=1}^n (\mathbf{h},g_i)(F,g_i))}{\E(\exp(\sum_{i=1}^n (\mathbf{h},g_i)(F,g_i)))} = \frac{\exp(\sum_{i=1}^n (\mathbf{h},g_i)(F,g_i))}{\exp(\sum_{i=1}^n (F,g_i)^2/2)}.
	\end{equation}
	
Now, for $\mathbf{h}$ and $\mathbf{h}+F$ to be mutually absolutely continuous, the family of random variables on the right hand side of \eqref{eq:approx_RN} must be uniformly integrable (in $n$). Indeed, they should be the conditional expectations, with respect to a family of sub $\sigma$-algebras, of the Radon--Nikodym derivative of $(\mathbf{h}+F)$ with respect to $\mathbf{h}$.
This family is \emph{not} uniformly integrable if $F\notin H_0^1(D)$, that is, $\sum_{i\ge 1} (F,g_i)^2 = \infty$. Hence we obtain the necessity of the condition  $F\in H_0^1(D)$ in the proposition.
	
%As a result, the stochastic processes $((h,\sqrt{\lambda_i}e_i))_{i\ge 1}$ and $((h+f,\sqrt{\lambda_i}e_i))_{i\ge 1}$ are mutually absolutely continuous if and only if the right hand side above converges in $L^1$. That is, if and only if $\sum_{i=1}^n \lambda_i (f,e_i)^2$ converges as $n\to \infty$, which holds if and only if $f\in H_0^1(D)$. When this does hold, the Radon-Nikodym derivative of the latter process with respect to the former is given by this $L^1$ limit, which is equal to
%$$\frac{e^{\sum_{i\ge 1} \lambda_i (h,e_i)(f,e_i)}}{e^{\sum_{i\ge 1} \lambda_i (f,e_i)^2}} = \frac{e^{(h,f)_\nabla}}{e^{(f,f)_\nabla}}.$$

%To conclude, we need to deduce the same result for the stochastic processes $(h,\rho)_{\rho\in \mathfrak{M}_0}$ and $(h+f,\rho)_{\rho\in \mathfrak{M}_0}$. It is clear that these cannot be mutually absolutely continuous if the stochastic processes $((h,\sqrt{\lambda_i}e_i))_{i\ge 1}$ and $((h+f,\sqrt{\lambda_i}e_i))_{i\ge 1}$ are not (since this corresponds to a restriction of the index set). Hence, we see that the condition $f\in H_0^1(D)$ in the proposition is necessary.

For the sufficiency, we observe that when $F\in H_0^1(D)$, the random variables on the right hand side of \eqref{eq:approx_RN} converge in $L^1(\P)$ to
$$\frac{\exp(\sum_{i\ge 1}  (\mathbf{h},g_i)(F,g_i))}{\exp(\sum_{i\ge 1} (F,g_i)^2/2)} = \frac{\exp((\mathbf{h},F)_\nabla)}{\exp((F,F)_\nabla/2)}$$
as $n\to \infty$. We also know by \cref{GF} that  whenever $\rho\in \mathfrak{M}_0$, $\sum_{i=1}^n  \lambda_i^{-1}(\rho,g_i)(\mathbf{h},g_i)$ converges to $(\mathbf{h},\rho)$ almost surely. This implies that for any $\rho_1,...,\rho_m\in \mathfrak{M}_0$ and any $\psi:\R^m\to \R$ continuous and bounded:
\[ \E \left(\psi((\mathbf{h}+F,\rho_1),...,(\mathbf{h}+F,\rho_m))\right)=\E\left( \frac{\exp((\mathbf{h},F)_\nabla)}{\exp((F,F)_\nabla/2)} \psi((\mathbf{h},\rho_1),...,(\mathbf{h},\rho_m))\right).\]
But this is exactly the statement that $\mathbf{h}+F$ is absolutely continuous with respect to $\mathbf{h}$, as a stochastic process indexed by $\mathfrak{M}_0$, with the desired Radon--Nikodym derivative. Since the inverse of the Radon--Nikodym derivative is also in $L^1$, we obtain the mutual absolute continuity. \end{proof}

\subsection{Markov property}
\label{SS:Markov}
We are now ready to state one of the main properties of the GFF, which is the (domain) Markov property. As in the discrete case, informally speaking, it states that conditionally on the values of $\mathbf{h}$ outside of a given subset $U$, the free field inside $U$ is obtained by harmonically extending $\mathbf{h}|_{D\setminus U}$ into $U$ and then adding  an independent GFF with Dirichlet boundary conditions in $U$. %to a function which is the harmonic extension inside $U$ of the values of $h$ on its boundary.
Note that in this case, however, it is not at all clear that such a harmonic extension is well defined. \ind{Markov property}

\begin{figure}
\begin{center}
\includegraphics[width = .4\textwidth]{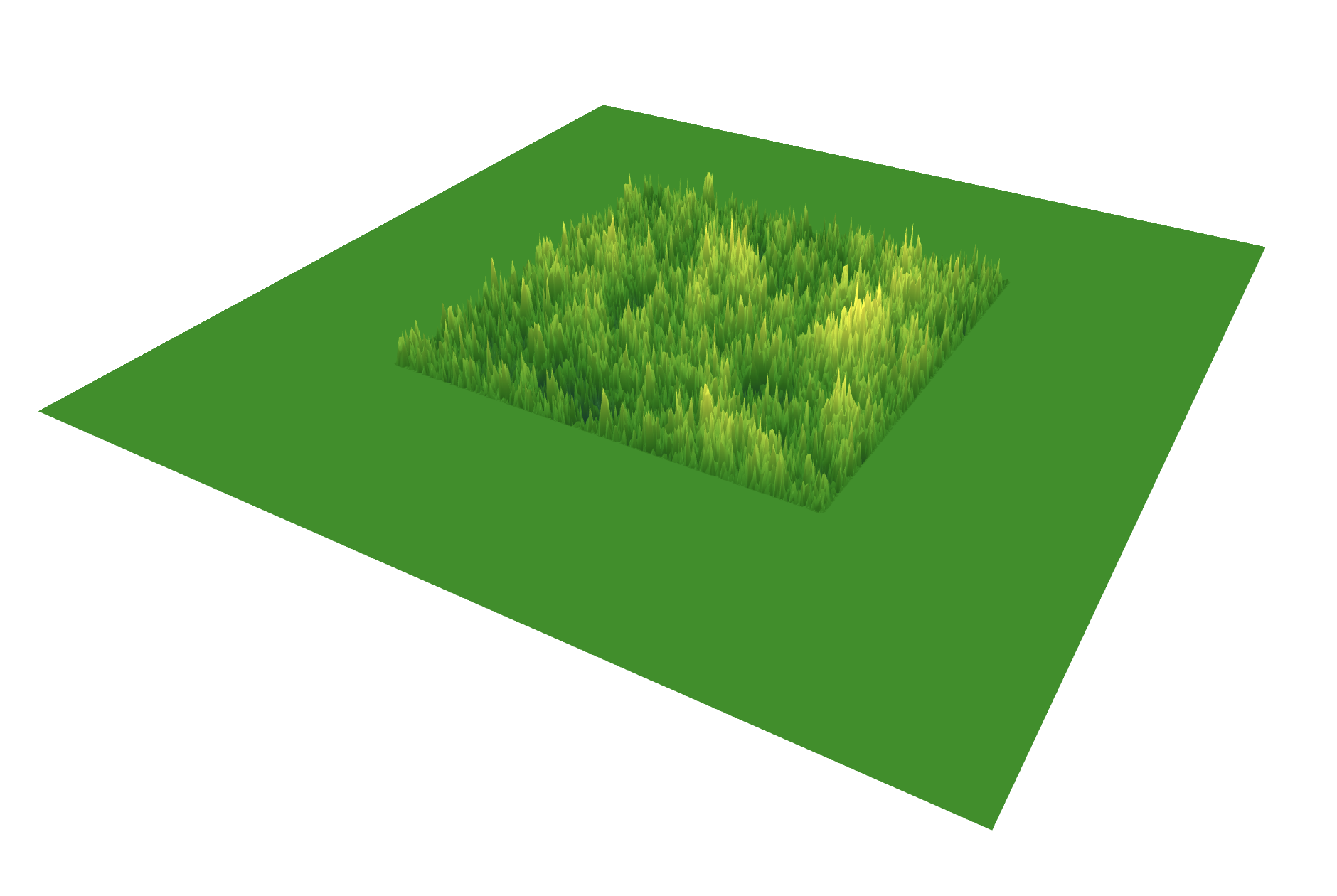}
\includegraphics[width = .4\textwidth]{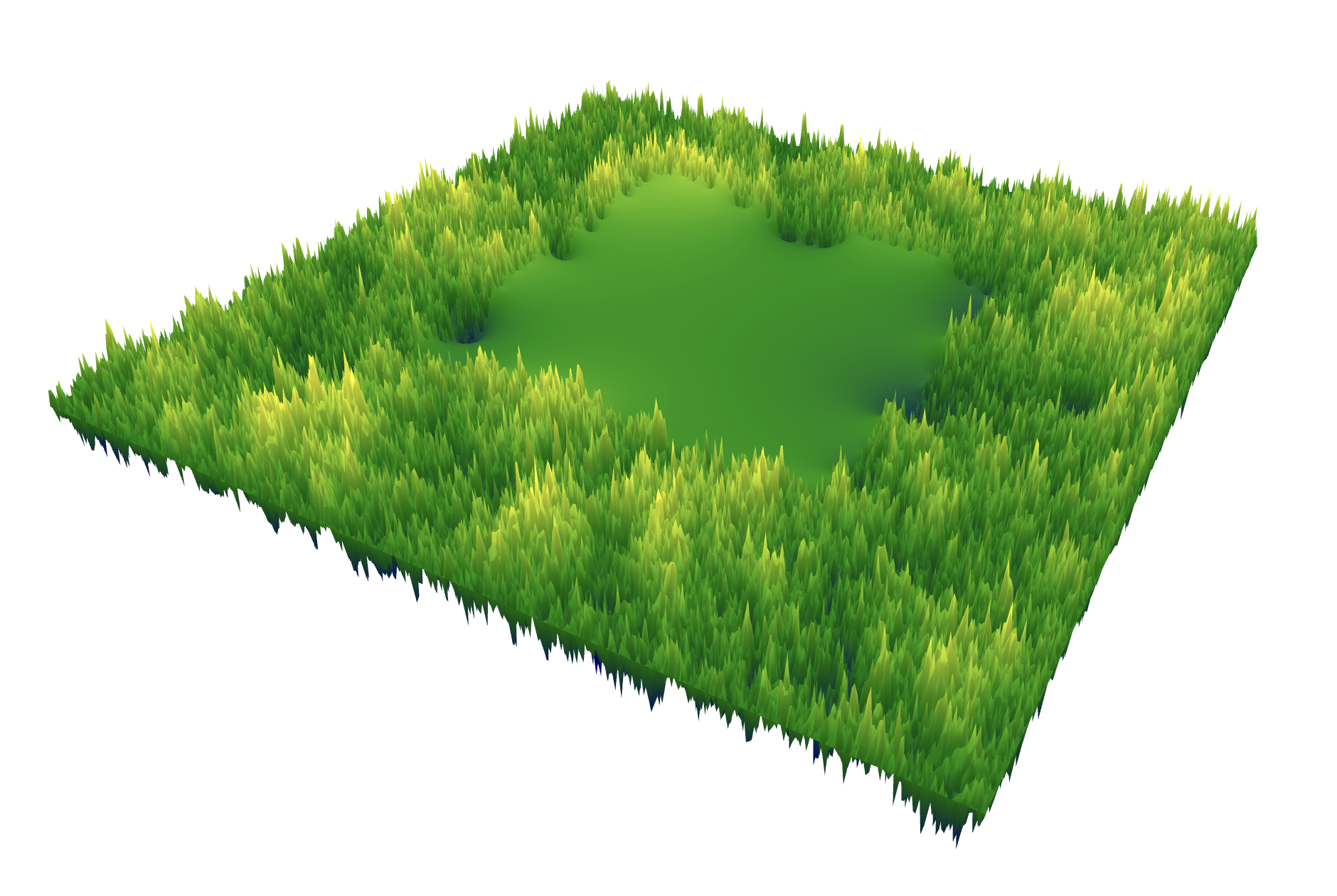}
\caption{The Markovian decomposition of the GFF: here $D$ is a square, and $U\subset D$ a slightly smaller square. The first graph shows $\mathbf{h}_0$, and the second shows $\varphi$. Their sum $\mathbf{h}$ is a GFF in $D$, shown in Figure \ref{F:dGFF}. Simulations by Oskar-Laurin Koiner. \label{F:Markov}
}
\end{center}
\end{figure}

\begin{theorem}[Markov property]\label{T:mp} Fix $U \subset D$ a regular subdomain. Let $\mathbf{h}$ be a GFF (with zero boundary conditions on $D$). Then we may write
$$
\mathbf{h} = \mathbf{h}_0 + \varphi,
$$
where:
\begin{enumerate}
\item $\mathbf{h}_0$ is a zero boundary condition GFF in $U$, and is zero outside of $U$;

\item $\varphi$ is harmonic in $U$; and

\item $\mathbf{h}_0$ and $\varphi$ are independent.
\end{enumerate}
\end{theorem}
This makes sense whether we view $\mathbf{h}$ as a random distribution or a stochastic process indexed by $\mathfrak{M}_0$. %Note also that we do not require $U$ to be connected. 
Note that since $\mathbf{h}_0 = 0$ on $U^c$, $\ph$ coincides with $\mathbf{h}$ on $U^c$. See Figure \ref{F:Markov} for an illustration.

\begin{cor}\label{C:H1loc}
	By \cref{R:conv_GFF_nonbounded}, when $D\subset \R^2$ is an arbitrary domain (that is, potentially unbounded) this Markov property implies that the random distribution $\mathbf{h}$ almost surely defines a random element of the local Sobolev space $H^{-1}_{\mathrm{loc}}(D)$.\footnote{$H^{-1}_{\mathrm{loc}}(D)$ is the space of distributions whose restriction to any $U\Subset D$ (that is, such that $\bar{U}$ is a compact subset of $D$) is an element of $H^{-1}_0(U)$.}\indN{{\bf Function spaces}! $H^{-1}_{\mathrm{loc}}(D)$;  distributions in element of $H^{-1}_0(U)$ for any $U\Subset D$}  In general dimension $d\ge 3$ it follows from the above Markov property that there is a version of the stochastic process $\mathbf{h}$ that is almost surely a random distribution (in fact, a random element of $H^{1-d/2-\eps}_{\mathrm{loc}}(D)$ for any $\eps>0$).
\end{cor}
\begin{proof}
The key point is the following Hilbertian decomposition:

\begin{lemma}\label{lem:markovprop_od}
Let $U$ be as in \cref{T:mp}. We have $$H_0^1(D) = \Supp(U) \oplus \Harm(U),$$
 where $\Harm(U)$\indN{{\bf Function spaces}! $\Harm(U)$; harmonic functions in $U$} consists of harmonic functions in $U$ (that is, elements of $H_0^1 (D)$ whose restriction to $U$ coincide with a harmonic function in $U$). % and $\Supp(U) \equiv H_0^1(U)$ consists of functions supported in $U$ (the closure with respect to the Dirichlet inner product of compactly supported smooth functions in $U$).
\end{lemma}

\begin{proof}
We first prove orthogonality. Let $f \in \Supp(U)\subset H_0^1(D)$. Then there exists $f_n \in \cD_0(U)\subset \cD_0(D)$ such that $f_n \to f $ in the $H_0^1(D)$ sense. Now let $g \in H_0^1(D)$ such that $g$ coincides in $U$ with a harmonic function. Note that
$$
(f_n, g)_\nabla = \int_D (\nabla f_n) \cdot (\nabla g) = \int_{U_n} (\nabla f_n) \cdot (\nabla g)
$$
where $U_n\subset U$ is chosen to be compactly contained in $U$, have smooth boundary, and contain the (closure of the) support of $f_n$ (in particular $\nabla f_n = 0$ outside of $U_n$).
Since $U_n$ is smooth, we can apply the Gauss--Green \ref{L:ipp} formula in $U_n$ with boundary term; because $g$ is non-zero on $\partial U_n$, this boundary term does need to be considered. However, the boundary term vanishes because $\partial g/ \partial n$ is a smooth function on $\partial U_n$ and $f_n = 0 $ on $\partial U_n$.

Therefore $(f_n, g)_\nabla = - \int_{U_n} f_n \Delta g$,
which is clearly 0 because in $U_n$, $\Delta g = 0$ and $f_n$ is a smooth function. This shows that $f_n$ and $g$ are orthogonal in $H_0^1(D)$. Then by taking a limit (since $f_n$ approximates $f$ in the $H_0^1(D)$ sense), $f$ and $g$ must also be orthogonal.

Now let us show that the sum of the two spaces spans $H_0^1(D)$. Let us suppose to begin with that $U$ is $C^1$ smooth. For $f \in H_0^1(D)$, let $f_0$ denote the orthogonal projection of $f$ onto $\Supp(U)$. Set $\varphi = f- f_0$: our aim is to show that $\varphi$ is harmonic in $U$.
\ind{Gauss--Green formula}
%A property of orthogonal projections gives that $\|f-f_0\|_\nabla^2 = \inf_{ g \in \Supp(U) } \| f- g\|_\nabla^2$. Now,
%\begin{align*}
% \inf_{ g \in \Supp(U) } \| f- g\|_\nabla^2 & = \inf_{\psi:V\to \R, \psi|_{U^c} = f|_{U^c} } \int_U |\nabla\psi |^2
%\end{align*}
%But it is easy to check that the minimiser, if the inf is attained (and it is attained, at a unique point, by properties of the orthogonal projection), must be attained at a harmonic function in $U$ (indeed, consider the Dirichlet energy of perturbations of the form $\psi + \eps u$, where $u$ is an arbitrary test function in $U$, and do a Taylor expansion as $\eps \to 0$ using the Gauss Green formula). Hence $\varphi$ is a harmonic function in $U$, and we are done.
Note that $\varphi $ is (by definition) orthogonal to $\Supp(U)$. Hence for any test function $\psi\in \cD_0(U)$, we have that $(\varphi, \psi)_\nabla = 0$. By the Gauss--Green formula (and since $U$ is $C^1$ smooth), we deduce that
$$
\int_D (\Delta \varphi) \psi = \int_U (\Delta \varphi) \psi = 0
$$
and hence $\Delta \varphi = 0$ as a distribution in $U$. Elliptic regularity arguments (going beyond the scope of these notes) show that a distribution which is harmonic in the sense of distributions must in fact be a smooth function, harmonic in the usual sense. Therefore $\varphi \in \Harm(U)$ and we are done.

If $U$ does not have $C^1$ boundary, let $(U_n)_{n\in \N}$ be a sequence of increasing open subsets of $U$ with $C^1$ boundaries, such that $\cup U_n = U$. For $f\in H_0^1(D)$, by the previous paragraph, we can write $f=f_0^n+\varphi^n$ for each $n\in \N$, where $f_0^n$ is the projection of $f$ onto $\Supp(U_n)$ and $\varphi^n\in \Harm(U_n)$. Then we just need to show that: (a) $f_0^n\to f_0$ as $n\to \infty$ for some $f_0\in \Supp(U)$; and (b) that $f-f_0$ is harmonic in $U$. For (a), we observe that $\Supp(U)=\overline{\cup_n \Supp(U_n)}$ (by definition of $\Supp(U)$ as the closure of $\cD_0(U)$ with respect to the Dirichlet inner product) and so the projections $f_0^n$ of $f$ onto $\Supp(U_n)$ converge, with respect to $\|\cdot \|_\nabla$, to $f_0\in \Supp(U)$. For (b), notice that by definition of $f_0$, $f-f_0$ is the limit of $\varphi_n$ as $n\to \infty$, with respect to $\|\cdot \|_\nabla$. In particular, it is clear that when restricted to any $U_n$, $f-f_0=\lim_n \varphi_n$ is harmonic in the distributional sense, and thus harmonic by elliptic regularity. Since this holds for any $n$, it follows that $f-f_0$ is harmonic in $U$.
%\nb{Actually orthogonality also needs to be proved in this case.}
\end{proof}

Having this decomposition in hand, we may deduce the Markov property in a rather straightforward way. Indeed, let $(f^0_n)_n$ be an orthonormal basis of $\Supp(U)$, and let $(\phi_n)_n$ be an orthonormal basis of $\Harm(U)$. For $((X_n, Y_n))_n$  an i.i.d. sequence of independent standard Gaussian random variables, set $\mathbf{h}_0 = \sum_n X_n f_n^0$ and $\varphi = \sum_n Y_n \phi_n$. Then the first series converges in $\cD_0'(D)$ since it is a series of a GFF in $U$. The sum of the two series gives $\mathbf{h}$ by construction, and so the second series also converges in the space of distributions. In the space of distributions, the limit of harmonic distributions must be harmonic as a distribution, and hence (by the same elliptic regularity arguments as above) a true harmonic function. This proves the theorem.
\end{proof}

\begin{rmk}It is worth pointing out an important message from the proof above: any orthogonal decomposition of $H_0^1(D)$ gives rise to a decomposition of the GFF into independent summands. \ind{Orthogonal decomposition of $H_0^1(D)$}
%\red{Exercise}
\end{rmk}

\noindent \textbf{Example}. When $d=1$. this is the statement that if $(b_s)_{s\in [0,1]}$ is a Brownian bridge from $0$ to $0$ and $[a,b]\subset [0,1]$, then conditionally on $(b_s)_{s\in [0,a]\cup[b,1]}$, the law of $(b_s)_{s\in [a,b]}$ is given by a the linear interpolation of $b_a$ and $b_b$, plus an independent Brownian bridge from $0$ to $0$ on $[a,b]$.

\begin{rmk}
	In the case when $D$ is an unbounded domain of $\R^d$ with $d\ne 2$, applying the Markov property in bounded subdomains shows that the GFF, viewed as a stochastic process with restricted index set $\cD_0(D)$, has a version that almost surely defines a distribution on $D$.
\end{rmk}

\subsection{Conformal invariance}

In the remainder of this chapter, we restrict ourselves to dimension $d=2$.

\medskip In this case the GFF possesses the important additional property of \textbf{conformal invariance}, which follows almost immediately from the construction in the previous section.
Indeed, a straightforward change of variable formula shows that the Dirichlet inner product is conformally invariant: if $\varphi: D \to D'$ is a conformal isomorphism, then
$$
\int_{D'} \nabla (f \circ \varphi^{-1}) \cdot \nabla (g \circ \varphi^{-1})  = \int_D \nabla f \cdot \nabla g.
$$
Consequently, if $(f_n)_n$ is an orthonormal basis of $H_0^1(D)$, then $(f_n \circ \varphi^{-1})_n$ defines an orthonormal basis of $H_0^1(D')$. (Watch out however, that eigenfunctions of $-\Delta$ are not conformally invariant in any sense).
{So by \cref{T:GFFseriesSob}:}

\begin{theorem}[Conformal invariance of the GFF]\label{thm:CIDGFF}
If $\mathbf{h}$ is a random distribution on $\cD_0'(D)$ with the law of the Gaussian free field on $D$, then the distribution $\mathbf{h}\circ \varphi^{-1}$, defined by setting $(\mathbf{h}\circ \varphi^{-1}, f)=(\mathbf{h},|\varphi'|^2( f\circ \varphi))$ for $f\in \cD_0'(D')$, has the law of a GFF on $D'$.
\end{theorem}
\ind{Conformal invariance! of Dirichlet GFF}

\medskip Recently, a kind of converse was shown in \cite{BPR20,BPR}: if a field $\mathbf{h}$ with zero boundary conditions satisfies conformal invariance and the domain Markov property, as well as a moment condition ($\E((\mathbf{h}, \phi)^{1+\eps}) < \infty$ for some $\eps>0$ and all $\phi \in \cD_0(D)$), then $\mathbf{h}$ must be a multiple of the Gaussian free field.  In fact, one can reduce the conformal invariance assumption to scale invariance, and obtain the result in all dimensions, \cite{APchar}. See \cite{BPR20,BPR,APchar} for details.

\subsection{Circle averages}

An important tool for studying the GFF is the process which describes its average values on small circles centred around a point $z\in D$. This is known as the \textbf{circle average process} around $z$.

More precisely, fix $z \in D$ and let $0<\eps< \dist(z, \partial D)$. Let $\rho_{z,\eps}$\indN{{\bf Miscellaneous}! $\rho_{z,\eps}$; uniform distribution on the circle of radius $\eps$ around $z$} denote the uniform distribution on the circle of radius $\eps$ around $z$. Note that $\rho_{z,\eps} \in \mathfrak{M}_0$. This follows from the fact that $G_0^D(x,y) \le - (2\pi)^{-1} \log |x- y| + O(1)$, and the fact, when we fix one of the variables $x$ on the circle, the integral over the circle of $-\log | x- y|$ with respect to $y$ is finite (just like the integral of $-\log r $ with respect to $r$ is finite in one dimension). More generally, this argument shows that the Lebesgue measure on any smooth curve is an element of $\mathfrak{M}_0$.

We set $\mathbf{h}_\eps(z) = (\mathbf{h}, \rho_{z,\eps}).$
The following theorem, is a consequence of the Kolmogorov--{\v C}entsov continuity theorem (a multidimensional generalisation of the more classical Kolmogorov continuity criterion), and  will not be proved here.  The interested reader is directed to Proposition 3.1 of \cite{DuplantierSheffield} for a proof. \ind{Circle average}

\begin{prop}[Circle average is jointly H\"older]
\label{T:contCA}
There exists a modification of $\mathbf{h}$ such that $(\mathbf{h}_{\eps}(z), z\in D, 0< \eps < \dist(z, \partial D))$ is  almost surely jointly H\"older continuous of order $\eta < 1/2$ on all compact subsets of $\{z\in D \text{ s.t. } 0< \eps < \dist(z, \partial D)\}$.
\end{prop}
In fact it can be shown that this version of the GFF is the same as the version which turns $h$ into a random distribution in \cref{GF}.
%\red{Exercise}
The reason circle averages are so useful is because of the following result.

\begin{theorem}[Circle average is a Brownian motion]
\label{T:ca}
Let $\mathbf{h}$ be a GFF on $D$. Fix $z \in D$ and let $0<\eps_0 < \dist(z, \partial D)$. For  $ t \ge t_0 = \log (1/\eps_0)$, set
$$
B_t = \sqrt{2\pi}\mathbf{h}_{e^{-t}}(z).
$$
Then $(B_t, t \ge t_0)$ has the law of a Brownian motion started from $B_{t_0}$: in other words, $(B_{t+ t_0} - B_{t_0}, t \ge 0)$ is a standard Brownian motion.
\end{theorem}

\begin{proof}
In order to avoid factors of $\sqrt{2\pi}$ everywhere, we use $h = \sqrt{2\pi} \mathbf{h}$ as defined in \eqref{Dnormalisation}, and call $h_\eps(z) = ( h, \rho_{z, \eps})$. The theorem statement is then equivalent to saying that $(B_t= h_{e^{-t}}, t \ge t_0)$ is a Brownian motion starting from $B_{t_0}$. Various proofs can be given. For instance, the covariance function can be computed explicitly (this is a good exercise)! Alternatively, we can use the Markov property of the GFF to see that $B_t$ must have stationary and independent increments. Indeed, suppose $\eps_1 > \eps_2$, and we condition on $h$ outside of $B(z, \eps_1)$. That is, we write $h = h^0 + \varphi$, where $\varphi$ is harmonic in $U  = B(z, \eps_1)$ and $h^0$ is a GFF in $U$ that is independent of $(h_\eps(z))_{\eps\ge \eps_1}$ (scaled in the same manner as \eqref{Dnormalisation}). Then $h_{\eps_2}(z)$ is the sum of two terms:  $h^0_{\eps_2}(z)$; and the circle average of $\varphi$ on $\partial B(z, \eps_2)$. By harmonicity of $\varphi$ the latter is nothing else than $h_{\eps_1}(z)$. This gives that the increment can be expressed as
$$
h_{\eps_2}(z) - h_{\eps_1}(z) = h^0_{\eps_2}(z)
$$
and hence, since $h^0$ is independent of $(h_\eps(z))_{\eps\ge \eps_1}$, the increments are independent.
Moreover, by applying the change of scale $w \mapsto (w-z) /\eps_1$, so that the outer circle is mapped to the unit circle, we see that the distribution of $h_{\eps_2}(z)-h_{\eps_1}(z)$ depends only on $r=\eps_2/\eps_1$. This means that they are also stationary.

To show from here that $h_{e^{-t}}(z)$ is a Brownian motion, it suffices to compute its variance. That is (by the Markov property), to check that if $h$ is a GFF in the unit disc $\D$ and $r<1$, then $h_r(0)$ has variance $-\log r$.

For this, let $\rho$ denote the uniform distribution on the circle $\partial (r\D)$ at distance $r$ from the origin, so that
\begin{equation}\label{varcomp}
\var (h_r(0)) = 2\pi \int_{\D^2} G_0^\D (x,y) \rho(\dd x) \rho(\dd y).
\end{equation}
The point is that by harmonicity of $G_0^\D(x, \cdot)$ in $\D \setminus\{x\}$ and the mean value property, the above integral is simply
\begin{equation}\label{varcomp2}
\var (h_r(0)) = 2\pi \int_{\D} G_0^\D (x,0) \rho(\dd x),
\end{equation}
which completes the proof since $G_0^\D(x,0) = - (2\pi)^{-1}\log |x| = - (2\pi)^{-1}\log r$ on $\partial (r\D)$.

To check \eqref{varcomp2} rigorously, first consider for a fixed $\eta>0$, the double integral
$$I_\eta =2\pi  \int_{\D^2} G_0^\D \big((1+\eta)x,y\big) \rho(\dd x) \rho(\dd y).$$
Then $I_\eta $ converges clearly to the right hand side of \eqref{varcomp} as $\eta \to 0$, and it is now rigorous to exploit the mean value property for the harmonic function $G_0^{\D} ((1+\eta)x, \cdot)$ in the entire ball $B(0,r)$ to deduce that
$$
I_\eta =2\pi  \int_{\D} G_0^{\D}\big((1+\eta)x, 0\big) \rho (\dd x).
$$
Letting $\eta \to 0$ proves \eqref{varcomp2}.

%So $h_\eps(z)$ has independent and stationary increments, and is Gaussian with mean zero and $\var(h_{e^{-t}}(z))=t$ for all $t$.  This concludes the proof.
%. One still needs to compute the variance of $(B_{t_0+ t } - B_{t_0})$, but by scale invariance it has to be a multiple of $t$. The constant, as before, has to be identified by explicitly computing it through the Green function. The fact that it is 1 boils down to our normalisation of the Green function.
\end{proof}
So, as we ``zoom in'' towards a point, the average values of the field oscillate like those of a Brownian motion. This gives us a very precise sense in which the field cannot be defined pointwise.

\subsection{Thick points}

An important notion in the study of Liouville quantum gravity is that of thick points of the Gaussian free field. Indeed, although these points are atypical from the point of view of Euclidean geometry, we will see that they are typical from the point of view of the associated quantum geometry. In order to be consistent with its applications in Gaussian multiplicative chaos and Liouville quantum gravity, we will once again here mostly work with the normalisation $h = \sqrt{2\pi} \mathbf{h}$ from \eqref{Dnormalisation}.

\begin{definition}Let $\mathbf{h}$ be a GFF in $D\subset \C$ open and simply connected, let $h = \sqrt{2\pi} \mathbf{h}$, and let $\alpha >0$. We say a point $z\in D$ is $\alpha$-thick if
$$
\liminf_{\eps \to 0} \frac{h_\eps(z)}{\log (1/\eps) } = \alpha.
$$
\end{definition}
In fact, the $\liminf$ in the definition could be replaced with a $\limsup$ or $\lim$. It is also clear by symmetry that the set of $(-\alpha)$-thick points with $\alpha>0$ has the same law as the set of $\alpha$-thick points; hence we restrict to the case $\alpha>0$ for simplicity.  

Note that a given point $z\in D$ is almost surely not thick: the typical value of $h_\eps(z)$ is of order $\sqrt{\log 1/\eps}$ since $h_\eps(z)$ is a Brownian motion at scale $\log 1/\eps$. At this stage, the most relevant result is the following fact due to Hu, Miller and Peres \cite{HMP} (though it was independently and earlier proved by Kahane in the context of his work on Gaussian multiplicative chaos). \ind{Thick points}

\begin{theorem}\label{T:thick}
Let $\cT_\alpha$\indN{{\bf Miscellaneous}! $\cT_\alpha$; $\alpha$-thick points of the Gaussian free field} denote the set of $\alpha$-thick points.
Then almost surely, the Hausdorff dimension\indN{{\bf Miscellaneous}! $\dd_H$; Hausdorff dimension} $\dd_H(\cT_\alpha)$ of $\cT_\alpha$ satisfies 
$$
\dd_H(\cT_\alpha) = (2 - \frac{\alpha^2}2)_+
$$
and $\cT_\alpha$ is almost surely empty if $\alpha > 2$.
\end{theorem}

\begin{proof}[Heuristics] The value of the dimension of $\cT_\alpha$ is easy to understand and to guess. Indeed, for a given $\eps>0$,
\begin{align*}
\P( h_\eps( z) \ge \alpha \log (1/\eps))  & = \P( \cN(0, \log (1/\eps) + O(1) ) \ge \alpha \log (1/\eps))
\\
&= \P( \cN(0,1)  \ge \alpha \sqrt{\log (1/\eps) + O(1)} ) \le \eps^{\alpha^2/2}
\end{align*}
using scaling and the standard bound
$ \P( X>t ) \le \text{const}\times  t ^{-1} e^{-t^2/2}$
for $X \sim \cN(0,1)$. Suppose without loss of generality that $D = (0,1)^2$ is the unit square. Then the expected number of squares of size $\eps$ such that the centre $z$ satisfies $h_\eps(z) \ge \alpha \log 1/\eps$ is bounded by $\eps^{-2 + \alpha^2/2}$. This suggests that the Minkowski dimension is less or equal to $2 - \alpha^2/2$ when $\alpha<2$ and that $\cT_\alpha$ is empty if $\alpha > 2$.
\end{proof}

\begin{proof}[Rigourous proof of upper bound.] We now turn the above heuristics into a rigorous proof that $\dd_H (\cT_\alpha) \le (2 - \alpha^2 / 2) \vee 0$, which follows closely the argument given in \cite{HMP}. The lower bound given in \cite{HMP} is more complicated, but we will obtain an elementary proof in the next chapter, via the Liouville measure: see Exercise \ref{Ex:thick} of Chapter 2.

To start the proof of the upper bound, we begin by stating an improvement of \cref{T:contCA}, which is Proposition 2.1 in \cite{HMP}. This is the circle average analogue of L\'evy's modulus of continuity for Brownian motion.

\begin{lemma}
  \label{L:Levymod}
Suppose $D$ is bounded with smooth boundary. %contains the unit disc and let $U = B(0, 1/2)$. \ellen{Why is this condition here? So that $h_r(z)$ will be defined when $r<1$ and $z\in U$? Still not necessarily true... Other than that is it necessary?}
Then there exists a version of the circle average process $(h_r(z))_{r<1,z\in D}$, such that for every $\eta<1/2, \zeta>0$ and $\eps>0$, there exists $M = M(\eta, \zeta, \eps)$ which is finite almost surely and such that
$$
| h_{r}(z) - h_{s}(w) | \le M \left( \log \frac1r\right)^{\zeta}\frac{(|z-w| + |r-s|)^{\eta}}{r^{\eta+ \eps}}
$$
holds for every $z,w \in D$ and for all $r,s \in (0,1)$ such that $r/s \in [1/2, 2]$ and $B(z,r), B(w,s)\subset D$.
\end{lemma}
See Proposition 2.1 in \cite{HMP} for a proof.

\medskip
Without loss of generality, we will now work in the case where $D$ is bounded with smooth boundary. This yields the proof in the general case by the domain Markov property, and since $\dd_H(\mathcal{T}_\alpha)=\lim_{n \to \infty}\dd_H(\mathcal{T}_\alpha \cap D_n)$ for a sequence of smooth, bounded domains $D_n$ with $\cup D_n=D$.

In this setting, the above lemma allows us to ``discretise'' the set of $\eps$ and $z$ on which it suffices to check thickness. More precisely, set $\eps>0$, $K >0$ and consider the sequence of scales $r_n = n^{-K}$. Fix $\zeta<1$, and $\eta<1/2$ arbitrarily (say $\zeta = 1/2, \eta = 1/4$), and let $M = M(\eta, \zeta, \eps)$ be as in the lemma. %(clearly we can assume that $D$ contains the unit disc without loss of generality, and it suffices to compute the Hausdorff dimension of thick points within $U$ \ellen{(Why? And why do we need this assumption - apart from because it's in the lemma above?)}
Then for any $z \in D$,
we have that if $r_{n+1} \le r \le r_n$,
\begin{align*}
  |h_r(z) - h_{r_n}(z)| & \le M K^\zeta (\log n)^\zeta \frac{(r_{n+1} - r_n)^\eta}{r_n^{\eta(1 + \eps)}}\\
  & \lesssim( \log n)^\zeta n^{K\eta (1+\eps) - (K+1) \eta} \lesssim( \log n )^\zeta
\end{align*}
if we choose $\eps = K^{-1}$. Thus any point $z \in D$ is in $\cT_\alpha$ if and only if
$$
\lim_{n \to \infty} \frac{h_{r_n}(z)}{\log 1/r_n} = \alpha.
$$
Now for any $n \ge 1$, let $\{z_{n,j}\}_{j} = D \cap  r_{n}^{1+\eps} \Z^2$ be a set of discrete points spaced by $r_n^{1+\eps}$ within $D$. Then if $z \in B(z_{n,j}, r_n^{1+\eps})$ we have, for the same reasons,
$$
|h_{r_n}(z) - h_{r_n}(z_{n,j})|  \lesssim( \log n)^\zeta.
$$
Thus for fixed $\delta>0$, we let 
$$
\cI_n = \{j: h_{r_n}(z_{n,j}) \ge (\alpha - \delta) \log (1/r_n)\}.
$$
Then for each $N \ge 1$, and each $\delta >0$,
$$
\cT'_\alpha = \bigcup_{n > N} \bigcup_{j \in \cI_n} B(z_{n,j}, r_n^{1+\eps})
$$
is a cover of $\cT_\alpha$. Consequently, if $\cH_q$ denotes $q$ dimensional Hausdorff measure for $q>0$,
$$
\E( \cH_q (\cT_\alpha)) \le \E\left(\sum_{n > N} \sum_{j \in \cI_n} \text{diam} B(z_{n,j}, r_n^{1+\eps})^{q}\right) \lesssim \sum_{n >N}r_n^{ - 2 -2\eps} r_n^{q(1+\eps)} \max_j \P( j \in \cI_n).
$$
For a fixed $n$ and a fixed $j$, as argued in the heuristics, $$\P( j \in \cI_n) \lesssim \exp( - \tfrac{ (\alpha - \delta)^2}{2} \log (1/r_n)) = r_n^{(\alpha - \delta)^2/2}$$
where the implied constants are uniform over $D$. We deduce
$$
\E( \cH_q (\cT_\alpha) ) \le \sum_{n > N} r_n^{- 2 -2\eps  + (\alpha - \delta)^2/2 + q(1+\eps)}.
$$
As $r_n = n^{-K}$ and $K$ can be chosen arbitrarily large, the right hand side tends to zero as $N \to \infty$ as soon as the exponent of $r_n$ in the above sum is positive, or, in other words, if $q$ is such that
$$
- 2 -2\eps+ (\alpha - \delta)^2/2 +q (1+\eps) >0.
$$
Thus we deduce that $\cH_q( \cT_\alpha) = 0 $ almost surely (and hence $\dd_H (\cT_\alpha) \le q$ whenever $q(1+\eps)> 2 + 2\eps- (\alpha - \delta)^2/2$. So
$$
\dd_H (\cT_\alpha) \le \frac{2 + 2\eps- (\alpha - \delta)^2/2}{1+\eps},
$$
almost surely. Since $\eps>0, \delta> 0$ are arbitrary, we deduce
$$
\dd_H(\cT_\alpha) \le 2 - \alpha^2/2,
$$
almost surely, as desired.
  \end{proof}
%A slightly more elaborate argument shows that in fact the Minkowski dimension of $\cT_\alpha$ is a.s. less than $2- \alpha^2 / 2$ and therefore so is, a.s., the Hausdorff dimension. To prove the lower bound, one uses a fairly standard second moment argument, see \cite{HMP} for details. A simple alternative new proof of this fact is also proposed in Exercise \ref{Ex:thick} of Chapter 2.

The value $\alpha =2$ corresponds informally to the maximum of the free field, and the study of the set $\cT_2$ is, informally at least, related to the study of extremes in a branching Brownian motion (see \cite{ABBS, ArguinBovierKistler}).

\subsection{Scaling limit of the discrete GFF }\label{SS:scalinglimit}

In this short section we briefly explain why the discrete GFF on appropriate sequences of planar graphs converges in the scaling limit to the continuum GFF. Before we give general arguments, let us point out a situation in which this is relatively straightforward to see. 

Let $D = (0,1)^2$ be the unit square, and $V_N = D \cap (\Z^2/N)$ be the portion of the square lattice (scaled to have mesh size $1/N$) that intersects $D$, and let $E_N$ be the edges of the whole square lattice scaled by $1/N$. Let $\partial_N$ denote the set of vertices $v \in V_N$ with at least one neighbour outside of $V_N$, which is the natural boundary of this graph. Let $h_N$ be the discrete Gaussian free field associated with $V_N,\partial_N$ (and with $q_{x,y} = 1$ for every pair of neighbouring vertices $x,y$ in $V_N$).  In order to discuss convergence to the continuum GFF, it is useful to extend the definition of $h_N$ to all of $\R^2$: namely, we extend $h_N$ to be constant on each face of the dual graph of $(V_N, E_N)$; that is, for $x \in V_N$, and $y \in (-1/(2N), +1/(2N)]^2$, we set
$h_N(y) = h_N(x).$ 

We then claim that for a fixed $k\ge 1$ and fixed test functions $\phi_1, \ldots, \phi_k \in \cD_0(D)$, the law of the vector $(h_N, \phi_i)_{i=1}^k$ converges (without scaling) as $N\to \infty$ to the law of $(h, \phi_i)_{i=1}^k$, where $h$ is a continuum GFF. In fact, we will check the following stronger convergence.
\begin{prop} \label{P:GFFsquare}Consider the above discrete GFF $h_N$ in the unit square. We then have the convergence in distribution: 
\begin{equation}\label{eq:hCVdistr}
h_N \to  h
\end{equation}
 as random variables on $H_0^{s}(D)$ for any $s<0$, where $h$ is a continuum GFF with zero boundary conditions on $D = (0,1)^2$.
\end{prop}

Note that the choice of normalisation is consistent from discrete to continuum, in that the discrete random walk associated to the graph $G_N = (V_N,E_N)$ converges, after speeding up time by a factor $N^2$, to a \emph{speed two} Brownian motion. 

\begin{proof}
For $k, m\ge 1$ let 
$$
f_{k,m}(x,y) =2\sin ({\pi k x}) \sin(\pi m y). 
$$
It is elementary that $f_{k,m}$ is an eigenfunction of $-\Delta$ in $D=(0,1)^2$ (with Dirichlet boundary conditions), corresponding to the eigenvalue $\lambda_{k,m} =\pi^2( k^2 + m^2)$, and has unit $L^2(D)$ norm; an elementary fact from Fourier analysis is that $(f_{k,m})_{k,m\ge 1}$ form an orthonormal basis of $L^2(D)$. 

Furthermore, on the unit square a minor miracle happens: namely, if $1\le k \le N$ and $1\le m \le N$ then $f_{k,m}$ is \emph{also} a discrete eigenfunction of the negative discrete Laplacian $- Q_N$, with associated eigenvalue 
$$
\lambda^N_{k,m} = 2-2\cos \left( \frac{\pi k}{N}\right) + 2- 2\cos \left( \frac{\pi m}{N}\right).
$$
In particular, letting $N \to \infty$ but keeping $k,m\ge 1$ fixed, we see that 
$$
\lambda^N_{k,m} \sim \frac{1}{N^2} \lambda_{k,m}.
$$
We denote by $f^N_{k,m}$ the eigenfunction $f_{k,m}$, normalised to have unit (discrete) $L^2$ norm, that is, 
$$
f_{k,m}^N (\cdot) = \frac1{c^N_{m,k}} f_{k,m} (\cdot) ; \text{ with } c^N_{k,m} =\left( \sum_{z\in V_N} f_{k,m}(z)^2\right)^{1/2}. 
$$
Clearly, as $N \to \infty$ with $k,m\ge 1$ fixed,
$$
(c^N_{k,m})^2  \sim4 N^2 \iint_D \sin^2(\pi k x) \sin^2 (\pi m y) \dd x \dd y = N^2. 
$$
In fact, using simple trigonometric identities we can check that for all $N\ge 1$ and all $1\le k,m \le N$ we have $c_{k,m}^N =N$ exactly.

The functions $f_{k,m}$ are linearly independent and thus $(\lambda_{k,m}^N)_{1\le k \le N, 1\le m \le N}$ give us all possible eigenvalues of $-Q_N$ (counted with multiplicity in case of repetition). 

We deduce, using Theorem \ref{T:discretefourier}, that the discrete GFF can be written as
$$
h_N (\cdot) = \sum_{k,m=1}^N \frac1{\sqrt{\lambda_{k,m}^N} } X_{k,m} f^N_{k,m} ( \cdot) =  \sum_{k,m=1}^N \frac1{c_{k,m}^N\sqrt{\lambda_{k,m}^N} } X_{k,m} f_{k,m} ( \cdot) ,
$$ 
where $(X_{k,m})_{1\le k,m \le N}$ are independent standard Gaussian random variables. By Theorem \ref{T:GFFseriesSob}, to deduce \eqref{eq:hCVdistr}, it remains to check that 

\begin{itemize}
	\item 
when $k,m\ge 1$ are fixed and $N\to \infty$, $c^N_{k,m} \sqrt{\lambda_{k,m}^N} \to \pi \sqrt{\lambda_{k,m}}$; and 
\item the expected $H_0^s$ square norm of the remainder series is controlled uniformly in $N$.
\end{itemize}
The first point is elementary given the above asymptotics. For the second point, we need to show that for all $\eps>0$ we can find $A\ge 1$ large but fixed (that is, independent of $N$) such that, 
\begin{equation}\label{Eq:convGFFboildown}
\E \left( \| \sum_{k=1 }^N \sum_{m=A}^N   \frac1{c_{k,m}^N\sqrt{\lambda_{k,m}^N} } X_{k,m} f_{k,m} ( \cdot) \|_{H_0^s}^2\right) \le \eps
\end{equation}
for all $N\ge 1$. 
Note that the left hand side above is equal to 
$$
\sum_{k=1}^N \sum_{m=A}^N \frac1{(c_{k,m}^N)^2 \lambda_{k,m}^N } \lambda_{k,m}^s
$$
by definition of the $H_0^s$ norm. 
To conclude, we simply observe that $1- \cos (x) \ge a x^2/2$ for all $x \in [0,1]$ and some $a>0$. Hence 
$$
\lambda_{k,m}^N \ge \frac{a}{N^2} \lambda_{k,m},
$$
and since $c^N_{k,m} = N$, we see that 
$$
\frac{1}{ c_{k,m}^N \lambda_{k,m}^N } \lambda_{k,m}^s \le C \lambda_{k,m}^{-1+s}
$$ 
for some constant $C>0$. We conclude that \eqref{Eq:convGFFboildown} holds as in the proof of Theorem \ref{T:GFFseriesSob}, that is, using Weyl's law. This proves \eqref{eq:hCVdistr}. 
\end{proof}

If we take a general bounded domain $D\subset \R^2$, the argument above can no longer be applied because there is no exact relation between the discrete and continuous eigenfunctions. A different argument is therefore needed. 

We fix a bounded domain $D \subset \R^2$. Let $(G_\delta)_{\delta >0}$ denote a sequence of undirected graphs (with weights on the edges) embedded in $D$. Denote their vertex sets by $v(G_\delta)$ and prespecified boundaries by $\partial_\delta \subset v(G_\delta)$. Let $\P^{\delta}_x$ denote the law of continuous time random walk on $G_\delta$ starting from some vertex $x$ of $G_\delta$, killed when it reaches the boundary $\partial_\delta$, and let $\E^\delta_x$ denote the associated expectation. Our main assumption is that the random walk under $\mathbb{P}^\delta_x$ converges to (speed two) Brownian motion as $\delta\to 0$, uniformly on compact time intervals and uniformly in space, in the sense that for any smooth test function $\phi \in \cD_0(D)$ we have 
\begin{equation}\label{e:CVassumption}
\left|\E^\delta_x [ \phi ( X_{s\delta^{-2} } )] - \E_x [ \phi( B_s)1_{\{ \tau > s \} } ]\right| \to 0,
\end{equation}
as $\delta\to 0$, uniformly over $s\in [0,T]$ for every $T>0$ and $x\in V(G_\delta)$. We also suppose that if $\tau_\delta$ is the time that the random walk under $\mathbb{P}_x^\delta$ first hits the boundary $\partial_\delta$, then $\delta^2\tau_\delta$ is uniformly integrable: that is, for every $\eps>0$ we can choose $K<\infty$ such that 
\begin{equation}
\E^\delta_x ( {\delta^2}{\tau_\delta} 1_{\{ \tau_\delta  \ge K \delta^{-2}\} }) \le \eps, \label{eq:uniformtightness}
\end{equation}
uniformly in the vertices $x$ of $G_\delta$. Finally, we assume that the vertices of $G_\delta$ have density asymptotically uniform, in the following sense: for any open set $A $ such that $A \subset D$, 
\begin{equation}\label{eq:densityone}
\frac{\# v ( G_\delta) \cap A}{\delta^{-2} } \to \text{Leb} (A)
\end{equation}
as $\delta \to 0$. 
%If $\alpha$ is a function on the vertices of $G_\delta$, extend it to all of $\R^2$ as follows. Fix a proper embedding of the dual graph $G_\delta^*$ of $G_\delta$. For any vertex $v$ of $G_\delta$, let $f$ be the face of $G_\delta^*$ that contains $v$, and set $ \alpha$ to be constant, equal to $\alpha(v)$, on the face $f$. 

\begin{theorem}
\label{T:GFFscalinglimit} Let $h_\delta$ denote the discrete GFF associated to the graph $G_\delta$ and $\partial_\delta$, as above, and suppose that \eqref{e:CVassumption},  \eqref{eq:uniformtightness} and \eqref{eq:densityone} hold. For a test function $\phi \in \cD_0(D)$, let $h_\delta (\phi) = \delta^2\sum_{x\in v(G_\delta)} h_\delta(x) \phi(x)$. Then for every $k\ge 1$, and for every set of test functions $\phi_1, \ldots, \phi_k \in \cD_0(D)$, we have
$$
(h_\delta(\phi_i))_{i=1}^k \to ( h, \phi_i)_{i=1}^k
$$
in distribution as $\delta \to 0$, where $h$ is a continuum GFF with zero boundary conditions in $D$.
\end{theorem}

\begin{proof}
Since the variables $h_\delta(\phi_i)$ are Gaussian and linear in $\phi_i$, it once again suffices to prove the statement for $k = 1$, in which case we write $\phi$ instead of $\phi_1$. Having fixed $\phi$, we observe that $h_\delta(\phi)$ is a centred Gaussian random variable with variance
\begin{equation}\label{eq:convergencegoal}
\sigma_\delta^2 = \delta^4 \sum_{x,y\in v(G_\delta)} G_\delta(x,y) \phi(x) \phi(y).
\end{equation}
In order to show that $h_\delta(\phi)$ converges to $(h, \phi)$ it suffices to check that 
$$\sigma_\delta^2 \to \sigma^2 = \iint_{D^2} G_0^D(x,y) \phi(x) \phi(y) \dd x \dd y$$ as $\delta\to 0$.
Our goal is therefore to show that the Green function for the random walk (in continuous time) on $G_\delta$ converges to the continuous Green function, in the integrated sense above. Under the sole assumption that random walk converges to Brownian motion, showing pointwise convergence of the Green functions is not completely straightforward; in fact if one wants any kind of uniformity in the arguments $x$ and $y$ this will typically be false close to the diagonal. Showing the integrated convergence of the Green function, which is what we require here, is fortunately much simpler. 

Indeed, fix $x \in v(G_\delta)$. Observe that by definition of the Green function as an occupation measure, 
$$
\delta^2 \sum_{y \in v(G_\delta)} G_\delta(x, y) \phi(y) = \delta^2 \E^\delta_x ( \int_0^{\tau_\delta} \phi(X_s) \dd s) 
$$
where $X$ is the continuous time random walk associated to $G_\delta$, as explained before the statement of the theorem.
We first change variables $s = u \delta^{-2}$ to get 
$$
\delta^2 \sum_{y \in v(G_\delta)} G_\delta(x, y) \phi(y) = \int_0^\infty \E^\delta_x ( \phi(X_{u \delta^{-2}}) 1_{\{ \tau_\delta > u \delta^{-2} \}} ) \dd u. 
$$
Fix $\eps>0$. Choose $K>0$ sufficiently large that $\int_K^\infty \P^\delta_x ( \tau_\delta > u \delta^{-2} ) \dd u \le \eps$ for all $x \in D$, which is possible since we assumed in \eqref{eq:uniformtightness} that $\delta^2\tau_\delta$ is uniformly integrable (uniformly in space). We may also assume without loss of generality that $\int_K^\infty \P( \tau> u ) \dd u \le \eps$ (where $\tau$ is the first hitting time of $\partial D$ by $B$), because $\sup_{x \in D}\E_x( \tau)< \infty$ as $D$ is bounded. We then observe, letting $M = \|\phi\|_\infty$, that
\begin{align*}
& \sup_{x\in v(G_\delta)} \left| \int_0^\infty \E^\delta_x ( \phi(X_{u \delta^{-2}}) 1_{\{ \tau_\delta > u \delta^{-2} \}} ) \dd u - \int_0^\infty \E_x( \phi(B_u) 1_{\{\tau > u\} } ) \dd u \right|\\
\le &  \int_0^K \sup_{x\in v(G_\delta)}\left | \E^\delta_x ( \phi(X_{u \delta^{-2}}) 1_{\{ \tau_\delta > u \delta^{-2} \}} ) - \E_x ( \phi(B_u) 1_{\{ \tau> u \}} )\right| \dd u + 2 M \eps 
\end{align*}
Since the position of the random walk on $G_\delta$ killed at $\partial_\delta$ converges uniformly on compact time intervals and uniformly in space to Brownian motion killed when leaving $D$ by \eqref{e:CVassumption}, we deduce that
$$
\limsup_{\delta \to 0} \sup_{x\in v(G_\delta)} \left| \int_0^\infty \E^\delta_x ( \phi(X_{u \delta^{-2}}) 1_{\{ \tau_\delta > u \delta^{-2} \}} ) \dd u - \int_0^\infty \E_x( \phi(B_u) 1_{\{\tau > u\} } ) \dd u \right| \le 2M \eps.
$$
Since $\eps>0$ is arbitrary, we deduce  (using \eqref{eq:Gformula}),
\begin{equation}
 \int_0^\infty \E^\delta_x ( \phi(X_{u \delta^{-2}}) 1_{\{ \tau_\delta > u \delta^{-2} \}} ) \dd u \to  \int_0^\infty \E_x( \phi(B_u) 1_{\{\tau > u\} })\dd u  = \int_D G_0^D (x,y) \phi(y) \dd y, 
 \label{eq:uniformconvergenceGreen}
 \end{equation}
as $\delta \to 0$, uniformly in $x \in v(G_\delta)$. 

It remains to sum over $x\in v(G_\delta)$. Since the above convergence is uniform, and the right hand side of \eqref{eq:uniformconvergenceGreen} is a continuous function of $x\in \bar D$, we deduce, using \eqref{eq:densityone}, that
$$
\delta^2 \sum_{x\in v(G_\delta)} \phi(x)  \int_0^\infty \E^\delta_x ( \phi(X_{u \delta^{-2}}) 1_{\{ \tau_\delta > u \delta^{-2} \}} ) \dd u \to  \iint_D G_0^D (x,y) \phi(x) \phi(y) \dd y \dd x, 
$$
as desired in \eqref{eq:convergencegoal}. This completes the proof of the theorem. 
\end{proof}

\begin{rmk}Let us conclude with some remarks on this theorem. 

\begin{enumerate}

\item When the area of each face $f$ surrounding a given $x \in v(G_\delta)$ is constant as a function of $x$ (and of order $\delta^2$) then the quantity $h_\delta (\phi) = \delta^2\sum_{x\in v(G_\delta)} h(x) \phi(x)$ may be viewed up to a multiple factor (coming from the area of each cell) as the integral of $h_\delta$ against the test function $\phi$, provided that we extend $h_\delta$ to all of $\R^2$ by setting it equal to $h(x)$ in the face $f$. In that case Theorem \ref{T:GFFscalinglimit} says that $h_\delta$, thus extended and viewed as a stochastic process indexed by $\cD_0(D)$, converges in the sense of finite dimensional distributions to a (multiple of) the continuum Gaussian free field. This applies in particular to any periodic lattice such as the square, triangular or hexagonal lattices.

\item In situations where a stronger convergence is desired (such as convergence in the Sobolev space $H^s_0(D)$ for some given $s<0$, as in Proposition \ref{P:GFFsquare}), the Rellich--Kondrachov embedding theorem is a useful criterion which can be used to establish relative compactness (and hence tightness) in such a Sobolev space. In particular, assuming that the boundary of $D$ is at least $C^1$, if $h_\delta$ is a family of random variables in $H^s_0(D)$ such that $\E( \|h_\delta\|_{H^{s'}_0}^2 ) \le C$ for some $s'>s$ and $C< \infty$ independent of $\delta$, then $(h_\delta)_{\delta>0}$ is tight in $H^s_0(D)$. 

This criterion is particularly simple to use in combination with Lemma \ref{L:M0Sobnorm} in order to show convergence in $H^{-1-\eps}$ for any $\eps>0$. Indeed, once we extend $h_\delta$ from the vertices $v(G_\delta)$ to a function defined on all $D$ (for instance we extend $h_\delta$ to be constant on each face, and suppose as above that each face has equal area) then by Lemma \ref{L:M0Sobnorm}, 
$$\|h_\delta\|_{H^{-1}}^2 = \iint_D G_D(x,y) h_\delta(x) h_\delta(y) \dd x \dd y .$$
Taking the expectation, and using similar estimates on the discrete Green function as the ones obtained in the proof of Theorem \ref{T:GFFscalinglimit}, it is not hard to see that 
$$
\E( \|h_\delta\|_{H^{-1}}^2 ) \to \text{const.} \iint G_D(x,y)^2 \dd x \dd y < \infty,
$$
with the above constant related to the area per face. 
Thus by the Rellich--Kondrachov criterion, $(h_\delta)_{\delta>0}$ is tight in the Sobolev space $H^{-1-\eps}$, for any $\eps>0$. By Theorem \ref{T:GFFscalinglimit}, the unique limit point is an appropriate multiple of the Gaussian free field. Hence convergence takes place in distribution in the space $H^{-1- \eps}$, for any $\eps>0$. 

\end{enumerate}

\end{rmk}
 
%We give below an example of such a result. Clearly the result is valid more generally (see the discussion after the proof). 

%The Green function associated to $-Q_N$ is simply

\subsection{Exercises}\label{sec:ex_dir_GFF}

\paragraph{Discrete GFF}
\nopagebreak

\begin{enumerate}[label=\thesection.\arabic*]

%\item Prove that the discrete Green function is symmetric, positive definite, and that $\Delta G(x, \cdot) = - \delta_{x}(\cdot)$.

\item Describe the GFF on a binary tree of depth $n$, where $\partial$ is the root of the tree.

\item Using an orthonormal basis of eigenfunctions for $- \hat Q$, show that the partition function $Z$ in \cref{P:gfflaw} is given by
$$
Z = \det(-\hat Q)^{-1/2}
$$

%\item Show that the law of a discrete GFF $(h(x))_{x \in V}$ is absolutely continuous with respect to $ d\mathbf{x} = \prod_{1\le i \le n} dx_i$, and the joint pdf is proportional to
%$$
%\P( h(\mathbf{x}) ) \propto \exp\left( -\frac14 \sum_{x \sim y }(h(x) - h(y))^2 \right) d\mathbf{x}
%$$
%\index{Dirichlet energy}

\item In this exercise we will show that the minimiser of the discrete Dirichlet energy is discrete harmonic. Fix $U \subset V$ and fix a function $g: V \setminus U \to \R$. Consider 
$$
\inf \{ \cE(f,f), \text{ over } f: V \to \R; f|_{V\setminus U} = g\},
$$
where $\cE(f,f)$ is defined in \eqref{E:Dirichletdiscrete}. 

(a) Show that the inf is attained at a function $f_0$.

(b) Show that $f_0$ is harmonic in $U$: that is, $Q f_0(x) = 0$ for all $x \in U$. To see this, it may be helpful to note that, for every function $\ph$ supported in $U$, and for every $\eps>0$, $\cE( f_0 + \eps \ph, f_0 + \eps \ph) \ge \cE(f_0, f_0)$, and to use the following integration by parts formula: if $u, v: V \to \R$ with $v$ supported on $U$, 
$$
\cE( u,v) =  - (Qu,v).
$$ 
where the inner product on the right hand side is defined in \eqref{l2discrete}.
%\item\label{ex:discretefourier} Prove that discrete GFF can be written as a weighted sum of eigenfunctions, as described in Remark \ref{rmk:discretefourier}.

\item Prove the spatial Markov property of the discrete GFF (Theorem \ref{T:markovdiscrete}). One way to do this is to consider the harmonic extension $\ph$ to $U$ of the boundary data (i.e. $h|_{U^c}$) and check that $h-\ph$ and $ \ph$ are jointly Gaussian vectors indexed by $U$, so the desired property follows by computing suitable covariances.

%\item Perhaps an exercise of the kind: conditional on a field of i.i.d. Gaussians deriving from a height function, then that height function is the GFF?

\end{enumerate}

\paragraph{Continuum GFF}

\begin{enumerate}[label=\thesection.\arabic*]

\setcounter{enumi}{4}

\item \label{Ex:green_half} Show that on the upper half plane,
$$
G_0^{\H}(x,y) = \frac1{2\pi}\log \left| \frac{x - \bar y}{x-y}\right|.
$$
Hint: use that $p_t^{\H}(x,y) = p_t(x,y) - p_t(x, \bar y)$ by symmetry,
and use the formula $e^{-a/t} - e^{- b /t} = t^{-1} \int_a^b e^{- x  /t} \dd x $.

Deduce the value of $G^\D_0(0, \cdot)$ on the unit disc.

\item Let $p_t(x,y)$ be the transition function of Brownian motion on the whole plane (with diffusivity 2). Show that $\int_0^1 p_t(x,y) \dd t =- (2\pi)^{-1} \log |x-y| + O(1)$ as $x \to y$. Then use this to argue that if $D$ is connected and bounded (for simplicity), then $G_0^D(x,y) = - (2\pi)^{-1}\log |x-y| + O(1)$ as $ x \to y$, recovering the third property of \cref{P:basicGFFc}.

\item \label{Green_charact} Let $D$ be a bounded domain and $z_0 \in D$. Suppose that $\phi(z)$ is harmonic in $D\setminus \{ z_0\}$ and that $$\phi(z) = - (2\pi)^{-1}(1+ o(1))\log | z-z_0| \text{ as }z \to z_0 \quad ; \quad  \phi(z)\to 0 \text{ as } z\to w\in \partial D.$$
     Show that $\phi(z) = G_0^D(z_0, z)$ for all $z \in D \setminus\{z_0\}$. (\emph{Hint: use the optional stopping theorem}.)

%\item Show that eigenfunctions of $- \Delta$ corresponding to different eigenvalues are orthogonal on $H_0^1(D)$.

%\item Compute the Green function in one dimension on $D = (0, \infty)$ and $D= (0,1)$ respectively. What are the corresponding Gaussian free fields? What is the spatial Markov property in this context?

\item Let $\mathbf{h}$ be a GFF in a domain $D$. Consider $\tilde{\mathbf{h}}_\eps(z)$, the average value of $\mathbf{h}$ on a square of side length $\eps$ centered at $z$. Let $\tilde h_\eps(z) = \sqrt{2\pi} \tilde{\mathbf{h}}_\eps(z)$. Is this a Brownian motion as a function of $t = \log 1/\eps$? If not, how can you modify it so that it becomes a Brownian motion? More generally, what about the average of the field on a scaled contour $\eps \lambda$, where $\lambda$ is a piecewise smooth loop (the so-called potato average...)?

\item \textbf{Radial decomposition.} Suppose $D = \D$ is the unit disc and $\mathbf{h}$ is a GFF in $D$. Show that $\mathbf{h}$ can be written as the sum
$$
\mathbf{h} = \mathbf{h}_{\text{rad}} + \mathbf{h}_{\text{circ}}
$$
where $\mathbf{h}_{\text{rad}}$ is a radially symmetric function, $\mathbf{h}_{\text{circ}}$ is a distribution with zero average on each disc, and the two parts are independent. Specify  the law of each of these two parts.
\ind{Radial decomposition}

\item Let $D$ be a proper simply connected domain and let $z \in D$.

    (a) Show that $$\log R(z;D) = \E_z(\log |B_T - z |)$$ where $T = \inf\{t>0: B_t \notin D\}$.
    (\emph{Hint: let $g$ be a map sending $D$ to $\D$ and $z_0$ to $0$. Let $\phi(z)=\frac{g(z)}{z-z_0}$ for $z \neq z_0$ and $\phi(z_0) = g'(z_0)$; and consider $\log | \phi|$.})

    (b) Deduce the following useful formula: let $D \subset \C$ be as above, let $U\subset D$ be a subdomain and for $z \in U$ let $\rho_z$ be the harmonic measure on $\partial U$ as seen from $z$. Then show that $\rho \in \mathfrak{M}_0$ and that
    $$
    \var (\mathbf{h}, \rho) = \frac1{2\pi}\log \frac{R(z;D)}{R(z,U)}.
    $$

   \item {Show that the constraints in \cref{GF:generald} uniquely identify $G^D$ when $d\ge 3$. For $x\in D$, defining $H_x(y)=(1/A_d)|x-y|^{2-d}$, let $h_x$ be the unique harmonic extension of $H_x|_{\partial D}$ into $D$. Show that the function $H(x,y)=H_x(y)-h_x(y)$, defined for $x\ne y \in D$, satisfies the constraints of \cref{GF:generald}. Deduce that $G^D=H$. Show this directly by proving that the transition probability $p^D_t(x,y)$ solves the heat equation in $D$.}
\end{enumerate}

\newpage

\section{Liouville measure}\label{S:Liouvillemeasure}

% !TEX root = master.tex

%\section{Liouville measure}\label{S:Liouvillemeasure}

In this chapter we fix $\gamma>0$ (the \textbf{coupling constant}) and introduce the Liouville measure. Informally speaking, this measure $\cM $ (depending on $\gamma$) takes the form
\ind{Coupling constant ($=\gamma$)}
\indN{{\bf Parameters}! $\gamma$; Coupling constant (GMC)}
\begin{equation}\label{Liouville_inf}
\cM(\dd z) = e^{\gamma h(z)} \dd z,
\end{equation} 
where $h = \sqrt{2\pi} \mathbf{h}$ is a GFF in two dimensions (normalised according to \eqref{Dnormalisation}). The scaling factor $\sqrt{2\pi}$ is introduced so that (formally) $\E[ h(x) h(y) ] = - \log |x-y| + O(1)$, that is, $h$ is logarithmically correlated. The construction will be generalised in  \cref{S:GMC} which is devoted to \textbf{Gaussian multiplicative chaos}, which are measures of the form \eqref{Liouville_inf} but for generic log-correlated fields Gaussian fields $h$. While the Gaussian free field in two dimensions is of course an example of such a field, so that Liouville measure really is just a particular case of the theory of Gaussian multiplicative chaos, some arguments specific to the GFF can be used to simplify the presentation and introduce relevant ideas in a clean way, without the need to introduce too much machinery. This is the reason why we have chosen to do the construction of Liouville measure (that is, in the case of the GFF) in this separate chapter.

\paragraph{Heuristics.} The informal definition \eqref{Liouville_inf} should be interpreted as follows. Some abstract Riemann surface has been parametrised, after Riemann uniformisation, by a domain of our choice -- perhaps the disc, assuming that it has a boundary, or perhaps the unit sphere in three dimensions if it doesn't. In this parametrisation, the conformal structure is preserved: that is, curves crossing at an angle $\theta$ at some point in the domain would also correspond to curves crossing at an angle $\theta$ in the original surface.
\ind{Riemann uniformisation}
However, in this parametrisation, the metric and the volume do not correspond to the ambient volume and metric of Euclidean space.
 Namely, a small element of volume $\dd z$ in the domain really corresponds to a small element of volume $e^{\gamma h(z)} \dd z$ in the original surface. Hence points where $h$ is \emph{very big} (for example, thick points) correspond in reality to relatively big portions of the surface; while points where $h$ is very low are points which correspond to small portions of the surface. The first points will tend to be typical from the point of view of sampling from the volume measure, while the second points will be where geodesics tend to travel.

\paragraph{Rigorous approach.} Let $D \subset \R^2$ be an open set and let $h$ be a Dirichlet (or zero boundary) GFF on $D$. When we try to give a precise meaning  \eqref{Liouville_inf}, we immediately run into a serious problem: the exponential of a distribution (such as $h$) is not \emph{a priori} defined. This corresponds to the fact that while $h$ is regular enough to be a \emph{distribution}, so small rough oscillations cancel each other when we average $h$ over macroscopic regions of space, these oscillations become highly magnified when we take the exponential and they can no longer cancel each other out. In fact, giving a meaning to \eqref{Liouville_inf} will require non-trivial work, and will be done via an approximation procedure, using
\begin{equation}
\label{mu_eps}
\cM_\eps(\dd z)  := e^{\gamma h_\eps(z)} \eps^{\gamma^2/2} \dd z,
\end{equation} 
for $\eps>0$, where $h_\eps(z)$ is a jointly continuous version of the circle average. (More general regularisations will be considered in \cref{S:GMC}). It is straightforward to see that $\cM_\eps$ is a (random) Radon measure on $D$ for every $\eps$. Our goal will be to prove the following theorem.

\begin{theorem}
\label{T:liouville}
Suppose $0\le\gamma<2$. If $D$ is bounded, then the random measure $\cM_\eps$ converges weakly almost surely to a
random measure $\cM$, the \emph{(bulk) Liouville measure}, along the subsequence $\eps = 2^{-k}$. $\cM$ almost surely has no atoms, and for  any $A \subset D$ open, we have $\cM(A)>0$ almost surely. In fact, $\E (\cM (A)) = \int_A R(z, D)^{\gamma^2/2} \dd z\in (0,\infty)$.
\end{theorem}
\ind{Liouville measure!Bulk}

We remind the reader that the notation $R(z, D)$ above stands for the conformal radius of $D$ seen from $z$.  That is, $R(z,D) = |f'(0)|$ where $f$ is (any) conformal isomorphism taking $\D$ to $D$ and 0 to $z$. \ind{Conformal radius} If $D$ is unbounded then weak convergence can be replaced by vague convergence with exactly the same proof.

In this form, the result is due  to Duplantier and Sheffield \cite{DuplantierSheffield}. It could also have been deduced from earlier work of Kahane \cite{Kahane85} who used a different approximation procedure, together with results of Robert and Vargas \cite{RobertVargas} showing universality of the limit with respect to the approximating procedure. (In fact, these two results would have given convergence in distribution of $\cM_\eps$ rather than in probability; and hence would not show that the limiting measure $\cM$ depends solely on the free field $h$. However, a strengthening of the arguments of Robert and Vargas due to Shamov \cite{Shamov} has recently yielded convergence in probability.) Earlier, H{\o}egh--Krohn \cite{hoeghkrohn} had introduced a similar model in the context of quantum field theory, and analysed it in the relatively easy $L^2$ phase when $0\le\gamma< \sqrt{2}$.  Here we will follow the elementary approach developed in \cite{BerestyckiGMC}, which works in the more general context of Gaussian multiplicative chaos (see \cref{S:GMC}), but with the simplifications that are allowed by taking the underlying field to be the GFF.

\subsection{Preliminaries}

Before we start the proof of \cref{T:liouville} we first observe that this is the right normalisation.

\begin{lemma}\label{L:var}
We have that $\var (h_\eps(x))  = \log (1/\eps) + \log R(x,D).$ As a consequence,
$$\E (\cM_\eps (A)) = \int_A R(z, D)^{\gamma^2/2} \dd z\in (0,\infty).$$
\end{lemma}

\ind{Conformal radius}
\ind{Circle average}

\begin{proof}
The proof is very similar to the argument in  \cref{T:ca} and is a good exercise. Fix $x \in D$.
%Consider $G_\eps(\cdot)  =  \Delta^{-1} (-2\pi \rho_{x, \eps}) (\cdot)=  2\pi \int_y G(\cdot, y) \rho_{x, \eps}(dy)$. Then by definition, $\Delta G_\eps(\cdot) = -2\pi \rho_{x, \eps}$ as a distribution. In particular, $G_\eps(\cdot)$ is harmonic in $B(z, \eps)$. Furthermore, note that, using integration by parts,
By definition,
\begin{align*}
\var (h_\eps(x)) &=  2\pi \Gamma( \rho_{x, \eps}) = 2\pi \int \rho_{x, \eps}(\dd z) \rho_{x, \eps}(\dd w) G_0^D(z,w).
\end{align*}
For a fixed $z$, $G_0^D(z, \cdot)$ is harmonic on $D \setminus\{z\}$ and so $\int \rho_{x, \eps}(\dd w) G_0^D(w,z) = G_0^D(x,z)$ by the mean value property and an approximation argument similar to \eqref{varcomp2}.
Therefore,
\begin{align*}
\var (h_\eps(x)) &=  2\pi \Gamma( \rho_{x, \eps}) = 2\pi \int \rho_{x, \eps}(\dd z) G_0^D(z,x).
\end{align*}
Now, observe that $2\pi G_0^D(x,\cdot) = - \log | x- \cdot | + \xi(\cdot)$ where $\xi(\cdot)$ is harmonic and $\xi(x) = \log R(x; D)$. Indeed let $\xi(\cdot)$ be the harmonic extension of $- \log |x - \cdot|$ on $\partial D$. Then $2\pi G_0^D(x, \cdot ) + \log| x - \cdot| - \xi (\cdot)$
has zero boundary values on $\partial D$, and is bounded and harmonic in $D\setminus\{x\}$. Hence it must be zero in all of $D$ by uniqueness of solutions to the Dirichlet problem among bounded functions (for example, by the optional stopping theorem). Note that $\xi(x) = \log R(x; D)$ by \eqref{E:cr}. Therefore, by harmonicity of $\xi$ and the mean value property,
$$
\var( h_\eps(x)) =   2\pi \int G_0^D(x, z) \rho_{x, \eps}(\dd z) = \log (1/\eps) + \xi(x)
$$
as desired.
\end{proof}

We now make a couple of remarks:

\begin{enumerate}

\item Not only is the expectation constant, but we have that for each fixed $z$, $ e^{\gamma h_\eps(z)} \eps^{\gamma^2/2}$ forms a martingale as a function of $\eps$. This is nothing but the exponential martingale of a Brownian motion.

\item However, the integral $\cM_\eps(A)$ is \textbf{not} a martingale. This is because the underlying filtration in which $ e^{\gamma h_\eps(z)} \eps^{\gamma^2/2}$ is a martingale depends on $z$. If we try to condition on $(h_\eps(z), z \in D)$, then this gives too much information, and we lose the martingale property.

\end{enumerate}

\subsection{Convergence and uniform integrability in the \texorpdfstring{$L^2$}{TEXT} phase}

The bulk of the proof consists in showing that for any fixed bounded Borel subset $S$ of $D$ (including possibly $D$ itself), we have that $\cM_\eps(S)$ converges almost surely along the subsequence $\eps = 2^{-k}$ to a non-degenerate limit. We will then explain in  \cref{S:weakconvergence}, using fairly general arguments, why this implies the almost sure weak convergence of the sequence of measures $\cM_\eps$ along the same subsequence.

 Let us now fix $S$ and set $I_\eps = \cM_\eps(S)$.  We first suppose that $\gamma\in [0,\sqrt{2})$. In this case, the so called $L^2$ phase, it is relatively easy to check the convergence (which actually holds in $L^2$), %check that $I_\eps$ is integrable. This is easy to check when $\gamma < \sqrt{2}$,
 but difficulties arise when $\gamma \in [\sqrt{2}, 2)$. (As luck would have it this coincides precisely with the phase which is interesting from the point of view of statistical physics on random planar maps).
%We start with the easy case when $\gamma < \sqrt{2}$; this is the so-called $L^2$ phase.

%\begin{lemma}
%\label{L2} $$
%\E(I_\eps^2) \le C \int_{S^2}  \frac{dxdy}{|x-y|^{\gamma^2}}
%$$
%The latter integral is finite if $\gamma^2 < 2$, so $I_\eps$ is bounded in $L^2$ and hence uniformly integrable if $\gamma < \sqrt{2}$.
%\end{lemma}
%
%\begin{comment}
%\begin{proof}
%By Fubini's theorem,
%\begin{align*}
%\E(I_\eps^2) & = \int_{S^2}  \eps^{\gamma^2} \E( e^{\gamma h_\eps(x) + \gamma h_\eps(y)}) dx dy \\
%& = \int_{S^2} \eps^{\gamma^2} \exp( \frac{\gamma^2}2 \var (h_\eps(x) + h_\eps(y)) ) dx dy\\
%& = \int_{S^2}  \exp( {\gamma^2} \cov (h_\eps(x) , h_\eps(y)) ) dx dy\\
%& \le O(1)  \int_{S^2}  \frac{dxdy}{|x-y|^{\gamma^2}},
%\end{align*}
%as desired.
%\end{proof}
%
%Let us now prove convergence, still in the $L^2$ phase where $\gamma  < \sqrt{2}$.
%\end{comment}

%Ideas from the previous part can then easily adapted to cover the case $\gamma \in [\sqrt{2}, 2)$ as well.

\begin{prop}\label{P:conv} If $\gamma\in [0,\sqrt{2})$ and $\eps>0$, $\delta = \eps/2$, then we have the estimate
$\E( (I_\eps - I_\delta)^2) \le C \eps^{2 - \gamma^2}$. In particular, $I_\eps$ is a Cauchy sequence in $L^2 (\P)$ and so converges to a limit in probability as $\eps\to 0$. Along the sequence $\eps = 2^{-k}$, this convergence occurs almost surely, and the limit is almost surely strictly positive.
\end{prop}
\begin{proof} For ease of notations, let $\bar h_\eps (z) = \gamma h_\eps(z) - (\gamma^2/2) \var ( h_\eps(z)) $, and let
$$\sigma(\dd z) = R(z, D)^{\gamma^2/2}{\dd z}.$$

The idea is to say that if we consider the Brownian motions $h_\eps(x)$ and $h_\eps(y)$ (viewed as functions of $\eps = e^{-t}$), then they are (approximately) identical until $\eps \le  |x-y|$, after which time they evolve (exactly) independently.

Observe that by Fubini's theorem,
\begin{align*}
\E( (I_\eps - I_\delta)^2) &  = \int_{S^2} \E\left(( e^{\bar h_\eps(x)} - e^{\bar h_\delta(x)}) ( e^{\bar h_\eps(y)} - e^{\bar h_\delta(y)})\right) \sigma(\dd x) \sigma(\dd y) \\
& = \int_{S^2} \E \left( e^{\bar h_\eps(x) + \bar h_\eps(y)} (1 - e^{\bar h_\delta(x) - \bar h_\eps(x)})(1 - e^{\bar h_\delta(y) - \bar h_\eps(y)})\right) \sigma(\dd x) \sigma(\dd y).
\end{align*}
By the Markov property, $\bar h_\eps(x)+\bar h_\eps(y)$, $\bar h_\eps(x) - \bar h_\delta(x)$ and $h_\eps(y) - h_\delta(y)$ are independent as soon as $| x - y | \ge 2\eps$.
Indeed, we can apply the Markov property in $U = B(x, \eps)$, which allows us to write $h = \tilde h + \varphi$ where $\varphi$ is harmonic in $U$ and $\tilde h$ is an independent GFF in $U$. Since $\tilde h$ is zero outside of $U$, the increment $\bar h _\delta (y) - \bar h_\eps(y)$ and the term $\bar h_{\eps}(y)+\bar h_\eps(x)$ depend only on $\varphi$, and are therefore independent of the increment $\bar h_\delta(x) - \bar h_\eps(x) $ (which, as noted in Theorem \ref{T:ca}, depends only on $\tilde h$.) Applying the same argument with $U=B(y,\eps)$ gives that  $\bar h _\delta (y) - \bar h_\eps(y)$ is independent of the pair $\{\bar h_\delta(x) - \bar h_\eps(x), \bar h_{\eps}(y)+\bar h_\eps(x)\}$. 

%Then the first of the three terms is measurable with respect to $h$ outside of $B(x,\eps)\cup B(y,\eps)$ (that is, depends only on $\varphi$), the second term depends only on $\tilde h$ in the ball $B(x, \eps)$, and the third term depends only on $\tilde h$ in $B(y, \eps)$. 

Hence if $|x - y | \ge 2 \eps$, we can factorise the expectation in the above integral as
\begin{align*}
& = \E( e^{\bar h_\eps(x) + \bar h_\eps(y)}) \E (1 - e^{\bar h_\delta(x) - \bar h_\eps(x)})\E (1 - e^{\bar h_\delta(y) - \bar h_\eps(y)})
%
%\{ ( 1 - 2^{- \gamma^2/2} e^{ - \gamma ( h_{\eps/2 } (x) - h _\eps(x))} ) ( 1 - 2^{- \gamma^2/2} e^{ - \gamma ( h_{\eps/2  }(y) - h _\eps(y))} ) | \cF \right\} \right]
\end{align*}
where both second and third terms are equal to zero, because of the pointwise martingale property.
%But the conditional expectation is in fact zero! Indeed,
%\begin{align*}
% \E\left\{ ( 1 - 2^{- \gamma^2/2} e^{ - \gamma ( h_{\eps/2 } (x) - h _\eps(x))} ) | \cF \right\} &= 1 - 2^{ -\gamma^2/2} \E(  e^{ - \gamma ( h_{\eps/2 } (x) - h _\eps(x))})  \\
% & = 1 - 2^{- \gamma^2/2} e^{\gamma^2/2 \log(2)} = 0
%\end{align*}
Therefore the expectation is just 0 as soon as $|x-y| > 2\eps$.

Also note that by the martingale property for a fixed $x$,
\begin{align*}
\E( (e^{\bar h_\eps(x)} - e^{\bar h_\delta(x)})^2 ) &= \E(  e^{2\bar h_\delta(x)} - e^{2\bar h_\eps(x)}  ) \\
& \le \E(  e^{2\bar h_\delta(x)} )  \le C \E(  e^{2\bar h_\eps(x)} )
\end{align*}
for some $C>0$.
Hence using Cauchy--Schwarz in the case where $|x- y | \le 2 \eps$,
\begin{align}
\E( (I_\eps - I_\delta)^2) & \le \int_{|x- y | \le 2 \eps}\sqrt{ \E( (e^{\bar h_\eps(x)} - e^{\bar h_\delta(x)})^2 ) \E( (e^{\bar h_\eps(y)} - e^{\bar h_\delta(y)})^2 ) } \sigma(\dd x) \sigma(\dd y) \nonumber \\
& \le C\int_{|x- y |\le 2 \eps} \sqrt{\E( e^{2 \bar h_\eps(x)}) \E( e^{2 \bar h_\eps(y)}) } \sigma(\dd x)\sigma( \dd y) \label{CauchySchwarz}\\
& \le  C \int_{|x- y |\le 2 \eps}  \eps^{\gamma^2} e^{\frac12 (2\gamma)^2 \log (1/\eps)}  \sigma(\dd x)\sigma( \dd y) \nonumber\\
& \le C \eps^{2 + \gamma^2  - 2 \gamma^2} = C\eps^{ 2 - \gamma^2}.\nonumber
\end{align}
Thus $I_\eps$ is a Cauchy sequence in $L^2(\P)$. To check almost sure convergence along the subsequence $\eps = 2^{-k}$, we just note that since $\gamma$ is assumed to be smaller than $\sqrt{2}$, the exponent $2-\gamma^2 $ is positive, and hence the sum $\sum_{k\ge 1} 2^{-k (2- \gamma^2)} < \infty$. The almost sure convergence thus follows from the Borel--Cantelli lemma. 

It remains to check that $\P(\lim_{\eps\to 0} I_\eps >0)=1$. We will appeal to Kolmogorov's $0-1$ law. We already know that $\P(\lim_{\eps\to 0} I_\eps>0)>0$, since $\E(\lim_{\eps\to 0} I_\eps)=\lim_{\eps\to 0} \E(I_\eps)>0$. Moreover, notice that if $(f_i)_{i\ge 1}$ is an orthonormal basis of $H^{1}_0(D)$, then  $\{h_\eps(x):x\in D, \eps>0\}$ and therefore $\lim_{\eps\to 0} I_\eps$, is a function of the sequence of coefficients $(h,f_i)_\nabla$. Now, we have seen that these coefficients are independent standard Gaussians, and it is clear that the event $\{\lim_{\eps\to 0}I_\eps>0\}$ is in the tail $\sigma$-algebra generated by the sequence (since this event is invariant under resampling any finite number of terms). Thus it has probability zero or one, and since the probability is positive, it must be one. This concludes the proof of the proposition.
\end{proof}

The moral of this proof is the following: while $I_\eps$ is not a martingale in $\eps$ (because there is no filtration common to all points $x$ such that $e^{\bar h_\eps(x)}$ forms a martingale), %even though it is a martingale for each fixed $x$),
we can use the pointwise martingales to estimate the second moment of the increment $I_\eps - I_\delta$. Only for points $x,y$ which are very close (of order $\eps$) do we get a non-trivial contribution.

We defer the proof of the general case $\gamma \in [0,2)$ until slightly later (see \cref{SS:general}); and for now show how convergence of masses $\cM_\eps(S)$ towards some nonnegative limit implies the almost sure weak convergence of the sequence of measures $\cM_\eps$.

\subsection{Weak convergence to Liouville measure}
\label{S:weakconvergence}

We now finish the proof of {the weak convergence in} \cref{T:liouville} (assuming convergence of masses of fixed bounded Borel subsets $S\subseteq D$ toward some nonnegative limit with probability one) by showing that the sequence of measures $\cM_\eps$ converges in probability for the weak topology towards a measure $\cM$. This measure will be defined by the limits of quantities of the form $\cM_\eps(S), $ where $S$ is a cube such that $\bar S \subset D$. These arguments are borrowed from \cite{BerestyckiGMC}.

%We start by proving that the total mass of the measures $\cM_\eps$ converges almost surely (along $\eps = 2^{-k}$, which we will not repeat). Let $D_n$ be an increasing sequence of domains such that $\cup_n D_n = D$. Let $\ell = \sup_n \cM(D_n)$. Let us show that $\cM_\eps(D) \to \ell$ almost surely as $\eps \to 0$. Note first that $\ell < \infty$ a.s., since by monotone convergence $\E(\ell) = \sup_n \E(\cM(D_n)) = \sigma (D) < \infty$ since $D$ is bounded.

%Let $\delta >0$. Then we can write $\cM_\eps(D) = \cM_\eps(D_n) + \cM_\eps(D \setminus D_n) $ for any $n\ge 1$. The first term converges almost surely as $\eps \to 0$ to $\cM(D_n)$ by the part of the theorem already proved, and for the second term, $ \E (\cM_\eps( D\setminus D_n)) = \sigma(D\setminus D_n) \le \delta^2$ for some $n$ sufficiently large depending only on $\delta$. Fixing that value of $n$, by Markov's inequality, $\P( \cM_\eps(D \setminus D_n) > \delta) \le \delta$, and for the same reason, $\P( |\ell - \cM(D_n) | > \delta) \le \delta$ as well. Then for $\eps$ small enough (depending only on on $n$ and $\delta$, and thus only on $\delta$) we see that $|\cM_\eps(D_n) - \cM(D_n) |\le \delta$ with probability at least $1- \delta$. Altogether, $ \P(| \cM_\eps(D) - \ell |  \ge 2 \delta ) \le 3\delta$. So $ \cM_\eps(D) \to \ell$ almost surely as $\eps \to 0$.

Note that since $\cM_\eps(D) $ converges almost surely,
we have that the measures $\cM_\eps$ are almost surely tight in the space of Borel measures on $\bar D$ with the topology of weak convergence (along the subsequence $\eps = 2^{-k}$, which we will not repeat). Let $\tilde \cM$ be any weak limit.

Let $\cA$ denote the $\pi$-system of subsets of $\R^2$ of
the form $ A = [x_1, y_1) \times [x_2,y_2)$ where $x_i, y_i \in \Q, i=1,2$ and such that $\bar A \subset D$, and note that the $\sigma$-algebra generated by $\cA$ is the Borel $\sigma$-field on $D$.
Observe that $\cM_\eps(A)$ converges almost surely to a limit (which we call $\cM (A)$) for any $A \in \cA$, by the part of the theorem which is already proved (or assumed in the case $\gamma \ge \sqrt{2}$). Observe that this convergence holds almost surely simultaneously for all $A \in \cA$, since $\cA$ is countable.

Let $A = [x_1, y_1) \times [x_2, y_2)\in \cA$. We first claim that
\begin{equation}\label{pm1}
\cM(A) = \sup_{x'_i,y'_i} \{\cM ( [x'_1, y'_1] \times [x'_2, y'_2])   \}
\end{equation}
where the sup is over all $x'_i, y'_i \in \Q$ with $ x'_i > x_i$ and $y'_i < y_i, 1 \le i \le 2$. Clearly the left hand side is almost surely greater or equal to the right hand side, but both sides have the same expectation by monotone convergence (for $\E$). Likewise, it is easy to check that
\begin{equation}\label{pm2}
\cM(A) = \inf_{x'_i, y'_i} \{\cM (  (x'_1, y'_1) \times(x'_2, y'_2))   \}
\end{equation}
where now the infimum is over all $x'_i, y'_i\in \Q$ with $x'_i < x_i$ and $y'_i > y_i, 1 \le i \le 2$.

We aim to check that $\tilde \cM(A) = \cM(A)$, which uniquely identifies  the weak limit $\tilde \cM$ and hence proves the desired weak convergence.

 %Then the normalised measures $\cM^\sharp_\eps : = \cM_\eps / \cM_\eps(D)$ are such that $\cM^\sharp_{\eps'_k} $ converges weakly to $\tilde \cM^\sharp := \tilde \cM / \tilde \cM (D) = \tilde \cM / \ell$ since weak convergence implies convergence of the total mass.

Note that by the portmanteau lemma, for any $A = [x_1, y_1) \times [x_2, y_2)$, and for any $x'_i, y'_i\in \Q$ with $x'_i < x_i$ and $y'_i > y_i, 1 \le i \le 2$, we have:
\begin{align*}
\tilde \cM(A) &\le \tilde \cM ((x'_1, y'_1) \times (x'_2, y'_2))\\
&  \le \liminf_{\eps \to 0} \cM_\eps ((x'_1, y'_1) \times (x'_2, y'_2))\\
& = \cM ((x'_1, y'_1) \times (x'_2, y'_2)).
\end{align*}
 (The portmanteau lemma is classically stated for probability measures, but there is no problem in using it here since we already know convergence of the total mass, and thus can equivalently work with the normalised measures $\cM_\eps / \cM_\eps(D)$).

Since the $x'_i, y'_i$ are arbitrary, taking the infimum over the admissible values and using \eqref{pm2} we get
$$
\tilde \cM(A) \le \cM(A).
$$
The converse inequality follows in the same manner, using \eqref{pm1}. We deduce that $\tilde \cM(A) = \cM(A) $, almost surely, as desired.
As already explained, this uniquely identifies the limit $\tilde \cM$. Hence $\cM_{\eps}$ converges almost surely, weakly, to $\cM$ on $D$. \qed

\subsection{The GFF viewed from a Liouville typical point}\label{sec:rooted_meas}

Let $h$ be a Gaussian free field on a domain $D$, with associated Liouville measure $\cM$ for some  $\gamma <2$. %Let $\cM$ be the associated Liouville measure.
An interesting question is the following: if $z$ is a random point sampled according to the Liouville measure, normalised to be a probability distribution (this is possible when $D$ is bounded), then what does $h$ look like near the point $z$? This gives rise to the concept of \emph{rooted measure} in the terminology of \cite{DuplantierSheffield} or to the Peyri\`ere measure in the terminology of Gaussian multiplicative chaos.
\ind{Rooted measure}
\ind{Peyri\`ere measure}
\ind{Gaussian multiplicative chaos}

We expect some atypical behaviour: after all, for any given fixed $z \in D$, $e^{\gamma h_\eps(z)} \eps^{\gamma^2/ 2}$ converges almost surely to 0, so the only reason $\cM$ could be non-trivial is if there are enough points on which $h$ is atypically big. Of course this leads us to suspect that $\cM$ is in some sense carried by certain thick points of the GFF. It remains to identify the level of thickness. As mentioned before, a simple back of the envelope calculation (made slightly more rigorous in the next result) suggests that these points should be $\gamma$-thick. As we will see, this is in fact a simple consequence of Girsanov's lemma: essentially, when we bias $h$ by $e^{\gamma h(z)}$, we shift the mean value of the field by $\gamma G_0^D(\cdot, z) = \gamma \log 1/|\cdot - z| + O(1)$, thereby resulting in a $\gamma$-thick point.

\begin{theorem}\label{T:rooted}
Suppose $D$ is bounded. Let $z$ be a point sampled according to the Liouville measure $\cM$, normalised to be a probability measure. Then, almost surely,
$$
\lim_{\eps \to 0} \frac{h_\eps(z)}{\log (1/\eps)} = \gamma.
$$
In other words, $z$ is almost surely a $\gamma$-thick point ($z\in \cT_\gamma$).
\end{theorem} \ind{Thick points}
When $D$ is not bounded we can simply say that $\cM ( \cT_\gamma^c) = 0,$ almost surely. In particular, $\cM$ is singular with respect to Lebesgue measure, almost surely.

\begin{proof}
The proof is elegant and at its core relatively simple, but requires a conceptual shift. This comes the following important but elementary lemma, which can be seen as a (completely elementary) version of Girsanov's theorem. Because it is Gaussian in nature rather relying on stochastic calculus, it is appropriate to also credit Cameron and Martin. \ind{Girsanov lemma}

\begin{lemma}[Tilting lemma / Girsanov / Cameron--Martin]\label{L:Girsanov} Let $X = (X_1, \ldots, X_n)$ be a Gaussian vector under the law $\P$, with mean $\mu$ and covariance matrix $V$. Let $\alpha \in \R^n$ and define a new probability measure by
$$
\frac{\dd \Q}{\dd \P} = \frac{e^{\langle \alpha, X\rangle}}{Z} ,
$$
where $Z = \E( e^{\langle \alpha, X \rangle} ) $ is a normalising constant. Then under $\Q$, $X$ is still a Gaussian vector, with covariance matrix $V$ and mean $\mu + V \alpha$.
\end{lemma}

It is worth rephrasing this lemma in plain words. Suppose we weigh the law of a Gaussian vector by some linear functional. Then the process remains Gaussian, with unchanged covariances, however the mean is shifted, and the new mean of the variable $X_i$ say, is
$$\mu'_i = \mu_i + \cov( X_i, \langle \alpha, X\rangle).$$
In other words, the mean is shifted by an amount which is simply the covariance of the quantity we are considering and what we are weighting by.

\begin{proof}
Assume for simplicity (and in fact without loss of generality) that $\mu = 0$. It is simple to check it with Laplace transforms: indeed if $\lambda \in \R^n$, then
\begin{align*}
\Q( e^{{\langle \lambda, X\rangle} }) & = \frac1{Z} \E( e^{ {\langle \lambda + \alpha, X\rangle} }) \\
& = \frac1{e^{\frac12 \langle \alpha, V \alpha \rangle}} e^{\frac12 \langle \alpha + \lambda , V( \alpha + \lambda )\rangle}\\
& = e^{\frac12 \langle \lambda, V \lambda   \rangle  +  \langle \lambda, V \alpha \rangle }
\end{align*}
The first term in the exponent $ \langle \lambda, V \lambda   \rangle  $ is the Gaussian term with variance $V$, while the second term $\langle \lambda, V \alpha \rangle$ shows that the mean is now $V\alpha$, as desired.
\end{proof}

Let $\P = \P(\dd h)$ be the law of the GFF, and let $Q_\eps$ denote the joint law on $(z,h)$ {defined by: %where $h$ has the law $\P$ and given $h$, $z$ is sampled proportionally to $\cM_\eps$. \ellen{(Is this really the right description? The marginal law of $h$ is not $\P$?)} That is,
\begin{equation}\label{eq:rootedChap2}
Q_\eps(\dd z, \dd h) = \frac1Z e^{\gamma h_\eps(z)} \eps^{\gamma^2/2} \dd z \P(\dd h).
\end{equation}
Here $Z$ is a normalising (non-random) constant depending solely on $\eps$. Note that the marginal law of $h$ is weighted by $\cM_\eps(D)$ under $Q_\eps$, and given $h$, the point $z$ is sampled proportionally to $\cM_\eps$.}

Also define $Q(\dd z, \dd h) = (\mathbb{E}(\cM_h(D))^{-1} \cM_h(\dd z) \P(\dd h)$ where by $\cM_h$ we mean the Liouville measure which is almost surely defined by $h$. Note that $Q_\eps$ converges to $Q$ weakly with respect to the product topology induced by the Euclidean metric for $z$ and the Sobolev $H^{-1}$ norm for $h$, say, or, if we prefer the point of view that $h$ is a stochastic process indexed by $\mathfrak{M}_0$, then the meaning of this convergence is with respect to the infinite product $D \times \R^\mathfrak{M}_0$: that is, for any fixed $m\ge 1$ and $\rho_1, \ldots, \rho_m \in \mathfrak{M}_0$, and any continuous bounded function $f$ on $D$,
$$
\E\left( (h, \rho_1) \ldots (h, \rho_m) \int f(z) \eps^{\gamma^2/2} e^{\gamma h_\eps(z)} \dd z\right) \to
\E\left( (h, \rho_1) \ldots (h, \rho_m) \int f(z) \cM_h(\dd z)\right).
$$
This can be verified exactly with the same argument which shows the weak convergence of the approximate Liouville measures. For simplicity we will keep the point of view of a stochastic process for the rest of the proof.

Recall that under the law $Q_\eps$, the marginal law of $h$ is simply
%\begin{equation}\label{E:biasing}
%Q_\eps(dh) = \frac1{Z} \cM_\eps(D) \P(dh)
%\end{equation}
that of a GFF biased by its total mass, so that in particular, $Z = \E( \cM_\eps(D))$ is (up to some small effects from the boundary, which we freely ignore from now on) equal to  $\int_D R(z, D)^{\gamma^2/2} \dd z$, and does not depend on $\eps$.
Furthermore, the marginal law of $z$ is
$$
Q_\eps(\dd z) = \frac1Z \dd z \E(  e^{\gamma h_\eps(z)} \eps^{\gamma^2/2} ) = \frac{\dd z}{Z} R(z, D)^{\gamma^2/2}.
$$
Here again, the law does not depend on $\eps$ and is nice, that is, absolutely continuous, with respect to Lebesgue measure. Finally, it is clear that under $Q_\eps$, given $h$, the conditional law of $z$ is just given by a sample from $\cM_\eps$.

We will simply reverse the procedure, and focus instead on the \emph{conditional distribution of $h$} given $z$. We start by explaining the argument without worrying about its formal justification, and add the justifications where needed afterwards.

Note that by definition,
$$
Q_\eps( \dd h |z) = \frac1{Z(z) } e^{\gamma h_\eps(z)} \eps^{\gamma^2/2} \P(\dd h),
$$
where $Z(z):=R(z,D)^{\gamma^2/2}$.
In other words, the law of the Gaussian field $h$ has been reweighted by an exponential linear functional. By Girsanov's lemma, we deduce that under $Q_\eps(\dd h |z)$, $h$ is a field with the same covariances as under $\P$, and \emph{non-zero} mean at point $w$ given by
$$\cov ( h(w), \gamma h_\eps(z) ) = \gamma \log (1/ | w- z|) + O(1).$$
(More rigorously, we apply Girsanov's lemma to the Gaussian stochastic process $(h, \rho)_{\rho \in \mathfrak{M}_0}$ and find that under $Q_\eps$, its covariance structure remains unchanged, while its mean has been shifted by $\cov ((h, \rho); \gamma h_\eps(z))$.)

 In the limit as $\eps \to 0$, this amounts to adding the function $\gamma G_0^D(\cdot, z)$ to the field $h(\cdot)$. We now argue that this must coincide with the law of $Q(\dd h|z)$. To see this, we use the previous paragraph to write for any $\eps>0$, and for any $m \ge 1$, $\rho_1, \ldots, \rho_m \in \mathfrak{M}_0$,  $\psi\in C_b(D)$:
\begin{multline*}
 \E_{Q^\eps}( (h, \rho_1) \ldots (h, \rho_m) \psi(z) )\! 
  = \! \int_D \! \dd z \,  \psi(z) R(z, D)^{\frac{\gamma^2}{2}}  \E_{h} ( \prod_{i=1}^m
 ((h,\rho_i)+\cov((h,\rho_i),\gamma h_\eps(z))) 
 %\ldots ((h,\rho_m)+\cov((h,\rho_m),\gamma h_\eps(z)))
 )\!.
\end{multline*}
Invoking the weak convergence of $Q_\eps$ to $Q$, we see that the left hand side of the above equality converges to $\E_{Q}( (h, \rho_1) \ldots (h, \rho_m) \psi(z) )$ as $\eps\to 0$. At the same time, an application of the dominated convergence theorem shows that the right hand side converges as $\eps\to 0$ to
\begin{equation}\label{condrootG}
 \int_D \dd z \psi(z) R(z, D)^{\frac{\gamma^2}{2}}  \E_{h} ( ( h + \gamma G_0^D(z, \cdot), \rho_1) \ldots (h+ \gamma G_0^D(z,\cdot), \rho_m)).
\end{equation}
Hence the law of $Q(\dd h|z)$ is as claimed.

To summarise, under $Q$ and given $z$, a logarithmic singularity of strength $\gamma$ has been introduced at the point $z$.
Hence we find that under $Q(\dd h|z)$, almost surely,
$$
\lim_{\delta \to 0} \frac{h_\delta(z)}{\log (1/\delta)} = \gamma,
$$
so $z \in \cT_\gamma$, almost surely as desired. In other words, $Q( \cM_h(\cT_\gamma^c) = 0) = 1. $

We conclude the proof of the theorem by observing that the marginal laws $Q(\dd h)$ and $\P(\dd h)$ are mutually absolutely continuous with respect to one another, so any property which holds almost surely under $Q$ holds also almost surely under $\P$. (This absolute continuity follows simply from the fact that $ \cM(S) \in (0, \infty), \P-$almost surely)
\end{proof}

\subsection{The full \texorpdfstring{$L^1$}{TEXT} phase}
\label{SS:general}

To address the difficulties that arise when $\gamma \ge \sqrt{2}$, we proceed as follows. Roughly, we claim that the second moment of $I_\eps$ blows up because of rare points which are \emph{too thick} and which do not contribute to the integral in an almost sure sense, but nevertheless inflate the value of the second moment. So we will remove these points by hand. To see which points to remove, we appeal the considerations of the previous section: this suggests that we should be safe to get rid of points that are strictly more than $\gamma$-thick. \ind{Thick points}

Let $\alpha > 0$ be fixed (it will be chosen $> \gamma$ and very close to $\gamma$ soon). We define a good event
$
G^\alpha_\eps(x) = \{ h_{\eps}(x) \le \alpha \log (1/\eps)  \},
$
for which the point $x$ is not too thick at scale $\eps$.
% Further let $\bar h_\eps(x) = \gamma h_\eps(x) - (\gamma^2/2 )  \log (1/\eps)$ to ease notations.

%\begin{lemma}[Ordinary points are not thick]
%\label{L:typicalthick} For any $\alpha >0$, we have that uniformly over $x\in S$, $\P(G^\alpha_\eps(x)) \ge 1- p(\eps_0)$ where the function $p$ may depend on $\alpha$ and for a fixed $\alpha > \gamma$, $p(\eps_0) \to 0$ as $\eps_0 \to 0$. \end{lemma}

%\begin{proof}
%This follows immediately from the fact that for a standard Brownian motion $B_t$, $B_t / t \to 0$ as $t \to \infty$, almost surely
%\end{proof}

\begin{lemma}[Liouville points are no more than $\gamma$-thick]  \label{L:Liouvillethick_GFF} For $\alpha > \gamma$ we have
$$
\E(e^{\bar h_\eps(x) } 1_{G^\alpha_\eps(x)} ) \ge 1- p(\eps)
$$
where the function $p$ may depend on $\alpha$ and for a fixed $\alpha > \gamma$, $p(\eps) \to 0$ as $\eps \to 0$, polynomially fast. The same estimate holds if $\bar h_\eps(x)$ is replaced with $\bar h_{\eps/2}(x)$.
\end{lemma}

\begin{proof}
Note that
$$
\E( e^{ \bar h_\eps(x)}  \indic{G^\alpha_\eps(x) }) = \tilde \P (  G^\alpha_\eps(x)), \text{ where } \frac{\dd \tilde \P}{\dd \P} =  e^{\bar h_\eps(x)}.
$$
By Girsanov's lemma, under $\tilde \P$, the process $X_s = h_{e^{-s}}(x)$ has the same covariance structure as under $\P$ and its mean is now $\gamma \cov(X_s, X_t) = \gamma s + O(1) $ for $ s\le t$. Hence it is a Brownian motion with drift $\gamma$, and the lemma follows from the fact that such a process does not exceed $\alpha t$ at time $t$ with high probability when $t$ is large (and the error probability is exponential in $t$, or polynomial in $\eps$, as desired).

Changing $\eps$ into $\eps /2$ means that the drift of $X_s$ is $\gamma s  + O(1)$ over a slightly larger interval of time, namely until time  $t + \log 2$. In particular the same argument as above shows that the same estimate holds for $\bar h_{\eps/2}(x)$ as well.
\end{proof}

We therefore see that points which are more than $\gamma$-thick do not contribute significantly to $I_\eps$ in expectation and can be safely removed. To this end, we fix $\alpha >\gamma$ and introduce:
\begin{equation}
J_\eps = \int_S  e^{ \bar h_\eps(x)} \mathbf{1}_{G_\eps(x)}\sigma(\dd x) ; \quad \hat{J}_{\eps/2}(x) = \int_S  e^{ \bar h_{\eps/2}(x)} \mathbf{1}_{G_\eps(x)}\sigma(\dd x)
\end{equation}
with $G_\eps(x) = G^\alpha_\eps(x)$, {and where we recall that $\sigma(\dd x) = R(x, D)^{\gamma^2/2} \dd x$}. Note that a consequence of \cref{L:Liouvillethick_GFF} is that
\begin{equation}
\label{I-J}
\E(|I_\eps - J_\eps|) \le p(\eps)|\sigma(S)|\to 0 \text{ and } \E ( |I_{\eps/2} - \hat J_{\eps/2}|) \le p(\eps) |\sigma(S)| \to 0
\end{equation}
as $\eps\to 0$.

\begin{lemma} \label{L:UI}
We have the estimate
$\E( (J_\eps - \hat{J}_{\eps/2})^2) \le \eps^{r}$ for some $r>0$.
In particular, $J_\eps$ is a Cauchy sequence in $L^2$. Along $\eps = 2^{-k}$, this convergence occurs almost surely.
\end{lemma}

\begin{proof} The proof of this lemma is virtually identical to that in the $L^2$ phase (see \cref{P:conv}). The key observation there was that if $|x- y | \ge 2\eps$, then the increments $h_\eps(x) - h_{\eps/2}(x)$ and $h_\eps(y) - h_{\eps/2}(y)$ are independent of each other, and in fact also of $\cF$: the $\sigma$-algebra generated by $h$ restricted to the complement of $B(x,\eps)\cup B(y,\eps)$. Since the events $G_\eps(x)$ and $G_\eps(y)$ are both measurable with respect to $\cF$, we may therefore deduce from \emph{that} proof (see \eqref{CauchySchwarz}) that
$$
\E( (J_\eps - \hat{J}_{\eps/2})^2) \le C \int_{|x-y|\le 2\eps}\sqrt{\E( e^{2 \bar h_\eps(x)}  \mathbf{1}_{G_\eps(x)} )  \E( e^{2 \bar h_\eps(y)}  \mathbf{1}_{G_\eps(y)} )} \sigma(\dd x)\sigma( \dd y).
$$
Now,
\begin{align*}
\E( e^{2 \bar h_\eps(x)}  \mathbf{1}_{G_\eps(x)} ) & \le \E( e^{2 \bar h_\eps(x)}  \mathbf{1}_{  \{ h_\eps(x) \le \alpha \log (1/\eps)  } )\\
& \le O(1) \eps^{-\gamma^2} \Q( h_\eps(x) \le \alpha \log 1/\eps)
\end{align*}
where by Girsanov's lemma, under the law $\Q$, $h_\eps(x)$ is a normal random variable with mean $2 \gamma \log (1/\eps) + O(1)$ and variance $\log 1/\eps + O(1)$. This means that
$$
 \Q( h_\eps(x) \le \alpha \log 1/\eps) \le O(1) \exp ( - \frac12 (2\gamma - \alpha)^2 \log 1/\eps)
$$
and hence
$$
\E( (J_\eps - \hat{J}_{\eps/2})^2)  \le O(1) \eps^{2 - \gamma^2} \eps^{ \frac12 (2\gamma - \alpha)^2}.
$$
Again, choosing $\alpha > \gamma$ sufficiently close to $\gamma$ ensures that the bound on the right hand side is at most $O(1) \eps^r$ for some $r>0$, as desired. This finishes the proof of the lemma. {It  also concludes the proof of  \cref{T:liouville} in the general case $\gamma<2$, by \eqref{I-J}, and recalling that $p(\eps)$ decays polynomially in $\eps$ for fixed $\alpha$, so we can apply Borel--Cantelli to get almost sure convergence along the sequence $\eps=2^{-k}$.}
\end{proof}

{\begin{proof}[Proof of Theorem \ref{T:liouville}] As a consequence of Lemma \ref{L:UI} and \eqref{I-J}, $I_\eps$ is a Cauchy sequence in $L^1 $ and so converges to a limit in probability. The almost sure convergence along the sequence $\eps =2^{-k}$ follows from the fact that $p(\eps)$ converges to zero polynomially fast by \cref{L:Liouvillethick_GFF} and the Borel--Cantelli lemma.  The limit of $I_\eps$ is almost surely strictly positive by the same $0-1$ argument as in the case $\gamma\le \sqrt{2}$, and the weak convergence in the sense of measures then follows by the argument detailed in Section \ref{S:weakconvergence}. The formula for the expectation of the limit in Theorem \ref{T:liouville} is a consequence of Lemma \ref{L:var}. Finally, the fact that $\cM$ almost surely has no atoms follows from Exercise \ref{Ex:thick}.
\end{proof}}

We note that the almost sure convergence over the entire range of $\eps$ (not just the dyadic values $\eps = 2^{-k}$)
was proved by Sheffield and Wang \cite{SheffieldWang}.

\paragraph{Understanding the phase transition for the Liouville measure.} %in Gaussian multiplicative chaos}

\ind{Gaussian multiplicative chaos}
The description of the Liouville measure viewed from a typical point, explained in Theorem \ref{T:rooted}, was established in the entire $L^1$ regime and so for all $0 < \gamma <2$. The fact that the Liouville measure $\cM = \cM_\gamma$ is supported on the $\gamma$-thick points, $\cT_\gamma$, is very helpful to get a clearer picture what changes when $\gamma=2$. %of Gaussian multiplicative chaos (Kahane's general theory of measures of the form $e^{\gamma X(z)}dz$ where $X$ is a log-correlated Gaussian field).
Indeed recall from Theorem \ref{T:thick} that $\dim( \cT_\gamma) = (2- \gamma^2/2)_+$, and $\cT_\gamma$ is empty if
$\gamma>2$. The point is that $\cM = \cM_\gamma$ does not degenerate when $\gamma<2$ \emph{because} there are thick points to support it. Once $\gamma >2$ there are no longer any thick points, and this makes it in some sense ``clear'' that any approximations to $\cM_\gamma$ must degenerate to the zero measure. 

When $\gamma = 2$ however, $\cT_\gamma$ is not empty, and there is therefore a hope to  construct a meaningful critical Liouville measure $\cM_2$. %corresponding to the \emph{critical Gaussian multiplicative chaos}.
Such a construction has indeed been carried out in two separate papers by Duplantier, Rhodes, Sheffield, and Vargas \cite{DRSV1, DRSV2}. However the normalisation must be done more carefully -- see these two papers for details, as well as the more recent works \cite{JunnilaSaksman, HRVdisk, Pow18chaos}.\ind{Gaussian multiplicative chaos!critical}

\subsection{Change of coordinate formula and conformal covariance}

\ind{Conformal covariance of Liouville measure}
Of course, it is natural to wonder in what way the conformal invariance of the GFF manifests itself at the level of the Liouville measure. As it turns out these measures are not simply conformally invariant. This is easy to believe intuitively, since the total mass of the Liouville measure has to do with total surface area (measured in quantum terms) enclosed in a domain, %(for example via circle packing), 
and so this must grow as the domain grows.

However, the measures are \textbf{conformally covariant}: that is, to relate their laws under conformal mappings one must include a correction term accounting for the inflation of the domain under the conformal map. This term is naturally proportional to the derivative of the conformal map.

To formulate the result, it is convenient to use the following notation. Suppose that $h$ is a given distribution -- perhaps a realisation of a GFF, but also perhaps one of its close relatives (for example, the GFF plus some smooth deterministic function) -- and suppose that its circle average process is well defined. Then we define $\cM_h$ to be the measure, if it exists, given by $\cM_h(\dd z) = \lim_{\eps \to 0} e^{\gamma h_\eps(z) } \eps^{\gamma^2/2}dz.$ Of course, if $h$ is just a GFF, then $\cM_h$ is nothing else but the measure we have constructed in the previous part. If $h$ can be written as $h = h_0 + \varphi$ where $\varphi$ is deterministic, $h_0$ is a GFF and $\e^{\gamma\varphi}\in L^1(\cM_{h_0})$, then $\cM_h (dz) = e^{\gamma \varphi(z)} \cdot \cM_{h_0}(dz)$ is absolutely continuous with respect to the Liouville measure $\cM_{h_0}$.

\begin{theorem}[Conformal covariance of Liouville measure] \label{T:ci}
Let $f:D \to D'$ be a conformal isomorphism, and let $h$ be a GFF in $D$. Then $h' = h \circ f^{-1}$ (where we define this image in the sense of distributions) is a GFF in $D'$, and
\begin{align*}
 \cM_h \circ f^{-1} & = \cM_{h \circ f^{-1} + Q \log |(f^{-1})'|}\\
 & = e^{\gamma Q  \log |(f^{-1})' |} \cM_{h'},
\end{align*}
where
$$
Q = \frac{\gamma}2 + \frac{2}{\gamma}.
$$\indN{{\bf Parameters}! $Q$; parameter in change of coordinates formula for LQG}
\end{theorem}

In other words, pushing forward the Liouville measure $\cM_h$ by the map $f$, we get a measure which is absolutely continuous (with density $|(f^{-1})'(z)|^{\gamma Q}$ at $z\in D'$) with respect to the Liouville measure on $D'$. The quantity $Q$ plays a very important role in the theory developed in the subsequent chapters. In physics it is known under the name of \textbf{background charge}. \ind{Background charge}

\mn \textbf{Informal proof.}
The idea behind this formula may be understood quite easily. Indeed, note that $\gamma Q = \gamma^2/ 2 + 2$. When we use the map $f$, a small circle of radius $\eps$ is mapped \emph{approximately} into a small circle of radius $\eps' = |f'(z) | \eps$ around $f(z)$. So $e^{\gamma h_\eps(z)} \eps^{\gamma^2/2} \dd z$ approximately corresponds to
$$
e^{\gamma h'_{|f'(z) | \eps} (z')} \eps^{\gamma^2/2} \frac{\dd z'}{|f'(z)|^2}
$$
by the usual change of variable formula. This can be rewritten as
$$
e^{\gamma h'_{\eps'} (z')} (\eps')^{\gamma^2/2} \frac{\dd z'}{|f'(z)|^{2 + \gamma^2/2}}
$$
Letting $\eps \to 0$ we get, at least heuristically speaking, the desired result.\qed

\begin{proof}[Proof of \cref{T:ci}.] Of course, the above heuristic is far from a proof, and the main reason is that $h_\eps(z)$ is not a very well behaved approximation of $h$ under conformal maps. It is better to instead work with a different approximation of the GFF, using an orthonormal basis of $H_0^1(D)$ as in \cref{ZGFFd}, which has the advantage of being conformally invariant.
	
	In view of this, we make the following definition: suppose $h = \sum_n X_n f_n$, where $X_n$ are i.i.d. standard normal random variables, and $f_n$ is an orthonormal basis of $H_0^1(D)$. Set $h^N(z) = \sum_{i=1}^N X_i f_i$ to be the truncated series, and define
$$
\cM^N(S) = \int_S\exp\left( \gamma h^N(z) - \frac{\gamma^2}{2} \var (h^N(z)) \right) \sigma(\dd z)
$$
where we recall that $\sigma(\dd z) = R(z, D)^{\gamma^2/2}\dd z$.
Note that $\cM^N(S)$ has the same expected value as $\cM(S)$. Furthermore, $\cM^N(S)$ is a non-negative martingale with respect to the filtration $(\cF_N)_{N}$ generated by $(X_N)_N$, so has an almost sure limit which we will call $\cM^*(S)$.
%Using uniform integrability of $\cM_\eps(S)$, and a few tricks such as Fatou's lemma, it is possible to show:
\begin{lemma}
\label{L:univ}
Almost surely, $\cM^*(S) = \cM(S)$.
\end{lemma}
%See exercises, and \cite[Lemma 5.2]{BerestyckiGMC} for details.

\begin{proof}
%Write $h = \sum_{n=1}^\infty X_n f_n$ where $X_n$ are i.i.d. standard Gaussian random variables and $f_n$ are an orthonormal basis of $H_0^1(D)$. Set $h^N = \sum_{n=1}^N X_n f_n$ to be the truncated series.
When we take the circle averages of the series we obtain
$$
h_\eps = h^N_\eps + h'_\eps
$$
where $h'_\eps$ is independent from $h^N$, and $h^N_\eps$ denotes the circle average of the function $h^N$.
Hence
$$
\eps^{\gamma^2/2} e^{\gamma h_\eps(z)} = e^{\gamma h^N_\eps(z)} \eps^{\gamma^2/2} e^{\gamma h'_\eps(z)}.
$$
Consequently, integrating over $S$ and taking the conditional expectation given $\cF_N$, we obtain that
$$
\cM^N_\eps(S): = \E( \cM_\eps(S) | \cF_N) =  \int_S\exp\left( \gamma h_\eps^N(z) - \frac{\gamma^2}{2} \var (h_\eps^N(z)) \right) \sigma(\dd z).
$$
When $\eps\to 0$, the right hand side converges to $\cM^N(S)$, since $h^N$ is a nice smooth function. Consequently,
$$
\cM^N(S) = \lim_{\eps \to 0} \E( \cM_\eps(S)|\cF_N).
$$
{Since $\cM_{\eps}(S)\to \cM(S)$ in $L^1$, we have $\cM^N (S)= \lim_{\eps\to 0} \E(\cM_\eps(S)|\cF_N)=\E(\cM(S)|\cF_N)$ and so by martingale convergence, $\cM^N(S)\to \cM(S)$ as $N\to \infty$. Hence $\cM(S)=\cM^*(S)$, as desired.}
\begin{comment}
By Fatou's lemma, the right hand side is greater or equal to $\E( \cM(S) | \cF_n)$.
Hence
$$
\cM^N(S) \ge \E(\cM(S) | \cF_N).
$$
We can take $N\to \infty$ and deduce from the martingale convergence theorem (and the fact that $\cM(S)$ is measurable with respect to $\cF_\infty$) that
$$
\cM^*(S) \ge \cM(S).
$$
However, from Fatou's lemma again, we know that
$$
\E( \cM^*(S)) \le \lim_{N\to \infty} \E( \cM^N(S)) = \int_D R(z,D)^{\gamma^2/2} dz.
$$
We also know (by uniform integrability of $\cM_\eps(S)$) that
$$
\E( \cM(S)) = \int_D R(z,D)^{\gamma^2/2} dz.
$$
Consequently we see that $\E(\cM^*(S)) \le \E( \cM(S))$. The only possibility is that $\E(\cM(S)) = \E(\cM^*(S)) $ and therefore also $\cM(S) = \cM^*(S)$, as desired.
\end{comment}
\end{proof}

To finish the proof of conformal covariance (\cref{T:ci}) we now simply recall that if $f_n$ is an orthonormal basis of $H_0^1(D)$ then $f_n \circ f^{-1}$ gives an orthonormal basis of $H_0^1(D')$. Hence if $h' = h \circ f^{-1}$, then its truncated series $h'_N$ can also simply be written as $h'_N = h^N \circ f^{-1}$. Thus, consider the measure $\cM^N$ and apply the map $f$. We obtain a measure $\tilde \cM'_N$ in $D'$ such that
\begin{align*}
\tilde \cM'_N(D') &= \int_{D'}   \exp \{ \gamma h^N (f^{-1}(z') )- \frac{\gamma^2}{2} \var (h^N(f^{-1}(z')))  \} R(f^{-1}(z'), D)^{\gamma^2/2} \frac{\dd z'}{|f'(f^{-1}(z'))|^{2}} \\
& = \int_{D'} \dd \cM'_N(z')  e^{( 2 + \gamma^2/2) \log|(f^{-1})'(z')|},
\end{align*}
where $\dd \cM'_N$ is the approximating measure to $\cM_{h'}$ in $D'$. (The second identity is justified by properties of the conformal radius). Letting $N\to \infty$, and recalling that $d\cM'_N$ converges to $d\cM_{h'}$ by the previous lemma, we obtain the desired statement of conformal covariance. This finishes the proof of \cref{T:ci}.
 \end{proof}

\subsection{Random surfaces}

The notion of \textbf{random surface} is a way of identifying Gaussian free field type distributions that give rise to different ``parametrisations'' of the same Liouville measure.
Essentially, we want to consider the surfaces encoded by $\cM_h$ and by $\cM_h \circ f^{-1}$ to be ``the same'' for any given conformal isomorphism $f:D \to D'$. By the conformal covariance formula (\cref{T:ci}) if $h$ is a GFF, we have 
\begin{equation}\label{E:coordinate_change}
\cM_h \circ f^{-1} = \cM_{ h'} \text{ almost surely, where } h'= h \circ f^{-1} + Q \log |(f^{-1})'|.
\end{equation}
 Thus we should think of $h$ and $h'$ as encoding the same quantum surface.

In fact, (when $h$ is a GFF) this equality holds almost surely for \emph{all} $D'$ and \textbf{all} conformal isomorphisms $f: D \to D'$ simultaneously. This result was proved by Sheffield and Wang in \cite{SheffieldWang}.

This motivates
the following definition, due to Duplantier and Sheffield \cite{DuplantierSheffield}. Define an equivalence relation on pairs $(D,h)$, consisting of a simply connected domain $D$ and an element $h$ of $\cD'(D)$, by declaring that \[ D_1 \sim D_2\] if there exists $f: D_1 \to D_2$ a conformal isomorphism such that
$$
h_2 = h_1 \circ f^{-1} + Q \log | ( f^{-1})'|.
$$
It is easy to see that this is an equivalence relation.

\begin{definition}\label{D:surface}
A (random) surface is a pair $(D,h)$ consisting of a domain and a (random) distribution $h \in \cD'(D)$, where the pair is considered modulo the above equivalence relation.
\end{definition}

Observe that this definition of (random) surface depends on the parameter $\gamma\ge 0$ of the Liouville measure (since $Q$ depends on $\gamma$).

\ind{Random surface}
\ind{Change of coordinates formula}

Interesting random surfaces arise, among other things, when we sample a point according to the Liouville measure (either in the bulk, or on the boundary for a free field with a non-trivial boundary behaviour, see later), and we `zoom in' near this point. Roughly speaking, these are the \emph{quantum cones} and \emph{quantum wedges} introduced by Sheffield in \cite{zipper}. A particular kind of wedge will be studied in a fair amount of detail later on in these notes (see \cref{T:wedge}).

\ind{Quantum wedges}
\ind{Quantum cones}

\subsection{Exercises}

\begin{enumerate}[label=\thesection.\arabic*]

\item Explain why \cref{L:UI} and \cref{L:Liouvillethick_GFF} imply uniform integrability of $\cM_\eps(S)$.

%\item \nb{Don't we do this already?} Let $\cM_n$ be a random sequence of Borel measures an open set $D\subset \R^d$. Suppose that for each open $S\subseteq D$, $\cM_n(S) $ converges in probability to a finite limit $\ell(S)$. Explain why $\ell$ defines uniquely a Borel measure $\cM$ and why $\cM_n$ converges in probability weakly to $\cM$.

%\item Prove the convergence of $\cM_\eps(S)$ towards a limit in the general case where $\gamma<2$. (Hint: just use the truncated integral $J_\eps$ introduced in the proof of uniform integrability, and check $L^2$ convergence).

\item Let $\cM $ be the Liouville measure with parameter $0\le \gamma <2$. Use uniform integrability and the Markov property of the GFF to show that $\cM(S)> 0 $ almost surely.

%\item Let $\cT(\alpha)$ be the set of $\alpha$-thick points. Is $\cT(\alpha)$ conformally invariant? (See \cite{HMP}).

\item (a) How would you normalise $e^{\gamma h_\eps(z)}$ if you were just aiming to define the Liouville measure on some {line} segment contained in $D$, or more generally a smooth simple curve in $D$? Show that with this normalisation you get a non-degenerate limit. 

(b) What is the conformal covariance property in this case?

\item \label{Ex:thick} Recall the events $G_\eps(z) = \{ h_\eps(z) \le \alpha \log 1/\eps\}$ % \text{ for all } r \in [\eps_0 , \eps]\}$ 
from the proof of uniform integrability of the Liouville measure in the general case. Show that {for any }
$0 \le d < 2 - \gamma^2/2<2$, 
$$
\E\left( \int_{S^2} \frac1{|x-y|^d} e^{\bar h_\eps(x)}\mathbf{1}_{ G_\eps(x)} \sigma(\dd x)\  e^{\bar h_\eps(y)} \mathbf{1}_{ G_\eps(y)} \sigma(\dd y) \right)\le C< \infty
$$
where $C$ does not depend on $\eps$. Deduce that
$$
\dim ( \cT_\gamma) \ge 2 - \gamma^2/ 2
$$
almost surely {and that $\cM$ almost surely has no atoms}. Conclude with a proof of  \cref{T:thick}.

%\item \label{Ex:rooted} Write carefully the continuity argument of $Q_\eps$ in the proof of \cref{T:rooted}. \nb{Is this still needed?} \ellen{I don't think so.} Under the rooted measure $Q(dh |z)$, what are the fluctuations of $h_\eps(z)$ as $\eps \to 0$? \ellen{Also this?}

%Show that if $h$ is a GFF and $z$ is sampled according to the Liouville measure, then given $z$, $h$ is absolutely continuous with respect to $\tilde h (\cdot)+ \gamma \log (1/|\cdot - z|)$, where $\tilde h$ is a Gaussian free field.

%\item Prove $\cM^*(S) = \cM(S)$.

%\item Conclude the proof of conformal covariance, assuming Lemma \ref{L:univ}.

\end{enumerate}

\newpage

\section{Gaussian multiplicative chaos}
\label{S:GMC}
% !TEX root = master.tex

\subsection{Motivation, background}

In the previous chapter we constructed the Liouville measure, which is an example of Gaussian multiplicative chaos. Gaussian multiplicative chaos is the theory, developed initially by Kahane in the 1980s, whose goal is the definition and study of random measures of the form \indN{{\bf Gaussian multiplicative chaos}! $\cM$; general Gaussian multiplicative chaos measure} 
$$
\cM (\dd z) = \exp ( \gamma h(z) - \frac{\gamma^2}{2} \E (h(z)^2)) \sigma(\dd z),
$$ \ind{Coupling constant ($=\gamma$)}
where $\gamma$ is a parameter (the \textbf{coupling constant})\indN{{\bf Parameters}! $\gamma$; Coupling constant (GMC)}, $h$ is a centred, logarithmically correlated Gaussian field, and $\sigma$ is a reference measure. In this chapter we will give a modern presentation of the general theory. There are two main reasons why we devote an entire chapter to this theory.
The first one, continuing on the theme of previous chapters, is because the tools we will develop in the process are very useful for the study of Liouville measure: for instance, they will enable us to describe the multifractal spectrum of Liouville measure, leading us to the KPZ relation\footnote{here KPZ stands for Khnizhnik--Polyakov--Zamolodchikov, and should not be confused with the Kardar--Parisi--Zhang equation.}, which was one of the initial motivations of the seminal work of Duplantier and Sheffield \cite{DuplantierSheffield}. The second, and possibly more important reason, is that
Gaussian multiplicative chaos has arisen in many natural models coming from  motivations beyond Liouville quantum gravity. For instance, Kahane's original motivation, following the pioneering ideas of Mandelbrot and in particular Kolmogorov, was the description of \textbf{turbulence} and especially the phenomenon of intermittency (see \cite{Frisch} for a classical book on the legacy of Kolmogorov in turbulence, and we refer, for instance, to \cite{ChevillardHDR} and \cite{CGRV} for a recent survey and article outlining the connection to Gaussian multiplicative chaos). In this case it is of course more natural to assume that the field lives in three (rather than two) dimensions. We have however already observed that the Gaussian free field is not logarithmically correlated except in dimension two: indeed the correlations are given by the Green function, which is a multiple of $|y-x|^{2-d}$ in dimensions $d \ge 3$. This is one of the reasons why developing a general theory (going beyond the two dimensional case of the Gaussian free field) is of great interest. Let us mention briefly a few further topics, where a connection to Gaussian multiplicative chaos has been observed.
\begin{itemize}
  \item In \textbf{random matrix theory}, Gaussian multiplicative chaos describes (or is conjectured to describe) the limiting behaviour of (powers of) the characteristic polynomial of a large random matrix drawn from many of the classical random matrix ensembles. See, for example, \cite{WebbCUE} followed by \cite{NSW} for the case of CUE, and \cite{BWW} for a general class of random Hermitian matrices including GUE. Other relevant works include (but are not limited to) \cite{ForkelKeating,Kivimae,LambertOstrovskySimm}.

  \item In number theory, Gaussian multiplicative chaos describes the (real part) of the \textbf{Riemann zeta function} on the critical line, when it is recentred at a random point. See \cite{SaksmanWebb}, and references therein for a long line of works leading to this result.

  \item In \textbf{mathematical finance}, Gaussian multiplicative chaos is used as a model for the stochastic volatility of a financial asset, following some highly influential works of Bacry, Muzy and Delour (\cite{BacryMuzyDelour}), Mandelbrot et al. \cite{Mandelbrot_finance}, and also Cont \cite{Cont}).

  \item In the study of \textbf{planar Brownian motion} a closely related theory has been developed by Jego \cite{jegoBMC}, \cite{jegoRW}, and A\"id\'ekon, Hu and Shi \cite{AHS}; both of these works generalise the older paper of Bass, Burdzy and Koshnevisan \cite{BBK} to the full $L^1$ regime. In the case of the random walk, see \cite{jegoRW}
 as well as \cite{AbeBiskup} (although this requires wiring the boundary); and
      also \cite{ABJL} in the context of the loop soup with explicit connections to Liouville measure.
\end{itemize}

The reader is also invited to consult the survey \cite{RhodesVargas_survey}, which contains many useful discussions, facts and references.

\subsection{Setup for Gaussian multiplicative chaos}\label{setup_gmc}
%{Now using $\mathbf{d}$ for the dimension of $\sigma$ and $d$ for the dimension of the space, to be consistent with the previous sections.}

 \dis{
 The next few sections of this chapter could be skipped by a reader interested only in the GMC measures associated to the Gaussian free field (that is, the Liouville measures). In this case the reader may wish to skip to  \cref{S:multifractal}, although tools such as Kahane's inequality (\cref{T:Kahane}) will be needed.}

%Following the construction we will present a set of useful tools for the study of Gaussian multiplicative chaos measures; in particular, Kahane's powerful \textbf{convexity inequality}. Finally, we will apply these tools to the study of moments of GMC and in particular derive the \textbf{multifractal spectrum of GMC} and discuss applications to the famous \textbf{KPZ formula}.

\label{SS:setup}
We first describe the setup in which we will be working.
We consider a more general setup than before and in particular for the rest of this chapter we do not assume we are working exclusively in two dimensions. Let $D \subset \R^d$ be a bounded domain.
{Consider a kernel $K(x,y)$ of the form
\begin{equation}\label{cov}
 K(x,y) = \log (|x-y|^{-1} ) + g(x,y)
\end{equation}
where $g$ is continuous over $\bar D \times \bar D$. 

{We will now define a Gaussian field whose pointwise covariance function is, in some sense given by \eqref{cov}. We will proceed in a way that is somewhat inspired by what we did in the case of the Gaussian free field (Chapter \ref{S:GFF}), defining first a space of measures $\mathfrak{M}$ against which we will be able to test the field, and viewing the field as a stochastic process indexed by this space of measures. However, we will face an immediate difficulty. Namely, if we define $\mathfrak{M}$ to be the space of bounded signed measures\footnote{A signed measure $\rho = \rho^+ - \rho^-$ is called bounded if $\rho^+$ and $\rho^-$ are finite.} such that $\iint |K(x,y) | \rho| (\dd x) |\rho|( \dd y)< \infty$, it is hard to show directly that this forms a vector space. (In the case of the GFF, this was a consequence of Lemma \ref{L:linearitycov}, which has no obvious analogue here.) In two dimensions, this would follow from the theory developed in that chapter, since integrability against $K$ is essentially equivalent to integrability agains $G_D$. However, we wish to develop the theory in a setting which is more general. Our approach will thus differ a little, and the index set of our stochastic process will first be defined as an abstract Hilbert space induced by $K$. While this works well, the downside is that it is not easy to check if a concrete bounded signed measure $\rho$ lies in this space. We will essentially check it for bounded measures which can not only integrate $K$ (or, equivalently, can integrate the logarithm), but can do so \emph{continuously}, see the energy condition \eqref{eq:assumptionrho} below.}

\medskip Let us now explain how to associate a stochastic process to the covariance function $K$ in \eqref{cov}. For smooth, compactly supported functions $f,g \in \cD(D)$ we write 
\begin{equation}\label{Kfg}
K(f,g) = \int K(x,y) f(x) g(y) \dd x \dd y.
\end{equation}
 Then $K$ is bilinear over the vector space $\cD(D)$. We assume that $K$ is positive definite in the sense that $K(f,f) \ge 0$ for all $f \in \cD(D)$ and $K(f,f) = 0$ implies $f = 0$. Thus $K(f,g)$ defines an inner product over $\cD(D)$ and $f\mapsto \| f\|_K := \sqrt{K(f,f)}$ defines a norm over $\cD(D)$. We let $\mathfrak{M}$ be the Hilbert space completion of $\cD(D)$ with respect to this norm. By definition, if $\rho\in \mathfrak{M}$ then there exists a sequence of smooth, compactly supported functions $f_n\in \cD(D)$ such that $f_n \to \rho$ with respect to the norm induced by $K$, i.e., $\|f_n - \rho, f_n -\rho\|_K^2 \to 0$. Note that if $K$ is the Green function with Dirichlet boundary conditions in a bounded domain $D \subset \R^2$ then $\mathfrak{M}{=H_0^{-1}(D)}$ {with inner product $(\cdot, \cdot)_{-1}$ and its restriction to signed measures} coincides with our earlier notation $\mathfrak{M}_0$.  
}

\indN{{\bf Function spaces}! $\mathfrak{M}$; completion of $\cD(D)$ with respect to the inner product induced by a covariance kernel $K$}

%Set \indN{{\bf Function spaces}! $\mathfrak{M}_+$; non-negative measures with finite energy with respect to a kernel $K$}
%$$
%\mathfrak{M}_+ =\{ \rho \text{ a nonnegative measure in $D$ such that:} \iint |K(x,y)| \rho(\dd x) \rho(\dd y) < \infty\},
%$$
% and set $\mathfrak{M}$\indN{{\bf Function spaces}! $\mathfrak{M}$; difference of two elements of $\mathfrak{M}_+$} to be the set of signed measures of the form  $ \rho = \rho_+ - \rho_-$, where $\rho_\pm \in \mathfrak{M}_+$. Note that $\mathfrak{M}$ contains all smooth compactly supported functions in $D$. 

Let $h$ be the centred Gaussian generalised function with covariance $K$. That is, we view $h$ as a stochastic process indexed by $\mathfrak{M}$, characterised by the two properties that: $(h, \rho)$ is linear in $\rho \in \mathfrak{M}$ in the sense that $(h, \alpha \rho_1 + \beta \rho_2 ) = \alpha (h, \rho_1) + \beta(h, \rho_2)$ almost surely for $\alpha,\beta\in \R$, $\rho_1,\rho_2\in \mathfrak{M}$; and for any $\rho \in \mathfrak{M}$,
$$
(h,\rho) \text{ is a centered Gaussian random variable with variance } \| \rho\|_K^2.
%\iint  K(x,y) \rho(\dd x)\rho(\dd y).
$$
{Equivalently, $((h, \rho))_{\rho \in \mathfrak{M}}$ is a centered Gaussian stochastic process with covariance
 \begin{equation}\label{covfunctionnew}
 (\rho, \rho')\in {\mathfrak{M}}^2 \mapsto \E [ (h,\rho)( h,{\rho'})] =  (\rho, \rho')_K
 \end{equation}
  (where the latter denotes the inner product associated to the Hilbert space $\mathfrak{M}$, associated to the norm $\| \cdot \|_K$). }
We will write $\int h(x) \rho(\dd x) $ or $(h,\rho)$ with an abuse of notation for the random variable $h_\rho$. {Note that this setup covers the case of a Gaussian free field in two dimensions with Dirichlet boundary conditions. }
%(more precisely, given a domain $D$, the restriction of the GFF with Dirichlet boundary conditions in $D$ to any subdomain $D'$ such that $ \bar D' \subset D$ satisfies the requirement \eqref{cov}). 
In fact, it also covers the case of the Gaussian free field with free or Neumann boundary conditions, see \cref{S:SIsurfaces}. %by changing $\gamma$ into $2\gamma$ if necessary. 
We extend the definition of $h$ outside of $D$ by setting $h|_{D^c} = 0$, so if $\rho$ {is, say, a measure} such that $\rho|_D \in \mathfrak{M}$, by definition $(h, \rho) = (h , \rho|_D)$. %{Why do we do this?}

As mentioned above, it is not straightforward to check if a concrete nonnegative measure lies in $\mathfrak{M}$. Let $\rho$ be a fixed non-negative finite Borel measure on $\R^d$ compactly supported in $D$, and suppose that $\rho$ satisfies the \emph{energy} condition
\begin{equation}
\label{eq:assumptionrho}
x\mapsto \cE_\rho(x) :=  \int \log (\frac1{|x-y|}) \rho(\dd y) \text{ is continuous on {$\bar B(0,R)$}},
\end{equation}
{where $R>0$ is large enough that $\text{Supp} (\rho) \subset B(0,R)$.} 
{In particular, $\cE_\rho$ is uniformly bounded on $\bar B(0,R)$.} It is easy to check that the condition \eqref{eq:assumptionrho} is satisfied whenever $\rho$ has an $L^{p}$ Lebesgue density supported in $D)$  for some ${p}>1$, but also in many other cases, for example, when $\rho$ is the uniform distribution on some circle contained in $D$. 

{For $\rho$ satisfying the energy condition \eqref{eq:assumptionrho} it is} not hard to see that $\rho$ can be identified with an element of $\mathfrak{M}$ as follows. Let us fix a smooth, symmetric nonnegative function $\psi$ supported in $B(0,1)$ such that $\int \psi (z) \dd z = 1$, and let $\rho_n (y) = \rho *\psi_n(y) = \int \rho (\dd z) \psi_n( y - z)$, where $\psi_n (z) = n^d \psi ( zn )$.  Then $\rho_n$ is a smooth, compactly supported function in $D$ for $n$ large enough. Furthermore, one can check that $\|\rho_m - \rho_n\|_K\to 0$ as $n, m \to \infty$, i.e., $(\rho_n)_{n\ge 0}$ is a Cauchy sequence with respect to the norm induced by the covariance kernel $K$. Indeed consider, for instance, $K(\rho_n, \rho_m)$. By definition of $\rho_n$ and $\rho_m$, using the symmetry of $\psi_n$ and Fubini's theorem {(in the nonnegative case, as $K$ is bounded below from \eqref{cov} and $\rho_n, \rho_m$ have finite mass)} 
%\ellen{$K$ isn't necessarily nonnegative, but I think Fubini's theorem is still OK by \eqref{eq:assumptionrho}}
\begin{align*}
K( \rho_n, \rho_m) & = \iint \rho_n(x) K(x,y) \rho_m(y) \dd x \dd y\\
%& = \iint \int \rho (\dd z) \psi_n(x-z) K (x,y)\int \rho( \dd w) \psi_m(y - w)  \dd x \dd y\\
& = \iint \rho (\dd z) \rho (\dd w)   \iint \psi_n(z- x)   K(x,y) \psi_m(w- y) \dd x \dd y \\
&= \E[ \iint   K (z + \frac{X}{n}, w + \frac{Y}{m} )\rho (\dd z ) \rho (\dd w) ] 
\end{align*}
where $X,Y$ are independent and distributed according to $\psi(z) \dd z$. Let us check that the above expression converges to $K( \rho, \rho)$. 
{First,} as $\rho$ is finite and the function {$g$} in the definition of $K$ is bounded and continuous, it suffices to replace $K(z + X/n, w+ Y/m)$ by $- \log | z- w - X/n/ -Y/m|$. Let us fix $z$ and consider the integral in $w$. This integral can be rewritten as $$
\int - \log ( |z +\frac{X}{n} - \frac{Y}{m} - w |) \rho ( \dd w) =  \cE_\rho ( z + \frac{X}{n} - \frac{Y}{m})$$ in the notation of \eqref{eq:assumptionrho}. Hence as $n, m\to \infty$ this converges to $\cE_\rho(z)$ by continuity of $\cE_\rho$ in $z$. 
Furthermore as $\cE_\rho$ is bounded and $\rho$ is a finite measure, {we have} in turn $\int\ \cE_\rho( z + X/n - Y/m) \rho (\dd z) \to \int \cE_\rho(z) \rho (\dd z) = K(\rho,\rho)$ by dominated convergence. Using boundedness again, we can take expectations by dominated convergence, {and we see that $K(\rho_n,\rho_m)\to K(\rho,\rho)$ as $n,m\to \infty$, as desired}. For the same reason, all four terms in the expansion of $K( \rho_n - \rho_m, \rho_n - \rho_m)$ converge to the same limit (i.e., $K( \rho, \rho)$) and thus $\rho_n - \rho_m$ is a Cauchy sequence with respect to the norm $\| \cdot \|_K$. {The sequence $(\rho_m)_{m\ge 1}$ therefore has a limit {in $\mathfrak{M}$} and this limit does not depend on the choice of smoothing (i.e., on the choice of $\psi$ subject to the above properties). Let us denote this limit by $\rho_{\mathfrak{M}}$. Note then that for two measures $\rho, \rho'$ satisfying \eqref{eq:assumptionrho} we have the following expression for the covariance function \eqref{covfunctionnew} of the field:
\begin{align}
 (\rho_{\mathfrak{M}}, \rho'_{\mathfrak{M}})_K = \lim_{n \to \infty} K( \rho_n, \rho'_n)  = K( \rho, \rho') := \iint \rho(\dd x ) K(x,y) \rho' (\dd y).\label{eq:extendedcov}
\end{align} 
%\ellen{*I think it's confusing to bring in $h$ above, I would delete the leftmost expression, since this is not what is used to justify that $\rho$ can be identified with $\rho_{\mathfrak{M}}$.} \nb{On the contrary this is what actually matters... I have tried to re-explain this better below.}
{Hence a measure satisfying \eqref{eq:assumptionrho} induces an element $\rho_{\mathfrak{M}}$ of the abstract Hilbert space $\mathfrak{M}$, and the corresponding covariance has the concrete expression \eqref{eq:extendedcov} in this case. As the map $\rho \mapsto \rho_{\mathfrak{M}}$ is linear, we will still write, with a similar abuse of notation as earlier, $(h, \rho)$ for the random variable $h_{\rho_{\mathfrak{M}}} = ( h, \rho_{\mathfrak{M}})$. To summarise, for measures $\rho$ and $\rho'$ satisfying \eqref{eq:assumptionrho}:
$$
\cov ( (h, \rho); (h, \rho')) = \iint \rho(\dd x ) K(x,y) \rho' (\dd y).
$$
This is the expression we will use in practice when constructing Gaussian multiplicative chaos below. }

%\nb{I will delete the two sentences below if you agree; the issue is that we have not proved that $\rho \mapsto \rho_{\mathfrak M}$ is eg injective so I don't want to use or imply things like ``we identify $\rho$ with $\rho_{\mathfrak{M}}$''.} 
%[[For such measures the associated \ellen{inner product, equivalently, covariance when tested against the field $h$} is simply $K(\rho, \rho')$, as in \eqref{Kfg} but with $f(x) \dd x$ and $g(y) \dd y$ replaced by $\rho(\dd x)$ and $\rho'(\dd y)$.
%For this reason we make a small abuse of notation we will still denote by $\rho = \rho_{\mathfrak{M}}$, so that we view $\rho$ directly as an element of $\mathfrak{M}$.]] 
%
%}

Let $\sigma$\indN{{\bf Gaussian multiplicative chaos}! $\sigma$; reference measure} be a Radon measure\footnote{in fact, on $\R^d$ every locally finite Borel measure is Radon, so it would suffice to assume that $\sigma$ is a locally finite Borel measure.} on $\bar  D$ of dimension at least  ${\mathfrak{d}}$\indN{{\bf Gaussian multiplicative chaos}! $\mathfrak{d}$; dimension of the reference measure} (where $0\le \mathfrak{d} \le d$), in the sense that,
\begin{equation}\label{dim}
\iint_{\ \bar D\times \bar D} \frac1{| x-y|^{{\mathfrak{d} - \eps}}}\sigma(\dd x) \sigma(\dd y) < \infty
\end{equation}
for all $\eps>0$ (so for example, if $\sigma$ is Lebesgue measure, then $d={\mathfrak{d}}$).
% In particular $\sigma$ is a finite measure since $\bar D$ is bounded.
 Note that ${\mathfrak{d}} \ge 0$ and may be equal to 0, but the statement of the theorem below will be empty in that case. In particular, we will only care about the case ${\mathfrak{d}}>0$, which prevents $\sigma$ from having any atoms. {Throughout this chapter, when $\mathfrak{d}>0$ we will fix a number $0<\mathbf{d}< \mathfrak{d}$, such that
\begin{equation}\label{dim2}
\iint_{\ \bar D\times \bar D} \frac1{| x-y|^{{\mathbf{d}}}}\sigma(\dd x) \sigma(\dd y) < \infty.
\end{equation}
}
\begin{rmk}
{Many results of this chapter (for example, Theorem \ref{T:conv}, \cref{T:finiteposmoments}) are stated under an assumption of a strict inequality involving $\mathbf{d}$; since $\mathbf{d}$ can be chosen arbitrarily close to $\mathfrak{d}$ the same results could be stated by replacing $\mathbf{d}$ by $\mathfrak{d}$.}
\end{rmk}
%Assume without loss of generality that $D$ contains the ball of radius 10 around the origin and let $S$ be the unit cube $(0,1)^k$.

%\nb{Suppose also that $\theta \in \mathfrak{M}$.} 

Now, fix $\theta$ a nonnegative, Borel measure on $\R^d$ which supported on the closed unit ball $\bar B(0,1)$ and such that $\theta(\R^d) = 1$. Suppose in addition that $x \mapsto \cE_\theta(x)$ is continuous (and hence uniformly bounded) on $\bar B(0,5)$, as in \eqref{eq:assumptionrho}. Thus {we may test} $h$ against $\theta$, as explained above.
%\begin{equation}
%\nb{
%\label{eq:assumptionrho}
%%\label{eq:assumptiontheta}
%x\mapsto \cE(x) :=  \int \log (1/|x-y|) \theta(\dd y) \text{ is continuous (hence uniformly bounded) on $\bar B(0,5)$}. 
%}
%\end{equation}
More generally, for $\eps>0$, set $\theta_\eps(\cdot)$\indN{{\bf Gaussian multiplicative chaos}! $\theta_\eps(\cdot)$; mollifier at scale $\varepsilon$} to be the image of the measure $\theta$ under the mapping $x \mapsto \eps x$, that is $\theta_\eps(A) =  \theta (A/\eps)$ for all Borel sets $A$. We view this as an approximation of the identity based on $\theta$ (and will sometimes write $\theta_\eps(x)\dd x$ for the measure $\theta_\eps(dx)$ with an abuse of notation). We also write $\theta_{x, \eps}(\cdot)$ for the measure $\theta_\eps$ translated by $x$. For $x \in D$, {it is now straightforward to see that $\theta_{x, \eps}$ also satisfies \eqref{eq:assumptionrho} and thus can be viewed as an element of $\mathfrak{M}$.}
%\footnote{It is tempting to only use $\theta \in \mathfrak{M}$ instead of \eqref{eq:assumptiontheta} as an assumption on $\theta$. However this leads to problems in the proofs below, as pointed out by an anonymous referee, even if we assume $\theta$ has a density. Consider the following instructive example in dimension $d=1$. Take $f(x) = (c/|x|) (1+ \log^+(1/|x|))^{-2}$ for $x \in (-1,1)$ and $f(x) = 0$ else. Then if $\theta(dx) = f(x) dx$, we have $\theta \in \mathfrak{M}$. However, it is not the case that \eqref{eq:roughcov} holds. In fact even the basic fact \eqref{covbound} does not hold, as the left hand side is unbounded as $t \to \infty$ (since $\int \log(1/|x|) f(x) dx = \infty$).\label{footnote}}
 %
 %by \eqref{eq:assumptiontheta}, the translated measure $\theta_{x,\eps}\in \mathfrak{M}$. 
 So we can define an $\eps$-regularisation of the field $h$ by setting for $\eps$ small,
\begin{equation}\label{convol}
h_\eps(x) = h \ast \theta_\eps(x) = \int h(y) \theta_\eps( x-y ) \dd y  = \int h(y) \theta_{x, \eps}(\dd y)\; , \; x \in D.
\end{equation}
One can check that $\var (h_\eps(x) - h_\eps(x')) \to 0 $ as $|x-x'| \to0$ for a fixed $\eps$, so there exists a version of the stochastic process $h$ such that $h_\eps(x)$ is almost surely a Borel measurable function of $x \in S$ (see for example Proposition 2.1.12 in \cite{GineNickl}).
%(With a bit more calculations and Kolmogorov's continuity criterion, it is in fact not hard to see that there exists a jointly continuous version of $h_\eps(x)$ as a function of $x \in S$ and $\eps>0$ small enough).
%\begin{remark} In certain situations one may wish to relax the assumption of smoothness on $\theta$. In fact this is not used at any point in the proof, except to guarantee that $X\ast \theta_\eps$ is well defined as a random variable. If this fact is known for other reasons (for example if $h$ is the two-dimensional GFF and $\theta$ is the uniform distribution on the unit circle, in which case $h_\eps$ is simply the usual circle average), then the theorems and arguments below go through without any change.
%\end{remark}
Hence for any Borel set $S \subset D$ and $\gamma\ge 0$ we may define
\indN{{\bf Gaussian multiplicative chaos}! $\cM_\eps$; approximation of $\mathcal{M}$ at spatial scale $\eps$}
\begin{equation}\label{regmeas}
\cM_\eps(S) = I_\eps = \int_S e^{\gamma h_\eps(z) - \frac{\gamma^2}{2}  \E( h_\eps(z)^2 ) }\sigma(\dd z).
\end{equation}
In the previous chapter where $h$ was the 2d Dirichlet Gaussian free field (multiplied by $\sqrt{2\pi}$), our choice for $\sigma$ was $\sigma(\dd z) = R(z, D)^{\gamma^2/2}\dd z$, and our choice for the measure $\theta$ was the uniform distribution on the unit circle (so $h_\eps(z)$ was the usual circle average process of $h$).
%However, the case where $\sigma$ is the occupation measure of an (independent) planar Brownian motion is also of interest as it is used for defining the Liouville Brownian motion \cite{GRV, NB}, the canonical diffusion process in Liouville quantum gravity. In this example, $\sigma$ is singular with respect to Lebesgue measure, yet $\sigma$ is of dimension two in the sense of \eqref{dim}.

\subsection{Construction of Gaussian multiplicative chaos}

With these definitions we can state the result that guarantees the existence of Gaussian multiplicative chaos. For simplicity we assume $D$ bounded, and let $S \subset D$ be a Borel subset (which may be equal to $D$ itself).

\begin{theorem}\label{T:conv}
Let $0 \le \gamma < \sqrt{2\mathbf{d}}$ {(equivalently, $0 \le \gamma < \sqrt{2 \mathfrak{d}}$)}. Then $\cM_\eps(S)$ converges in probability and in $L^1(\P)$ to a limit $\cM(S)$. The random variable $\cM(S)$ does not depend on the choice of the regularising kernel $\theta$ subject to the above assumptions. Furthermore, the collection $(\cM(S))_{S\subset D}$ defines a Borel measure $\cM$ on $D$, and $\cM_\eps$ converges in probability towards $\cM$ for the topology of weak convergence of measures on $D$. {Finally, the measure $\cM$ almost surely has no atoms.}
\end{theorem}

\dis{In later chapters, we will sometimes also use the notation $\cM_h$ or $\cM_h^\gamma$ to indicate the dependence of $\cM$ on the underlying field $h$ or the field $h$ and the parameter $\gamma$.}
\indN{{\bf Gaussian multiplicative chaos}! $\cM_h$; Gaussian multiplicative chaos associated with a field $h$}
\indN{{\bf Gaussian multiplicative chaos}! $\cM_h^\gamma$; Gaussian multiplicative chaos associated with a field $h$ and parameter $\gamma$}

%{Maybe too much work, but this might be a place to clarify that we get almost sure\ convergence along a specific subsequence (this is not written up yet), and actually can deduce convergence almost surely\ as $\eps\to 0$ (as I wrote up in the sketch for our old project with Juhan). }

Let us assume without loss of generality that $\mathbf{d}>0$, so that $\sigma$ has no atoms.

\medskip As before, the main idea will be to pick $\alpha > \gamma$ and
consider the normalised measure $e^{\gamma h_\eps(x)}\dd x$, but \emph{restricted to good points}; that is, points that are not too thick.
We will check that the $L^1$ contribution of bad points is negligible (essentially by the above Cameron--Martin--Girsanov observation), while the remaining part is shown to remain bounded and in fact convergent in $L^2(\P)$. The key will be to take a good and slightly more subtle definition of the notion of \emph{good points}, that makes the relevant $L^2$ computation very simple.

In \cite{BerestyckiGMC}, uniqueness of the limit was obtained by comparing to a different approximation of the field, arising from the Karhuhen--Loeve expansion of $h$. This gives another approximation of the measure which turns out to be a martingale, and hence also has a limit. \cite{BerestyckiGMC} then showed that the two measures must agree, thereby deducing uniqueness. Here we present a slightly simpler argument partly based on a remark made Hubert Lacoin (private communication).

\subsubsection{Uniform integrability}
The goal of this section will be to prove:

\begin{prop}\label{P:UI}
$I_\eps$ is uniformly integrable.
\end{prop}

\begin{proof}
Let $\alpha > 0$ be fixed (it will be chosen $> \gamma$ and very close to $\gamma$ soon). It is helpful in this proof to introduce another regularisation kernel $\theta^*$ satisfying the same assumptions as $\theta$, namely \eqref{eq:assumptionrho}. We then denote by $h^*_r (z) = h\star \theta_{z, r}^*$.  
For $s \in S$, and $n \in \Z$, set 
$$
E_n(x) = \{ h^*_{e^{-n}}(x) \le \alpha n \}. 
$$

%We will use the following notation in the rest of the article: for $r >0$ we define
%\begin{equation}\label{rbar}
%\bar r = e^{  \lceil \log r \rceil} = \inf \{ e^k: k \in \Z, e^k > r \}
%\end{equation}
% to be the closest upper $e$-adic approximation of $r$. 
We define a \textbf{good event}
\begin{equation}\label{good}
G^\alpha_\eps(x) = \bigcap_{n=n_0}^{n(\eps)} E_n(x),
%\{ h_{\bar r}(x) \le \alpha \log (1/\bar r) \text{ for all $r \in [\eps, \eps_0] $} \}
\end{equation}
where $e^{-n_0} = \eps_0  \le 1$ for instance, and $ n (\eps) = \lceil \log (1/\eps)\rceil$. This is the good event that the point $x$ is never too thick (with respect to the regularisation kernel $\theta^*$) up to scale $\eps$. Mostly, we think of the case where $\theta = \theta^*$, but in the proof of uniqueness (independence of the limit with respect to the choice of the regularisation kernel) it will be useful to define the good event in such a way that it does not depend on $\theta$, hence our choice to let it depend on some other fixed arbitrary regularisation kernel $\theta^*$.   
 Further let $\bar h_\eps(x) = \gamma h_\eps(x) - (\gamma^2/2 )  \E( h_\eps(x)^2)$ to ease notations.

\begin{lemma}[Ordinary points are not thick]
\label{L:typicalthick} For any $\alpha >0$, we have that uniformly over $x\in S$, $\P(G^\alpha_\eps(x)) \ge 1- p(\eps_0)$ where the function $p$ may depend on $\alpha$ and for a fixed $\alpha > 0$, $p(\eps_0) \to 0$ as $\eps_0 \to 0$. \end{lemma}

\begin{proof}
Set $X^*_t = h^*_\eps(x)$ for $\eps = e^{- t}$. Then a direct computation {using the covariance formula \eqref{eq:extendedcov}} 
and \eqref{cov} (see below in \cref{L:cov}, and more precisely \eqref{eq:roughcov}), implies that
\begin{equation}\label{covbound}
|\cov(X^*_s, X^*_t) - s\wedge t| \le O(1),
\end{equation}
 where the implicit constant is uniform.
In particular $\var (X^*_t ) = t + O(1)$.

Note that for each $k\ge 1$, $\P( X^*_k \ge \alpha k /2) \le e^{- \alpha^2 k^2/ ( 8\var (X^*_k))} $ which decays exponentially in $k$ by the above, and so is smaller than $Ce^{- \lambda k}$ for some $\lambda>0$ (we may take $\lambda = \alpha^2/16$ for instance).
Hence
$$
\P( \exists k \ge k_0:|X^*_k | \ge \alpha k) \le \sum_{k \ge k_0} Ce^{- \lambda k }
$$
We call $p(\eps_0) $ to be the right hand side of the above for $k_0 = \lceil-\log(\eps_0)\rceil $ which can be made arbitrarily small by picking $\eps_0$ small enough.
This proves the lemma.
\end{proof}

\begin{lemma}[Liouville points are no more than $\gamma$-thick]  \label{L:Liouvillethick} For $\alpha > \gamma$ we have
$$
\E(e^{\bar h_\eps(x) } 1_{G^\alpha_\eps(x)} ) \ge 1- p(\eps_0).
$$
\end{lemma}

\begin{proof}
Note that
$$
\E( e^{ \bar h_\eps(x)}  \indic{G^\alpha_\eps(x) }) = \tilde \P (  G^\alpha_\eps(x)), \text{ where } \frac{\dd \tilde \P}{\dd \P} =  e^{\bar h_\eps(x)}.
$$
Analogously to $X_t^*$, set $X_t := h_\eps(x))$ for $\eps = e^{- t}$. By the Cameron--Martin--Girsanov lemma, under $\tilde \P$, the process $(X^*_s)_{- \log \eps_0 \le s\le t }$ has the same covariance structure as under $\P$ and its mean is now $\gamma \cov(X^*_s, X_t) = \gamma s + O(1) $ for $ s\le t$ (see again below in \cref{L:cov}, and more precisely \eqref{eq:roughcov}). Hence
$$
\tilde \P( G^\alpha_\eps(x))  \ge \P (  G^{\alpha - \gamma }_\eps(x) ) \ge 1-p(\eps_0)
$$
by \cref{L:typicalthick} since $\alpha > \gamma$.
\end{proof}

We thus see that points which are more than $\gamma$-thick do not contribute significantly to $I_\eps$ in expectation and can therefore be safely removed. We therefore fix $\alpha >\gamma$ and introduce:
\begin{equation}
J_\eps = \int_S  e^{ \bar h_\eps(z)} \indic{G_\eps(z)}\sigma(\dd z)
\end{equation}
with $G_\eps(x) = G^\alpha_\eps(x)$. We will show that $J_\eps$ is uniformly integrable from which the result follows.

Before we embark on the main argument of the proof, we record here for ease of reference an elementary estimate on the covariance structure of $h_\eps(x)$. Roughly speaking, the role of the first estimate \eqref{eq:roughcov} is to bound from above (up to an unimportant constant of the form $e^{O(1)}$) the contribution to $\E(J_\eps^2)$ coming from points $x,y$ that are close to each other. That will suffice to prove uniform integrability. The role of the finer estimate \eqref{eq:finecov} is to get a more precise estimate to the contribution to $\E( J_\eps^2)$ coming from points $x,y$ which are macroscopically far away, which we will be able to assume thanks to \eqref{eq:roughcov}. This time the error in the covariance up to an additive term $o(1)$ will translate into an error up to a factor $e^{o(1)} = 1+ o(1)$ in the estimation of this contribution. In turn this will imply convergence.

\begin{lemma}
  \label{L:cov} We have the following estimate:
  \begin{equation}\label{eq:roughcov}
  \cov( h_\eps(x), h^*_r(y)) = \log 1/(|x-y|\vee r \vee \eps) + O(1).
  \end{equation}
  Moreover, if $\eta>0$ and $|x-y| \ge \eta$, then
  \begin{equation}\label{eq:finecov}
  \cov( h_\eps(x), h^*_\delta(y)) = \log (1/|x-y|) + g(x,y) + o(1)
  \end{equation}
  where $o(1)$ tends to 0 as $\delta, \eps \to 0$, uniformly in $|x- y | \ge \eta$.
\end{lemma}

Note that, as $\theta$ and $\theta^*$ are completely arbitrary subject to \eqref{eq:assumptionrho}, the same estimates hold with $h_\eps(x)$ replaced by $h_\eps^*(x)$ and/or $h^*_r(y)$ replaced by $h_r(y)$.

\begin{proof}
  We start with the proof of \eqref{eq:roughcov}. Assume without loss of generality that $\eps \le r$. Note that {by \eqref{eq:extendedcov}},
  \begin{align}
  \cov( h_\eps(x), h^*_r(y)) &= \iint K(z,w) \theta_{x, \eps}(\dd w) \theta^*_{y,r}(\dd z)\nonumber\\
  & = \iint -\log (| w-z |)  \theta_{x,\eps}(\dd w) \theta^*_{y,r}(\dd z) + O(1)\label{roughcov1}
  \end{align}
  We consider the following cases: (a) $ r \le |x-y |/3$, and (b) $r \ge | x-y | /3$.

  In case (a),  $|x-y| \le \eps + |w-z| +r \le 2 r + |w-z| \le (2/3) |x-y| + |w-z|$ by the triangle inequality, so $|w-z| \ge (1/3) | x-y|$ and we get
   $$
    \cov( h_\eps(x), h^*_r(y)) \le - \log |x-y| + O(1)
   $$
    as desired in this case.

  The second case (b) is when $r \ge |x-y|/3$. Then by translation and scaling so that $B(y,r)$ becomes $B(0,1)$, the right hand side of \eqref{roughcov1} is equal to
  $$
  \log (1/r) + \iint -\log |w-z| \theta_{\frac{x-y}{r}, \frac{\eps}{r}} (\dd w) \theta^*(\dd z)
  $$
Conditioning on $w$ (which is necessarily in $\bar B(0,4)$ under the assumptions of case (b)), we see that by the assumption \eqref{eq:assumptionrho} on $\theta^*$, the second term is bounded by $O(1)$, uniformly, so that
$$
    \cov( h_\eps(x), h^*_r(y)) \le - \log r + O(1)
   $$
   as desired in this case. This proves \eqref{eq:roughcov}.
  %Possibly subcase of (b) is when $|x-y | \le \eps$.

  The proof of \eqref{eq:finecov} is similar but simpler. Indeed, we get (as in \eqref{roughcov1}),
  \begin{equation}\label{finecov1}
  \cov(h_\eps(x), h^*_\delta(y)) = \iint - \log | w-z| \theta_{x,\eps}(\dd w) \theta^*_{y, \delta}(\dd z) + g(x,y) + o(1)
  \end{equation}
  where the $o(1)$ term tends to 0 as $\eps, \delta \to 0$, coming from the continuity of $g$, and hence is uniform in $x,y$ (not even assuming $|x-y | \ge \eta$). Now note that
  $$
  \big| \log | w- z | - \log |x- y | \big| \le  \frac{4 \max(\eps, \delta)}{|x-y|}
  $$
  as soon as $\max(\eps, \delta ) \le \eta/4 \le |x-y |/4$. Therefore the right hand side of \eqref{finecov1} is $- \log | x- y |  +  g(x,y) + O( \max (\eps, \delta)) + o(1) $ when $|x-y | \ge \eta$, which proves the claim \eqref{eq:finecov}.
\end{proof}

\begin{lemma}
\label{L:UI2}
For $\alpha>\gamma$ sufficiently close to $\gamma$, $J_\eps$ is bounded in $L^2(\P)$ and hence uniformly integrable.
\end{lemma}

\begin{proof} By Fubini's theorem,
\begin{align}
\E (J_\eps^2) & = \int_{S\times S} \E( e^{\bar h_\eps(x) + \bar h_\eps(y)}\indic{G_\eps(x) \cap G_\eps(y)})  \sigma(\dd x )\sigma(\dd y) \nonumber \\
& = \int_{S\times S}   e^{\gamma^2 \cov (h_\eps(x), h_\eps(y))}
\tilde \P ( G_\eps(x)\cap  G_\eps(y)) \sigma(\dd x)\sigma( \dd y) \label{eq:twopointcorrel}
\end{align}
where $\tilde \P$ is a new probability measure obtained by the Radon--Nikodym derivative
$$
\frac{\dd \tilde \P}{\dd \P} = \frac{e^{\bar  h_\eps(x) + \bar h_\eps(y)}}{\E(e^{\bar  h_\eps(x) + \bar h_\eps (y)}) }.
$$
\begin{comment}
Observe that for $|x - y| \le 3 \eps$ say, by Cauchy--Schwarz and \eqref{covbound},  $ \cov (h_\eps(x), h_\eps(y)) \le \log(1/\eps) + g(x,x) + o(1)$. On the other hand, for $|x- y | \ge 3 \eps$,
\begin{align}
\cov (h_\eps(x), h_\eps(y)) &= g(x,y) + o(1)  + \iint \log( |w-z  |^{-1} ) \theta_\eps(x - w) \theta_\eps(y - z) dw dz \nonumber \\
&  = \log(1/|x-y| ) + g(x,y) + o(1). \nonumber
\end{align}
\end{comment}
Note that since $\sigma$ has no atoms, we may assume that $x \neq y$.
% By \cref{L:cov} (more precisely by \eqref{eq:roughcov})
%\begin{equation}\label{twopointcorr}
%\cov(h_\eps(x), h_\eps(y) ) =  - \log(|x-y| \vee \eps) + g(x,y) + O(1).
%\end{equation}
Also, if $\eps \le e^{-1} \eps_0$ and $|x- y | \le e^{-1} \eps_0  $ (else we bound the probability below by one), we have
$$
\tilde \P ( G_\eps(x) \cap  G_\eps(y)) \le \tilde \P ( h^*_{r} (x) \le \alpha \log 1/r)
$$
where
\begin{equation}\label{eps'}
r = e^{-n}, \text{where } n = \lceil \log (\frac{1}{\eps \vee |x-y|}) \rceil.
\end{equation}
Furthermore, by Cameron--Martin--Girsanov, under $\tilde \P$ we have that $h^*_r(x)$ has the same variance as before (therefore $\log 1/r + O(1)$) and a mean given by
\begin{equation}\label{twopointcorr2}
\cov_{\P} ( h^*_{r}(x), \gamma h_\eps(x) + \gamma h_\eps(y)) = 2\gamma \log 1/r + O(1),
\end{equation}
again by \cref{L:cov} (more precisely, by \eqref{eq:roughcov}).
%%%%%%%%
%
%%%%%%%%%
Consequently,
\begin{align}
 \tilde \P ( h^*_{r} (x) \le \alpha \log 1/r) & = \P( \cN( 2\gamma \log (1/r), \log 1/r) \le \alpha \log (1/r) + O(1)) \nonumber \\
& \le \exp ( - \frac12 (2\gamma - \alpha)^2 (\log (1/r) + O(1))) = O(1) r^{(2\gamma - \alpha)^2/2}.\label{corr_naive}
\end{align}
We deduce
\begin{align}
\E(J_\eps^2) & \le O(1) \int_{S\times S} |(x-y )\vee \eps|^{ ( 2\gamma - \alpha)^2/2-\gamma^2}  \sigma(\dd x )\sigma(\dd y).\label{J2naive}
\end{align}
(We will get a better approximation in the next section).
Clearly by \eqref{dim} this is bounded if
$$
( 2\gamma - \alpha)^2/2-\gamma^2 > - \mathbf{d}
$$
and since $\alpha$ can be chosen arbitrarily close to $\gamma$ this is possible if
\begin{equation}\label{UI}
\mathbf{d} - \gamma^2 /2 >0 \text{ or } \gamma < \sqrt{2\mathbf{d}}.
\end{equation}
This proves the lemma.
\end{proof}

To finish the proof of \cref{P:UI}, observe that $I_\eps = J_\eps + J'_\eps$.
We have $\E(J'_\eps) \le p(\eps_0)$ by  \cref{L:Liouvillethick},  and for a fixed $\eps_0$, $J_\eps$ is bounded in $L^2$ (uniformly in $\eps$). Hence $I_\eps$ is uniformly integrable.
\end{proof}

\subsubsection{Convergence}

As before, since $\E(J'_\eps)$ can be made arbitrarily small by choosing $\eps_0$ sufficiently small, it suffices to show that $J_\eps$ converges in probability and in $L^1$. In fact we will show that it converges in $L^2$, from which convergence will follow. To do this we will show that $(J_\eps)_\eps$ forms a Cauchy sequence in $L^2$, and we start by writing
\begin{equation}\label{basic}
\E( (J_\eps - J_\delta)^2 ) = \E(J_\eps^2 ) + \E ( J_\delta^2) - 2 \E(J_\eps J_\delta).
\end{equation}
Our basic approach is thus to estimate, better than before, $\E(J_\eps^2)$ from above and $\E(J_\eps J_\delta)$ from below. Essentially, the idea is that for $x,y$ which are at a small but macroscopic distance, we can identify the limiting distribution of $(h^*_{r}(x), h^*_{r}(y))_{r\le \eps_0}$ under the distribution $\P$ biased by $e^{\bar h_\eps(x) + \bar h_\delta(y)}$. On the other hand when $x,y$ are closer than that we know from the previous section that the contribution is essentially negligible.

Let us write, similar to \eqref{eq:twopointcorrel},
$$
\E(J_\eps J_\delta) = \int_{S^2}  e^{\gamma^2  \cov( h_\eps(x), h_\delta(y))} \tilde \P( G_\eps(x)\cap G_\delta(y)) \sigma(\dd x)\sigma( \dd y),
$$
where $\dd \tilde \P (\cdot)  = \tfrac{e^{ h_\eps(x) +  h_\delta(y)}}{ \E (e^{ h_\eps(x) +  h_\delta(y)}) } \dd \P (\cdot) .$ We fix $\eta>0$ arbitrarily for now, and suppose that $|x - y | \ge \eta$. 

\medskip Observe that for any fixed $\eps_1 \le \eps_0$, as $\eps \to 0$, and uniformly over $x,y \in S$ with $|x-y| \ge \eta$ and $r \ge \eps_1$, if $z \in\{x,y\}$, 
\begin{equation}\label{shift1}
\cov( h^*_{r}(z), h_\eps(x) ) \to \int_D K(w,x) \theta^*_{r}(w-z) \dd w.
\end{equation}
Likewise, uniformly over $x, y \in S$ such that $|x-y|\ge \eta$, and over $r \ge \eps_1$, if $z\in \{x,y\}$, as $\delta \to 0$:
\begin{equation}\label{shift2}
\cov(h^*_{r}(z), h_\delta(y)) \to \int_D K(w,y) \theta^*_{r}(w-z) \dd w.
\end{equation}
%\ellen{I don't see the difference between these two equations??} 
(Note that both right hand sides of \eqref{shift1} and \eqref{shift2} are finite by \eqref{eq:roughcov}.) Consequently, by Cameron--Martin--Girsanov, as $\eps \to 0, \delta \to 0$, the joint law of the processes $(h^*_{r}(x), h^*_{r}(y))_{r\le \eps_0}$ under $\tilde \P$ converges to a joint distribution
$(\tilde h^*_{r}(x), \tilde h^*_{r}(y))_{r \le \eps_0}$ whose covariance is unchanged and 
whose mean is given by 
\begin{equation}\label{eq:addshift}
\E( \tilde h^*_r(z) ) = \gamma [\int_D K(w,x) \theta^*_{r}(w-z) \dd w + \int_D K(w,y) \theta^*_{r}(w-z) \dd w], \quad z \in \{x,y\}, r\le \eps_0.
\end{equation}
This convergence is for the weak convergence on compacts of $r \in (0, \eps_0]$, and is uniform in $|x-y| \ge \eta$.

The core of the argument will be to prove that $\tilde \P(G_\eps(x) \cap G_\delta(y))$ converges as $\eps, \delta \to 0$, uniformly over $|x-y| \ge \eta$, to the analogous probabilities for $( \tilde h^*_r(x), \tilde h^*_r(y))$. To state this more precisely, let us introduce some notation. Define the event $\tilde E_n(z)$ (for $z \in \{x,y\}$) in a way analogous to the event $E_n(z)$: 
$$
\tilde E_n(z) = \{ \tilde h^*_{e^{-n}}(z) \le\alpha n \}; \quad n \in \Z,
$$
and consider the corresponding good event for the field $\tilde h$,
$$
\tilde G(z) = \bigcap_{n=n_0}^\infty \tilde E_n(z). 
$$

\begin{lemma}
\label{L:prooftwopoints} As $\eps \to 0, \delta \to 0$, uniformly over $|x- y | \ge \eta$, 
 \begin{equation}\label{corr_better}
\tilde \P( G_\eps(x) \cap G_\delta(y)) \to g_\alpha(x,y): = \P( \tilde G (x) \cap  \tilde G(y)) .
\end{equation}
\end{lemma}

\begin{proof}
This will follow from the fact that the joint law of the processes $(h^*_{r}(x), h^*_{r}(y))_{r\le \eps_0}$ under $\tilde \P$ converges on compact sets of $(0, \eps_0]$ to the joint distribution of 
$(\tilde h^*_{r}(x), \tilde h^*_{r}(y))_{r \le \eps_0}$, but requires an argument since the good events $\tilde G(x), \tilde G(y)$ do not depend only on the behaviour of $(h^*_r(x), h^*_r(y))$ in some fixed compact of $(0, \eps_0]$. 

Let us start with the upper bound for \eqref{corr_better}. For any $n_1 > n_0$, and any $\eps, \delta$ such that $n(\eps) , n(\delta) \ge n_1$,
\begin{align*}
\tilde \P( G_\eps(x) \cap G_\delta(y))& 
%= \tilde \P ( (\bigcap_{n=n_0}^{n(\eps)} E_n (x)) \cap ( \bigcap_{n=n_0}^{n(\delta)}E_n (y)))
 \le \tilde \P ( \bigcap_{n=n_0}^{n_1} E_n (x) \cap E_n(y)) 
\end{align*}
so using the convergence on compact sets, 
\begin{align*}
\limsup_{\delta, \eps\to 0} \tilde \P( G_\eps(x) \cap G_\delta(y)) &\le  \P ( \bigcap_{n=n_0}^{n_1} \tilde E_n (x) \cap \tilde E_n(y))
\end{align*}
and since $n_1>n_0$ is arbitrary we get
$$
\limsup_{\eps, \delta \to 0} \tilde \P( G_\eps(x) \cap G_\delta(y)) \le \P ( \tilde G(x) \cap \tilde G(y)).
$$
Let us now turn to the lower bound. 
The key is to observe that the constraint $ h^*_{e^{-n}}(z) \le \alpha n $ defining $E_n(z)$ is essentially guaranteed for large $n$ under $\tilde \P$, whether $z$ is $x$ or $y$. This is because $|x-y|\ge \eta$ so $x$ and $y$ are well separated (note that this would be false however if $x$ and $y$ would not be separated: for instance if $x = y$, we would typically expect $h^*_{e^{-n}} (x)$ to be of the order of $2\gamma n$, which is $> \alpha n$.) 
More precisely, applying the same argument as in Lemma \ref{L:Liouvillethick}, no matter how small $a>0$ is, we can find $n_1= n_1( \eta, \alpha, \gamma,a)$, but independent of $\eps$ or $\delta$, such that for all $\eps, \delta>0$,
under $\tilde \P$, 
\begin{equation}\label{eq:22}
\tilde \P(\bigcup_{n=n_1}^{ n(\delta) \vee n (\eps)} \{h^*_{e^{-n}} (z) \ge \alpha n\} )\le a. 
\end{equation}
Indeed, consider $z = x$ first. Lemma \ref{L:Liouvillethick} analyses the effect of biasing by $e^{\bar h_\eps(x)}$. To analyse the effect of further biasing by $e^{\bar h_\delta(y)}$, observe that since $x$ and $y$ are separated by a distance at least $\eta$, the resulting additional shift in the mean of $ h^*_{e^{-n}}(x)$, which is explicitly given by $\gamma \cov (h^*_{e^{-n}}(x), h_\delta(y) ) $, is at most $O(1)$ by \eqref{eq:roughcov}. The same observation applies to $z = y$. 

Therefore, using \eqref{eq:22},
$$
\tilde \P( G_\eps(x) \cap G_\delta(y) ) \ge \tilde \P ( \bigcap_{n=n_0}^{n_1} E_n(x) \cap E_n(y) ) - 2a. 
$$
We can take a limit as $\delta, \eps \to 0$ using convergence on compacts to deduce 
\begin{align*}
\liminf_{\eps \to 0} \tilde \P( G_\eps(x) \cap G_\delta(y) ) &\ge  \P ( \bigcap_{n=n_0}^{n_1}\tilde E_n(x) \cap \tilde E_n(y) ) - 2a\\
& \ge \P ( \bigcap_{n=n_0}^{\infty}\tilde E_n(x) \cap \tilde E_n(y) ) - 2a.
\end{align*}
This completes the proof of \eqref{corr_better} because $a>0$ was arbitrary. 
\end{proof}

Using this lemma we can easily conclude:

\begin{lemma}
\label{L:limsup}
We have
$$
\limsup_{\eps\to 0} \E( J_\eps^2) \le \int_{S\times S} e^{\gamma^2 g(x,y)} \frac1{|x-y|^{\gamma^2}}  g_\alpha(x,y)\sigma(\dd x)\sigma( \dd y)
$$
where $g_\alpha(x,y)$ is a non-negative function depending on $\alpha, \eps_0$ and $\gamma$ such that the above integral is finite.
\end{lemma}
\begin{proof}
We fix $\eta >0$ arbitrarily small (in particular, $\eta$ may and will be smaller than $e^{-1}\eps_0$). If $|x-y| \le \eta$ we use the same bound as in \eqref{J2naive}. The contribution coming from the part $|x-y| \le \eta$ can thus be bounded, uniformly in $\eps$, by $f(\eta) $ (where $f(\eta) \to 0$ as $\eta \to 0$ and the precise order of magnitude of $f(\eta)$ is determined by \eqref{dim}, and is at most polynomial in $\eta$). 
%We thus focus on the contribution coming from $|x- y| \ge \eta$.

On the other hand, for $|x-y| \ge \eta$, taking $\delta = \eps$ in \eqref{L:prooftwopoints}, 
%Indeed, by \eqref{eq:roughcov}, under $\tilde \P$, the drifts of $h_r(x)$ and of $h_r(y)$ (with $ r \ge \eps$) are each $\gamma \log (1/r) + O(1)$ where the $O(1)$ term is uniform in $|x-y | \ge \eta$. Because of this, up to an error in $\tilde \P$ probability that is arbitrarily small (uniformly in $|x-y | \ge \eta$), the events $G_\eps(x), G_\eps(y)$ as well as $\tilde G(x), \tilde G(y)$ depend only on the ``macroscopic" behaviour of $h_r(x)$ and $h_r(y)$; that is, depend only on $(h_r(x), h_r(y))_{r \ge \eps_1}$ for some $\eps_1$.
%Indeed the upper bound is obvious from the weak convergence, and the lower bound follows from the obvious coupling between the law of $(h_{r}(x), h_{r}(y))_{r\le \eps_0}$ under $\tilde \P$ and its limit $(\tilde h_r(x), \tilde h_r(y))_{r\le \eps_0}$, as well as the convergence of the means in \eqref{shift1} and \eqref{shift2}. \note{In fact we need the convergence to be uniform in $r$, which I think is the case when $|x- y| \ge \eta$.}
after applying \cref{L:cov} (more specifically \eqref{eq:finecov}) for the pointwise limit and \eqref{eq:roughcov} for the use of dominated convergence):
\begin{equation}\label{convJ2}
\int_{S^2; |x- y|\ge \eta} \!\!\!\!\!\!\!\!\!\!\!\!\!\!\!\!{e^{\gamma^2 \cov( h_\eps(x), h_\eps(y))}} \tilde \P( G_\eps(x), G_\eps(y))\sigma( \dd x ) \sigma( \dd y )\to \int_{S^2 ; |x- y | \ge \eta} \frac{e^{ \gamma^2 g(x,y)}}{|x-y|^{\gamma^2}} g_\alpha (x,y) \sigma(\dd x) \sigma(\dd y).
\end{equation}
Since we already know that the piece of the integral coming from $|x- y | \le \eta$ contributes at most $f(\eta) \to 0$ when $\eta \to 0$, it remains to check that the integral on the right hand side of \eqref{convJ2} remains finite as $\eta \to 0$. But we have already seen in \eqref{corr_naive} that for $|x- y | \le \eps_0/3$, $\tilde \P( G_\eps(x) \cap  G_\eps(y)) \le O(1) |x-y|^{( 2\gamma - \alpha)^2/2-\gamma^2}$; hence this inequality must also hold for $g_\alpha (x,y)$. Hence the result follows as in \eqref{UI}.
\end{proof}

\begin{lemma}
\label{L:liminf}
We have
$$
\liminf_{\eps,\delta \to 0} \E( J_\eps J_\delta ) \ge \int_{S\times S} e^{\gamma^2 g(x,y)} \frac1{|x-y|^{\gamma^2}}  g_\alpha(x,y) \sigma(\dd x) \sigma( \dd y).
$$
\end{lemma}

\begin{proof}
In fact, the proof is almost exactly the same as in \cref{L:limsup}, and is even easier, since we get a lower bound by restricting ourselves to $|x- y |  \ge \eta$. We deduce immediately from Lemma \ref{L:prooftwopoints} that
$$
\liminf_{\eps,\delta \to 0} \E( J_\eps J_\delta ) \ge \int_{S^2; | x- y  | \ge \eta} e^{\gamma^2 g(x,y)} \frac1{|x-y|^{\gamma^2}}  g_\alpha(x,y) \sigma(\dd x) \sigma( \dd y).
$$
Since $\eta$ is arbitrary, the result follows.
\end{proof}
\begin{proof}[Proof of convergence in \cref{T:conv}]
Using \eqref{basic} together with \cref{L:limsup,L:liminf}, we see that $J_\eps$ is a Cauchy sequence in $L^2$ for any $\eps_0>0$.
Combining with \cref{L:Liouvillethick}, it therefore follows that $I_\eps$ is a Cauchy sequence in $L^1$ and hence converges in $L^1$ (and also in probability) to a limit $I = \cM(S)$. The proof of weak convergence follows by the argument in \cref{S:weakconvergence}. {We leave the proof that $\cM$ has no atoms as an exercise: see Exercise \ref{Ex:thick} or Exercise \ref{Ex:atom} for another proof.}
 \end{proof}

\begin{rmk}
Note that $\lim_{\eps\to 0} \E(J_\eps^2)$ depends on the regularisation $\theta$,
even though, as we will see next, $\lim_{\eps \to 0} I_\eps$ does not.
\end{rmk}

\paragraph{Proof of uniqueness in \cref{T:conv}.} To prove uniqueness, we take $\tilde \theta$ another nonnegative Borel measure on $\R^d$ supported in $\bar B(0,1)$ and satisfying \eqref{eq:assumptionrho}. Let $\tilde h_\delta (x) = h * \tilde \theta_\delta(x)$, and let $\tilde J_\delta$ be defined as $J_\delta$ but with $\tilde \theta$ instead of $\theta$: that is,
$$
\tilde J_\delta  = \int_S e^{\gamma \tilde h_\delta(z) - (\gamma^2/2) \E( \tilde h_\delta(z)^2)} \mathbf{1}_{\{  G_\delta(z)\}} \sigma(\dd z)
$$
where we use the same good event $ G_\delta(z)$ as in \eqref{good}, based on the regularisation kernel $\theta^*$. Then the argument of \cref{L:liminf} can be used to show that the same conclusion holds for $J_\delta$ replace by $\tilde J_\delta$: that is,
$$
\liminf_{\eps,\delta \to 0} \E( J_\eps \tilde J_\delta ) \ge \int_{S\times S} e^{\gamma^2 g(x,y)} \frac1{|x-y|^{\gamma^2}}  g_\alpha(x,y) \sigma(\dd x) \sigma( \dd y).
$$
Hence we deduce $\lim_{\eps \to 0, \delta \to 0} \E( (J_\eps - \tilde J_\delta)^2) = 0$ and this implies that the limits associated with $\theta$ and $\tilde \theta$ are almost surely the same.

%{Have used implicitly later on (GFF case) that (any restriction of) $\cM$ or $\nu$ is (locally) measurable with respect to the realisation of the field as a distribution. I think from this section it is clear in the bulk case. In the boundary case I think we may need some extra clarification. Probably better later on after Neumann GFF actually introduced...}

\begin{rmk}\label{R:extensionGMC}
The proof of Theorem \ref{T:conv} given above extends without difficulty to a variety of settings going somewhat beyond the stated assumptions. In such cases, the input that is crucially required for the argument to extend without major modifications is Lemma \ref{L:cov} which controls the correlations of the regularised Gaussian field. An example would be the setting of a Gaussian, logarithmicall correlated field on a Riemannian manifold. 
\end{rmk}

\subsection{Shamov's approach to Gaussian multiplicative chaos}\label{S:shamov}

An alternative {and powerful} viewpoint on Gaussian multiplicative chaos was also developed in Shamov \cite{Shamov}. It is closely related to the generalisation of ``rooted measures'' for the GFF: see \cref{sec:rooted_meas}. In what follows $h$ will be a centred Gaussian field with logarithmically diverging covariance kernel $K$ as in \eqref{setup_gmc} (although the original paper \cite{Shamov} works in a more general setting).

Before stating the result, let us make an observation about changes of measure for the field $h$. If $\rho \in \mathfrak{M}$ we write $K\rho$ for the linear operator $\rho'\mapsto K(\rho,\rho')$ on $\mathfrak{M}$. Note that if $\rho\in \cD_0(D)$ we have
\begin{equation}\label{Trho}
K\rho(x)=\int_D K(x,y)\, \rho(\dd y).
\end{equation} 
%\nb{I think we should call it $K\rho$, not $T\rho$. I also assume that $\rho$ is smooth, else this doesn't work.}

% 
Applying Girsanov's Lemma (\cref{L:Girsanov}), we see that if $\rho \in \mathfrak{M}$ then the field $h+K\rho$ (viewed as a stochastic process indexed by $\rho'\in \mathfrak{M}$) is absolutely continuous with respect to $h$, with associated Radon--Nikodym derivative 
 \begin{equation} 
 \label{eq:CM_general}\frac{\exp((h,\rho))}{\exp(\frac{1}{2}{(\rho,\rho)_K})}.
 \end{equation} 
%\nb{Related: should this be $\|\rho\|_K$ in the denominator? We can only identify this with $(\rho,T\rho)$ under the energy condition at this stage.} 
%
 Note the connection with \cref{sec:ac} in the case of the zero boundary GFF: when $\rho\in \mathfrak{M}_0$ ($\mathfrak{M}_0$ corresponding to the zero boundary condition Green function) then $(h,\rho)=(h, F )_\nabla$, where $F$ is defined by $-\Delta F = 2\pi \rho$ and is an element of $H_0^1(D)$. By \eqref{eqn:frho} this is exactly the statement that $F=K\rho$, and  the above expression is equal to $\exp(( h, F)_\nabla)/\exp(\frac{1}{2} (F, F )_\nabla)$ as in \cref{lem:CMGFF}. See \cite{Arunotes} for more on Shamov's approach when the field is the planar Dirichlet GFF.

\begin{definition}[Shamov's definition of GMC]\label{def:gmc_shamov}
	Let $h$ be as above and $\sigma$ as in \eqref{dim}. Let $\gamma\in (0,2)$. A measure $\cM^\gamma$ is a $\gamma$-multiplicative chaos measure for $h$, with background measure $\sigma$ if:
	\begin{itemize}
		\item $\cM^\gamma$ is measurable with respect to $h$ as a stochastic process indexed by $\mathfrak{M}$ (note that this  allows us to write $\cM^{\gamma}(\dd x)=\cM^{\gamma}(h,\dd x)$);
		\item $\E(\cM^\gamma(S))=\sigma(S)$ for all Borel sets $S\subset D$;
		\item For every fixed (deterministic) %deterministic Borel measurable function 
		function $\rho\in \mathcal{D}_0(D)$ (i.e. $\rho$ is smooth and compactly supported), if $\xi=K\rho$ is given by \eqref{Trho} then %such that $\xi(x) =T\rho(x)$ $\sigma$-almost everywhere with $\rho\in \mathfrak{M}$,
			\begin{equation}\label{eq:shamov1} \cM^\gamma(h+\xi,\dd x) = \exp(\gamma \xi(x)) \cM^\gamma(h,\dd x) \text{ almost surely.} \end{equation}
	\end{itemize}
\end{definition}

We use the notation $\mathcal{M}^\gamma$ above to distinguish it from $\mathcal{M}$ in the previous sections. However, we will see just below that $\cM^\gamma$ exists, and in fact must be equal to $\mathcal{M}$.

%Note that although $\xi$ is only defined almost everywhere with respect to $\sigma$ (for example when the field is a GFF with Dirichlet boundary conditions in $D$, then $\xi$ will only be an element of $H_0^1(D)$), the measure $\exp(\gamma \xi(x))\cM^{\gamma}(h,\dd x)$ still makes sense unambiguously.

%\nb{Actually I am a bit uncertain about this. After all when we define $\xi = T\rho$ we have selected a particular a.e. representative, haven't we ? So perhaps the condition should be for every $\xi = T\rho$, $\sigma-$a.e.? This way this translates exactly into $\xi \in H_0^1$ for the GFF? }

%Indeed, the assumption that $\E(\cM^{\gamma})=\sigma$ implies that if {$\tilde \xi$} is such that $\xi(x)={\tilde \xi(x)}$ for $\sigma$-almost every $x$, then by Fubini's theorem

%$$
%\E( \int_S  1_{\{\xi(x) \neq \tilde \xi (x)\}}\exp( {\gamma \xi(x)}) \cM^{\gamma} (h,\dd x))=0.
%$$
%It follows that $ \int_A  1_{\{\xi(x) \neq \tilde \xi (x)\}}\exp( {\gamma \xi(x)}) \cM^{\gamma} (h,\dd x) = 0 $ almost surely simultaneously for all Borel sets $A \subset S$.
%This implies that on an event of probability one,
%$$\int_A \exp( {\gamma \tilde \xi(x)})\cM^{\gamma} (h,\dd x)=\int_A \exp(\gamma \xi(x))\cM^{\gamma}(h,\dd x)$$ for all $A \subset S$. Hence the measures $\cM^\gamma(h+ \xi, \dd x) $ and $\cM^\gamma( h+ \tilde \xi, \dd x)$ agree with probability one, and so are unambiguously defined.
%\nb{I think a $\pi$-system argument is needed to say that this holds for all sets simultaneously. Also, do we want to say that when $h$ is the GFF it suffices to check it for smooth functions $\xi$?}

\begin{theorem}[Shamov, \cite{Shamov}]\label{thm:shamov}
	Assume the setup of \cref{def:gmc_shamov}. Then a multipli\-cative chaos measure for $h$ with parameter $\gamma$ and background measure $\sigma$ exists. Moreover, it is unique.
\end{theorem}

We note that the uniqueness part of \cref{thm:shamov} may be particularly useful if one wants to identify some limit as being a GMC measure, since the conditions are in many contexts relatively easy to check. %{Actually, these conditions can be slightly weakened so as to restrict $\xi$ to an appropriately dense subspace; for instance, in the case where $h$ is the GFF with Dirichlet boundary conditions, it suffices to know \eqref{eq:shamov1} for smooth functions with compact support, see \cref{R:shamovdense}.}

\begin{rmk}
As we will see in the proof below, the condition \eqref{eq:shamov1} ensures that the effect of weighting the law of the field by $\cM^\gamma(D)$ is to add the singularity $\gamma K(x,\cdot)$ to the field at a point $x$ chosen from $\cM^{\gamma}$, a property which we will shall see amounts to Girsanov's transform for the field reweighted by the mass of $\cM^\gamma(D)$. So essentially, Shamov's approach characterises the GMC measure as a certain Radon--Nikodym derivative for the field.
\end{rmk}

\begin{proof}

(i). \emph{Proof of existence.}
For the existence part we will show that the GMC measure constructed in \cref{T:conv} does satisfy the stated conditions. The first two properties are an obvious consequence of the construction: in particular, since the GMC measure $\cM$ is obtained as a limit in probability as $\eps \to 0$, with respect to the weak convergence of measures, of a sequence of measures $\cM_\eps$ that are obviously measurable with respect to $h$, so is their limit $\cM$. The third property requires some additional arguments.  {First, since {$\rho \in \cD_0(D)$}, by \eqref{eq:CM_general} we have that} $h + \xi = h + K\rho $ is absolutely continuous with respect to $h$. {In particular, it makes sense to apply Theorem \ref{T:conv} to $\tilde h = h +\xi$ and we have $\cM^\gamma(h+\xi,\dd x)=\lim_{\eps\to 0} \eps^{\gamma^2/2} e^{\gamma(h_\eps+\xi_\eps)(x)} \dd x=\lim_{\eps\to 0} e^{\gamma \xi_\eps(x)}\cM_\eps(\dd x)$ and $\cM^\gamma(h,\dd x)=\lim_{\eps\to 0} e^{\gamma h_\eps(x)} \dd x= \lim_{\eps\to 0} \cM_\eps(\dd x)$, where the limits are {in probability} and $f\mapsto f_\eps{= f * \theta_\eps}$ represents smoothing with some fixed mollifier at scale $\eps$. Hence, it suffices to show that $e^{\gamma \xi_\eps}\cM_\eps\to e^{\gamma \xi} \cM$ {in probability}, with respect to the topology of weak convergence as $\eps\to 0$. This follows by the triangle inequality since $e^{\gamma \xi_\eps}\to e^{\gamma \xi}$ uniformly on compacts of $D$ as $\eps\to 0$, {as $\xi = K\rho$ is continuous on $\bar D$, by a simple application of dominated convergence (recall that $\rho \in \cD_0(D)$)}. }

%(b) Let $\tilde h = h + \xi$. Then $\tilde h_\eps = h_\eps + \xi_\eps$. Thus $\tilde \cM_\eps (\dd x)= e^{\gamma \xi_\eps(x)} \cM_\eps(\dd x)$. We need to check that $\xi_\eps (x) \to \xi(x)$ (I know a reason why $\xi_\eps (x) \to \xi(x)$ almost everywhere wrt $\rho$ , but additional arguments are needed to conclude that $e^{\gamma \xi_\eps(x)} \cM_\eps(\dd x)$ converges to $e^{\gamma \xi(x)} \cM(\dd x)$.

\medskip 

%(In particular, given the uniqueness of Theorem \ref{thm:shamov}, it follows that the measures of  \cref{thm:shamov} and  \cref{T:conv} are the same). It remains to prove the uniqueness.

(ii). \emph{Proof of uniqueness.}
	Suppose that a measure $\cM^{\gamma}$ satisfying the constraints of \cref{def:gmc_shamov} exists. We will consider the probability measure (often called the \textbf{rooted measure}, and already encountered in Chapter \ref{S:Liouvillemeasure} in the context of the Liouville measure associated to the Dirichlet GFF in \eqref{eq:rootedChap2}): \ind{Rooted measure}
	\begin{equation}
	Q(\dd x, \dd h)=\frac{ \cM^{\gamma}(h,\dd x)}{\E(\cM^{\gamma}(h,D))} \P(\dd h).
	\end{equation}
	Note that $\int_D \int Q ( \dd x, \dd h) = \E_h ( \int_D \cM^\gamma(h, \dd x) ) / \E_h ( \cM^\gamma(h, D)) = 1$ so that $Q$ is a probability law on pairs $(x, h)$. 
	 We will show that:
	 \begin{enumerate}
\item 	  Under $Q$, the marginal law of $x$ has density proportional to $\sigma(\dd x)$, 
\item Given $x$, the conditional law of the field (viewed as a stochastic process indexed by $\mathfrak{M}$) 
is that of $h$ plus the deterministic function $\gamma K(x,\cdot)$. 
\end{enumerate}
Observe that these two properties completely characterise the law $Q$ and thus, by disintegration, the conditional law $Q( \dd x | h)$ of $x$ given $h$ under $Q$. On the other hand, the definition of $Q$ means that this conditional law is exactly $\cM^{\gamma}(h,dx)$ and so we have identified $\cM^{\gamma}$ uniquely (note that this doesn't identify only the law of $\cM^\gamma$ but really the joint law of $h$ and $\cM^\gamma$).
	
	To show the claim concerning $Q$, it is enough to prove that the $Q$ marginal law of $x$ is equal to $\sigma(\dd x)/\sigma(D)$, and that for any $\rho_1,\cdots, \rho_m \in \mathfrak{M}$ and $a_1, \ldots, a_m\in \R$ the $Q$ conditional law of $(a_1(h,\rho_1)+\ldots+a_m(h,\rho_m))$ given $x$ is a normal random variable with the correct mean and covariance.  In other words (using linearity of $h$ on the space $\mathfrak{M}$) it suffices to show that for any $g\in L^1(\sigma)$ on $D$, and $\rho\in \mathfrak{M}$
	\begin{equation}\label{E:Q}
	\E_Q(\e^{(h,\rho)} g(x))  =  \E( \int_D \e^{(h+\gamma K(x,\cdot),\rho)} g(x) \frac{\sigma(\dd x)}{\sigma(D)}). \end{equation}
	{In fact, it follows from the definition of $\mathfrak{M}$ that for $\rho\in \mathfrak{M}$, $(h,\rho)$ is the limit {in probability} (under $\mathbb{P}$ but therefore also under $Q$ by absolute continuity) of $(h,\rho_n)$ as $n\to \infty$, for some appropriate sequence $\rho_n$ with $\rho_n \in \cD_0(D)$ for all $n$. To characterise the law of $h$ under $Q$ it therefore  suffices to show \eqref{E:Q} for $\rho\in \cD_0(D)$.}
	
{Let us conclude the proof by showing this.}	Note that by Fubini's theorem, when $\rho \in \cD_0(D)$ the right hand side of \eqref{E:Q} is equal to
$$ \sigma(D)^{-1} \int_D \e^{\frac{1}{2}\var((h,\rho))+\gamma \int K(x,y)\rho(\dd y) } g(x)\sigma(\dd x)
 = \sigma(D)^{-1} \int_D \e^{\frac{1}{2}(\rho, \rho)_K+\gamma K\rho(x) } g(x) \sigma(\dd x)$$ (recalling the notation $K\rho$ in \eqref{Trho}).
	Furthermore, the left hand side of \eqref{E:Q}  (using the assumption that $\E(\cM^\gamma(h,D))=\sigma(D)$ and the definition of $Q$) is equal to
	$$ \sigma(D)^{-1} \E(\int_D \e^{(h,\rho)} g(x) \cM^{\gamma}(h,\dd x)).$$
%    So the real aim is to show that
%\begin{equation}
%\E(\int_D\e^{(h,\rho)}g(x) \cM_{\gamma}(h,dx)) = \int_D \e^{\frac{1}{2}(T\rho, \rho)+\gamma T\rho(x) } g(x) \sigma(dx).
%\end{equation}	
However, using the observation  \eqref{eq:CM_general} and {the property} \eqref{eq:shamov1}, we have
\begin{align}
\E(\int_D\e^{(h,\rho)}g(x) \cM^{\gamma}(h,\dd x)) & = \E( \int_D \e^{\tfrac{1}{2}(\rho,\rho)_K}  g(x) \cM^{\gamma}(h+K\rho, \dd x)) \nonumber \\ & = \E(\int_D \e^{\tfrac{1}{2}(\rho,\rho)_K}\e^{\gamma K\rho(x)} g(x) \cM^{\gamma}(h,\dd x)) \nonumber \\
& = \int_D \e^{\frac{1}{2}(\rho, \rho)_K+\gamma K\rho(x) } g(x) \sigma(\dd x),\label{eq:shamovproof}
\end{align}
where in the last line we again used the assumption that $\E(\cM^{\gamma}(h,\dd x))=\sigma(\dd x)$. Dividing by $\sigma(D)$ this is the same as the right hand side of \eqref{E:Q}, so we get the desired result.
\end{proof}

\subsection{Rooted measures and Girsanov lemma for GMC}

We now return to the notation $\cM$ where we suppress the dependence on $\gamma$.

Let $h$ be as in Section \ref{SS:setup} and let $\sigma$ be as in \eqref{dim}. Closely related to the previous theorem (and in particular the rooted measure appearing in its proof) is a description of the law of $h$ after reweighting by $\cM(D) / \sigma(D)$. In the case of the two dimensional Gaussian free field with Dirichlet boundary conditions, this has already been described in \eqref{condrootG}, which is a consequence of Lemma \ref{L:Girsanov}. The result in this case is that the law of the field, when biased by the total mass $\cM(D)$, can simply be described by first sampling a point $z$ with an appropriate (deterministic) law (corresponding to the appropriate multiple of $\sigma(\dd z) = R(z, D)^{\gamma^2/2} \dd z$), and then adding to $h$ a function of the form $\gamma G_D(z, \cdot)$. Since $G_D$ is nothing but the covariance of the field, it is easy to guess that such a description generalises to the broader Gaussian multiplicative chaos framework. This is indeed what the next theorem shows, which is essentially a reformulation of the work done in Theorem \ref{thm:shamov}.
\ind{Girsanov lemma!for GMC}

\begin{theorem}[Girsanov's lemma for GMC]\label{T:Girsanov}
  Let $h$ be as in Section \ref{SS:setup} and $\sigma $ as in \eqref{dim}, $\cM$ the $\gamma$-multiplicative chaos measure for $h$ with reference measure $\sigma$. Then for any $\rho \in \mathfrak{M}$, and any non-negative Borel function $g$ on $D$,
  $$
  \E[e^{ (h, \rho)} \int_Dg(x) \cM(\dd x) ] = \int_D \sigma(\dd x) g(x) \E[ e^{ ( h + \gamma K(x, \cdot), \rho)} ] .
  $$
\end{theorem}

\begin{proof}
  We note that this could be proved using the same argument
   explained in \eqref{condrootG}. However, given Theorem \ref{thm:shamov}, it is simpler to proceed as follows. We recall that $\cM$ is a Gaussian multiplicative chaos for $h$ in the sense of Shamov (Definition \ref{def:gmc_shamov}), {as explained in Theorem \ref{thm:shamov}}. Therefore, it satisfies \eqref{E:Q}, as shown in Theorem \ref{thm:shamov}. This implies the result.
\end{proof}

Since $\rho \in \mathfrak{M}$ is arbitrary and the law of $(h, \rho)$ for arbitrary $\rho$ characterises the law of a Gaussian additive process $((h, \rho))_{\rho \in \mathfrak{M}}$ uniquely, or (alternatively) applying the argument used in the proof of Theorem \ref{thm:shamov} since we already know that $\cM$ is also a GMC in the sense of Definition \ref{def:gmc_shamov},
we deduce: %(taking $g = 1$):

\begin{cor}
  \label{C:Girsanov}
 Define a law $\Q$ on pairs $(x,h)$ through the formula
 $$
 \Q (\dd x, \dd h) = \frac1{\sigma(D)} \cM_h(\dd x) \P(\dd h),
 $$
 where $\cM_h$ is the GMC measure almost surely associated with $h$. Then $\Q$ is indeed a probability measure (the \textbf{rooted} or Peyrière measure associated to $h$) and satisfies: \ind{Rooted measure}
\begin{enumerate} 
\item The marginal law of $h$ is given by $ \Q(\dd h) = (\cM_h(D)/ \sigma(D)) \P (\dd h)$. (In particular, the law of $h$ under $\Q$ is absolutely continuous with respect to $\P$.)

\item Given $h$, the point $x\in D$ is sampled according to the normalised GMC measure associated to $h$, $\cM_h(\cdot) /\cM_h(D)$. 

\item The marginal law of $x$ is $\Q ( \dd x) = \sigma(\dd x)/\sigma(D)=:\hat{\sigma} ( \dd x)$.

\item Finally, the law of $(x,h)$ under $\Q$ is the same as the law of $(x, h + \gamma K(x, \cdot))$, under $ \hat \sigma (\dd x) \otimes \P(\dd h)$.  

\end{enumerate}

 \end{cor}
 
%\ellen{ (Here and in what follows, $ \tfrac{\sigma}{\sigma(D)}\otimes \P=:  \hat{\sigma}\otimes \P$  denotes the product measure under which $z$ is sampled according to $\hat \sigma $ and $h$ is sampled according to $\P$, independently of $z$).}

\begin{cor}
	  \label{C:Girsanov2}
{Under $ \hat{\sigma} \otimes \P$, $\int_D e^{\gamma^2 K(x,y)}\cM_h(\dd y)<\infty$ almost surely. Furthermore, the law of the measure $e^{\gamma^2K(x,y)}\cM_h(\dd y)$ under $ \hat{\sigma} \otimes \P$ is equal to the $\Q$-law of $\cM_h(\dd y)$.}
\end{cor}

\begin{proof}[Proof of Corollary \ref{C:Girsanov2}]
{If we show the second point, then the first one follows automatically since $\Q$ is absolutely continuous with respect to $\P$ and thus $\cM(D)<\infty$ almost surely under $\Q$. We therefore concentrate on the second point. To this end we note that} since $\Q$ is absolutely continuous with respect to $\hat{\sigma} \otimes \P$, the measure $\cM_h$ almost surely does not have any atoms under $\Q$. In particular, it almost surely does not have an atom at the marked point $x$ and therefore, 
	\begin{equation}\label{eq:noatomX}
\cM_h(D\setminus B_\eta(x))\uparrow \cM_h(D)<\infty 
\end{equation} 
almost surely  (under $\Q$) as $\eta\downarrow 0$. 
{On the other hand, by monotone convergence, as $\eta \to 0$, 
$$\int_{D\setminus B_\eta(x) } e^{\gamma^2 K(x, y)} \cM_h( \dd y ) \uparrow  \int_{D} e^{\gamma^2 K(x, y)} \cM_h( \dd y ).$$ 
Therefore it suffices to show that for each $\eta>0$, the law of $\cM_h$ restricted to $D\setminus B_\eta(x)$ under $\Q$ coincides with the law of $e^{\gamma^2K(x,y)}\cM_h(\dd y)$ restricted to $D\setminus B_\eta(x)^c$ under $\hat \sigma \otimes \P$.}

{For this, we first note that $K(x, \cdot)$ is bounded and continuous on $\bar D\setminus B_\eta(x)$. Let $\cM_\eps = \cM_{h, \eps}$ denote the approximation to the Gaussian multiplicative chaos measure in Theorem \ref{T:conv} associated to the field $h$. Since $\cM_{\eps, h} \Rightarrow \cM_h$ in probability as $\eps \to 0$ (where $\Rightarrow$ denotes weak convergence), we deduce immediately that 
\begin{equation}\label{eq:weakconvergenceRN}
e^{\gamma^2 K(x, y)} \cM_{\eps,h}(\dd y) \Rightarrow e^{\gamma^2 K(x, y)} \cM_h(\dd y) \text{ on } D\setminus B_\eta(x)
\end{equation}
in probability (under $\hat  \sigma \otimes \P$). On the other hand, when $\eps>0$ is fixed, it is obvious that 
$$
\cM_{h + \gamma K(x, \cdot) , \eps} (\dd y) = e^{\gamma^2 K_\eps (x, y)} \cM_{\eps, h} ( \dd y),
$$
where $K_\eps(x,\cdot) = K(x, \cdot) * \theta_\eps$ is the $\eps$-regularisation of $K(x, \cdot)$.  By the triangle inequality, and the uniform convergence of $K_\eps(x, \cdot)$ to $K(x, \cdot)$ on $\bar D\setminus B_\eta(x)$, and the fact that $\E( \cM_h(\dd y)) = \sigma(\dd y)$ under $\P$, we deduce from \eqref{eq:weakconvergenceRN} that 
$$
\cM_{h + \gamma K(x, \cdot) , \eps} \Rightarrow e^{\gamma^2 K(x, y)} \cM_h(\dd y)  \text{ on } D \setminus B_\eta(x)
$$ 
weakly in probability under $\hat \sigma \otimes \P$. Using  Corollary \ref{C:Girsanov} (item 4.), the left hand side however has the same law as $\cM_{h, \eps}$ under $\Q$. Taking limits in distribution as $\eps \to 0$, we deduce that the restriction of $\cM_{h}$ to $D \setminus B_\eta(x)$ under $\Q$ has the same law as the restriction of $e^{\gamma^2 K(x, y)} \cM_h(\dd y) $ to $D \setminus B_\eta(x)$ under $\hat \sigma \otimes \P$, as desired.
}
\end{proof}

As will be illustrated below, Girsanov's theorem (either Theorem \ref{T:Girsanov}, Corollary \ref{C:Girsanov} or {Corollary \ref{C:Girsanov2}}) is the basis of many calculations for GMC.

\subsection{Kahane's convexity inequality}
We now present a fundamental tool in the study of Gaussian multiplicative chaos, which is Kahane's convexity inequality. Essentially, this is an inequality that will allow us to ``compare'' the GMC measures associated with two slightly different fields.  Such comparison arguments are very useful in order to do scaling arguments and so compute moments and multifractal spectra, which is our next goal. This inequality was actually crucial to Kahane's construction of Gaussian multiplicative chaos \cite{Kahane85}, although modern approaches such as the one presented just above (coming from \cite{BerestyckiGMC}) do not rely on this.

\ind{Gaussian multiplicative chaos}

\medskip More precisely, the content of Kahane's inequality is to say that given a \textbf{convex} function $f$, and two centred Gaussian fields $X = (X_s)_{s\in T}$ and $Y = (Y_s)_{s\in T}$ with covariances $\Sigma_X$ and $\Sigma_Y$ such that $\Sigma_X(s,t) \le \Sigma_Y(s,t)$ pointwise, we have
$$
\E(f ( \cM_X (D)) \le \E( f(\cM_Y(D)))
$$
for $\cM_X$, $\cM_Y$ the GMC measures associated with $X$ and $Y$. The precise statement of the inequality comes in different flavours depending on what one is willing to assume about $f$ and the fields. A statement first appeared in \cite{Kahane86}, which had an elegant proof but relied on the extra assumption that $f$ is increasing. As we will see this assumption is crucially violated for us (for example, in the proof of \cref{T:moments} we will use $f(x) = -x^q$ with $q<1$, so $f$ is convex but decreasing). The assumption of increasing $f$ is removed in \cite{Kahane85}, whose proof we will follow roughly here.

\ind{Kahane's convexity inequality}
\begin{theorem}[Kahane's convexity inequality] \label{T:Kahane} %Let $(X_i)_i, (Y_i)_i$ be two Gaussian vectors such that:
%\begin{equation}
%\E( X_i X_j) \le \E(Y_i Y_j) \text{ for all } i,j.
%\end{equation}
%Then for all non-negative vectors $(p_i)_i$, and all smooth convex functions $f:(0, \infty) \to \R$, with at most polynomial growth at $0$ and $\infty$, we have
%$$
%\E\left( f\left( \sum_{i=1}^n p_i e^{X_i - \frac12 \E(X_i^2) } \right) \right) \le
%\E\left( f\left( \sum_{i=1}^n p_i e^{Y_i - \frac12 \E(Y_i^2) } \right) \right).
%$$
Suppose that $D\subset \R^d$ is bounded and that $(X(x))_{x\in D}, (Y(x))_{x\in D}$ are almost surely  continuous centred Gaussian fields with $$K_X(x,y):=\mathbb{E}(X(x)X(y))\le \E(Y(x)Y(y))=:K_Y(x,y) \text{ for all } x,y\in D.$$ Assume that $f:(0,\infty)\to \R$ is convex with at most polynomial growth at $0$ and $\infty$, and $\sigma$ is a Radon measure as in \eqref{dim}. Then
$$ \E \left( f \left(\int_{D} \e^{X(x)-\frac{1}{2}\E(X(x)^2)}\sigma(\dd x) \right) \right) \le \E \left( f \left(\int_{D} \e^{Y(x)-\frac{1}{2}\E(Y(x)^2)}\sigma(\dd x) \right) \right).$$
\end{theorem}

%\begin{remark}
%  A word of caution. In various texts (including Kahane's original paper \cite{Kahane86}) the inequality is derived from a version of Slepian's inequality stated for functions $F$ which is only assumed to have second directional derivatives that are positive in the Schwartz sense -- however it is easy to find counterexamples.
%\end{remark}

%{Put Gaussian IBP somewhere?}
%\\

\begin{proof}
The proof is closely related to a Gaussian Integration by Parts formula (see for example \cite{Zeitouni_LN}). \ind{Integration by parts!Gaussian}
Define, for $t \in [0,1]$:
$$
Z_t = \sqrt{1-t} X + \sqrt{t} Y.
$$
Thus $Z_0 = X$ and $Z_1 = Y$. Since the fields $X$ and $Y$ are assumed to be continuous, the maxima and minima of $X$ and $Y$ on $D$ have sub-Gaussian tails by Borell's inequality (see for example \cite[Theorem 2]{Zeitouni_LN}). This means that if $f$ is as in the statement of theorem, we have that 
$$
h(t):=\E\left(f(\int_D Q_t(x) \, \sigma(\dd x))\right):=\E\left(f\left(\int_D \e^{Z_t(x)-\frac{1}{2}\E((Z_t(x))^2)} \, \sigma (\dd x)\right)\right)
$$
is well defined for all $t \in [0,1]$. 

Suppose first that $f$ is smooth. This means we can actually differentiate the above expression and obtain that
$$ \frac{\dd h}{\dd t} \!= \!\frac{1}{2} \E\left( f'\left(\int_D Q_t(x) \, \sigma (\dd x)\right) \int_D\!\sigma (\dd y)  \left(\frac{-X(y)}{\sqrt{1-t}}+\frac{Y(y)}{\sqrt{t}}+K_X(y,y)-K_Y(y,y)\right) Q_t(y)  \right)\!. 
$$
Here we have ``differentiated under the integral sign'' twice (once for the derivative of the integral $\int_D Q_t(x) \sigma(\dd x)$ and once for the expectation) which is permitted since $\int_D Q_t(x) \sigma(\dd x)$ has sub-Gaussian tails and $f$ has at most polynomial growth at $0$ and $\infty$ (since $f'$ is increasing this means that $f'$ also has at most polynomial growth at $0$ and $\infty$).

Consequently, by Fubini's theorem, it suffices to show that for any fixed $y$ :
\begin{equation}\label{eqn:pos_der}
\E\left( \left(\frac{-X(y)}{\sqrt{1-t}}+\frac{Y(y)}{\sqrt{t}}+K_X(y,y)-K_Y(y,y)\right)\,  Q_t(y) \, f'\left(\int_D Q_t(x) \, \sigma (\dd x)\right)    \right)\ge 0.
\end{equation}
Indeed, this  then implies that $h$ is increasing and so $h(0)=\E(f(\int_D \e^{X(x)-(1/2)\E(X(x)^2)}\sigma(\dd x)))$ is less than or equal to $h(1)=\E(f(\int_D \e^{Y(x)-(1/2)\E(Y(x)^2)}\sigma(\dd x)))$, as desired.

To show \eqref{eqn:pos_der}, we fix $y$ and write $$U_t(y):=\frac{-X(y)}{\sqrt{1-t}}+\frac{Y(y)}{\sqrt{t}},$$
so that $U_t(y)$ is the time derivative of the interpolation $Z_t(y)$.
Note that $\E(U_t(y)Z_t(x))=K_Y(x,y)-K_X(x,y)\ge 0$ for all $x$. This means that we can decompose $$Z_t(x)=A_t(x)U_t(y)+V_t(x)$$ for each $x\in D$, where $A_t(x)=(K_Y(x,y)-K_X(x,y))/\E(U_t(y)^2)\ge 0$ and $V_t(x)$ is centred, Gaussian and independent of $U_t(y)$. This corresponds to writing the conditional law of $Z_t(x)$ given $U_t(y)$. Let us rewrite the expectation in \eqref{eqn:pos_der} in terms of $U_t(y)$ and $V_t(y)$. To start with, we decompose \begin{equation}\label{eqn:dec}
Q_t(x) = \e^{A_t(x)U_t(y)-\frac{1}{2}A_t(x)^2\E(U_t(y)^2)} \e^{V_t(x)-\frac{1}{2}\E(V_t(x)^2)}
\end{equation}
for each $x\in D$.
%
%
%(The property that $A_t(x) \ge 0$ says that $Z_t(x)$ is positively correlated with $U_t(y)$. This is to be expected since by bilinearity, this covariance is up to factor of two the time-derivative of $\cov (Z_t(x), Z_t(y))$. The latter is an increasing function of $t$: by assumption, this increases from $K_X(x,y)$ to $K_Y(x,y)$.)
%Note that $A_t(x)$ and the law of $V_t(x)$ do depend on $y$ but since we have fixed $y$ we omit this dependence from the notation.
Thus applying \eqref{eqn:dec} with $x=y$, the expectation in \eqref{eqn:pos_der} can be rewritten as
$$\mathbb{E}\left( \left( U_t(y)-A_t(y)\E(U_t(y)^2)\right) \e^{A_t(y)U_t(y)-\frac{1}{2}A_t(y)^2\E(U_t(y)^2)} \e^{V_t(y)-\frac{1}{2}\E(V_t(y)^2)} f'\left( \int_D \, Q_t(x) \sigma (\dd x)\right)\right).  $$
Now, in order to write this an expectation involving the single Gaussian random variable $U_t(y)$, we consider the conditional expectation (now expanding $Q_t(x)$ as in \eqref{eqn:dec} for clarity):%of this random variable given $U_t(y)$.
%Looking at the decomposition \eqref{eqn:dec}, since $U_t(y)$ is independent of $V_t(x)$ for each $x \in D$ (and thus, by Gaussianity, of $(V_t(x), x \in D)$), and since $A_t(x) \ge 0$, we see that
$$
\E\left(\e^{V_t(y)-\frac{1}{2}\E(V_t(y)^2)} f'\left( \left.\int_D \, \e^{A_t(x)U_t(y)-\frac{1}{2}A_t(x)^2\E(U_t(y)^2)} \e^{V_t(x)-\frac{1}{2}\E(V_t(x)^2)} \sigma (\dd x)\right) \, \right| \, U_t(y)\right).
$$
Since $U_t(y)$ is independent of $V_t(x)$ for each $x \in D$ (and thus, by Gaussianity, of $(V_t(x), x \in D)$), and since $A_t(x) \ge 0$ and $f'$ is increasing,  we see that the above conditional expectation is an almost surely increasing function of $U_t(y)$. Hence \eqref{eqn:pos_der} can be written as
$$ \mathbb{E}\left( g(U_t(y))\left( U_t(y)-A_t(y)\E(U_t(y)^2)\right) \e^{A_t(y)U_t(y)-\frac{1}{2}A_t(y)^2\E(U_t(y)^2)} \right), $$
where $g$ is an increasing function.
 Approximating $g$ by a positive linear combination of step functions and writing $a=A_t(y)$, $\sigma^2=\E(U_t(y)^2)$ it therefore suffices to prove that
$$\int_x^\infty e^{-z^2/2\sigma^2}(z-a\sigma^2)\e^{az-\frac{a^2\sigma^2}{2}} \, \dd z \ge 0$$
for any $x\in \R$.

If $x\ge a\sigma^2$ then the above clearly holds by positivity of the integrand. On the other hand, if $x\le a\sigma^2$ then the integral is greater than
$$\int_{-\infty}^\infty e^{-z^2/2\sigma^2}(z-a\sigma^2)\e^{az-\frac{a^2\sigma^2}{2}} \, \dd z = \frac{\dd }{\dd a} \int_{-\infty}^\infty e^{-z^2/2\sigma^2}\e^{az-\frac{a^2\sigma^2}{2}}=\frac{\dd }{\dd a}(1)=0.$$ This concludes the proof when $f$ is smooth.

In the general case of a convex function $f$, we approximate $f$ by smooth convex functions $f_n \to f$ pointwise with (uniform) polynomial growth at zero and infinity, and then apply dominated convergence, using again the fact that $\sup_x X(x)$ has Gaussian tails; such approximations are easily obtained by approximating the weak derivative of $f$ (a measure) by smooth functions via convolution. 
\end{proof}

\subsection{Scale invariant fields}
\label{SS:scale invariant field}

When we apply Kahane's convexity inequality we will want to compare our Gaussian field with an auxiliary Gaussian field enjoying an exact scaling relation. In this section we explain a modification, due to Rhodes and Vargas (\cite{RhodesVargasMRM}) of a construction due to Bacry and Muzy (\cite{BacryMuzy}), that will give us the desired scale invariant field. (In the case of the two dimensional GFF the Markov property gives a close analogue but would lead to extra technicalities.) The main result of this section is Theorem \ref{C:scaling2}, which does not seem to appear in the literature in this form.

\subsubsection{One dimensional cone construction}
\begin{figure}
\begin{center}
  \includegraphics{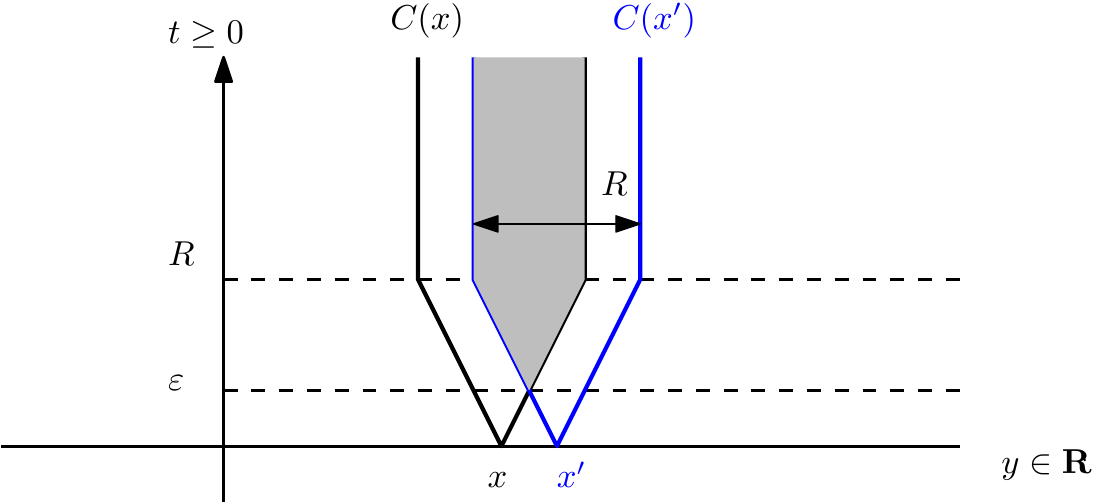}
  \end{center}
  \caption{The truncated cones in the construction of the scale invariant auxiliary field. The covariance of the field at $(x,x')$ is obtained by integrating $\dd y \dd t / t^{2}$ in the shaded area.}
\end{figure}
We first explain the construction we will use in one dimension where things are easier.
 Fix $0 < \eps <R$ and for $x \in \R$, consider the \textbf{truncated cone} $C_{\eps, R}(x)$ in $\R^{2}$ given by
\begin{equation}\label{D:conegmc}
C(x) = C_{R}(x) = \{ z=(y,t) \in \R \times [0, \infty):  |y - x | \le (t \wedge R)/2\}.
\end{equation} 
where $| y - x | $ denotes Euclidean norm in $\R$. Define a kernel
$$
c_{\eps, R} (x, x') = \int_{ y \in \R} \int_{t \in [ \eps, \infty)}  \mathbf{1}_{\{(y,t) \in C(x) \cap C(x')\}} \frac{\dd y \dd t}{t^{2}}.
$$ 
Note that since the domain of integration has been truncated at $t=\eps$, the integral is finite.
We claim that $c_{\eps, R}$ is non-negative definite and so can be used to defined a Gaussian field $X_{\eps, R}$ on $\R$ whose covariance is given by $c_{\eps, R}$. Indeed, for any $n \ge 1$, for any $x_1, \ldots, x_n \in \R^d$ and $\lambda_1, \ldots, \lambda_n \in \R$,
\begin{align}\label{posdef1}
  \sum_{i,j=1}^n \lambda_i \lambda_j c_{\eps, R} (x_i, x_j) & = \sum_{i,j=1}^n \lambda_i \lambda_j \int_{\R} \int_\eps^\infty \mathbf{1}_{\{(y,t) \in C(x_i)\}} \mathbf{1}_{\{(y,t) \in C(x_j)\}} \frac{\dd y \dd t}{t^{2}} \nonumber \\
  & = \int_{\R^d} \int_\eps^\infty  \left( \sum_{i=1}^n \lambda_i \mathbf{1}_{\{(y,t) \in C(x_i)\}}\right)^2 \frac{\dd y \dd t}{ t^{2}} \ge 0.
\end{align}
As the covariance kernel $c_{\eps,R}$ is a nice continuous function of $x, x'$ we can check (again using for example Proposition 2.1.12 in \cite{GineNickl}) that there exists a centred Gaussian field $X_{\eps,R}$ whose covariance is given by $c_{\eps , R}$ and which is almost surely Borel measurable as a function on $\R$.

\begin{rmk}\label{R:posdef}
This computation showing that $c_{\eps, R}$ is non-negative definite works because the covariance is defined to be of the form $c(x, x') = \int_S f_x(z) f_{x'} (z) \nu(\dd z)$, for some fixed function $f_x$ of $z$ associated to each $x \in \R$, where the integral can be on some arbitrary space $S$ with measure $\nu$. Here the space $S$ is $\R \times (0, \infty)$, $\nu (dz) = \mathbf{1}_{\{ t \ge \eps\}}dy dt / t^2$, and $f_x(z) = \mathbf{1}_{z \in C_{\eps, R}(x)}$. We will use other choices when considering the higher dimensional case.
\end{rmk}

The key property of $X_{\eps,R}$ (and the reason for introducing it) is that its covariance can be computed exactly.
This not only shows that the field is logarithmically correlated, but enjoys an exact scaling relation, as follows.

\begin{lemma} \label{L:covariance1}
Define the function
\begin{equation}\label{D:g}
g_{\eps, R} (x) =
\begin{cases}
  \log_+ (R/|x|) &  \text{ if } |x | \ge \eps \\
  \log ( R/ \eps) + 1 - (|x |/\eps) &   \text{ if } |x | \le \eps,
\end{cases}
\end{equation}
where $\log_+(x)=\log(x)\vee 0$.
Then for all $x, y \in \R$,
%For all $x, y \in \D$ we have
\begin{equation}\label{conelog1}
c_{\eps,R}(x,x') = g_{\eps, R} (x - x').
\end{equation}
In particular, $c_{\eps, R} (x, x') = \log (R / (|x-x'|\vee \eps)) + O(1)$, where the $O(1)$ term does not depend on $x, x'$, $\eps$ or $R$ (and is in fact bounded between $0$ and 1).
\end{lemma}

\begin{proof}
By translation invariance and symmetry we can assume that $x'= 0$ and $x >0$. If $x \ge R$ there is nothing to prove, so assume first that $\eps \le x \le R$. Then the two cones first intersect at height $x \ge \eps$. Moreover the width of the intersection of these cones at height $t \ge x$ is $(t-x) \wedge (R-x)$, so
\begin{align*}
c_{\eps, R} (0,x) &= \int_x^R (t-x) \frac{\dd t}{t^2} + \int_R^\infty (R- x) \frac{\dd t}{t^2}\\
& = \log (R/x) - x(\frac1{x} + \frac{1}{R}) + \frac{(R-x)}{R} \\
& = \log (R/x)
\end{align*}
as desired. When $x \le \eps$, the computation is almost the same, but the lower bound of integration for the first integral is $\eps$ instead of $x$, which gives the desired result.
\end{proof}

We now explain why this implies a scaling property. We fix the value $R$ of truncation and write $X_\eps $ for $X_{\eps, R}$ (often we will choose $R=1$ and write $Y_\eps$ for $X_{\eps, 1}$).
%For this we fix the value $R$ (corresponding to the large scale truncation), and let $Y_\eps = X_{\eps, 1}$.

\begin{cor}
 \label{C:scaling1}
 For $\lambda<1$,
$$
(X_{\lambda \eps}(\lambda x))_{x \in B(0,R/2)} =_d (\Omega_\lambda + X_\eps(x))_{x \in B(0,R/2)},
$$
where $\Omega_\lambda$ is an independent centred Gaussian random variable with variance $\log(1/\lambda)$.
 \end{cor}

\begin{proof}
One directly checks that for all $x, x' \in \R$ such that $| x - x'| \le R$ (and so also $| x- x'| \le R/ \lambda$ automatically),
$$
c_{\lambda\eps, R}(\lambda x, \lambda y) =  c_{\eps,R}(x,y) + \log (1/ \lambda)
$$
and hence the result follows.
%(this can be seen either by looking at the exact value of the function $g_{\eps , R}$ computed in the previous lemma or by noting that $dy dt / t^2$ is scale invariant in $\R^2$. Hence, since the dependence in $R$ only appears in the logarithmic term,
%$$
%c_{\lambda\eps, R} (\lambda x, \lambda y)  = \log (1/\lambda) + c_{\eps,R}(x,y)
%$$
%{Should $R=1$ in this proof? Also don't we need $x,x'\in B(0,1/2)$ to have $|x-x'|\le 1$? If $|x-x'|> 1$ but $|\lambda x - \lambda x'|< 1$ this is a problem right?  }
\end{proof}

\subsubsection{Higher dimensional construction}\label{S:hdc}

The one dimensional Bacry--Muzy construction presented above is beautiful and simple but does not trivially generalise to more than one dimension. This is because if one considers truncated cones in $\R^{d+1}$ (instead of $\R^2$) and integrates with respect to the scale invariant measure $\dd y \dd t /t^{d+1}$, the volume of the intersection of two truncated cones based at $x$ and $x'$ does not lead to nice formulae which yield scale invariance in the sense of \cref{C:scaling1} (see \cite{Chainais} for an article where this model is nevertheless studied).

To overcome this problem we follow (in a slightly simplified setting) a very nice construction proposed by Rhodes and Vargas \cite{RhodesVargasMRM} in which the exact one dimensional computation of Bacry and Muzy can be exploited to give a field in any number of dimensions satisfying both logarithmic correlations and exact scaling relations. The basic idea is to define the cones on $\R^d$ based at $x$ and $x' \in \R^d$ by first applying a random rotation in order to preserve isotropy, and then applying the one dimensional construction to the first coordinates of $x$ and $x'$.

To be more precise, let $d \ge 1$ and consider $\mathcal{R}$ the orthogonal group of $\R^d$: that is, $d$ dimensional matrices $M$ such that $MM^t = I$. %(In two dimensions, $\mathcal{R}$ is nothing but the set of rotations). 
Let $\Sigma$ denote Haar measure on $\mathcal{R}$ normalised to be a probability distribution. 
%(Thus in two dimensions $\Sigma$ can be identified with the uniform distribution on the unit circle).
\begin{comment}
We consider a Gaussian random measure $\cM$ on $\mathcal{R} \times \R \times (0, \infty)$ with intensity $d \sigma \otimes dy  dt / t^2 $. That is, $\cM$ satisfies:

\begin{enumerate}
  \item For any disjoint Borel sets $A_1, A_2, \ldots$, the random variables $\cM(A_n)$ are independent and
      $$
      \cM( \bigcup_{n=1}^\infty A_n) = \sum_{n=1}^\infty  \cM (A_n), almost surely
      $$
  \item if $A$ is a Borel set then $\cM(A)$ is a Gaussian random variable with mean zero and variance $\sigma^2$ where
      $$
      \sigma^2 =
      \int_{\mathcal R \times \R \times (0, \infty)}\mathbf{1}_A  d\sigma \otimes \frac{dy dt }{ t^2}.$$
\end{enumerate}

This is analogous to the more familiar notion of a Poisson random measure with given intensity but there is no problem in extending this notion to Gaussian random measure (and indeed to random measures based on any given infinitely divisible distribution, as is done in \cite{BacryMuzy} and \cite{RhodesVargasMRM}).
\end{comment}
If $\rho \in \cR, x \in \R^d $, let $\rho x$ denote the vector of $\R^d$ obtained by applying the isometry $\rho$ to $x$, and let $(\rho x)_1$ denote its first coordinate. Define the \textbf{cone like} set $\mathbf{C}_{R} (x)$ as follows:
$$
\mathbf{C}_{R}(x) : = \{( \rho, t, y) \in \mathcal{R} \times \R \times (0, \infty): (t,y) \in C_{R} ((\rho x)_1) \}.
$$
where if $z \in \R$, $C_R(z)$ is the truncated cone of \eqref{D:conegmc}. Thus for any given $\rho$, we first apply $\rho$ to $x$ and consider the truncated cone (in two dimensions) based on the first coordinate of $\rho x$.
As in the previous section we define a field through its covariance kernel
$$
\mathbf{c}_{\eps, R} (x, x') = \int_{\mathcal{R} \times \R \times (0, \infty)} \mathbf{1}_{\{( \rho, t, y) \in \mathbf{C}_R(x) \cap \mathbf{C}_R (x') \} } \mathbf{1}_{\{t \ge \eps\}} \dd  \Sigma (\rho) \otimes \frac{\dd y \dd t}{t^2}.
$$
We note that this is non-negative definite for the same reasons as \eqref{posdef1} (see especially \cref{R:posdef}). Hence, as before we can consider an almost surely Borel measurable function $x \in \R^d \mapsto X_{\eps, R} (x) \in \R$ which is a translation-invariant, centered Gaussian field with $\mathbf{c}_{\eps, R}$ as its covariance kernel.

\begin{theorem}
 \label{C:scaling2} Fix any $R>0$. The restriction of the field $(X_{\eps, R})_{x\in \R^d}$ to $x \in B(0, R/2)$, viewed as a function of $\eps>0$ and $x\in B(0, R/2)$, is scale-invariant in the following sense:
 for any $\lambda<1$,
\begin{equation}\label{scaleinv}
(X_{\lambda \eps,R}(\lambda x))_{x \in B(0,R/2)} =_d (\Omega_\lambda + X_{\eps,R}(x))_{x \in B(0,R/2)},
\end{equation}
where $\Omega_\lambda$ is an independent centred Gaussian random variable with variance $\log(1/\lambda)$. Furthermore, its covariance function $\mathbf{c}_{\eps, R}$ satisfies
\begin{equation}\label{scale_estimate}
\mathbf{c}_{\eps, R} (x, x' ) = \log (\frac1{\|x-x'\| \vee \eps}) + O(1),
\end{equation}
uniformly over $x, x' \in B(0, R/2)$, where the implicit constant $O(1)$ above depends only on the dimension $d \ge 1$.
 \end{theorem}

\begin{proof} We start by noticing that we have the following exact expression for the covariance. Recall the function $g_{\eps, R}(t)$ for $t \in \R$ from \cref{L:covariance1}:
$$
g_{\eps, R} (t) =
\begin{cases}
  \log_+ (R/|t|) & \text{ if } |t | \ge \eps \\
  \log ( R/ \eps) + 1 - (|t |/\eps) & \text{ if } |t | \le \eps.
\end{cases}
$$
Using Fubini's theorem, and since $g_{\eps, R}$ gives the covariance in the one dimensional case (\cref{L:covariance1}), we have:
\begin{equation}  \label{L:covariance2}
 \mathbf{c}_{\eps, R} (x, x') = \int_{\rho \in \mathcal{R}} g_{\eps, R} ( (\rho x)_1 - (\rho x')_1 ) \dd \Sigma(\rho).
\end{equation}
The scale invariance in \eqref{scaleinv} then follows easily.
%Fix a bounded continuous function $\phi: \R^d \to [0, \infty)$  which is a bounded \emph{positive definite function} and such that $\phi(x) = 0$ if $|x|\ge 1$. For instance, we choose $\phi(x) = {(1- |x|^\cM)_+}$ for an appropriate value of $\cM$ depending on the dimension: in dimension $d=1$, $\cM = 1$ works (as can be checked readily by computing the Fourier transform) and $\cM = 1/2$ works in dimension $d=2$. This is the so-called Kuttner--Golubov problem. \red{CHANGE}
%See the discussion in Example 2.3 in \cite{RobertVargas} to check that this is indeed a positive definite function. Define an auxiliary centred Gaussian random field $(X_\eps(x))_{x\in \R^d}$  by specifying its covariance
%$$
%c_\eps(x,y): = \E(X_\eps(x) X_\eps(y))   = \log_+\left(\frac1{|x- y | \vee \eps} \right)+ \phi\left( \frac{|y-x|}{\eps}\right).
%$$
%As before, we fix the value $R$ (corresponding to the large scale truncation), and write $X_\eps $ for $ X_{\eps, R}$. Again, typically we will choose $R=1$ in which case we write $Y_\eps$ for $X_{\eps, 1}$.
Indeed, if $x, x' \in \R^d$ are such that $| x - x'| \le R$ (and so also $| x- x'| \le R/ \lambda$ automatically), then note that
$$
g_{\lambda\eps, R}(\lambda x - \lambda y) =  g_{\eps,R}(x - y) + \log (1/ \lambda)
$$
which, as already noticed in \cref{C:scaling1}, immediately implies \eqref{scaleinv}.

Let us now turn to the proof of \eqref{scale_estimate}.
Since $g_{\eps, R}(t) = \log_+ (R/ (|t|\vee \eps)) + O(1)$, we have
\begin{equation}\label{scale_est_int}
\mathbf{c}_{\eps, R}(x, x') = \int_{\cR} \log ( \frac{R}{|\rho (x- x')_1| \vee \eps }) \dd  \Sigma(\rho) + O(1).
\end{equation}
Since $| \rho(x-x')_1| \le | \rho (x- x') |= \|x - x'\|$ for any $\rho \in \cR$, we immediately get the lower bound
\begin{equation}\label{scale_estimate_LB}
\mathbf{c}_{\eps, R} (x,x') \ge \log ( \frac{R}{ \|x-x'\| \vee \eps}) + O(1).
\end{equation}
To get a bound in the other direction, we observe that for a fixed vector $u \in B(0, R/2)$, the distribution of $\rho u$ under the Haar measure $\dd \Sigma(\rho)$ is uniform on the sphere of radius $\|u\|$. Its first coordinate $(\rho u)_1$ therefore has an absolutely continuous distribution with respect to $(1/\|u\|)$ times (one dimensional) Lebesgue measure. As a consequence, if
$$
\cR_k(u) =  \{ \rho \in \cR: | (\rho u)_1| \in  [2^{-(k+1)}\|u\|, 2^{ - k } \|u\|] \}
$$
then
\begin{equation}\label{HaarLeb}
  \Sigma ( \cR_{k} (u) ) \le O(  2^{ - k}) ,
\end{equation}
where the implicit constant depends only on the dimension $d \ge 1$. We note that the right hand side does not depend on $u$ since the quantity on the left hand side is clearly scale (and rotation) invariant.  Therefore, from \eqref{scale_est_int} with $x-x'=u$, and since $\eps \ge 2^{-k-1} \eps$, 
\begin{align*}
  \mathbf{c}_{\eps, R}(x, x') & = O(1) + \sum_{k \ge 0} \int_{\cR_k(u)} \log ( \frac{R}{ |(\rho u)_1 | \vee \eps}) \dd  \Sigma(\rho) \\
  & \le O(1) + \sum_{k \ge 0} \int_{\cR_k(u)} \log ( \frac{R}{\|2^{-k-1} u\| \vee \eps}) \dd \Sigma(\rho) \\
  & \le O(1) + \sum_{k \ge 0} \int_{\cR_k(u)} \log ( \frac{R}{\|u\| \vee \eps}) \dd  \Sigma(\rho) + O(k) \Sigma (\cR_k(u))\\
  & \le O(1) + \log ( \frac{R}{\|u\| \vee \eps})  + \sum_{k\ge 0} O( k2^{-k})
\end{align*}
where we have used \eqref{HaarLeb} in the last line. This proves \eqref{scale_estimate}.
%{Same comment as for the 1d analogue.}
\end{proof}

\begin{rmk}\label{rmk:logPD}The covariance kernel takes a particularly nice form in a fixed neighbourhood of a given point when $\eps \to 0$. Indeed, note that if $x \in B(0, R)$ and $|(\rho x)_1 |\ge  \eps$, then writing $x = \|x\| e_x$ where $e_x$ is the unit vector in the direction of $x$, we have (letting $e_1$ denote the unit vector in the first coordinate),
$$ g_{\eps,R} ((\rho x)_1 ) = \log( R/ \langle \rho x, e_1 \rangle) =  \log (R / \|x\|) + \log( R/ \langle \rho e_x, e_1 \rangle) .$$
When we integrate against $\dd \Sigma$, we can take advantage of rotational symmetry to note that
$$
C = \int_{\rho \in \mathcal{R}}\log( R/ \langle \rho e_x, e_1 \rangle) \dd  \Sigma (\rho)
$$
does not in fact depend on $x$.

Therefore for \emph{any} $ x \in B(0, R)$,
$$
\lim_{\eps \to 0} \mathbf{c}_{\eps, R} (x, 0) = \log (R/\| x \|) + C.
$$
It follows from this observation that in $B(0, R)$ that the function $x \mapsto \log (R/\| x\|) + C$ is positive definite in $B(0, R)$. We can get rid of the constant $C$ by changing the value of $R$, and so we deduce that
$$
x \mapsto K(x) := \log (R/ \|x\|) \text{ is positive definite in a small neighbourhood of 0},
$$
a fact which appears to have been first proved for all dimensions in \cite{RhodesVargasMRM}. Note that the size of this neighbourhood does depend on the dimension $d$.
\end{rmk}

\begin{rmk}
	It was shown in \cite{Decomplogfields} that if the continuous term $g$ from the decomposition \eqref{cov} of $K$ is an element of $H_{\mathrm{loc}}^{d+\eps}(D\times D)$ for some $\eps>0$, then $K$ is locally positive definite on $D$. This provides an alternative justification that $K(x)=\log(R/\|x\|)$ is locally positive definite on $\R^d$.
\end{rmk}

\begin{rmk}
By contrast, note that $ x \mapsto \tilde K(x) = \log_+ (R/ \|x\|) $ is positive definite in the \emph{whole space} if and only if $d \le 3$:  see Section 5.2 of \cite{RobertVargas} for a nice proof based on Fourier transform.

  When $d =1$ or $d=2$ one can also show that $\tilde K(x)$ is not only positive definite but of $\sigma$-positive type in the sense of Kahane: that is, it is a sum $\tilde K(x) = \sum_{n=1}^\infty K_n(x)$ where the summands $K_n$ are not only positive definite functions, but also pointwise non-negative ($K_n(x) \ge 0$). When $d=3$ it is an open question to determine whether $\tilde K(x)$ is $\sigma$-positive.
\end{rmk}

\subsection{Multifractal spectrum}
\label{S:multifractal}

We now explain how Kahane's convexity inequality can be used to obtain various estimates on the moments of the mass of small balls, and in turn to the multifractal spectrum of Gaussian multiplicative chaos. We take $h,\theta$ as in \cref{setup_gmc}, and we assume that $d=\mathbf{d}$ and the reference measure $\sigma$ is Lebesgue measure for simplicity.

\begin{theorem}[Scaling relation for Gaussian multiplicative chaos]\label{T:moments} Let $\gamma \in (0,\sqrt{2d}).$ Let $B(r)$ be a ball of radius $r$ in the domain $D$. Then uniformly over all such balls, and for any $q \in \R$ (including $q<0$) such that $ \cM_\eps(B(0,1))^q $ is uniformly integrable in $\eps$,
\begin{equation}\label{MF}
\E( \cM( B(r))^q) \asymp r^{( d+ \gamma^2/2)q  - \gamma^2 q^2/2},
\end{equation}
where $a_r\sim b_r$ if $C^{-1}a_r\le b_r \le C a_r$ for some constant $C$ depending only on $\sup_{\bar D} |g|$, $q$, and $\gamma$.
The function
\begin{equation}\label{E:xi}
\xi(q) = q(d + \gamma^2/2) - \gamma^2 q^2 / 2
\end{equation}
is called the \textbf{multifractal spectrum} of Gaussian multiplicative chaos. \indN{{\bf Gaussian multiplicative chaos}! $\xi(\cdot)$; multifractal spectrum function of Gaussian multiplicative chaos}
\end{theorem}

\ind{Liouville measure!Scaling relation}

\ind{Liouville measure!Multifractal spectrum}

\begin{rmk}
In the next section, we will see that the assumption on $q$ is equivalent to
$$
q< \frac{2d}{\gamma^2}.
$$
At this stage we already know it at least for $0 \le q <1$.
\end{rmk}

\begin{rmk} \textbf{What is a multifractal spectrum?} The above theorem characterises the multifractal spectrum of Gaussian multiplicative chaos. To explain the terminology, it is useful to consider the opposite case of a \emph{monofractal} object. For instance, Brownian motion is a monofractal because its behaviour is (to first order at least) described by a single fractal exponent, $\alpha = 1/2$. One way to say this is to observe that for all $q$
$$\E(|B_t|^q) \asymp t^{q/2}.$$
(A variety of exponents can however be obtained by considering logarithmic corrections, see for example \cite{MortersPeres}). The monofractality of Brownian motion is thus expressed through the fact that its moments have a power law behaviour where the exponent is \emph{linear} in the order of the moment $q$. By contrast, note that the function $\xi$ in \cref{T:moments} is \textbf{non-linear}, which is indicative of multifractal behaviour. That is, several fractal exponents (in fact, a whole spectrum of exponents) are needed to characterise the first order behaviour of Gaussian multiplicative chaos. Roughly speaking, the multifractal formalism developed among others in \cite{Falconer} is what allows the data of a non-linear function such as the right hand side of \eqref{MF} to be translated into a knowledge about the various fractal exponents and their relative importance.
% {What does "first order" mean in the above?}
\end{rmk}

\begin{proof}[Proof of \cref{T:moments}]
%Recall that we seek to establish: for $q \in [0,1]$,
%\begin{equation}\label{MF}
%\E(\cM(B(r))^q) \asymp r^{(2 + \gamma^2/2)q - q^2 \gamma^2/2}.
%\end{equation}
Set $R=1$ and let $Y_\eps = X_{\eps, 1}$ denote the scale invariant field constructed in Theorem \ref{C:scaling2}. As hinted previously, the idea will be to compare $h_\eps$ to the scale invariant field $Y_\eps$. Note that by the estimate \eqref{scale_estimate} in Theorem \ref{C:scaling2} on the one hand, and \eqref{eq:roughcov} on the other hand,
there exist constants $a, b>0 $ such that
\begin{equation}\label{encadr-cov}
\mathbf{c}_{\eps}(x,y) - a \le \E(  h_\eps(x)  h_\eps(y)) \le \mathbf{c}_\eps(x,y) + b
\end{equation}
where $\mathbf{c}_\eps = \mathbf{c}_{\eps,1}$ is the covariance function for $Y$.
As a result it will be possible to estimate the moments of $\cM(B(r))$ up to constants by computing those of $\tilde \cM(B(0,r))$, where $\tilde \cM$ is the chaos measure associated to $Y$. %{(Do we actually want to do some translation and use $\tilde \cM(B(0,r))$ for the comparison? In the below there is lots of switching between $B(r)$ and $B(0,r)$. I guess we just need to be consistent).} %obtained from $X$ in a similar way as $\cM$ is obtained from $h$:
More precisely, from \eqref{encadr-cov} and Kahane's convexity inequality (applied to the fields $h_\eps$ and $Y_\eps + \cN(0, a)$ in one direction and to the fields $Y_\eps$ and $h_\eps + \cN(0, b)$ in the other direction, with the function $f$ taken to be the concave or convex function $x \mapsto x^q$), we get:
\begin{equation}\label{encadrKahane}
\E ( (\cM_{\eps}(S))^q ) \asymp \E( (\tilde \cM_\eps(S))^q )
\end{equation}
for $S\subset D$, where the implicit constants depend only on $a, b$ and $q \in \R $, and not on $S$ or $\eps$, and where
$$
\tilde \cM_\eps(z) = \exp (  \gamma Y_\eps(z) - (\gamma^2/2) \E( Y_\eps(z)^2) ) \dd z.
$$
It therefore suffices (also making use of the translation invariance of $Y$) to study the moments of $\tilde \cM_\eps (B(0,r))$.

We now turn to the proof of \eqref{MF}.
%Define the normalised field to be
%$\tilde  X_\eps(x) = \gamma X_\eps(x) - (\gamma^2/2) \E( X_\eps(x)^2) = \gamma X_\eps(x) + (\gamma^2/2) \log \eps + C $, for some constant $C$, and
Note that $\E( Y_\eps(x)^2) = \log (1/\eps) + O(1)$.
Fix $\eps>0$, and $\lambda = r<1$.
Then
\begin{align*}  \tilde \cM_{r\eps}(B(0,r)) &\asymp   \int_{B(0,r)}   e^{\gamma Y_{r \eps}(z) } (r\eps)^{\gamma^2/2} \dd z\\
 & =   r^{d + \gamma^2/2} \int_{B(0,1)} e^{\gamma Y_{r\eps} (rw)} \eps^{\gamma^2/2} \dd w
\end{align*}
by the change of variables $z = rw$. Hence by  \cref{C:scaling2},
\begin{equation}\label{prescaling}
\tilde  \cM_{r\eps}(B(0,r)) \asymp  r^{d + \gamma^2/2} e^{\gamma \Omega_r}  \tilde \cM'_\eps(B(0,1))
\end{equation}
where $\tilde \cM'$ is a copy of $\tilde \cM$ and $\Omega_r$ is an independent $\cN(0, \log (1/r))$ random variable.  Raising to the power $q$, taking expectations and using \cref{encadrKahane}, we get:
\begin{align}
\E( \cM_{r\eps}(B(r)^q) ) & \asymp \E( \tilde \cM_{r\eps}(B(0,r)^q) )\nonumber \\
& = r^{q (d + \gamma^2/2)} \E( e^{\gamma q \Omega_r}) \E( \tilde \cM_\eps(B(0,1))^q) \nonumber \\
& \asymp r^{ \xi(q)} \E( \cM_\eps(B(0,1)^q))\label{scalingeps}
\end{align}
where
$$
\xi(q) = q(d + \gamma^2/2) - \gamma^2 q^2 / 2
$$
is the {multifractal spectrum} from the theorem statement.
Suppose now that $q$ is such that $\cM_\eps(B(0,1))^q$ is uniformly integrable. Then
$$
\E( \cM(B(r))^q) \asymp  r^{\xi(q)},
$$
as desired.
\end{proof}

\subsection{Positive moments of Gaussian multiplicative chaos (Lebesgue case)}
\label{S:positivemomentsGMC}
We continue our study of GMC initiated above in $\R^d$ with the reference measure $\sigma$ taken to be the Lebesgue measure, for a logarithmically correlated field $h$ satisfying the general assumptions of  \cref{SS:setup}. Let $\cM$ be the associated GMC measure.
The goal of this section will be to prove the following theorem on its moments. (See Section \ref{S:posgen} for similar results where $\sigma$ is allowed to be more general than Lebesgue measure).
%First we will need to make an assumption on the underlying reference measure $\sigma$ that its support is of Minkowski dimension $\le d$. We first recall the definition of Minkowski (or box counting) dimension:

%\ind{Minkowski dimension} Let $\cS_m$ denote the $m$th level dyadic covering of the domain $\R^k$ by cubes $S_i, i \in \cS_m$ of side-length $2^{-m}$. The $2^{-m}$-box counting number of $S$ is, by definition:
%$$
%N(S; 2^{-m}) = \sum_{i \in \cS_m} \indic{S_i \cap S \neq \emptyset}
%$$
%The (Euclidean) Minkowski dimension of $A$ is then defined by
%\begin{equation}\label{Minkowskidef}
%D_M (A) = \inf \{ \delta: \limsup_{m \to \infty} 2^{- m \delta} N(A, 2^{-m})  < \infty\}.
%\end{equation}

\begin{theorem}
  \label{T:finiteposmoments} Let $S \subset D$ be bounded and open, and suppose that $\sigma(\dd x) = \dd x$ is the Lebesgue measure on $\R^d$. Let $\gamma\in (0,2)$ and $q >0$. Then $\E(\cM(S)^q) < \infty$ if
  \begin{equation}\label{E:finiteposmoments}
  q < \frac{2d}{\gamma^2}.
  \end{equation}
\end{theorem}

In fact, the theorem shows that $ (\cM_\eps(S))^q$ is uniformly integrable in $\eps$, so that \cref{T:moments} applies to this range of values of $q$.

\medskip Before starting the proof of this theorem, we note that %$q$ is allowed to be negative: hence $\cM(S)$ has finite moments of all negative orders, which is to say that GMC has a very light tail at zero. By contrast, the tail at $\infty$ has a power law:
from \cref{T:finiteposmoments},
\begin{equation}\label{E:tail}
\P( \cM(S) > t) \le t^{ - 2d/\gamma^2 + o(1)}; \ \ t \to \infty.
\end{equation}
In fact, much more precise information is known about the tail at $\infty$: a lower bound matching this upper bound can be obtained so that it becomes an equality. In fact, the $o(1)$ term in the exponent can also be removed and a constant identified: in the case of the two dimensional GFF this was done by Rhodes and Vargas \cite{RhodesVargas_tail}, and the universality of this behaviour (including the calculation of the constant itself) was shown subsequently in a paper by Mo-Dick Wong \cite{Wong_tail}.

%We give two proofs of this important result. The first one (which is ``classical") is when the reference measure $\sigma$ has a nice density with respect to Lebesgue measure, (in particular, the dimension of the measure $d$ coincides with the dimension of the ambient space $\R^k$ and so $k=d$). The second one, more general, will be given below and deals with the general setup announced in the theorem.

\begin{proof}Note that we already know uniform integrability of
$ \cM_\eps(S)$ (\cref{T:conv}) so we can assume $q>1$. For simplicity (and without loss of generality) we assume that $S$ is the unit cube in $\R^d$. Let $\cS_m$ denote the $m$th level dyadic covering of the domain $\R^d$ by cubes $S_i, i \in \cS_m$ of sidelength $2^{-m}$.   Given $q<2d/\gamma^2$, we define $n=n(q) \ge 2$ such that $n-1 < q \le n$. We will show by \textbf{induction} on $n$ that
$$
M_\eps := \E(\cM_\eps(S)^q)
$$
is uniformly bounded.
%{Is the induction argument actually just for the off-diagonal terms?}

  Let us consider the case $n=2$ first. We first subdivide the cubes of $\cS_m$ into $2^d$ disjoint groups so that no two cubes within any given group touch (including at the boundary); thus any two cubes within a given group are at distance at least $2^{-m}$ from one another. The reader should convince themselves that this is actually possible (it is a generalisation of the usual checkerboard pattern for $\Z^2$). Let $\cS'_m$ denote one of these $2^d$ groups of cubes of sidelength $2^{-m}$.
  %By the triangle inequality in $L^q$, it will suffice to check that
 % $$
 % \E[\Big(\sum_{i \in \cS'_m} \cM(S \cap S_i)\Big)^q] <\infty
 %  $$
 %  if $q < 2d/\gamma^2$.

We will now take advantage of some convexity properties, using the fact that $q/2 \le 1$ (recall that $n=2$ and $n-1<q\le n$ by definition). We write, for given $m$,
  \begin{align}
    \Big(\sum_{i \in \cS'_m} \cM_\eps(S \cap S_i)\Big)^q & = \left( \sum_{i,j \in \cS'_m} \cM_\eps(S_{i}) \cM_\eps(S_{j})\right)^{q/2} \nonumber \\
    & \le \sum_{i,j\in \cS'_m} \cM_\eps(S_i)^{q/2} \cM_\eps(S_j)^{q/2}, \label{squaredec}
  \end{align}
  where we have used the elementary fact that $(x+y)^{\alpha} \le x^\alpha + y^\alpha$ if $x,y >0$ and $\alpha \in (0,1)$. (This is easily proven by writing $(x+y)^{\alpha} - x^\alpha=\int_x^{x+y} \alpha t^{\alpha -1}\dd t\le \int_0^y \alpha t^{\alpha-1} \, \dd t = y^\alpha$, since the integrand $\alpha t^{\alpha-1}$ is decreasing in $t$).

  We consider the on-diagonal and off-diagonal terms in \eqref{squaredec} separately.
%  $$
%\E [\Big(\sum_{i \in \cS'_m} \cM_\eps(S \cap S_i)\Big)^q]
%  $$
  We start with the on-diagonal terms (the estimate in this case works for general $q>0$ so is not restricted to the case $n=2$):

  \begin{lemma}\label{L:ondiagbase}
  Assume the set up of  \cref{T:finiteposmoments}. Then there exist a constant $c_q$ such that for all sufficiently large $m$, and for all $\eps>0$,
    \begin{equation}
    \E (\Big(\sum_{i \in \cS'_m} \cM_\eps(S_i)^q\Big))
 \le c_q 2^{dm - \xi(q) m} \E ( \cM_{\eps2^m} (S)^q). \label{E:ondiagbase}
    \end{equation}
  \end{lemma}

  \begin{proof}
By \eqref{scalingeps}, applied with $r = 2^{-m}$, we have for each $i \in \cS'_m$,
$$
\E( \cM_{\eps} (S_i)^q )\le c_q 2^{- \xi(q) m}\E( (\cM_{\eps2^m} (S))^q ) .
$$
Since there are at most $2^{dm}$ terms in this sum, we deduce the lemma.
  \end{proof}

  For the off-diagonal terms, we simply observe that in the case where the two indices are distinct:

  \begin{lemma}\label{L:offdiagbase}   Assume the set up of  \cref{T:finiteposmoments}. Then for any fixed $m$ and $q <2$, there exists a constant $C_{m,q}$ independent of $\eps$ such that
  $$
 \E \big( \sum_{i\neq j \in \cS'_m } \cM_\eps(S_{i})^{q/2} \cM_\eps(S_{j})^{q/2}\big) \le C_{m,q}.
  $$
  \end{lemma}

  \begin{proof}Note that by Jensen's inequality (since $q/2 \le 1$),
  $$
\E\big(  \cM_\eps(S_i)^{q/2} \cM_\eps(S_j)^{q/2}\big) \le \E \big( \cM_\eps(S_i) \cM_\eps(S_j) \big)^{q/2}
    $$
    for all $i \neq j\in \cS'_m$.
    The expectation can easily be computed, and we have for some constant $c_m$,
    $$
    \lim_{\eps \to 0} \E\big(  \cM_\eps(S_i)\cM_\eps(S_j)\big)  \le \int_{x \in S_i, y \in S_j} e^{\gamma^2 K(x,y)} \dd x \dd y \le c_{m}<\infty
    $$
    since the squares $S_i$ and $S_j$ are at distance at least $2^{-m}$ from one another. Taking the $q$th power and summing over all terms $i\neq j \in \cS'_m$ gives the lemma.
  \end{proof}

  We put these two lemmas together as follows. First, note that for $q< 2d / \gamma^2$, $2d - \xi(q) <0$. We can therefore choose $m$ large enough that $c_q 2^{dm - \xi(q) m} < 1/(2^d)^{q+1}$, where $c_q$ is as in  \cref{L:ondiagbase}. From \eqref{squaredec} we obtain
  \begin{equation}\label{eq:2^dterms}
  \E\big( \big(\sum_{i \in \cS'_m} \cM_\eps(S_i)\big)^q \big) \le \frac1{(2^d)^{q+1}} \E( \cM_{2^{m} \eps}(S)^q ) + C_{m,q},
  \end{equation}
  where $C_{m,q}$ comes from \cref{L:offdiagbase}.
  Adding the contributions from all $2^d$ groups (and using the fact that $(x_1+\ldots + x_{2^d})^q \le (2^d)^{q-1} (x_1^q +  \ldots + x_{2^d}^q)$ by convexity), and bounding each of the $2^d$ terms by \eqref{eq:2^dterms},
  $$
  \E\big( \big(\sum_{i \in \cS_m} \cM_\eps(S_i)\big)^q \big) \le \frac1{2^d} \E(\cM_{2^{m} \eps}(S)^q ) + (2^d)^{q}C_{m,q}.
  $$
  Therefore, recalling that $M_\eps = \E(\cM_\eps(S)^q)$, we have
  $$
  M_\eps \le \frac1{2^d} M_{2^m \eps} + 2^{dq}C_{m,q}.
  $$
  Taking the sup over $\eps > \eps_0$, and since $2^m \eps \ge \eps$, we get
  $$
  \sup_{\eps > \eps_0}M_\eps  \le \frac1{2^d} \sup_{\eps > \eps_0} M_\eps + 2^{dq}C_{m,q}
  $$
  and hence
  $$
  \sup_{\eps > \eps_0}M_\eps  \le \frac{2^{dq}}{1- 1/2^d} C_{m,q}.
  $$
  We conclude proof for $q\in (1,2]$ that is, $n=2$, by letting $\eps_0\to 0$ and Fatou's lemma.
  %{Monotone convergence?}

  We now consider the general case, which is in fact very similar to when $n=2$. We use the fact that $q/n \le 1$ and thus, arguing as in \eqref{squaredec},
   \begin{align}
  (\sum_{i \in \cS'_m} \cM_\eps(S \cap S_i) )^q & \le
    \sum_{i_1, \ldots, i_n \in \cS'_m} \cM_\eps(S_{i_1})^{q/n} \ldots \cM_\eps(S_{i_n})^{q/n} \label{squaredec2}
  \end{align}
As before, we consider the on-diagonal (when all indices are equal) and off-diagonal terms separately.
%  $$
%\E [\Big(\sum_{i \in \cS'_m} \cM_\eps(S \cap S_i)\Big)^q]
%  $$
  The on-diagonal terms were already estimated in  \cref{L:ondiagbase}, and we have the same upper bound \eqref{E:ondiagbase} for all sufficiently large $m$ and all $\eps>0$. %where we obtained that for all sufficiently large $m$, and for all $\eps>0$,
 %   \begin{equation}\label{E:scaling_it}
%    \E [\Big(\sum_{i \in \cS'_m} \cM_\eps(S_i)^q\Big)]
 %\le c_q 2^{dm - \xi(q) m} \E [ \cM_{\eps2^m} (S)^q].
  %  \end{equation}
%for some constant $c_q$ by scaling.
For the off-diagonal terms, we obtain the following estimate.

  \begin{lemma}\label{L:offdiag_it}   Assume the set up of  \cref{T:finiteposmoments}. Then for any fixed $m$ and $q <2d/\gamma^2 $, there exists a constant $C_{m,q}$ independent of $\eps$ such that
  $$
 \E \big( \sum_{i_1, \ldots, i_n \in \cS'_m: i_1 \neq i_2 } \cM_\eps(S_{i_1})^{q/n} \ldots \cM_\eps(S_{i_n})^{q/n}\big) \le C_{m,q}.
  $$
  \end{lemma}

  \begin{proof}Note that by Jensen's inequality (since $q/n \le 1$), if $i_1 \neq i_2 \in \cS'_m$,
  $$
\E\big(  \cM_\eps(S_{i_1})^{q/n} \ldots \cM_\eps(S_{i_n})^{q/n}\big) \le \E \big( \cM_\eps(S_{i_1})\ldots  \cM_\eps(S_{i_n}) \big)^{q/n}
    $$
    As before, this expectation can be computed exactly. To begin with, we rewrite the index set $\{i_1, \ldots, i_n\}$ in a way that takes into account which indices are equal and which are distinct. Thus let $\{i_1, \ldots, i_n\} =\{j_1 , \ldots, j_p\}$ where the $j_k$ are pairwise distinct and $2\le p\le n$ (since at least two indices are distinct). Call $m_k$ the multiplicity of $j_k$ in $\{i_1, \ldots, i_n\}$, that is, the number of times $j_k$ is present in that set, so that $m_1 + \ldots + m_p = n$ (with $m_k \ge 1$ by assumption).
%    Hence
%    $$
 %   \cM_\eps(S_{i_1})^{q/n} \ldots \cM_\eps(S_{i_n})^{q/n} = \cM_\eps(S_{j_1})^{qm_1/n} \ldots \cM_\eps(S_{j_p})^{qm_p/n}
 %   $$
Then
    $$
    \E\big(  \cM_\eps(S_{i_1})\ldots \cM_\eps(S_{i_n})\big)  =
    \int_{x_1 \in S_{i_1}} \ldots \int_{x_n \in S_{i_n}} e^{(\gamma^2/2) \sum_{1\le k \neq \ell \le n} K_\eps(x_{k},x_{\ell})} \dd x_{1} \ldots \dd x_{n}.
    $$
When $x_k \in S_{i_k}, x_\ell \in S_{i_\ell}$ and $S_{i_k} \neq S_{i_\ell}$, the term $K(x_k, x_\ell) = - \log |x_k - x_\ell| + O(1)$ is bounded above by a constant $c_m$ since the cubes are separated by a minimum distance of $2^{-m}$.
Hence
\begin{align*}
    \E\big(  \cM_\eps(S_{i_1})\ldots \cM_\eps(S_{i_n})\big)
 &\le c'_m \prod_{k=1}^p \int_{S_{j_k}} \ldots \int_{S_{j_k}} e^{(\gamma^2/2) \sum_{1\le k \neq \ell \le m_k } K_\eps(x_{k},x_{\ell})} \dd x_1 \ldots \dd x_{m_k}\\
 & = c'_m \prod_{k=1}^p \E ( \big(\cM_\eps(S_{j_k})\big)^{m_k})
\end{align*}
Now, since $m_k \le n-1$ (as there are at least two distinct indices in the set $\{i_1, \ldots, i_n\}$), and since $S_{j_k} \subset S$, we have that
$$\E ( \big(\cM_\eps(S_{j_k})\big)^{m_k}) \le \E( \cM_\eps (S)^{n-1})$$
which, by the induction hypothesis, is uniformly bounded in $\eps$, by a constant depending only on $m$ and $q$. This concludes the proof of the lemma.
\end{proof}
Putting together \eqref{E:ondiagbase} and  \cref{L:offdiag_it}, we conclude the proof that $M_\eps$ is uniformly bounded for arbitrary $q < 2d/\gamma^2$, as in the case $q < 2 \wedge (2d/\gamma^2)$. This finishes the proof of  \cref{T:finiteposmoments}.
  \end{proof}

 We complement \cref{T:finiteposmoments} with two results. The first one shows that the condition $q< 2d/\gamma^2$ is sharp for the finiteness of the moment of order $q>0$. The second will show a partial result in the general framework of Gaussian multiplicative chaos with respect to a $d$ dimensional reference measure $\sigma$ (that is, satisfying \eqref{dim}). We start with the first result.

\begin{prop}
  \label{P:inf_pos_mom}
    Assume the set up of  \cref{T:finiteposmoments} (in particular, that $\sigma(\dd x) = \dd x$ is the Lebesgue measure on $\R^d$). Let $q > 2d/\gamma^2$. Then
  $$\E(\cM(S)^q) = \infty.$$
\end{prop}

\begin{proof}
  The proof argues by contradiction, and has the same flavour as \cref{T:finiteposmoments} but is much simpler (essentially, we can ignore the off-diagonal term). Suppose that for some $q> 2d/\gamma^2$, $\E(\cM(S)^q) < \infty$. By Kahane's inequality, there is no loss of generality in assuming that the Gaussian field $h$ is in fact an exactly scale invariant field $X$ satisfying  \cref{C:scaling2}. Then for any cube $S_i$ of sidelength $2^{-m}$, by \eqref{scalingeps} (or more precisely \eqref{prescaling}),
  $$
  \E ( (\cM(S_i))^q ) \asymp 2^{- m \xi(q)} \E( (\cM(S)  )^q).
  $$
  On the other hand, keeping the same notations as in the proof of \cref{T:finiteposmoments}, and since $(x+y)^q \ge x^q + y^q$ for $q>1$ and $x,y > 0$,
  $$
  \cM(S)^q \ge \sum_{i \in \cS_m} \cM(S_i)^q.
  $$
  Hence, taking expectations,
  $$
  \E(\cM(S)^q) \gtrsim 2^{dm - \xi(q)m} \E((\cM(S))^q)
  $$
  However, when $q> 2d/\gamma^2$, we have that $d - \xi(q) >0$. Since $m$ is arbitrary and the implicit constant does not depend on $m$, we get the desired contradiction.
\end{proof}

\subsection{Degeneracy of multiplicative chaos for $\gamma > \sqrt{2d}$}

As in the previous Section \ref{S:positivemomentsGMC} we continue to consider the Lebesgue case, so $\sigma(\dd x) = \dd x$ and $d= \mathbf{d}$, while 
the logarithmically correlated field $h$ satisfies the general assumptions of  \cref{SS:setup}. In this short section we give a (non-optimal) bound showing that if $\gamma > \sqrt{2d}$ then the measures $\cM_\eps$ converge to zero and so degenerate, in contrast to the case $\gamma < \sqrt{2d}$ studied in Theorem \ref{T:conv}. This shows that the value $\gamma = \sqrt{2d}$ is indeed a natural (sharp) boundary for the validity of that result. 

At the value $\gamma = \sqrt{2d}$ itself, it is still possible with additional renormalisation to define a \emph{critical} Gaussian multiplicative chaos which is a measurable function of the field $h$. See \cite{DRSV2, DRSV1} as well as \cite{JunnilaSaksman, Pow18chaos} for a construction, and see \cite{Powell} for a survey. In the case $\gamma > \sqrt{2d}$, even if we normalise more carefully than above, any limit would only be in law rather than in probability or almost surely; and thus will not be a measurable function of the field $h$. This is partly related to the phenomenon known as \emph{freezing} in the statistical physics literature. See in particular \cite{MadauleRhodesVargas, BarralJinRhodesVargas} and references therein for a more precise discussion related to supercritical Gaussian multiplicative chaos. \ind{Freezing} \ind{Gaussian multiplicative chaos!supercritical}
\ind{Gaussian multiplicative chaos!critical}

\begin{theorem}\label{T:GMCsuper_poly} Let $D \subset \R^d$ be a bounded domain, let $\sigma (\dd x) = \dd x$ be the Lebesgue measure and let $h$ be a  logarithmically correlated field $h$ satisfying the assumptions of  \cref{SS:setup}. Let $\gamma> \sqrt{2d}$ and let $\cM_\eps(\dd x) = \exp({\gamma h_\eps(x)  - \tfrac{\gamma^2}{2} \var [ h_\eps(x) ] }) \dd x$, as in \eqref{regmeas}. Then for any compactly supported $A \subset D$, $\cM_\eps(A)\to 0$ in probability as $\eps \to 0$. More precisely, for any $0< \alpha<1$, there exists $C>0$
 depending only on $A, D, \gamma, d$ and the law of $h$ (but not on $\eps$), such that 
\begin{equation}
\label{eq:GMCsuper_poly}
\E[ (\cM_\eps(A) )^\alpha ]\le C \eps^{\xi(\alpha) -d}, 
\end{equation}
where $\xi(\alpha) =  \alpha ( d + \gamma^2/2) - \gamma^2\alpha^2/2$. In particular, for $\alpha=\sqrt{2d}/\gamma\in (0,1)$ we have for $C$ depending on %As a consequence, we can find $0< \alpha<1$ depending only on 
	$A, D, \gamma, d$ and the law of $h$ (but not on $\eps$),  that
\begin{equation}
\label{eq:GMCsuper}
\E[ (\cM_\eps(A) )^\alpha ]^{1/\alpha} \le C \eps^{\psi }, \quad \text{where } 
\psi=(\frac{\gamma}{\sqrt{2}}-\sqrt{d})^2>0.
\end{equation}
 In particular, $\cM_\eps$ converges weakly to the null measure in probability. 
\end{theorem}

\begin{rmk}\label{R:Freezing}
The bound we obtain in \eqref{eq:GMCsuper} shows that $\cM_\eps(A)$ decays polynomially fast to zero, and is close to  optimal except for polylogarithmic terms: indeed, 
Madaule, Rhodes and Vargas showed in \cite[Corollary 2.3]{MadauleRhodesVargas} that
$$ (\log (1/\eps))^{\tfrac{3\gamma}{2\sqrt{2d}}} \eps^{-(\tfrac{\gamma}{\sqrt{2}} - \sqrt{d})^2} \cM_\eps(A)$$ 
converges -- necessarily in law, as mentioned above -- to a nondegenerate limit. 
The analysis in \cite{MadauleRhodesVargas} comes from a fine understanding of the behaviour of the maximum of the field; this connection to the behaviour of the maximum of the regularised field is at the heart of the freezing transition alluded to above. This is however not entirely captured by the rather rough scaling argument which is given below.
\end{rmk}
\ind{Freezing}

\begin{proof} 
By Kahane's inequality (Theorem \ref{T:Kahane}), we may without loss of generality assume that $h_\eps = X_{\eps, 10}$ (say) is the scale invariant field defined on all of $\R^d$ constructed in Theorem \ref{C:scaling2}. %is exactly scale-invariant, and is thus defined on all of $\R^d$. 
By increasing $A$ slightly we can then assume that $A$ is a cube, which we may (still without loss of generality) take to be the unit cube.
Since $\E [ h_\eps(x)^2 ] = - \log (1/\eps) + O(1)$, there is also no loss of generality in proving the result for the measure $\cM_\eps(\dd x) = \eps^{\gamma^2/2} e^{\gamma h_\eps(x)} \dd x$ instead of \eqref{regmeas}. (This turns out to have some advantages in the proof below). 

% Fix a positive integer $m$ and let $r = 2^{-m}.$ 
Let $m$ be such that $r=2^{-m}$ satisfies $2^{-m-1}\le \eps \le 2^{-m}$. We split the cube $A$ into a union of $r^{-d} = 2^{md}$ cubes of sidelength $r =2^{-m}$. We take these cubes to be open and non-overlapping and call the $C_i, i \in I$. Then (since the Lebesgue measure of the boundary of these cubes equals zero), 
$$
\cM_\eps (A) = \sum_{i\in I} \cM_\eps( C_i). 
$$   
As noted in the proof of Theorem \ref{T:finiteposmoments}, the function $x>0\mapsto x^\alpha$ is subadditive when $\alpha <1$. Hence we see that 
$$
\cM_\eps(A)^\alpha \le \sum_{i\in I} \cM_\eps(C_i)^\alpha. 
$$
Now, arguing as in \eqref{scalingeps} using the exact scale-invariance and translation invariance of $h$, we see that
$$
 \E [  \cM_\eps(C_i)^\alpha] =  r^{\xi(\alpha)}  \E [ \cM_{\eps /r} (A)^\alpha ],
$$
where $\xi (\alpha) = \alpha ( d + \gamma^2/2) - \gamma^2\alpha^2/2$. (Note that this is an equality rather than an identity valid only up to constants, as is the case in \eqref{scalingeps}, thanks to our choice of normalisation of $\cM_\eps$ which is slightly different from \eqref{regmeas}). We deduce that 
\begin{equation}\label{eq:scalingexactMeps}
\E [ \cM_\eps(A)^\alpha ] \le r^{\xi (\alpha) } r^{-d} \E [ \cM_{\eps /r} (A)^\alpha ] . 
\end{equation}
 %Let us check that since $\gamma>\sqrt{2d}$ we can find $\alpha <1$ such that $\xi(\alpha ) - d >0$. 
% At $\alpha = 1$, the value of $\xi(\alpha ) - d$ is $\xi(1) - d = 0$, whereas its derivative is $$\xi'(1) = d+ \gamma^2/2 - \gamma^2 = d- \gamma^2/2 < 0.$$ 
% Thus, fix any $0<\alpha<1$ such that $\xi(\alpha) - d >0$. Iterating \eqref{eq:scalingexactMeps} $n$ times, where $n= \lfloor \log (1/\eps) / \log (1/r)\rfloor$,  we see that 
% $$
 %\E [ \cM_\eps(A)^\alpha ] \le r^{n ( \xi ( \alpha ) - d) } \E [ \cM_{\eps/r^n} (A)^\alpha ].
% $$
Now,  by Jensen's inequality, $ \E [ \cM_{\eps/{r}} (A)^\alpha ] \le \E[ \cM_{\eps/{r}} (A) ]^\alpha \le C$, so that (since ${r} \asymp \eps$),
$$
 \E [ \cM_\eps(A)^\alpha ]\le C \eps^{\xi (\alpha) - d}.
 $$
As desired, this shows that the left hand side decays at least polynomially fast when $\alpha$ is such that $\xi(\alpha)-d>0$ (note that this is true for some range of $\alpha$ since at $\alpha = 1$ the value of $\xi(\alpha ) - d$ is $\xi(1) - d = 0$, whereas its derivative is $\xi'(1) = d+ \gamma^2/2 - \gamma^2 = d- \gamma^2/2 < 0$). The optimal decay for $\E[ (\cM_\eps(A) )^\alpha ]^{1/\alpha}$ is found by choosing $\alpha$ to maximise the function $\alpha \in [0,1] \mapsto (\xi (\alpha) - d)/\alpha$. This is achieved when %$\xi'(\alpha ) = 0$, i.e., when $d+ \gamma^2/2 - \alpha \gamma^2 = 0$ or, in other words, when 
	$\alpha =\sqrt{2d}/\gamma$, and then
$$
(\xi(\alpha ) - d)/\alpha = (\frac{\gamma}{\sqrt{2}}-\sqrt{d})^2, 
$$
so \eqref{eq:GMCsuper} follows.
\end{proof}

\subsection{Positive moments for general reference measures}
\label{S:posgen}
We now introduce the second result complementing  \cref{T:finiteposmoments}, which is an extension of  \cref{T:finiteposmoments} to the general setup of Gaussian multiplicative chaos relative to a $\mathbf{d}$-dimensional reference measure $\sigma$ (satisfying \eqref{dim}). In order to not make the exposition too cumbersome, we limit the proof to the case where $q < (2\mathbf{d}/\gamma^2) \wedge 2$ (hence, at least in the $L^1$ regime where $\gamma \in [\sqrt{\mathbf{d}}, \sqrt{2\mathbf{d}})$, there is no loss of generality at all).

Before doing so it may be useful to explain where the previous proof breaks down if $\sigma$ is not Lebesgue measure. The main issue lies in the scaling argument of  \cref{L:ondiagbase}; when we consider a cube of $S_i$ of $\cS'_m$ (sidelength $2^{-m}$), blowing this up by a factor $2^m$ will of course still produce a cube of unit sidelength, but the Gaussian multiplicative chaos is now with respect to a reference measure which is no longer $\sigma$, but instead reflects the local behaviour of $\sigma$ near the cube $S_i$. For very inhomogeneous fractals this behaviour could be wildly different, and so the inequality in that Lemma has no reason to hold true.

Instead we will need a different approach that accounts for the possible inhomogeneities of the fractal supporting the reference measure $\sigma$. The proof below comes from work (written roughly in parallel with this part of the book) in \cite{BerestyckiSheffieldSun}. It is based on Girsanov's lemma (Theorem \ref{T:Girsanov}).

\begin{prop}
  \label{P:finiteposmoments_gen} Let $S \subset D$ be bounded and open, and suppose that the reference measure $\sigma$ satisfies the dimensionality condition \eqref{dim}. Then if $0 <q < 2 \wedge (2\mathbf{d} / \gamma^2)$,
  \begin{equation}\label{E:finiteposmoments_gen}
  \E( \cM(S)^q) < \infty
  \end{equation}
  and moreover, $\cM_\eps(S)$ converges to $\cM(S)$ in $L^q$.
\end{prop}

\begin{proof}
Again, we can assume without loss of generality that $q>1$. Set $\delta=q-1 \in (0,1)$. Write
$$
\E (\cM({S})^{q} ) = \E(\cM(S) \cM(S)^\delta) = {\E(\cM(S))}\Q( \cM(S)^\delta)={\sigma(S)\E^*(\cM(S)^\delta)}
$$
where $\Q$ denotes the law of the field biased by $\cM(S)$.  Using Girsanov's lemma for Gaussian multiplicative chaos, \cref{T:Girsanov} {(and Corollary \ref{C:Girsanov2})}, we can rewrite this as
$$
\Q(\cM(S)^\delta) = \int_S \sigma(\dd x) \E\left( \Big(\int_S e^{ \gamma^2 K(x,y)  } \cM(\dd y)\Big)^\delta\right).
$$
Next, for each $n\ge 0$, let $A_n(x)$ denote the annulus at distance between $2^{-n}$ and $2^{-n-1}$ from $x$; that is, $A_n(x) = \{ y: |y-x| \in [2^{-n-1}, 2^{-n})\} $. Then using the fact that
$K(x,y)  = - \log |x-y| + O(1)$ and the inequality $(a_1 + \ldots + a_n)^\delta\le a_1^\delta + \ldots + a_n^\delta$ for $\delta <1$ and $a_i >0$,
we see that
\begin{align*}
\Q(\cM(S)^\delta)  & \le  C\sum_{n=0}^\infty \int_S \sigma(\dd x)  \E \left( \Big(\int_{A_n(x)} e^{-\gamma^2 \log |x-y| } \cM(\dd y) \Big)^\delta\right) \\
& \le C \sum_{n=0}^\infty 2^{n \gamma^2\delta}  \int_S \sigma(\dd x) \E ( \cM(A_n(x))^\delta).
%& \le C \sum_{n} 2^{n \gamma^2 \delta} \int_A\sigma(dx) \sigma(A_n(x))^\delta
\end{align*}
%Now by possibly losing a constant, we can assume the underlying field is exactly scale invariant around $x$ by Kahane's inequality,
Now fix $x$, and consider a field $X$ which is exactly scale invariant around $x$ as in \cref{S:hdc}. Hence for any $\lambda<1$,
$$
(X(x+\lambda z))_{z \in S}=(\tilde X(z)+\Omega_{-\log \lambda})_{z \in S},
$$
where $\tilde X$ has the same law as $X$ and $\Omega_{r}$ is a Gaussian with variance $r$ independent of $\tilde X$. Write $X_\eps$ for the field truncated at level $\eps$, as in \cref{S:hdc}.
%Let $\cM_\eps (dx) = e^{\gamma X_\eps(x) - (\gamma^2/2) \E ( X_\eps(x)^2) } dx \asymp e^{\gamma X_\eps(x)} \eps^{\gamma^2/2} dx$, and let

Set $\lambda = 2^{-n} \le 1$. Denote by $\sigma_{\lambda, x} (dz)$ the image measure of $\sigma$ under the map $y=x+\lambda z\mapsto z$ (so that the total mass $\sigma_{\lambda, x} (A_1(0)) = \sigma(A_n(x))$). By applying \cref{C:scaling2} and changing variables $y\mapsto z$, we obtain:
\begin{align*}
%\E [\cM_{\eps 2^{-n} } (A_n(x))^\delta] &=
\E \left(\left( \int_{A_n(x)} e^{\gamma X_{\lambda \eps} (y)} (\lambda \eps)^{\tfrac{\gamma^2}{2}} \sigma(\dd y) \right)^\delta \right)&\le  \lambda^{\tfrac{\gamma^2\delta }{2}} \E \left( \left( \int_{{A_1(0)}} e^{\gamma X_{ {\eps}} ({x+}\lambda z)} \eps^{\tfrac{\gamma^2}{2}} \sigma_{\lambda, x} (\dd z) \right)^\delta\right)\\
 & =   \lambda^{\tfrac{\gamma^2\delta }{2}} \E \left( e^{\delta \gamma \Omega_{-\log \lambda}} \left( \int_{A_1(0)} e^{\gamma  \tilde{X}_\eps(z)} \eps^{\tfrac{\gamma^2}{2}} \sigma_{\lambda, x} (\dd z)\right)^\delta \right) \\
& = \lambda^{\tfrac{\gamma^2 \delta }{2}} e^{-\tfrac{\delta^2 \gamma^2 \log (\lambda)}{2}} \E \left(  \left( \int_{A_1(0)} e^{\gamma {\tilde{X}_\eps(z)}} \eps^{\tfrac{\gamma^2}{2}} \sigma_{\lambda, x} (\dd z)\right)^\delta \right)\\
& \le \lambda^{\tfrac{\gamma^2\delta}{2} - \tfrac{\delta^2 \gamma^2 }{2} }  \E \left(  \left( \int_{A_1(0)} e^{\gamma  \tilde{X}_\eps(z)} \eps^{\tfrac{\gamma^2}{2}} \sigma_{\lambda, x} (\dd z)\right) \right)^\delta,
\end{align*}
where the last inequality is by Jensen's inequality since $\delta < 1$.

By Kahane's inequality (\cref{T:Kahane}) and comparing to the trivial field, there exists an absolute constant $C>0$ such that
\begin{align*}
\E(\cM_{\eps 2^{-n} } (A_n(x))^\delta) & \le C \lambda^{\gamma^2\delta/2 - \delta^2 \gamma^2 /2 }  \sigma(A_n(x))^\delta.
\end{align*}
Letting $\eps\to 0$, we get that for any ${n\ge 0}$,
$$
 \E(\cM(A_n(x))^\delta) \le C 2^{- n (\gamma^2\delta/2 - \delta^2 \gamma^2 /2 )} \sigma(A_n(x))^\delta.
$$
We deduce that
$$
\Q(\cM(S)^\delta) \le C \sum_n 2^{n (\gamma^2 \delta/2 + \delta^2 \gamma^2 /2 ) } \int \sigma(\dd x) \sigma(A_n(x))^\delta.
$$
To estimate the last integral, let $\bar y$ and $\bar x$ be two independent points distributed according to $\sigma(\cdot \cap S)/ \sigma(S)$. Then note that
\begin{align*}
 \sigma(A_n(x))&  \le \sigma(S)\P (|\bar y- x |  \le 2^{-n})
 \end{align*}
 so that by Jensen's inequality again (as $\delta < 1$),
 \begin{align*}
  \int_S\sigma(\dd x) \sigma(A_n(x))^\delta & \le  \int_S \sigma(\dd x) \sigma(S)^\delta \P \Big( | \bar y - \bar x | \le 2^{-n} \Big| \bar x = x\Big)^\delta  \\
  &\le
  \sigma(S)^{\delta+1} \E ( \P\Big( |\bar x-\bar y| <  2^{-n} \Big| \bar x \Big)^\delta)\\
  & \le \sigma(S)^{\delta+1}\P( |\bar x-\bar y| \le 2^{-n} )^\delta\\
  & \le  \sigma(S)^{\delta+1} \E ( |\bar x-\bar y|^{-{\mathbf{d}}})^\delta 2^{-n{\mathbf{d}} \delta}.
 \end{align*}
by Markov's inequality.
%{Should it be $\sigma(D)^{\delta+1}$ on the RHS? Because we need to apply Jensen to a probability measure?}
Now $\delta < 2\mathbf{d}/\gamma^2 -1$ implies that $\gamma^2 \delta/2 + \delta^2 \gamma^2 /2-\mathbf{d}\delta=\delta(  \gamma^2/2 - \mathbf{d}  + \gamma^2\delta/2) <0$.
Putting everything together, we can find $c=c(\mathbf{d},\gamma,\delta)$ such that
$$
\E(\cM(S)^q)=\sigma(S)\Q(\cM(S)^\delta) \le  c(\delta) \E(|\bar x-\bar y|^{-\mathbf{d}})^\delta <\infty %c(\delta) \cE(\sigma)^\delta . % \sum_n 2^{n( \gamma^2 \delta/2 + \delta^2 \gamma^2 /2 ) - n d\delta }
$$
{by \eqref{dim2}.} This concludes the proof.
%\nb{Because we need the $\mathbf{d}$-dimensional energy to be finite, I changed the definition of $\mathbf{d}$ in \eqref{dim2}. This matches the definition of \cite{BerestyckiSheffieldSun}.}
\end{proof}

\subsection{Negative moments of Gaussian multiplicative chaos}

We now turn our attention to negative moments of the chaos measures. We will first show in Proposition \ref{P:small_neg_mom} that in the general set up, $\cM(S)$ admits moments of order $q\in [q_0,0]$ for some $q_0<0$. This proof is based on a similar argument appearing in \cite{GHSS}. We will then explain how to bootstrap this to get existence of all negative moments (see \cref{T:negmom}). Note that, in particular, this implies strict positivity of the measures with probability one. %This strategy should adapt to show the same result for all dimensions, as stated in \cite{RobertVargas}, but we omit this here as the argument seems less straightforward {*change as appropriate}.

We work in the general setting: $\sigma$ is a Radon measure with dimension at least $\mathbf{d}$ ($0< \mathbf{d}\le d$) and  $0 \le \gamma<\sqrt{2\mathbf{d}}$. The first ingredient we will need concerns the $\beta$-dimensional energy of the measure $\cM$.

\begin{lemma}\label{L:quantumenergy}
	Assume that $\sigma(D)>0$. Suppose that $0 \le \beta<\mathbf{d}\vee \sqrt{2\mathbf{d}}\gamma$, and $x$ is a point chosen from the measure $\sigma(\dd x)$ in $D$ (normalised to be a probability measure), independently of the field. Then $$\int_{D} |x-y|^{-\beta} \cM(\dd y)<\infty$$ almost surely and in fact, has finite $r$th moment for $r>0$ small enough.
\end{lemma}

\begin{proof}
If $\sqrt{2 \mathbf{d}} \gamma \le \mathbf{d}$ then $\beta < \mathbf{d}$, in which case
this energy will have finite expectation directly by assumption \eqref{dim}. So let us assume that $\sqrt{2 \mathbf{d}} \gamma > \mathbf{d}$, and thus $\beta < \sqrt{2 \mathbf{d}} \gamma$.
This means that for $r>0$ small enough we will have  $1>1-r>1/2 \vee \beta^2/(2\mathbf{d}\gamma^2)$.
	For such an $r$, we bound
	\begin{align*}
	\E\left((\int_D |x-y|^{-\beta} \cM(\dd y))^r\right) & \leq C  \E\left( \int_D \left(\int_D \e^{\beta K(x,y) } \cM(\dd y)\right)^r \frac{\sigma(\dd x)}{\sigma(D)} \right)  \\ & \leq C \E (\cM(D)^r \cM_{\gamma^{-1}\beta}(D))
	\end{align*}
	for some finite $C$, where the last inequality follows from Girsanov's lemma {(Theorem \ref{T:Girsanov} and Corollary \ref{C:Girsanov2})} and $\cM_{\gamma^{-1}\beta}$ is the chaos measure of the field with parameter $\gamma^{-1}\beta$ rather than $\gamma$ (note that by our assumptions on the parameters, we have $\gamma^{-1} \beta < \sqrt{2 \mathbf{d}}$, so this chaos is indeed well defined and non-trivial).
	
	Now by H\"{o}lder's inequality with $p=r^{-1}$ and  $q=(1-r)^{-1}$, the above is less than or equal to
	$$\E(\cM(D))^r \E(\cM_{\gamma^{-1}
	\beta}(D)^{\frac{1}{1-r}})^{1-r}.$$
By \cref{P:finiteposmoments_gen}, this is finite as long as $(1-r)^{-1}\le 2 \wedge 2\mathbf{d}/(\gamma^{-1}\beta)^2$, which is exactly our assumption on $r$.
\end{proof}

\begin{cor}\label{C:sing_mass_small}
	Take the same set up as in \cref{L:quantumenergy}. Then there exists $M$ large enough, depending only on $\gamma$ and $\mathbf{d}$, such that $$\P(E_s):=\P \left( \int_{B(x,s^{-M})} \e^{\gamma^2 K(x,y)} \cM(\dd y) \le \frac{1}{s}\right) \ge \frac{1}{2} $$ for all $s$ sufficiently large.
\end{cor}

\begin{proof}
	By the assumptions in Section \ref{SS:setup} on $K$, we know that $\e^{\gamma^2 K(x,y)}\le c |x-y|^{-\gamma^2}$ for some $c<\infty$. Writing $|x-y|^{-\gamma^2}=|x-y|^{-(\gamma^2+2/M)}|x-y|^{2/M}$, we therefore have that $$\e^{\gamma^2 K(x,y)}\le cs^{-2} |x-y|^{-(\gamma^2+2/M)} \text{ for all } y\in B(x,s^{-M}).$$ Hence
	$$\P \left( \int_{B(x,s^{-M})} \e^{\gamma^2 K(x,y)} \cM(\dd y) \le \frac{1}{s}\right)\ge \P\left(\int_D|x-y|^{-(\gamma^2+2/M)} \cM(\dd y)\le \frac{s}{c}\right).$$
If $M$ is such that $\beta := \gamma^2+2/M<\mathbf{d}\vee \sqrt{2\mathbf{d}} \gamma$ (it is always possible to choose $M$ in this manner, consider separately the cases $\gamma \le \sqrt{\mathbf{d}}$ and $\mathbf{d} < \gamma < \sqrt{2 \mathbf{d}}$), then by \cref{L:quantumenergy}, the right hand side converges monotonically to $1$ as $s\to \infty$.	
\end{proof}

From here the key observation is the following. If we write $\Q$ for the field biased by $\cM(D)/\sigma(D)$ as before, then for $s>0$, $$\Q(\exp(-s \cM(D)))=\sigma(D)^{-1}\E(\cM(D)\exp(-s\cM(D)))\le \frac{e^{-1}}{\sigma(D)s}$$
simply because $x\e^{-sx}\le e^{-1}/s$ for all positive $x,s$. This says that under $\Q$, $\cM(D)$ is unlikely to be too small. Of course we would actually like such a statement under $\mathbb{P}$. The trouble is that the field has an extra log singularity under $\Q$, and so it could be this that saves $\cM(D)$ from being very small. The work now is essentially to rule this out using \cref{C:sing_mass_small}.

To do this, we first claim that if $E_s$ is the event in \cref{C:sing_mass_small}, then
	\begin{equation}
	\label{eqn:laplace_event} \E(\exp(-c s^{1+ M\gamma^2}\cM(D))\mathbf{1}_{E_s})\le \frac{C}{s\sigma(D)}\end{equation}
	for some $c,C,s_0<\infty$ %(depending only on $\gamma, \mathbf{d}$)
	and all $s\ge s_0$, where these constants depend on $\mathbf{d}, \gamma$ and also $K$. Indeed, by Girsanov's lemma again (Theorem \ref{C:Girsanov}),
	$$\Q(\exp(-s\cM(D)))=\E(\exp(-s \int_D e^{\gamma^2K(x,y)}\cM(\dd y)))$$
	where under $\P$, $x$ is a point chosen according to $\sigma$, independently of the field. Moreover, on the event $E_s$,
$$s\int_D \e^{\gamma^2 K(x,y)} \cM(\dd y)\le 1+{c}s^{1+ M\gamma^2} \cM(D)
$$ for some $c<\infty$. This implies that 
$$\exp(-{c}s^{1+ M\gamma^2}\cM(D))\mathbf{1}_{E_s}\le e^{-1} \exp(-s \int_D e^{\gamma^2K(x,y)}\cM(\dd y)).$$ 
Taking expectations, we get
\begin{align*}
\E (\exp(-{c}s^{1+ M\gamma^2}\cM(D))\mathbf{1}_{E_s})& \lesssim \E [\exp(-s \int_D e^{\gamma^2K(x,y)}\cM(\dd y))]\\
& = \Q[ (\exp(-s\cM(D)))]\\
& = \frac1{s \sigma(D)} \E [ {s\cM(D)}\exp ( - s \cM(D))]
\end{align*}
 and so \eqref{eqn:laplace_event} holds.
	
	Note that if it weren't for the indicator function in \eqref{eqn:laplace_event}, this would imply that $\cM(D)$ has \emph{some} finite negative moments {(using the identity $y^{-p}=\Gamma(p)^{-1} \int_0^\infty t^{p-1}\exp(-ty) \, \dd t$ for $p>0$, this will be detailed below)}. On the other hand, we have shown in \cref{C:sing_mass_small} that the event in the indicator function is rather likely. Putting these ideas together more precisely, we obtain the following.
	
	\begin{prop}\label{P:small_neg_mom}
		Assume that $\sigma(D)>0$ and $0 \le \gamma<\sqrt{2\mathbf{d}}$. For some $q_0<0$, depending only on $\gamma$ and $\mathbf{d}$, it holds that $\E(\cM(D)^{q_0})<\infty$. %In fact, there exists $c>0$ such that $\E(\cM_\eps(D)^{q_0})<c$ for all $\eps>0$ small enough. 
	\end{prop}

\begin{proof} %For this proof it will help to work with an exactly scale invariant field as in \cref{S:hdc}. %Since $\sigma(D)>0$ by assumption, we may find a point $x$ and a ball $B(x,r)$ around $x$ with $\sigma(B(x,r))>0$. Note that it suffices to show that $\E(\cM(B(x,r)^q))<\infty$ for some $q_0$ and all $q\in [q_0,0]$, since $\cM(B(x,r))\le \cM(D)$.
%	By Kahane's inequality (Theorem \ref{T:Kahane}, because $x\to x^{q}$ for $q<0$ is convex) we can equivalently prove the proposition with such a field. {*More details} So from now on, we assume that setting.
	
Let us first observe that, without loss of generality, we may assume that $K(x,y)\ge 0$ for all $x,y\in D$. Indeed, we can always find some $D'\subset D$ with $\sigma(D')>0$ and $K(x,y)\ge 0$ for all $x,y\in D'$ (since $K$ diverges logarithmically near the diagonal), and then it clearly suffices to show that $\E(\cM(D')^{q_0})<\infty$. Note that $\sigma$ also has dimension at least  $\mathbf{d}$ when restricted to $D'$.

The advantage of assuming this setup, is that we can make use of the following tool (see for example, \cite{FKG}): \ind{FKG inequality}
	
	\begin{theorem}[FKG inequality]\label{T:FKG}
		Let $(Z(x))_{x\in U}$ be an almost surely continuous centred Gaussian field on $U\subset\R^d$ with $\E(Z(x)Z(y))\ge 0$ for all $x,y\in U$. Then, if f, g are two bounded, increasing measurable functions,
		$$\E \big( f((Z(x))_{x\in U})g((Z(x))_{x\in U})\big) \ge \E \big( f((Z(x))_{x\in U})\big)  \E \big( g((Z(x))_{x\in U})\big).$$
	\end{theorem}
	%Since $\cM(D)>\cM(D')$ for $D'\subset D$, we may assume that $D$ is a small enough neighbourhood of the origin that $K(x,y)=\log_+(|x-y|)$ is positive definite in $D$ {(ref)}. Since the function $x\mapsto x^{q}$ for $q<0$ is convex, we may use Kahane's convexity inequality to reduce ourselves to dealing with this covariance kernel. {(more details)} In fact, the only property we will use is that $K(x,y)\ge 0$ for all $x,y$. This ensures {(FKG...)}
	
	To apply this, we need to work with continuous fields, so let us consider the regularised field $h_\eps$ of \eqref{convol} and regularised measure $\cM_{\eps}$ of \eqref{regmeas}, and denote by $E_s^\eps$ the event $E_s$ of Corollary \ref{C:sing_mass_small} with $\cM$ replaced by $\cM_{\eps}$ (and we still define $E_s^\eps$ in terms of a point that is sampled independently of the field and with probability proportional to $\sigma(\dd x)$). Since $h_\eps$ is almost surely continuous and the functions $\mathbf{1}_{E_s^\eps}$ and $\exp(-cs^{1+M\gamma^2}\cM_\eps(D))$ are both bounded, decreasing functions of the field $h_\eps$, we can apply Theorem \ref{T:FKG} to see that
	$$\E(\exp(-cs^{1+ M\gamma^2}\cM_\eps(D))\mathbf{1}_{E^\eps_s})\ge \E(\exp(-cs^{1+M\gamma^2}\cM_\eps(D)))\P(E^\eps_s)$$
	for all $\eps$. (Recall that $\E$ is over the field as well as the independent random point $x$, so we actually apply the FKG inequality conditionally given $x$, then note that the first term in the right hand does not depend on $x$). By dominated convergence, we therefore obtain that
	$$\E(\exp(-{c}s^{1+ M\gamma^2}\cM(D))\mathbf{1}_{E_s})\ge \E(\exp(-{c}s^{1+ M\gamma^2}\cM(D)))\P(E_s).$$
	Hence by \cref{C:sing_mass_small} and \eqref{eqn:laplace_event}, for $M$ large enough (depending only on $\gamma, \mathbf{d}$) and $s_0$ large enough (depending on $\gamma, \mathbf{d}$ \emph{and} $K$):
		$$
\E(\exp(-{c}s^{1+ M\gamma^2}\cM(D)))\le \frac{2C}{s\sigma(D)} \quad \quad \forall s\ge s_0,
$$
or to put it another way, for some $t_0$ sufficiently large,
\begin{equation}\label{laplace_useful}
\E(\exp(-t \cM(D)))\le \frac{2C}{{(t/c)}^{1/(1+ M\gamma^2)}\sigma(D)} \quad \quad \forall t\ge t_0.
\end{equation}
Finally, since $y^{-p}=\Gamma(p)^{-1} \int_0^\infty t^{p-1}\exp(-ty) \, \dd t$ for $p>0$, this implies that
\begin{align*}
\E(\cM(D)^{-p}) & = \int_0^\infty t^{p-1} \E( e^{ - t \cM (D)}) \dd t \\
& \lesssim 1 + \int_{t_0}^\infty t^{p-1 - 1/(1+  M \gamma^2)} \dd t.
\end{align*}
The integral in the right hand side is finite as soon as $p< 1/(1+ M\gamma^2)$, and so for such values of $p$ we get $\E(\cM(D)^{-p})<\infty$. Note that this only depends on $M,\gamma$, so since $M$ is a function of $\gamma$ and $\mathbf{d}$, the obtained $q_0$ also depends only on $\gamma$ and $\mathbf{d}$. This completes the proof of Proposition \ref{P:small_neg_mom}. %{*Write proof of $\varepsilon$ case*}
\end{proof}

Now we explain how to extend this to all negative moments, using an iteration procedure. This idea first appeared in the setting of multiplicative cascade measures (a toy model for multiplicative chaos) in \cite{molchan}.

\begin{theorem}\label{T:negmom}
Suppose that $\sigma(D)>0$ and $0 \le \gamma < \sqrt{2 \mathbf{d}}$. Then $$\E(\cM(D)^q)<\infty$$ for all $q<0$. %In fact, for any $q<0$ there exists a constant $c_q$ such that $\mathbb{E}(\cM_\eps(D)^q)\le c_q$ for all $\eps>0$ small enough. 
\end{theorem}

We emphasise that we need only our standing assumptions on the measure $\sigma$ and the field $h$ from Section \ref{SS:setup} here (as long as $\sigma(D)>0$), and that there are no restrictions on $d$ or $\mathbf{d}>0$.

\begin{proof}
To begin with, note that since $\sigma$ does not have any atoms, we can find two distinct points $x_1, x_2$ in the support of $\sigma$. Therefore we can find open sets $D_1$ and $D_2$ such that $x_1 \in D_1, x_2 \in D_2$, $\bar D_1 \cap  \bar D_2 = \emptyset $ and $\sigma(D_1) \sigma(D_2)>0$. Furthermore, by the assumption on $K$ (more precisely, the continuity of $g$), we may assume that $K(x,y)\le C$ whenever $x\in D_1$, $y\in D_2$.

The key point is that by \cref{P:small_neg_mom}, there exists $q_0<0$ such that $\E(\cM(D_i)^{q})<\infty$ for all $q\in [q_0,0]$ and $i=1,2$. Indeed we have seen that $q_0$ depends only on $\mathbf{d},\gamma$, as long as the base measure has strictly positive mass.

The idea to make use of this is to note the trivial bound $\cM(D)\ge \cM(D_1)+\cM(D_2)$, and then apply the AM-GM inequality to see that $\cM(D)\ge \sqrt{\cM(D_1)\cM(D_2)}$. This gives that
\[ \E(\cM(D)^{q})\le \E(\cM(D_1)^{q/2}\cM(D_2)^{q/2}) \]
for $q<0$. If $\cM(D_1)$ and $\cM(D_2)$ were independent, we could factorise the right hand side and choose $q=2q_0$, therefore showing that negative moments exist with orders in the larger interval $[2q_0,0]$. We could then iterate to get all negative moments.

The problem of course is that they are not actually independent. To get around this we will use the assumption that $K(x,y)\le C$ for $x\in D_1, y\in D_2$, together with Kahane's inequality (Theorem \ref{T:Kahane}).

More precisely, let us denote our field restricted to $D_1\cup D_2$ by $X$. Let us also define a Gaussian field $Y$ on $D_1\cup D_2$ by setting it equal to $Y_1+Y_2+Z$ where: $Y_1,Y_2$ are independent; $Y_1$ has the law of $X|_{D_1}$ on $D_1$ and is $0$ on $D_2$; $Y_2$ has the law of $X|_{D_2}$ on $D_2$ and is $0$ on $D_1$; and $Z$ is an independent normal random variable with variance $C$. Then the covariance kernel of $Y$ dominates (pointwise) the covariance kernel of $X$. Since polynomials of negative order are convex, we can apply Kahane's inequality (\cref{T:Kahane}) (and a limiting argument so that we can compare the respective GMC measures) to obtain that \begin{align*}\E((\cM(D_1)+\cM(D_2))^q) &\le \E((\cM_Y(D_1)+\cM_Y(D_2))^q) \\ & \le \E(\cM_Y(D_1)^{q/2}\cM_Y(D_2)^{q/2}) \\ & =\E(\e^{\frac{q}{2}(\gamma Z-\frac{\gamma^2}{2}C)})\E(\cM_Y(D_1)^{q/2})\E(\cM_Y(D_2)^{q/2}),\end{align*}
where we have applied AM-GM in the second line. By construction, if $q\in[2q_0,0]$, the right hand side is finite. So we obtain that $\E(\cM(D)^q)<\infty$ for all $q\in[2q_0,0]$. Repeating the argument one obtains the existence of any negative moment.% {*Write proof of $\varepsilon$ case*}
 \end{proof}

%\nb{Can we say the negative moments of the approximations are uniformly bounded in $\eps$? This is useful in LCFT.  I think it follows from a martingale approximation + Jensen which implies monotonicity of the negative moments.}

%{*Add something about lower tail being lognormal and ref. exercises}

Since $\cM(D)$ has finite negative moments of all orders (as shown by the previous theorem), we deduce that the tail at zero of $\cM(D)$ decays faster than any polynomial. It is natural to wonder whether the decay can be characterised precisely. A lognormal upper bound for this decay (meaning, $\P( \cM(D) \le \delta) \le \exp ( - c (\log 1/\delta)^2)$ for some $c>0$) was first established in \cite{DuplantierSheffield}, see also \cite{Arunotes}. In some one dimensional cases of Gaussian multiplicative chaos, the exact law of the total mass is in fact known (this was obtained by Remy \cite{Remy}, proving a well known conjecture of Fyodorov and Bouchaud \cite{FB}). In exercise \eqref{E:log}, we propose a lognormal lower bound valid in great generality.

\subsection{KPZ theorem}
\label{S:KPZ}

In this section, we consider the Gaussian free field with zero boundary conditions in a domain $D \subset \R^2$. The KPZ formula relates the ``quantum'' and ``Euclidean'' sizes of a given set $A$, which is either deterministic, or random but independent of the field. This often has a particularly natural interpretation  in the context of discrete random planar maps and critical exponents; see \cref{S:exponents}.
Concrete examples are given in \cref{S:maps}.

We will first formulate the KPZ theorem using the framework of Rhodes and Vargas \cite{RhodesVargas}. This article appeared simultaneously with (and independently from) the paper by Duplantier and Sheffield \cite{DuplantierSheffield}. The results of these two papers are similar in spirit, but the version we present here is a bit easier to state, and in fact stronger. The formulation (and sketch of proof) corresponding to \cite{DuplantierSheffield} will be given in \cref{S:DS}. We will also include a version, due to Aru \cite{Aru}. Although this last statement is weaker, its proof is completely straightforward given our earlier work.

We first introduce the notion of \emph{scaling exponent} of a set $A$ (in the sense of \cite{RhodesVargas}), starting with the Euclidean version. Let $A \subset D$ be a fixed Borel set and write $\dd_H(A)$ for the (Euclidean) Hausdorff dimension of $A$. Since $0 \le \dd_H(A) \le 2$, we may write
\begin{equation}\label{scalexponent}
\dd_H(A) = 2 (1- x),
\end{equation}
for $x \in [0,1]$. The number $x$ is called the \textbf{(Euclidean) scaling exponent} of $A$.
%The scaling exponent of $A$, tells us how likely it is, if we pick a point $z$ uniformly in $\D$, that $z$ falls within distance $\eps$ of the set $A$.
\ind{Scaling exponents}

We now define the quantum analogue of the scaling exponent. Let $$C_\delta(A):=\inf\{\sum_i \cM(B(x_i,r_i))^\delta:   \{B(x_i,r_i)\}_i \text{ is a cover of } A\} ,$$
so that $C_\delta(A)$ can be viewed as a (multiple) of the quantum Hausdorff content of $A$. We now define 
$$\dd_{H,\gamma}(A)=\inf\{\delta>0: C_\delta(A)=0\} \in [0,1]$$
and call $d_{H, \gamma} (A)$ its ``quantum Hausdorff dimension''. Finally, we define the \textbf{quantum scaling exponent} $\Delta$ by 
$$
\Delta = 1- \dd_{H, \gamma}(A).
$$
The terms ``quantum Hausdorff dimension" and content should perhaps be qualified, for the following reasons.

\begin{enumerate}

\item Although it does not feature in these notes, a random metric associated with $e^{\gamma h}$ ($h$ a GFF in $D$) has recently been constructed in a series of works culminating with \cite{DDDF,GM,GM2}.
The Hausdorff dimension $d_\gamma$ of $D$ equipped with this random metric is currently {unknown}, except for the special case \{$\gamma = \sqrt{8/3}$, $d_\gamma = 4$\}. The general bound $d_\gamma>2$ is also known, as well as more precise estimates: see \cite{DingGwynne,PGbounds}.

\item Under this random metric, the actual value of the Hausdorff dimension of $A
\subset D$ is then given by
$$
{d_\gamma} (1- \Delta).
$$
Again it always holds that $\Delta \in [0,1]$,  and note the analogy with \eqref{scalexponent}.

\item Recently, a metric version of the KPZ formula has been obtained by Gwynne and Pfeffer \cite{GwynnePfeffer}; more details concerning the relation between scaling exponent and Hausdorff dimension can be found there.
\end{enumerate}

\begin{rmk}
There is no consensus (even in the physics literature) about the value of {$d_\gamma$}. {Until recently it seemed that the prediction $$d_\gamma = 1+ \frac{\gamma^2}{4} + \sqrt{(1+ \frac{\gamma^2}4)^2 + \gamma^2}$$ by Watabiki \cite{Watabiki} had a reasonable chance of being correct, but it has now been proved false -- at least for small $\gamma$ \cite{DGWatabiki}.}
Simulations are notoriously difficult because of large fluctuations. As mentioned earlier, the only value that is known rigorously is when $\gamma=\sqrt{8/3}$. In this case the metric space is described by the Brownian map \cite{Miermont,LeGall,LQGandTBMIII} and the Hausdorff dimension is equal to 4.
\ind{Watabiki formula}
\end{rmk}

%\item In the definition of $\Delta$, the probability $\P( B^\delta(z) \cap A   \neq \emptyset)$ is averaged over everything: $z$, $h$ and $A$ itself if it is random. So it is really an annealed probability.

%\end{enumerate}

We are now ready to state the KPZ theorem in this setup.

\ind{KPZ relation}

\begin{theorem}[{Almost sure Hausdorff KPZ formula}]\label{T:KPZas}
{Suppose that $A$ is deterministic and that $\gamma\in (0,2)$. Then, almost surely it holds that \begin{equation*}
	x = \frac{\gamma^2}{4} \Delta^2 + (1- \frac{\gamma^2}4) \Delta.
	\end{equation*}}	
		\end{theorem}

We will not prove this result and refer to \cite{RhodesVargas} for details. (We will, however, soon see the proof of a closely related result due to Aru \cite{Aru}). We make a few observations.
\begin{enumerate}

\item $x= 0,1$ if and only if $\Delta = 0,1$.

\item This is a quadratic relation with positive discriminant so can be inverted.

\item In the particular but important case of uniform random planar map scaling limits (see \cref{S:maps}), $\gamma  = \sqrt{8/3}$ and so the relation is
\begin{equation}
x = \frac23 \Delta^2 + \frac13 \Delta.
\end{equation}
\end{enumerate}

As we have already mentioned, various forms of the KPZ relation
have now been proved; the above statement comes from the work of Rhodes and Vargas \cite{RhodesVargas}.
Other versions can be found in the works of Aru \cite{Aru}, Duplantier and Sheffield \cite{DuplantierSheffield} which will both be discussed later in this chapter. See
also works of Gwynne and Pfeffer \cite{GwynnePfeffer}
for a KPZ relation in the sense of metric (Hausdorff) dimensions; Gwynne, Holden and Miller
\cite{GwynneHoldenMiller} for an effective KPZ formula which can be
used rigorously for determining a number of SLE exponents, and Berestycki, Garban, Rhodes and Vargas \cite{HKPZ} for a KPZ
relation formulated using the Liouville heat kernel.
%TO ADD (The transition probabilities of Liouville Brownian motion will be discussed
%later in \cref{S:LBM}).

   %However, almost sure versions also exist -- let us mention, in particular, the work of Rhodes and Vargas \cite{RhodesVargas} who proved a theorem which is in some ways stronger than the above (it is more robust, and deals with an almost sure notion of dimension) and which appeared simultaneously to the paper of Duplantier and Sheffield \cite{DuplantierSheffield} from which the above theorem is taken.
% More recently, a version of the KPZ formula was formulated and proved using the so-called Liouville heat kernel \cite{HKPZ}, thereby eliminating the need to reference the underlying Euclidean structure in the notion of dimension.

 %{We will conclude this chapter with a proof of Aru's result \cite{Aru}, and a sketch proof of \cref{T:KPZ}. The proof of \cref{T:KPZas} is omitted: see the paper \cite{RhodesVargas} for more details (this actually works for a much more general class of random measures).}

\subsubsection{Proof in the case of expected Minkowski dimension}

We now state Aru's version of the KPZ formula \cite{Aru} which, as already mentioned, has a straightforward proof given our earlier work.
%. The second is close to part of the argument by Duplantier and Sheffield \cite{KPZ}. Both are instructive in their own right.
%Aru's proof relies on the formula giving the power law spectrum of Liouville measures, which we state now but establish later in Section \ref{S:multifractal}. \red{change text here}
%\ind{Scaling relation for Liouville measure}
%\ind{Multifractal spectrum of Liouville measure}
%We defer the proof of this proposition until the end of the section, and content ourselves in saying it is simply the result of approximate scaling relation between $h_r(z)$ and $h_{r/2}(z)$ say -- boundary conditions of the GFF make this scaling relations only approximative and hence introduce a few tedious complications.
%Also we point out that the fact the exponent in the right hand side is not linear in $q$ is an indication that the geometry has an interesting \emph{multifractal structure}: by definition this means that the geometry is not defined by a single fractal exponent (compare for instance with $\E(|B_t|^q)$ which is proportional to $t^{q/2}$).
\ind{Minkowski dimension}
This statement uses an alternative notion of fractal dimension: Minkowski dimension rather than Hausdorff.

We will only state the case $d=2$ of this result, even though the arguments generalise easily to arbitrary dimensions.
We again use the notation $\cS_n$ for the $n$th level dyadic covering of $D$ by squares $S_i, i \in \cS_n$ of sidelength $2^{-n}$.
For $\delta>0$, the (Euclidean) {$(\delta,2^{-n})$}-Minkowski content of $A$ is defined by
$$
M_\delta(A; 2^{-n}) = \sum_{i \in \cS_n} \indic{S_i \cap A \neq \emptyset} \Leb (S_i)^{\delta},
$$
and the (Euclidean) Minkowski dimension\indN{{\bf Miscellaneous}! $\dd_{\hspace{.05cm}M}$; Minkowski dimension} (fraction) of $A$ is then
$$
\dd_{\hspace{.05cm}M} (A) = \inf \{ \delta: \limsup_{n \to \infty} M_\delta(A, 2^{-n}) < \infty\}.
$$
Note that since we used $\Leb(S_i)$ in the definition of the Minkowski content rather than the more standard sidelength $2^{-n}$ of $S_i$, the above quantity $\dd_{\hspace{.05cm}M}$ is in $[0,1]$ and is related to the more standard notion of Minkowski dimension $D_M$  through the identity $\dd_{\hspace{.05cm}M} = D_M / 2$. Finally, we define the \textbf{Minkowski scaling exponent} $$x_M = 1- \dd_{\hspace{.05cm}M}.$$

On the quantum side, we define
$$
M^\gamma_\delta(A,2^{-n}) =   \sum_{i \in \cS_n} \indic{S_i \cap A \neq \emptyset}\ \cM (S_i)^{\delta},
$$
 and the quantum expected Minkowski dimension by
 $$
 q_M = \inf \{ \delta: \limsup_{n \to \infty} \E( M^\gamma_\delta(A, 2^{-n})) < \infty\}.
 $$
The quantum Minkowski scaling exponent is then set to be
 $$
 \Delta_M = 1 - q_M.
 $$
 The KPZ relation for the Minkowski scaling exponents is then $x_M =( \gamma^2/4 )\Delta_M^2 + (1- \gamma^2/4) \Delta_M$ (formally this is the same as the relation in  \cref{T:KPZ}). Equivalently, this can be rephrased as follows.
 \begin{prop}[Expected Minkowski KPZ, \cite{Aru}] \label{P:Aru} Suppose $\bar A$ lies at a positive distance from $\partial D$
 and that $A$ is bounded. Then
 \begin{equation}\label{Mink}
 \dd _M =(1+ \gamma^2/4) q_M - \gamma^2 q_M^2/4.
 \end{equation}
\end{prop}
\begin{proof}
{First recall \cref{T:moments} from earlier in this chapter, which implies that if $0\le q\le 1$, then
$$\E(\cM(B(r))^q)\asymp r^{(d+\gamma^2/2)q-\gamma^2q^2/2}$$ for balls $B(r)$ of Euclidean radius $r$ lying strictly within the domain $D$.}

Fix $\mathrm{d} \in (0,1)$ and let $q$ be such that
$
\mathrm{d}  = (1+ \gamma^2/4) q - q^2\gamma^2/4
$
and note that $q \in (0,1)$.
%Then since $d_M + \delta > d_M$, for $n$ sufficiently large we have that
%$$
%\sum_{i \in \cS_n} \indic{S_i \cap A \neq \emptyset } \Leb(S_i)^{d_M + \delta} < 1.
%$$
%Consequently, when we go the quantum side, by the scaling relation of Proposition \ref{P:moments},
Therefore,
$$
\E( \sum_{i \in \cS_n} \indic{S_i \cap A \neq \emptyset } \cM(S_i)^{q})   \asymp \sum_{i \in \cS_n} \indic{S_i \cap A \neq \emptyset } \Leb(S_i)^{\mathrm{d}}
$$
and consequently the limsup of the left hand side is infinite if and only if the limsup of the right hand side is infinite.
In other words, $\dd_{\hspace{.05cm}M}$ and $q_M$ satisfy \eqref{Mink}. %(Note here, we are ignoring boundary effects).
\end{proof}

\subsubsection{Duplantier--Sheffield's KPZ theorem}
\label{S:DS}

We end this chapter with a short description of Duplantier and Sheffield's definitions of scaling exponents, as well as a sketch
of proof of the resulting KPZ formula \cite{DuplantierSheffield}. (The statement is a bit weaker than \cref{T:KPZas}, since
the notions of scaling exponents are slightly harder to use, and the formula holds only in expectation as opposed to almost surely).

In this section, the (Euclidean) scaling exponent of $A
\subset D$ is the limit, if it exists, defined by
$$x' = \lim_{\eps \to 0}\frac{\log \P( B(z, \eps) \cap A  \neq\emptyset) }{\log (\eps^2)},
%=  \lim_{\eps \to 0}\frac{\log \Leb (A_\eps) }{\log (\eps^2)}
$$
where $\P$ is the joint law of $A$ (if it is random) and a point $z$ chosen proportionally to Lebesgue measure in $D$. We will assume that $D$ is bounded.
%{Is $\P$ here the joint probability measure on $A$ and $z$ - we should define it. If $A$ is not a subset of $\D$, how do we define the probability measure on $z$?}

We need to make a few comments about this definition.

\begin{enumerate}

\item First, this is equivalent to saying that the volume of $A_\eps$, the Euclidean $\eps$-neighbourhood of $A$, decays like $\eps^{2x'}$. In other words, $A$ can be covered with $\eps^{-(2 - 2x')}$ balls of radius $\eps$, and hence typically the Hausdorff dimension of $A$ is simply
$$
\dd_H(A) = 2 - 2x' = 2(1-x'),
$$
consistent with our earlier definition of Euclidean scaling exponent. In particular, note that $x'\in [0,1]$ always; $x'=0$
means that $A$ is practically the full space, $x'=1$ means it is practically empty.

\item In the definition we divide by $\log (\eps^2)$, because $\eps^2$ is the volume (with respect to the Euclidean geometry on $\R^2$) of a ball of radius $\eps$. In the quantum world, we would need to replace this by the Liouville area of a ball of radius $\eps$ -- see below.

\end{enumerate}

%We therefore make the following definition -- to be taken with a pinch of salt, as we shall allow ourselves to consider slightly different variants. Since we do not have a metric to speak about, we resort to the following expedient:
\ind{Isothermal ball}
The quantum analogue of this is the following.
For $z \in D$, we denote by $B^\delta(z)$ the quantum ball of mass $\delta$: that is, the Euclidean ball centred at $z$ whose radius is chosen so that its Liouville area is precisely $\delta$. (In \cite{DuplantierSheffield}, this is called the \emph{isothermal} ball of mass $\delta$ at $z$).
The quantum scaling exponent of $A\subset D$ is then the limit, if it exists, defined by
$$\Delta' = \lim_{\delta \to 0}\frac{\log  \P( B^\delta (z) \cap A \neq \emptyset )}{\log ( \delta)},$$
where $z$ is sampled from the Liouville measure $\cM$ normalised to be a probability distribution.
%where $A^\delta $
%is the quantum $\delta$-enlargement of $A$; that is, $A^\delta = \{z \in \D: B^\delta(z) \cap A \neq \emptyset\}$.
%\end{definition}

\begin{theorem}[{Expected} Hausdorff KPZ formula] \label{T:KPZ}Suppose $A$ is independent of the GFF, $\gamma\in (0,2)$, and $D$ is bounded. Then if $A$ has Euclidean scaling exponent $x'$,
 it has quantum scaling exponent $\Delta'$, where $x'$ and $\Delta'$ are related by the formula
\begin{equation}
x' = \frac{\gamma^2}{4} (\Delta')^2 + (1- \frac{\gamma^2}4) \Delta'.
\end{equation}
\end{theorem}

We will now sketch the argument used by Duplantier and Sheffield to prove \cref{T:KPZ}, since it is interesting in its own right and gives a somewhat different perspective
(in particular, it shows that the KPZ formula can be seen as a large deviation probability for Brownian motion).

\medskip \noindent\textbf{Informal description of the idea of the proof.}
We wish to evaluate the probability $\P( B^\delta(z) \cap A \neq \emptyset)$, where $z$ is a point
sampled from the Liouville measure, and $B^\delta$ is the Euclidean ball of Liouville mass $\delta$ around $z$. Of course the event that this ball intersects $A$ is rather unlikely, since the ball is small. But it can happen for two reasons. The first one is simply that $z$ lands very close (in the Euclidean sense) to $A$ -- this has a cost governed by the Euclidean scaling exponent of $A$, by definition, since we may think of $z$ as being sampled from the Lebesgue measure and then sampling the Gaussian free field given $z$, as in the description of the rooted measure \cref{sec:rooted_meas}. However, it is more economical for $z$ to land relatively further away from $A$, and instead require that the ball of quantum mass $\delta$ have a bigger than usual radius. As the quantum mass of the ball of radius $r$ around $z$ is essentially governed by the size of the circle average $h_r(z)$, which behaves like a Brownian motion plus some drift, we find ourselves computing a large deviation probability for a Brownian motion.
The KPZ formula is hence nothing else but the large deviation function for Brownian motion.

\begin{proof}[Sketch of proof of \cref{T:KPZ}] Now we turn the informal idea above into more concrete mathematics, except for two approximations that we will not justify. Suppose $z$ is sampled according to the Liouville measure $\cM$. Then we know from \cref{T:Girsanov} (see also \eqref{condrootG}) that the
{joint law of the point $z$ and the free field} is absolutely continuous with respect to {a point $z$ sampled from Lebesgue measure, together with the field $h^0(\cdot) + \gamma \log| \cdot - z| + O(1)$, where $h^0$ is a GFF that is independent of $z$}. (See \cref{sec:rooted_meas}). Hence the mass of the ball of radius $\eps$ about $z$ is approximately given by
\begin{align}\label{approx}
\cM(B(z,\eps)) & \approx  \eps^{\gamma^2/2} e^{\gamma h_\eps(z)}\times  \eps^2 \nonumber \\
 & \asymp e^{\gamma ( h^0_\eps(z)  + \gamma \log 1/\eps)} \eps^{2+ \gamma^2/2} \nonumber  \\
 & = \eps^{2 - \gamma^2/2} e^{\gamma h^0_\eps(z)} .
\end{align}
%(Note that $2 + \gamma^2 /2 = \gamma Q$, where $Q = \gamma /2 + 2 /\gamma$.) Call $r_\delta$ the first radius such that the above is less than $\delta$.
It takes some time to justify rigorously the approximation in \eqref{approx}, but the idea is that the field $h_\eps$
fluctuates on a spatial scale of size roughly $\eps$. Hence we are not making a big error by pretending that $h_\eps$ is constant on $B(z,\eps)$,
equal to $h_\eps(z)$. In a way, making this precise is the most technical part of the paper \cite{DuplantierSheffield}.
%{Also the fact that $h$ is only absolutely continuous w.r.t. that law is a problem?}
 We will not go through the arguments which do so, and instead we will see how, assuming it, one is naturally led to the KPZ relation.

Working on an exponential scale (which is more natural for circle averages) and writing $B_t = h^0_{e^{-t}}(z)$, we find that
$$
\log \cM(B(z,\e^{-t})) \approx \gamma B_t - (2 - \gamma^2/2) t .
$$
We are interested in the maximum radius $\eps$ such that $\cM(B(z,\eps)) $ will be approximately $\delta$: this will give us the Euclidean radius of the quantum ball of mass $\delta$ around $z$. So let
\begin{align*}
T_\delta &= \inf \{ t \ge 0: \gamma B_t - (2 - \gamma^2/2) t \le \log \delta \}\\
& = \inf \{ t \ge 0: B_t + (\frac{2}{\gamma} - \frac{\gamma}{2}) t \ge \frac{\log (1/\delta)}{\gamma} \}.
\end{align*}
where the second equality is in distribution. Note that since $\gamma < 2$ the drift is positive, and hence $T_\delta < \infty$ almost surely.

Now, recall that if $\eps>0$ is fixed, the probability that $z$ will fall within (Euclidean) distance $\eps$ of $A$ is
approximately $\eps^{2x'}$. Hence, applying this with $\eps  = e^{- T_\delta}$ the probability that $B^\delta(z)$ intersects $A$ is, approximately, given by
$$
\P ( B^\delta(z) \cap A \neq \emptyset) \approx \E( \exp( -2x' T_\delta)).
$$
This is the second approximation that we will not seek to justify fully. Consequently, we deduce that
$$
\Delta' = \lim_{\delta \to 0} \frac{\log \E( \exp( -2x' T_\delta))}{\log \delta}.
$$
For $\beta >0$, consider the martingale
$$
M_t = \exp ( \beta B_t - \beta^2 t / 2),
$$
and apply the optional stopping theorem at the time $T_\delta$ (note that this is justified). Then we get, letting $a = 2/\gamma - \gamma /2 $, that
$$
1 = \exp ( \beta  \frac{\log (1/\delta)}{\gamma} ) \E( \exp ( - (a\beta + \beta^2 / 2) T_\delta )).
$$
Finally set $2x' = a \beta + \beta^2/ 2$, so that $\E( \exp ( -2x' T_A)) = \delta^{\beta /\gamma}$. In other words,
$\Delta' = \beta/ \gamma$,
where $\beta$ is such that $2x' = a \beta + \beta^2/ 2$. Equivalently, $\beta = \gamma \Delta'$, and
$$
2x' = (\frac2{\gamma} - \frac{\gamma}2) \gamma \Delta' + \frac{\gamma^2}{2} (\Delta')^2.
$$
This is exactly the KPZ relation.
\end{proof}

\subsubsection{Applications of KPZ to critical exponents}
\label{S:exponents}
\ind{Critical exponents}

\dis{This section explains in a non-rigorous manner how the KPZ relation can be used to compute critical exponents in some models of statistical mechanics in two dimensions. This section can be skipped on a first reading, and is only relevant later in connection with the end of  \cref{S:maps}. This section also assumes basic familiarity with the notion of random planar maps and the conjectures %described in  \cref{S:maps}
	related to their conformal embeddings, see \cref{sec:RPM_LQG}.}

At the discrete level, the KPZ formula can be interpreted as follows.
Consider a random planar map $M$ of size $N$ (where `size' refers indistinctly to the number of faces, vertices or edges). Suppose that a certain subset $A$ within $M$ has a size $|A| \approx N^{1-\Delta}$, so that $A$ is ``fractal-like". We have in mind a set $A$ which is defined conditionally independently given the map, and of course depends on $N$ (but we do not indicate this in the notation). For instance, $A$ could be the set of double points of a random walk on the map run until its cover time, or the set of pivotal edges for percolation on the map with respect to some macroscopic event.
% {Can we add some kind of explanation of what sort of object $A$ should be? It should be something that is naturally defined in a random planar map or a Euclidean lattice of any size? And that should converge along with the planar map to something that is independent of the limiting LQG measure?}
We may also consider the Euclidean analogue $A'$ of $A$ within a Euclidean box of area $N$  (and thus of side length $n = \sqrt{N}$). Namely, $A'$ is the set that one obtains when the map $M$ is exactly this subset of the square lattice.
In this case we again expect $A' $ to be fractal-like, and so $|A'| \approx N^{1-x} = n^{2 - 2x}$ for some $x\in [0,1]$. If $A'$ has a scaling limit then this $x$ is nothing but its Euclidean scaling exponent (indeed, the discrete size of $A'$ is essentially the number of balls of a fixed radius required to cover a scaled version of it). Likewise, if $A$ has a scaling limit then $\Delta$ is nothing but its quantum scaling exponent.

Hence the KPZ relation suggests that
$x, \Delta$ should be related by
$$
x = \frac{\gamma^2}{4} \Delta^2 + (1 - \frac{\gamma^2}{4}) \Delta
$$
in the limit as $N \to \infty$. Here $\gamma$ refers to the universality class
of the map; this assumes that $A$ is (when embedded suitably in the plane) independent of the field $h$ which represents the embedding of the map in the limit.

In particular, observe that the approximate (Euclidean) Hausdorff dimension of $A'$ is then $ 2-2x$, consistent with our definitions. See \cref{S:maps} for concrete examples, where this is used, for instance, to guess the loop-erased random walk exponent.

\subsection{Exercises}

\begin{enumerate}[label=\thesection.\arabic*]

\item By considering the set of thick points or otherwise, argue that the KPZ relation does not need to hold if the set $A$ is allowed to depend on the free field; for instance show that we can have $\Delta'=0$, while $x'> 0$. 
This type of example was first considered by \cite{Aru} who also considered the case of flow lines associated to the GFF.

%\item {Let $h$ be a centred Gaussian field with covariance $K$ as in \cref{setup_gmc}. Show that the Gaussian multiplicative chaos measures constructed in \cref{T:conv} satisfy \eqref{eq:shamov2}. By considering the analogue of the measure $Q$ from \cref{sec:rooted_meas}, and using the disintegration theorem, show that there is only one coupling $(h,\cM_\gamma)$ such that: \begin{itemize}	
%		\item $\cM_\gamma$ is measurable with respect to $h$ as a stochastic process;
%		%\item  $\E(\cM_\gamma(S))=\sigma(S)$ for all $S\subset D$; and
%		\item \eqref{eq:shamov2} holds.
%	\end{itemize}  }

%\nb{Would like to know how this works...!}

\item Suppose $K(x,y) \ge 0$ for all $x,y \in D$. Let $A \subset D$, and let $q \in [0,1]$. Show that $\E( \cM(A)^q)$ is a non-increasing function of $\gamma \in {[0,\sqrt{2\mathbf{d}})}$. (Hint: use Kahane's inequality).

\item {Let $A\subset D$. For $\gamma<{\sqrt{\mathbf{d}}}$, show that $\cM(A)$ admits a continuous modification in $\gamma$. (Hint: use the Kolmogorov continuity criterion.) }

\item \label{Ex:atom} (a) Use the scaling invariance properties developed in the proof of the multifractal spectrum to show that $\cM$ almost surely has  no atoms. %(You should admit that $\E( \cM(S)^q) < \infty$ for some $q>1$.)

(b) Give an alternative proof, using the energy estimate in Exercise \ref{Ex:thick} of \cref{S:Liouvillemeasure}.
\ind{Non atomicity of Liouville measure}

\item This exercise gives a flavour of Kahane's original pioneering argument for the construction of GMC in \cite{Kahane85}. Suppose that $K$ is a  covariance kernel of the form \eqref{setup_gmc}, that can be written in the form  $$K(x,y)=\sum_{n=1}^\infty K_n(x,y)$$ for all $x\ne y$ in $D\subset \R^d$, where for each $n$, $K_n:D\times D\to \R$ is positive definite and satisfies $K_n(x,y)\ge 0$ for all $x,y\in D$. Such a covariance kernel was called \textbf{$\sigma$-positive}  by Kahane. Show that there exists a sequence of centred Gaussian fields $(h^n)_{n\ge 1}$ such that the fields $(h^n-h^{n-1})_{n\ge 1}$ are independent centred Gaussian fields with covariances $K_n$ for each $n$. Let $\sigma$ be a reference Radon measure satisfying \eqref{dim} for some $\mathbf{d}>0$. For $0 \le \gamma< \sqrt{2\mathbf{d}}$, we use this decomposition to construct a natural sequence of `chaos approximations' $\cM_n$ by setting
    $$
    \cM_n(A) = \int_A \exp\{ \gamma h_n(x) - \tfrac{\gamma^2}{2} \E( h_n(x)^2)\} \sigma(\dd x),
    $$
    for any Borel set $A$.
    Prove that $\cM_n(A)$ has an almost sure limit $\cM(A)$ as $n\to \infty$ which defines a random measure.
    %\nb{In this exercise we will admit that $\E( \cM(A)) =|A|$.}
	
	Suppose we are given two $\sigma$-positive decompositions for $K$, say $$K(x,y) = \sum_{n=1}^\infty K_n(x,y) = \sum_{n=1}^\infty K'_n(x,y),$$ and let $\cM$ and $\cM' $ be the associated chaos measures constructed above. Using Kahane's inequality (and without using Theorem \ref{T:conv}), show that for any Borel set $A$, $\E( \cM(A)^q) \le \E (\cM'(A)^q)$ for $q \in (0,1)$ (note that this argument does not require knowing that either $\cM$ or $\cM'$ is non-zero). Deduce that the laws of  $\cM$ and $\cM'$ (as random measures) are identical. This is Kahane's theorem on uniqueness of GMC; Kahane's inequality \cite{Kahane86} was discovered for the purpose of this proof.

\item {We now take the same setup as above, and assume the result of \cref{T:conv}. Show that in the case $\gamma<\sqrt{2\mathbf{d}}$, the limit $\cM$ constructed above agrees with the GMC measure of Theorem \ref{T:conv}.}

%\item  \nb{Also, I think the nice part is to show the independence of the (law of the) limit using Kahane's inequality: By taking more summands, we get a kernel $K' \ge K$ so all fractional moments must match! }{*To add}.

%{\emph{\textbf{Remark.} These were the types of field that Kahane originally considered in \cite{Kahane85}. In this paper, he showed that the limit measure is non-trivial iff $\gamma< \sqrt{2d}$ and that it is independent of the choice of decomposition $K=\sum_n K_n$. }}

\item \label {E:log}{If $K$ is as in \eqref{setup_gmc}, define the linear operator $T$ on $L^2(D)$ by setting $$Tf(x)=\int_D K(x,y) f(y) \, \dd y$$ for each $f\in L^2(D)$. When $D$ is bounded, one can show using standard operator theory that there exists an orthonormal basis $\{f_k\}_{k\ge 0}$ of $L^2(D)$, made up of eigenfunctions for $T$. The ordering can be chosen so that the associated eigenvalues $\{\lambda_k^{-1}\}_{k\ge 0}$ satisfy  $0\le \lambda_1\le \lambda_2\le \lambda_3 , \ldots$. 
%Moreover, with this ordering, the first eigenfunction $f_1$ will be strictly positive in $D$. 
}

\begin{itemize}
		\item[(a)] { Show that for each $x$ in $D$, $$\sum_{k =0}^n \lambda_k^{-1} f_k(x) f_k(\cdot) \to K(x,\cdot) \text{ in } L^2(D)$$ as $n\to \infty$. Let $h$ be the centred Gaussian field with covariance $K$. By considering the joint law of $\{\lambda_k^{1/2}(h,f_k)\}_{k\ge 1}$, show that for any smooth compactly supported test function $f$ on $D$, if $h^n:=\sum_{k=0}^n (h,f_k) f_k$, then $(h^n,f)$
converges almost surely and in $L^2(\P)$  to $(h,f)$ as $n\to \infty$.}

{\emph{\textbf{Remark.} This decomposition of $h$ is known as the \textbf{Karhuhen--Loeve expansion}.}}

\item[(b)] {Show further that for $\gamma\ge 0$, the sequence of measures defined by
$$\cM^n(S):= \int_S\exp(\gamma h^n(z) - \gamma^2/2 \var(h^n(z))\, \dd z \quad \quad S\subset D, n\ge 0$$
has an almost sure limit with respect to the topology of weak convergence of measures. When $\gamma<\sqrt{2\mathbf{d}}$, show that $\cM^n(S)$ is a uniformly integrable family for any fixed $S$. Use this to show that $\lim_n \cM^n(S)$ agrees almost surely with $\cM_\gamma(S)$, where $\cM_\gamma$ is the GMC measure for $h$ constructed in \cref{T:conv}.}
\item[(c)] Suppose that $f_1$ is non-negative and bounded. {Use \eqref{laplace_useful} to show that for $\delta>0$, $\P(\cM_\gamma(D)\le \delta)\ge c\P(Z\le \delta)$ where $Z$ is an appropriately chosen lognormal random variable and $c>0$ does not depend on $\delta$.}

\end{itemize}

\end{enumerate}

\newpage
\section{Statistical physics on random planar maps}
\label{S:maps}

\subsection{Fortuin--Kasteleyn weighted random planar maps}

\def\m{\mathbf{m}}
\def\t{\mathbf{t}}
\def \bm {\boldsymbol{m}}
\def \bt {\boldsymbol{t}}
\def \aH {\textsf{H}}
\def \ah {\textsf{h}}
\def \aC {\textsf{C}}
\def \ac {\textsf{c}}
\def \aF {\textsf{F}}
\def \ae {\mathsf{e}}
\def \len {\textsf{Len}}
\def \Area {\textsf{Area}}

In this chapter we change our focus from the continuous to the discrete, and describe the model of random planar maps weighted by self dual Fortuin--Kasteleyn percolation. As we will see, these maps can be thought of as canonical discretisations of Liouville quantum gravity (but there are in fact many other models of planar maps which are believed to be related to Liouville quantum gravity).

We proceed as follows.  We first recall the notion of planar map and \textbf{decorated planar map} before defining a probability measure on such maps (maps decorated by self dual FK loops).
In \cref{sec:RPM_LQG}, we discuss aspects of the conjectured connection between this model of planar maps and \textbf{Liouville quantum gravity}. In \cref{sec:bij_tree} we focus on the case where the decoration is a spanning tree. Here we describe in detail a powerful \textbf{bijection} due (independently) to Mullin, Bernardi and Sheffield,  between tree decorated maps and pairs of independent, positive random walk excursions (equivalently, two dimensional random walk excursions in the positive quadrant). This bijection is a convenient way to approach the question of scaling limits. We use it in \cref{sec:LERWexp} to compute the (quantum) scaling exponent of the \textbf{loop-erased random walk} (LERW). Using the KPZ relation of \cref{S:KPZ}, we find that it agrees with various known properties of LERW on the square lattice, including the Hausdorff dimension of its scaling limit SLE$_2$. In \cref{S:Sheffield_bij}, we discuss Sheffield's bijection, which is a generalisation of the aforementioned bijection to decorations which are no longer spanning trees but \textbf{densely packed loop configurations}. Again, this bijection is from decorated maps to pairs of excursions in a suitable sense. In this case, however, the excursions are far from independent; this has an interpretation in terms of a \textbf{discrete mating of trees} which will be described in the continuum in Chapter \ref{S:MOT}. This description is used in \cref{sec:local} to show the existence of an infinite volume local limit. A scaling limit result is discussed which, roughly speaking, shows that the limiting trees are correlated infinite CRTs.

\ind{Planar maps!Definition}
\ind{Planar maps!Dual map}
\ind{Planar maps! Decorated}
\paragraph{Planar map, dual map.} Recall that a \textbf{planar map} is a proper embedding of a
(multi) graph with a finite number of edges in the plane $\C \cup \{\infty\}$ (viewed as the Riemann sphere), which is viewed up to
orientation preserving homeomorphisms from the sphere to itself.   Let
$\bm_n$ be a map with $n$ edges and $\bt_n$ be {a subgraph spanning \emph{all} of its vertices}. We call the pair
$(\bm_n, \bt_n)$ a (subgraph) \textbf{decorated map}. Let $\bm_n^\dagger$\indN{{\bf Planar maps}! $m^\dagger$; dual map of a map $m$} denote
the \emph{dual map} of
$\bm_n$. Recall that the vertices of the dual map  correspond to the
faces of $\bm_n$ and two vertices in the dual map are adjacent if and
only if their corresponding faces are adjacent to a common edge in the
primal map. Every edge $e$ in the primal map corresponds to an edge $e^\dagger$\indN{{\bf Planar maps}! $e^\dagger$; dual edge of an edge $e$} in the
dual map which joins the vertices corresponding to the two faces
adjacent to $e$. The
dual subgraph $\bt_n^\dagger$ is the graph formed by the subset of edges $\{e^\dagger:
e\notin \bt_n\}$ and all dual vertices. We fix an edge in the map $\bm_n$, {to which we also assign an orientation,} and define it to
be the
root edge. With an abuse of notation, we will still write $\bm_n$ for the rooted map; and we let $\cM_n$ be the set of maps with $n$ edges together with one distinguished edge called the root.\indN{{\bf Planar maps}! $\cM_n$; maps with $n$ edges and one distinguished root edge}

\ind{Planar maps!Loops}
\ind{Planar maps!Refinement edges}
\paragraph{Canonical triangulation.} Given a subgraph decorated map $(\bm_n,\bt_n)$ with $\bm_n\in \cM_n$ and a subgraph $\bt_n$ of $\bm_n$, one can associate to it a set of
loops where in some sense each loop forms the interface between two clusters (connected components) associated to $\bt_n$ and its planar dual. To be more precise, let us say that two vertices $x$ and $y$ of $\bm_n$ are in the same component if we can travel from $x$ to $y$ using edges in $\bt_n$; by convention a vertex $x$ is always in its own component (hence that component will consist only of $x$ if $x$ is not covered by $\bt_n$). We can use the same definition to talk about clusters on the planar map $\bm_n^\dagger$ dual to $\bm_n$ and the dual configuration of edges $\bt_n^\dagger$; the loops then separate primal and dual clusters. To define these loops more precisely, we will need to discuss not only the dual planar map but also a couple of related maps that can be constructed from superposing the primal and dual maps.

We first consider an auxiliary map which we call the \textbf{Tutte map}, and which is formed
by joining the dual vertices in every face of $\bm_n$ with the primal vertices
incident to that face. We call these edges \textbf{refinement
  edges} (drawn in green in Figure \ref{fig:clusters}). Thus the vertex set of the Tutte map consists of all primal and dual vertices, but note that its edge set does not contain any of the original edges of $\bm_n$ or its dual. It is easy to check that this Tutte map is a quadrangulation, meaning each face has exactly four (refinement) edges surrounding it. Each of the original edges of $\bm_n$ or $\bm_n^\dagger$ is a diagonal of one of these quadrangles. In other words, every
edge in $\bm_n$ corresponds to a quadrangle in the Tutte map; this quadrangle can be viewed as the union of two triangles, one on either side of the edge.
% \ellen{(so the refinement doesn't include the edges of $\bm_n$ itself?)} formed by the union of the two triangles
%incident to its two sides.

In fact this construction defines a
bijection between maps with $n$ edges and quadrangulations with $n$
faces, sometimes called the \textbf{Tutte bijection}.
\ind{Tutte bijection}
\begin{figure}
\centering{
\textbf{a.} \includegraphics[scale=.4]{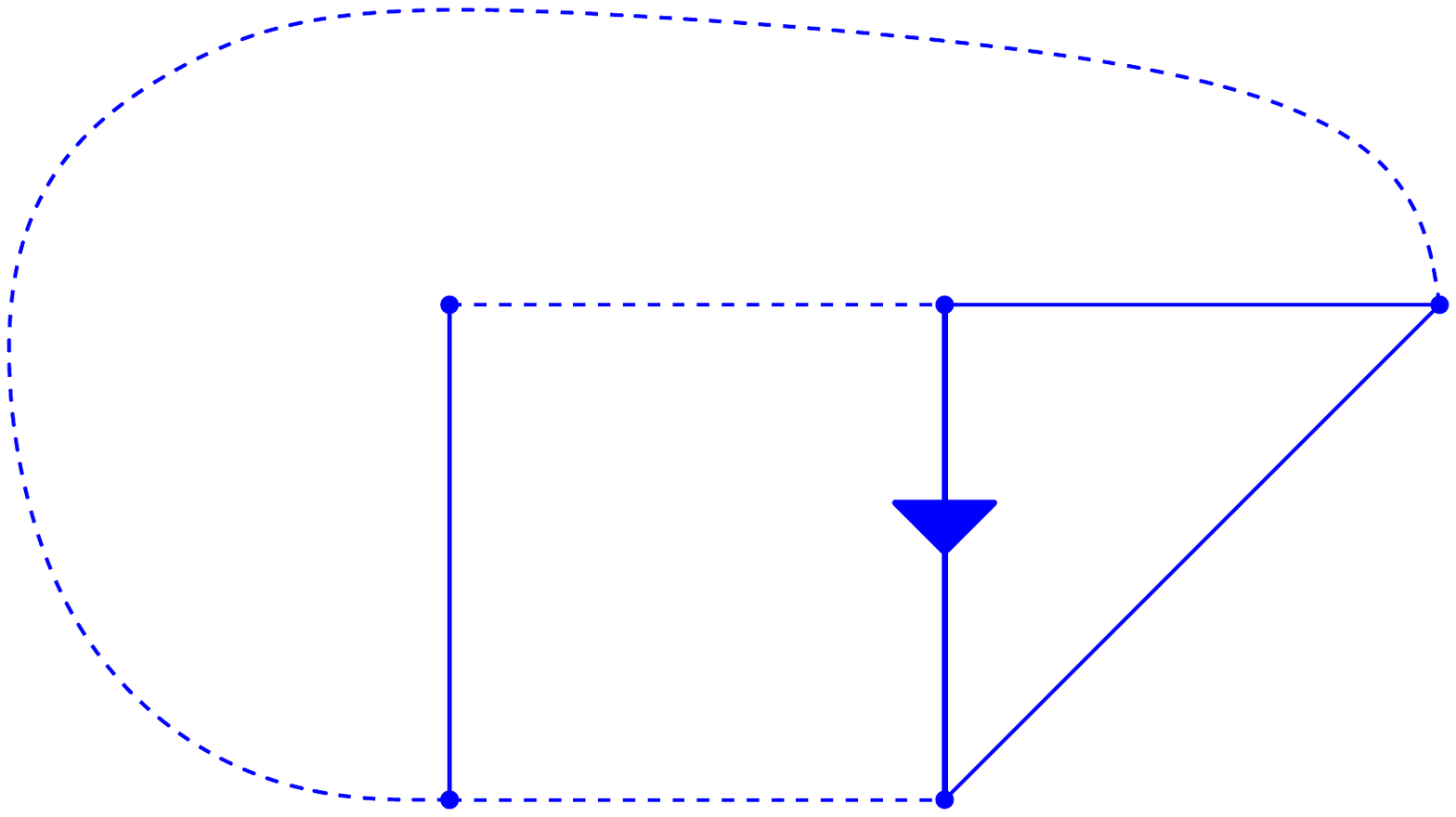} \quad \quad
\textbf{b.} \includegraphics[scale=.4]{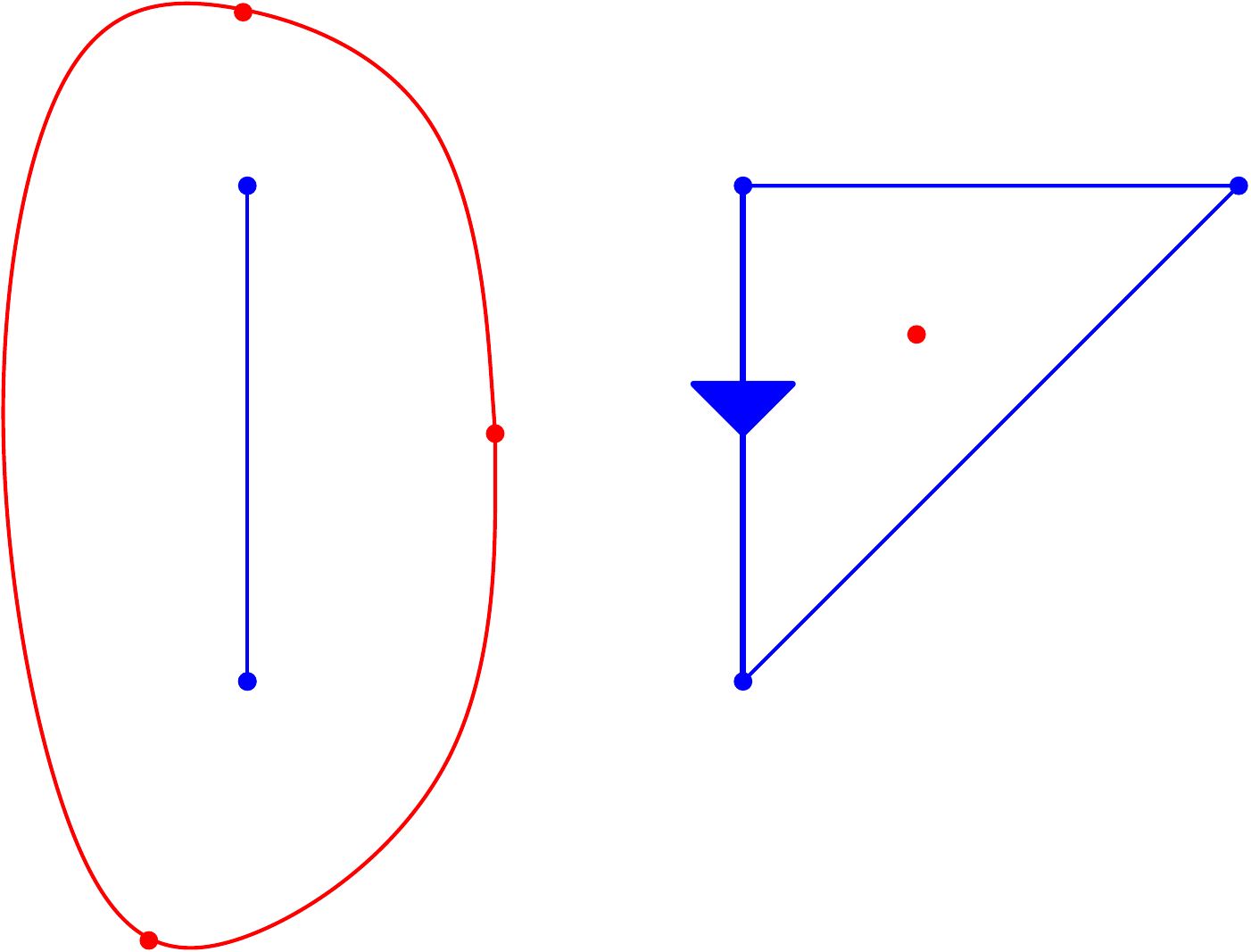}
}

\vspace{.5cm}
\centering{
\textbf{c.} \includegraphics[scale=.4]{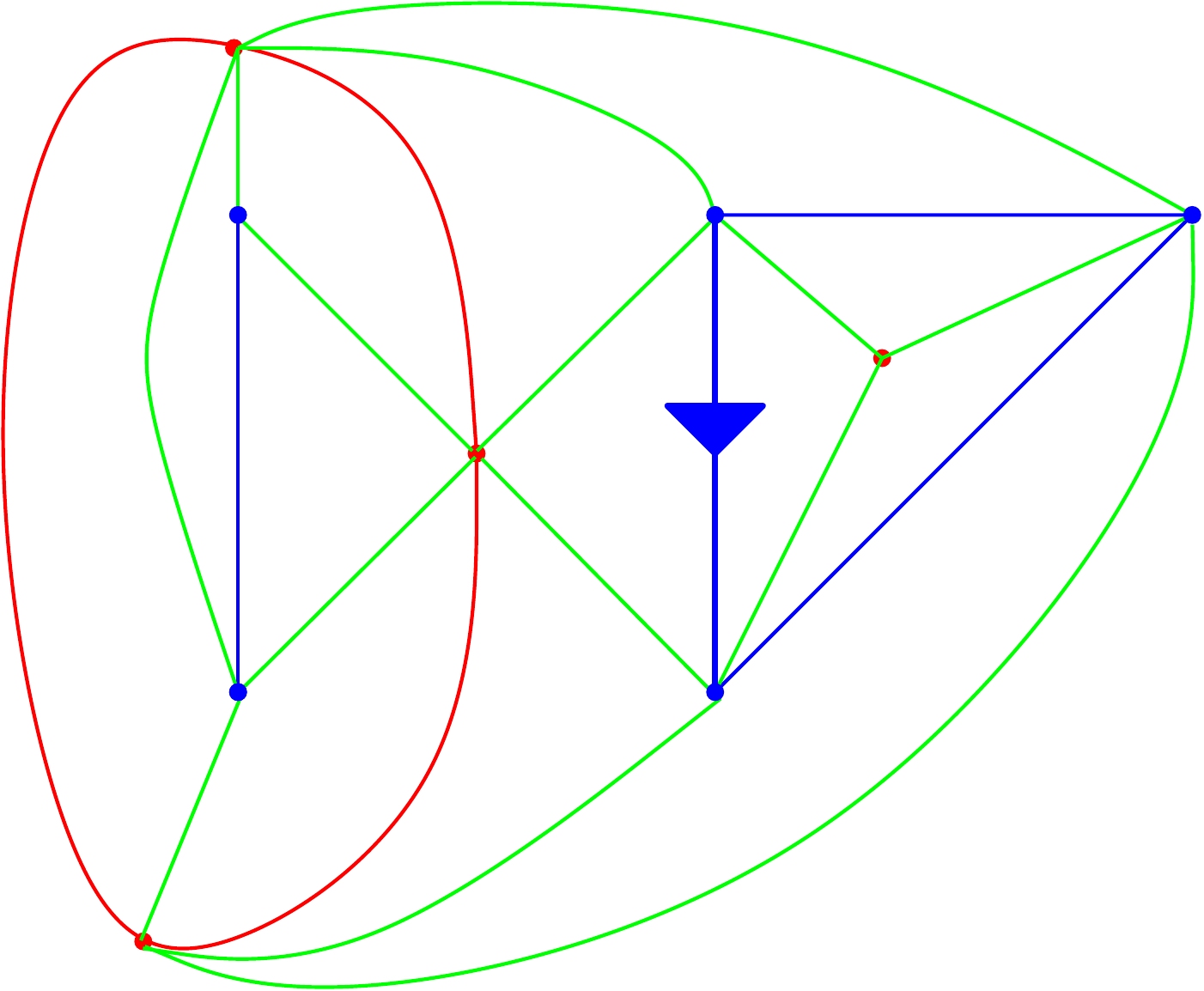}\quad \quad
\textbf{d.} \includegraphics[scale=.4]{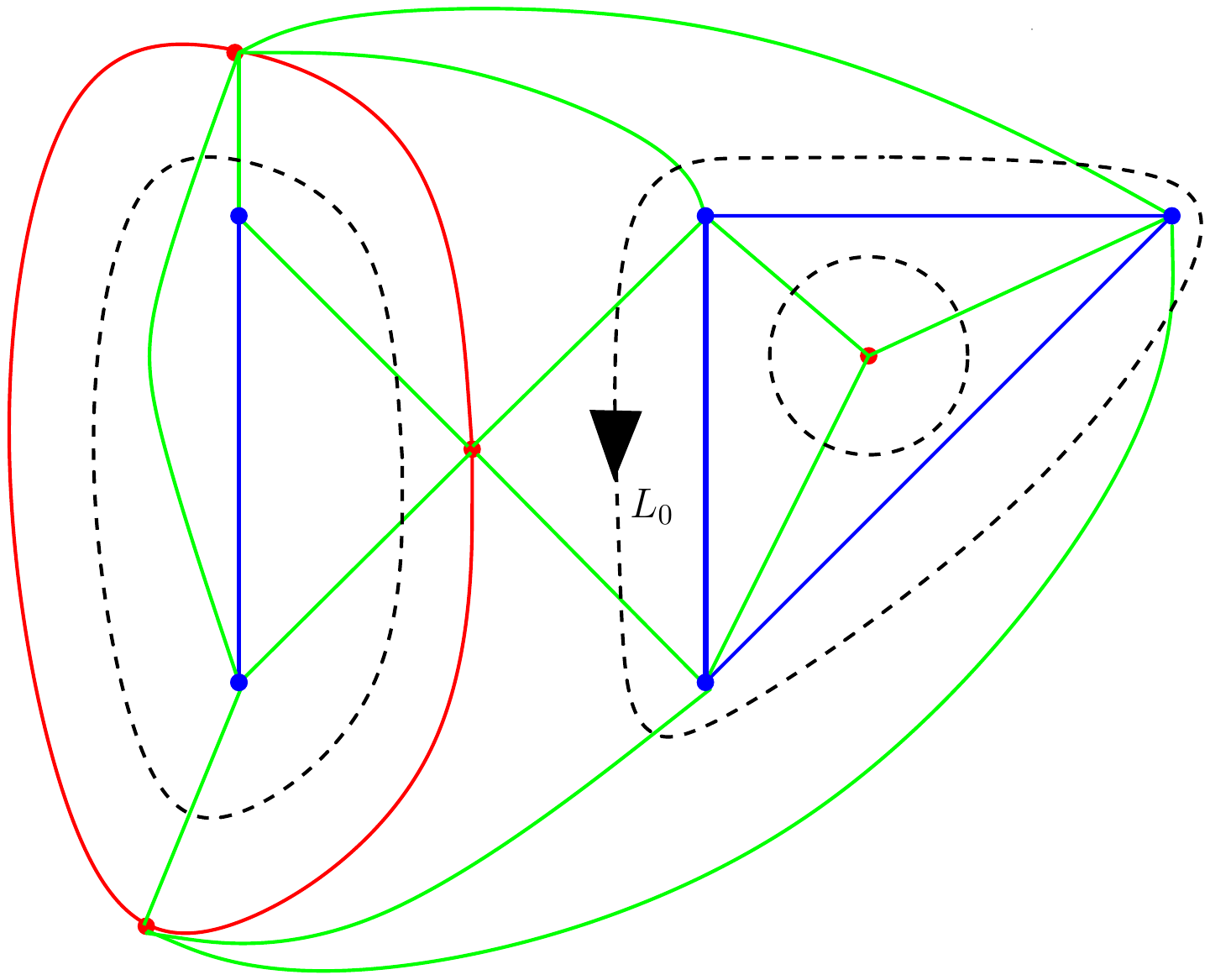}
}
\caption{A map $\bm$ decorated with loops associated to a set of open edges $\bt$.
{\textbf{a.}} The map is
  in blue, with solid open edges and dashed closed edges. {\textbf{b.}} Open clusters and corresponding open dual clusters are shown in blue and red. {\textbf{c.}} Every dual vertex is joined to its adjacent primal vertices by a green edge. This results in a refined map $\bar \bm$ which is a triangulation. {\textbf{d.}} The primal and dual open clusters are separated by loops, which are drawn in black and are dashed. Each loop is identified with the set of triangles through which it passes: note that it crosses each triangle in the set exactly once.
  {The oriented root edge of the map is indicated with a blue arrow in subfigures \textbf{a}, \textbf{b} and \textbf{c}. The loop $L_0$ is marked with an arrow in subfigure \textbf{d}, and the arrow indicates the orientation of the loop, parallel to the orientation of the root edge.} %The roots of these maps (edge or triangle) are indicated with arrows.
%
 %
 % with full lines for closed edges and dotted lines
 % for open edges. The dual map is in red, with only the open dual dges
 % drawn. The green edges form the associated quadrangulation along
 % with the primal and dual edges form the ($\bar m $ in \cref{sec:critical-fk-model}). The percolation clusters and dual clusters are separated by loops, in dashed purple, that cross every quandrangulation edge once.
 }
\label{fig:clusters}
\end{figure}

Given a subgraph decorated map $(\bm_n,\bt_n)$ define the \textbf{refinement} map $\bar
\bm_n$\indN{{\bf Planar maps}! $\bar{m}$; refinement map of a map $m$} to be formed by the
union of $\bt_n,\bt_n^\dagger$ and the
refinement edges; note that its vertex set is the same as the Tutte map, that is, every primal and dual vertex of $\bm_n$. The addition of $\bt_n$ and $\bt_n^\dagger$ makes the refinement map a triangulation: indeed, every quadrangle from the Tutte map has been split into two (either with a diagonal from $\bt_n$ or from $\bt_n^\dagger$).
%\ellen{Does this define $\bm_n$? What is $\bt_n$?)}
The root edge of $\bm_n$ induces a \textbf{root
  triangle} on the refinement map, which is taken to be the triangle immediately to the right of the root edge of $\bm_n$.
  %It is easy to see that such a map is a
%triangulation: every face in the refinement of $\bm_n$ is divided into
%two triangles either by a primal edge in $\bt_n$ or a dual edge in
%$\bt_n^\dagger$.

Note that every triangle consists of two refinement edges and one edge from either $\bt_n$ (primal edge) or $\bt^\dagger_n$ (dual edge).
For future reference, we call such a triangle in $\bar \bm_n$ a \textbf{primal triangle} or \textbf{dual triangle} respectively (see \cref{F:triangle}).

\paragraph{Loops.} Finally, given $(\bm_n, \bt_n)$ we can define the loops induced by $\bt_n$ as follows. For each connected component $C$ of either $\bt_n$ or $\bt_n^\dagger$, we draw a loop surrounding it (meaning a closed curve in the complement of $C$ in the sphere; the complement contains two components, and by convention we draw it in the ``exterior'' one that contains the point on the sphere designated to be $\infty$; note that even in the case where the connected component $C$ is reduced to a single vertex there is still a loop surrounding it which separates the sphere in two components). If this loop is drawn sufficiently close to $C$ it identifies a unique collection of triangles that are adjacent to $C$ (in the sense that they share at least a vertex with it). We view the component $C$ itself as an open cluster for a percolation configuration either on $\bm_n$ or its dual, and will use the word ``cluster'' interchangeably from now on. %\nb{There are two components in the complement, neither of which is distinguished on the sphere... should we specify that this is the ``exterior" one wrt a point chosen to be $\infty$?}

In what follows, one should visualise the loop of $C$ as being a closed curve drawn sufficiently close to $C$ in its complement, as above. However for precision, we will actually identify the loop with the collection of triangles through which it passes. See Figure \ref{fig:clusters} for an illustration. In this way, each loop is simply a collection of triangles ``separating'' a primal connected component of $\bt_n$ from a dual connected component in $\bt_n^\dagger$, or vice versa. Note that the set of loops is ``space
filling'' in the sense that every triangle of the refined map is contained in a loop. {We denote by $L_0$ the loop that is associated with the root triangle. It comes with a natural orientation, which is parallel to the orientation of the root edge of $M_n$.}

Also, given the Tutte map and the collection of closed curves described above,
one can recover the spanning subgraph $\bt_n$ (hence also $\bt^\dagger_n$) that generates it. Let $\ell(\bm_n,\bt_n)$ denote the number of loops
corresponding to a configuration $(\bm_n,\bt_n)$. Note that this is equal to the number of clusters in $\bt_n$ plus the number of clusters in $\bt_n^\dagger$ minus one; indeed, each new cluster generates a new loop.

\ind{Planar maps!Fortuin--Kasteleyn model}
\paragraph{Fortuin--Kasteleyn model.} The particular distribution on planar maps that we will now consider
was introduced in \cite{Sheffield_burger}. Let $q\ge 0$ and let $n \ge
1$: we will define a random map $M_n\in \cM_n$  decorated with a (random) subset
$T_n$  of edges. %The map $T_n$
 {As in the deterministic setting}, this induces a dual collection of edges $T_n^\dagger$
on the dual map of $M$ (see \cref{fig:clusters}). The law of {$(M_n,T_n)$ is defined by declaring that for any fixed planar map} $\bm$ with
$n$ edges, and $\bt$ a given subset of edges of $\bm$, %the probability to
%pick a particular $(\m,\t)$ is, by definition, proportional to

\begin{equation}\label{FK}
\P( M_n = \bm, T_n= \bt) \propto \sqrt{q}^\ell, \quad { \ell=\ell(\bm,\bt).}
\end{equation}
Recall from above that $\ell$ is the (total) number of loops separating primal and
dual clusters in $(\bm,\bt)$.

Equivalently, {the map $M_n$ is chosen with probability proportional to the ``partition function'' of the self dual Fortuin--Kasteleyn model on it, and} given the map $M_n$, the collection of edges $T_n$ is then sampled from this %self-dual
Fortuin--Kasteleyn model. This is in
turn closely related to the critical $q$-state Potts model, see
\cite{BKW}.
{Note that $M_n$ is actually a rooted map (as all of our maps are) and with this definition, the root edge of the map and its orientation are chosen uniformly at random (given the unrooted version).} See Figure \ref{F:maps} for simulations of $(M_n,T_n)$ at different values of $q$. \\ %Accordingly, the map $M_n$ is chosen with probability
%proportional to the partition function of the Fortuin-Kasteleyn
%model on it.

\ind{FK model}
\begin{figure}
\begin{center}
\includegraphics[width = .4\textwidth]{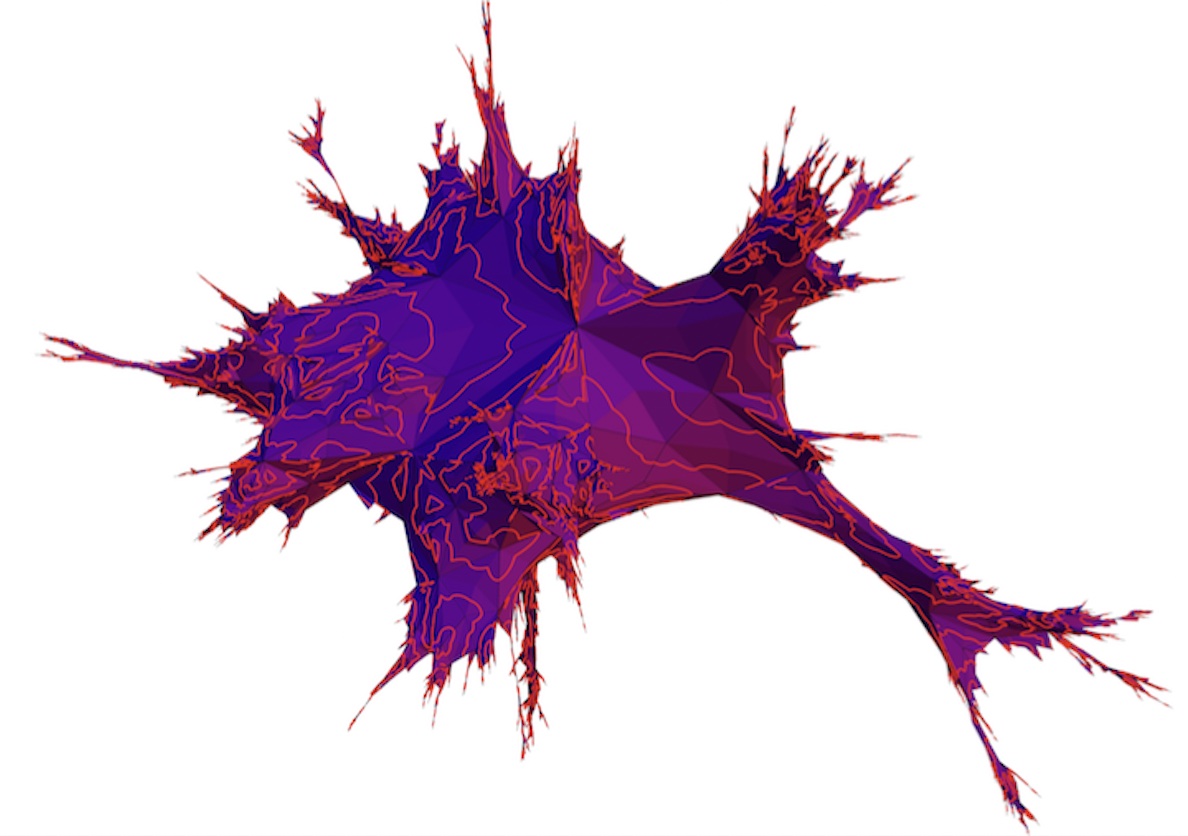}
\includegraphics[width = .4\textwidth]{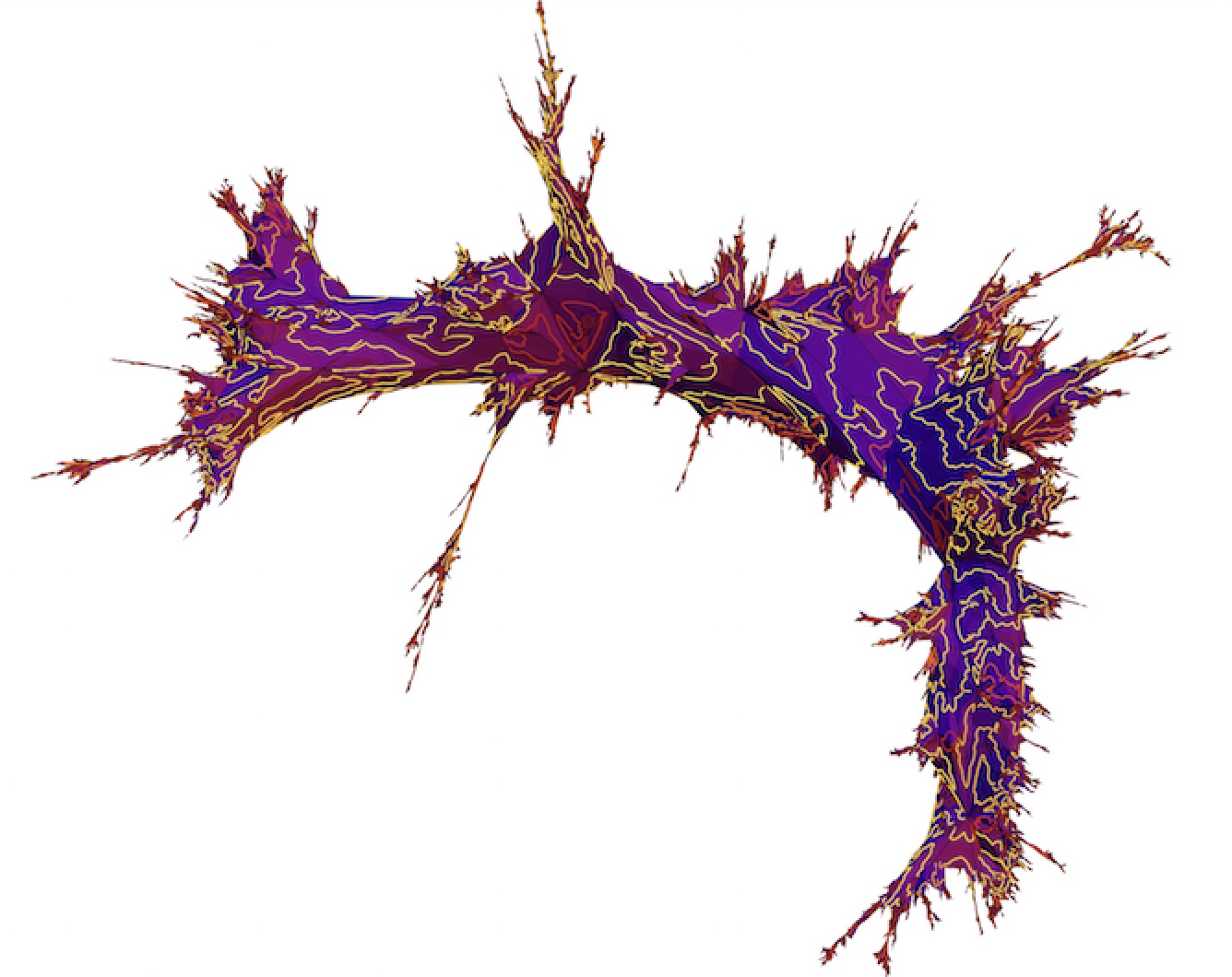}
\includegraphics[width = .6\textwidth]{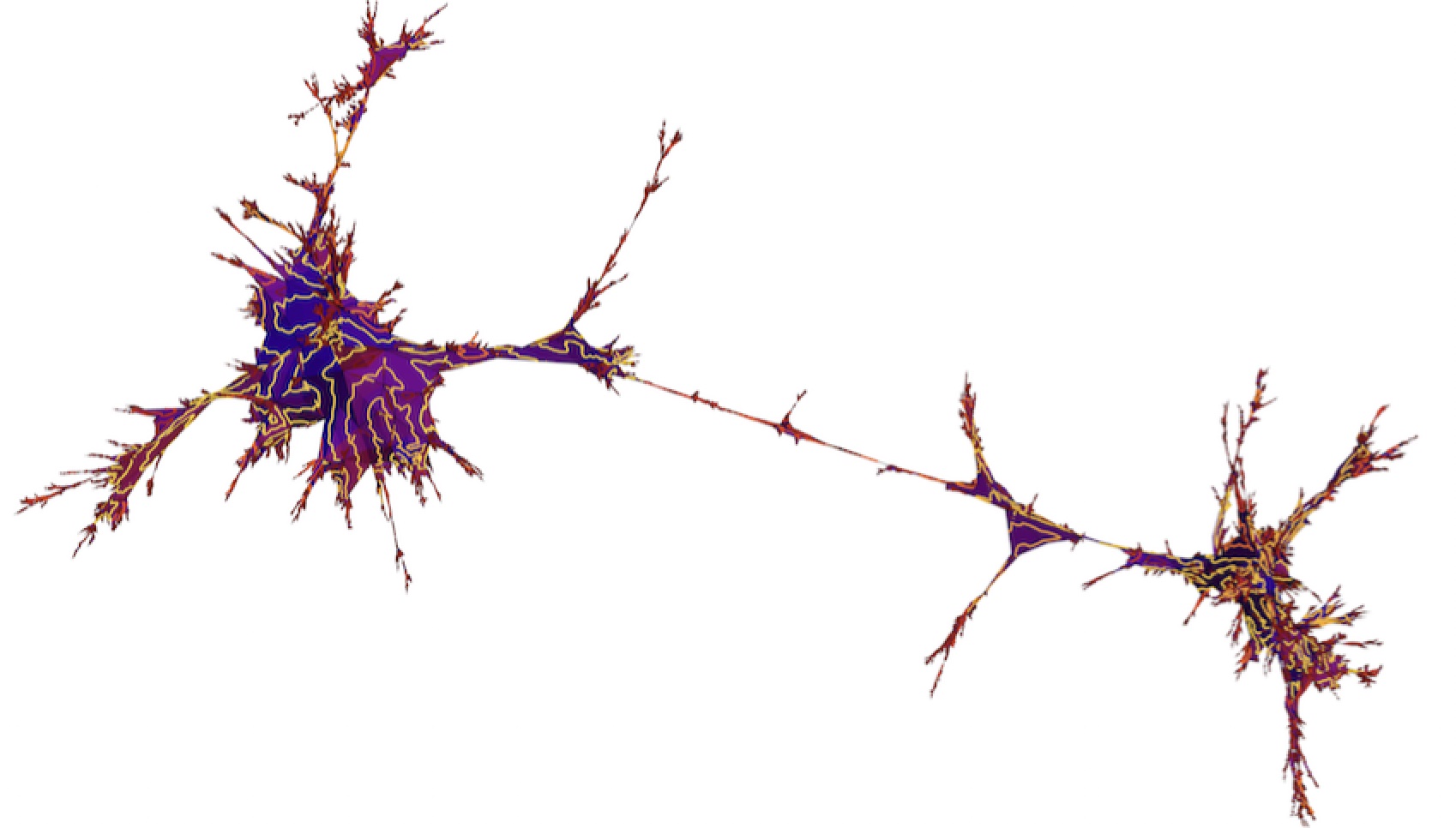}
\caption{A map weighted by the FK model with $q=0.5$, $q = 2$ (corresponding to the FK-Ising model) and $q=9$ respectively, together with some of their loops. Simulation by J. Bettinelli and B. Laslier. When $q>4$ it is believed that the maps become tree-like, and the limiting metric space should be Aldous' continuum random tree.\label{F:maps}}
\end{center}
\ind{Continuous Random Tree (CRT)}
\end{figure}

%\ellen{(Should this discussion belong in the next subsection - maybe combined with the part where it's discussed again?)}

\paragraph{Uniform random planar maps.} Observe that when $q=1$, the FK model \eqref{FK} has the property that the map $M_n$ is chosen \textbf{uniformly} at random among the set $\cM_n$ all of (rooted) maps with $n$ edges, because the total number of possible configurations for $\bt_n$ is $2^n$ \emph{independently} of $\bm_n$. Furthermore, given $M_n = \bm_n$, $T_n$ is chosen uniformly at random from the $2^n$ possibilities: this corresponds to each edge being present (open) with probability $1/2$, independently of one another. In other words,  $T_n$ corresponds to bond percolation with parameter $1/2$ given the map $M_n$. This is in fact the critical parameter for this percolation model, as shown in the work of Angel \cite{Angel}.

\ind{Planar maps!Uniform case}
\ind{Planar maps!Percolation}

The case of a uniformly chosen planar map in $\cM_n$ is one in which remarkably detailed
information is known about its structure. In
particular, a landmark result due to Miermont \cite{Miermont} and Le
Gall \cite{LeGall} shows that, viewed as a metric space and rescaling
edge lengths to be $n^{-1/4}$, the random map converges to a  multiple
of a certain universal random metric space known as the
\textbf{Brownian map}. (In fact, the results of Miermont and Le Gall
apply respectively to uniform quadrangulations with $n$ faces and to
$p$-angulations for $p=3$ or $p$ even, whereas the convergence result
concerning uniform planar maps in $\cM_n$ was established a bit
later by Bettinelli, Jacob and Miermont
\cite{betscaling}). Critical percolation on a related
half plane version of the maps has been analysed in a work of
Angel and Curien \cite{AngelCurien}, while information on the full
plane percolation model was later obtained by Curien and
Kortchemski \cite{CK13}. Related works on loop models (sometimes
rigorous, sometimes not) appear in \cite{guionnet12,borot12,eynard95,b11, BBG,CCM19}.

%\paragraph{Ising.} Likewise, when $q = 2$, this corresponds to the FK-Ising model, and we find that $\kappa =  16/3$. Here again, this is known rigorously \ellen{(what is known rigourously?)} by results of Chelkak and Smirnov. \ind{Planar maps!Ising}

One reason for the particular choice of the FK model in \eqref{FK} is the belief that for $q<4$, after Riemann uniformisation, a large sample of
such a map closely approximates a \emph{Liouville quantum gravity} surface. We will try to summarise this conjecture in the next subsection.

\begin{comment}
 described by an area measure of the form $e^{\gamma h(z)}dz$, where the parameter $\gamma\in(0,2)$ is believed to be related to the parameter $q$ of \eqref{FK} by the relation
\begin{equation}
\label{q_gamma}
q = 2 + 2\cos\left(\frac{8\pi}{\kappa'}\right); \quad \gamma = \sqrt{\frac{16}{\kappa'}} .
\end{equation}
Note that when $q \in (0,4)$ we have that $\kappa' \in (4,8)$ so that
it is necessary to generate the Liouville quantum gravity with the
associated dual parameter $\kappa  = 16/\kappa' \in (0,4)$. This
ensures that $\gamma = \sqrt{\kappa } \in (0,2)$.
\end{comment}

\subsection{Conjectured connection with Liouville quantum gravity}\label{sec:RPM_LQG}

%\ellen{(I think we could reorganise this section. Just to be a bit clearer concerning what is a conjecture about limits of loops, limits of measures, or both. Also there is some repetition, which is fine, but should probably emphasise what is repetition and what is not.)}

The distribution \eqref{FK} gives us a natural family of distributions on planar maps (indexed by the parameter $q \ge 0)$. As already mentioned, in this model, the weight of a particular map $\bm\in \cM_n$ is proportional to the \emph{partition function} $Z(\bm,q)$ of the critical FK model on the map. %$\bm\in \cM_n$ (the set of planar maps with $n$ edges).
\begin{comment}
Now, in the Euclidean world, scaling limits of critical FK models are widely believed to be related to collection of loops known as CLE$_\kappa$, where $\kappa$ is related to the parameter $q$ via the relation
\begin{equation}
\label{q_kappa}
q = 2 + 2\cos\left(\frac{8\pi}{\kappa'}\right).%; \quad \gamma = \sqrt{\frac{16}{\kappa'}} .
\end{equation}
\ellen{(note this is the same as formula above?)}
For instance, when $q=1$, the FK model reduces to percolation and one finds $\kappa = 6$, which is of course now rigorously established thanks to Stas Smirnov's proof of conformal invariance of critical percolation. \ellen{(Do we want to add a reference?)}
\end{comment}

\paragraph{Conformal Embedding.} Suppose that $q < 4$ in what follows. It is strongly believed that in the limit $n\to \infty$, the geometry of such maps are related to Liouville quantum gravity with parameter $\gamma$, where
\begin{equation}
\label{q_gamma}
q = 2 + 2\cos\left(\frac{\gamma^2\pi}{2}\right).
%\quad \gamma = \sqrt{\frac{16}{\kappa'}} .
\end{equation}

(Note that this equation has no real solution if $q >4$.)

%atisfying $\gamma^2 = 16/\kappa'$. In SLE theory, the value $16/\kappa$ is special, as it transforms an SLE into a dual SLE: roughly, the exterior boundary of an SLE$_\kappa$ is given by a form of SLE$_{16/\kappa}$ when $\kappa \ge 4$. (Note that the self dual point $\kappa = 4$ corresponds to $\gamma = 2$, which we have already seen is critical for Liouville quantum gravity.)
To be more precise about this, one must relate the world of planar maps to the world of Liouville quantum gravity by specifying a ``natural'' embedding of the maps into the plane. There are various ways to do this, and a couple of the simplest are as follows.

\begin{figure}
\begin{center}
\includegraphics[scale=.2]{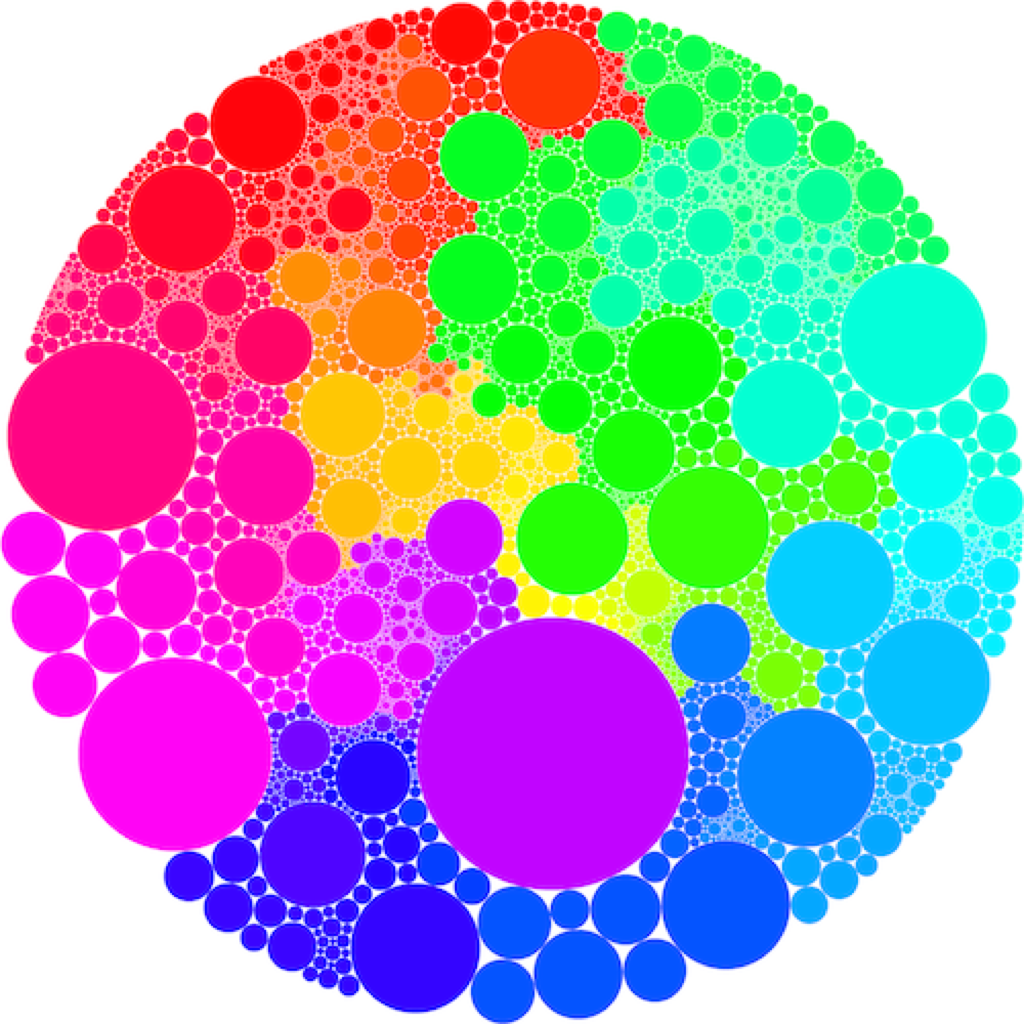}
\caption{Circle packing of a uniform random planar map. Simulation by Jason Miller.}
\label{F:CP}
\end{center}
\end{figure}

\begin{itemize}
	\item
 \textbf{Via the circle packing theorem.}
 By a theorem of Koebe--Andreev--Thurston (see the book by K. Stephenson \cite{Stephenson} for a comprehensive introduction), any planar map can be represented as a circle packing. A circle packing is a collection of circles in the plane such that any two of the corresponding discs either are tangent to one another, or do not overlap.
   In the circle packing representation, the vertices of the map are given by the centres of the circles, and the edges correspond to tangent circles. See \cref{F:CP}. Each circle packing representation of a map gives an embedding in the plane, and when the map is a simple triangulation, this embedding is unique up to M\"obius transformations.
\ind{Circle packing}

\item \textbf{Via the uniformisation theorem.} In this approach, a given map is viewed as a Riemann surface by declaring that each face of degree $p$ is a regular $p$-gon of unit area, endowed with the standard metric, and specifying the charts near a vertex in the natural manner. This Riemann surface can then be embedded into the disc (say) by the uniformisation theorem (which is a generalisation of the Riemann mapping theorem from subsets of $\C$ to arbitrary Riemann surfaces). \ind{Riemann uniformisation}
\end{itemize}

These embeddings are essentially unique up to M\"obius transforms (in the first case, we can circle pack the refinement map ${\bar \bm_n}$ instead of ${\bm_n}$). The choice of M\"obius transform can be fixed by requiring, for instance, that the root edge is mapped to $(0,1)$.

Once an embedding has been chosen, a natural object to study is the measure $\mu_n$ in the plane which puts mass $1/N$ ($N$ being the number of vertices in $M_n$) at the position of each embedded vertex. The conjecture alluded to above says that in the limit as $n \to \infty$, if $M_n$ is sampled from \eqref{FK}, then this measure $\mu_n$ should converge to $\gamma$-Liouville quantum gravity. More precisely, if $\gamma$ and $q$ are related by \eqref{q_gamma}, it should converge in distribution for the topology of weak convergence, to a \emph{variant} of the Liouville measure $\mu_\gamma$ (this variant will be specified, for example, in Chapter \ref{S:LCFT}).

{\begin{rmk}\label{R:uniformmaps}
		Note that when $q=1$, which we have already discussed is the case of uniformly chosen random planar maps, we have $\cos(\gamma^2 \pi/2)=-1/2$, that is, $\gamma=\sqrt{8/3}$. Consequently, the limit of a (conformally embedded) uniformly chosen map should be related to Liouville quantum gravity with this parameter. This has been  verified for a slightly different type of conformal embedding called the \textbf{Cardy embedding} in a recent breakthrough of Holden and Sun \cite{HS19}.
\ind{Cardy embedding}
	\end{rmk}}

\paragraph{Loops and CLE.} The loops induced by the FK model \eqref{FK} may be viewed as a decoration on the map. Indeed as we have already mentioned, given the map, they are the cluster boundaries of a self dual FK percolation model on it with parameter $q$. It is therefore natural to wonder about their geometry in the scaling limit, after embeddings of the type discussed above. The widely shared belief is that they converge to so called \textbf{conformal loop ensembles} CLE$_{\kappa'}$ where the parameter $\kappa'$ is given by
\begin{equation} \label{q_kappa}
\kappa' = \frac{16}{\gamma^2}; \text{ and thus }
q = 2 + 2\cos\left(\frac{8\pi}{\kappa'}\right).
%\quad \gamma = \sqrt{\frac{16}{\kappa'}} .
\end{equation}\indN{{\bf Parameters}!$\kappa'$; dual parameter value of $\kappa\in (0,4]$, $\kappa'=16/\kappa$}
\indN{{\bf Parameters}! $q$; FK model parameter}
In fact, one can also study the self dual FK percolation model and its associated loops on a Euclidean lattice, and the same belief is held. That is, these loops are also conjectured to converge to CLE$_{\kappa'}$ in the scaling limit, where the relationship between $q$ and $\kappa'$ is the same as in \eqref{q_kappa}. The fact that these two conjectures are the same should heuristically be considered as a consequence of conformal invariance. That is, if the scaling limit of FK loops is conformally invariant, it should be independent of the underlying metric: only their conformal type should matter.

For instance, we have already noticed that when $q=1$, the associated FK model is just bond percolation. In this case we already know (at least in the case of the triangular lattice) that the scaling limit of the associated loops is given by CLE with parameter $\kappa' = 6$ (\cite{Smirnov}, \cite{CamiaNewman}). This is consistent with the value $\gamma = \sqrt{8/3}$ being the Liouville quantum gravity parameter for the scaling limit of uniform planar maps, as described in Remark \ref{R:uniformmaps}.

Likewise, for $q=2$ the associated FK model is the FK representation of the critical Ising model. It was proven in \cite{KemppainenSmirnovII} (see also \cite{IsingSLE} for interfaces and \cite{BenoistHongler} for Ising loops)  that the scaling limit of these loops is given by CLE$_{16/3}$. The associated parameter $\gamma$ is thus $\gamma = \sqrt{3}$.

A small summary of these values is provided in the table below.

%\ellen{\textbf{Suggestion (??)} \begin{itemize} \item discuss (expected/known) convergence of FK interfaces on Euclidean lattices to SLE \item  discuss that it is also expected to hold on random lattices \item observe this leads one to expect a connection between SLE$_\kappa$ and $\gamma$-LQG with $\gamma^2=\kappa$ \item refer to section 6 where this is shown to be the case
%\item describe informally what will be shown in the next section of this chapter and how this corresponds to continuum mating-of-trees, demonstrating $\gamma$-LQG/SLE$_\kappa$ link. \end{itemize} }

\begin{center}\begin{tabular}{|c||c|c|c|} 
	\hline
	\textbf{ FK Model \eqref{FK}} & $q$ & $\gamma$ &  $\kappa'$ \\
	\hline
	\hline
	General $q\in [0,4)$ & $2+2\cos(\gamma^2\pi/2)$ & $\gamma\in [\sqrt{2},2)$ & $16/\gamma^2\in (4, 8]$ \\
	\hline
	Uniform map + critical bond percolation  &  $1$ &  $\sqrt{8/3}$   & $6$ \\
	\hline
	Spanning tree decorated map & $ 0$ & $\sqrt{2}$  & 8 \\
	\hline
	Critical Ising decorated map & $2$ & $\sqrt{3}$  & $16/3$ \\
	\hline
%\ellen{(?)}	$q$-Potts model & $q=3,4,...$ & - & - & - \\	
%\hline
\end{tabular} \label{Table:values}
\end{center}
\ind{FK model}

%Either way, it is believed that in the $n\to \infty$ limit, after appropriate rescaling, the empirical distribution on vertices converges %of the triangulation converges, after appropriate rescaling,
%to a version of Liouville quantum gravity with
%\begin{equation}
%\label{q_kappa}
%q = 2 + 2\cos\left(\frac{8\pi}{\kappa'}\right); \quad \gamma = \sqrt{\frac{16}{\kappa'}} .
%\end{equation}

%the parameter $\gamma = \sqrt{16/\kappa}$, and $\kappa$ is associated $q$ by $q = 2 + 2\cos(8\pi/\kappa)$.

%There are now several results making this connection precise at various levels: see \cite{BLR, DuplantierMillerSheffield, FKstory1}. \ellen{(Should update.)}

%\begin{rmk}
%Since the partition function of percolation is constant, a map weighted by the percolation model ($\kappa = 6$) is the same as a uniformly chosen random map. Consequently, the limit of a uniformly chosen map should be related to Liouville quantum gravity with parameter $\sqrt{16/\kappa} = \sqrt{16/6} = \sqrt{8/3}$.
%\end{rmk}

%\ellen{(Could add the examples we removed from Chapter 3 here?)}

\subsection{Mullin--Bernardi--Sheffield's bijection in the case of spanning trees}
\label{sec:bij_tree}

%We briefly recall the Hamburger--Cheeseburger bijection due to
%Sheffield (see also related works by Mullin \cite{mullin67} and
%Bernardi \cite{bernardi07,bernardi08}).

We will now discuss the case where the map {$M_n\in \cM_n$} is chosen with probability proportional to the number of spanning trees {it admits}. Here a spanning tree is a collection of unoriented edges, covering every vertex, and which contains no cycle. (By contrast, in Section \ref{sec:LERWexp} we will also encounter spanning trees for which a designated vertex called the root has been singled out; in which case one may view the edges in the spanning tree as being oriented towards the root). 
In other words, {for any (rooted) map $\bm_n\in \cM_n$ with $n$ edges and $\bt_n$ a set of edges on it}
\begin{equation}\label{maptree}
{\P(M_n=\bm_n, T_n=\bt_n) \propto \indic{\bt_n \text{ is a spanning tree on }\bm_n}.}
\end{equation}
 This can be understood as the limit when $q \to 0^+$ of the Fortuin--Kasteleyn model discussed above in \eqref{FK}, since in this limit the model concentrates on configurations where $\ell = 0$, equivalently, $\bt_n$ is a tree. In fact it is immediate in this case that given {$M_n=\bm_n$}, $\bt_n$ is a uniform spanning tree (UST) on $\bm_n$. We will discuss a powerful bijection due to Mullin \cite{mullin67} and Bernardi \cite{bernardi07,bernardi08} which is key to the study of such planar maps. This bijection is actually a particular case of a bijection due to Sheffield, which is sometimes called the ``hamburger--cheeseburger'' bijection. Sheffield's bijection can be used for arbitrary $q\ge 0$, however the case $q=0$ of trees is considerably simpler and so we discuss it first. (We will use the language of Sheffield, in order to prepare for the more general case later.) Although the hamburger--cheeseburger bijection is the only example we will treat in detail here, we mention that there are other powerful bijections of a similar flavour that can be used to connect random planar map models to Liouville quantum gravity and SLE: see for example \cite{BHS,schnyderwoods,GKMW,KMSW}.

\begin{figure}
\begin{center}
\includegraphics[width=.2\textwidth]{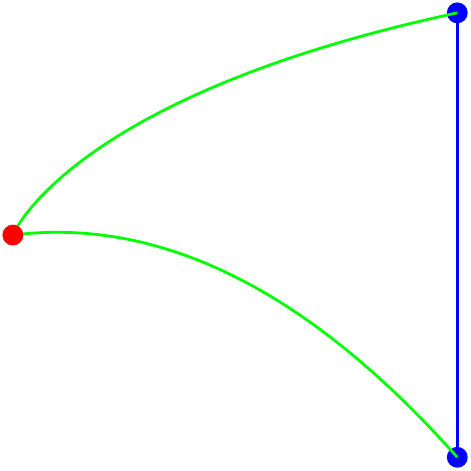} \quad \quad %or \quad
\quad
\includegraphics[width=.2\textwidth]{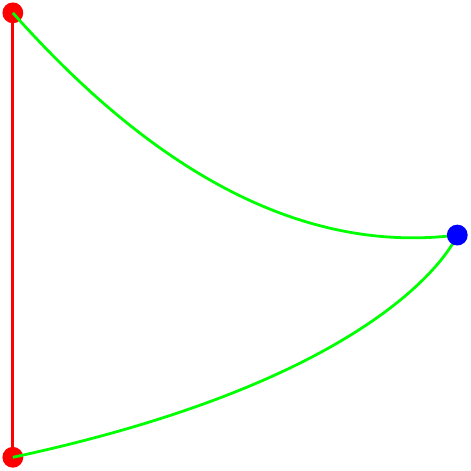}\quad \quad
%\end{center}
%\begin{center}
\includegraphics[width=.4\textwidth]{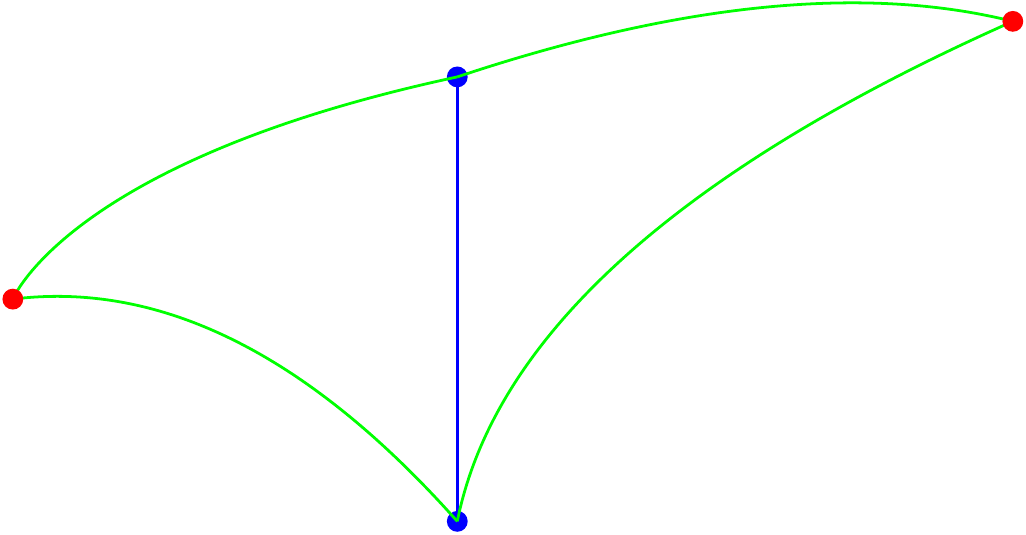}
\caption{Refined or green edges split the map and its dual into primal and dual triangles. Each primal triangle sits opposite another primal triangle, resulting in a primal quadrangle as above. }
\label{F:triangle}
\end{center}
\end{figure}

\medskip {To describe the $q=0$ hamburger--cheeseburger bijection, we first fix a deterministic pair $(\bm_n,\bt_n)$ as above (with an oriented root edge chosen for $\bm_n$ and $\bt_n$ a spanning tree on $\bm_n$) -- see Figure \ref{F:maptree} -- and describe how to associate with it a certain sequence of letters corresponding to ``hamburgers'' and ``cheeseburgers''}. Recall that adding refinement edges to a map splits it into triangles of exactly two types: primal triangles (meaning two refined edges and one primal edge) or dual triangles (meaning two refined edges and one dual edge). For ease of reference, primal triangles will be associated with hamburgers, and dual triangles with cheeseburgers. Note that for a primal edge in a primal triangle, the triangle opposite that edge is obviously a primal triangle too. Hence it is better to think of the map as being split into quadrangles with either a primal or dual diagonal, % one diagonal is primal or dual
as illustrated in \cref{F:triangle}.

\medskip We will reveal the map, triangle by triangle, by exploring it along a space-filling (in the sense that it visits every triangle once) path. When we do this, we will
keep track of the first time that the path enters a given quadrangle by saying
that {either a hamburger or a cheeseburger is produced, depending on whether the quadrangle is primal or dual.}
%a hamburger (or cheeseburger if it is a dual quadrangle) is produced at
%that point.
Later on, when the path comes back to the quadrangle for the second
and final time, we will say that the burger has been eaten. We will use the
letters $\ah, \ac$ to indicate that a hamburger or cheeseburger has been
produced and we will use the letters $\aH, \aC$ to indicate that a burger
has been eaten (or \emph{ordered} and eaten immediately). So in this description we will have one letter for every triangle.
\ind{Sheffield's bijection}
\ind{Hamburgers / cheeseburgers| see {Sheffield's bijection} }

\medskip It remains to specify in what order are the triangles visited;
equivalently, to describe the space-filling path. In the case {that we consider now, where the
decoration $\bt_n$ consists of a single spanning tree,} %(corresponding to $q=0$ as
%we will see later)
the path is simply the contour path going around the
tree {(starting from the root)}, that is, the unique loop $L_0$ separating the primal and dual spanning trees, with its orientation inherited from that of the root edge of $\bm_n$. Hence in this case, we can associate to $(\bm_n,\bt_n)$  a sequence $w$ (or \textbf{word}) made up of $M$
letters in the alphabet $\Theta = \{\ah, \ac, \aH, \aC\}$. We will see below that subject to certain natural conditions on the word $w$, this map is actually a bijection.
%\ellen{(Should we reference Figure 9 here? See also comment below. Or do you mean that the rooted map describes a sequence of $N$ letters?)}
 %We also claim that it is always the case that

\begin{figure}
	\centering
	\includegraphics[scale = 0.8]{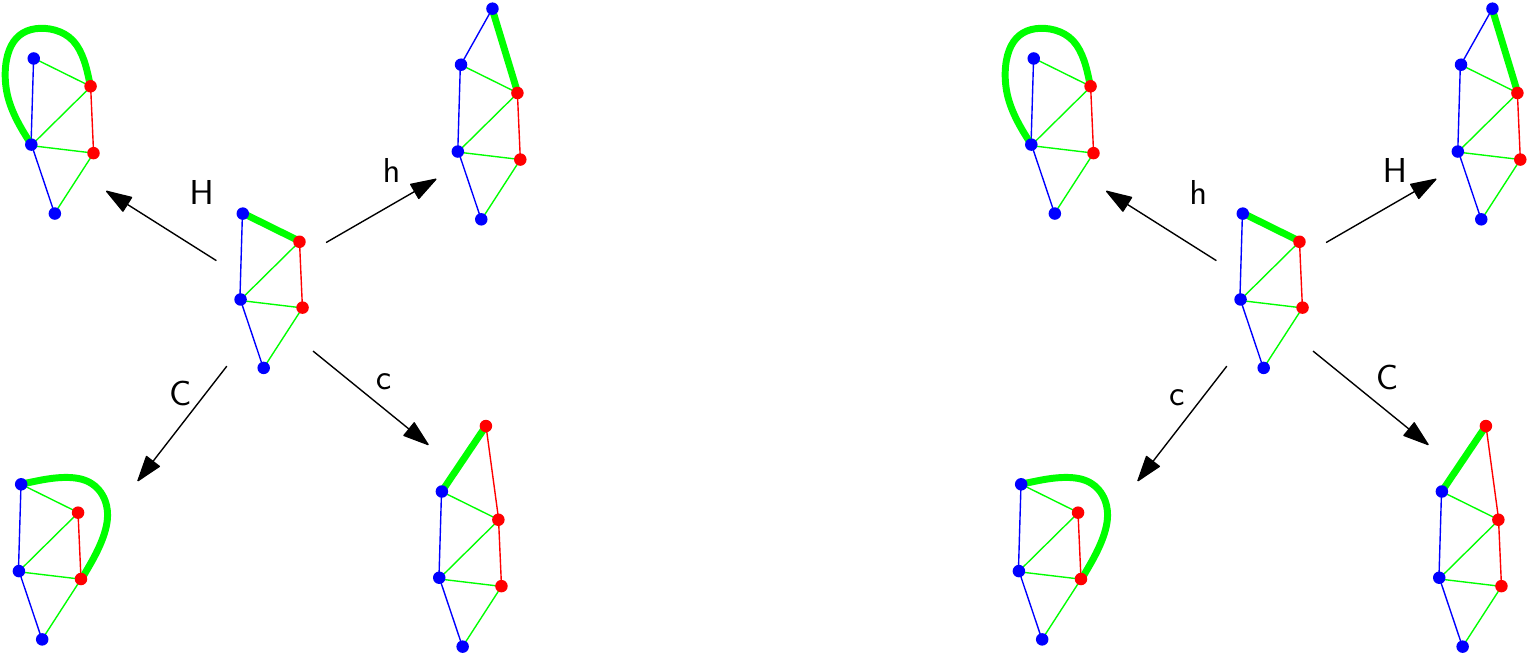}
	\caption{From symbols to map. The current position of the interface (or last discovered refined edge) is indicated with a bold line.
		Left: reading the word sequence from left to right or \emph{into the future}. The map in the centre is formed from the
		symbol sequence $\ah\ac\ah$.
		%The operation of gluing triangles depending on the next symbol is shown.
		Right: The corresponding operation when we go from right
		to left (or into the \emph{past}); this is useful for instance when taking a local limit, see Section \ref{sec:local}. The map in the centre now corresponds to the
		symbol sequence $\aH\aC\aH$.
%\ellen{*I moved this figure to where it is first referenced.}
}
	\label{fig:map_from_word}
\end{figure}

{Observe that we always have} $M = 2n$. To see why, recall that there is one letter for every triangle, so $M$ is the total number of triangles. Moreover, each triangle can be identified with an edge (or in fact half an edge, because each edge is visited once when the burger is produced and once when it is eaten), and so
$$
M = 2(E(\bt_n) + E( \bt^\dagger_n)) = 2 ( V(\bt_n) -1  + V(\bt^\dagger_n) -1).
$$
%Now, the number of vertices of $\bt_n$ is the same as $\bm_n$
Now $V(\bt_n)=V(\bm_n)$, and  %the number of vertices of $\bt^\dagger_n$ is the same as the number of vertices in the dual map,
$V(\bt^\dagger_n)=V(\bm_n^\dagger)=F(\bm_n)$. %which is equal to the number of faces of $\bm_n$.
This gives that
\begin{equation}\label{Euler}
M = 2 (V(\bm_n) + F(\bm_n) - 2),
\end{equation}
and applying Euler's formula together with the fact $E(\bm_n)=n$, we find that $M = 2n$. Alternatively note directly that $E(\bt_n) + E( \bt^\dagger_n) = E(\bm_n) = n$ since each edge of $\bm_n$ corresponds to an edge that is either open in $\bt_n$ or $\bt^\dagger_n$.
\ind{Euler's formula}

\begin{figure}
\textbf{a.} \includegraphics[width = .4\textwidth]{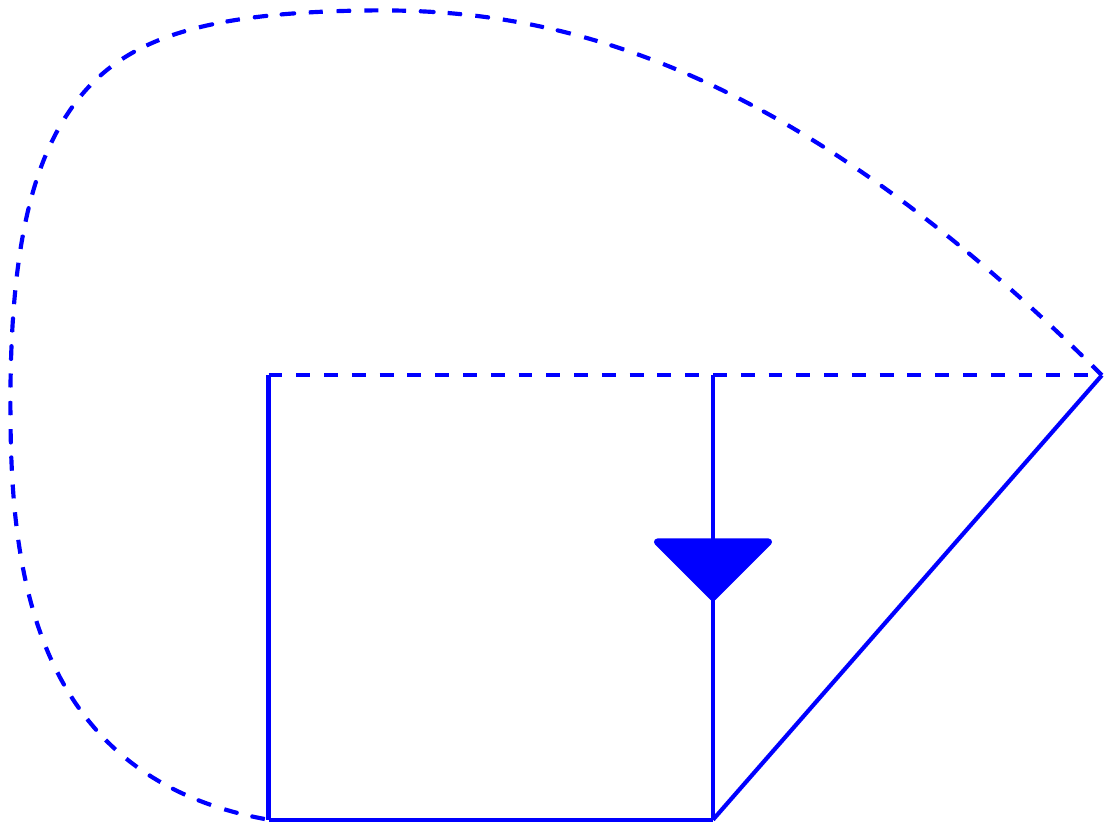} \quad \quad
\textbf{b.}  \includegraphics[width = .4\textwidth]{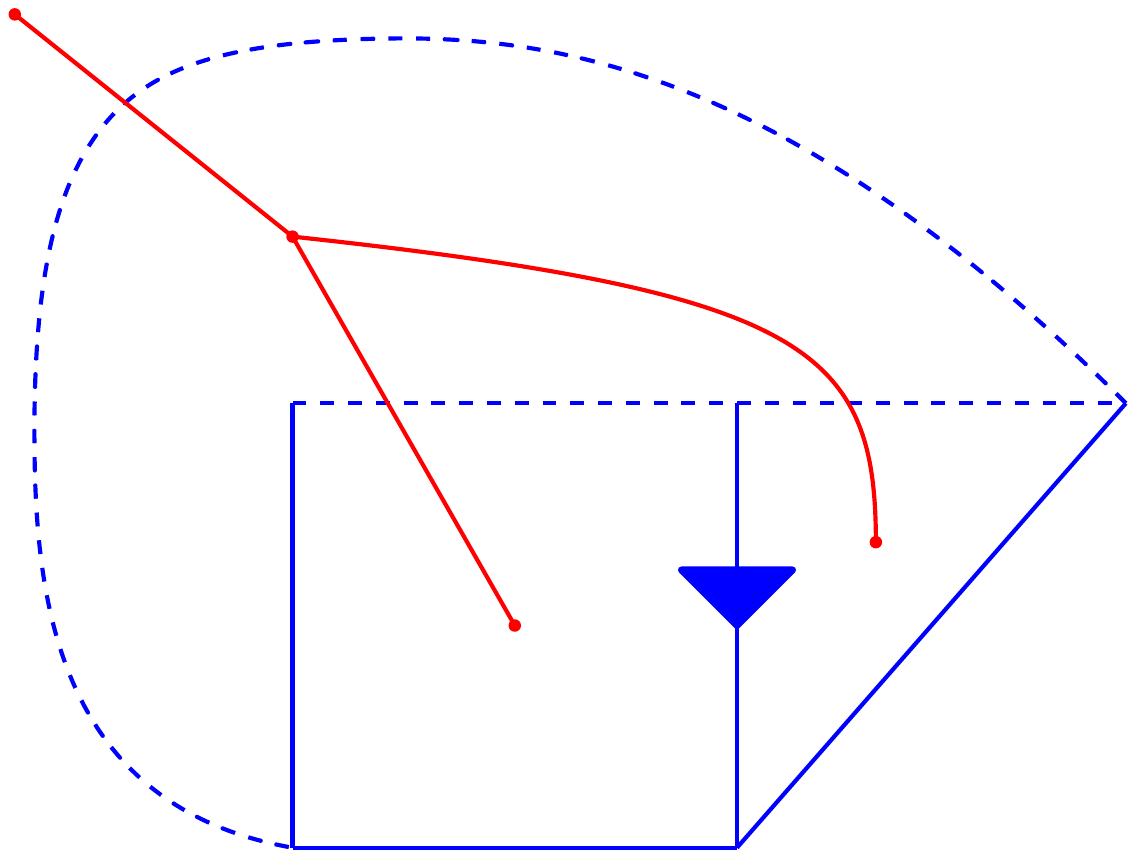}\quad \quad

\vspace{.5cm}

\textbf{c.} \includegraphics[width = .4\textwidth]{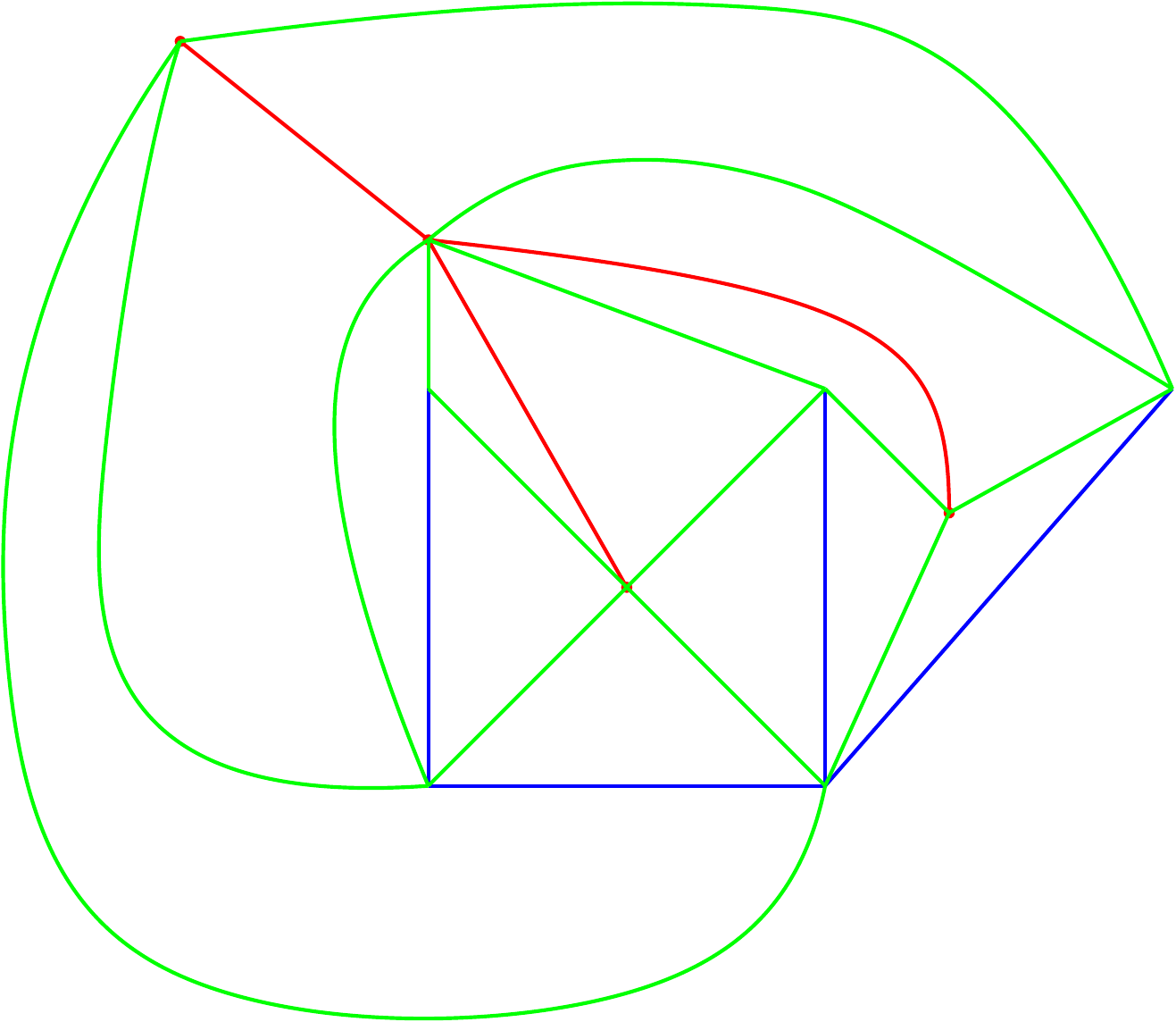}\quad \quad
\textbf{d.}  \includegraphics[width = .4\textwidth]{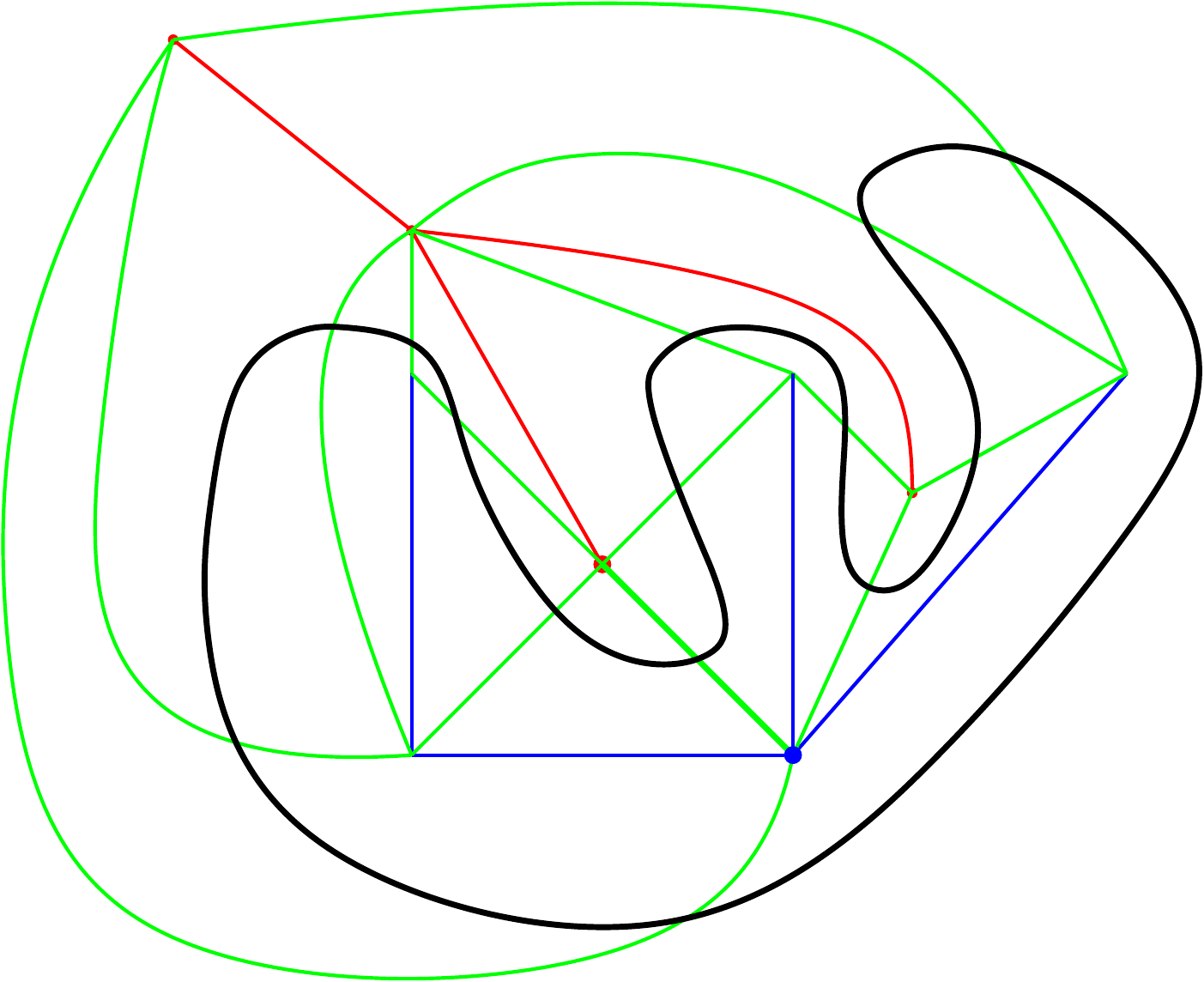}
\caption{ \textbf{a:} a map with a spanning tree. \textbf{b:} Spanning tree and dual tree. \textbf{c:} Refinement edges. \textbf{d:} Loop separating the primal and dual spanning trees, to which a root (refined) edge has been added in bold.\label{F:maptree}}
%\caption{The word hamburger-cheeseburger sequence associated to this tree is: h c c H h C c H h C h C H H }
  %
  %
\end{figure}

{To summarise, given $(\bm_n,\bt_n)$ a rooted, spanning tree decorated map with $n$ edges, we can uniquely define a word $w$ of length $2n$ in the letters $\{\ah,\ac,\aH,\aC\}$. Observe further that under the reduction rules
$$ \overline{\ac \aC } = \overline{\ah \aH  }= \emptyset, \; \overline{\ac\aH} = \overline{\aH \ac} \text{ and } \overline{\ah \aC} = \overline{\aC \ah},
$$
we have $\bar{w}=\emptyset$ (here $\bar{w}$ denotes the reduction of the word $w$). This corresponds to the fact that every burger produced is eaten, and every food order corresponds to a burger that was produced before. Subject to the condition $\bar w = \emptyset$, it is easy to see that the map $(\bm_n, \bt_n) \mapsto w$ is a bijection. See, for example, Figure \ref{fig:map_from_word} for elements of a proof by picture.}

{Now we go a step further, and associate to this word $w$ a pair} %the quantities
$(X_k, Y_k)_{1\le k \le 2n}$, which count the number of hamburgers and cheeseburgers respectively in the stack at any given time $1 \le k\le 2n$ (that is, the number of hamburgers or cheeseburgers which have been produced prior to time $k$ but {get eaten }%whose matches come
	after time $k$). Note that $(X, Y)$ is a process which starts from the origin at time $k = 0$, and ends at the origin at time $k = 2n$.
Moreover, by construction $X$ and $Y$ both stay non-negative throughout. We call a process $(X_k,Y_k)_{0 \le k \le 2n}$ satisfying these properties a \textbf{discrete excursion} (in the quarter plane).
\ind{Discrete excursion}
{So at this point, we have associated with any $(\bm_n, \bt_n)$ as above, a unique discrete excursion $(X,Y)$ of length $2n$.}

%It is obvious that $(X,Y)_{0 \le k \le 2n}$ identify uniquely  the word $w$ encoding the decorated map $(\bm_n, \bt_n)$. \ellen{(I removed this first sentence? Isn't it saying the same thing as the next?)}
Conversely, given such a process $(X,Y)$ we can associate to it a word $w$ in the letters of $\Theta$ such that $(X,Y)$ is the net burger count of $w$. Obviously $w$ reduces to $\emptyset$ and so, as we have seen above, this word $w$ specifies a unique pair $(\bm_n, \bt_n)$.
%\ellen{(Explain why $(X,Y)$ therefore specifies a unique rooted spanning tree decorated map $(\bm_n, \bt_n)$ with $n$ edges? Reference Figure 9?)}

\begin{figure}
\textbf{a.} \includegraphics[width = .4\textwidth]{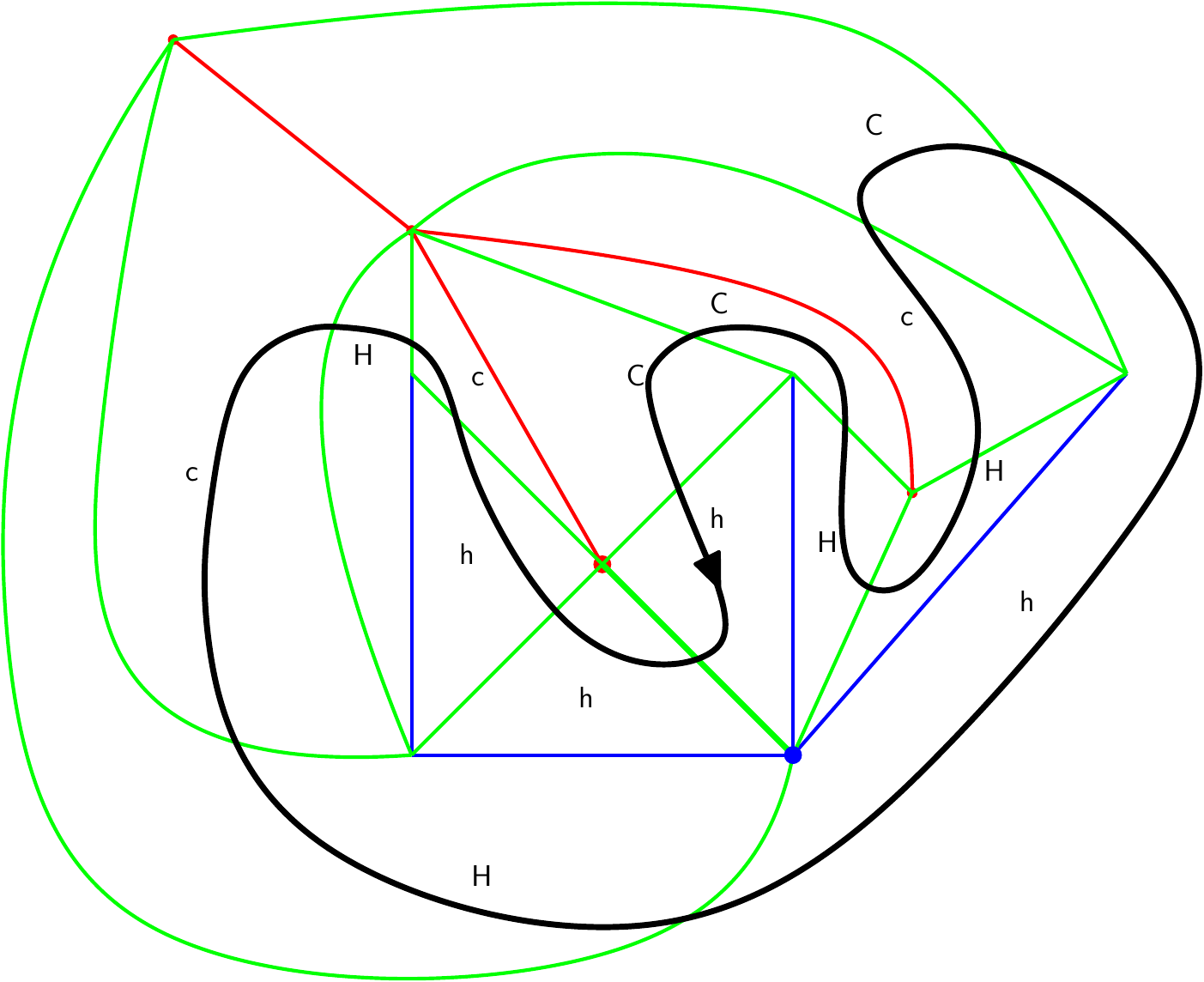} \quad 
\textbf{b.} \includegraphics[width = .48\textwidth]{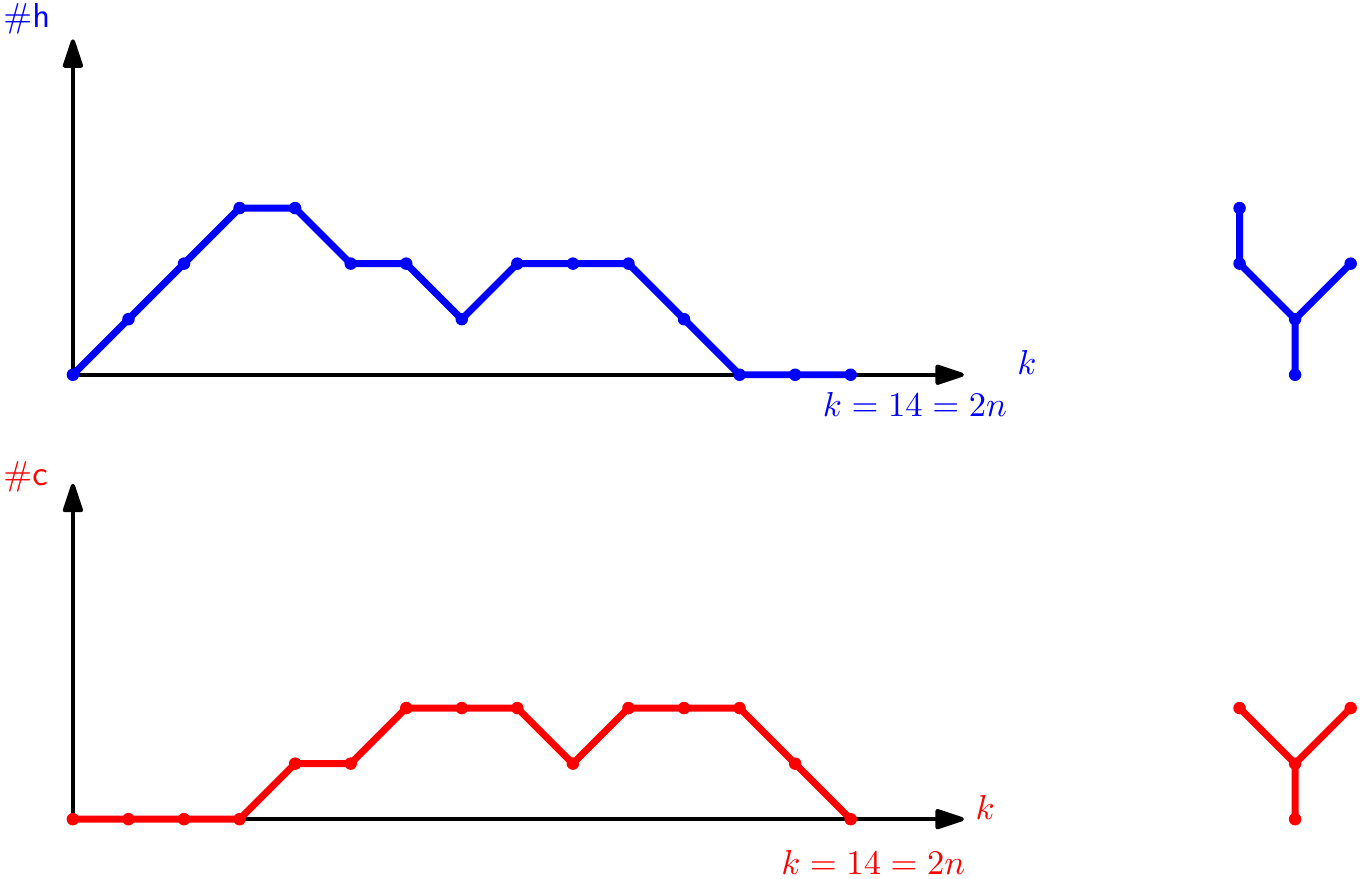}
\caption{\textbf{a:} The word associated to $(\bm_n, \bt_n)$ is: {$w=\ah \ah \ah \ac \aH \ac\aH \ah \aC \ac \aH \aH \aC \aC$} 
	\textbf{b:} The hamburger and cheeseburger counts, as well as the trees encoded by these excursions (which are identical to the primal and dual spanning trees, respectively).
\label{F:maptree2}}
\end{figure}

Another property which is easy to check (and easily seen on Figure \ref{F:maptree2}) is that the excursions $X$ and $Y$ encode the spanning tree $\bt_n$ and dual spanning tree $\bt^\dagger_n$ in the sense that they are (after removing steps where $X$, respectively $Y$, remain constant) the contour functions of these trees. More precisely, at a given time $k$, $X_k$ denotes the height in the tree (distance to the root) of the last vertex discovered prior to time $k$. 

\begin{rmk}\label{R:queues}
It may be useful to recast the above connections in the language of \textbf{queues}, where customers are being served one at the time. More precisely, a queue (in discrete time) is a process where at each unit of time either a new customer arrives, or a customer at the front of the queue is being served and leaves the queue forever. Any queue can be equivalently described by a tree $\bt$ or an excursion $X$. Indeed, a tree structure $\bt$ can be defined from the queue, by declaring that any customer $c$ arriving during the service of a customer $c'$ is a child of $c'$. An excursion $X$ can be defined by simply counting the queue length at each time. Note that $X$ is nothing  but the contour function of the tree $\bt$ (meaning the discrete process which measures the height of the tree $\bt$ as it goes around it in depth-first order; see \cite{LeGallsurvey} for much more about this). In our case, the tree $\bt$ is simply either the primal spanning tree on the map or its dual.
\end{rmk}

{When $(M_n,T_n)$ are \emph{random} and sampled according to \eqref{maptree}, the corresponding random excursion $(X,Y)$ is clearly chosen uniformly from the set of all possibilities.} It  therefore follows from classical results of Durrett, Iglehart and Miller \cite{excursion} that as $n \to \infty$, $$
\frac1{\sqrt{n}} (X_{\lfloor 2n t\rfloor } , Y_{\lfloor 2nt \rfloor})_{0 \le t \le 1} \to (e_t, e'_t)_{0 \le t \le 1}
$$
where $e, e'$ are independent  Brownian (one dimensional) excursions (that is, the pair $(e,e')$ is Brownian excursion in the quarter plane), for example, in the Skorokhod sense (alternatively for the topology of uniform convergence if the paths are linearly interpolated instead of piecewise constant as above). This property implies (see for example Lemma 2.4 in Le Gall's comprehensive survey \cite{LeGallsurvey}) that, in the Gromov--Hausdorff sense, the primal and dual spanning trees converge after rescaling the distances by a factor $n^{-1/2}$, to a pair of independent \textbf{Continuous Random Trees} (CRTs) \cite{Aldous}.

We summarise our findings, in the case of UST weighted random planar maps, in the following theorem.

\begin{theorem}
  \label{T:ust} The set of {(rooted)} spanning tree decorated maps $(\bm_n, \bt_n)$ {with $n$ edges} are in bijection with excursions $(X_k, Y_k)_{0 \le k \le 2n}$ in the quarter plane. When $(M_n,T_n)$ is random and distributed according to \eqref{maptree}, the pair of trees $(T_n,T_n^\dagger)$ converges, for the Gromov--Hausdorff topology and after scaling distances (in each tree) by a factor $n^{-1/2}$, to %Scaled by $\sqrt{n}$, the pair of trees $(\bt_n, \bt^*_n)$ converges to
  a pair of independent Continuous Random Trees (CRTs).
\end{theorem}
\ind{Continuous Random Tree (CRT)}

Note that the map %$\bm_n$
{$M_n$} itself can then be thought of as a gluing of two discrete trees (that is, the primal and dual spanning trees, which are glued along the space-filling path). In the scaling limit, this pair of trees becomes a pair of independent CRTs. As it turns out, the procedure of gluing these two trees has a continuum analogue, which is described in the work of Duplantier, Miller and Sheffield \cite{DuplantierMillerSheffield}. This is the \textbf{mating of trees} approach to LQG, and is an extremely powerful and fruitful  point of view that we will describe in more detail later on. 
\ind{Mating of trees (discrete)}

\subsection{The loop-erased random walk exponent}
\label{sec:LERWexp}

A loop-erased random walk (or LERW for short) is the process that one obtains when erasing the loops chronologically as they appear on a simple random walk trajectory. More precisely, fix a vertex $x$ in a locally finite graph $G$ and a subset $U$ of vertices, and suppose that the hitting time $H_U < \infty$, $\P_x$-almost surely, where $\P_x$ denotes the law of simple random walk $(X_n)_{n\ge 0}$ on $G$ starting from $x$, and $H_U=\inf \{ n \ge 0: X_n \in U\}$ is the hitting time of $U$ for that walk.

\begin{definition} \label{D:LERW}
A \textbf{loop-erased random walk} from $x$ to $U$ is the process obtained from $(X_n)_{0\le n \le H_U}$ by chronologically erasing the loops from $X$. More precisely, the loop-erasure $\beta = (\beta_0, \ldots, \beta_\ell)$ is defined inductively as follows: $\beta_0 = x$. If $\beta_n \in U$ then $n = \ell$, else $\beta_{n+1} = X_L$, where $L = 1+ \max \{ m \le H_U: X_m = \beta_n\}.$
\end{definition}

Somewhat remarkably, the loops can also be erased antichronologically, and this does not change the resulting distribution:

\begin{lemma}\label{L:lerwreversal}
Let $X$ be a random walk starting from $x$, stopped at the time $ H = H_U$ when it hits $U$. Let $\beta$ denote the loop-erasure of $X$, and let $\gamma$ denote the loop-erasure of the time reversal $\hat X = ( X_H, X_{H-1}, \ldots, X_0)$. Then $\gamma$ has the same law as the time reversal of $\beta$. 
%
%The same result holds when we condition on $X_H = u \in U$. 
\end{lemma}  

This allows us to speak unambiguously of the law of the loop-erasure of a (portion of) simple random walk. 

There is a well known and very deep connection between uniform spanning trees and loop-erased random walks, which was discovered by Wilson \cite{Wilson_alg} (see also Propp and Wilson \cite{ProppWilson}), and may be used to  efficiently simulate such trees. This relation is known as \textbf{Wilson's algorithm}; see Chapter 4 of \cite{LyonsPeres} for a thorough discussion. This relation extends to weighted graphs (provided that one replaces the uniform distribution on the set of spanning trees by a natural Gibbs distribution).
However for simplicity we will not discuss it here and continue to assume that our graphs are unweighted. 
(As we will see in the proof below, this connection is more naturally expressed when we think of our trees as being oriented towards a designated vertex called the root; beware however that this is in contrast with our definition of spanning trees in the previous sections as being unoriented.) 

Here we will only need the following result, which may be seen as a straightforward consequence of Wilson's algorithm, but which was first discovered by Pemantle \cite{Pemantle} (prior to \cite{Wilson_alg}). We state and prove it here, since the proof is short and rather beautiful.

\begin{theorem}
  \label{T:pemantle}
  Let $G$ be an arbitrary finite (connected, unoriented) graph $G$, and let $T$ be a uniform spanning tree on $G$, that is, a uniformly chosen subset of unoriented edges which is acyclic and spanning.
  Let $x,y$ be any fixed vertices in $G$. Then the unique branch of $T$ between $x$ and $y$ has the same distribution as (the trace of) a loop-erased random walk from $x$ run until it hits $y$.
\end{theorem}

 \begin{proof}[Sketch of proof.]
 A \emph{rooted spanning tree} is just a pair $(t,x)$ where $t$ is an (unrooted) spanning tree and $x$ a fixed vertex of $G$. Alternatively we can view the rooted tree $(t,x)$ as an oriented tree, where all the edges of $t$ are oriented towards the root $x$. This is also known as an arborescence, that is, a rooted, directed acyclic graph in which each vertex except the root has a unique oriented edge leading out of it. We will sketch the proof of the theorem for the measure on rooted trees given by 
 \begin{equation}\label{eq:lawtree}
 \P(  (T,X) =  (t,x)) \propto \pi(x) = \deg(x)
 \end{equation}
    and by describing the law of the branch between $y$ and $X$, conditional on $\{X = x\}$; the stated result then follows easily.

   For a possibly infinite path $\gamma = (\gamma_0, \gamma_1, \ldots)$ on the vertices $V$ of the graph, let $T(\gamma)$ be the set of oriented edges of the form $(\gamma_{H(w) -1}, w)$ where $w \neq \gamma_0$, and $H(w) = \inf\{n\ge 0: \gamma_n = w\}$ is the first hitting time of $w$ by the path $\gamma$. In other words, in $T(\gamma)$ we keep the edge $(\gamma_n, \gamma_{n+1})$ if and only if $\gamma_{n+1}$ has not been previously visited. It is obvious that this generates an acyclic graph, and if the path visits every vertex (which will be almost surely the case when $\gamma$ has the law of a random walk, by recurrence of $G$) then $T(\gamma)$ is a spanning tree rooted at $o = \gamma_0$.
Note also that in $T(\gamma)$, the unique branch connecting the root $o$ and a given vertex $w$ can be described by chronologically erasing the loops of the time reversed path $(\gamma_{H(w)}, \gamma_{H(w) -1 }, \ldots, \gamma_0)$. 
   
   %rooted tree obtained by retaining only the edges $(\gamma_j, \gamma_{j+1})$ which do not close a loop: that is, keep this edge if and only if there is no $i < j$ such that $\gamma_i = \gamma_j$ (the root is the starting point $\gamma_0$ of the path). 
   
Appealing to \cref{L:lerwreversal}, in order to conclude the proof of the theorem it suffices to show that when $\gamma_0$ is chosen according to the stationary distribution of the graph $G$ (that is, proportionally to the degree of a vertex) and $\gamma=(\gamma_0, \gamma_1, \ldots, )$ is a random walk starting from $\gamma_0$, then the law of $(T(\gamma), \gamma_0)$ is the one in \eqref{eq:lawtree}.
 
 To do this, suppose that $(X_n)_{n \in \Z}$ is a bi-infinite stationary random walk on $G$ (constructed, for example, using Kolmogorov's extension theorem), so that $X_0$ is distributed according to its equilibrium measure $\pi$, and let $\gamma^{(n)} = (X_{-n}, X_{-n+1}, \ldots)$ be the path started from $\hat X_n := X_{-n}$. Then the claim is that $(T(\gamma^{(n)}), \hat X_{n})$ defines a certain Markov chain on rooted spanning trees. Indeed a straightforward computation using the definition of conditional probability and the fact that $\pi$ is a reversible measure for $X$ show that $\hat X_n $ is itself a Markov chain, with transition probabilities
 $$
 \hat p(x,y) = p(y,x) \frac{\pi(y)}{\pi(x)}
 $$
 where $p$ are the original transitions of the simple random walk on $G$. In particular, given $(T(\gamma^{(n)}) , \hat X_n)=(t,x)$, the probability that $\hat X_{n+1} = y$ is equal to $\hat p(x,y)$. Moreover, given $(T(\gamma^{(n)}) , \hat X_n, \hat X_{n+1})=(t,x,y)$, $T(\gamma^{(n+1)})$ is obtained deterministically from $t$ by adding the edge $e = (x,y)$ (which creates a cycle) and removing from $t$ the unique outgoing edge from $y$. The next state of the chain corresponds to $(T(\gamma^{(n+1)}), y)$. (This corresponds to what Lyons and Peres \cite{LyonsPeres} describe as the forward procedure, but applied to the time reversed chain $\hat X$). 
 
 It is not hard to see that this Markov chain on the space of rooted trees is irreducible. A calculation (involving the so called ``backward procedure'', see Section 4.4 in \cite{LyonsPeres}) shows further that the unique stationary distribution of this chain is proportional to 
 $$
 \P ( ( T,X) = (t,x)) \propto \psi((t,x)) :=  \prod_{\ve } p( \ve), 
 $$ 
 where the product is over all the oriented edges $\ve$ of the rooted tree $(t,x)$. (We warn the reader however that $\psi$ is not in general a reversible measure for this chain.) Furthermore, it is a classical fact, known as the \textbf{Markov chain tree theorem}, that for general weighted graphs, $\psi$ is in fact itself proportional to the invariant measure of the associated random walk. In the case which occupies us here where the graph is unweighted (that is, all edges have unit weight), this is particularly easy to see: indeed, since $t$ is spanning and every vertex except the root has exactly one outgoing oriented edge, 
 $$
 \prod_{\ve} p(\ve) = \prod_{v\ne x} \frac{1}{\deg(v)} =\frac{\deg(x)}{\prod_{v\in G} \pi(v)} \propto \pi(x);
 $$
   that is, the invariant distribution is given exactly by \eqref{eq:lawtree}. On the other hand, since  $(\gamma^{(n)},\hat{X}_n)_{n\in \Z}$ is stationary, so must be $(T(\gamma^{(n)}),\hat{X}_n)_{n\in \Z}$ (because $T(\gamma^{(n)})$ is a deterministic function of $\gamma^{(n)}$). In particular, $(T(\gamma^{(0)}),\hat{X}_0)$ has distribution \eqref{eq:lawtree}, as required.
    \end{proof}

\begin{rmk}
The fact that $T(\gamma^{(0)})$ has the law \eqref{eq:lawtree} is closely related to the algorithm of Aldous \cite{AldousUST} and Broder \cite{Broder} for sampling a uniform spanning tree. As mentioned in \cite{LyonsPeres}, both authors credit Persi Diaconis for discussions. This algorithm was initially used in the study of the Uniform Spanning Tree (notably by Pemantle \cite{Pemantle}) before Wilson's algorithm (\cite{Wilson_alg, ProppWilson}) became available. 

However, we note that in fact, Wilson's algorithm to generate a full UST is a simple consequence of Theorem \ref{T:pemantle}: indeed, to generate the tree $T$ we can first sample the branch connecting a fixed vertex $x$ to the boundary using a loop-erased random walk by Theorem \ref{T:pemantle}. The conditional law of the rest of the tree is then a uniform spanning tree on a modified graph where this branch has been wired to make a single vertex and become part of the boundary. Applying Theorem \ref{T:pemantle} recursively in this manner then gives Wilson's algorithm. 

Of course, direct (and relatively short) proofs of this algorithm exist. See, in particular, \cite[Chapter 4.1]{LyonsPeres} for a proof close to the original spirit of \cite{Wilson_alg}, \cite{LawlerLimic} for a proof using loop measures, and finally see \cite[Chapter 2.1]{WWnotes} for a proof based on the Green function and the discrete Laplacian. 
\end{rmk}

We may deduce from \cref{T:ust} the following result about the loop-erased random walk.
\ind{Wilson's algorithm}

\begin{theorem}\label{T:lerw}
  Let {$(M_n, T_n)$} be chosen as in \eqref{maptree} and let $x, y$ be two vertices chosen independently and uniformly on $M_n$. Let $(\Lambda_k)_{0 \le k \le \xi_n}$ be a LERW starting from $x$, run until the random time $\xi_n$ when it first hits $y$. Then
  $$
  \frac{\xi_n}{\sqrt{n}} \to \xi_\infty
  $$
  in distribution, where $\xi_\infty$ is a random variable that has a non-degenerate distribution (in the sense that $\xi_\infty \in (0, \infty)$ almost surely.
\end{theorem}

\begin{proof}
Let $(X_k,Y_k)_{1\le k \le 2n}$ be the pair of excursions which describes the map $(M_n, T_n)$. Then note that $\xi_n$ may be identified with the ``tree distance'' $X(J_1)+X(J_2)-2\min_{j\in [J_1,J_2]}X(j)$ where $J_1,J_2$ are uniformly (and independently) chosen between $1$ and $2n$. As a consequence, Theorem \ref{T:lerw} holds with $\xi_\infty=e(U_1)+e(U_2)-\inf_{u\in [U_1,U_2]} e(u)$,  where $e$ is a Brownian excursion $e$ and $U_1,U_2$ are chosen uniformly and independently from $(0,1)$.
\end{proof}

\begin{rmk}
In fact, as was already observed by Aldous \cite{Aldous}, the continuum random tree is invariant ``under rerooting'', that is, moving to the root to a uniformly chosen position. As a consequence, the law of the random variable $\xi_\infty$ above may be more simply written as $e(U)$, where $U$ is a uniform random variable on $(0,1)$. In fact, as noted in \cite{Aldous}, this can be derived directly from a simple path transformation of the Brownian excursion. See also \cite{DuquesneLeGall_invariance} for a discussion in the more general context of L\'evy trees. 
\end{rmk}

\ind{Loop-Erased Random Walk}
\ind{Critical exponents!Loop-Erased Random Walk}

\paragraph{Scaling exponent of LERW.} We now explain how the above result can be used to compute an exponent for the loop-erased random walk. Let $\Lambda = \{\Lambda_0, \ldots, \Lambda_{\xi_n}\}$ denote the loop-erasure of a random walk on $M_n$, run from a uniformly chosen vertex $x$ until the hitting time of another uniformly chosen vertex $y$, as above. Then $\Lambda$ may be viewed as an independent random ``fractal'' set on $M_n$, whose size is $|\Lambda| = \xi_n = n^{1/2 + o(1)}$ by Theorem \ref{T:lerw}. Since $M_n$ has $n^{1+ o(1)}$ vertices (indeed it has by definition $n$ edges, and the degree distribution of a given vertex is known to be very concentrated), this means that $\Lambda$ has a \textbf{quantum scaling exponent} given by
$$
\Delta = 1/2
$$
(recall our discussion from Section \ref{S:exponents}). We can therefore (at least informally) use the \textbf{KPZ relation} \ind{KPZ relation} to compute the equivalent exponent for a loop-erased random walk on the square lattice. To do so, we must first find the correct value of $\gamma$: the constant in front of the GFF which describes the scaling limit of the conformally embedded planar map $M_n$. This is given by the relation \eqref{q_gamma} when $q=0$ (which, as explained at the beginning of this section, indeed corresponds to the uniform spanning tree weighted map model of \eqref{maptree}). Plugging $q=0$ yields $$\gamma = \sqrt{2}.$$
Note that this is consistent with the conjecture (known to be true on the square lattice by results of \cite{LSW})  that the interface separating a uniform spanning tree from its dual, converges in the scaling limit to an SLE curve with parameter $\kappa'= 8$.

Therefore, the \textbf{Euclidean scaling exponent} $x$ of the loop-erased random walk should satisfy
$$
x = \frac{\gamma^2}{4} \Delta^2 + (1 - \frac{\gamma^2}{4}) \Delta = 3/8.
$$
In particular, we conclude that in the scaling limit, a loop-erased random walk on the square lattice has dimension
$$
d_{\text{Hausdorff}} = 2- 2x = 5/4.
$$
This is in accordance with Beffara's formula \cite{Beffara} for the dimension of SLE: indeed, in the scaling limit, LERW is known to converge to an SLE$_\kappa$ curve with $\kappa = 2$. This is closely related to the above mentioned scaling limit result for the UST, due to Lawler, Schramm and Werner \cite{LSW}, and is also proved in \cite{LSW}.
Beffara's result \cite{Beffara} states that the Hausdorff dimension of SLE$_\kappa$ is  $(1+ \kappa/ 8 ) \wedge 2$. In the case $\kappa=2$ this is exactly $5/4$, as above.

In fact, this exponent for LERW had earlier been derived by Kenyon in a remarkable paper \cite{KenyonLERW}, building on his earlier work on the dimer model and the Gaussian free field \cite{Kenyon_GFF}.

\subsection{Sheffield's bijection in the general case}
\label{S:Sheffield_bij}

We now describe the situation when $\bar{\bm}_n\in \cM_n$ but the collection of edges $\bt_n$ is arbitrary (that is, not necessarily a tree), which is more delicate. {Note that in the case of spanning trees there was only one loop present, but now there will generally be more than one.} These loops are \textbf{densely packed} in the sense that every triangle is part of some loop, as illustrated in Figure \ref{fig:clusters}. Indeed, each triangle consists of an edge of some type and a vertex of the opposite type, so must contain a loop separating the two associated clusters. In this case we will see that we can still define a canonical space-filling interface (that is, a curve which visits every single triangle exactly once). We will now describe this curve (see also Figure \ref{fig:bij}).

Recall that $L_0$ is the loop containing the root triangle of the map $\bar{\bm}_n$, oriented parallel to the orientation of the root edge of $\bm_n$. We view $L_0$ as an oriented collection of adjacent triangles (the triangles traversed by the loop). In general, $L_0$ does not cover every triangle of $\bar{\bm}_n$, and we may consider the connected components $C_1, \ldots, C_k$ which are obtained by removing all the triangles of $L_0$. Note that $L_0$ is adjacent to each of these components, in the sense that for each $1 \le i \le k$, it contains a triangle that is opposite a triangle in $C_i$. For each $i$, let $T_i$ be the \emph{last} (with respect to the orientation of the loop and its starting point) triangle that is adjacent to $C_i$. The triangle opposite $T_i$ is in $C_i$ and together they form a quadrangle. In order to explore all of the map and not just $L_0$, we will first modify the map by \emph{flipping} the diagonal of this quadrangle, for every $1 \le i \le k$. It can be seen that having done so, we have reduced the number of loops on the map by exactly $k$ (each such flipping has the effect of merging two loops). We may then iterate this procedure until there is only a single loop left, the loop $L_0$ (which now fills the whole map). This loop separates primal and dual clusters of the modified map, in the sense that it has only primal clusters on one side, and dual clusters on the other (we will see below that these clusters are in fact spanning trees).

%The idea is that the space-filing path starts, as before, by going around the component of the root edge, that is exploring the loop $L_0$ of the root.  However, {this will no longer explore the whole map.} %we also need to explore the rest of the map.
%To {explore the rest of the map, we} consider the last triangle which shares an edge with the inside of the loop (so, roughly -- but not necessarily exactly -- at the time when we are about to close the loop $L_0$). This triangle is part of a quadrangle where only one of the diagonals is drawn (either primal or dual diagonal). Then we change the map by flipping the diagonal of this quadrangle. This allows us to continue the exploration, essentially the loop we were about to close.

%This can be iterated; a bit of thought shows that this results in a space-filling path which visit every quadrangle exactly twice, going around some virtual tree.
So we now have a canonical space-filling path which allows us to explore the map as in Section \ref{sec:bij_tree}. As before, we can describe the type of triangles we see in this exploration using the symbols $\ah, \ac, \aH, \aC$. When we explore a triangle corresponding a flipped quadrangle for the first time, we record its type (either $\ah, \ac$) according to its type after having flipping the edge. However, when we visit its opposite triangle we record the fact that this is a special edge (which must be flipped to recover the original map) by the symbol $\aF$. The letter $\aF$ stands for ``flexible'' or ``freshest'' order. (We will see below a more precise interpretation in terms of queues, or hamburgers and cheeseburgers.) In this way, we may associate to the decorated map $(\bm_n, \bt_n)$ a list $w$ of $2n$ symbols $w=(X_i)_{1\le i \le 2n}$ taking values in the alphabet $\Theta = \{ \ah, \ac, \aH, \aC, \aF\}$.

\begin{figure}
\centering{
%%\includegraphics[width=0.6\textwidth]{example}
%1.\includegraphics[scale=.4]{map_closed} \quad \quad
%2. \includegraphics[scale=.4]{map_dual}
%}
%
%\vspace{.5cm}
%\centering{
%3. \includegraphics[scale=.4]{map_tutte}\quad \quad
\textbf{a.}\includegraphics[scale=.4]{map_loops} \quad \quad
\textbf{b.} \includegraphics[scale=.4]{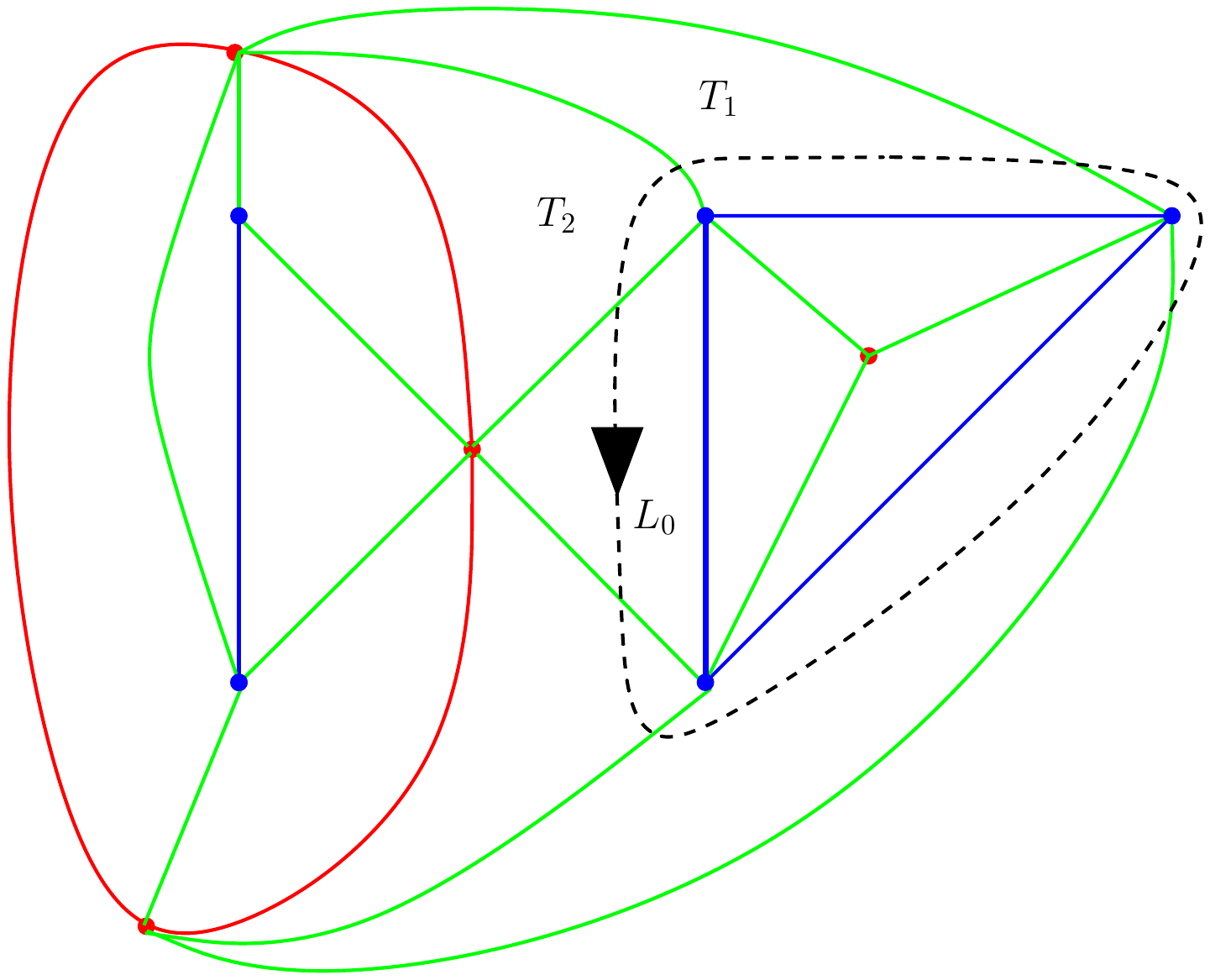}}

\vspace{.5cm}
\centering{
\textbf{c.}\includegraphics[scale=.4]{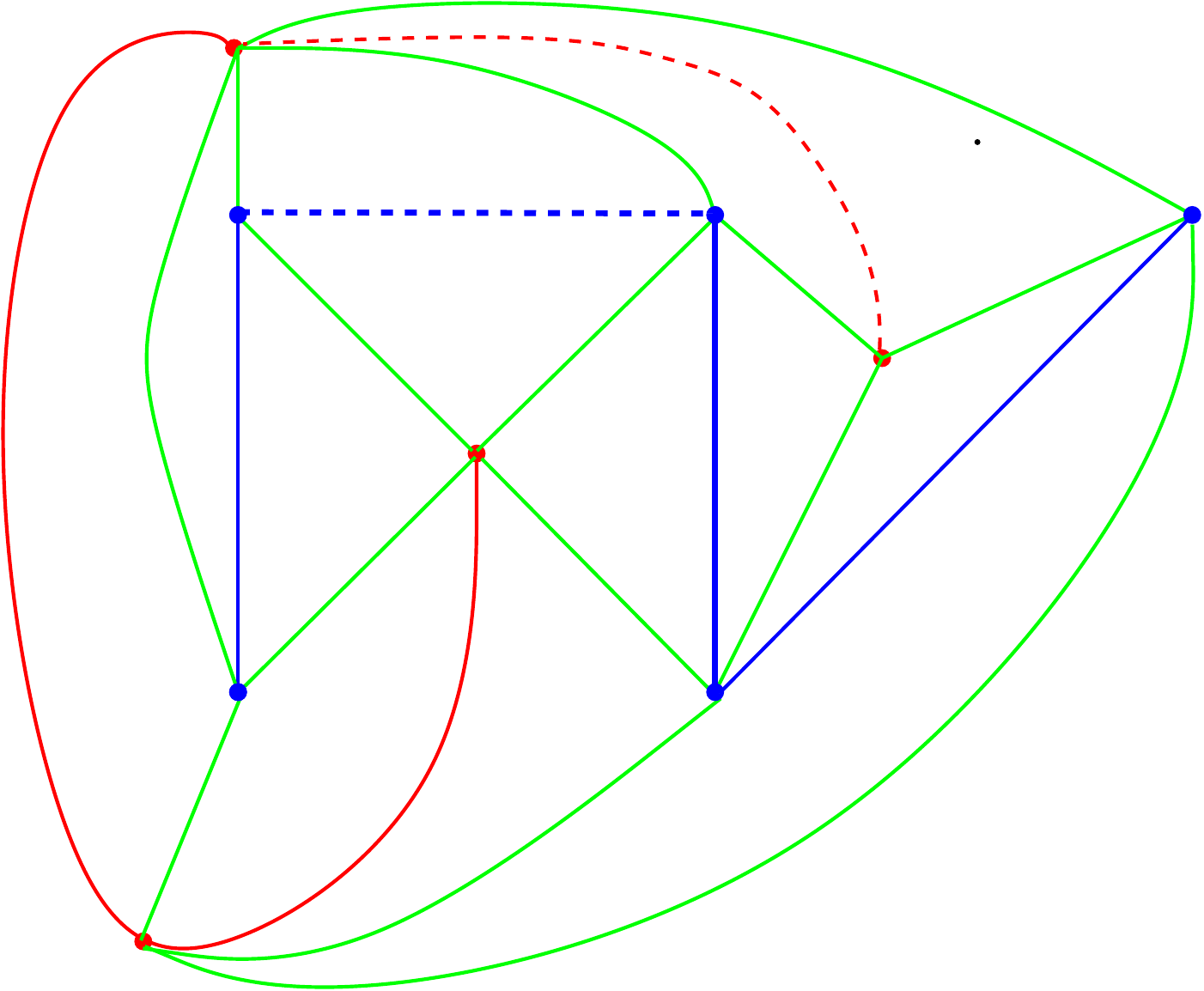}\quad \quad
\textbf{d.} \includegraphics[scale=.4]{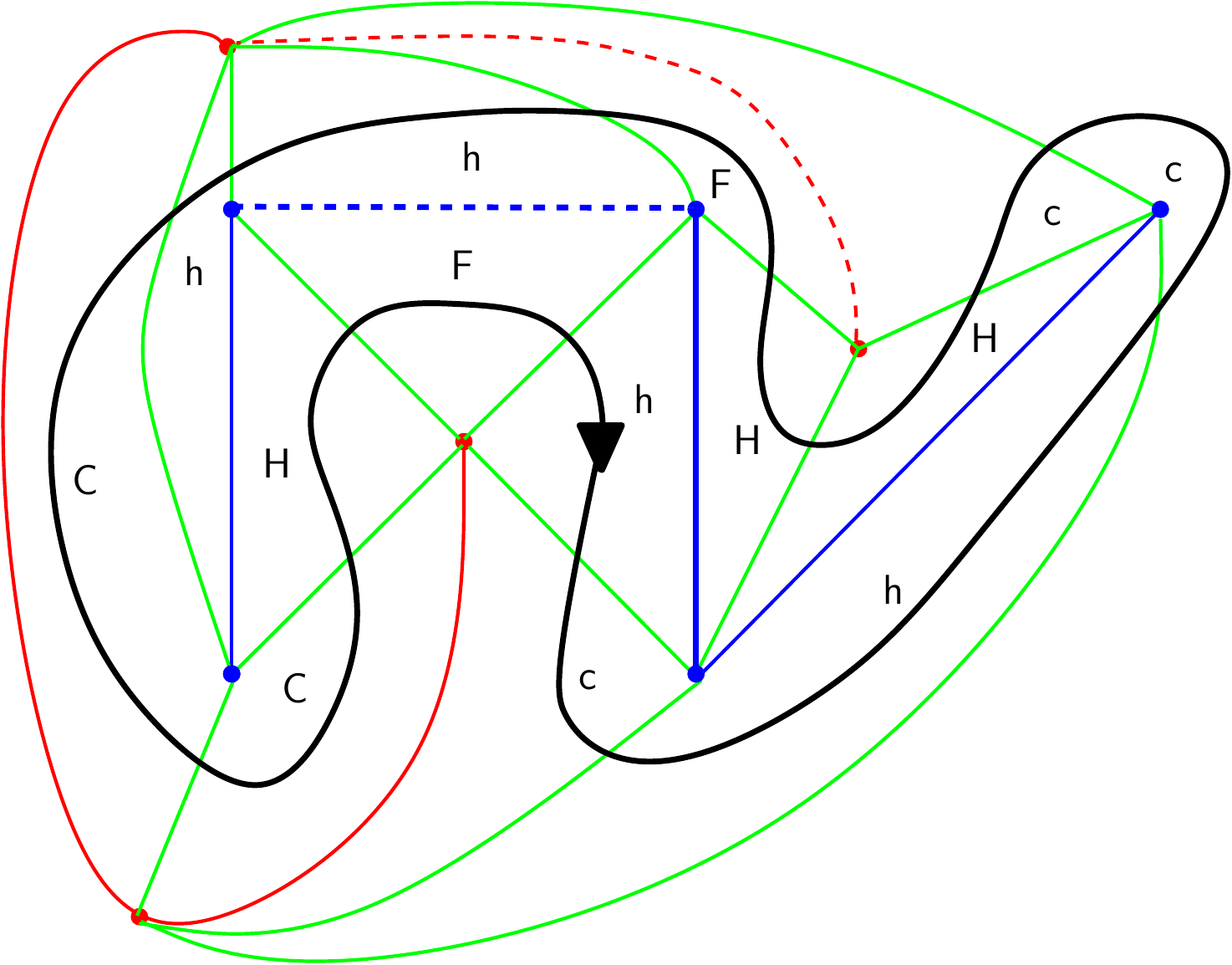}}
\caption{{Generating a word from a  decorated map in the general case. \textbf{a.} The decorated map is as in \cref{fig:clusters}, with the (oriented) root loop $L_0$. \textbf{b.} The complement of $L_0$ consists of two components, $C_1$ and $C_2$.
%In this example, $C_1$ is the collection of three triangles surrounded by the blue loop towards the right of the picture, and $C_2$ is the collection of five triangles surrounded by the red loop towards the left.
$T_1$ and $T_2$ are the \emph{last} triangles visited by the loop $L_0$ that share an edge with a triangle in $C_1$ and $C_2$ respectively. \textbf{c.} We flip the diagonals of the quadrangles associated with $T_1$ and $T_2$.
%The triangles $A_1$ and $A_2$ are each associated with a quadrilateral, and we flip the diagonal in each of these quadrilaterals. So here, the quadrilateral of $A_1$ changes from primal to dual, and the quadrilateral of $A_2$ changes from dual to primal.
\textbf{d.} We obtain a single space-filling loop (drawn in black).
%Once the edges have been flipped the blue edges form a tree and there is a single loop (marked in black) that goes around this tree anti-clockwise as in \cref{sec:bij_tree}.
To this path we can again associate a word in  \{\ah, \ac, \aH, \aC\}. However, we also record the second visit to a flipped quadrangle by replacing the symbol $\aC$ or $\aH$ by the symbol \aF. The word here is thus \ah \ac \ah \ac \ac \aH \aH \aF \ah \ah \aC\aC \aH \aF.  } Note the non-obvious fact that after flipping, the primal and dual clusters have become trees.}
%
 %
 % with full lines for closed edges and dotted lines
 % for open edges. The dual map is in red, with only the open dual dges
 % drawn. The green edges form the associated quadrangulation along
 % with the primal and dual edges form the ($\bar m $ in \cref{sec:critical-fk-model}). The percolation clusters and dual clusters are separated by loops, in dashed purple, that cross every quandrangulation edge once.
\label{fig:bij}
\end{figure}

We will see below the properties of this word (essentially, it reduces to $\emptyset$ with the appropriate definition of reduction when there is an $\aF$) and that the map from $(\bm_n, \bt_n)$ to $w$, subject to this constraint, is a bijection. For now, we make the important observation
 that each loop corresponds to a unique symbol $\aF$, except for the loop through the root.

\paragraph{Inventory accumulation.}
Recall that we can interpret an element in $\{\ah,\ac,\aH,\aC\}^{2n}$ as a
last in, first out inventory accumulation process in a burger factory with two types of product:
hamburgers and cheeseburgers. Think of a sequence of events,
occurring once per unit time, in
which either a burger is produced (either ham or cheese) or there is
an order of a burger (either ham or cheese). The burgers are put in a
single \textbf{stack} and every time there is an order of a certain
type of burger, the freshest burger in the stack of the corresponding
type is removed. The symbol $\ah $ (resp. $\ac$) corresponds to a ham
(resp. cheese) burger production and the symbol $\aH$ (resp. $\aC$)
corresponds to a ham (resp. cheese) burger order.

 The inventory interpretation of the symbol $\aF$ is the following:
this corresponds to a customer
demanding the freshest or the topmost burger in the stack, irrespective
of the type. In particular, whether an $\aF$ symbol corresponds to a
hamburger or a cheeseburger order depends on the topmost burger type
at the time of the order. Thus overall, we can think of the inventory
process as a sequence of symbols in $\Theta$ with the following
reduction rules
\begin{itemize}
\item $\overline{\ac \aC }= \overline{\ac \aF }= \overline{\ah \aH }= \overline{\ah \aF} = \emptyset$,
\item $\overline{\ac\aH }= \overline{\aH \ac}$ and $\overline{\ah \aC }= \overline{    \aC \ah}$.
\end{itemize}
Given a sequence of symbols $w$, we denote by $\bar w$ the reduced
word formed via the above reduction rule.

\paragraph{Reversing the construction.} Given a sequence $w$ of symbols from $\Theta$, such that $\bar w =
\emptyset$, we can construct a decorated map $(\bm_n,\bt_n)$
as follows. First, we convert all the $\aF$ symbols to either an $\aH$ or
a $\aC$ symbol depending on its order type. Then, we construct a spanning
tree decorated map as described in Section \ref{sec:bij_tree}
(see in particular \cref{fig:map_from_word}). The condition $\bar w = \emptyset$ ensures
that we can do this. To obtain the original loop decorated map, we simply flip the type of
every quadrangle which has one of the triangles corresponding to an $\aF$
symbol. That is, if a quadrangle formed by primal triangles has one of
its triangles coming from an $\aF$ symbol, then we replace the primal map
edge in that quadrangle by the corresponding dual edge and vice
versa. The interface is now divided into several loops (and the number
of loops is exactly one more than the number of $\aF$ symbols). In particular:
\begin{theorem}[Sheffield, \cite{Sheffield_burger}]\label{T:sheffieldbij}
The map $(\bm_n, \bt_n) \mapsto w$ (subject to $\bar{w} = \emptyset$) is a bijection.
\end{theorem}

\paragraph{Two canonical spanning trees.} It is not obvious but true that after flipping, the corresponding primal and dual decorations of the map have become two mutually dual spanning trees. One way to see this is as follows: observe that after flipping, we have (as already argued) a single space-filling loop which separates primal and dual clusters of the resulting modified map. These clusters are of course spanning, and they cannot contain non-trivial cycles, else the loop would either not be space-filling or consist of multiple loops. Therefore, we can again think of $M_n$ as a gluing of two spanning trees, which are glued along the space-filling path (that is, along their contour functions). Again, this perspective is a crucial intuition which guides the \textbf{mating of trees approach} to Liouville quantum gravity \cite{DuplantierMillerSheffield}.
We will survey this later on (see in particular Section \ref{SS:MOTdescription}). \ind{Mating of trees}

\paragraph{Generating FK-weighted maps.} A remarkable consequence of Theorem \ref{T:sheffieldbij} is the following simple way of generating a random planar map from the FK model \eqref{FK}. Fix $p \in [0,1/2)$, which will be suitably chosen (as a function of $q$) below in \eqref{pq}. Let $(X_1,\ldots, X_{2n})\in (\Theta)^{2n}$ be i.i.d.\ with the
following law
\begin{equation}
  \label{eq:3}
  \P( \ac) = \P( \ah ) = \frac14, \P( \aC ) = \P ( \aH ) = \frac{1-p}{4}, \P( \aF ) = \frac p2,
\end{equation}
conditioned on $\overline{X_1,\ldots, X_{2n}} = \emptyset$.
\ind{FK model}

Let {$(M_n, T_n)$} be the random associated decorated map (via the bijection described above). Then observe that since $n$ hamburgers and cheeseburgers must be produced, and since $\#\aH + \#\aC = n - \# \aF$,
\begin{align}
  \label{eq:2}
  {\P((M_n,T_n)=(\bm_n,\bt_n))} & = \left(\frac14\right)^n  \left( \frac{1-p}4\right)^{\#\aH + \#\aC} \left(\frac{p}{2}\right)^{\#\aF}\nonumber  \\
  & \propto \left(\frac{2p}{1-p}\right)^{\# \aF} = \left(\frac{2p}{1-p}\right)^{\# \ell(\bm_n,\bt_n)-1}
\end{align}
Thus we see that {$(M_n,T_n)$} is a realisation of the critical FK-weighted cluster random map model with
\begin{equation}\label{pq}
\sqrt{q}
= \frac{2p}{(1-p)}.
\end{equation}
 Notice that $p \in [0,1/2)$ corresponds to
$q =
[0,4)$. From now on we fix the value of $p$ and $q$ in this regime. Recall that $q=4$ is believed to be a critical value for many properties of the map; indeed later on we will later show that a  phase transition occurs at $p=1/2$ ($q=4$) for the geometry of the map. Intuitively, it is perhaps not surprising that the value $p=1/2$ marks a distinction from the point of view of inventory accumulation.

\subsection{Infinite volume limit}\label{sec:local}
The following theorem due to Sheffield \cite{Sheffield_burger}, and made more precise later by Chen
\cite{chen}, shows that the decorated map $(M_n,T_n)$ has a local limit
as $n \to \infty$ in the local topology. Roughly two {(decorated)} maps
are close in the local topology if the finite maps (and their decorations) near a large
neighbourhood of the root are isomorphic as decorated maps.

\begin{theorem}[\cite{Sheffield_burger,chen}]
  \label{thm:Local}
Fix $p \in [0,\frac12)$. We have
\[
(M_n,T_n) \xrightarrow[n \to \infty]{(d)}(M,T)
\]
{with respect to the local topology, where $(M,T)$ is the unique infinite decorated map associated with a bi-infinite i.i.d.\ sequence of symbols $(X_n)_{n \in \Z}$ having law \eqref{eq:3}.
}
\end{theorem}

\begin{proof}[Sketch of proof]
We now give the idea behind the proof of \cref{thm:Local}. Let
$X_1,\ldots,X_{2n}$ be i.i.d.\ with law given by \eqref{eq:3},
and denote by $E_{2n}$ the event that
$\overline{X_1\ldots X_{2n}} = \emptyset$.

A key step  is to show the following.
\begin{lemma}[\cite{Sheffield_burger,FKstory2}]\label{L:reduction}
Let $X_1, \ldots, X_{2n}$ be i.i.d. with law  \eqref{eq:3}. Then $\P( E_{2n} ) $ decays subexponentially in $n$, that is,
$ \log \P( E_{2n}) / n \to 0$ as $n \to \infty$.
\end{lemma}

We will not prove this statement (although we will later come back to it and explain it informally). Instead we explain how Theorem \ref{thm:Local} follows.

Notice that uniformly selecting a symbol $1\le I \le 2n$ corresponds to selecting a
uniform triangle in $(\bar M_n, T_n)$, which in turn corresponds to %a unique root triangle in $\bar M_n$ (or equivalently, a unique
a unique oriented edge in $M_n$ (recall that $\bar M_n$ denotes the refinement map associated to $M_n$). Because of invariance of the decorated map $(M_n, T_n)$ under re-rooting, we claim that it suffices to check the convergence in distribution of a large neighbourhood of the triangle corresponding to $X_I$ in $\bar M_n$.

Let $r>0$. We will first show that for any fixed word $w$ of length $2r+1$ in the alphabet $\Theta$,
\begin{equation}\label{localroot}
\P ( X_{I-r}\ldots X_{I+r} = w \big| E_{2n}) \to \P(w) : = \P( X_{-r} \ldots X_r = w),
\end{equation}
where on the left hand side the addition of indices has to be interpreted cyclically within $\{1, \ldots, 2n\}$, and on the right hand side, $(X_n)_{n\in \mathbb{Z}}$ is the random bi-infinite word whose law is described in \cref{thm:Local}.

To see \eqref{localroot}, observe that the conditional probability on the left hand side (conditionally given the entire sequence $X= (X_{1}, \ldots X_{2n})$ satisfying $E_{2n}$, and averaging just over $I$),
is equal to $f +o(1)$ as $n \to \infty$, where $f$ is the fraction of occurrences of $w$ in $X$, that is, $f = (2n)^{-1} \sum_{i=r+1}^{2n-2r-1} 1_{\{ X_{i-r} , \ldots X_{ i+r} = w\}}$, and the $o(1)$ term is uniform, accounting for boundary effects. Hence it suffices to check that $\E( f | E_{2n}) \to \P( w)$. To do this, for arbitrary $\eps>0$ we define $A_n = \{ | f - \P(w)| \le  \eps\}$, and write
\begin{align*}
\E(f | E_{2n}) & =  \E (f 1_{A_n} | E_{2n} ) + \E (f 1_{A_n^c} | E_{2n} ).
\end{align*}
Now the first term $\E (f 1_{A_n} | E_{2n} )$ is equal to $(\P(w) + O(\eps)) \P( A_n |E_{2n})) $, while the second term satisfies
$$
\E (f 1_{A_n^c} | E_{2n} ) \le \P( A_n^c | E_{2n} ) \le \frac{\P(A_n^c)}{\P(E_{2n})}.
$$
However, $\P( A_n^c) \to 0$ \emph{exponentially fast} as $n \to \infty$, by basic large deviation estimates (Cramer's theorem). This means that $\E (f 1_{A_n^c} | E_{2n} )$ converges to zero by Lemma \ref{L:reduction}, and also %(again exponentially fast).
that $\P(A_{n} | E_{2n}) \to 1$ as $n\to \infty$. We can conclude that $\E(f1_{A_n})$ and therefore $\E(f | E_{2n} )$ converges to  $\P(w)$ as $n \to \infty$, which proves \eqref{localroot}.

To conclude the theorem, it remains to show that convergence of the symbols locally around a letter implies local convergence of the maps. This is a consequence of Exercise \ref{Ex:symbolmaps}; see also Figure \ref{fig:map_from_word}.
\end{proof}

One important feature related to Theorem \ref{thm:Local} is that every symbol
in the i.i.d.\ sequence $\{X_i\}_{i \in \Z}$ has an almost sure unique
\textbf{match}, meaning that every burger order is fulfilled (it corresponds to a burger that was produced at a finite time before), and every burger that is produced is
consumed at some finite later time, both with probability $1$; see \cite[Proposition 3.2]{Sheffield_burger}. In the language of maps, this is equivalent to saying that the map $M$ has no edge ``to infinity''. For future reference, let $\varphi(i)$ denote the match of
the $i$th symbol. Notice that $\varphi:\Z \mapsto \Z$ defines an
involution on the integers.

\subsection{Scaling limit of the two canonical trees}

We now state (without proof) one of the main results of Sheffield \cite{Sheffield_burger}, which gives a scaling limit result for the geometry of the infinite volume map $(M,T)$ defined in \cref{thm:Local}. Recall that $(M,T)$ is completely described by a doubly infinite sequence $(X_n)_{n \in \Z}$ of i.i.d symbols in the alphabet $\Theta$, having law \eqref{eq:3}. Associated to such a sequence we can define two processes $(H_n)_{n \in \Z}$ and $(C_n)_{n \in \Z}$ which count the respective number of hamburgers and cheeseburgers present in the queue at time $n\in \Z$ (of course, we convert the flexible $\aF$ orders into their appropriate values to count the numbers of hamburgers and cheeseburgers in the queue at time $n$). These numbers are defined relative to time 0, so $(H_0,C_0)=(0,0)$. In other words, let $\tilde w = (\tilde X_n)_{n \in \Z}$ denote the infinite word obtained from $w = (X_n)_{n\in \Z}$ by transforming the $\aF$ symbols into their actual values $\aH$ and $\aC$, and let 
$$
H_n  = \begin{cases}
\sum_{i=1}^n 1_{\{\tilde X_i = \ah\} } - 1_{\{\tilde X_i = \aH\}}  & \text{ if } n >0\\
\sum_{i=n}^{-1} 1_{\{\tilde X_i = \ah\} } - 1_{\{\tilde X_i = \aH\}} & \text{ if } n < 0;
\end{cases}
$$
similarly for $C_n$.

This scaling limit is most conveniently phrased as a scaling limit for $H= (H_n)_{n \in \Z}$ and $C=(C_n)_{n \in \Z}$ (although the statement of Sheffield \cite{Sheffield_burger} concerns instead $H+C$ and the discrepancy $H-C$). We first state the result and then make some comments on its significance below.

\begin{theorem} \label{T:scalinglimitFK}
  Let $p \in [0,1]$, and let $C, H$ be as above. Then
  $$
  \left( \frac{H_{\lfloor nt\rfloor }}{\sqrt{n}} , \frac{C_{\lfloor nt\rfloor }}{\sqrt{n}}\right)_{-1 \le t \le 1} \to (L_t, R_t)_{-1 \le t \le 1}
  $$
in distribution as $ n \to \infty $ for the topology of uniform convergence,
where $(L_t, R_t)_{t \in \R}$ is a two-sided Brownian motion in $\R^2$,  starting from 0 and having covariance matrix given by
$$
\var (L_t) = \var(R_t) = \frac{1+ \alpha}{4} |t|\quad ; \quad \cov (L_t, R_t) = \frac{1- \alpha}{4} | t|
$$
and
$$
\alpha = \max(1- 2p, 0) .
$$
\end{theorem}

See \cite[Theorem 2.5]{Sheffield_burger} for a proof.
We now make a few important remarks about this statement.

\begin{itemize}
  \item This scaling limit result should be thought of as saying something about the large scale geometry of the map $(M,T)$ or,  equivalently, what it looks like after scaling down by a large factor. However, what this actually means is not \emph{a priori} obvious: really, the theorem only says that the pair of trees converges to correlated (infinite) CRTs. This is a (relatively weak) notion of convergence which has been called \textbf{peanosphere topology}; see more about this in Chapter \ref{S:MOT}. In particular, it does not say anything about convergence of the metric on $M$. \ind{Peanosphere convergence}

  \item Notice that when $p \ge 1/2$ (corresponding to $q \ge 4$ in terms of the FK model \eqref{FK}, see \eqref{pq}) we have $\alpha = 0$, so $L_t = R_t$ for all $t \in \R$. This is because the proportion of $\aF$ orders is large enough that there can be no discrepancy in the scaling limit between hamburgers and cheeseburgers.

  \item However, when $p \le 1/2$ (corresponding to $q \le 4$), the correlation between $L$ and $R$ is non-trivial. When $p=0$ (corresponding to $q=0$) they are actually independent. This last case should be compared with the case of spanning tree weighted maps (Theorem \ref{T:ust}). In general, this suggests that the scaling limit of the map $(M,T)$, if it exists, can be viewed as a gluing of two (possibly correlated) infinite CRTs; meaning that their contour (or alternatively their height) functions are described by a two-sided infinite Brownian motion (rather than a Brownian excursion of duration one). This fact is made rigorous (and will be discussed later on in Section \ref{SS:MOTdescription}) in the \textbf{mating of trees approach to LQG} of  \ind{Continuous Random Tree (CRT)} \cite{DuplantierMillerSheffield}.\ind{Mating of trees} Note in particular that in the case $q \ge 4$, the two corresponding trees are identical, meaning that the map should degenerate to a CRT in the scaling limit. This is in contrast with the case $q<4$, where the limit maps are expected to be homeomorphic to the sphere almost surely.

  \item $H_n, C_n$ also have a geometric interpretation, as the boundary lengths at time $n$ on the left and right hand sides of the space-filling interface (relative to time 0).

\end{itemize}

\begin{figure}
\begin{center}
  \includegraphics[scale=.5]{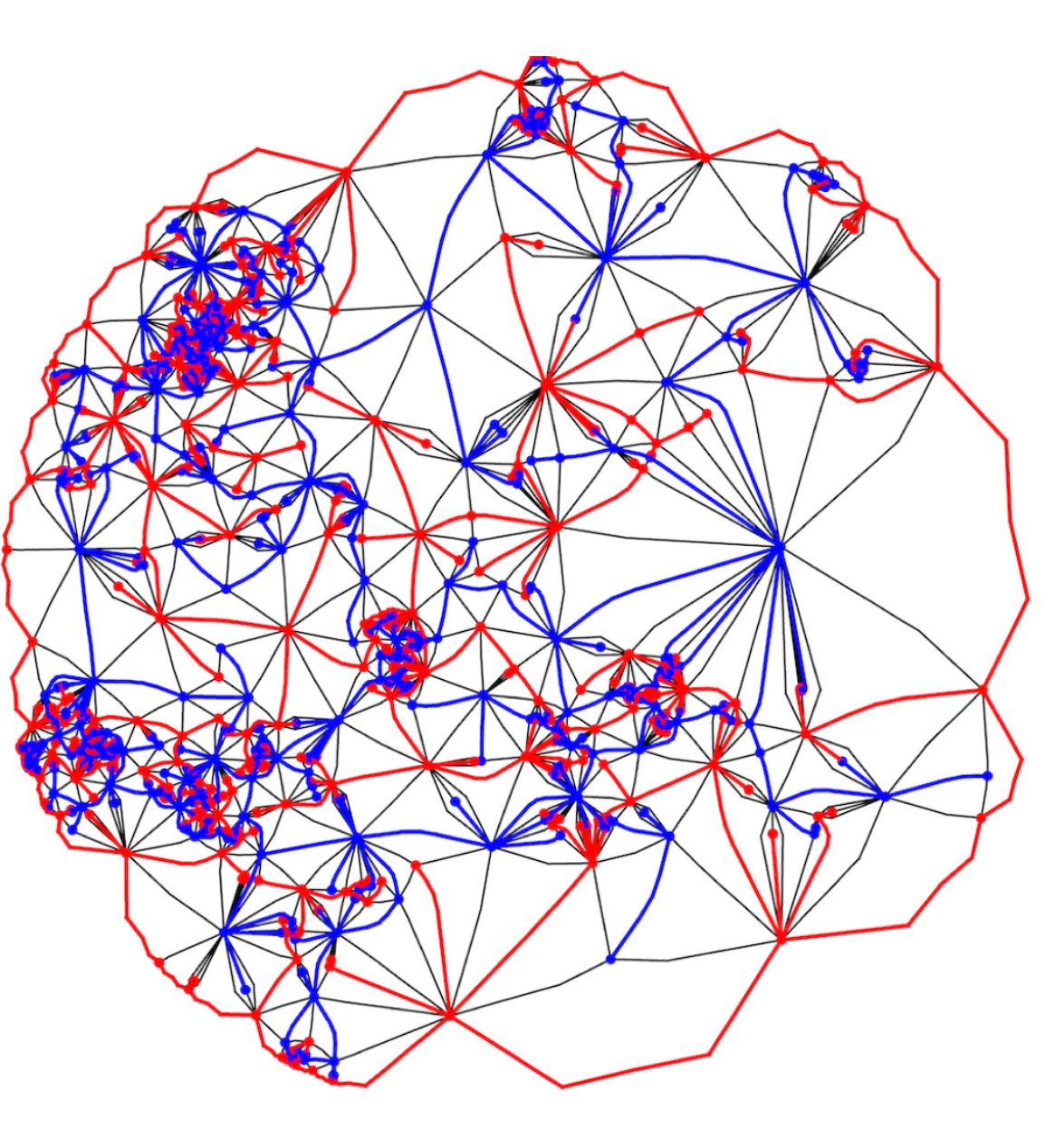}
  \end{center}
  \caption{A random planar map with law \eqref{FK} for $q=1$ (uniform case), generated using Sheffield's bijection of \cref{T:sheffieldbij}. The map has been embedded using circle packing. Shown in blue and red are the primal and dual spanning trees. In the infinite volume limit and then in the scaling limit,  Theorem \ref{T:scalinglimitFK} shows that these trees become correlated infinite CRTs. Simulation by Jason Miller.}\label{F:matingtrees}
\end{figure}

\newpage 

\subsection{Exponents associated with FK-weighted random planar maps}

\ind{FK model}

\dis{In this short section, some critical exponents of random planar maps are computed heuristically. This section can be skipped on a first reading, as none of those results are needed later on. }

It is possible to use Theorem \ref{T:scalinglimitFK} to obtain very  precise information on the geometry of loops on the map $(M,T)$. In particular, it is possible to check that large loops have statistics that coincide with those of CLE$_{\kappa'}$, where the value of $\kappa'$ is related to $q \in (0,4)$ via  \eqref{q_kappa}, thereby giving credence to the general conjectures formulated in Section \ref{sec:RPM_LQG}. This line of reasoning has been pursued very successfully in a string of papers by Gwynne, Mao and Sun \cite{FKstory1, FKstory2, FKstory3}. We will present here a slightly less precise (but easier to state) result proved in \cite{BLR}. Let $(X_i)_{i\in \Z}$ denote the symbols encoding $(M,T)$, and let us condition on the event $X_0 = \aF$. This $\aF$ symbol is associated to a loop in $M$ (which by definition goes through the triangle encoded by $X_0$). Let $\mathsf{L}$ denote its length (the number of triangles through which the loop passes) and $\mathsf{A}$ its area (number of triangles surrounded by it).   
 
%$L_0$ denote the loop containing the root triangle in the infinite map $(M, T)$; here as usual we identify a loop with the set of triangles through which it passes.
%Let $\mathsf{L}$ denote its length (the number of triangles it passes through) and let $\mathsf{A}$ denote its area, meaning the number of
%triangles surrounded by it (where the ``inside'' of the loop is the connected component of the loop which does not contain infinity).

Let
\begin{equation}\label{p0}
p_0 = \frac{\pi}{4 \arccos \left( \frac{\sqrt{2 - \sqrt{q}}}{2}
  \right) } =\frac {\kappa'}{8} \in (1/2, 1),
  \end{equation}
where $q$ and $\kappa'$ are related as in \eqref{q_kappa}. The following is the main result in \cite{BLR}.

\begin{theorem}
\label{T:typ}
Let $0< q< 4$. The random variables $\mathsf L$ and $\mathsf A$
satisfy
\begin{equation}\label{Tboundary}
  \P( \mathsf L> k) = k^{-1/p_0  + o(1)},
\end{equation}
and
\begin{equation}\label{Tarea}
 \P( \mathsf A > k) =  k^{-1 + o(1)}
\end{equation}
as $k\to \infty$.
\end{theorem}

As noted in \cite{BLR}, the laws of $\mathsf{L}$ and $\mathsf{A}$ correspond respectively to the limits of the length and area of a \emph{uniformly} chosen loop in the finite decorated planar map $(M_n, T_n)$ as $n \to \infty$. (By contrast, if we consider without any conditioning the length and area of the loop going through the triangle encoded by $X_0$, this would lead to different exponents, due to a size-biasing effect.)
%In other words, there is no ``biasing'' by the length of the loop when we condition on $X_0 = \aF$. 

Results in \cite{FKstory1, FKstory2, FKstory3} are analogous and more precise, in particular showing regular variation of the tail at infinity. (As a consequence, the sum of loop lengths and areas, in the order that they are discovered by the space-filling path, can be shown to converge after rescaling to a stable L\'evy process with appropriate exponent).

A particular consequence of Theorem \ref{T:typ} is that we expect the longest loop in the map $M_n$ to have size roughly $n^{p_0 + o(1)}$; that is,
\begin{equation} \label{E:max}
  \max_{\ell \in (M_n, T_n)}  |\ell | = n^{p_0 + o(1)}
\end{equation}
as $n\to \infty$.
%Indeed, first note that the loop $L_0$ containing the origin is biased by its length (since it contains 0) and so for a uniformly chosen loop $L$, we expect the length exponent to satisfy
%$$
%\P( |L| \ge k )= k^{- 1/p_0 + o(1)}.
%$$
Heuristically, to derive \eqref{E:max}, one then observes that $M_n$ contains order $n$ loops whose lengths are roughly i.i.d. with tail exponent $\alpha = 1/p_0$. The maximum value of this sequence of lengths is then easily shown to be of order $n^{1/\alpha + o(1)} = n^{p_0 + o(1)}$.

We will not prove Theorem \ref{T:typ}, but we will discuss in Exercise \ref{Ex:kpz} an interesting application using the KPZ formula.
These exponents are obtained (both in \cite{BLR} and \cite{FKstory1, FKstory2, FKstory3}) through a connection with a random walk in a cone. A simple setting, where it is easier to see this connection, is in the following result.

\begin{prop}[\cite{FKstory2}]
  \label{P:ext} Let $0< q< 4$, and let $E_{2n}$ be the event that the word $w= X_1 \ldots X_{2n}$ reduces to $\bar w = \emptyset$. Then
  $$
  \P(E_{2n}) = n^{-2p_0 -1 + o(1)} = n^{ -1 - \kappa'/ 4 + o(1)},
  $$
  as $n \to \infty$. In particular, $\P(E_{2n})$ decays subexponentially.
\end{prop}

\begin{proof}[Sketch of proof] We give a rough idea of where this exponent comes from, as it allows us to illustrate the connection to random walk in a cone, as mentioned above. A rigorous proof of this result may be found in \cite{FKstory2}.

%To avoid issues of size-biasing, instead of the loop $L_0$ we will condition on the event that the match of the symbol at time 0 is a fresh order $\aF$ (corresponding to the closure of a loop which under this conditioning started at time 0), and describe the statistics of the resulting loop $L$ (this is in some sense a ``typical loop''). In other words, if $\varphi (0) = n$ then we condition on $n>0$ and $X_{n} = \aF$. Let $E$ be this event.

The first step is to describe $E_{2n}$ in terms of the burger count processes $H$ and $C$ of Theorem \ref{T:scalinglimitFK}.
%Without loss of generality we may assume that at time 0 a hamburger is produced so $X_0 = \ah$ (matched to $\aF$ at time $n$). Now consider the word $w = X_0 \ldots X_n$ and observe that $\bar w  = \aC \ldots \aC$ can only consist of cheeseburger orders, of some length $\ell$. \ellen{This is not correct **} This length $\ell$ corresponds exactly to the length of the loop \ellen{** The length of the loop is the number of maximal sub-excursions + number of triangles not in a sub-excursion minus one ** } whereas the total duration $n$ of the loop is the number of triangles enclosed by the loop, that is, its area. \ellen{I think this second statement is correct.}
In particular, we note that the event $E_{2n}$ is equivalent to the conditions
\begin{itemize}
  \item $C_i, H_i \ge 0$ for $0 \le i \le 2n$; and

  \item $C_{2n} = 0$, $H_{2n} =0$
\end{itemize}
on $H$ and $C$.
Indeed, the first condition holds since if at some point $1\le k\le 2n$ the burger count $C$ or $H$ becomes negative, this must be because of an order whose match in the bi-infinite sequence $(X_k)_{k \in \Z}$ was in the past, that is, $\varphi(k) <0$. Therefore, the event $E_{2n}$ is equivalent to the process $Z_k = (C_k,H_k)_{1\le k \le 2n}$ being an excursion in the top right quadrant of the $(C,H)$ plane, starting and ending at the origin.

This probability may be computed approximately (or rather, heuristically here) using Theorem \ref{T:scalinglimitFK}. To do this it is useful to apply first a linear map of the $(C,H)$ plane so as to deal with independent Brownian coordinates in the limit. More precisely, we apply the linear map $\Lambda$ defined by
$$
\Lambda =(1/\sigma)\left(
\begin{array}{cc}
1 & \cos(\theta_0)\\
0  & \sin (\theta_0)
\end{array}
\right),
$$
where $\theta_0= \pi /(2p_0) = 4\pi/\kappa' = 2 \arctan( \sqrt{1/(1-
  2p)})$ and $\sigma^2 = (1-p)/2$. A direct but
tedious computation shows that $\Lambda (L_t, R_t)$ is indeed a
standard planar Brownian motion. (The computation is easier to do by
reverting to the original formulation of \cref{T:scalinglimitFK} in
\cite{Sheffield_burger}, where it is shown that $C + H$ and $(C-H)/\sqrt{1-2p}$ converge to a standard planar Brownian motion).
\begin{figure}
\begin{center}
\includegraphics[scale=.9]{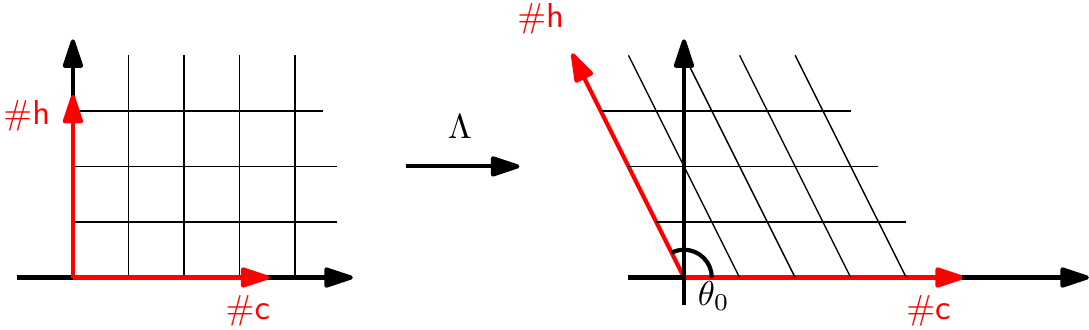}
\caption{The coordinate transformation. In these new axis, the burger counts $H$ and $C$ become independent Brownian motions; the event $E_{2n}$ then corresponds to $\Lambda(Z)$ making an excursion of duration $2n$ in the cone $\cC(\theta_0)$ of angle $\theta_0 = \pi / (2 p_0) = 4\pi / \kappa'$, starting and ending at its apex. }
\label{F:coord}
\end{center}
\end{figure}
Note that the top right quadrant transforms under $\Lambda$, see \cref{F:coord}, into the cone $\cC(\theta_0)$ of angle $\theta_0$ with apex at zero.

%Because of Theorem \ref{T:scalinglimitFK}, the question (heuristically) reduces to the following simpler question. Let $B$ be a
%standard
We therefore consider an analogous question for two dimensional Brownian motion. Namely, let $B$ be a standard planar Brownian motion, starting from some point $z \in \cC(\theta_0)$ with $|z|=1$. Let $T$ be the first time that $B$
leaves $\cC(\theta_0)$. Then from Theorem \ref{T:scalinglimitFK} it is reasonable to guess that
$$
\P(E_{2n}) \approx \P_z( T>t; |B_t| \le 1), \text{ with }
t = n^{1 + o(1)}.$$
(Indeed, note that if $T> t$ and $|B_t| \le 1$ then the Brownian motion is likely to exit the cone soon after time $t$ and not far from the apex. This intuition is for instance made rigorous in \cite{BLR} and \cite{FKstory1, FKstory2, FKstory3}.)

For this we first claim that
\begin{equation}\label{coneBM}
\P(T > t )= t^{- p_0+o(1)}
\end{equation}
as $t \to \infty$. To see why this is the case, consider the conformal isomorphism $z \mapsto z^{\pi/\theta_0}$, which sends the cone $\cC(\theta_0)$ to the upper half plane. In the upper half plane, the function $z \mapsto  \Im (z)$ ($\Im(z)$ being the imaginary part of $z$) is harmonic with zero boundary condition, and so in the cone, the function
$$
z\mapsto g(z) : = r^{\pi/\theta_0} \sin \left(\frac{\pi
    \theta}{\theta_0}\right) ; \quad  z\in \cC(\theta_0),
$$
is also harmonic. Applying the optional stopping theorem at time $ t$ to the martingale $M_t:=g(B_{t\wedge T})$, the only non-zero contribution to $M_{ t }$ comes from the event $T>t$. On the other, conditionally on $T>t$, $B_t$ is likely to be at distance $\sqrt{t}$ from the origin, in which case $M_{t} \approx t^{\pi/(2 \theta_0)} = t^{p_0}$. It is not hard to deduce \eqref{coneBM}.

We now claim that the desired probability satisfies
\begin{equation}
   \P_z( T>t; |B_t| \le 1) = t^{ -2p_0 -1 + o(1)} \text{ as } t\to \infty.
\end{equation}
To see this, we split the interval $[0,t]$ into three intervals of equal length $t/3$. In order for the event on the left hand side to be satisfied, three things must happen during these three intervals.
\begin{itemize}
\item Over the interval $[0, t/3]$, $B$ must not leave the cone. This has probability $t^{- p_0 + o(1)}$ by \eqref{coneBM}.

\item At the other extreme, if we reverse the direction of time, we also have a Brownian motion started close to the tip of the cone that must not leave the cone for time $t/3$. Again, this has probability $t^{- p_0 + o(1)}$.

\item Finally, given the behaviour of the process over $[0, t/3]$ and $[2t/3, t]$, the process must go from $B_{t/3}$ to $B_{2t/3}$ in the time interval $[t/3, 2t/3]$, and stay inside the cone. The latter requirement actually has probability bounded away from zero (because $B_{t/3}$ and $B_{2t/3}$ are typically far away from the boundary of the cone), so it remains to compute the probability to transition between these two endpoints. However this is roughly of order $t^{-1+o(1)}$, since we are dealing with a Brownian motion in dimension two.

\end{itemize}

Altogether, we obtain that $\P_z( T>t; |B_t| \le 1) = t^{ -2p_0 -1 + o(1)},$ as desired.

\end{proof}

\subsection{Exercises}

\begin{enumerate}[label=\thesection.\arabic*]

\item \label{Ex:symbolmaps} This exercise follows the arguments of Chen  \cite{chen} and gives a very nice concrete construction of the planar map associated to a word. 

Let $x_{1}, \ldots, x_{2n}$ be a sequence of $2n $ letters in the alphabet $\Theta = \{\ac, \ah, \aC, \aH, \aF\}$ and suppose that the corresponding word $ w = x_{1} \ldots x_{2n}$ reduces to $\bar w = \emptyset$. For each $1 \le i \le 2n$, denote by $\varphi (i)$ the unique match of $i$: meaning that if $i$ corresponds to production of a specific burger, then $\varphi(i)$ is the unique time at which this burger is consumed, and vice versa.

    Let us draw a map as follows. Start with the line segment (drawn in the complex plane) having vertices $1, \ldots, 2n$ and horizontal nearest neighbour edges. Draw an arc between $i$ and $\varphi(i)$ for each $1\le i \le 2n$; this arc is drawn in the upper half plane for a hamburger, and in the lower half plane for a cheeseburger.

    (a) Show that the arcs can be drawn in a planar way (so they don't cross one another); in other words, if $n_1 < n_2 < n_3 < n_4$ it is not possible that $\ph(n_1) = n_3$ and $\ph(n_2) = n_4$, unless $x_{n_1} \neq x_{n_2}$.

    (b) Add an additional edge between $1$ and $2n$ in the upper half plane, above every other edge, and call the resulting map  $\mathbf{A}$. Check that $\mathbf{A}$ has $2n$ vertices and $3n$ edges. 
    
    (c) Consider the planar dual $\mathbf{\Delta}$ of $\mathbf{A}$, and show this is a triangulation. 
    
    (d) Color in blue the edges of $\mathbf{\Delta}$ that stay in the upper half plane, and in red those that stay in the lower half plane. Show that the set of blue edges and the set of red edges define two trees. Let $\mathbf{Q}$ be the set of remaining edges in $\mathbf{\Delta}$, and colour them green. Show that $\mathbf{Q}$ is a quadrangulation. 
      
    (e) Explain how the map $\mathbf{\Delta}$ is related to the triangulation $\mathbf{m}_n$ encoded by Sheffield's bijection, and show that the straight line segment from $1$ to $2n$ together with the additional edge linking the two extreme vertices in $\mathbf{A}$ corresponds to the space-filling loop in Sheffield's bijection.
    
  (f) {Deduce that local convergence of maps is equivalent to local convergence of the symbols encoding them via Sheffield's bijection, as claimed in \cref{thm:Local}.}

\item \textbf{The reduced walk.} {Consider the infinite decorated planar map $(M,T)$ of \cref{thm:Local}, and let $(X_n)_{n \in \Z}$ denote the bi-infinite sequence of symbols encoding it via Sheffield's bijection. Let us assume that $q>0$ or equivalently $p>0$, where $p$ and $q$ are related via \eqref{pq} and $p$ is the proportion of $\aF$ symbols. Define a \emph{backward} exploration process $(c_n,h_n)_{n \ge 0}$ of the map, which keeps track of the number of
$\aC$ and $\aH$ in the reduced word, as follows. Let $(c_0,h_0 ) =
  (0,0)$. Suppose we have performed $n$ steps of the exploration and defined
$c_n,h_n$ and in this process, we have revealed the 
symbols $(X_{-m}, \ldots, X_0)$.
We inductively define the following.
\begin{itemize}
\item If $X_{-m-1}$ is a $\aC$  (resp. $\aH$), define
$(c_{n+1},h_{n+1}) = (c_n,h_n)   + (1,0)$
(resp. $(c_n,h_n)   + (0,1)$).
\item If $X_{-m-1}$ a $\ac$ (resp. $\ah$), $(c_{n+1},h_{n+1}) = (c_n,h_n)
  + (-1,0)$ (resp. $(c_n,h_n) + (0,-1)$).
\item If $X_{-m-1}$ is $\aF$, then we explore $X_{-m-2},X_{-m-3}\ldots$ until we find the  match of $X_{-m-1}$.
  % Notice that the reduced word
  %$\cR_{n+1}=\overline{X_{\ph( -m-1)} \ldots X_{-m-1}}$ is either of the form
%$\aC \aC \ldots \aC $ or
 % $\aH \aH \ldots \aH$ depending on whether the match of the $\aF$
 % is a $\ah$ or $\ac$ respectively. Either happens with equal
 % probability by symmetry.
  Let $|\cR_{n+1}|$ denote
  the number of symbols in the reduced word
  $\cR_{n+1} =\overline{X_{\varphi( -m-1)} \ldots X_{-m-1}}$. Show that $\cR_{n+1}$ contains only order symbols of one type. If $\cR_{n+1}$ consists
of $\aH$ symbols,
  define $(c_{n+1},h_{n+1}) = (c_n,h_n) + (0,|\cR_{n+1}|)$. Otherwise, if $\cR_{n+1}$ consists of $\aC$ symbols
define $(c_{n+1},h_{n+1}) = (c_n,h_n) + (|\cR_{n+1}|,0)$.
\end{itemize}
Show that the walk
$(c_n,h_n)_{n \ge 0}$ is a sum of \emph{independent} and identically distributed random variables. Note that this is in contrast to Theorem \ref{T:scalinglimitFK}. It can be shown that these random variables are in fact centred when $q\le4$ (see \cite{Sheffield_burger}). }

\item \textbf{Bubbles.} {Consider the infinite decorated planar map $(M,T)$ of \cref{thm:Local}, and let $(X_n)_{n \in \Z}$ denote the bi-infinite sequence of symbols encoding it via Sheffield's bijection. Let us assume that $q>0$ or equivalently $p>0$, where $p$ and $q$ are related via \eqref{pq} and $p$ is the proportion of $\aF$ symbols. Let us condition on the event $X_0 = \aF$. Let $\varphi(0)\le 0$ denote the match of this symbol. The word $w = X_{\varphi(0)}\ldots X_0$ encodes a finite planar map, called the \textbf{bubble} or envelope of the map at 0. This bubble corresponds to a finite number of loops of $(M,T)$ (note that this is in general more than a single loop of $(M, T)$ containing the root triangle, as there can be other $\aF$ symbols in $w$). This notion was pivotal in \cite{BLR} where it was used to derive critical exponents of Theorem \ref{T:typ}. This exercise gives one of the main steps in the derivation of this theorem.}

    {(a) Assume without loss of generality that $X_{\varphi(0)} = \ah$. Give a description of the reduced word $\bar w$. By considering the random length $N = |\varphi(0)|$ of $w$ and the random length $K$ of the reduced word $\bar w$, describe the event $\{N = n, K=k\}$ in terms of a certain cone excursion for the reverse two dimensional walk $(C_{-i}, H_{-i})_{0 \le i \le n}$. Explain why $N$ is the area of the bubble and $K$ the length of its outer boundary.}

    {(b) Arguing at the same level of rigour as in Proposition \ref{P:ext}, show that there are exponents $p_{\text{area}}$ and $p_{\text{boundary}}$ such that
    $$
    \P( N \ge n) = n^{- p_{\text{area}} + o(1)}; \quad  \P( K \ge k) = k^{- p_{\text{boundary}} + o(1)}
    $$
    where  $p_{\text{boundary}}/2 = p_{\text{area}} = p_0$, and $$
    p_0 = \frac{\pi}{4 \arccos \left( \frac{\sqrt{2 - \sqrt{q}}}{2}
  \right)}  =\frac {\kappa'}{8} \in (1/2, 1),
$$
was defined in \eqref{p0}.}

{\emph{The next three exercises use exponents derived in this chapter together with the KPZ formulas of the previous chapter to give predictions (in some cases proved through other methods) about the value of certain critical exponents associated with random fractals which can be defined without any reference to random planar maps.}}

\item \label{Ex:kpz}
Use \eqref{E:max}, the KPZ relation, and the relation
$$q = 2 + 2 \cos(8\pi / \kappa')$$
between $q$ and $\kappa'$ from \eqref{q_kappa},
to recover {(non-rigorously)} that the dimension of SLE$_{\kappa'}$ is $1+ \kappa'/8$ for $\kappa'\in (4,8)$.
\ind{Critical exponents!FK percolation}

\item Consider a simple random walk on the (infinite) uniform random planar map $G$, that is, take $G$ to be the infinite volume of the FK-weighted maps for $q=1$ defined in Theorem \ref{thm:Local}, and let $(X_n)_{n\ge 0}$ be a simple random walk on $G$ starting from the root. If $n \ge 1$, a pioneer point for the walk $(X_1 , \ldots X_n)$ is a point $x$ such that $x$ is visited at some time $m \le n$ and is on the boundary of the unbounded  component of $G\setminus\{X_1, \ldots, X_m\}$.
A beautiful theorem of Benjamini and Curien \cite{BenjaminiCurien} shows that when such a simple random walk 
first exits a ball of radius $R$, it has had $\approx R^3$ pioneer points (technically this result is only proved when $G$ is the so called Uniform Infinite Planar Quadrangulation, although it is also believed to hold for infinite planar maps within the same universality class such as the one considered above).

Analogously, for $(B_s)_{s\ge 0}$ a planar Brownian motion, we define the set {$\cP_t$ for given $t>0$ to be all points of the form} $B_s$ {for some $0 \le s \le t$}, such that $B_s$ is on the ``frontier'' at time $s$ ({where by frontier we mean} the boundary of the unbounded component of the complement of $B[0,s]$).

{Using a {(non-rigorous)} KPZ-type argument,} derive  {the dimension of the} Brownian pioneer points {$\cP_t$ for any fixed $t \ge 0$}.  (The answer is 7/4, as rigorously proved in a famous paper of Lawler, Schramm and Werner  \cite{LSW43} using SLE techniques).
\ind{Pioneer points}

\item Consider a simple random walk $(X_n)$ on the infinite local limit of FK-weighted planar maps {(as in \cref{thm:Local})}, starting from the root. Try to argue using the KPZ relation {(again without being fully rigorous)}, that {the graph distance between $X_n$ and $X_0$} must be approximately equal to $ n^{1/D}$ where $D$ is the dimension of the space. (Hint: the range of Brownian motion must satisfy $\Delta = 0$; more precisely, by the time a random walk leaves a ball of radius $R$, it has visited of order $R^2/\log R$ vertices with high probability). In particular, on the UIPT, one conjectures that {this distance is} $\approx n^{1/4}$.
    This has now been proven rigorously in \cite{GwynneHutchcroft} and \cite{GwynneMiller_rw}.
\ind{Subdiffusivity}

%\item Another proof of uniformity. \ellen{*?*}

\end{enumerate}

\newpage

\section{Introduction to Liouville conformal field theory}\label{S:LCFT}
% !TEX root = master.tex

%Introduction to Liouville CFT

\def\hsg{\tilde{h}^{\Sp,g}}
\def\MM{\cM^g}

In this chapter we present a short introduction to the theory initiated in the pioneering paper of David, Kupiainen, Rhodes and Vargas \cite{DKRV}, which we will refer to as Liouville conformal field theory, or Liouville CFT for short. We use this in order to avoid confusion with the SLE based theory developed in Chapters \ref{S:SIsurfaces} onwards, for which we choose to stick with the label of Liouville quantum gravity.

The main objectives of the two theories are similar (that is to say, making rigorous sense of Polyakov's conformal theory of quantum gravity), and indeed we will see concrete statements connecting these two approaches in Chapter \ref{S:SIsurfaces}. %see \ref{}.
Nevertheless they are entirely independent, and can be read in whichever order one chooses. In particular the Liouville CFT we are about to present does not require knowledge of SLE (it depends only on Gaussian multiplicative chaos theory). It also presents the advantage of being closer to the original formulation of Polyakov.

We will start with some heuristics and then move on to a rigorous definition motivated by these heuristics, staying for simplicity in the context of the Riemann sphere. (See \cite{GRVpolyakov} for an extension to more general Riemann surfaces; this extension is highly nontrivial due to the need to choose the conformal class at random with a suitable law, in contrast to the case of the sphere where this is not necessary.) We will then prove the existence of the theory (which is to say, the finiteness of some observables subject to the so called \textbf{Seiberg bounds}). An absolutely remarkable feature of the theory is that it is in some sense integrable or exactly solvable. We will show a simple result which hints at this integrability: the $k$ point function of the theory can be computed as a negative fractional moment of Gaussian multiplicative chaos. We conclude with a brief overview of some recent developments, including a short discussion of \textbf{conformal bootstrap} (\cite{GKRV}) and the proof by Kupiainen, Rhodes and Vargas \cite{DOZZ} of the celebrated \textbf{DOZZ formula}. \ind{Seiberg bounds}
\ind{Conformal bootstrap}
\ind{DOZZ formula}

\subsection{Preliminary background}

\subsubsection{Quantum and conformal field theory}

It is helpful to begin with a brief and \emph{very} informal overview of some underlying notions which help put Polyakov's proposal in context. A \textbf{statistical field theory} (also known as a \textbf{Euclidean field theory}) is, very roughly speaking, a random field $(\ph(x))_{x \in \R^d}$, or collection of such fields, defined in the continuum space $\R^d$ (or some region $D \subset \R^d$). \ind{Quantum field theory} A probabilist might intuitively think about the scaling limits of discrete fields naturally arising in statistical mechanics; for example, the magnetisation field in the Ising model, at or away from the critical points (this field counts the sum of all Ising spins in a given region). As this example suggests, one should not expect the ``statistical fields'' to be defined pointwise; instead, like the GFF they should be understood as random distributions. Physicists typically describe such fields via their $k$ point \textbf{correlation functions}:
$$
(x_1, \ldots, x_k) \mapsto \E[ \ph(x_1) \ldots \ph(x_k)].
$$
Although the field $\ph$ is typically not pointwise defined, such correlation functions \emph{are} typically well defined. For instance, in the case of the Gaussian free field, they can be computed from the knowledge of the two point function,  a multiple of the Green function, and \textbf{Wick's rule} which expresses the $k$ point functions of Gaussian fields in terms of their two point functions. (Note however that the two point function does not in general determine the $k$ point function.) \ind{Correlation functions} \ind{Wick's rule} 

Another subtlety is that, in many cases of interest, the underlying measure $\P$ with respect to which the above correlations are computed is in fact not a probability distribution but rather a positive measure which may well have infinite mass. For this reason, the correlations will usually be written as $\langle  \ph(x_1) \ldots \ph(x_k) \rangle$ rather than as expectations.
Furthermore, the quantities that are actually of interest to physicists are analytic continuations of these correlation functions in terms of the underlying parameters defining the model (for example, inverse temperature).
Indeed the resulting quantities can be interpreted in terms of \textbf{quantum field theory}. Roughly speaking, the statistical field theory described above corresponds to a quantum field theory via what is known as a 
 \textbf{Wick rotation}: essentially, multiplying one of the spatial coordinates (the `time' coordinate) of a quantum field theory by $i$ allows us to go from the quantum theory to a real valued, and indeed positive, measure, which describes the statistical field theory. See \cite{Mussardo} for an account of exactly solvable models of statistical field theory.

\medskip A (Euclidean) \textbf{Conformal Field Theory} is a particular case of statistical field theory, in which the theory is required to satisfy certain additional invariance properties under conformal mapping, often referred to as \textbf{conformal symmetries}. Note that this makes sense even in dimensions greater or equal to three, in which case conformal maps are simply diffeomorphisms that preserve angles locally. The central objects in conformal field theory are a family of \textbf{primary fields} denoted by $\{\psi_\alpha\}_{\alpha \in A}$. For instance, in the case of the Ising model, the primary fields are given by $\{1, \sigma, \cE\}$ where $\sigma$ is the spin field (the scaling limit of the sum of Ising spins in a given region) and $\cE$ is the energy field (the energy of an edge $e = (x,y)$ is the contribution $\sigma_x \sigma_y$ to the total energy of the configuration, and $\cE$ gives the scaling limit of the sum of these energies in a given region). In Liouville conformal field theory, the primary fields $\psi_\alpha$ will, roughly speaking, be given by $\psi_\alpha (z) = e^{\alpha h (z)}$ (suitably interpreted), and $h$ will be ``sampled'' from an infinite measure which is related to the law of a Gaussian free field. 

In a conformal field theory, these primary fields can be multiplied with one another in some formal sense, and this allows us to talk about correlation functions $\langle \psi_{\alpha_1 } (z_1) \ldots \psi_{\alpha_k} (z_k)\rangle$, as above. The first assumption of conformal symmetry is that, whenever $f$ is a M\"obius map (that is, a conformal isomorphism from the underlying domain $D$ in which theory is defined to itself) 
\begin{equation}\label{E:globalsymmetry}
\langle \psi_{\alpha_1} ( f (z_1) ) \ldots \psi_{\alpha_k} (f (z_k)) \rangle  =\left( \prod_{i=1}^k | f'(z_i)|^{-2\Delta_{\alpha_i}}\right) \langle \psi_{\alpha_1} (z_1) \ldots \psi_{\alpha_k} (z_k)\rangle,
 \end{equation}
 for some numbers $\Delta_\alpha, \alpha \in A$\indN{{\bf Liouville CFT}! $\Delta_\alpha$; conformal weights} called the \textbf{conformal weights} \ind{Conformal weights} associated to the primary fields. The assumption \eqref{E:globalsymmetry} describes a global symmetry condition as it imposes a constraint on how the correlation functions change under the application of a globally defined conformal isomorphism on $D$. Two dimensional conformal field theories also satisfy a more local kind of symmetry condition. There are several viewpoints that may be used to formulate these more local symmetries.  One way is via the so called \textbf{Virasoro algebra}. This is beyond the scope of the present chapter, but roughly speaking, the Virasoro algebra is generated by a family of operators $(L_n)_{n \in \Z}$ 
together with a central element that 
%which can be identified with a number $c \in \R$, 
commutes with every $L_n$ and so is in the centre of the algebra. More generally the $L_n$ satisfy certain commutation relations involving the central element.
%\ellen{I think the algebra has a central element, and the central charge is something is defined when we have a given representation of the Virasoro algebra};
%\ellen{this seems a bit meaningless before talking about the representation?} \nb{I agree, but cannot come up with something better}). 
Infinitesimal conformal symmetries are enforced by requiring that there is a representation of this algebra (often but not always unitary)
on a vector space containing the primary fields (together with the so called descendant fields). In this representation the central element can be identified with a number $c \in \R$ called \textbf{the central charge}. Moreover, it makes sense to ``apply'' $L_n$ to a primary field $\psi_\alpha$, and it is worth noting that the operators $L_n$ associated to the levels $n = -1, 0, 1$ correspond in some informal sense to M\"obius maps, so that this is indeed a generalisation of \eqref{E:globalsymmetry}. While such a rigorous description has recently been announced for Liouville conformal field theory (see \cite{BGKRV}) we will not pursue this here.

 %\ellen{not sure this is true?} \nb{that's what I understood from the discussion with Eveliina}.  
 
\medskip  Another possibility (note that it is not \emph{a priori} obvious whether the two descriptions are equivalent, and we do not claim this) is via the so called \textbf{Weyl invariance} (or more precisely in our case \textbf{Weyl anomaly}) property. To state this it is necessary to enrich the problem by considering the theory, with respect to which the correlation functions $\langle \psi_{\alpha_1 } (z_1) \ldots \psi_{\alpha_k} (z_k)\rangle$ are computed, as being defined on a manifold $M$ instead of a domain $D$ and suppose that $M$ is endowed with a background metric $g$. 
When we do so, for every metric $g$ on $M$ we should have an associated collection of correlation functions, which we denote $\langle \psi_{\alpha_1 } (z_1) \ldots \psi_{\alpha_k} (z_k)\rangle_g$. To get a feeling for what this might correspond to in the case of the Ising model, say, consider the following toy example: let $D$ be a domain and $U \subset D$ be a fixed subdomain, and take $g $ to be twice the Euclidean metric in $U$, and the Euclidean metric in the complement $ D \setminus U$. The corresponding correlation function should describe the scaling limits of Ising correlations for graphs in which the density of vertices in $U$ is twice as large as that in $D\setminus U$.

With these notations, let us now
 describe 
 %%
 %%
\begin{comment}
 the simpler notion of \textbf{diffeomorphism invariance} before we describe Weyl invariance. 
Let $F: M \to M$ be a diffeomorphism, and let $F_* g$ denote the pushforward of the metric $g$ under $F$ (that is, $F_* g (x)$ is the metric whose value at the point $F(x)$ is, in matrix notations, $ (DF^{-1})^T g(x) DF^{-1}$, $DF$ is the Jacobian matrix of the diffeomorphism of $F$, and $DF^{-1}$ is its inverse). Diffeomorphism invariance requires 
\begin{equation}\label{E:diffeo_inv}
\langle \psi_{\alpha_1 } (z_1) \ldots \psi_{\alpha_k} (z_k)\rangle_{g}  =  \langle \psi_{\alpha_1 } (F(z_1)) \ldots \psi_{\alpha_k} (F(z_k))\rangle_{F_*g}. 
\end{equation}
 Physically this amounts to requiring the theory to be independent of the choice of coordinates. Now let us describe 
\end{comment}
%%
%%
the Weyl invariance property. If $g$ is a metric, and $\rho: M \to \R$ is a smooth function, we obtain a conformally equivalent metric $\tilde g$ by setting $\tilde g = e^\rho g$: that is, the angle of the curves on $M$ are locally the same under $g$ and $\tilde g$, and the distances in $\tilde g$ are locally multiplied by $e^\rho$. Such a rescaling of the background metric is sometimes known as a \textbf{Weyl transformation}. 
Then, Weyl invariance would be the identity
 \begin{equation}\label{E:Weyl_inv}
 \langle \psi_{\alpha_1 } (z_1) \ldots \psi_{\alpha_k} (z_k)\rangle_{e^{\rho}g}  =  \langle \psi_{\alpha_1 } (z_1) \ldots \psi_{\alpha_k} (z_k)\rangle_g. 
 \end{equation}
 %Note the similarity between \eqref{E:diffeo_inv} and \eqref{E:Weyl_inv}. Both require the correlation functions to be invariant under some natural changes of the metric, but neither implies the other. 
 However, while Weyl invariance is a natural requirement for conformal theories describing classical physics, in \emph{quantum} conformal field theories this is not the case; instead, one has the \textbf{Weyl anomaly} 
 \begin{equation}\label{E:Weyl_anomaly}
 \langle \psi_{\alpha_1 } (z_1) \ldots \psi_{\alpha_k} (z_k)\rangle_{e^{\rho}g}  =  e^{\frac{c}{96 \pi} 
% \int_M  (| \nabla^g \rho |^2 + 4 \rho ) v_g 
A( \rho , g)  
 } \langle \psi_{\alpha_1 } (z_1) \ldots \psi_{\alpha_k} (z_k)\rangle_g, 
 \end{equation}
 where the anomaly term $A(\rho, g)$ is defined by 
 $$
A( \rho, g) =   \int_M  (| \nabla^g \rho |^2 + 2R_g \rho ) v_g .
 $$
(Sometimes, the Weyl anomaly formula is expressed slightly differently, see for instance Remark \ref{R:Weyl_correlation}.) Here $c$ is the \textbf{central charge} of the theory, $v_g$ is the volume form on $M$ associated to the metric $g$, $R_g$ is the scalar curvature, and $\nabla^g \rho$ denotes the gradient of $\rho$ computed in the metric $g$. The Weyl anomaly formula \eqref{E:Weyl_anomaly} replaces \eqref{E:Weyl_inv} and allows us to consider arbitrary rescalings of the metric; property \eqref{E:Weyl_anomaly} captures the desired ``local'' conformal transformations mentioned earlier. In the context of Liouville conformal field theory we will be able to prove the Weyl anomaly formula (see Theorem \ref{T:Weyl}).
 %for smooth observables; combining it with Theorem \ref{T:LfieldMobius} would show the form \eqref{E:Weyl_anomaly} at least in the case where $e^\rho g $ is the pushforward $m_*g$ of the metric $g$ by a M\"obius map $m$; note that in that case the conformal anomaly term $A( \rho, g) = 0$
In particular, Theorem \ref{T:Weyl} identifies the central charge of the theory. 
(We also note that from this point of view it is natural to consider infinitesimal deformations of the metric, that is, when $\rho = \rho_\delta = \delta \hat \rho$ for some fixed smooth function $\hat \rho:M \to \R$, so that $e^\rho g = (1 + \delta \hat \rho + o(\delta) ) g $; the corresponding change in the correlation functions would involve to the first order a quantity known as the \textbf{stress energy tensor} of the theory).

%namely, they obey a change of variable formula in which the correlation function is multiplied by a given power of the (modulus of the) derivative of the conformal map. 
%\ellen{Weyl invariance was what Joona was talking about in "classical field theory''. Probably want to say that we will not have Weyl invariance in our case but rather a Weyl anomaly; I think that's generally the case in CFT?}

\medskip Conformal field theory grew in the 1980s after it was observed (or rather predicted) by Polyakov that such conformal symmetries arise when the underlying statistical mechanics models are taken at their critical point \cite{Poly70,BPZ,BPZ2, FQS}.
 Furthermore, it turns out that at least in two dimensions, adding this requirement of conformal symmetries to a natural list of axioms for quantum field theory (as introduced by Osterwalder and Schrader \cite{OSaxioms1,OSaxioms2}), drastically impacts the space of solutions to these axioms. This leads to a classification of conformal field theories at least in the case of unitary representations. \ind{Central charge}

\subsubsection{Polyakov action}

%\ellen{Link more clearly to previous section...}

Having discussed the general context of quantum and conformal field theories, we now turn our attention to the specific case of Liouville conformal field theory, which will occupy us in this chapter. In order to assist the reader we begin with rather general considerations on statistical mechanics. 
In physics, a \textbf{Hamiltonian} $H$ is a function which assigns the energy $H(\sigma)$ to a configuration $\sigma$. In statistical physics, we are used to the idea of sampling a configuration $\sigma$ according to the \textbf{Gibbs measure} (with respect to an underlying reference measure denoted by $\text{d} \sigma$), namely
\begin{equation}\label{E:gibbs}
\P(\sigma) \propto \exp ( - \beta H(\sigma)) \text{d}\sigma.
\end{equation}
 Here $\beta\ge 0$ is a parameter playing the role of the inverse temperature of the system.

\ind{Gibbs measure}
\ind{Hamiltonian}
\ind{Action}
\ind{Polyakov action}

An \textbf{action} is an energy integrated against time: it represents the amount of energy needed to bring the system from one configuration to another. For a two dimensional field $\ph:\R^2 \to \R$ (where as above we view the field as a random distribution, so that $\ph$ is not really pointwise defined), the \textbf{Polyakov  action} $S(\ph)$ associated to the field $\ph$ can be thought of directly as the energy of the configuration $\ph$, so that for a probabilist used to statistical mechanics, there is no difference between the Hamiltonian of the system (the energy of the configuration $\ph$) and the action $S(\ph)$. The reason for this apparently confusing terminology is that in this 2d model of quantum gravity, one should remember that one of the two dimensions is space and the other is time. Thus by specifying the energy $S(\ph)$ we have already integrated against time and are thus properly dealing with an action. We will keep this convention and refer to $S(\ph)$ as the Polyakov action, but it should simply be thought of as the energy of the configuration $\ph$. We are now ready to give an expression (which for the moment is purely formal) for this Polyakov action on the sphere. \medskip 

To describe the action, we first need to fix a Riemannian metric $g$ on the two-sphere $\mathbb{S}=\{x\in \R^3: |x|\le 1\}$\indN{{\bf Geometry}! $\Sp$; unit two-sphere}. For computational purposes, it will often be easier to consider the pushforward of $g$ under a conformal isomorphism 
\begin{equation}
\label{E:psi} 
\psi:\Sp\to \hat{\C}
\end{equation} 
from $\Sp$ to the extended complex plane $\hat \C:=\C\cup \{\infty\}$.
\indN{{\bf Geometry}! $\hat{\mathbb{C}}$; extended complex plane $\mathbb{C}\cup\{\infty\}$}
%\indN{{\bf Geometry}! $\hat \C$; extended complex plane $\C\cup \{\infty\}$}
From now on, we will assume that the map $\psi$ in \eqref{E:psi} has been fixed; for example, we could take it to be stereographic projection. 
We will write $\hat{g}(z)$ for the pushforward of the metric $g$ on $\Sp$ to $\hat\C$, which we identify with a non-negative function on $\C$. So, a small region of area $\eps$ around the fixed point $z \in \C$ represents a region on $\mathbb{S}^2$ of area approximately $\hat{g}(z) \eps$ as $\eps \to 0$, while the distance between two points on the sphere is obtained by minimising the integral of $\sqrt{\hat{g}(z)}$ along paths between the two corresponding points on $\hat\C$.

We will be particularly interested in the spherical metric $g_0$ on $\Sp$\indN{{\bf Geometry}! $g_0$; spherical metric on $\mathbb{S}$}, which corresponds on $\hat \C$ to the function
\begin{equation}\label{spheremetric}
\hat g_0(z) = \frac{4}{(1 + |z|^2)^2}.
\end{equation}
\indN{{\bf Geometry}! $\hat g_0$; spherical metric on $\hat{\mathbb{C}}$}
For instance one can check that $\int_{\C} \hat g_0(z) \dd z = 4 \pi$, as required for the area of the unit sphere. In fact, without loss of generality, in what follows we will consider only metrics $g$ on $\Sp$, conformally equivalent to $g_0$: this means that on $\hat\C$, $\hat g$ must take the form
\begin{equation}\label{E:metric}
\hat g (z) = e^{\rho(z)} \hat g_0(z); \ \ z \in \C,
\end{equation}
with $\rho$ a twice differentiable function on $\C$ with finite limit at infinity such that $\int_{\C} |\nabla \rho|^2 <\infty$. We call $v_g$ the associated volume form on $\Sp$ and $v_{\hat g}$ the associated volume form on $\hat \C$.\indN{{\bf Geometry}! $v_g$; volume form associated with a metric $g$}

From now on, we also assume that the parameter $\gamma \in (0,2)$ is fixed and let
\begin{equation}\label{Qpolyakov}
Q = \frac{\gamma}{2} + \frac{2}{\gamma},
\end{equation}
which is the value first encountered in the change of coordinate formula for Liouville measure and the definition of random surfaces (see Theorem \ref{T:ci} and Definition \ref{D:surface} respectively). Finally, we let $\mu > 0$ denote a constant (the \textbf{cosmological constant}) whose value -- apart from the important fact that it is positive -- will not be of any relevance in the following.
	
	 With these notations, Polyakov's ansatz is to define the action as follows:
%\ellen{*check edit*}
\begin{equation}
  \label{E:Polyakov_informal}
  S(\ph) = \frac1{4\pi} \int_{\C} \Big[|\nabla^g \ph(z)|^2 + R_g Q \ph(z)  + 4 \pi \mu e^{\gamma \ph(z)}  \Big] v_g( \dd z),
\end{equation}\indN{{\bf Liouville CFT}! $S(\ph)$; Polyakov action}\ind{Polyakov action}\ind{Action!Polyakov action}
 where $R_g$\indN{{\bf Geometry}! $R_g$; scalar curvature associated to $g$} is the scalar curvature associated to $g$. On $\hat \C$,  the scalar curvature can be written explicitly as
\begin{equation}
\label{E:scalar}
R_{\hat g} (z) = - \frac1{\hat g(z)} \Delta \log \hat g(z); \, z\in \C.
\end{equation}
\ind{Scalar curvature}
The theory we are about to discuss is slightly simpler when the scalar curvature $R_g(z)$ is a constant.  In particular this includes the spherical metric $\hat g_0$, for which  $R_{\hat g_0} \equiv 2$ (this can be seen by expressing the Laplacian in polar coordinates; we leave this as an exercise). We also note that in general, due to the Gauss--Bonnet theorem (see for example, \cite{RG}),
  $$
  \int_\C R_{\hat g} (z) v_{\hat g} (\dd z) = 8 \pi.
  $$ \ind{Gauss--Bonnet theorem}
Returning to \eqref{E:Polyakov_informal}, we call the reader's attention to the exponential term $e^{\gamma \ph(z)}$ which is of course a priori not well defined for a generic distribution, but can be made sense of as in Chapter \ref{S:GMC} via Gaussian multiplicative chaos provided that $\ph$ is a logarithmically correlated Gaussian field.

Given the action $S(\varphi)$, by analogy with \eqref{E:gibbs}, one is led to formally define the associated Gibbs measure 
\begin{equation}\label{PolyakovMeasure}
\mathbf{P}(\ph) = \exp ( - S( \ph) ) \text{D} \ph,
\end{equation}
on a for now unspecified space of generalised functions defined on $\Sp$ (or $\hat \C$). The temperature has been set to 1 for simplicity; other choices lead to an equivalent theory since 
\[ \beta S(\varphi; \gamma, \mu) = S(\sqrt{\beta}\varphi, \frac{\gamma}{\sqrt{\beta}}, \beta \mu)\]
for $\beta>0$. The crucial thing to notice in \eqref{PolyakovMeasure} is that the choice of the reference measure $\text{D} \ph$, which should be heuristically viewed as a kind of Lebesgue (uniform) measure over the space of fields on $\mathbb{S}$ (or $\hat\C$), is \emph{not} specified precisely. 

In this chapter we will detail how \cite{DKRV} nevertheless succeeded in assigning a meaning to this Gibbs measure $\mathbf{P}$, which we will refer to in the rest of this chapter as the \textbf{Polyakov measure}.
\ind{Polyakov measure}
Note that this will in fact be an \emph{infinite measure} (in particular, not a probability measure). When we integrate this against an observable $F$ we will more typically write $\int F(\ph) \mathbf{P}(\dd \ph)=\langle F \rangle$ in agreement with the physics convention. We will compute these ``expectations'' for particular choices of $F$, and these will define the correlation functions of the theory. Informally these $F$ will be of the form $F(\ph)=\exp(\alpha \ph(z))$ and products thereof.

\subsection{Spherical GFF}

\ind{GFF!on the sphere}

\ind{GFF!on a compact manifold}

A key idea of \cite{DKRV} is to give meaning to the Polyakov measure in \eqref{PolyakovMeasure}  by combining the term $\exp( - \int |\nabla^g \ph(z)|^2 v_g(\dd z))$ and the reference measure $\text{D}\ph$, and then suitably \emph{reweighting} the resulting measure to account for the remaining terms on the right hand side of \eqref{E:Polyakov_informal}. In view of Theorem \ref{T:dGFF_Dirichlet} it is natural to want to interpret the product $\exp( - \int |\nabla^g \ph(z)|^2 v_g(\dd z) ) \text{D} \ph$ as the law of a Gaussian free field. However, in the absence of a boundary on which to impose boundary conditions, one has to choose a suitable version of the Gaussian free field, for which there is no obvious candidate at first sight.
In \cite{DKRV}, David, Kupiainen, Rhodes and Vargas made a simple but ingenious proposal, which consists of two steps. The first step requires us to define the Gaussian free field with zero average on the sphere $\mathbb{S}$ (or \textbf{spherical GFF} for short). In fact, one can do this on any compact surface $(\Sigma,g)$, in which case we speak of the \textbf{zero average GFF on $\Sigma$} with respect to the metric $g$. We will explain the construction in this generality since it is not more difficult. The details are reminiscent of other Gaussian free fields discussed in the book (see Chapter \ref{S:GFF} and the discussion of the \textbf{Neumann GFF} that will appear in Chapter \ref{S:SIsurfaces}). The reader who is keen to get on with the rigorous definition of Polyakov measure is encouraged to skip straight to Section \ref{SSS:polyakovdef}, where the second step is described.

\ind{GFF!Zero average}

\subsubsection{Laplacian on a compact manifold}
\medskip Let $(\Sigma, g)$ be a connected compact surface (i.e, a two dimensional, closed, connected and bounded Riemannian manifold), where $g$ refers to the Riemannian metric on $\Sigma$; in particular, $\Sigma$ has no boundary. The Riemannian structure induces a Laplace operator which we will denote by
{$\Delta^{\Sigma,g}$}\indN{{\bf Geometry}! $\Delta^{\Sigma,g}$; Laplace operator on Riemannian manifold $(\Sigma,g)$}. For instance, if $(\Sigma,g)=(\Sp,g)$ and $\hat g$ is the metric  obtained after pushing forward $g$ via the fixed conformal isomorphism $\psi:\Sp\to \hat\C$ of \eqref{E:psi}, then  for smooth $f$ defined on $\Sp$ we simply have 
$$
\Delta^{\Sp, g}f(z)   =\frac1{\hat g(\psi(z))} \Delta [f\circ \psi^{-1}](\psi(z))\, ; \quad z\in \Sp,,
$$
where $\Delta=\frac{\partial^2}{\partial x^2}+\frac{\partial^2}{\partial y^2}$ is the usual Laplacian on $\C$.
In other words, the  Brownian motion $X$ on $\Sp$ {with respect to $g$ (that is, the diffusion with infinitesimal generator $\Delta^{\Sp,g}$)} can be obtained by performing a time change to the standard Euclidean Brownian motion on $\mathbb{C}$
\begin{equation}\label{clock_sphere}
X_t = B_{F^{-1} (t)} ; \quad  F(t) = \int_0^t \hat g(B_s) \dd s.
\end{equation}
and mapping back to $\Sp$ via the inverse of $\psi$.
(A similar recipe, properly interpreted, may be used to define \emph{Liouville Brownian motion}, see \cite{diffrag} and \cite{GRV}.)
\ind{Liouville Brownian motion}

It can be seen that on a compact connected surface, the negative Laplacian $-\Delta^{\Sigma,g}$ has a discrete spectrum, which we denote by $$0 = \lambda_0 < \lambda_1 \le \lambda_2 \ldots \uparrow +\infty,$$
with each distinct eigenvalue repeated according to its multiplicity. By the Sturm--Liouville decomposition (see, for example, \cite[VI.1]{Chavel}) we can assume that the corresponding eigenfunctions $e_n$ form an orthonormal basis of $L^2(\Sigma; v_g)$. %, where we write $v_g$ for the volume measure on $\Sigma$ associated with $g$.  
Note that this has bounded total mass, since $\Sigma$ is bounded. The eigenvalue $\lambda_0 = 0$ is associated to the constant eigenfunction $e_0 = 1/\sqrt{ v_g(\Sigma)}$, corresponding to the fact that the Brownian motion on $\Sigma$ with respect to $g$ converges to the uniform distribution.

\medskip We will give several equivalent definitions of the GFF with zero $v_g$ average in $\Sigma$. The first one is as a random series. Before we give this definition, we briefly introduce the function space in which this series will converge, which is a variant of the Sobolev space $H^{s}(D)$ discussed in Chapter \ref{S:GFF}. %when introducing the GFF with Dirichlet boundary conditions.

A distribution on $(\Sigma,g)$ is simply a continuous linear functional on \emph{test functions}, where test functions are simply smooth functions on $(\Sigma, g)$ (since the space is bounded). Here continuity refers to the usual topology on test functions, meaning uniform convergence of derivatives of all orders. If $f$ is a distribution on $(\Sigma, g)$ and $\phi$ is a test function, we write $$(f, \phi)_g =\text{``} \int_\Sigma f(x)\phi(x) v_g(\dd x) \, \text{''} $$ for the action of $f$ on $\phi$. This is exactly the $L^2(\Sigma,g)$ inner product\indN{{\bf Inner products}! $(\cdot, \cdot)_g$; $L^2$ inner product on Riemannian manifold $(\Sigma,g)$} when $f$ is itself in $L^2(\Sigma,g)$ and in particular depends on the metric $g$. %Here we include the subscript since if we view a smooth function $f$ as a distribution then its action will correspond to the right hand side above, which in particular depends on $g$.
Note that the smooth function 1 is a valid test function, so that the total integral on $\Sigma$ of a distribution $f$ is well defined; we will write
$$
(f,1)_g = : v_g(f)
$$
for this integral and refer to it as the average of $f$ on $(\Sigma,g)$.
When this integral is equal to zero we say that the distribution has \textbf{zero average}. For $s \in \R$, we define $H^s (\Sigma,g)$ to be the space of zero average distributions $f$ on $\Sigma$ such that
$$
\sum_{n \ge 1} (f, e_n)_g^2 \lambda_n^s < \infty.
$$\indN{{\bf Function spaces}! $H^s(\Sigma,g)$; Sobolev space of index $s$ on  $(\Sigma,g)$}
Note that $e_n$ is smooth for every $n \ge 1$ so is a valid test function. It is not hard to see that the left hand side defines a Hilbert space norm $\|\cdot \|_{H^s(\Sigma,g)}$ on zero average distributions, and that convergence in that Hilbert space implies convergence in the sense of distributions. 
% In fact, the Sobolev embedding theorem tells us that for any $\alpha\in (0,1)$ if $s$ is large enough,
%\begin{equation}
%	\label{eq:sob_embedding}
%	\|f\|_{H^s(\Sigma,g)}\ge \sup_{x\ne y\in \Sigma}\frac{|f(x)-f(y)|}{d_g(x,y)^\alpha}.
%\end{equation}
%A proof can be found in \cite[Theorem 2.20]{Aubin}.
%In particular if $f\in H^s(\Sigma, g)$ for $s$ large enough then, after possible modification on a set of $v_g$-measure zero, $f$ defines a (H\"{o}lder) continuous function.
Note also that changing the metric $g$ to a conformally equivalent one as in \eqref{E:metric} leads to the same Sobolev spaces $\{H^s(\Sigma,g)\}_{s\in \R}$ in the sense that if $f\in H^s(\Sigma,\tilde{g})$, then $f-\bar{v}_g(f)\in H^s(\Sigma,g)$.

Let us also record that we have the Gauss--Green formula on $(\Sigma,g)$: that is, for any twice continuously differentiable functions $u,w$ on $(\Sigma,g)$ with zero average,
\begin{equation}
	\label{E:GG_manifold}
	\int_\Sigma u(x) \Delta^{\Sigma,g} w(x) v_g(\dd x)= \int_\Sigma w(x) \Delta^{\Sigma,g} u(x) v_g(\dd x).
\end{equation}
See for example,
\cite[(29); Chapter 1]{Aubin}.

\subsubsection{Definition of the zero average GFF on \texorpdfstring{$(\Sigma, g)$}{TEXT}}\label{S:ZAGFF}

We can now introduce the GFF with zero $v_g$ average on $\Sigma$:
\begin{definition}
\label{D:sphereGFF}
Let $(X_n)_{n\ge 1}$ denote a sequence of i.i.d. standard normal random variables. The GFF with zero average on $\Sigma$ is the random distribution $\mathbf{h}^{\Sigma, g}$ on $\Sigma$ obtained from the series
$$
\mathbf{h}^{\Sigma, g} = \sum_{n=1}^\infty \frac{X_n} {\sqrt{\lambda_n}} e_n ,
$$
which converges almost surely in any of the spaces $H^{s}(\Sigma,g)$ with $s<0$, and hence in the space of (zero average) distributions.
As usual we will also write $h^{\Sigma,g} = \sqrt{2\pi} \mathbf{h}^{\Sigma,g}$.
\end{definition}

The convergence of the series in this definition follows from Weyl's law (see \cite[VI.4, page 155]{Chavel}) in a manner similar to Lemma \ref{lem:weyl}.\ind{Weyl law} An argument similar to \eqref{eqn:exp_fs} also shows that the convergence of this series would also hold if we replaced $\{\lambda_n^{-1/2}e_n\}_{n\ge 1}$ with any orthonormal basis of $H^1(\Sigma, g)$. Furthermore, if $f \in C^\infty(\Sigma,g) $  has zero $v_g$ average, then
\begin{equation}\label{varianceGFFcentered}
\var ((\textbf{h}^{\Sigma,g}, f)_g) = \sum_{n=1}^\infty \frac{(f, e_n)_{L^2(\Sigma,g)}^2}{\lambda_n} = \| f\|^2_{H^{-1}(\Sigma,g)}
\end{equation}
where the right hand side indeed does not depend on the basis.

Observe  that, as in Theorem \ref{T:ito}, this allows us to alternatively define $\textbf{h}^{\Sigma,g}$ to be a stochastic process indexed by $ H^{-1}(\Sigma,g)$. More precisely:
$$\{ (\textbf{h}^{\Sigma,g},f)_g\}_{f\in H^{-1}(\Sigma,g)}$$
defines a Gaussian stochastic process indexed by $H^{-1}(\Sigma,g)$, with $\mathbb{E}((\mathbf{h}^{\Sigma,g},f)_g)=0$ for all $f$ and
$$\mathrm{cov}((\mathbf{h}^{\Sigma,g},f)_g (h^{\Sigma,g},\tilde{f})_g)=(f,\tilde{f})_{H^{-1}(\Sigma,g)}=\sum_{n=1}^\infty \frac{(f,e_n)_g(\tilde{f},e_n)_g}{\lambda_n}; \quad f,\tilde{f}\in H^{-1}(\Sigma,g).$$
By \eqref{varianceGFFcentered} and the polarisation identity, the restriction of this stochastic process to $f\in C^\infty(\Sigma,g)$ agrees with the definition of $\textbf{h}^{\Sigma,g}$ as a zero average distribution in  \cref{D:sphereGFF}.

\medskip It will also be useful to have an description for the above covariance structure in terms of a covariance kernel. This will be the  \textbf{centred Green function} of the $g-$Brownian motion on $\Sigma$, that is, the Brownian motion on $\Sigma$ with respect to $g$ (see \cref{L:varianceGFFcenteredGreen} below). To define the centered Green function, let $(X_t, t \ge 0)$ denote  this Brownian motion (so its generator is $\Delta^{\Sigma,g}$), and let $p^{\Sigma,g}_t (x, y)$ denote the \textbf{heat kernel}, which is characterised (for a fixed $x\in \Sigma$ and as a function of $y$), as the density of the law of $X_t$ when started from $x$ with respect to $v_g$. A standard fact (see \cite[VI.1, equation (13)]{Chavel}) is that we may decompose $p_t^\Sigma(x,y)$ according to the orthonormal basis of eigenfunctions $e_n $ of $L^2(\Sigma;v_g)$ as
\ind{Green function!on a compact surface}
\begin{equation}\label{E:evdecomposition}
p^{\Sigma,g}_t(x,y) = \sum_{n=0}^\infty e^{-\lambda_n t} e_n(x) e_n(y).
\end{equation}
\indN{{\bf Brownian motion}! $p^{\Sigma,g}_t(\cdot,\cdot)$; heat kernel on a Riemannian manifold $(\Sigma,g)$}
This series converges absolutely and uniformly on $\Sigma \times \Sigma$ for any $t>0$ (see again \cite[VI.1]{Chavel}). In fact, see \cite[Remark 10.15 and (10.51)]{Grigoryan}, we have that for all $n\ge 0$,
\begin{equation} \label{E:ptconv}  \sup_{t>t_0}\sup_{x,y\in \Sigma} e^{\frac{\lambda_{n+1}}{2}t_0}| p_t^{\Sigma, g}(x,y)-\sum_{k=0}^n e^{-\lambda_k t} e_k(x) e_k(y)| \le C(t_0) \end{equation}
where the constant $C(t_0)$ depends only on $t_0$ and in particular not on $n$.
As a consequence, we can define the associated Green function with zero average
$$
G^{\Sigma,g}(x,y) : = \int_0^\infty [p^{\Sigma,g}_t(x,y) - \frac{1}{v_g(\Sigma)} ] \dd t.
$$ \indN{{\bf Green functions}! $G^{\Sigma,g}(\cdot,\cdot)$; Green function with zero average on $(\Sigma,g)$}
When $x\ne y$, $p^{\Sigma,g}_t(x,y)$ is bounded as $t\downarrow 0$,  and so \eqref{E:ptconv} and \eqref{E:evdecomposition} imply that the above integral converges. In fact, they ensure that \begin{equation}\label{Green_ev}
	G^{\Sigma,g}(x,y) = \sum_{n=1}^\infty \frac1{\lambda_n} e_n(x) e_n(y),
\end{equation}
with the convergence of the series holding pointwise. Observe that if $x\ne y$ then by its definition as the transition density of the $g$-Brownian motion on $\Sigma$, $p^{\Sigma,g}_t(\cdot,y)$ is uniformly bounded as $t\to 0$ on a neighbourhood of $x$. This together with \eqref{E:ptconv} implies that $G^{\Sigma,g}(\cdot, y)$ is continuous at $x$. In other words, $G^{\Sigma,g}$ is continuous away from the diagonal. %Note also that, since we have subtracted $1/v_g(\Sigma)$, the average $v_g(G^{\Sigma,g}(\cdot, y))=v_g (G^{\Sigma,g}(x, \cdot))$ is indeed equal to $0$ for every $x \in \Sigma$.

Note also that \eqref{E:ptconv} implies the upper bound $|G^{\Sigma,g}(x,y)|\le \int_0^1 p_t^{\Sigma,g}(x,y) \dd t + C$ where $C$ is a constant not depending on $x,y\in \Sigma$. In particular, if $f$ is uniformly bounded on $\Sigma$, then $\int_\Sigma G^{\Sigma,g}(x,y)f(y)v_g(\dd y)$ converges absolutely and
\begin{equation}\label{E:G_dist_cont} \int_\Sigma |G^{\Sigma,g}(x,y)f(y)| v_g(\dd y) \le \int_0^1 \E(|f(B_t)|) \dd t + C \mathrm{vol}(G) \sup_\Sigma|f| \le C' \sup_\Sigma |f|.\end{equation}
This means that $G^{\Sigma,g}(x,\cdot)$ defines a distribution on $(\Sigma,g)$ for every $x$. As should be expected from \eqref{Green_ev}, we have
$$v_g(G^{\Sigma,g}(x,\cdot))=\int_\Sigma G^{\Sigma,g}(x,y)v_g(\dd y)=0$$
(that is, $G^{\Sigma,g}(x,\cdot)$ does have zero average), and
$$\int_{\Sigma} G^{\Sigma,g}(x,y)e_n(x) v_g(\dd y) = \frac{e_n(x)}{\lambda_n} \text{ for all } n\ge 1 $$
which can both be justified rigorously using \eqref{E:ptconv}. As a result, $G^{\Sigma,g}$ is an ``inverse'' of (minus) the Laplacian in the following sense:

Let $\bar v_{ g}( \dd x) = v_{ g} (\dd x) / v_{ g} ( \Sigma)$ be the normalised volume measure associated to the metric $ g$ on $\Sigma$, and for $f: \Sigma\to \R$ a bounded measurable function, let
 $$\bar v_{ g}(f) : = \frac1{v_{ g}(\Sigma)}\int_\Sigma f(x) v_{ g}(\dd x). %= \frac1{v_{ g}(\Sigma)}\int_\Sigma f(x) \hat g(x) \dd x.
 $$

\begin{lemma}\label{L:gsg_lap}
	If $\phi$ is a smooth function on $(\Sigma,g)$, then for every $x\in \Sigma$,
	 \begin{equation}
		\label{E:gsg_lap}
		\int_\Sigma G^{\Sigma,g}(x,y) \Delta^{\Sigma,g}\phi (y) v_g(\dd y)= -\phi(x)+\bar{v}_g(\phi).
	\end{equation}
\end{lemma}

\begin{rmk}
	This result should not be surprising, since with respect to the basis $\{e_n\}_{n\ge 1}$ of zero average functions, $-\Delta$ is a ``diagonal'' operator with diagonal entries $(\lambda_n)_{n\ge 1}$ and by \eqref{Green_ev}, the Green function can also be viewed as a diagonal operator with diagonal entries $(\lambda_n^{-1})_{n\ge 1}$. \end{rmk}

\begin{proof}
	We make use of the Sobolev embedding theorem for compact manifolds,  \cite[Theorem 2.20]{Aubin}, which implies in particular that for $f$ smooth, the series $\sum_{n=0}^\infty (f,e_n)_{L^2(\Sigma,g)}  e_n $ converges uniformly to $f$ on $\Sigma$. This applies directly to $\phi$ as in the statement of the lemma, and also to $\Delta^{\Sigma,g}\phi$, where $(\Delta^{\Sigma,g}\phi,e_n)_{L^2(\Sigma,g)}=-\lambda_n (\phi,e_n)_{L^2(\Sigma,g)}$ for each $n\ge 0$ by \eqref{E:GG_manifold}. This means (using \eqref{E:G_dist_cont}) that the left hand side of \eqref{E:gsg_lap} is equal to $$-\sum_{n=1}^\infty (G^{\Sigma,g}(x,\cdot),e_n)_{L^2(\Sigma,g)} \lambda_n (\phi,e_n)_{L^2(\Sigma,g)}=-\sum_{n=1}^\infty (\phi,e_n)_{L^2(\Sigma,g)} \phi(x) = -\phi(x)+\bar{v}_g(\phi)$$ as required.
\end{proof}

\begin{rmk}\label{R:gsg_lap}
In fact, \eqref{E:gsg_lap} can be extended, by approximation, to the case where $\phi$ is only assumed to be twice continuously differentiable on $(\Sigma, g)$. See \cite[Proposition 4.14]{Aubin}.
\end{rmk}

Let us now finally use the above to relate this Green function to the zero average GFF $\mathbf{h}^{\Sigma,g}$: 
\begin{lemma}\label{L:varianceGFFcenteredGreen}
	For a smooth function $f$ on $(\Sigma,g)$,
\begin{equation}
\label{E:varianceGFFcenteredGreen}
\var ((\textbf{h}^{\Sigma,g}, f)_g) = \int_\Sigma \int_\Sigma f(x) G^{\Sigma,g}(x,y) f(y) v_g (\dd x) v_g(\dd y).
%\text{``}\iint_{\Sigma^2} G^{\Sigma,g}(x,y) f(x) f(y) v_g(\dd x) v_g(\dd y)\text{''}
\end{equation}
 In other words, the Green function is the covariance kernel of $\mathbf{h}^{\Sigma,g}$. \end{lemma}
\begin{proof}
	By the Sobolev embedding theorem again, we have that for $f$ smooth with zero $v_g$ average, the series  $$\sum_{n=1}^\infty \lambda_n^{-1} e_n(x) (f,e_n)_{L^2(\Sigma,g)}$$ converges to a smooth $\phi$ with $\Delta^{\Sigma,g} \phi = f$. Hence $\int_\Sigma G^{\Sigma,g}(x,y) f(y) v_g(\dd y)=\phi(x)$ and thus since $(\phi,e_n)_{L^2(\Sigma,g)}=\lambda_n^{-1}(f,e_n)_{L^2(v_g)}$ for each $n$:
	$$\int_\Sigma \int_\Sigma f(x) G^{\Sigma,g}(x,y) f(y) v_g (\dd x) v_g(\dd y) = \|f\|_{H^{-1} (\Sigma,g)}^2.
	$$
	Combining with \eqref{varianceGFFcentered}, we reach the conclusion.
\end{proof}

\subsubsection{The spherical case}
Let us now specialise to the case where $\Sigma= \Sp$ is the sphere, and the metric is still as in \eqref{E:metric}. 
In this case, we can easily observe how $G^{\Sp,g}$ changes when we apply a M\"{o}bius transformation to $\Sp$.
\begin{lemma}\label{L:Gmob}
	Suppose that $m:\Sp \to \Sp$ is a M\"{o}bius transformation, and let $m_* g$ \indN{{\bf Miscellaneous}! $m_*g$; pushforward of a metric $g$ by a map $m$}be the pushforward of $g$ under $m$. Then
	\[G^{\Sp,m_* g}(x,y)=G^{\Sp,g}(m^{-1}(x),m^{-1}(y))\]
	for all $x\ne y\in \Sp$.
\end{lemma}
\begin{proof}
	It follows from a straightforward calculation that if $(e_n)_{n\ge 0}$ are an orthonormal basis of $L^2(\Sigma, v_g)$, such that $\Delta^{\Sigma,g}e_n=-\lambda_n e_n$ for $n\ge 1$, then $(e_n\circ m^{-1})_{n\ge 0}$ are an orthonormal basis of $L^2(\Sigma, v_{m_*g})$, such that $\Delta^{\Sigma,m_*g}(e_n\circ m^{-1})=-\lambda_n (e_n\circ m^{-1})$ for $n\ge 1$. The result then follows from \eqref{Green_ev} (it is easy to check that the series definition \eqref{Green_ev} cannot depend on the choice of orthonormal basis of eigenfunctions, since each of the eigenspaces is finite dimensional).
\end{proof}

Using \eqref{E:varianceGFFcenteredGreen}, this implies that $h^{\Sp,g}$ transforms in the following way:
\begin{cor}\label{C:hs_mobius}
	Suppose that $m:\Sp\to \Sp$ is a M\"{o}bius transformation. Then 
	\[
	(\mathbf{h}^{\Sp,g},f\circ m)_g=(\mathbf{h}^{\Sp,m_*g},f)_{m_* g}
	\]
	for $f\in C^\infty(\Sp)$. In other words, if 
	we define a zero average distribution $\mathbf{h}^{\Sp,g}\circ m^{-1}$ on $(\Sp,m_*g)$ by setting $({\mathbf{h}}^{\Sp,g}\circ m^{-1}, f)_{m_*g}=({\mathbf{h}}^{\Sp,g},f\circ m)_g$ for $f\in C^\infty(\Sp)$. Then we have that $${\mathbf{h}}^{\Sp,g}\circ m^{-1} \overset{(\mathrm{law})}{=} {\mathbf{h}}^{\Sp,m_* g}$$ as zero average distributions on $(\Sp,m_*g)$.
\end{cor}

As mentioned previously, in order to do explicit computations in this case, we will want to reparametrise $\Sp$ by the extended complex plane $\hat \C$. Recall from the conversation around \eqref{E:psi} that $\psi:\Sp\to \hat{\C}$ is a fixed conformal isomorphism (for example, stereographic projection) and we write (with an abuse of notation) $\hat g(z)$ for the pushforward of a metric $g$ on $\Sp$ under $\psi$. Then by the same proof as for \cref{L:Gmob}, we have that
\begin{equation}\label{E:hChS} ({\mathbf{h}}^{\hat\C,\hat g},f)_{\hat g}=({\mathbf{h}}^{\Sp,g},f\circ \psi)_g \end{equation}
for all smooth functions $f$ on $\hat{\C}$. In particular, 
%\[ \var\big((h^{\hat\C,g},f)_g\big) = \int_{\C\times \C} G^{\Sp,g}(\psi^{-1}(x),\psi^{-1}(y))f(x)g(y) v_g(\dd x) v_g(\dd y)\]
%for all such $f$. That is, 
\begin{equation}\label{E:GCGS} G^{\hat \C,\hat g}(x,y)=G^{\Sp,g}(\psi^{-1}(x),\psi^{-1}(y)) \text{ for } x \ne y\in \hat{\C}. \end{equation}

\begin{rmk}\label{R:hs_mobius}
$G^{\hat\C,\hat g}$ and $\mathbf{h}^{\hat \C,\hat g}$ therefore satisfy the same transformation rules under M\"{o}bius maps. If $m:\hat{\C}\to \hat{\C}$ is a M\"{o}bius transformation, that is, of the form $m(z)=(az+b)/(cz+d)$ with $ad-bc=1$, then
\[ 
G^{\hat\C, m_*\hat g}(x,y)=G^{\hat \C, \hat g}(m^{-1}(x),m^{-1}(y)) \, \quad x\ne y\in \C
\]
and 
\[
(\mathbf{h}^{\hat\C,\hat{g}},f\circ m)_{\hat{g}}=(\mathbf{h}^{\hat \C,m_*\hat g},f)_{m^*\hat{g}} \, \quad f\in C^\infty(\hat{\C})
\]
\end{rmk}
\ind{Green function!on the sphere}
 
 \begin{lemma}
 \label{L:sphereWholePlaneGFF} 
 %Let $\bar v_{\hat g}( \dd x) = v_{\hat g} (\dd x) / v_{\hat g} ( \hat \C)$ be the normalised volume measure associated to the metric $\hat g$ on $\hat{\C}$, and for $f: \hat \C\to \R$ a bounded measurable function, let
 %$$\bar v_{\hat g}(f) : = \frac1{v_{\hat g}(\hat \C)}\int_\C f(x) v_{\hat g}(\dd x) = \frac1{v_{\hat g}(\hat \C)}\int_\C f(x) \hat g(x) \dd x.$$
 %\index[ops]{$\bar{v}_{g}(f)$; average of a distribution on $(\Sigma,g)$ with respect to the volume form $v_g$} 
 We have, for all $x\ne y \in \C$,
 \begin{equation}\label{E:gSequiv}
 G^{\hat\C,\hat g} (x, y) = \frac1{2\pi} \Big[ - \log (|x-y|)
 + \bar v_{\hat g} \big( \log (|x-\cdot|)\big) +  \bar v_{\hat g} \big( \log (|y-\cdot|)\big)  - \theta_{\hat g}\Big],
 \end{equation}
 where
 \begin{equation}\label{E:thetag} \theta_{\hat g} = \iint_{\C^2}\log |x-y| \bar v_{\hat g}(\dd x) \bar v_{\hat g}(\dd y).\end{equation}
 \end{lemma}

In addition, we write
\[
G^{\hat{\C},\hat g}(x,\infty):=\lim_{R\to \infty} G^{\hat{\C},\hat g}(x,R)
\]
which is also well defined for $x\in \C$, by \eqref{E:GCGS} and since $G^{\Sp,g}(\psi^{-1}(x),y)$ is continuous at the pole $y=\psi^{-1}(\infty)$.

\begin{rmk}
	Note that $G^{\hat\C,\hat g}$ is \emph{not} harmonic away from the diagonal (in contrast to, say, the case of the Dirichlet Green function). Indeed, it follows from \eqref{E:gsg_lap} that for fixed $x$, $\Delta G^{\hat \C,\hat g}(x,\cdot)=-\delta_x+\bar v_{\hat g}$ as a distribution. This is necessary, due to the requirement that $G^{\hat\C,\hat g}$ has zero average with respect to $v_{\hat g}$.
\end{rmk}

\begin{proof}
Fix $x \in \C$ and define the function $y \mapsto F(y)$ as the difference between the left hand the right hand side of \eqref{E:gSequiv}.
Note that at the moment it is not clear whether $F$ is defined at $x\ne y$, but we can still view it as a distribution on $\C$.
We will show below that $F$  is harmonic in the sense of distributions on $\C$; that is,
%We start with the boundedness. {In fact, since we will show harmonicity (and in particular boundedness on compacts) of $F$, we only need to justify boundedness at infinity.} It is not hard to see that $\bar v_g (\log | y - \cdot|)  = \log |y| + O(1)$ as $y \to \infty$. This implies that the function on the right hand side of \eqref{E:gSequiv} is bounded at infinity. On the other hand, the boundeness of $G^{\Sp,g}(x,\cdot)$ at $\infty$ follows from Lemma \ref{L:Gmob}.
%It is also immediate that $\bar v_g ( F) = 0$. It therefore only remains to check that $F$ is harmonic on $\C$, % (for the spherical Laplacian $\Delta^\Sp$ or equivalently for the Euclidean Laplacian $\Delta$). In order to appeal to elliptic regularity on the plane rather than on the sphere, it is slightly more convenient (if perhaps less intrinsic) to do the calculation with $\Delta$.
%and by elliptic regularity, for this it is enough to show that $\Delta F=0$ in the sense of distributions on $\C$. That is, it suffices to show
for any compactly supported smooth test function $f$ on $\C$:
\begin{equation}\label{E:harmonic}
	({F},\Delta f)=\int_\C F(y) \Delta f(y) \dd y = 0.
\end{equation}
Admitting \eqref{E:harmonic}, it then follows from elliptic regularity that $F$ can be extended continuously to the point $x$ and with this extension is a classically harmonic function on $\C$. It is also immediate that $\bar v_{\hat g} ( F) = 0$, and $F$ is bounded at infinity (the boundedness holds since $\bar v_{\hat g} (\log | y - \cdot|)  = \log |y| + O(1)$ as $y \to \infty$ and since $G^{\hat\C,\hat g}(x,\cdot)=G^{\Sp,g}(\psi^{-1}(x),\psi^{-1}(\cdot))$ is bounded at infinity by definition). This means that $F$ must be identically zero, which completes the proof of the lemma.

So it remains to prove \eqref{E:harmonic}. For this, we first observe directly from \eqref{E:gsg_lap} that
\begin{equation}\label{E:greensphere1}
	\int_\C G^{\hat\C,\hat g}(x,y) \Delta f (y) \dd y = \int_\Sp G^{\Sp,g}(x,y) \Delta^{\Sp,g} f(y) v_g(\dd y) = -f(x) + \bar v_g(f).
\end{equation}
On the other hand, we know that
\begin{equation}\label{E:logharm}
-\frac1{2\pi} \int_\C \log|x-y| \Delta f(y) \dd y = - f(x)
\end{equation}
since $\Delta \tfrac1{2\pi}\log|x-\cdot|=\delta_x(\cdot)$ in the distributional sense. We are left to compute the integral
$$ \frac1{2\pi} \int_{\C}  (\bar v_{\hat g} \big( \log (|x-\cdot|)\big) +  \bar v_{\hat g} \big( \log (|y-\cdot|)\big)  - \theta_{\hat g}) \Delta f(y) \, \dd y =   \frac1{2\pi}  \int_\C \bar v_g \big( \log|y-\cdot|\big) \Delta f(y) \, \dd y $$ (with the equality coming from the fact that two of the terms in the integral on the left hand side do not depend on $y$, and the Gauss--Green formula). Here we can appeal to Fubini, since $f$ is smooth and compactly supported, and write this as
\begin{equation}
 \frac1{2\pi} \int_\C \int_\C \log(|y-w|)\Delta f(y) \dd y\,  \bar{v}_{\hat g}(\dd w) = \int_\C f(w) \bar v_{\hat g} (\dd w) = \bar v_{\hat g}(f)
 \label{E:greensphere2}
\end{equation}
by \eqref{E:logharm}. Combining \eqref{E:greensphere1} and \eqref{E:greensphere2}, we conclude that
$$(\tilde{F},\Delta f)= f(x)-\bar{v}_{\hat g}(f)-f(x)+\bar{v}_{\hat g}(f)=0$$
as required.
\end{proof}

As a corollary, we obtain the following expression for the circle average variance. Recall that $h^{\Sp,g} = \sqrt{2\pi} \mathbf{h}^{\Sp,g}$ (similarly $h^{\hat\C,\hat g}=\sqrt{2\pi} \mathbf{h}^{\hat\C,\hat g }$) and that from \eqref{E:varianceGFFcenteredGreen} and the expression of the spherical Green function $G^{\hat\C,g}$, circle averages of $h^{\hat \C,\hat g}$ are well defined. Let $h_\eps^{\hat\C,\hat g}$ denote the circle average of $h^{\hat\C, \hat g}$ at (Euclidean) distance $\eps$.

\begin{cor}
\label{C:spherecircleaverage}
As $\eps \to 0$, we have
\begin{equation}\label{E:spherecircleaveragevariance}
\var ( h_\eps^{\hat\C,\hat g}(z)) = \log \frac1\eps + v(z) + o(1);
\end{equation}
where $v(z) = 2 v_g ( \log |z - \cdot|) - \theta_g$.
\end{cor}

%When $g = \hat g$ (so $\rho = 0$ in \eqref{E:metric}), we simply write $G^{\Sp}$ instead of $G^{\Sp, \hat g}$.
The following explicit formula for the spherical metric Green function, $G^{\hat\C,\hat g_0}$ will also be useful.

\begin{lemma}
\label{L:covariance_hsphere}
We have
\begin{equation}\label{E:greensphere}
G^{\hat\C,\hat g_0} (x, y) = \frac1{2\pi}\big(-\log |x-y| - \frac14( \log \hat g_0(x) + \log \hat g_0(y) )  + \log 2 - 1/2\big)
\end{equation}
%\ellen{*do we want this second bit later, when we've defined the circle average and got rid of the $2\pi$s?*}
for $x\ne y \in \C$, and $G^{\hat\C,\hat g_0}(x,\infty)=(1/2\pi)(-(1/4)\log \hat{g}_0(x)+(1/2)\log 2-1/2)$ for $x\in \C$. Consequently, as $\eps \to 0$, uniformly in $z \in\C$,
\begin{equation}
\var (h_\eps^{\hat\C,\hat g_0}(z)) = \log \frac1{\eps} +\hat v(z)+ o(1)
\end{equation}
where
\begin{equation}\label{constantterm}
\hat v(z) = -  \frac12 \log \hat g_0(z) + \log 2 - \frac12 .
\end{equation}\
\end{lemma}

\begin{proof}
	Recall that by \eqref{E:greensphere2}, $\Delta \bar{v}_{\hat g_0}(\log|x-\cdot|)= 2\pi \bar{v}_{\hat g_0}(\cdot)$ in the sense of distributions on $\C$.
On the other hand, since 
	$$
	2\hat g_0(z) = R_{\hat g_0} \hat g_0(z) = - \Delta \log \hat g_0(z),
	$$
	by definition of the scalar curvature in \eqref{E:scalar} (and recalling that $R_{\hat g_0}=2$), for any smooth compactly supported $f$ on $\C$ we have  $$-\frac{1}{4} \int_\C \Delta f(y) \log \hat{g}_0(y) \dd y = -\frac{1}{4} \int_\C f(y) \Delta(\log \hat{g}_0(y)) \dd y = \frac{1}{2} \int_\C f(y) \hat{g}_0(y) \dd y = 2\pi \bar{v}_{\hat g_0}(f) $$ due to Gauss--Green. This implies that $\bar{v}_{\hat g_0}(\log|x-\cdot|)+\tfrac{1}{4}\log \hat{g}_0(x) $ is harmonic in $\C$, and the rest of the proof amounts to computing constants.

First, it is straightforward to check that $\bar{v}_{\hat g_0}(\log|y-\cdot|)+\tfrac{1}{4}\log(\hat g_0(y))\to \tfrac{1}{2}\log 2$ as $x\to \infty$, so that  by harmonicity $\bar{v}_{\hat g_0}(\log|y-\cdot|)+\tfrac{1}{4}\log(\hat g_0(y))\equiv \tfrac{1}{2}\log 2$. Also note that  $\bar{v}_{\hat g_0}(\log|x-\cdot|)$ has $\bar{v}_{\hat g_0}$ average $\theta_{\hat g_0}$ (by definition), and $-\frac{1}{4}\log \hat{g_0}$ has $\bar{v}_{\hat g_0}$ average $\tfrac{1}{2}-\tfrac{1}{2}\log 2$ (as can be shown by switching to polar coordinates and doing the integral explicitly). It therefore follows that $\theta_{\hat g_0}=1/2$ and we obtain the result.
\end{proof}

\begin{rmk} \label{R:constantcurvature}
More generally, if the curvature $R_{\hat g}$ is constant, then
$$ \bar{v}_{\hat g}(\log|x-\cdot|)+\frac{1}{2R_g}\log\hat{g}(x) \equiv \theta_{\hat{g}}+\frac{1}{2R_{\hat g}}\bar{v}_{\hat g}(\log(\hat g))
$$ so that 
$$
G^{\hat{\C},\hat{g}}(x,y)=\frac{1}{2\pi}\left(-\log|x-y|-\frac{1}{2R_{\hat{g}}}\log \hat{g}(x)-\frac{1}{2R_{\hat{g}}}\log \hat{g}(y)+c_{\hat g}\right)
$$
and		
$$
v(z) = -\frac1{R_{\hat g}} \log \hat{g}(z)  + c_{\hat g},
$$
where 
$$
c_{\hat g}=\theta_g+\frac{1}{R_{\hat g}}\bar{v}_{\hat{g}}(\log(\hat{g})). 
$$ 
\end{rmk}

The following special case will come in useful later on, when studying the behaviour of Liouville theory under M\"{o}bius transformations. The proof uses similar reasoning to above, and we leave it as a guided exercise.
\begin{lemma}\label{L:thetaRcmg}
	When $\hat{g}=m_*\hat{g}_0$, with $m$ a M\"{o}bius transform of $\hat{\C}$, we have 
	\begin{equation}
		\label{eq:thetamg}
		\theta_{m_*\hat g_0}=-\frac{1}{2}\bar{v}_{m_*\hat{g}_0}(\log(m_*(\hat g_0)))+\log(2)-\theta_{\hat g}
	\end{equation}
and
	\begin{equation}
	\label{eq:Rcmg}
	R_{m_*\hat{g}_0}\equiv R_{\hat g_0} \equiv 2 \quad ; \quad c_{m_*\hat g_0} = c_{\hat g_0}=\log(2)-\frac{1}{2}.
\end{equation} 
\end{lemma}
\begin{proof}
See	Exercise \ref{Ex:mobius}.
\end{proof}

Let us conclude this section by noting that $\mathbf{h}^{\Sp,g}$ for different metrics $g$ on $\Sp$ are simply recentrings of the same field. More precisely:

\begin{lemma} \label{L:hsnorms}
Suppose that $g_1,g_2$ are two metrics on $\Sp$ as in \eqref{E:metric}. Then
$${\mathbf{h}}^{\Sp,g_1}-\bar{v}_{g_2}({\mathbf{h}}^{\Sp,g_1}) \overset{(\mathrm{law})}{=} {\mathbf{h}}^{\Sp,g_2}$$
as zero average distributions on $(\Sp,g_2)$. Equivalently, 
$${\mathbf{h}}^{\hat \C,\hat g_1}-\bar{v}_{\hat g_2}({\mathbf{h}}^{\hat \C,\hat g_1}) \overset{(\mathrm{law})}{=} {\mathbf{h}}^{\hat \C,\hat g_2}.$$
\end{lemma}

\begin{proof}
	It is enough to prove the second statement. This can either be verified using the explicit expression for the covariance in Lemma \ref{L:sphereWholePlaneGFF}, or using the fact that $\var(h^{\hat \C,\hat g_2},f)_{\hat g_2}=\|f\|^2_{H^{-1}(\hat \C, \hat g_2)}$ for all smooth functions $f$ on $\C$ (similarly with $g_2$ replaced by $g_1$). We leave this to the reader as Exercise \ref{ex:shift}.
\end{proof}

\subsubsection{GMC on the Riemann sphere}

\ind{Gaussian multiplicative chaos}

Let $(\Sp, g)$ denote the sphere with a metric $g$ assumed to be conformally equivalent to the standard round metric $g_0$, and let $h=\sqrt{2\pi} \mathbf{h}^{\Sp, g}$ denote the (rescaled) Gaussian free field on $\Sp$ with zero average with respect to $v_g$. Associated with $h $ there is a notion of Gaussian multiplicative chaos $\cM_{h;g}$, formally described by
$$
\cM_{h;g} (\dd x) = e^{\gamma h(x) } v_g(\dd x),
$$
and understood rigorously as 
\begin{equation}\label{eq:GMCsphere}
\cM_{h;g} (\dd x) = \lim_{\eps \to 0} \eps^{\gamma^2/2} e^{\gamma h_\eps (x)  } v_g(\dd x),
\end{equation}
where $h_\eps(x)$ denotes the mollification of $h$ at scale $\eps$ with respect to the underlying metric  $g$. 
%We introduce the notation ${\cM}$ here to avoid confusion with the notation $\cM_h^\gamma$ which we occasionally use to emphasise the dependence on $\gamma$ when describing the (Euclidean) Gaussian multiplicative chaos of a field $h$). 
Existence of the above limit requires a slight extension of the theory presented in Chapter \ref{S:GMC} and more specifically Theorem \ref{T:conv}, since the latter deals with only logarithmically correlated fields on $\R^d$. Rather than go through the necessary adjustments (which however do not present any difficulty, see Remark \ref{R:extensionGMC}), it is equivalent and slightly simpler for our purposes to discuss the pushforward of $\cM_{h;g}$ to the extended complex plane $\hat \C$ (see also \cite{Cercle, DSHKS} for a discussion of GMC measures arising naturally from higher dimensional extensions of Liouville CFT). 

Hence, we rigorously define $\cM_{h;g}$ by conformally mapping $h$ to $\hat \C$ using the fixed conformal isomorphism $\psi:\Sp \to \hat \C$ (see \eqref{E:psi}), and then pulling back the resulting measure associated with $\hat g$. That is, if $\hat h  = h \circ \psi^{-1}$, we define for Borel sets $A\subset \Sp$, 
 \begin{equation}\label{E:LCFTshift}
\cM_{h;g} (A) = \lim_{\eps \to 0} \cM_{\hat h+ ({Q}/2) \log \hat g ;\eps}( \psi(A)):= \lim_{\eps\to 0} \int_{\psi(A)} \eps^{\gamma^2/2} e^{\gamma (\hat h +\tfrac{Q}2 \log \hat g)_\eps(z) } \dd z,
 \end{equation}\indN{{\bf Gaussian multiplicative chaos}! ${\mathcal{M}}_{h;g}(A)$; Gaussian multiplicative chaos of a field $h$ on a Riemannian manifold $(\Sigma,g)$}
where $\hat g$ denotes the pushforward of $g$ by $\psi$.
%(which corresponds to \eqref{spheremetric} when $g=g_0$). 
The subscript $\eps$ on the right hand side here denotes a mollified version of the field $\hat h + (Q/2) \log \hat g$ at (Euclidean) radius $\eps$. %(so this corresponds to the theory of \cref{S:GMC} with base measure $\sigma(\dd z)=\lim_{\eps\to 0} \eps^{\gamma^2/2}e^{\var((\hat h + (Q/2) \log \hat g)_\eps(z))}$). 
Notice that the field $\hat h$ needs to be shifted by $(Q/2) \log \hat g$, as specified by the change of coordinate formula: see Theorem \ref{T:ci}, and note that $|(\psi^{-1})'(z)| = \sqrt{\hat g(z)}$ for $ z\in \hat \C$.

For instance, when $g=g_0$, with this choice of normalisation, we have (using \eqref{constantterm}), 
\begin{align*}
\E ( \cM_{h;g_0}(\Sp)) &= \lim_{\eps \to 0} \eps^{\gamma^2/2} \int_{\hat \C} e^{\tfrac{\gamma^2}{2} \var ( \hat h_\eps(z) ) + \tfrac{\gamma Q}{2} \log \hat g_0(z) } \dd z\\
& = \int_{\hat \C}  e^{\tfrac{\gamma^2}{2} \hat v(z) + (1+ \tfrac{\gamma^2}{4} ) \log \hat g_0 (z) } \dd z\\
& = \int_{\C} e^{\tfrac{\gamma^2}{2} ( \log 2 - 1/2) } \hat g_0(z) \dd z\\
& = e^{\tfrac{\gamma^2}{2} ( \log 2 - 1/2) }  v_{g_0}(\Sp).
\end{align*}
One can check (also using \eqref{constantterm}) that this agrees with the limit of the expectation in the expression in \eqref{eq:GMCsphere}. 

The appearance of the (non-universal) constant $\log 2 - 1/2$ here is a consequence of our choice of normalisation for the GFF, which is required to have zero average. In the theory below, the choice of this additive constant does not play a role. 
%The reason for this exact definition, that is, the shift $\tfrac{Q}{2}\log g$ will be further explained in Section \ref{SS:LCFTshift}.

\subsection{Defining the Polyakov measure}

%\subsubsection{Definition}
\label{SSS:polyakovdef}

\medskip \textbf{Step 1: GFF on the sphere with mean zero.}   We fix a metric $g$ conformally equivalent to $ g_0$ as in \eqref{E:metric}. We then consider the scaled version $h^{\Sp,g}  = \sqrt{2\pi} \h^{\Sp,g}$ and, as usual, write $h^{\hat \C, \hat g} $ for the same field parametrised by the extended complex plane, that is, $h^{\hat \C, \hat g}  = h^{\Sp, g} \circ \psi^{-1}$ as in \eqref{E:hChS} where $\psi: \Sp \to \hat \C$ is the conformal isomorphism (for example, stereographic projection) that we fixed in \eqref{E:psi}. 
\ind{GFF!on the sphere}

\medskip \noindent \textbf{Step 2: the Lebesgue shift.} As observed earlier, there is nothing canonical about normalising the field to have zero average. In fact, the expression for $ \exp ( - \int_\C |\nabla^g \ph(z) |^2 v_g ( \dd z))$ is clearly invariant under shifting the field $\ph$  by an additive constant.  %The solution in these chapters was to view the fields as being defined modulo an additive constant. 
The heart of the construction of \cite{DKRV} is to start with $h^{\Sp,g}$, and then add a constant $c$ ``distributed'' according to Lebesgue measure on $\R$, the unique Radon measure (up to normalisation) on $\R$ invariant under additive shift. We then define the resulting ``law'' $\lambda$ to be $\exp ( - \int_\C |\nabla^g \ph|^2 v_g( \dd z) ) \text{D} \ph$. Before we explain precisely what this means, we mention that we write ``law'' in the above (and below) in quotation marks since it is in fact a measure of infinite total mass (because Lebesgue measure on $\R$ has infinite total mass). Let us be more specific. We define $H^{-1}(\Sp)$ to be the subspace of distributions on $(\Sp,g)$ of the form $\{\ph+c \, ; \, \ph\in H^{-1}(\Sp,g), c\in \R\}$ equipped with the topology induced by the natural product topology on $H^{-1}(\Sp,g) \times \R$. (It is easy to see that, unlike $H^{-1} (\Sp,g)$, this space is independent of the choice of $g=e^\rho g_0$ as in \eqref{E:metric}, hence we drop it from the notation). \indN{{\bf Function spaces}! $H^{-1}(\Sp)$; distributions of the form $\{\ph+c \, ; \, \ph\in H^{-1}(\Sp,g_0), c\in \R\}$} Then we define a measure $\lambda$ on $H^{-1} (\Sp)$, by setting
$$\lambda(A)=\int_{c\in \R} \mathbb{P}(h^{\Sp,g}+c\in A) \dd c $$
for an arbitrary Borel set $A \subset H^{-1} (\Sp)$.\footnote{In what follows we use the generic notation   $\mathbb{P}$  (and associated expectation $\mathbb{E}$) for the law of a field, for example $h^{\Sp,g}$ or $h^{\hat\C,\hat g}$, when the particular law in question is implicit from the notation.}
%More informally,
%$$
%  \lambda(\dd\ph) = \int_{c \in \R} \P( h^{\Sp,g} + c = \ph) \, \dd c
%$$
Note that for a non-negative Borel functional $F$ on $H^{-1} (\Sp)$ (an ``observable'') we have
\begin{equation}\label{E:lambdaintegral}
  \int_{\ph \in H^{-1} (\Sp)} F(\ph) \lambda (\text{d} \ph)= \int_{c \in \R} \E \Big[ F( h^{\Sp,g} + c) \Big] \text{d}c.
\end{equation}

This measure $\lambda$ allows us to assign a meaning to the term $\exp ( -\frac{1}{4\pi} \int |\nabla^g \ph|^2 v_g(\dd z) ) \text{D} \ph $ in the Polyakov measure \eqref{PolyakovMeasure}. Namely, we set
$$
\exp\Big( - \frac{1}{4\pi} \int |\nabla^g \ph|^2 v_g(\dd z) \Big) \text{D} \ph: = \lambda (\text{d} \ph)
$$
The fact that the total mass of the measure $\lambda$ is infinite is consistent with the fact that we expect the left hand side to be invariant under shifting $\ph$ by an arbitrary constant $c$.

%\red{Until the end of this Subsection I tried to use $\hat g$ on $\hat \C$ and $g$ on $\Sp$ generic.}
\medskip The measure $\lambda$ can also be pushed forward by $\psi$ to a measure $\hat \lambda$ on $H^{-1}(\hat \C)=\{\varphi+c; \varphi\in H^{-1}(\hat{\C},\hat g_0), c\in \R\}$ (which is a subspace of $H^{-1}_{\mathrm{loc}}(\C)$ as defined in \cref{C:H1loc}).\indN{{\bf Function spaces}! $H^{-1}(\hat{\mathbb{C}})$; distributions of the form $\{\ph+c \, ; \, \ph\in H^{-1}(\hat{\mathbb{C}},\hat{g}_0), c\in \R\}$} Similarly to $\lambda$, $\hat \lambda$ is then the ``law'' of $h^{\hat \C,  \hat g} +c$ with $c$ distributed according to Lebesgue measure, and  $h^{\hat \C, \hat g} = h^{\Sp,g} \circ \psi^{-1}$ as in \eqref{E:hChS}.

\medskip \textbf{Definition of the  Polyakov measure(s).} These two ingredients, the mean zero spherical GFF and its Lebesgue shift, allow us to give a rigorous definition of the Polyakov measure $\mathbf{P}$ of \eqref{PolyakovMeasure}. As before, we will allow ourselves to consider two closely related versions $\mathbf{P}$ and $\widehat{\mathbf{P}}$, depending on whether we want to consider the fields on $\Sp$ or on $\hat \C$.

%In fact we will give two closely related definitions. The first, which we call the spherical Polyakov measure, is more intrinsic and closer to the original Polyakov formulation, but a bit harder to compute with. In the second, we view the field as being defined on $\hat\C$ rather than $\Sp$, and in this case must add a shift to the field (related to the metric $g$), whose meaning will be discussed further in Section \ref{SS:LCFTshift}. 
%It is this second definition which we will actually work with, and we simply call it Polyakov measure. 
%This leads to nicer formulae, and although initially it looks further away from the original Polyakov formulation, there are intrinsic reasons to add in this metric shift. We will attempt to present one such justification. 

%Let us start with the definition of the \textbf{spherical Polyakov measure} $\mathbf{P}'$.  
Simply put, the (spherical) Polyakov measure $\mathbf{P}$ corresponds to reweighting $\lambda(\dd\ph)$ by the remainder of the terms in $S(\ph)$, 
%that is $$\exp\big(-\frac{1}{4\pi} \int_\Sp (Q R_g\ph +4\pi \mu e^{\gamma \ph}) v_g\big).$$
%Formally we will define $\mathbf{P}$ through its integral against functionals $F$ (which may be thought of as observables, hence the corresponding integrals are simply ``expectations'' of observables). 
%We will denote these either by $\mathbf{P} ( F)$ or $\mathbf{P}_g(F)$ to emphasise the dependence on the underlying metric $g$.
% or, in agreement with physics conventions, $\langle F \rangle_g$, where the subscript $g$ is used in order to emphasise the dependence on the metric $g$.
%
% which corresponds to a simply reweighting $\lambda(d\ph)$ by the remainder of the terms in $S(\ph)$, that is $$\exp\big(-\frac{1}{4\pi} \int_\C (Q R_g\ph +4\pi \mu e^{\gamma \ph}) v_g\big).$$ We will call this the \textbf{unshifted Polyakov measure}.
%We will then deduce how this transforms under a conformal change of metric, and observe that if we want a covariant transformation, we should actually modify the definition of $\mathbf{P}'$ by shifting the field slightly. This will be our final definition of the Polyakov measure $\mathbf{P}$ (see Definition \ref{D:Polyakovmeasure}).
%
%Let us start with $\mathbf{P}'$. 
%If $F$ is a non-negative Borel functional on $H^{-1} (\Sp)$ then we define the $\mathbf{P}$ expectation of $F$ (with respect to the metric $g$) by
 \begin{equation}\label{E:Polyakovdef'}
    \mathbf{P}(\dd \ph )=\mathbf{P}_g(\dd \ph): = %\int_{\ph \in H^{-1} (\Sp)% } F\left(\ph\right)  
    \exp \Big( - \frac{Q}{4\pi} (\ph,R_g)_g  -   \mu \cM_{\ph;g} (\Sp) \Big) \lambda (\dd \ph) .
 \end{equation} 
 Here we recall that $\cM_{\ph;g}( \Sp)$ is the total mass of the Gaussian multiplicative chaos measure with parameter $\gamma$ associated to $\ph = h^{\Sp,g} + c $, as defined in \eqref{eq:GMCsphere}. 
Thus for a non-negative functional $F$ on $H^{-1} ( \Sp)$, we have  
 \begin{align*}
 \mathbf{P}_g(F)  &=  \int_{\ph \in H^{-1} (\Sp) } F\left(\ph\right)  
    \exp \Big( - \frac{Q}{4\pi} (\ph,R_g)_g  -   \mu \cM_{\ph;g} (\Sp) \Big) \lambda (\dd \ph) \\
    & = \int_{-\infty}^{\infty} \E [F( h^{\Sp, g} + c)  \exp \Big( - \frac{Q}{4\pi} (h^{\Sp, g} + c,R_g)_g  -   \mu e^{\gamma c}\cM_{h^{\Sp, g};g} (\Sp) \Big) ]\dd c,
 \end{align*}
 where the expectation above is over the law of $h^{\Sp,g}$.
 
%Note that $\MM_\varphi(\Sp)$ is a measurable function of $\varphi$. Indeed, $h^{\hat\C,g}$ is a log-correlated field, so the theory of Chapter \ref{S:GMC} tells us that for any fixed $c$, 
%\begin{equation} \label{E:eboundMC} \MM_{h^{\Sp,g}+c}(\C)= e^{\gamma c} \lim_{\eps\to 0} \MM_{h^{\Sp,g} +(Q/2)\log g;\eps}(\C)  \end{equation}
%where the limit holds in probability and also, by \cref{C:spherecircleaverage}, in $L^1(\P)$. In particular $\mathbb{E}( \MM_{h^{\Sp,g}+c}(\C))=e^{\gamma c}\mathbb{E}( \MM_{h^{\Sp,g}}(\C))<\infty$.
%(In fact, one first needs to approximate $\C$ by a bounded an increasing sequence of bounded domains to apply the theory of that chapter, but this is straightforward).

Concretely, for computations it is more convenient to work on the extended complex plane. If we want to express \eqref{E:Polyakovdef'} in terms of an expectation over $h^{\hat \C,\hat g}$, a straightforward rewriting gives the following expression:
\begin{equation}\label{E:Pgshiftplane}
\int_{c\in \R} \E\left( F( h^{\hat \C, \hat g}\circ \psi^{-1} + c)  
\exp ( - \frac{Q}{4\pi} \int_\C R_{g} (h^{\hat \C, \hat g} + c) v_{\hat g} (\dd z)  - \mu e^{\gamma c} \cM_{h^{\hat \C, \hat g} + \tfrac{Q}2 \log \hat g}(\hat \C) ) 
\right) \dd c 
\end{equation}
where the equality 
\[\mu e^{\gamma c}\cM_{h^{\Sp, g};g} (\Sp)=\mu e^{\gamma c} \cM_{h^{\hat \C, \hat g} + \tfrac{Q}2 \log \hat g}(\hat \C)
\] is a consequence of \eqref{E:LCFTshift}.
The above is, however, not quite the right definition for the law $\widehat{\mathbf{P}}$ of the Polyakov measure in $\hat \C$, since we have only partly taken into account the change of coordinates from $\Sp$ to $\hat \C$. Instead, we define $\widehat{\mathbf{P}}_{\hat g}$ to be the ``law'' of $\ph \circ \psi^{-1} + (Q/2) \log \hat g$ under $\mathbf{P}_g$: 

\begin{definition}
	\label{D:Polyakovmeasure}
If $F$ is now a non-negative Borel functional on $H^{-1} ( \hat \C)$% \cong H^{-1} ( \C)$
, we set
%Now we come to the second definition, that we will actually take to be the Polyakov measure $\mathbf{P}$. The expression is very closely related to \eqref{E:Pgshiftplane}, but it will be a measure on fields on $\hat \C$ rather than $\Sp$ (that is, we will take our observable $F$ to be a function of the field on $\hat \C$ directly) and we will also shift the field by $(Q/2) \log g$. That is, we define
\begin{align}\label{metricshift}
& \widehat{\mathbf{P}}_{\hat g}(F)=
\mathbf{P}_g\left( F (( \cdot + \tfrac{Q}2 \log \hat g) \circ \psi) \right)\\
& = \int_{\R} \E\left[ F\left( h^{\hat \C,\hat g} + \tfrac{Q}2 \log \hat g + c \right)\exp\Big( -\tfrac{Q}{4\pi} (h^{\hat\C ,\hat g} + c ,R_{\hat g})_{\hat g}   - \mu e^{\gamma c} \cM_{h^{\hat \C ,\hat g} + \tfrac{Q}2 \log \hat g} (\C) \Big) \right] \dd c.\nonumber
\end{align}
	We will write $\langle F \rangle_{\hat g}$ for $\widehat{\mathbf{P}}_{\hat g}(F)$ in what follows, in agreement with physics conventions. 
\end{definition} 
\indN{{\bf Liouville CFT}! $\widehat{\mathbf{P}}_{\hat g}$; Polyakov measure}\indN{{\bf Liouville CFT}! $\langle \cdot \rangle_{\hat g}$; expectation with respect to the Polyakov measure}

%We refer to the shift of $\ph$  by $(Q/2) \log g$ in the observable $F$ as the \textbf{metric shift} (we will comment on this below). 
Note that $(c, R_{\hat g})_{\hat g} = 8\pi c $ by Gauss--Bonnet, hence 
 $$\exp( -\frac{Q}{4\pi} (c,R_{\hat g})_{\hat g} ) = \exp ( - 2Qc).$$ 
 Combining with \eqref{metricshift} we reach the following explicit definition of the Polyakov measure $\widehat{\mathbf{P}}_{\hat g}$: 
%	If $F$ is a non-negative Borel functional on $H^{-1} (\C)$ then we define the \textbf{Polyakov expectation} of $F$ (with respect to the metric $g$) by 
	%$$\mathbf{P}(F)=\mathbf{P}_g(F)=\langle F\rangle_g:=\mathbf{P}'_g(F_g)$$ so that
	\begin{equation*}
	\widehat{\mathbf{P}}_{\hat g}(F)\!\! =\! \! \int_{\R}\! \E\!\left[ F\!\left( h^{\hat \C,\hat g} \!+\! \tfrac{Q}2 \log \hat g \!+\! c \right)\exp\!\Big( \!-\tfrac{Q}{4\pi} (h^{\hat\C ,\hat g},R_{\hat g})_{\hat g} -2Qc  - \mu e^{\gamma c} \cM_{h^{\hat \C ,\hat g} + \tfrac{Q}2 \log \hat g} (\C) \Big) \right] \dd c. 
	\end{equation*}
%Note that the curvature term can be simplified: by Gauss--Bonnet, $(c, R_g)_g = 8\pi $, so that $$\exp( -\frac{Q}{4\pi} (c,R_g)_g ) = \exp ( - 2Qc).$$ Furthermore, 
Observe that when $g$ (or equivalently $\hat g$) has constant curvature, the term
$\exp( -\frac{Q}{4\pi} (h^{\hat \C,\hat g},R_{\hat g})_{\hat g})$ also disappears, since $h^{\hat \C,\hat g}$ has zero average with respect to $v_{\hat g}$ by definition. Thus in this case we get a particularly simple expression, which we will use repeatedly below:
\begin{equation}\label{E:Fgnice}
	\langle F\rangle_{\hat g} =   \int_{c \in \R} \E\left[ F\left( h^{\hat \C,\hat g} + \tfrac{Q}2 \log \hat g + c \right)\exp\Big( -2Qc   - \mu e^{\gamma c} \cM_{h^{\hat \C,\hat g} + \tfrac{Q}2 \log \hat g} (\C) \Big) \right] \dd c.
\end{equation}

\subsection{Weyl anomaly formula}
\label{SS:Weyl}
We have seen that the expression for $\langle F\rangle_{\hat g}$ simplifies considerably when $\hat g$ (or $g$ on $\Sp$) has constant curvature. When $R_{\hat g}$ is not constant, we have the additional term $\exp(-\frac{Q}{4\pi}(h^{\Sp,\hat{g}}+c,R_{\hat g})_{\hat g})$ in the expectation. 
The effect of this extra term is to further tilt the law of the field. 
%In fact, the form of this tilt is much simpler if we write $h^{\Sp,g}=h^{\Sp,\hat{g}}-v_g(h^{\Sp,\hat{g}})$ immediately in \eqref{E:Polyakovdef_int}.
It turns out the effect of this tilt can be described exactly, thanks to Girsanov's lemma.

\begin{theorem}[Weyl anomaly]\label{T:Weyl}
	Let $g$ be a metric on $\Sp$, with pushforward $\hat g$ to $\hat\C$ by $\psi$ of the form $\hat g(z) = e^{\rho(z)} \hat g_0(z)$ for $\rho$ as in \eqref{E:metric}.
	Then for each non-negative Borel function $F$ on $H^{-1}(\hat\C)$, we have
	\begin{align}\label{E:Weyl}
		\langle F \rangle_{\hat g} &  =  \exp\left(\frac{6Q^2}{96\pi} \int_\C [ |\nabla^{\hat g_0}\rho(x)|^2+ 4\rho (x) ]v_{\hat g_0}(\dd x)  \right)\, \langle F \rangle_{\hat g_0}.
	\end{align}
\end{theorem}
\ind{Weyl anomaly Formula}

See \cref{C:Weyl_correlation} and \cref{R:Weyl_correlation} for a Weyl anomaly formula valid for correlation functions. 

In the above expression, $\nabla^{\hat g_0}$ is the gradient operator in the metric ${\hat g_0}$. The only thing that is needed about this operator is the fact that $$\int_\C |\nabla^{\hat g_0}\rho(x)|^2 v_{\hat g_0}(\dd x)=-\int_\C \rho(x) \tfrac{1}{\hat g_0(x)} \Delta \rho(x) v_{\hat g_0}(\dd x)$$ (the Gauss--Green identity on $(\hat \C,\hat g_0)$).

\begin{comment}\begin{rmk}
	Notice that the expression \eqref{E:Weyl} differs from the formula \eqref{E:Weyl_anomaly} which is stated for the correlations of primary fields, $\langle \psi_{\alpha_1}(z_1)\dots \psi_{\alpha_k}(z_k)\rangle$, rather than for $\langle F \rangle$ with $F$ a Borel function on $H^{-1}(\C)$. In the case of primary fields there is an additional term $\exp(-\sum_{i=1}^k \Delta_{\alpha_i} \rho(z_i))$ on the right hand side of the anomaly formula. This is because, as we will see in the next section, $\langle \psi_{\alpha_1}(z_1)\dots \psi_{\alpha_k}(z_k)\rangle$ formally corresponds to $\langle F \rangle$ with $F(\varphi)=\exp(\alpha_1\varphi(z_1)+\dots+\alpha_k\varphi(z_k))$, but this $F$ is not a Borel function on $H^{-1}(\C)$ and needs to defined via regularisation. This regularisation procedure is the same as that used when defining Gaussian multiplicative chaos, and the additional term $\exp(-\sum_{i=1}^k \Delta_{\alpha_i} \rho(z_i))$ is analogous to the covariant term appearing in \cref{T:ci}. 
\end{rmk}
\end{comment}

\begin{rmk}
In \cite{DKRV} their definition for $\langle F \rangle_{\hat g}$ includes an extra multiplicative factor $\exp ( \frac1{96 \pi} \int_{\C}  [ \nabla^{\hat g_0}\rho(x)|^2+ 4\rho (x) ]v_{\hat g_0}(\dd x))$ which corresponds to the partition function of the Gaussian free field in \eqref{E:lambdaintegral}. This leads to the anomaly
	\begin{align}\label{E:Weylstandard}
		\langle F \rangle_{\hat g}  = \exp \left( \frac{c_L}{96\pi} \int_{\C}  [ \nabla^{\hat g_0}\rho(x)|^2 \rho(x)+ 2R_{\hat g_0}\rho (x) ]v_{\hat g_0}(\dd x)\right) \langle F\rangle_{\hat g_0}
	\end{align}
	with $$c_L=1+6Q^2.$$ 
	This is the more classical Weyl anomaly formula for a conformal field theory with \textbf{Central charge} $c_L$. Note that $c_L \in (25, \infty)$.
\end{rmk}
\ind{Central charge}

%This also allows us to express $\mathbf{P}_g'$ in terms of $\mathbf{P}'_{\hat g}$, and leads to the following lemma:

\begin{comment}
\begin{lemma}\label{L:Weylbad}
	Let $g$ be a metric on $\Sigma$, not necessarily of constant curvature, with $g = e^\rho \hat g$ as in \eqref{E:metric}.
	%The Polyakov measure associated to the $g$ is then chosen to assign
	Then for each non-negative Borel function $F$ on $H^{-1}(\C)$, we have
	\begin{align}\label{E:Weylbad}
		\mathbf{P}'_g(F) &  =  \exp\left(Q^2\bar{v}_{\hat{g}}(\rho)-\frac{Q^2}{16\pi}\int_\C \rho(x) \Delta^{\Sp,\hat g}\rho(x) v_{\hat{g}}(\dd x) \right)\, \mathbf{P}'_{\hat g}(F_\rho)
	\end{align}
where $F_\rho(\ph)=F(\ph-\frac{Q}{2}\rho)$.
\end{lemma}
\end{comment}

\begin{proof}
	Recalling Lemma \ref{L:hsnorms} (applied to $h^{\hat \C,\hat{g}}$), and applying the change of variables $c\mapsto c-\bar{v}_{\hat g}(h^{\hat\C,\hat{g}_0})$ in the definition of the Polyakov measure, we first rewrite $\mathbf{P}_g(F) = \langle F \rangle_g$ as 
	\begin{equation*}  \int \!\E\! \left[ \exp(-\frac{Q}{4\pi} (h^{\hat \C,\hat{g}_0},R_{\hat g})_{\hat g})  F\left( h^{\hat\C,\hat{g}_0} +   \tfrac{Q}{2}\log \hat g + c \right)\!  \exp\! \Big(\! \! - \! \!2Q c \! - \! \! \mu e^{\gamma c}\cM_{h^{\hat\C,\hat g_0} + (Q/2)\log \hat g} (\C) \Big)\right]\! \dd c. 
	\end{equation*}
	By Girsanov, the effect of the term $\exp(-\frac{Q}{4\pi}  (h^{\hat\C,\hat{g}_0},R_{\hat g})_{\hat g})$ is to shift the field $h^{\hat\C,\hat{g}_0}$ by
	\begin{equation}\label{E:tilt_R} -\frac{Q}{4\pi}\int_\C 2\pi G^{\hat\C,\hat{g}_0}(x,y)R_{\hat g}(y) {\hat g}(y)\, \dd y\end{equation}
and to multiply the whole expression by
	\begin{equation}
		\label{E:tilt_R_norm} \mathbb{E}\left[\exp(-\frac{Q}{4\pi}  (h^{\hat\C,\hat{g}_0},R_{\hat g})_{\hat g})\right].
	\end{equation}
	In fact, both of these expressions simplify quite nicely. Recalling that
	$$R_{\hat g} \hat g = - \Delta \log \hat g = -\Delta \rho - \Delta \log \hat{g}_0 = -\Delta \rho + R_{\hat{g}_0}\hat{g}_0,$$
	we have
	\begin{align} -\frac{Q}{2}\int_\C G^{\hat\C,\hat{g}_0}(x,y)R_{\hat g}(y) {\hat g}(y) \dd y & = -\frac{Q}{2}\int_\C G^{\hat\C,\hat{g}_0}(x,y)\big(R_{\hat{g}_0}(y)\hat{g}_0(y)-\Delta \rho(y)\big) \dd y\nonumber \\
		& = \frac{Q}{2}\int_\C G^{\hat\C,\hat{g}_0}(x,y)\Delta^{\hat\C,\hat{g}_0} \rho(y) v_{\hat{g}_0}(\dd y) \nonumber \\
		& =-\frac{Q}{2}(\rho(x)-\bar{v}_{\hat{g}_0}(\rho)) \end{align}
	where the second line follows because $R_{\hat{g}_0}=2$ and $v_{\hat{g}_0}(G^{\hat\C,\hat{g}_0}(x,\cdot)) = 0$, while the third line follows from Remark \ref{R:gsg_lap}. For this we used the assumption that $\rho$ is twice continuously differentiable.
	Similarly,
	\begin{align}
		\label{E:varianceweyl}
			\mathbb{E}\left[\exp(-\frac{Q}{4\pi}  (h^{\hat\C,\hat{g}_0},R_{\hat g})_{\hat g})\right] & = \exp(\frac{Q^2}{16\pi} \int_\C R_{\hat g}(x) \hat{g}(x) (\rho(x)-{\bar{v}}_{\hat{g}_0}(\rho)) \, \dd x) \nonumber \\
		& = \exp\big(\frac{Q^2}{16\pi}\int_\C (R_{\hat{g}_0}\hat{g}_0(x)-\Delta \rho(x))(\rho(x)-\bar{v}_{\hat{g}_0}(\rho)) \dd x\big) \nonumber \\
		& = \exp\big(-\frac{Q^2}{16\pi}\int_{\C}\rho(x) \Delta^{\hat\C,\hat g_0}\rho(x) v_{\hat{g}_0}(\dd x) \big).  \end{align}
	where the last line follows since $$\int_\C R_{\hat g_0} \hat{g}_0(x) \rho(x) \dd x = 2 {v}_{\hat g_0}(\rho)=8\pi \bar{v}_{\hat g_0}(\rho) = \int_\C R_{\hat g_0} \hat{g}_0(x) \bar{v}_{\hat g_0}(\rho) \dd x$$ and $\int_\C \Delta \rho(x) \, \dd x = \int_\C \Delta^{\hat\C,\hat{g}_0}\rho(x) v_{\hat g_0}(\dd x)=0$ by \eqref{E:GG_manifold}.
	
	Notice that subtracting $-\frac{Q}{2}\rho$ from the field in our expression for $\mathbf{P}_g(F)$ has exactly the effect of turning $h^{\hat\C,\hat{g}_0}+\tfrac{Q}{2}\log \hat g$ into $h^{\hat\C,\hat{g}_0}-\tfrac{Q}{2}\log \hat{g}_0$. Combined with a further change of variables $c\mapsto c+\frac{Q}{2}\bar{v}_{\hat{g}_0}(\rho)$ in the integral we reach the conclusion:
	\begin{equation} \langle F \rangle_{\hat g} = \langle F \rangle_{\hat g_0} \exp\big(Q^2\bar{v}_{\hat{g}_0}(\rho)-\frac{Q^2}{16\pi}\int_\C \rho(x) \Delta^{\hat\C, \hat g_0}\rho(x) v_{\hat{g}_0}(\dd x) \big), \end{equation}
where the anomaly term can be rewritten as in the statement of the Theorem.
%
	%
	\begin{comment}To complete the proof of the theorem it remains to check that
	\begin{equation} \label{E:cgcg} c_g = c_{\hat g} \exp\big( -2Q^2 \bar{v}_{\hat g}(\rho) + \frac{Q^2}{8\pi} \int_\C \rho(x) \Delta^{\Sp, \hat{g}}\rho(x) v_{\hat g}(\dd x)\big). \end{equation}
	For this we use the same trick as above, writing
	\begin{align*}
		c_g & = \exp \left( - \frac{Q^2}{8\pi} \int_\C R_g(x)  \log g(x) v_g(\dd x) \right) \\
		& = \exp \left( - \frac{Q^2}{8\pi} \int_\C (R_{\hat g}(x)\hat g(x)-\Delta \rho(x)) ( \log {\hat g}(x) + \rho(x)) \dd x \right) \\
		& = c_{\hat g}  \exp\big(  \frac{Q^2}{8\pi} \int_\C \rho(x) \Delta^{\Sp, \hat{g}}\rho(x) v_{\hat g}(\dd x)-Q^2 \bar{v}_{\hat g}(\rho)\big) \exp\big( \frac{Q^2}{8\pi}\int_\C \Delta^{\Sp,\hat{g}} \rho(x) \log \hat{g}(x) v_{\hat g}(\dd x) \big)
	\end{align*}
	where $\int_\C \Delta^{\Sp,\hat{g}} \rho(x) \log \hat{g}(x) v_{\hat g}(\dd x) = \int_\C \rho(x) \Delta^{\Sp,{\hat g}}\log \hat{g}(x) v_{\hat g}(\dd x) = - \int_\C R_{\hat{g}} \rho(x) v_{\hat g}(dx)=8\pi \bar{v}_{\hat g}(\rho)$ by \eqref{E:GG_manifold}. Substituting this into the above gives \eqref{E:cgcg}.
	\end{comment}
	%
\end{proof}

\begin{lemma}[Weyl Anomaly for M\"obius transforms]
	\label{L:weylmobius}
When $\hat{g}=m_*\hat{g}_0$, with $m$ a M\"{o}bius transform of $\hat{\C}$, we have 
\[ \langle F \rangle_{m_*\hat g_0} = \langle F \rangle_{\hat g_0}\]
for all non-negative Borel functions $F$.
\end{lemma}
\begin{proof}
	Recall from \cref{L:thetaRcmg} that $R_{m^*\hat{g}_0}\equiv R_{\hat g_0}\equiv 2$, and by definition of $\rho$ with $\hat g=m_*\hat g_0$, 
	$$
	-\Delta \rho(x)=R_{m_*\hat{g}_0}m_*\hat g_0(x)-R_{\hat g_0}\hat g_0(x)=2m_*\hat g_0(x)-2\hat g_0(x).
	$$
	We therefore have 
	\begin{equation}\label{E:anomalymobius1}
- \frac{Q^2}{16\pi}\int_\C \rho(x)\Delta \rho(x) \dd x  = \frac{Q^2}{2}\left(\bar{v}_{m_*\hat g_0}(\rho)-\bar{v}_{\hat g_0}(\rho)\right)
	\end{equation}
while also 
	\begin{equation}\label{E:anomalymobius2}
	- \frac{Q^2}{16\pi}\int_\C \rho(x)\Delta \rho(x) \dd x  = 2Q^2 \int_{\hat \C}\int_{\hat \C} (2\pi G^{\hat \C,\hat g_0})(x,y)\bar{v}_{m_*\hat g_0}(\dd x)\bar{v}_{m_*\hat g_0}(\dd y)
\end{equation}
by \eqref{E:varianceweyl}.
Now, on the one hand, by \eqref{E:gSequiv}, we have
\[
2\pi G^{\hat \C, m_*\hat g_0}(x,y)= -\log|x-y|+\bar{v}_{m_*\hat g_0}(\log|x-\cdot|)+\bar{v}_{m_*\hat g_0}(\log|x-\cdot|)-\theta_{m_*\hat g_0}.
\]
On the other, since $2\pi G^{\hat \C, m_*\hat g_0}$ and $2\pi G^{\hat \C, \hat g_0}$ are the variances of $h^{\hat \C,m_*\hat g_0}$ and $h^{\hat \C,\hat g_0}$ respectively, and we know by \cref{L:hsnorms} that  $h^{\hat \C,m_*\hat g_0}$ is equal in distribution to $h^{\hat \C, \hat{g}_0}-\bar{v}_{m_*\hat g_0}(h^{\hat\C,\hat g_0})$, we have 
\begin{multline*}
2\pi G^{\hat\C,m_*\hat g_0}(x,y)= 2\pi G^{\hat\C,\hat g_0}(x,y)-\int_\C 2\pi G^{\hat\C,\hat g_0}(x,y) \bar{v}_{m_*\hat g_0}(\dd x)-\int_\C 2\pi G^{\hat\C,\hat g_0}(x,y) \bar{v}_{m_*\hat g_0}(\dd x) \\ +\int_{\hat \C}\int_{\hat \C} (2\pi G^{\hat \C,\hat g_0})(x,y)\bar{v}_{m_*\hat g_0}(\dd x)\bar{v}_{m_*\hat g_0}(\dd y).
\end{multline*}
Using that $2\pi G^{\hat \C, \hat g_0}=-\log|x-y|-(1/4)(\log \hat g_0(x)+\log \hat g_0(y))+\log(2)-\theta_{\hat g_0}$ and equating the two expressions for $2\pi G^{\hat \C, m_* \hat g_0}$ above, we are left with the equality
\begin{align*}
\int_{\hat \C}\int_{\hat \C} (2\pi G^{\hat \C,\hat g_0})(x,y)\bar{v}_{m_*\hat g_0}(\dd x)\bar{v}_{m_*\hat g_0}(\dd y) & =  -\theta_{m_*\hat g_0}+\log(2)-\theta_{\hat g_0}-\frac{1}{2}\bar{v}_{m_*\hat g_0}(\log \hat{g}_0)\\
& = \frac{1}{2}\bar{v}_{m_*\hat{g}_0}(\log(m_*(\hat g_0)))-\frac{1}{2}\bar{v}_{m_*\hat g_0}(\log \hat{g}_0) \\
& = \frac{1}{2}\bar{v}_{m^*\hat g_0}(\rho),
\end{align*}
where the second equality follows from the expression \eqref{eq:thetamg} relating $\theta_{m_*\hat g_0}$ and $\theta_{\hat g_0}$.
Combining this with \eqref{E:anomalymobius1} and \eqref{E:anomalymobius2} we deduce that 
$$ 
\frac{Q^2}{2}\left(\bar{v}_{m_*\hat g_0}(\rho)-\bar{v}_{\hat g_0}(\rho)\right)=Q^2\bar{v}_{m_*\hat g_0}(\rho)
$$
and so $\bar{v}_{m_*\hat g_0}(\rho)=-\bar{v}_{\hat g_0}(\rho)$. We conclude that the anomaly term 
$$
\exp\big(Q^2\bar{v}_{\hat{g}_0}(\rho)-\frac{Q^2}{16\pi}\int_\C \rho(x) \Delta^{\hat\C, \hat g_0}\rho(x) v_{\hat{g}_0}(\dd x) \big)=\exp\big(Q^2\bar{v}_{\hat g_0}(\rho)-Q^2\bar{v}_{\hat g_0}(\rho)\big)=1,
$$
which completes the proof.

\end{proof}

\subsection{Convergence of correlation functions within Seiberg bounds}

\textbf{In this section we drop the superscripts $\hat \C,\hat g$ for ease of notation.} In particular $h=h^{\hat \C,\hat g}$. \medskip

%Definition \ref{D:Polyakovmeasure} has the advantage that it stays explicitly close to Polyakov's original formulation. However, 
It is not immediately obvious for which observables $F$ can we say that the associated expectation $\langle F \rangle_{\hat g}$ is finite. Let us consider the simplest case where $F = 1$ and $\hat g=\hat{g}_0$. Then recall that by \cref{D:Polyakovmeasure} of the Polyakov measure we have
\begin{equation}\label{divPartitionfunction}
\langle 1 \rangle_{\hat g_0}  = \int_{c \in \R} \E \Big[ \exp \Big( - 2Qc - \mu e^{\gamma c} \cM_{h+\tfrac{Q}{2}\log \hat{g}} (\C)\Big) \Big] \text{d} c.
\end{equation}
The two possible divergences we need to worry about are at $c \to \infty$ and $c \to - \infty$. The first one (when $c \to \infty$) is not a problem since $Q>0$ and the GMC mass is also strictly positive, so that overall the integrand decays (doubly) exponentially as $c \to \infty$, hence the integral converges for large $c$. The second limit however is divergent: indeed, when $c \to - \infty$, the exponential term is  $\exp ( - 2Qc + o(1))$ which blows up exponentially. This implies $\langle 1 \rangle_{\hat g} = \infty$.

\medskip It turns out that we get a convergent expectation if we choose for our observable the natural \textbf{correlation functions} of the model (denoted by $V = V_{\alpha_1, \ldots, \alpha_k} (\mathbf{z})$), that is, informally,
\begin{equation}\label{eq:Vepsilon1}
V(\ph) = e^{\alpha_1 \ph(z_1) + \ldots + \alpha_k \ph(z_k)}
\end{equation}
\indN{{\bf Liouville CFT}! $V_{\alpha_1,\dots, \alpha_k}(\mathbf{z})$; vertex operator}
where $\alpha_1, \ldots, \alpha_k$ are real numbers and $\mathbf{z} = (z_1, \ldots, z_k) \in \C^k$ with $z_i\ne z_j$ for $i\ne j$.  In physics language, $e^{\alpha_i \ph (z_i)}$ is a \textbf{vertex operator} and $z_i$ is called an \textbf{insertion}.

\ind{Correlation functions}
\ind{Vertex operator}
\ind{Insertion}

As usual, since $\ph$ is a distribution, it is not entirely clear what one means by $e^{\alpha_i \ph (z_i)}$, so in order to even speak about 
$\langle F \rangle_{\hat g}$ some regularisation is necessary. Let $h_\eps$ denote the circle average of $h$.
Define 
\begin{equation} \label{eq:Vepsilon2}V_\eps(\ph)=\prod_{i=1}^k \eps^{\alpha_i^2/2}e^{\alpha_i \ph_\eps(z_i)},  
\end{equation} 
so that
\begin{equation*}
\langle V_\eps \rangle_{\hat g}   =  \int_{c \in \R} \E \left[ \prod_{i=1}^k\eps^{\alpha_i^2/2}  e^{\alpha_i \left(h_\eps(z_i) +  \frac{Q}{2}\log {\hat g}(z_i)  + c\right) }\exp \Big( - 2Q c  - \mu e^{\gamma c} \cM_{h+\tfrac{Q}{2}\log {\hat g}} (\C) \Big) \right] \text{d} c.
\end{equation*}
 We will attempt to define $\langle V \rangle_{\hat g}$ by taking a limit of $\langle V_\eps\rangle_{\hat g}$ as $\eps\to 0$.

%and discussed in Remark \ref{R:LCFTshift}) regularised at scale $\eps$, that is,
%\begin{align}\label{E:shiftreg}
%\MM_{h ; \eps}(\C)&  = \int_{\C} \eps^{\gamma^2/2}\exp \left( \gamma \left[h_\eps(z) + \tfrac{Q}2 \log g(z)\right] \right) \dd z.
%%& =\int_{\C} \exp \left[ \gamma h_\eps(z)  - \frac{\gamma^2}2 \E( h_\eps(z)^2)  \right] v_g(\dd z)
%\end{align}
%We also point out that when we specialise to $g = \hat g$, we have the identity:
%\begin{equation}
%\label{E:shiftreg2}
%\cM^{\hat g}_{h ; \eps}(\C)  =\int_{\C} \exp \left[ \gamma h_\eps(z)  - \frac{\gamma^2}2 \E( h_\eps(z)^2)  \right] v_g(\dd z).
%\end{equation}
%In particular, as $\eps \to 0$, this converges to a random variable in $(0, \infty)$ almost surely. This last property (but not the exact identity \eqref{E:shiftreg2}) remains true for a general metric $g$ subject to \eqref{E:metric}.

Note that if $\sum_{i=1}^k \alpha_i > 2Q$, the problem near $c =- \infty$ leading to the divergence of \eqref{divPartitionfunction} should disappear. On the other hand, no new problem is created at $c= \infty$ because the decay of the integrand is doubly exponential in that region. This suggests that there is a chance that $\langle V \rangle_{\hat g}=\lim_{\eps\to 0} \langle V_\eps \rangle_{\hat g}$ may be finite if $\sum_{i=1}^k \alpha_i > 2Q$. On the other hand if any of the $\alpha_i$ is too large then the expectation can also explode as we are adding a logarithmic singularity of strength $\alpha_i$ to the field. One might naively guess that the maximal allowed value for $\alpha_i$ could be $\alpha_i = \gamma$ (corresponding to a Liouville typical point) or perhaps $\alpha_i=2$ (corresponding to the maximal thickness of any point $h$). Surprisingly, the maximal allowed value is in fact strictly larger, namely it suffices to require $\alpha_i <Q$. That the expectation is convergent and non-zero under these two conditions, collectively known as the Seiberg bounds, is the content of the next theorem and one of the main results of \cite{DKRV}.
With these notations, the main theorem of this section is the following:

\begin{theorem}
  \label{T:Seibergbounds}
  Suppose $\alpha_1, \ldots, \alpha_k \in \R$ satisfy
  \begin{equation}\label{E:Seiberg1}
  \sum_{i=1}^k \alpha_i > 2Q
  \end{equation}
and $z_1,\dots, z_k\in \C$ are distinct.
  Then $\langle V_\eps \rangle_{\hat g}< \infty.$ Suppose furthermore
  \begin{equation}\label{E:Seiberg2}
    \alpha_i < Q \quad \text{ for } i = 1, \ldots, k.
  \end{equation}
  Then $\langle V_\eps \rangle_{\hat g}$ converges to a limit $\langle V \rangle = \langle V \rangle_{\hat g}$ in $(0, \infty)$ as $\eps \to 0$.
\end{theorem}\indN{{\bf Liouville CFT}! $\langle V \rangle_{\hat g}$; correlation function}

Before giving the proof of this theorem, we make a few comments. The two bounds \eqref{E:Seiberg1} and \eqref{E:Seiberg2} are known collectively as the Seiberg bounds. They are known to be optimal in the sense that if either of these bounds fail, then either $\langle V_\eps \rangle_{\hat g} = \infty$ or it converges to zero as $\eps \to 0$. Note that for these two bounds to be simultaneously satisfied, it is necessary that $k \ge 3$. The necessity of these \textbf{three insertions} will be discussed below both from a geometric and probabilistic perspective. \ind{Seiberg bounds} \ind{Insertion}

\medskip We now start the proof of this theorem.

\begin{proof} Without loss of generality, by the Weyl anomaly (Theorem \ref{T:Weyl}) we take $\hat{g}=\hat{g}_0$.
Our first task will be to re-express $\langle V_\eps \rangle_{\hat g_0}$ via Girsanov's theorem (Lemma \ref{L:Girsanov}). To do this we will need to view each term $e^{\alpha h_\eps(z_i)} \eps^{\alpha_i^2/2}$ as an exponential biasing of the Gaussian field $h$. Notice however that the normalising factor, namely $\eps^{\alpha_i^2/2}$, is not quite equal to $\E[e^{\alpha_i h_\eps( z_i)}]^{-1}$, so we must account for this using Lemma \ref{L:covariance_hsphere}.

The result is that we may write
\begin{align}
& \langle V_\eps \rangle_{\hat{g}_0} = \int_{c \in \R} \E \left[ \prod_{i=1}^k\eps^{\alpha_i^2/2}  e^{\alpha_i (h_\eps(z_i) +  \tfrac{Q}{2} (\log \hat{g}_0)_\eps (z_i)+ c) }\exp \Big( - 2Q c  - \mu e^{\gamma c}\cM_{h+\tfrac{Q}{2}\log \hat{g}_0}(\C) \Big) \right] \text{d} c \nonumber \\
& = e^{C_\eps(z_1, \ldots, z_k) }   \prod_i \hat{g}_0 (z_i)^{\tfrac{\alpha_i}{2}(Q-\tfrac{\alpha_i}{2})}  \int_{c \in \R}\E \left[  \exp \Big( (\sum_{i=1}^k \alpha_i - 2Q) c  - \mu e^{\gamma c} \cM_{\hat h^\eps +\tfrac{Q}{2}\log \hat{g}_0} (\C)\Big) \right] \text{d} c \label{LCFTgirsanov}
\end{align}
where the field $\hat h^\eps$ is obtained from $h$ by applying the Girsanov shift,
\begin{equation}\label{eq:shiftGirsanov_seiberg}
\hat h^\eps (\cdot) = h (\cdot) + \sum_{i=1}^{k} \alpha_i \int_0^{2\pi} 2\pi G^{\hat\C,\hat g_0} (z_i + \eps e^{i \theta}, \cdot) \frac{\text{d} \theta}{2\pi}.
\end{equation}
(The factor $2\pi$ in front of the Green function is due to our normalisation: $h = \sqrt{2\pi} \mathbf{h}$). The normalising constant $C_\eps$ above satisfies
 $$
C_\eps(z_1, \ldots, z_k) : =\sum_{i=1}^k \sum_{j>i} 2\pi \alpha_i \alpha_j  G^{\hat\C,\hat g_0} (z_i, z_j) + \sum_{i=1}^k \frac{\alpha_i^2}{2} (\log (2)-\frac{1}{2}) +o(1).
$$
Since this normalising factor has a well behaved limit as $\eps \to 0$, which we call $C(\mathbf z)$, it suffices to consider the integral term in \eqref{LCFTgirsanov}. Writing $$Z_\eps := \cM_{\hat h^\eps+\tfrac{Q}{2}\log \hat{g}_0}(\C) ,$$ applying Fubini's theorem, and changing variables $u = e^{\gamma c} Z_\eps$, so that $\dd u =  \gamma e^{\gamma c } Z_\eps \dd c = \gamma u \dd c$ we obtain that
\begin{align}
\langle V_\eps \rangle_{\hat g_0}& \sim e^{C(\mathbf{z})} \prod_i \hat{g}_0 (z_i)^{\tfrac{\alpha_i}{2}(Q-\tfrac{\alpha_i}{2})} \, \E \left[  \int_{c \in \R} \exp \Big( (\sum_{i=1}^k \alpha_i - 2Q) c  - \mu e^{\gamma c}Z_\eps\Big)  \text{d} c \right] \nonumber \\
 & \sim e^{C(\mathbf{z})} \prod_i \hat{g}_0 (z_i)^{\tfrac{\alpha_i}{2}(Q-\tfrac{\alpha_i}{2})} \, \E \left[  \int_{u >0}
 \Big( \frac{u}{ Z_\eps} \Big)^{\frac{\sum_i \alpha_i - 2Q}{\gamma}}
 e^{- \mu u} \frac{\dd u }{\gamma u }
 \right] \nonumber \\
 & \sim \frac{e^{C(\mathbf{z})}}{\gamma } \prod_i \hat{g}_0 (z_i)^{\tfrac{\alpha_i}{2}(Q-\tfrac{\alpha_i}{2})} \int_{u >0}   {u}^{\frac{\sum_i \alpha_i - 2Q}{\gamma} -1 } e^{-\mu u} \dd u \cdot \E\left[Z_\eps^{ - \frac{\sum_i \alpha_i - 2Q}{\gamma}} \right].\label{LCFTgamma}
\end{align}
Above we use the Landau notation $a_\eps \sim b_\eps$ to mean that the ratio $a_\eps / b_\eps \to 1$ as $\eps\to 0$.
The integral over $u$ does not depend on $\eps$ (in fact, it is nothing but the Gamma function evaluated at $\frac{\sum_i \alpha_i - 2Q}{\gamma}  $). Hence, the proof of the theorem eventually boils down to proving the following lemma:
\begin{lemma}
\label{L:negmomentsLCFT}
Let $s = \frac{\sum_i \alpha_i - 2Q}{\gamma} > 0 $. Then the limit
$$
\lim_{\eps \to 0} \E( Z_\eps^{-s}) = :\E(Z_0^{-s})
$$
exists and lies in $(0,\infty)$.
\end{lemma}

\begin{proof}
Recall that
\begin{align*}
Z_\eps & = \lim_{\delta\to 0} \int_{\C} \delta^{\gamma^2/2}e^{\gamma[\hat h^\eps_\delta (z) + \tfrac{Q}2 (\log \hat g_0)_{\delta}(z)]} \dd z\\
& = \lim_{\delta\to 0} \int_\C e^{\gamma (H_\eps)_\delta(z)} \delta^{\gamma^2/2}e^{\gamma[h_\delta (z) + \tfrac{Q}2 (\log \hat g_0)_{\delta}(z)]} \dd z
%& = \int_{z\in \C} e^{\gamma H_\eps (z)} \MM_{h;\eps} (\dd z)
\end{align*}
where the subscript $\delta$ is used everywhere to denote the circle average at radius $\delta$, and the convergence is in probability. Here $$
H_\eps (z)  = \sum_{i=1}^{k} \alpha_i \int_0^{2\pi} 2\pi G^{\hat\C,\hat g_0} (z_i + \eps e^{i \theta}, z) \frac{\text{d} \theta}{2\pi}.
$$
Since $H_\eps$ is a smooth function of $z$ for fixed $\eps$, $(H_\eps)_\delta(z)\to H_\eps(z)$  uniformly in $z$ as $\delta\to 0$. Together with the fact that $H_\eps$ is uniformly bounded and that $\cM_{h+(Q/2)\log \hat g_0}(\C)$ is a limit in $L^1(\P)$ of $ \int_{\C} \delta^{\gamma^2/2}e^{\gamma[h_\delta (z) + (Q/2) (\log \hat{g}_0)_{\delta}(z)]} \dd z$, this implies that
$$Z_\eps= \lim_{\delta\to 0} \int_\C e^{\gamma (H_\eps)_\delta(z)} \delta^{\gamma^2/2}e^{\gamma[h_\delta (z) + \tfrac{Q}2 (\log \hat g_0)_{\delta}(z)]} \dd z=\int_\C e^{\gamma H_\eps (z)} \cM_{h+(Q/2)\log \hat{g}_0} (\dd z),$$
where the limit holds in $L^1(\mathbb{P})$ and in probability.
In particular, $Z_\eps$ has finite expectation for each $\eps>0$. \medskip
%is the circle-average of the Girsanov shift and, as in \eqref{E:shiftreg}, $\MM_{h;\eps} (\dd z) = \eps^{\gamma^2/2} e^{\gamma[ h_\eps(z) + (Q/2) \log g(z)] } dz$ is the $\eps$-approximation of Liouville measure associated with $h = h^{\Sp,g}$, shifted by $(Q/2) \log g$.
%Note first that $\MM_{h;\delta}$ converges weakly to the Liouville measure $\MM_h$ whose total expected mass is finite, as observed \ref{XXX}.
%\begin{align}
%\E[ \cM_h (\C)] &=\lim_{\eps\to 0} \eps^{\gamma^2/2}  \int_\C  \E [ e^{\gamma h_\eps (z)}] e^{\tfrac{\gamma Q}2 \log g(z)}dz \nonumber \\
%& = e^{( \log 2 - 1/2) \gamma^2/2 } \int_{z \in \C} g(z) dz < \infty. \label{totalmassLCFT}
%\end{align}

\noindent Before proceeding with the proof, let us make a couple of remarks.
\begin{itemize}
\item It is clear that $H_\eps (z)$ converges to a function $H (z) = \sum_{i=1}^k \alpha_i 2\pi G^{\hat\C,\hat g_0}(z_i, z)$ as $\eps\to 0$. It is therefore natural to expect that $Z_\eps$ converges (in probability, say) to
$$
Z_0: = \int_{\C} e^{\gamma H(z)} \cM_{h+\tfrac{Q}{2}\log \hat{g}_0} (\dd z) .
$$
This will indeed follow from our proof.
\item Most of the proof consists in checking that $Z_0$ is in fact finite almost surely under the second Seiberg bound \eqref{E:Seiberg2} (that is, $\alpha_i < Q$). This makes the negative moment $\E(Z_0^{-s})$ strictly positive.
This is however far from obvious: for instance, the expectation of $Z_0$ actually blows up if one of the $\alpha_i$ satisfies $\gamma \alpha_i >2$ (which is allowed since $Q>2/\gamma$). Nevertheless, $Z_0$ remains finite almost surely, even though its expectation is infinite. The fact that $Z_0$ remains finite under the Seiberg condition \eqref{E:Seiberg2} is instead a consequence of a scaling argument, as we will now see.
\end{itemize}

We divide the proof in several steps. The first step is the easy upper bound on $\E[(Z_\eps)^{-s}]$ (corresponding to a lower bound on $Z_\eps$).

\medskip \textbf{Step 1.} For any $q>0$ there exists $C=C_s>0$ such that $\E[ (Z_\eps)^{-q}] \le C$ for all $\eps>0$. This is indeed easy to see, since if we consider any bounded set $B$ (which does not even need to stay disjoint from the insertions $z_i, 1\le i \le k$), then
\begin{align*}
Z_\eps &\ge \int_B e^{\gamma H_\eps (z) } \cM_{h+(Q/2)\log \hat{g}_0}(\dd z)\\
& \ge C \cM_{h+(Q/2)\log \hat{g}_0}(B)
\end{align*}
where $C$ is a uniform (in $\eps$ and $z \in B$) lower bound on $e^{\gamma H_\eps(z)}$. Taking the negative moment of order $-q<0$, the claimed upper bound therefore follows from Theorem
\ref{T:negmom}. %\nb{Actually what is needed here is that the approximation are unif bounded, so we need to say something there. DKRV are not very precise here. }

\medskip \textbf{Step 2.} We fix $r>0$ and let $A_r = \cup_{i=1}^k B(z_i, r)$. We decompose $Z_\eps$ according to whether $z \in A_r$ or not. Thus we write
\begin{align*}
Z_\eps  &= \int_{A_r} e^{\gamma H_\eps (z)} \cM_{h+(Q/2)\log \hat{g}_0} (\dd z) + \int_{A_r^c} e^{\gamma H_\eps (z)} \cM_{h+(Q/2)\log \hat{g}_0} (\dd z)\\
& =: Z_{r, \eps} + Z_{r, \eps}^c.
\end{align*}
In this second step we show that for fixed $r$, the term $Z_{r,\eps}^c$ corresponding to points far away from the insertions is well behaved and has a limit in probability. Note that $H_\eps \to H$ uniformly on $A_r^c$, and that $\cM_{h+(Q/2)\log \hat{g}_0} (\dd z)$ is a measure of uniformly bounded total expectation (as already observed, it is a limit in $L^1(\P)$). Hence
$$
\int_{A_r^c} |e^{\gamma H_\eps(z)} - e^{\gamma H(z)} |\cM_{h+(Q/2)\log \hat{g}_0} (\dd z) \to 0
$$
in $L^1(\P)$ and in probability.  %Furthermore,
%$$\int_{A_r^c} e^{\gamma H(z)} \MM_{h;\eps} (\dd z) \to \int_{A_r^c} e^{\gamma H(z)} \MM_h (
%\dd z )$$
%in probability, by weak convergence in probability of $\MM_{h;\eps}$ to $\MM_h$ and the easily seen fact that $H$ is bounded: this itself follows from the fact that $G^{\Sp,g} (\cdot, \cdot)$ is bounded away from the diagonal, see Lemma \ref{L:sphereWholePlaneGFF}.
We deduce that $Z_{r, \eps}^c \to Z_{r, 0}^c: = \int_{A_r^c} e^{\gamma H(z)}\cM_{h+(Q/2)\log \hat{g}}(\dd z)$ in probability as $\eps\to 0$ for each fixed $r$.

\medskip \textbf{Step 3.}  In this third step we show that the term $Z_{r,\eps}$ corresponding to the points close to the insertions does not blow up if $\alpha_i <Q$ for all $1\le i \le k$. More precisely, we will show that there is a function $C(r)>0$ such that $C(r) \to 0$ as $r\to 0$, and such that for sufficiently small $p>0$ and all $\eps>0$,
\begin{equation}\label{E:goalStep3LCFT}
\E [ (Z_{r, \eps})^p] \le C(r).
\end{equation}
This is the most technical part of the proof.
To begin with, we observe that by subadditivity of $x \mapsto x^p$ for $p<1$, it suffices to prove the result for $k = 1$ insertions. Without loss of generality we take $z_1 = 0$, and we write $\alpha= \alpha_1$ for the power associated with the corresponding insertion. Recall that $\alpha<Q$ as per \eqref{E:Seiberg2}.

Since $Z_{r,\eps}^p$ is decreasing with $r$ for fixed $\eps$, we may also assume without loss of generality that $r = e^{- k_0}$ with $k_0\in \N$. We then decompose the ball $B(0,r)$ into the disjoint annuli $B_k = B(0, e^{-k} ) \setminus B(0, e^{- k -1})$ so that $B(0, r) = \cup_{k \ge k_0} B_k$. Note that for any fixed $k \ge k_0$ one has
%$$
%e^{\gamma H(z)} \asymp |z|^{- \alpha \gamma} \asymp e^{k \alpha \gamma} ;
%$$
$$e^{\gamma H_\eps(z)} \lesssim e^{k\alpha \gamma}$$
for $z\in B_k$, where the implied constant is uniform over $z\in B_k$ and $\eps>0$.
It follows that
\begin{align*}
\E [ (Z_{r, \eps})^p] &\lesssim \sum_{k = k_0}^{\lceil \log (1/\eps)\rceil} e^{ k \alpha \gamma p} \E [ (\cM_{h+(Q/2)\log \hat{g}_0} (B_k))^p]\\
& \lesssim  \sum_{k \ge k_0} e^{ k \alpha \gamma p} e^{- k \xi (p)},
\end{align*}
where $\xi(p) = p (2 + \gamma^2/2) -p^2 \gamma^2/2 $ is the multifractal spectrum function. Here we have used Theorem \ref{T:moments} and more precisely \eqref{scalingeps}, together with the obvious
fact that for $p<1$, by Jensen's inequality,
$$\E[(\cM_{h+(Q/2)\log \hat{g}_0} (B(0,1)))^p] \le \E [ \cM_{h+(Q/2)\log \hat{g}_0} (\C) ]^p    \lesssim 1.$$

As a consequence, the claim \eqref{E:goalStep3LCFT} follows if we can find $0<p<1$ sufficiently small so that
$$
\alpha \gamma p - \xi (p) <0.
$$
Linearising $\xi(p)$ around $p =0$, it suffices that
\begin{equation}\label{E:proofStep3}
\alpha \gamma - (2+ \gamma^2/2) <0.
\end{equation}
But from the Seiberg bound \eqref{E:Seiberg2}, since $\alpha< Q = (\gamma/2 + 2 / \gamma)$, we see that $\alpha Q < 2 + \gamma^2/2$ so that \eqref{E:proofStep3} is fulfilled.

\medskip \textbf{Step 4.} We now conclude the proof of Lemma \ref{L:negmomentsLCFT}. Recall from Step 2 that the limit in probability $Z_{r,0}^c = \lim_{\eps \to 0} Z_{r, \eps}^c $ satisfies $Z_{r,0}^c = \int_{A_r^c } e^{\gamma H(z)} \cM_{h+(Q/2)\log \hat{g}_0} (\dd z)$ and is thus monotone increasing in $r$. We can therefore set $$Z_0 = \lim_{r \to 0} Z_{r, 0}^c,$$
where the above limit is in the almost sure sense.
By Markov's inequality and \eqref{E:goalStep3LCFT} of Step 3, we also have that
 $$\P( Z_{r, \eps} > C(r)^{p/2} ) \le  C(r)^{p/2} \to 0$$
 as $r\to 0$, uniformly in $\eps$.
Using the triangle inequality and taking $r$ small and then $\eps$ small, we deduce that $$Z_\eps \to Z_0 \text{ in probability as } \eps\to 0.$$ Moreover, taking the difference, we see that $Z_{r, \eps}  = Z_\eps - Z_{r,\eps}^c \to Z_{r, 0} := Z_0-Z_{r,0}^c$ %\int_{A_r} e^{\gamma H(z) } \MM_h(\dd z)$
in probability for each $r>0$.

\medskip From Step 1, we see that for our fixed $s:=\gamma^{-1}(\sum \alpha_i-2Q)$, $\E( (Z_\eps)^{-2s}) $ is uniformly bounded in $\eps$, so that $(Z_\eps)^{-s}$ is uniformly integrable and hence $\E((Z_\eps)^{-s} ) \to \E((Z_0)^{-s})<\infty $ as $\eps \to 0$. %Step 1 again shows that this limit is necessarily finite.
To show that the limit is also non-zero, it suffices to show that $Z_0<\infty$ (so that $(Z_0)^{-s} >0$ and hence its expectation is also strictly positive). For this we take $r =1$, and write $Z_0 = Z_{1,0} + Z_{1,0}^c$ in the notation of Step 2. The second term has finite expectation and is therefore finite almost surely. The first term has finite $p$th moment for sufficiently small $p>0$ by \eqref{E:goalStep3LCFT} and Fatou's lemma. It is therefore also finite almost surely. We deduce that $Z_0 <\infty$, which concludes the proof of Lemma \ref{L:negmomentsLCFT}.
\end{proof}
Plugging Lemma \ref{L:negmomentsLCFT} into \eqref{LCFTgamma}, we conclude the proof of Theorem \ref{T:Seibergbounds}.
\end{proof}

\begin{rmk}
It can be shown that the Seiberg bounds are sharp in the following sense. If $\sum_i \alpha_i \le 2Q$, then $\langle V_\eps \rangle_{\hat g_0} = \infty$, while if  $\max_i \alpha_i \ge Q$, then the random variable $Z_0$ in Lemma \ref{L:negmomentsLCFT} is almost surely infinite, so that its negative moment of order $s$ is zero, hence $\langle V_\eps \rangle_{\hat g_0} \to 0$.  See (3.17) in \cite{DKRV}. Thus, the Seiberg bounds are a \emph{necessary and sufficient condition} for the correlation functions to be well defined (at least without further normalisation).
\end{rmk}
%\nb{To-do: possibly, add the term $Q/2 \log g$ in the definition of the correlation function for the curvature term.}

In the course of the proof, we obtained a very important expression for the value of the correlation function $\langle V\rangle_{\hat g_0}=\lim_{\eps\to 0} \langle V_\eps \rangle_{\hat g_0}$. This shows that the correlation function of the model (which in a few moments we will view as the partition function of a random field) can be computed exactly as some fractional negative moment of a Gaussian multiplicative moment and the Gamma function, and is the first hint of the remarkable \textbf{integrability} of the model. It is worth restating this expression as a separate corollary.

\begin{cor} \label{C:Vformula}  Suppose that $z_1, \ldots, z_k \in \C$ are distinct, and $\alpha_1, \ldots, \alpha_k$ satisfy the Seiberg bounds \eqref{E:Seiberg1} and \eqref{E:Seiberg2}. Write $V = V_{\alpha_1, \ldots, \alpha_k} (\mathbf{z})$ for the associated correlation function.
Set \begin{equation}\label{E:tildeh} \tilde{h}^{\hat\C}=h^{\hat\C,\hat g_0}+\tfrac{Q}{2}\log \hat{g}_0 + 2\pi \sum_{i=1}^k \alpha_i  G^{\hat\C,\hat g_0}(\cdot, z_i);\end{equation} 
\begin{equation} 
C_\alpha(\mathbf{z})=  2\pi \sum_{i=1}^k\sum_{j=i+1}^k {\alpha_i\alpha_j} G^{\hat\C,\hat g_0}(z_i,z_j)+\sum_{i=1}^k \alpha_i^2 c_{\hat g_0}  ;  
\end{equation} 
where we recall that $c_{\hat g_0}=\log(2)-1/2$; and 
\begin{equation}\label{E:conformalweights} \Delta_{\alpha}=\tfrac{\alpha}{2}(Q-\tfrac{\alpha}{2})
\quad ; \quad s=\frac{\sum_{i=1}^k \alpha_i -2Q}{\gamma}. 
\end{equation} \indN{{\bf Liouville CFT}! $\Delta_\alpha$; conformal weights} \ind{Conformal weights}
Then we have
	\begin{equation}\label{E:Vgform} \langle V \rangle_{\hat g_0} = \gamma^{-1} e^{C_\alpha(\mathbf{z})} \prod_i \hat{g}_0(z_i)^{\Delta_{\alpha_i}}  	\Gamma \left(s; \mu\right)
	\E\Big( \cM_{\tilde{h}^{\hat\C}}(\C)^{-s} \Big) ,
	 \end{equation}
	where $\Gamma(s; \mu) = \int_0^\infty u^{s-1} e^{-\mu u} \dd u$\indN{{\bf Miscellaneous}! $\Gamma(s,\mu)$; the Gamma function with parameters $s$ and $\mu$} is the Gamma function as before. Here, following the general notation in the chapter, $\cM_{h}(\C)=\lim_{\eps\to 0} \int_\C \eps^{\gamma^2/2}e^{\gamma h_\eps(z)} \dd z$.

  More generally, if $F$ is a non-negative measurable functional on $H^{-1} (\hat \C)$, 
we have  	\begin{align}\langle V F \rangle_{\hat g_0} & := \lim_{\eps\to 0} \langle V_\eps F \rangle_{\hat g_0}  = \gamma^{-1} e^{C_\alpha(\mathbf{z})} \prod_i \hat{g}_0(z_i)^{\Delta_{\alpha_i}}  \nonumber \\ & \times 
	\int_{u>0} \E\big( F(\tilde{h}^{\hat \C}  + \tfrac{\log u}{\gamma} -\tfrac{ \log \cM_{\tilde{h}^{\hat \C}}(\C)}{\gamma}) \cM_{\tilde{h}^{\hat \C}}(\C)^{-s}\big) u^{s-1}e^{-\mu u} \, \dd u .\label{E:VF} \end{align}
\end{cor}

\begin{rmk}\label{R:Vformula}
	One can check that the proof of \cref{T:Seibergbounds} goes through unchanged when $\hat{g}_0$ is replaced by a metric $\hat g$ as in \eqref{E:metric} with constant scalar curvature. This results in analogous explicit expressions for $\langle V \rangle_{\hat g}$ and $\langle F V \rangle_{\hat g}$ (the constant $c_{\hat g_0}$ being replaced by its general expression $c_{\hat g}$ from \cref{R:constantcurvature}). This is in particular the case when $\hat g$ is the pushforward of $\hat g_0$ by a M\"obius map $m$, and note furthermore that in this case $c_{\hat g} = c_{\hat g_0}$ (see \eqref{L:thetaRcmg}).
\end{rmk}

	The Weyl anomaly formula of Theorem \ref{T:Weyl} was written for general ``smooth'' observable $F $ on $H^{-1} (\hat\C)$. But it is immediate to deduce from it a formula which includes the correlation functions. 
	
	\begin{cor}\label{C:Weyl_correlation} Let $V = V_{\alpha_1, \ldots, \alpha_k} (\mathbf{z})$ be the correlation functions as above; and let $F$ be an arbitrary non-negative measurable functional on $H^{-1} (\hat \C)$, and let $\hat g = e^{\rho} \hat g_0$ where $\rho$ is a twice differentiable function on $\C$ with a finite limit at infinity and $\int_{\C} |\nabla \rho(z) |^2 \dd z < \infty$, as in \eqref{E:metric}. Then 
	\[
	\langle V F \rangle_{\hat g}=\exp\left(\frac{6Q^2}{96\pi}\int_{\C} [|\nabla^{\hat g_0}\rho(x)|^2+2R_{\hat g_0}\rho(x)]v_{\hat g_0}(\dd x)\right) \;  \langle V  F\rangle_{\hat g_0}.
	\]
	\end{cor}

\begin{proof}
By \cref{T:Weyl} we have $\langle V_\eps F \rangle_{\hat g}=\exp(\frac{6Q^2}{96\pi}\int_{\C} [|\nabla^{\hat g_0}\rho(x)|^2+2R_{\hat g_0}\rho(x)]v_{\hat g_0}(\dd x)) \langle V_\eps F\rangle_{\hat g_0}$ for every $\eps>0$, so the result follows by taking $\eps\to 0$.
\end{proof}

\begin{rmk} \label{R:Weyl_correlation}
The above Weyl anomaly formula is partly a consequence of how we chose to define the correlation functions $V$ in \eqref{eq:Vepsilon1}, since implicitly the chosen regularisation in \eqref{eq:Vepsilon2} is given in terms of the Euclidean metric rather than the intrinsic metric $g$. If instead one thinks of choosing circle averages with respect to the metric $g$, a more intrinsic definition of $\langle V_\eps\rangle_{\hat g}$ would be
\[
\langle V_\eps \rangle_{\hat g} = \int \E \left[ \prod_{i=1}^k ( \sqrt{ \hat g(z_i)} \eps)^{\alpha_i^2/2} e^{ \alpha_i ( h_\eps(z_i) + c ) } \exp ( - 2 Q c - \mu e^{\gamma c} \cM_{h + \tfrac{Q}{2} \log \hat g } (\hat \C) ) \right] \dd c.
\]
Note that compared to  \eqref{eq:Vepsilon1} there is an extra factor $ \sqrt{ \hat g(z_i)}$ in front of the normalising factor $\eps$ and that the term $\tfrac{Q}{2} \log \hat g(z_i)$ is \emph{not} included in the first exponential term. 
This would lead to a slightly different version of the Weyl anomaly formula: namely, 
\[
	\langle V F \rangle_{\hat g}=\exp\left(\frac{6Q^2}{96\pi} A( \rho, \hat g_0) - \sum_{i=1}^k \Delta_{\alpha_i} \rho(z_i) \right) \;  \langle V  F\rangle_{\hat g_0}.
\]
where $A( \rho, \hat g_0) =  \int_{\C} [|\nabla^{\hat g_0}\rho(x)|^2+2R_{\hat g_0}\rho(x)]v_{\hat g_0}(\dd x)$ is as in \cref{C:Weyl_correlation}. This version of the Weyl anomaly formula is for instance the one that is used in \cite[(1.3)]{GKRV_segal}.
\end{rmk}

As a consequence of \cref{R:Vformula} and \cref{C:Weyl_correlation}, together with \cref{L:Gmob} and \cref{C:hs_mobius}, we obtain the following theorem describing how $\langle V \rangle_{\hat{g}_0}$ changes when the insertions $\{z_i\}_i$ are transformed using a M\"{o}bius map. This transform is described in \cite{DKRV} as a version of the KPZ formula, cf. \cref{T:KPZ}. 

\begin{theorem}\label{T:LfieldMobius}
	Suppose that $m:\hat\C \to \hat \C$ is a M\"{o}bius transformation of the Riemann sphere, and $\alpha_i, z_i$ are as in \cref{C:Vformula}. Then 
	\begin{equation}\label{E:LFieldMobius}
		\langle V_{\alpha_1,\dots, \alpha_k} (m(\mathbf{z})) \rangle_{\hat{g}_0}=\prod_i |m'(z_i)|^{-2\Delta_{\alpha_i}}\langle V_{\alpha_1,\dots, \alpha_k} (\mathbf{z}) \rangle_{\hat{g}_0}
	\end{equation}
where $m(\mathbf{z})=(m(z_1),\dots, m(z_k))$. Moreover, if $F$ is a non-negative functional on $H^{-1} (\hat \C)$, 
	\begin{equation}\label{E:LFieldMobiusF}
	\langle F V_{\alpha_1,\dots, \alpha_k} (m(\mathbf{z})) \rangle_{\hat{g}_0}=\prod_i |m'(z_i)|^{-2\Delta_{\alpha_i}}\langle F_m V_{\alpha_1,\dots, \alpha_k} (\mathbf{z}) \rangle_{\hat{g}_0}
\end{equation}
where $F_m(h):=F(h\circ m^{-1} + Q\log |(m^{-1})'|)$ for $h\in H^{-1}(\hat\C)$.
\end{theorem}

\begin{proof}
%	See exercise \ref{Ex:mobius}.
	By setting $F=1$, it suffices to prove the second statement. The Weyl anomaly formula \cref{L:weylmobius} gives that 
\[
	\langle FV_{\alpha_1,\dots, \alpha_k} (m(\mathbf{z})) \rangle_{\hat{g}_0}  = 	\langle FV_{\alpha_1,\dots, \alpha_k} (m(\mathbf{z})) \rangle_{m_*\hat{g}_0}.
\]
On the other hand, using Remark \ref{R:Vformula} and using that $c_{m_*\hat g_0}=c_{\hat g_0} = \log(2)-1/2$, we see that 
\begin{multline*}\langle FV_{\alpha_1,\dots, \alpha_k} (m(\mathbf{z})) \rangle_{m_*\hat{g}_0} = \gamma^{-1} e^{2\pi \sum_{i=1}^k \sum_{j=i+1}^k \alpha_i\alpha_j G^{\hat \C,m_*\hat g_0}(m(z_i),m(z_j))+\sum_i \alpha_i^2(\log(2)-1/2)} \nonumber \\  \times  \prod_i (m_*\hat{g}_0)(m(z_i))^{\Delta_{\alpha_i}} 
	\int_{u>0} \E\big( F(\tilde{h}  + \tfrac{\log u}{\gamma} -\tfrac{ \log \cM_{\tilde{h}}(\C)}{\gamma}) \cM_{\tilde{h}}(\C)^{-s}\big) u^{s-1}e^{-\mu u} \, \dd u  \end{multline*}
where 
$$
\tilde{h}=h^{\hat\C,m_*\hat g_0}+\tfrac{Q}{2}\log m_*\hat{g}_0 + 2\pi \sum_{i=1}^k \alpha_i  G^{\hat\C,m_*\hat g_0}(\cdot, m(z_i)).
$$
Recall \cref{R:hs_mobius}, which says that 
$$h^{\hat \C,m_*\hat g_0} \overset{d}{=} h^{\hat \C,\hat g_0}\circ m^{-1}; \text{ and }G^{\hat \C,m_*\hat g_0}(m(z),m(w))=G^{\hat{\C},\hat g_0}(z,w)\text{ for $z\ne w$ in $\hat\C$.} $$ 
 Also using the explicit form $m_*\hat{g}_0(x)=\hat{g}_0(m^{-1}(x))|(m^{-1})'(x)|^2$, we therefore have that 
$$
\tilde{h}\overset{(d)}{=}\tilde{h}^{\hat \C}\circ m^{-1} +Q\log (|(m^{-1})'|)
$$
and 
$$2\pi \sum_{i=1}^k \sum_{j=i+1}^k \alpha_i\alpha_j G^{\hat \C,m_*\hat g_0}(m(z_i),m(z_j))+\sum_i \alpha_i^2(\log(2)-1/2)=C_\alpha(\mathbf{z})$$
where $\tilde{h}^{\hat \C}$ and $C_{\alpha}(\mathbf{z})$ are as in \cref{C:Vformula}. Finally, observe that 
\begin{align*}
\prod_i (m_*\hat{g}_0)(m(z_i))^{\Delta_{\alpha_i}} &= \prod_i (\hat g_0(z_i))^{\Delta_{\alpha_i}}\prod_i |(m^{-1})'(m(z_i))|^{2\Delta \alpha_i}\\\
& = \prod_i (\hat g_0(z_i))^{\Delta_{\alpha_i}}\prod_i |m'(z_i)|^{-2\Delta \alpha_i},
\end{align*}
which yields
\[
\langle FV_{\alpha_1,\dots, \alpha_k} (m(\mathbf{z})) \rangle_{\hat{g}_0}  = 	\langle FV_{\alpha_1,\dots, \alpha_k} (m(\mathbf{z})) \rangle_{m_*\hat{g}_0}=\prod_i |m'(z_i)|^{-2\Delta \alpha_i}\langle F_m V_{\alpha_1,\dots,\alpha_k}(\mathbf{z})\rangle_{\hat g_0},
\]
as desired.
\end{proof}

\begin{cor}
Set $k = 3$ and suppose $\alpha_1,\alpha_2$ and $\alpha_3 $ satisfy the Seiberg bounds. 
Then there exist a constant $C(\alpha_1, \alpha_2, \alpha_3)>0$ called the \textbf{structure constant} such that for any $z_1, \ldots, z_3 \in \hat \C$ 
$$
\langle V_{\alpha_1,\alpha_2, \alpha_3} (\mathbf{z}) \rangle_{\hat{g}_0} = C(\alpha_1, \alpha_2, \alpha_3) |z_1 - z_2|^{2\Delta_{1,2}} | z_2 - z_3|^{2\Delta_{2,3}} | z_3 - z_1|^{2\Delta_{1,3}}
$$
where $\Delta_{1,2} = \Delta_{\alpha_3} - \Delta_{\alpha_1} - \Delta_{\alpha_2}$, $\Delta_{2,3} = \Delta_{\alpha_1} - \Delta_{\alpha_2} - \Delta_{\alpha_3}$, and $\Delta_{1,3} = \Delta_{\alpha_2} - \Delta_{\alpha_1} - \Delta_{\alpha_3}$.
\end{cor}

\begin{proof}
Since $k =3$ we can find a M\"obius map $m$ sending $z_1 \mapsto 0$, $z_2 \mapsto 1$, $z_3 \mapsto \infty$. The map $m$ has an explicit form, namely 
$$
m(z) = \frac{z-z_1}{z- z_3} \frac{z_2 - z_3}{ z_2 - z_1}. 
$$ 
Note that then 
$$
m'(z)= \frac{z_1 - z_3}{ (z-z_3)^2} \frac{z_2 - z_3}{z_2 - z_1}. 
$$
Hence 
$$
m'(z_1) = \frac{z_2 - z_3}{(z_2- z_1) (z_1 - z_3)}; m'(z_2) = \frac{z_1 - z_3}{(z_2- z_1) (z_2 - z_3)}. 
$$
while 
$$
m'(z_3) = \infty; \text{ and }m'(z) \sim \frac{(z_1 - z_3) (z_2 - z_3)}{z_2 - z_1} \times \frac1{(z-z_3)^2}
$$
as $z \to z_3$. Let us evaluate \eqref{E:LFieldMobius} at $\mathbf{z} = z_1, z_2, z$ with $z \to z_3$. 
Then writing $\Delta_i = \Delta_{\alpha_i}$, we have
\begin{equation*}
	 \frac{\langle V_{\alpha_1,\alpha_2, \alpha_3} (0,1, m(z) ) \rangle_{\hat{g}_0}}{ \langle V_{\alpha_1,\alpha_2,\alpha_3}(z_1,z_2,z)\rangle_{\hat g_0}}\! \sim \!  |z- z_3|^{4\Delta_3} |z_1 - z_2|^{2(\Delta_1 + \Delta_2 + \Delta_3)} |z_2 - z_3|^{2(- \Delta_1 -  \Delta_3 +\Delta_2) } |z_1 - z_3|^{2(\Delta_1 - \Delta_2 -\Delta_3)}
\end{equation*}
as $z\to z_3$.
From this we learn that, writing $y = m(z)$,
$$
\lim_{y\to \infty} |y|^{4 \Delta_3} \langle V_{\alpha_1,\alpha_2, \alpha_3} (0,1, y ) \rangle_{\hat{g}_0}  = C(\alpha_1, \alpha_2, \alpha_3)
$$
exists, and equals 
$$
\langle V_{\alpha_1,\alpha_2, \alpha_3} (\mathbf{z}) \rangle_{\hat{g}_0} |z_1 - z_2|^{ - 2\Delta_{1,2}} | z_2 - z_3|^{- 2\Delta_{2,3}} | z_3 - z_1|^{-2\Delta_{1,3}}.
$$
This concludes the proof. 
\end{proof}

\begin{definition}[Correlation functions with an insertion at $\infty$]
	\label{R:Vangleinfinity} Generalising the argument in the proof above we see that for $\alpha_1,\dots,\alpha_k$ satisfying \eqref{E:Seiberg1} and \eqref{E:Seiberg2}, 
\[
\lim_{y\to \infty} |y|^{4\alpha_k} \langle V_{\alpha}(z_1,\dots, z_{k-1},y) \rangle_{\hat{g}_0}=:\langle V_\alpha(z_1,\dots, z_{k-1},\infty)\rangle_{\hat{g}_0}
\]
exists.\end{definition}

 It then follows from taking a limit as $y\to \infty$ in \cref{C:Vformula}, with $z_k=y$, that for any non-negative Borel function on $H^{-1}(\hat \C)$:
	\begin{align} & \langle V_\alpha(z_1,\dots, z_{k-1},\infty) F \rangle_{\hat g_0} = c \gamma^{-1} e^{C_{\alpha_1,\dots, \alpha_{k-1}}(z_1,\dots, z_k)} \prod_{i=1}^{k-1} \hat{g}_0(z_i)^{\Delta_{\alpha_i}-\tfrac{\alpha_i\alpha_k}{4}}  \nonumber \\ & \times 
	\int_{u>0} \E\big( F(\tilde{h}^{\hat \C}  + \tfrac{\log u}{\gamma} -\tfrac{ \log \cM_{\tilde{h}^{\hat \C}}(\C)}{\gamma}) \cM_{\tilde{h}^{\hat \C}}(\C)^{-s}\big) u^{s-1}e^{-\mu u} \, \dd u .\label{E:VFinfinity} \end{align}
with $c=c(\alpha)$ depending only on $\alpha$ (and not $F$ or $z_1,\dots, z_{k-1}$), $s=\gamma^{-1}(\sum_{i=1}^k \alpha_i -2Q)$, and
\begin{equation}\label{E:tildehinfinity} \tilde{h}^{\hat\C}=h^{\hat\C,\hat g_0}+(\tfrac{Q}{2}-\tfrac{\alpha_k}{4})\log \hat{g}_0 + 2\pi \sum_{i=1}^{k-1} \alpha_i  G^{\hat\C,\hat g_0}(\cdot, z_i).
\end{equation}
Moreover, if $m(z)=az+b$ is a M\"{o}bius transformation fixing $\infty$, then 
\begin{equation}\label{eq:mob_infinity}
	\langle V_\alpha(m(z_1),\dots, m(z_{k-1}),\infty)\rangle_{\hat g_0}=\prod_{i=1}^{k-1} |m'(z_i)|^{-2\Delta_{\alpha_i}} a^{2\Delta_{\alpha_k}} \langle V_{\alpha}(z_1,\dots, z_{k-1},\infty) \rangle_{\hat{g}_0}.
\end{equation}

\subsection{An alternative choice of background metric}
\label{SSS:zerocircleg}

It is common in the literature on spherical Liouville CFT (for example in \cite{RhodesVargas_tail,Vargas_LN,DOZZ}) to define the correlations starting with the field $h^\fc$, the whole plane GFF with zero average on the unit circle, in place of $\hat{h}^{\hat\C,\hat g_0}$. One can, for instance, define $h^\fc$ by setting $h^\fc:=h^{\hat\C,\hat g_0}-(h^{\hat \C,\hat g_0},\rho_1)$ where $\rho_1$ is uniform measure on the unit circle. 

Let us begin by computing the covariance $G^\fc$ of $h^\fc$. 

\begin{lemma}\label{L:Ginfty}
	\begin{equation}\label{E:Ginfty} 2\pi G^{\fc}(x,y)=-\log|x-y|+\log(|x|\vee 1)+\log(|y|\vee 1)
	\end{equation}
	for $x\ne y\in \C$. 
\end{lemma}\indN{{\bf Green functions}! $G^{\fc}$; Green function for whole plane GFF with zero average on the unit circle}
\begin{proof}From the definition of $h^\fc$, we have
$$
 G^\fc (x,y)
 \! = \!G^{\hat \C, \hat g_0 }(x,y) -\! \int \! G^{\hat \C,\hat g_0}(x,z)\rho_1(\dd z) - \! \int \! G^{\hat \C,\hat g_0}(y,z)\rho_1(\dd z) +\! \iint \!G^{\hat \C, \hat g_0} (z,w) \rho_1(\dd z) \rho_1(\dd w). 
$$
Recall that a formula for $G^{\hat \C,\hat g_0}(x,z)$ is provided in \eqref{E:greensphere}. Furthermore, we claim that $\log \hat g_0(z) =1$ if $|z| =1$, and for $x\in \C$,
\begin{equation}\label{eq:int_log}
\int \log |x- z| \rho_1( \dd z) = \log(|x| \vee 1).  
\end{equation}
To justify \eqref{eq:int_log} we consider the cases $|x|>1$ and $|x|< 1$ separately (the case $|x|=1$ follows by dominated convergence and continuity). When $|x| > 1$, \eqref{eq:int_log} is straightforward by harmonicity of $\log |x- \cdot|$ in $B(0,1)$. When $|x| < 1$, we note that for $z\in\partial B(0,1)$,
$$
|x-z| = | x - z | |\bar z| = |1- \bar z x | = | 1- \bar x z|, 
$$
so that 
$$
\int \log | x - z| \rho_1( \dd z) = \int \log | 1 - \bar   x z | \rho_1( \dd z). 
$$
As $z$ varies across the unit circle, $ 1 - \bar x z$ varies across a circle centred at 1 of radius $|x| <1$. Using harmonicity of the log function we deduce that the right hand side is 0, which proves \eqref{eq:int_log}. 

Together with \eqref{E:greensphere}
this immediately implies that
\begin{equation}\label{eq:covcircle}
2\pi \int G^{\hat \C,\hat g_0}(x,y)\rho_1(\dd y) = -\log(|x|\vee 1)-\tfrac14\log \hat g_0(x) + \log(2)-1/2,\quad  \text{ for } x\in \C.
\end{equation}
The lemma follows. 
\end{proof} 

Recalling \eqref{E:gSequiv}, $h^{\fc}$ formally corresponds to the GFF with zero average for the metric $\hat g = \hat g_{\fc}$, $\hat{g}_\fc(z)=(|z|\vee 1)^{-4}$. (This metric corresponds to gluing two copies of a unit Euclidean disc along their boundaries, but this fact will not be used below.) As this metric is not of the form $e^\rho \hat{g}_0$ with $\rho$ twice differentiable, it does not quite fit into the earlier framework. 

Nonetheless if we set 
\begin{equation}\label{E:correlationsfc}
	\langle F\rangle_{\fc} =   \int_{c \in \R} \E\left[ F\left( h^{\fc} -2Q \log (|\cdot|\vee 1) + c \right)\exp\Big( -2Qc   - \mu e^{\gamma c} \cM_{h^{\fc}-2Q\log(|\cdot|\vee 1)} (\C) \Big) \right] \dd c,
\end{equation}
analogously to \eqref{E:Fgnice}, then we can prove the following.

\begin{lemma}
	\begin{equation}\label{E:Fg0isFgc}
		\langle F \rangle_{\fc}=e^{-2Q^2(\log(2)-1/2)}\langle F \rangle_{\hat g_0}
	\end{equation}
	for all  non-negative Borel functions $F$ on $H^{-1}(\hat \C)$. Moreover, for $V=V_\alpha(\mathbf{z})$ as in \cref{T:Seibergbounds}, 
	\[\langle VF \rangle_\fc :=\lim_{\eps\to 0} \langle V_\eps F \rangle_\fc\]
	exists and can be explicitly expressed as 
	\begin{multline} \label{E:VFc}
		\langle VF \rangle_\fc = \gamma^{-1} \prod_{i=1}^k (|z_i|\vee 1)^{-4\Delta_{\alpha_i}+\alpha_i \sum_{l\ne i} \alpha_l} \prod_{j=i+1}^k |z_i-z_j|^{-\alpha_i\alpha_j} \notag \\ 
		\times \int_{u>0} \mathbb{E}(F(\tilde{h}^\fc+\gamma^{-1}(\log u - \log \cM_{\tilde h^\fc}(\C)))\cM_{\tilde h^\fc}(\C)^{-s}) u^{s-1}e^{-\mu u} \dd u,    
	\end{multline} 
	with $s=\gamma^{-1}(\sum \alpha_j - 2Q)$ and 
	$
	\tilde{h}^\fc=h^\fc -2Q\log(|\cdot|\vee 1)+\sum \alpha_i 2\pi G^{\fc}(\cdot,z_i).
	$
	Equivalently, 
	$
	\tilde{h}^{\fc}=h^\fc+(\sum \alpha_i -2Q)\log(|\cdot|\vee 1)-\sum \alpha_i \log|\cdot-z_i|+\sum \alpha_i \log(|z_i|\vee 1).
	$
\end{lemma}

\begin{proof}
We start by proving \eqref{E:Fg0isFgc} using Girsanov's theorem. Sine $h^{\fc}$ is equal in law to $h^{\hat\C,\hat g_0}-(h^{\hat\C,\hat g_0},\rho_1)$ we can write 
\begin{multline*}
	\langle F\rangle_{\fc} =   \int_{c \in \R} \E \Big[ F\left( h^{\hat\C,\hat g_0}-(h^{\hat\C,\hat g_0},\rho_1) -2Q \log (|\cdot|\vee 1) + c \right) \\ \exp\Big( -2Qc   - \mu e^{\gamma c} \cM_{h^{\hat\C,\hat g_0}-(h^{\hat\C,\hat g_0},\rho_1)-2Q\log(|\cdot|\vee 1)} (\C) \Big) \Big] \dd c
\end{multline*}
which by applying the change of variables $\hat c=\hat c- (h^{\hat\C,\hat g_0},\rho_1)$, is equal to 
\begin{multline*}
 \int_{\hat c \in \R} \E \Big[\! F \left(\! h^{\hat\C,\hat g_0}\!-\!2Q \log (|\cdot|\vee 1) + \hat c \right) \\
  \exp\Big(-2Q(h^{\hat\C,\hat g_0},\rho_1) -2Q\hat c   - \mu e^{\gamma \hat c} \cM_{h^{\hat\C,\hat g_0}-2Q\log(|\cdot|\vee 1)} (\C) \Big) \Big]\! \dd \hat c.
\end{multline*}
Defining $\tilde{\P}$ by 
\[ \frac{\dd\tilde{\P}}{\dd\P}:=\frac{\exp(-2Q(h^{\hat\C,\hat g_0},\rho_1))}{\mathbb{E}(\exp(-2Q(h^{\hat\C,\hat g_0},\rho_1)))}=\frac{\exp(-2Q(h^{\hat\C,\hat g_0},\rho_1))}{\exp(2Q^2(\log(2)-1/2)}\]
(recall \eqref{eq:covcircle}) we obtain that
\begin{multline*}
\langle F\rangle_{\fc} = e^{2Q^2(\log(2)-1/2)}	\int_{\hat c \in \R} \tilde{\E} \Big[ F\left( h^{\hat\C,\hat g_0}-2Q \log (|\cdot|\vee 1) + \hat c \right)  \\
\exp\Big(-2Q\hat c   - \mu e^{\gamma \hat c} \cM_{h^{\hat\C,\hat g_0}-2Q\log(|\cdot|\vee 1)} (\C) \Big) \Big] \dd \hat c.
\end{multline*}
By Girsanov's theorem, the law of $h^{\hat \C,\hat g_0}$ under $\tilde{\P}$ is the same as under $\P$, except with mean shifted by 
\[-2Q\int 2\pi G^{\hat \C,\hat g_0}(\cdot,y)\rho_1(\dd y)=2Q\log(|\cdot|\vee 1)+(Q/2)\log \hat g_0 -2Q(\log(2)-\tfrac{1}{2})\]
(again using \eqref{eq:covcircle}). This yields that
\begin{multline*}
	\langle F\rangle_{\fc} = e^{2Q^2(\log(2)-\tfrac{1}{2})}	
\int_{\hat c \in \R} {\E} \Big[ F\left( h^{\hat\C,\hat g_0}+\tfrac{Q}{2}\log \hat g_0 + (\hat c-2Q(\log(2)-\tfrac{1}{2}))\right) \times
\\
 \exp\!\Big(\!-2Q(\hat c-2Q(\log(2)-\tfrac{1}{2})) +4Q^2(\log(2)-\tfrac{1}{2})  - \mu e^{\gamma (\hat c-2Q(\log(2)-\tfrac{1}{2}))} \cM_{h^{\hat\C,\hat g_0}+\tfrac{Q}{2}\log \hat g_0} (\C) \Big) \! \Big]\! \dd \hat c.
\end{multline*}
Performing one final change of variables $c=\hat c- 2Q(\log(2)-\tfrac{1}{2})$, we obtain that 
\begin{multline*}\langle F \rangle_{\fc}=e^{-2Q^2(\log(2)-\tfrac{1}{2})}\int_{c \in \R} {\E} \Big[ F\left( h^{\hat\C,\hat g_0}+\tfrac{Q}{2}\log \hat g_0 + c \right) \\
\exp\Big(-2Qc  - \mu e^{\gamma c} \cM_{h^{\hat\C,\hat g_0}+\tfrac{Q}{2}\log \hat g_0} (\C) \Big) \Big] \dd  c,
\end{multline*}
which by \eqref{E:Fgnice} is exactly \eqref{E:Fg0isFgc}.

The explicit expression for 	$\langle VF \rangle_\fc :=\lim_{\eps\to 0} \langle V_\eps F \rangle_\fc$  follows from exactly the exact same argument as in the $\hat{g}_0$ case (\cref{T:Seibergbounds} and \cref{C:Vformula}). In summary:
\begin{itemize}
	\item 
The field $\tilde{h}^{\fc}$ is obtained by shifting $h^{\fc}-2Q\log(|\cdot|\vee 1)$ by $\sum_i \alpha_i (2\pi G^{\fc}(\cdot, z_i))$; this shift arises from Girsanov's theorem exactly as in \eqref{LCFTgirsanov}. $\tilde{h}^{\fc}$ is analogous to $\tilde{h}^{\hat\C}$ in \cref{C:Vformula}.
\item There is a compensation term $\lim_{\eps\to 0}\eps^{\alpha_i^2/2}\exp(\tfrac{1}{2}\var(\sum_i \alpha_i h^{\fc}_\eps(z_i))$ coming from Girsanov's theorem, as in \eqref{LCFTgirsanov}. Using the expression for $2\pi G^{\fc}$, this is equal to 
\[ \exp(\frac{1}{2}\sum_i \alpha_i^2 (2\log(|z_i|\vee 1)) - \sum_{i=1}^k \sum_{j=i+1}^k \alpha_i\alpha_j \log|z_i-z_j| +\sum_{i=1}^k \sum_{j\ne i} \alpha_i \alpha_k \log(|z_i|\vee 1))\]
which can be rewritten as 
$\prod_i (|z_i|\vee 1)^{\alpha_i^2-\alpha_i \sum_{j\ne i}\alpha_j} \prod_i\prod_{j=i+1}^k |z_i-z_j|^{-\alpha_i\alpha_j}$. There is also a term $e^{-\sum 2 \alpha_i Q\log(|z_i|\vee 1)}$ arising from the insertions, which is equal to $\prod_i (|z_i|\vee 1)^{-2Q\alpha_i}$. Putting these together gives the factor 
\[\prod_{i=1}^k (|z_i|\vee 1)^{-4\Delta_{\alpha_i}+\alpha_i \sum_{l\ne i} \alpha_l} \prod_{j=i+1}^k |z_i-z_j|^{-\alpha_i\alpha_j}\]
(analogous to $e^{C_\alpha(\mathbf{z})}\prod \hat g_0(z_i)^{\Delta_{\alpha_i}}$ in \cref{C:Vformula}).

\end{itemize} 

\end{proof}

We also have a similar expression when one of the insertions is at $\infty$ (we assume this is the $k$th and final insertion for simplicity). Recall \cref{R:Vangleinfinity}.

\begin{rmk}[Correlations with an insertion at infinity]
Suppose that $z_1,\dots, z_{k-1}$ are distinct and $\alpha_1,\dots, \alpha_k$ are as in \cref{T:Seibergbounds}. Then
	\begin{equation*} 
		\langle V_\alpha(z_1,\dots, z_{k-1},\infty) F \rangle_{\hat g_0}   =e^{2Q^2(\log(2)-1/2)}\langle V_\alpha(z_1,\dots, z_{k-1},\infty) F \rangle_\fc \end{equation*}
	is equal to 
	\begin{multline}  c' \gamma^{-1} \prod_{i=1}^{k-1} (|z_i|\vee 1)^{-4\Delta_{\alpha_i}+\alpha_i \sum_{l\ne i} \alpha_l}\prod_{j=i+1}^{k-1} |z_i-z_j|^{-\alpha_i\alpha_j} \notag \\  
		\times \int_{u>0}  \mathbb{E}(F(\tilde{h}^\fc+\gamma^{-1}(\log u - \log \cM_{\tilde h^\fc}(\C)))\cM_{\tilde h^\fc}(\C)^{-s}) u^{s-1}e^{-\mu u} \dd u  \label{E:VFcinf}\end{multline}
with $c'=c'(\alpha)$ depending only on $\alpha$ and 
\[ \tilde{h}^\fc=h^\fc+(\sum_{i=1}^k \alpha_i -2Q)\log(|\cdot|\vee 1)-\sum_{i\ne k} \alpha_i \log|\cdot-z_i|+\sum_{i\ne k} \alpha_i \log(|z_i|\vee 1).\]
\end{rmk}

\subsection{Geometric and probabilistic interpretation of Seiberg bounds}

\label{SS:geom_interpretation}

\ind{Gibbs measure}

The following discussion is intended to guide intuition but is not meant to be fully  rigorous. 

\medskip Assuming that the Seiberg bounds hold, the finite partition function $\langle V\rangle_{\hat g}$ in Theorem \ref{T:Seibergbounds} implicitly defines a random field (that is, sampled from a probability distribution) that we will soon call the \textbf{Liouville field}, see Definition \ref{D:Liouville}.  It is believed, as will be detailed more precisely in Remark \ref{R:mapconjecture}, that this Liouville field and its multiplicative chaos describe the scaling limit of suitably (that is, conformally) embedded random planar maps; the precise formulation of this conjecture goes back to \cite{DKRV}. This gives a \textbf{probabilistic justification} as to why insertions are necessary to define correlation functions, and why we need at least three of them. Indeed, conformal embeddings into the sphere are typically only unique up to M\"{o}bius transforms, so that in order to get a well defined, unique embedding, it is necessary to choose the embedded position of three vertices of the map in advance. It is natural to choose these three vertices uniformly at random on the planar map; in that case note that (because of Girsanov's theorem) their associated insertion weights should correspond to $\alpha_i = \gamma$. Conversely, if we take $k =3$ and $\alpha_i = \gamma$, the Seiberg bounds can only be satisfied when
$$
3 \gamma > 2 Q
$$ 
or equivalently 
$$
\gamma > \sqrt{2}.
$$
We remind the reader that the range $\gamma \in [ \sqrt{2}, 2)$ is exactly the range of values that one should obtain in the scaling limit of FK-decorated planar maps, see Section \ref{sec:RPM_LQG}.
\ind{FK model}

To appreciate the necessity of the insertions in order to get a finite partition function from a \textbf{geometric} point of view, it is useful to pause the exposition of the theory and to make a few heuristic considerations. 
Since the probability distribution associated to $\langle V \rangle_{\hat g}$ should formally be a Gibbs measure as in \eqref{E:gibbs}, it is intuitively useful to view this field as a random perturbation around the \textbf{ground state} of the theory, that is, the state $\ph$ of minimal energy, particularly when $\gamma \to 0$ and the field is essentially deterministic. To begin with, one might wonder what the ground state corresponding to the Polyakov action \eqref{E:Polyakov_informal} looks like without any insertion. Let us simply study the associated variational problem; that is, let $\ph$ be a minimiser of $S(\ph)$, let $f$ be an arbitrary test function, let $\eps>0$ and consider the action $S( \ph + \eps f)$. Then 
\begin{align*}
S( \ph + \eps f) & = \frac1{4\pi}\int \Big[|\nabla^{\hat g}(\ph + \eps f) |^2 + R_{\hat g} Q ( \ph + \eps f) + 4 \pi \mu e^{\gamma ( \ph + \eps f)} \Big] v_{\hat g}(\dd z)\\
& = \frac1{4\pi}\int \Big[|\nabla^{\hat g} \ph  |^2 + R_{\hat g} Q  \ph  + 4 \pi  \mu e^{\gamma \ph } \Big] v_{\hat g}(\dd z) + \\
& \quad \quad +  \frac{\eps}{4\pi}\int \Big[ 2\langle \nabla^{\hat g} \ph, \nabla^{\hat g} f \rangle +  Q R_{\hat g} f + 4\pi \mu \gamma e^{\gamma \ph} f \Big] v_{\hat g} (\dd z) + o(\eps)
\end{align*}
so that, by the Gauss--Green formula, since $\ph$ is a minimiser,
$$
\frac1{4\pi} \int \Big [ - 2 \Delta^{\hat g} \ph +  Q R_{\hat g} + 4 \pi \mu \gamma e^{\gamma \ph} \Big] f v_{\hat g} (\dd z) = 0.
$$
Because $f$ is arbitrary we deduce (first in the sense of distributions and then in the pointwise sense using elliptic arguments, which are not required here since this is anyway entirely heuristic), that
\begin{equation}\label{E:Liouville_eq_gamma}
\Delta^{\hat g} \ph = \frac{Q}2 R_{\hat g} + 2 \pi \mu \gamma e^{\gamma \ph}.
\end{equation}
Now let $u =\gamma \ph$, and consider the situation as $\gamma \to 0$ and $\mu \gamma^2 $ is kept constant (in particular, the cosmological constant $\mu$ tends to infinity, this is the so called semiclassical limit). Then we get an equation of the form
\begin{equation} \label{E:Liouville_eq}
\Delta^{\hat g} u = R_{\hat g} - K e^u ;\quad \quad  \text{ where } K = - 2 \pi \mu \gamma^2 <0.
\end{equation}
\ind{Liouville equation}
This equation is called \textbf{Liouville's equation}, and arises when searching for metrics $\tilde g$ conformally equivalent to $\hat g$ and of constant curvature $K$. Indeed, let us write $\tilde g = e^\rho  \hat g$ and suppose $R_{\tilde g} = K$. Since $R_{\tilde g} = - \Delta^{\tilde g} \log \tilde g $, this becomes
\begin{align*}
 - \Delta^{\tilde g} \log {\tilde g} & = K\\
 - e^{ - \rho} \Delta^{\hat g} \log (  \hat g e^\rho) & = K\\
 -\Delta^{ \hat g} \log  \hat g  - \Delta^{ \hat g} \rho &= K e^\rho\\
 R_{\hat g} - \Delta^{\hat g} \rho & = K e^\rho,
\end{align*}
which is the same as \eqref{E:Liouville_eq} with $\rho = u$ and $K = - 2\pi \mu \gamma^2$ (recall that $\mu \gamma^2$ is chosen to be a constant). Thus the Polyakov action is formally minimised by a function $\rho$ corresponding to a metric $\tilde g= e^\rho \hat g$ which has constant negative curvature $K$. But of course this is impossible on the sphere, in view of the Gauss--Bonnet theorem, which implies that the integral of the curvature should be $8 \pi$.

Formally, adding insertions in the computation of $\langle V \rangle$ can be thought of as changing the Polyakov action; the minimiser then satisfies 
$$
\frac1{4\pi} \int \Big [ - 2 \Delta^{\hat g}\ph +  Q R_{\hat g} + 4 \pi \mu \gamma e^{\gamma \ph} \Big] f v_{\hat g} (\dd z)  - \sum_{i=1}^k \alpha_i  f(z_i)= 0.
$$
or in other words,
\begin{equation}\label{E:Liouville_eq_gamma_insertions}
\Delta^{\hat g} \ph = \frac{Q}2 R_{\hat g} + 2 \pi \mu \gamma e^{\gamma \ph} - 2\pi  \sum_{i=1}^k \alpha_i \delta_{\{z_i\}}.
\end{equation}
instead of \eqref{E:Liouville_eq_gamma}, where $\delta_{\{z_i\}}$ is the Dirac mass on $\Sp$ (with respect to the underlying metric $\hat g$). Scaling the weights $\alpha_i$ by defining new weights $\tilde \alpha_i = \gamma \alpha_i$, we get that the new weights satisfy the rescaled Seiberg bounds:
$$
\sum_{i=1}^k \tilde \alpha_i > 2 \gamma Q, \tilde \alpha_i < \gamma Q
$$
which as $\gamma \to 0$ becomes
\begin{equation}
\label{rescaledSeiberg}
\sum_{i=1}^k \tilde \alpha_i > 4, \tilde \alpha_i < 2.
\end{equation}
Then with these rescaled weights, setting $u = \gamma \ph$ in \eqref{E:Liouville_eq_gamma_insertions} and letting $\gamma \to 0$ with $\mu \gamma^2$ constant as above (and $\tilde \alpha_i = \gamma \alpha_i$ fixed), the equation satisfied by $u$ becomes
\begin{equation}\label{E:Liouville_eq_insertions}
\Delta^{\hat g} u = R_{\hat g} + K e^{u} - 2\pi  \sum_{i=1}^k \tilde \alpha_i \delta_{\{z_i\}}.
\end{equation}
This modified form of Liouville's equation describes a metric $\tilde g = e^u \hat g$ such that $\tilde g$ has constant curvature $K = - 2 \pi \mu \gamma^2 $ away from the points $\{z_i\}_i$, but conical singularities at each of the $z_i$. That is, the metric is locally of the form $1/ | z_i - z|^{\tilde \alpha_i}$ as $z \to z_i$. In this setting there is no obstruction from Gauss--Bonnet to the existence of such metrics. 

\medskip This connection is made precise by Lacoin, Rhodes and Vargas \cite{LRV,LRVnew}. Indeed, the authors prove that in the limit $\gamma \to 0$ (with $\mu\gamma^2$ and $\gamma \alpha_i$ kept fixed as above) the associated normalised Liouville field concentrates near the solution of the modified Liouville equation \eqref{E:Liouville_eq_insertions}. Furthermore, they show that the fluctuations are asymptotically given by a massive Gaussian free field, and obtain a large deviation theorem where the rate function is given by the Polyakov action (shifted by the insertions); as expected this rate function is thus zero at the minimiser.

\subsection{Liouville fields}

As mentioned above, the finiteness of the partition function $\langle V \rangle_{\hat g}$ when the Seiberg bounds are satisfied, allows us not only to define ``expectations'' of observables such as $\langle VF \rangle_{\hat g} $ in Corollary \ref{C:Vformula}, but actually random \emph{fields} sampled from the associated probability distribution. 

\begin{definition}
  We define the Liouville field (associated to insertions $\mathbf{z} = (z_1, \ldots, z_k)$ and parameters $\alpha = (\alpha_1, \ldots, \alpha_k)$ satisfying the Seiberg bounds \eqref{E:Seiberg1} and \eqref{E:Seiberg2}) to be the random field $\haz$
  %\indN{{\bf Liouville CFT}! $\haz$; Liouville field on $\hat{\C}$} %law $\mathbb{P}^L_{\alpha, \mathbf{z}}$\
  in $H^{-1}(\hat\C)$, such that for any observable $F$,
  $$
  \E[F(\haz)] : = \frac{\langle FV_{\alpha}(\mathbf{z}) \rangle_{\hat g}}{\langle V_{\alpha}(\mathbf{z}) \rangle_{\hat g}}.
  $$ 
  This law does not depend on the choice of the metric ${\hat g}$.
% \index[functions]{$\P^L_\az$; law of the Liouville field}\index[functions]{$\E^L_{\az}$; expectation associated with the Liouville field}
\label{D:Liouville}
\end{definition} \ind{Liouville field}

The remarkable independence of this law from the metric ${\hat g}$ is a direct consequence of the Weyl anomaly formula (see Theorem \ref{T:Weyl} and Remark \ref{R:Vformula}). In what follows we will always work with the definition starting from the spherical metric $\hat g_0$ for simplicity. 

\begin{rmk}\label{R:Liouvillefieldc}
	By \eqref{E:correlationsfc} we also have $\mathbb{E}[F(\haz)]=\frac{\langle F V_{\alpha}(\mathbf{z})\rangle_\fc}{\langle V_{\alpha}(\mathbf{z}) \rangle_\fc}$, where $\langle \cdot \rangle_\fc$ is as defined in \cref{SSS:zerocircleg}.
\end{rmk}

\cref{T:LfieldMobius} also implies that the field transforms in the following way under M\"{o}bius transformations of the sphere. 

\begin{cor}\label{C:LfieldMobius}
Suppose that $m:\hat\C\to \hat \C$ is a M\"{o}bius transform of the Riemann sphere, and $\{\alpha_i,z_i\}$ are as in \cref{D:Liouville}. Then 
\begin{equation*}
	h^L_{\alpha,\mathbf{z}} \overset{(\mathrm{law})}{=} h^L_{\alpha,m(\mathbf{z})}\circ m + Q \log |m'|
\end{equation*}
where $m(\mathbf{z})=(m(z_1),\dots, m(z_k))$.
\end{cor}

Note that the law of $\haz$ also depends on $\gamma$, but we omit this from the notation (as with everything previously in this chapter).  A natural next question is to identify the law of the Liouville field in a way that is more explicit than the definition. We will not do this right away, but in the end we will get a very nice description by conditioning on the area; this will give us the \emph{unit volume Liouville sphere} in Section \ref{SS:LS}. For now, we will first make a simple (but surprising) observation about the law of the total mass of the multiplicative chaos measure associated to the Liouville field $h^L_{\alpha, \mathbf{z}}$.  
\begin{lemma}\label{L:Liouvillegamma}
Suppose that $\{\alpha_i,z_i\}$ are as in \cref{D:Liouville}. Then
$$\cM_{\haz}(\C)\sim \Gamma(s,\mu); \quad s = \tfrac{\sum_i \alpha_i-2Q}{\gamma}>0.$$
That is, $\cM_{\haz}(\C)$ has density proportional to
$u^{s-1}e^{-\mu u}\indic{u>0}$ with respect to Lebesgue measure on $\R$.
\end{lemma}
\begin{proof}
Recall the definition $\tilde{h}^{\hat\C}=h^{\hat\C,\hat g_0}+\tfrac{Q}{2}\log \hat{g_0} + \sum_{i=1}^k \alpha_i G^{\hat\C,\hat g_0}(\cdot, z_i)$.	For $A\subset \R_+$, we have that 
	\begin{align*}
		\P(\cM_{\haz}(\C)\in A)& = \frac{\langle V_\az \mathbf{1}_{\{\cM_{\haz}(\C)\in A\}}\rangle_{\hat g_0}}{\la V_\az \ra_{\hat g_0}} \\ 
		& = \frac{	\int_{u>0} \E\big( \mathbf{1}_{\{\cM_{h}(\C)\in A\}} \big(\cM_{\tilde h^{\hat\C}}(\C)\big)^{-s}\big)u^{s-1}e^{-\mu u} \, \dd u}{	\int_{u>0} \E\big( \big(\cM_{\tilde h^{\hat\C}}(\C)\big)^{-s}\big)u^{s-1}e^{-\mu u} \, \dd u}
	\end{align*}
where the second line follows from \eqref{E:VF} with $h=\tilde{h}^{\hat\C}+\gamma^{-1}\log u - \gamma^{-1}\cM_{\tilde{h}^{\hat\C} + (Q/2)\log \hat g_0}(\C)$. Notice however that $\cM_h(\C)=u$ by definition, so that the above becomes 
\begin{equation*}
	\frac{	\int_{u\in A} \E\big( \big(\cM)_{\tilde{h}^{\hat\C}}(\C)\big)^{-s}\big)u^{s-1}e^{-\mu u} \, \dd u}{	\int_{u>0} \E\big( \big(\cM_{\tilde{h}^{\hat\C}}(\C)\big)^{-s}\big)u^{s-1}e^{-\mu u} \, \dd u} = \frac{	\int_{u\in A} u^{s-1}e^{-\mu u} \, \dd u}{	\int_{u>0} u^{s-1}e^{-\mu u} \, \dd u}
	\end{equation*}
as required.
\end{proof}

\begin{rmk}\label{R:mapconjecture}
This Gamma law is precisely what one expects to get for the limiting distribution of the total area of a random planar map when it is chosen according to an appropriate Boltzmann--Gibbs measure (that is, with a random number of vertices and edges). Let us explain more precisely what we mean by this. In a celebrated work, Tutte showed that the number of planar maps with $n$ edges and $k$ designated roots grows like $C e^{n \beta} n^{- 7/2 + k}$, where $C>0$ and $\beta>0$ are two (essentially) unimportant constants; here $\beta = \log 12$. If we want to embed this conformally using for example circle packing, it is natural to take $k =3$ (so we have $k$ designated points that can be mapped to three fixed locations on the Riemann sphere), so we can rewrite this as $C e^{n\beta } n^{ 1/2 - 1}$. If we assign each map with $n$ edges a weight equal to $e^{ - n  \beta(1+ \eps\mu)}$ 
(this is the slightly subcritical Boltzmann--Gibbs law, which samples very large maps when $\eps$ is small) then we see from Tutte's formula that we should have in the limit as $\eps \to 0$, after suitable rescaling, a distribution for the area which is proportional to $e^{-\mu u} u^{1/2-1}$, that is, a Gamma$(s, \mu)$ law with the parameter $s = 1/2$. This matches the value that one obtains for $\gamma = \sqrt{8/3}$, $k = 3$, $\alpha_ i = \gamma$. Indeed, in that case the formula in Lemma \ref{L:Liouvillegamma} gives
$$
s = \frac{3 \gamma - 2Q}{\gamma} = 3 - (1 + 4/\gamma^2) = 3 - 5/2 =1/2,
$$
as desired. This is no mere coincidence, and analogous results should hold for more general planar maps weighted by the $O(n)$ model or FK weighted planar maps as considered in Chapter \ref{S:maps}. This led \cite{DKRV} to formulate the precise conjecture that after conformally embedding these maps with the three roots sent to some fixed points $\mathbf{z} = (z_1, z_2, z_3)$ of the Riemann sphere, the uniform measure on vertices of the map converges to the Gaussian multiplicative chaos measure associated to the Liouville field $\haz$ with $\alpha = (\gamma , \gamma, \gamma)$. 
\end{rmk}

\ind{Boltzmann--Gibbs distribution}

\ind{FK model}

\ind{Circle packing} 

\ind{Conformal embedding}

\subsection{Unit volume Liouville sphere}\label{SS:LS}

In order to describe the Liouville field $\haz$ defined in the previous subsection, the next step is to identify the law of $\haz$ conditional on the total area. Remarkably, the result does \textbf{not} depend on the actual area except for a (conditionally) deterministic shift corresponding to the area itself. The resulting field will be called the \textbf{unit volume Liouville sphere}. 

\begin{prop}\label{P:unitLiouville} Let $\{\alpha_i,z_i\}$ be as in \cref{D:Liouville}. Then
$\cM_{\haz}(\C)$ and $\haz-\cM_{\haz}(\C)$ are independent, and the law of	$\haz-\cM_{\haz}(\C)$ is equal to that of $$
	\tilde{h}^{\hat\C}- \gamma^{-1}\log(\cM_{\tilde{h}^{\hat\C}}(\C)) $$ 
	weighted by $\big(\cM_{\tilde{h}^{\hat\C}}(\C)\big)^{-s}$, where $\tilde{h}^{\hat\C}= h^{\hat\C,\hat{g}_0}+\tfrac{Q}{2}\log \hat{g}_0+\sum_{i=1}^k\alpha_i G^{\hat\C,\hat g_0}(z_i,\cdot)$.
\end{prop}

\ind{Liouville field! Unit volume Liouville sphere}

\begin{rmk}
Note that this law \emph{does} depend on $(\alpha_i, z_i)_i$ because $s = \gamma^{-1} (\sum_i \alpha_i-2Q )>0$, and because the field  $\tilde{h}^{\hat\C}$ also depends on them. 
%Also, since the law of the Liouville field does not depend on the metric $g$ conformally equivalent to $\hat g$, the law of $\tilde{h}^{\Sp,g}+\tfrac{Q}{2} \log g - \gamma^{-1}\log(\MM_{\tilde{h}^{\Sp,g}}(\C)) $ weighted by $\MM_{\tilde{h}^{\Sp,g}}(\C)^{-s}$ does not depend on the metric. Later we will extend this slightly to more general $g$ for which the associated GFF has zero average on, say, the unit circle. This will be useful when we check equivalence between the unit volume sphere constructed here and that introduced in the work of Duplantier, Miller and Sheffield \cite{DuplantierMillerSheffield} which will be discussed in Chapter \ref{C:MOT}.  
\end{rmk}

\begin{proof}
	Let $F$ be a non-negative Borel measurable function on $H^{-1}(\hat\C)$ and let $A$ be a Borel subset of $[0,\infty)$. Then just as in the proof of Lemma \ref{L:Liouvillegamma} (using the fact that $\cM_h(\C) = u$ when $h = \tilde{h}^{\hat\C}- \gamma^{-1} {\cM_{\tilde{h}^{\hat \C}}(\C)}+\gamma^{-1} \log u$), 
	we have
	\begin{align*}
		& \E_{\az}^L(F(\haz-\gamma^{-1}\log(\cM_{\haz}(\C)))\indic{\cM_{\haz}(\C)\in A}) \\ & = \frac{	\int_{u\in A} \E\big( F(\tilde{h}^{\hat\C} -\gamma^{-1}\log{\cM_{\tilde{h}^{\hat\C}}(\C)}) \big(\cM_{\tilde{h}^{\hat\C}}(\C)\big)^{-s}\big)u^{s-1}e^{-\mu u} \, \dd u}{	\int_{u>0} \E\big( \big(\cM_{\tilde{h}^{\hat\C}}(\C)\big)^{-s}\big)u^{s-1}e^{-\mu u} \, \dd u} \\
		& = \frac{\E\big( F(\tilde{h}^{\hat\C}-{\gamma^{-1}}\log \cM_{\tilde{h}^{\hat\C}}(\C)) \big(\cM_{\tilde{h}^{\hat\C}}(\C)\big)^{-s}  \big) }{\E\big( \big(\cM_{\tilde{h}^{\hat\C}}(\C)\big)^{-s}\big)}\frac{\int_{u\in A} u^{s-1}e^{-\mu u}\dd u}{\int_{u>0} u^{s-1}e^{-\mu u} \dd u} \\
		& = \frac{\E\big( F(\tilde{h}^{\hat\C} -\gamma^{-1}\log \cM_{\tilde{h}^{\hat\C}}(\C)) \big(\cM_{\tilde{h}^{\hat\C}}(\C)\big)^{-s}\big)}{\E\big( \big(\cM_{\tilde{h}^{\hat\C}}(\C)\big)^{-s}\big)} \P_{\az}^L(\cM_{\haz}(\C)\in A)
	\end{align*}
This immediately yields the statement of the proposition.
\end{proof}

\begin{definition}[Unit volume Liouville sphere]\label{D:unitLiouville}
The unit volume Liouville sphere $h_{\az}^{L,1}$ %\indN{{\bf Liouville CFT}! $h_{\az}^{L,1}$; unit volume Liouville sphere on $\hat{\C}$}  
is the random field in $H^{-1}(\hat \C)$ whose law is that of
$$
 \tilde{h}^{\hat\C}-\gamma^{-1} \log(\cM_{\tilde{h}^{\hat\C}}(\C)) \text{ weighted by } (\cM_{\tilde{h}^{\hat\C}}(\C))^{-s},
$$
where $\tilde{h}^{\hat\C}= h^{\hat\C,\hat g_0}+\tfrac{Q}{2}\log \hat{g_0}+\sum_{i=1}^k\alpha_i G^{\hat\C,\hat g_0}(z_i,\cdot)$. 
\end{definition}

\begin{rmk}[Extended Seiberg bounds]\label{R:extendedSB}
		It can be shown that $\mathbb{E}(\cM_{\tilde{h}^{\hat\C}}(\C)^{-s})<\infty$ if (and only if) $\alpha_i<Q$ for all $i$ and $Q-\tfrac12\sum \alpha_i <\tfrac{2}{\gamma} \wedge \min_i (Q-\alpha_i)$. See \cite[Lemma 3.10]{DKRV}. Therefore, a ``unit volume Liouville field'' with the law described in \cref{D:unitLiouville} can still be defined in this extended setting.
	\end{rmk}

Note that this law does not depend on the cosmological constant $\mu$. Furthermore, by \cref{P:unitLiouville}, the law is unaffected if we replace $\hat{g}_0$ with any metric $\hat{g}=e^\rho \hat{g}_0$ as in \eqref{E:metric}.

\begin{rmk}\label{R:dis_LF} Recalling Lemma \ref{L:Liouvillegamma}, we also obtain the following decomposition of the Liouville field $\haz$:
$$
\haz = h_{\az}^{L,1} + \frac1\gamma \log X,
$$
where $h_{\az}^{L,1}$ is the unit volume Liouville sphere, and $X$ is an independent random variable with the Gamma$(s; \mu)$ distribution. 
\end{rmk}

\begin{rmk}[The case $\mu=0$]
	\label{R:mu0}
	Define the infinite measure 
	
\begin{align*} m_\az^{L}(F) & := \lim_{\mu\to 0} \langle V_{\alpha}(\mathbf{z}) F \rangle_{\hat{g}_0} \\
	& = \gamma^{-1} e^{C_\alpha(\mathbf{z})} \prod_i \hat{g}_0(z_i)^{\Delta_{\alpha_i}}  
	\int_{u>0} \E\big( F(\tilde{h}^{\hat \C}  + \tfrac{\log u}{\gamma} -\tfrac{ \log \cM_{\tilde{h}^{\hat \C}}(\C)}{\gamma}) \cM_{\tilde{h}^{\hat \C}}(\C)^{-s}\big) u^{s-1} \, \dd u  .\end{align*}	 
%\indN{{\bf Liouville CFT}! $m_{\az}^L$; infinite Liouville field measure} 
%	& =  e^{C_\alpha(\mathbf{z})} \prod_i \hat{g}_0(z_i)^{\Delta_{\alpha_i}}  \int_\R e^{(\sum\alpha_i-2Q)c}\, \mathbb{E}(F(\tilde{h}^{\hat\C}+c)) \, dc \end{align*}
%(where the last line follows from a change of variables $u=\cM_h e^{\gamma c}$, and is included here for future reference). 
This infinite measure will play a role in the identification of the unit volume Liouville field with the ``unit volume quantum sphere'' which will be introduced in \cref{S:SIsurfaces}.  
 Then we can write
\[ m_\az^{L}(F)=\gamma^{-1} e^{C_\alpha(\mathbf{z})} \prod_i \hat{g}_0(z_i)^{\Delta_{\alpha_i}}  \int_{u>0} u^{s-1} \mathbb{E}(F(h_\az^{L,1} + \gamma^{-1} \log u)) \dd u\]
by \cref{D:unitLiouville} (of $h_{\az}^{L,1}$).
In other words, as in the $\mu>0$ case, we can disintegrate the infinite measure $m_{\az}^{L}$ on $H^{-1}(\hat \C)$ with respect to the total GMC mass of the field. The marginal of the mass is proportional to $u^{s-1} \dd u$ and the law of the field conditioned to have mass $u$ is simply that of the unit volume Liouville sphere $h_{\az}^{L,1}$ plus the constant $\gamma^{-1} \log(u)$. 

\end{rmk}

\begin{rmk}\label{R:3ptexplicit}
	It will be useful later on to express $m_{\alpha,\mathbf{z}}^L$ in terms of the field $h^\fc$ with average zero on the unit circle, defined in \cref{SSS:zerocircleg}, and $\mathbf{z}=(0,\infty,z)$ for $z\ne 0$, $z\in \C$. Namely, using \eqref{E:Fg0isFgc}, we obtain that 
	\begin{align}
	 m_{\alpha,\mathbf{z}}^L(F) & =	 C 
	(|z|\vee 1)^{-4\Delta_{\alpha_3}+\alpha_3(\alpha_1+\alpha_2)}|z|^{-\alpha_1\alpha_3} \times \nonumber \\
	 & \quad  \quad  \int_{u>0}  \mathbb{E}(F(\tilde{h}^\fc+\gamma^{-1}(\log u - \log \cM_{\tilde h^\fc}(\C)))\cM_{\tilde{h}^\fc}(\C)^{-s}) u^{s-1} \dd u  
	 \end{align}
	with $C$ depending only on $\alpha_1,\alpha_2,\alpha_3$ and 
	\[ \tilde{h}^\fc=h^\fc+(\alpha_1+\alpha_2-2Q)\log(|\cdot|\vee 1)-\alpha_1 \log|\cdot| +2\pi \alpha_3 G^\fc(z,\cdot).\]
\end{rmk}

\subsection{Some integrability results}

We have so far discussed the way the Liouville correlation functions evolve under global geometric deformations. However a key step in the development of conformal field theory was accomplished in a celebrated paper of Beliavin, Polyakov and Zamolodchikov \cite{BPZ} in which ``infinitesimal'' geometric deformations were considered and shown to lead to differential identities for the correlation functions. Unlike the identities such as the KPZ identity of Theorem \ref{T:LfieldMobius}, which expresses global invariance of the correlations under M\"obius maps and thus with three degrees of freedom, we get as a result of these considerations an infinite hierarchy of equations (one for each number of insertions, that is, number of points where the correlation functions are being evaluated) and hence infinitely many  degrees of freedom for the parameters of these equations. Comparing to the Virasoro point of view on CFT which was briefly alluded to at the start of this chapter, this is analogous to the fact that the Virasoro algebra contains not only the operators $L_{-1}, L_0, $ and $L_1$ but more generally the infinite dimensional family of operators $\{L_n\}_{n\in \Z}$. 

A first set of identities is obtained by considering the behaviour of the correlation functions $\langle V_{\alpha_1, \ldots, \alpha_k} ( \mathbf{z})\rangle_g $ under infinitesimal deformations of the metric $g \to g + \eps f$. Taking a derivative in the Weyl anomaly formula (Theorem \ref{T:Weyl}) would lead to a field called the \textbf{stress energy tensor} $T(z)$. Correlations between $T(z)$ and the insertion operators are shown to satisfy, as a function of $z$, two families of differential equations known as the \textbf{Ward identities}.
\ind{Stress energy tensor} \ind{Ward identities}
We will not enter into details here except to refer the interested reader to \cite{KRV_BPZ} where this is discussed in detail and furthermore rigorously. Instead we state here the so called \textbf{BPZ equations}.\ind{BPZ equations}

Recall from \eqref{E:conformalweights} that for $\alpha>0$, the conformal weight $\Delta_{\alpha}$ of the operator $V_\alpha$ is given by $\Delta_\alpha=\tfrac{\alpha}{2}(Q-\tfrac{\alpha}{2})$.

\begin{theorem}[Theorem 2.2 in \cite{KRV_BPZ}]
\label{T:BPZ}
Fix $\alpha \in \{ -\tfrac{\gamma}2, -\tfrac{2}{\gamma}\},$ and suppose $k\ge 2$ and $\alpha_1, \ldots, \alpha_k$ satisfy $\sum_{i=1}^k \alpha_i + \alpha > 2Q$, $\alpha_i< Q$ for $1\le i \le k$.  
 Then 
\begin{equation}
\left( 
\frac1{\alpha^2} \partial_z^2 
+ \sum_{i=1}^k \frac{\Delta_{\alpha_i}}{(z-z_i)^2} 
+ \sum_{i=1}^k \frac{1}{z-z_i} \partial_{z_i} 
\right)
\langle V_\alpha(z) \prod_{i=1}^k V_{\alpha_i}(z_i) \rangle_{\hat g} = 0. 
\label{eq:BPZ}
\end{equation}
\end{theorem}

\begin{proof}[Outline of proof]
The proof is very technical and we will only give an extremely rough summary here; of course, readers are once again referred to \cite{KRV_BPZ} for details. 
A key step in the proof of \eqref{eq:BPZ} is the following identity. Write 
$$
G(x;\mathbf{z}) = \langle V_\gamma (x) V_{\alpha_1} (z_1) \ldots V_{\alpha_k} (z_k)\rangle_{\hat g}. 
$$
Note that this is not quite the correlation function appearing in the left hand side of \eqref{eq:BPZ} as here the ``weight'' of the insertion is $\gamma$, whereas in \eqref{eq:BPZ} it is $\alpha \in \{ - \tfrac{2}{\gamma}, - \tfrac{\gamma}2\}$. Write also $G(\mathbf{z})$ for $\langle V_{\alpha_1, \ldots, \alpha_k} ( \mathbf{z})\rangle_{\hat g}$. 

\ind{Integration by parts!Gaussian}
Then using Gaussian integration by parts (already mentioned in the proof of Kahane's convexity inequality in Theorem \ref{T:Kahane}) and plenty of careful estimates (which require ingenious tricks) one can check that for every fixed $1\le i \le k$,
\begin{equation}\label{eq:preBPZ}
\partial_{z_i} G(\mathbf{z}) = -\frac12 \sum_{j\neq i} \frac{\alpha_i \alpha_j}{z-z_j}  G(\mathbf{z})  + \frac{\alpha \mu \gamma}{2} \int_{\C} G(x; \mathbf{z} ) \dd x. 
\end{equation}
This corresponds to (3.27) in \cite{KRV_BPZ}. (A priori it is not even clear that the integral on the right hand side is finite, but this could be deduced with some work from Corollary \ref{C:Vformula}, cf. the proof of Proposition 5.5 in Section 6.8 of \cite{KRV_BPZ}). 

In a second step we may apply the formula with $\mathbf{z}$ replaced by $(z, \mathbf{z})$ and $\alpha_1, \ldots, \alpha_k$ replaced by $\alpha, \alpha_1, \ldots, \alpha_k$. We can then differentiate this identity with respect to $z$ a second time using \eqref{eq:preBPZ} itself to identify the derivatives in the right hand side. After a long calculation and some remarkable cancellations when $\alpha \in \{-\tfrac{2}{\gamma} , \tfrac{\gamma}2\}$, the authors of \cite{KRV_BPZ} end up with the identity \eqref{eq:BPZ}. 
\end{proof}

The proof of the BPZ equations in \cite{KRV_BPZ} is a major step in the proof of the celebrated \textbf{DOZZ formula} (named after Dorn, Otto, Zamolodchikov and Zamolodchikov) which gives an explicit formula for the \textbf{structure constant} $C( \alpha_1, \alpha_2, \alpha_3)$ determining the three point correlation function, and its proof in \cite{DOZZ} is a landmark of Liouville conformal field theory. 
Write 
$$
\ell(z) = \frac{\Gamma(z)}{\Gamma(1-z)},
$$
and furthermore define the special Upsilon function \indN{{\bf Miscellaneous}!$\Upsilon$; special Upsilon function} by 
$$
\log \Upsilon_{\gamma/2} (z) = \int_0^\infty \left( ( \tfrac{Q}2 - z)^2 e^{-t} - \frac{  \sinh^2 (( \tfrac{Q}2 - z) \tfrac{t}{2})}
{\sinh( \tfrac{t\gamma}{4}) \sinh ( \tfrac{t}{\gamma} )} \right)\frac{ \dd t}{t}; \quad 0 < \Re (z) < Q;
$$
which can be analytically continued to $\C$ (this is by no means obvious and in fact follows from functional identities satisfied by the function).

\begin{theorem}
\label{T:DOZZ}
For any $\alpha_1, \alpha_2, \alpha_3$ satisfying the Seiberg bounds, setting $\bar \alpha = \alpha_1 + \ldots + \alpha_3$ 
$$
C_\gamma(\alpha_1, \alpha_2, \alpha_3) = \left( \pi \mu \ell( \tfrac{\gamma^2}{4}) (\tfrac{\gamma}{2})^{2- \gamma^2/2}\right)^{\tfrac{2Q - \bar\alpha}{\gamma}} \!\!
\frac{
\Upsilon'_{\gamma/2}(0) \Upsilon_{\gamma/2} (\alpha_1) \Upsilon_{\gamma/2}(\alpha_2) \Upsilon_{\gamma/2} (\alpha_3) }
{
\Upsilon_{\gamma/2}( \tfrac{\bar \alpha - 2Q}{2} ) \Upsilon_{\gamma/2} ( \tfrac{\bar \alpha}2  - \alpha_1)  \Upsilon_{\gamma/2} ( \tfrac{\bar \alpha}2  - \alpha_2) \Upsilon_{\gamma/2} ( \tfrac{\bar \alpha}2  - \alpha_3)
}
$$
\end{theorem}

Given the BPZ equations, a relatively short sketch of the main arguments can be found in the Section 5 of the lecture notes \cite{Vargas_LN}. 

\medskip The \textbf{conformal bootstrap}, recently proved by Guillarmou, Kupiainen, Rhodes and Vargas \cite{GKRV}, allows one to express correlation functions of order $n+1$ in terms of those of order $n$. In combination with the above DOZZ formula (Theorem \ref{T:DOZZ}), 
this gives exact formulae for correlation functions of \emph{all} orders.

\ind{Conformal bootstrap}

\subsection{Exercises}

\begin{enumerate}[label=\thesection.\arabic*]
	
\item By using spherical polar coordinates, show that the spherical metric $\hat{g}$ satisfies $$\int_\C \hat g_0(z)  \dd z = 4\pi  \text{ and } R_{\hat g_0}(z)=-\tfrac{\Delta \log \hat{g}_0(z)}{\hat{g}_0(z)} \equiv 2.$$

\item \label{ex:shift} Prove \cref{L:hsnorms}, using the fact that for general $g$ as in \eqref{E:metric}, $$\var((h^{\hat \C,g},f)_g)=\|f\|^2_{H^{-1}(\hat \C,g)}$$ for all $f\in H^{-1}(\hat \C,g)$ with $v_g$ average zero. 

\item Let $\{\alpha_i\}_{i=1}^k$ satisfy the Seiberg bounds and $\mathbf{z}=(z_1,\dots, z_k)$ be fixed. Denote $V_\mu:=V_{\alpha_1,\dots, \alpha_k, \mathbf{z}}$ when the cosmological constant is equal to $\mu>0$. Using \cref{C:Vformula}, show that 
$$\langle V_\mu \rangle_{\hat g}= \mu^{\tfrac{2Q-\sum_i \alpha_i}{\gamma}} \langle V_1 \rangle_{\hat g}. $$

\item Give a proof of \eqref{E:VF}.

\item \label{Ex:mobius} %The goal of this exercise is to prove the \textbf{KPZ relations}, \cref{T:LfieldMobius} and \cref{C:LfieldMobius}. Suppose that $m:\hat\C\to \hat \C$ is a M\"{o}bius transformation of the Riemann sphere, $z_1,\cdots, z_k$ are points on $\C$ and $\alpha_1,\cdots, \alpha_k$ satisfy the Seiberg bounds \eqref{E:Seiberg1}, \eqref{E:Seiberg2}.
Suppose that $m:\hat \C \to \hat \C$ is a M\"{o}bius transformation, and
 $$
 m_*\hat{g}_0(z)=\hat{g}_0(m^{-1}(z))|(m^{-1})'(z)|^2,
 $$ 
 that is, viewed as metrics, $m_*\hat g_0$ is the pushforward of $\hat g_0$  by the map $m$.

(a) Recalling the definition $R_{\hat g}=-(1/\hat g)\Delta \log(\hat g)$, show that $R_{m_*\hat{g}_0}\equiv 2$.

(b) Using that $$
G^{\hat{\C},\hat{g}}(x,y)=\frac{1}{2\pi}\left(-\log|x-y|-\frac{1}{2R_{\hat{g}}}\log \hat{g}(x)-\frac{1}{2R_{\hat{g}}}\log \hat{g}(y)+c_{\hat g}\right)
$$
for $\hat{g}$ with constant curvature, and M\"obius invariance of the Green function, 
that is $
G^{\hat\C, m_*\hat{g}_0}(x,y)=G^{\hat \C, \hat{g}_0}(m^{-1}(x),m^{-1}(y))$ for $x\ne y \in \C$,
deduce that $c_{m_*\hat{g}_0}=c_{\hat{g}_0}$. Hint: it may be helpful to write $m(z)$ in the explicit form $(az+b)/(cz+d)$.

(c) Finally, using \cref{R:constantcurvature}, show that 	
$$
\theta_{m_*\hat g_0}=-\frac{1}{2}\bar{v}_{m_*\hat{g}_0}(\log(m_*(\hat g_0)))+\log(2)-\theta_{\hat g}.
$$

\item Write down a general formula for $\langle V\rangle_g$ when $g$ does not have constant scalar curvature but is in the same conformal class as $\hat{g}_0$. Deduce that the law of the Liouville field does \emph{not} depend on the choice of $g$.
\end{enumerate}

%
%\begin{comment}\subsection{Weyl conformal anomaly}
%\label{SS:Weyl}

%\subsection{Probabilistic and geometric interpretation of the Seiberg bounds}

%
%\subsection{To do}
%
%(1) Write the probabilistic interpretation of Seiberg bounds.
%
%(2) Liouville field.
%
%(3) Law of total area. Law of field given total area. (Only explain but not do how it will fit with quantum sphere).
%
%(4) Reflection coeff.

%\input{LCFT_ellen}

\newpage

\section{Gaussian free field with Neumann boundary conditions}\label{S:NGFF}
% !TEX root = master.tex
\def\bbh{\bar{\mathbf{h}}}
%In the next chapter, we will introduce some variants of the planar Gaussian free field which define random surfaces that are scale-invariant in a precise sense: ``zooming in'' near a given point of the surface produces the same surface in distribution.

%We first need to set up the scene correctly. This will require working with a slightly 
So far in this book, we have encountered: 
\begin{itemize}
	\item Gaussian free fields on graphs (Section \ref{SS:DGFF});
	\item Gaussian free fields with zero boundary conditions on proper, regular domains of $\R^d$ (Chapter \ref{S:GFF}); and
	\item Gaussian free fields on compact surfaces (Section \ref{S:ZAGFF}).
\end{itemize}
 The purpose of this chapter is to introduce a different version of the GFF on simply connected domains of $\mathbb{C}$, but now with non-zero boundary conditions. This object will be the so called \textbf{Neumann} or \textbf{free boundary} GFF. It is the basic building block for constructing the special ``scale invariant quantum surfaces'' that will be the focus of Chapter \ref{S:SIsurfaces}.
 
 In general if we wish to add boundary data to a GFF it is natural to simply add a function that is harmonic in the domain (though it can have relatively wild behaviour on the boundary).
We will seek to impose \textbf{Neumann boundary conditions}. Recall that for a smooth function, this means that the normal derivative of the function vanishes along the boundary (if the domain is smooth).
Of course for an object as rough as the GFF it is a priori unclear what this condition should mean. Indeed, we will see that the resulting object is actually the same as when we don't impose any conditions at all (which is why the field can also be called a free boundary GFF, as is done for example in the papers \cite{zipper} and  \cite{DuplantierMillerSheffield}). Indeed, note that in the discrete, a random walk on a graph with Neumann/``reflecting'' boundary conditions or no/``free'' boundary conditions are by definition the same thing (and both converge to reflecting Brownian motion, whose generator is $\tfrac12$ the Laplace operator with Neumann boundary conditions).

\paragraph{Outlook} Let $D$ be a proper, simply connected domain of $\C$.
We will first show how to define the Neumann GFF as a random distribution on $D$, just as in  \cref{ZGFFd} for the Dirichlet GFF and \cref{S:ZAGFF} for the GFF on a Riemann surface. This allows for a straightforward deduction of several nice properties, which is why we present this point of view first. 
In  \cref{NGFF:sp} we will then go on to show that the Neumann GFF can be defined as a stochastic process (as in the Dirichlet case), and that this object coincides with the random distribution defined here when its index set is restricted appropriately. In the penultimate section of this chapter we will discuss some further variants of the Gaussian free field, and how they relate to one another, and conclude in the final section with an analysis of \emph{boundary} Gaussian multiplicative chaos.

\paragraph{Warning} One technical complication when working with the Neumann GFF, {compared to the Dirichlet case}, is that it is really only defined up to a global additive constant. This corresponds to the fact that if one tries to extend the Dirichlet inner product $(\cdot, \cdot)_\nabla$ to test functions that are not necessarily compactly supported in $D$, it is no longer an inner product. Indeed, functions that are constant on the domain will have zero Dirichlet norm. Alternatively (as we will see later) one can think of the additive constant as arising from the fact that the Green function with Neumann boundary conditions is \emph{not} canonically defined (or equivalently, that Brownian motion reflected on the boundary of $D$ is recurrent). 

Note that this complication was already present for the GFF on a Riemann surface, see \cref{S:ZAGFF}. In that setting we fixed the additive constant by requiring the field to have zero average with respect to the Riemannian volume form.
In this chapter it will be useful to have access to both of the following viewpoints.
\begin{enumerate}
	\item We can view the Neumann GFF as a \textbf{distribution modulo constants} (two distributions are equivalent if their difference is a constant function). Equivalently, a distribution modulo constants can be defined as a continuous linear functional on test functions whose integral is required to be zero. 
	\item We can specify a 
\textbf{particular representative} of the Neumann GFF's 
	\textbf{equivalence class modulo constants} (for example by requiring that the average of the field over a specific region is zero). We will then speak of ``fixing the additive constant''. Note that while this point of view may appear to be more concrete, fixing the additive constant for the free field in this way actually causes it to lose some useful properties, such as conformal invariance.
\end{enumerate} {When using the Neumann GFF, we will therefore always need to be careful to say whether we consider the modulo constants version, or a version that has had the constant fixed in a particular way.} 

\subsection{The Neumann GFF as a random distribution}\label{sec:NGFF_dist}

Let $\bar{\cD}(D)$ be the space of $C^\infty$ functions in ${D}$ with $(f,f)_\nabla<\infty$ (``finite Dirichlet energy'')\indN{{\bf Function spaces}! $\bar{\cD}(D)$; smooth functions in $D$ with finite Dirichlet energy, considered modulo constants}, defined \textbf{modulo constants}. That is, two functions are equivalent if their difference is a constant function. Note that these functions are 
\emph{not} assumed to have compact support in $D$. It is clear that on this space, $(\cdot, \cdot)_\nabla$ really is an inner product.
Hence we can define $\bar{H}^1(D)$ to be the Hilbert space closure of $\bar{\cD}(D)$ with respect to $(\cdot, \cdot)_\nabla$.\indN{{\bf Function spaces}! $\bar{H}^1(D)$; Hilbert space closure of $\bar{\cD}(D)$ with respect to $(\cdot, \cdot)_\nabla$}

We define a distribution modulo constants to be a continuous linear functional on the space of test functions $f\in \mathcal{D}_0(D)$ such that $\int_D f(x)\, \dd x =0$, and denote the set of such test functions by $\tilde \cD_0(D)$. \indN{{\bf Function spaces}! $\tilde{\mathcal{D}}_0(D)$; test functions $f\in \mathcal{D}_0(D)$ with total integral zero} 
We write $\bar{\cD}_0'(D)$ for the space of distributions modulo constants, \indN{{\bf Function spaces}! $\bar{\cD}_0'(D)$; distributions modulo constants on $D$} and equip it with the topology of weak-$\star$ convergence. That is, a sequence $T_n$ of distributions modulo constants converges to a distribution $T$ if and only if $(T_n,f)\to (T,f)$ for any test function $f\in \tilde \cD_0(D)$.

\begin{rmk}[Notation]
	In this section we will use the general notation $\,\bar{\cdot}\,$ to refer to spaces of objects or objects defined modulo constants, and the notation $\,\tilde{\cdot}\,$ for spaces of objects or objects with zero average over $D$.
\end{rmk}

As in Section \ref{SS:randomdistributions}, a random variable $X$ defined on a probability space $(\Omega, \cF, \P)$ and taking values in the space of distributions modulo constants, is simply a function $X: \Omega \to \bar{\cD}'_0(D)$ which is measurable with respect to the Borel $\sigma$-field on $\bar\cD'_0(D)$ induced by the weak-$*$ topology. Arguing as in Lemma \ref{L:measurability}, we see that convergence of a sequence of random variables $X_n \in \bar{\cD}_0(D)$ is a measurable event. Thus, it makes sense to ask about almost sure convergence of such sequences.

\ind{GFF!Neumann, random series}
We now give the definition of the Neumann GFF as a random element of $\bar{\cD}_0'(D)$.

\begin{theorem}
	\label{TD:NeumannGFF2}
	Let $\{\bar{f}_j\}_{j\ge 1}$ be \emph{any} orthonormal basis of $\bar H^1(D)$, and $\{X_j\}_{j\ge 1}$ be a sequence of independent $\mathcal{N}(0,1)$ random variables. Then the random series
	\begin{equation}\label{NeumannGFF2}
		\bar{\bf h}_n:= \sum_1^n X_j \bar f_j
	\end{equation}
	converges almost surely in the space of distributions modulo constants.
	Moreover, the law of the limit $\bar{\mathbf{h}}=\bar{\mathbf{h}}^D$ does not depend on the choice of orthonormal basis $\{\bar{f}_j\}_j$, and can be written as the sum of a Dirichlet boundary condition GFF on $D$ and an independent harmonic function modulo constants. 
\end{theorem}

\begin{definition}[Neumann GFF as a distribution modulo constants]
	\label{D:NeumannGFF2}
	We define the Neumann GFF $\bar{\mathbf{h}}$ to be the random distribution modulo constants constructed in Theorem \ref{TD:NeumannGFF2}. 
\end{definition}

\begin{rmk} (\emph{Neumann boundary conditions})  Suppose that $D=\D$. In defining $\bar H^1(\D)$ we started from the space $\bar{\cD}(\D)$ of smooth functions (modulo constants) on $\D$ with no restriction on their boundary conditions. However, we could equally have started with the space of smooth functions (modulo constants) with Neumann boundary conditions, and ended up with the same space $\bar H^1(\D)$ after taking the closure with respect to $(\cdot, \cdot)_\nabla$. Indeed, there exists an orthonormal basis of $L^2(\D)$ made up of eigenfunctions of the Laplacian with Neumann boundary conditions (see for example \cite[Theorem 8.5.2]{Jost}). Then omitting the first eigenfunction (which has eigenvalue 0) and dividing the rest by the square roots of their respective eigenvalues and considering them modulo constants, provides an orthonormal basis of $\bar{H}^1(\D)$. 		
	Thus, one can think of the Neumann GFF as either having no imposed (``free'') boundary conditions, or as having Neumann boundary conditions.
	
	The connection with Neumann boundary conditions will also become more apparent when we define the Neumann GFF as a stochastic process. Indeed, we will see that its covariance function is given by a Green function in the domain, with Neumann instead of Dirichlet boundary conditions. As already mentioned,  in the discrete, a random walk on a graph with Neumann/``reflecting'' boundary conditions or no/``free'' boundary conditions are really one and the same thing. So the discrete Green's function will be the same if either free or Neumann boundary conditions are imposed. \ind{GFF! Free boundary conditions}
\end{rmk}

\begin{comment}\begin{rmk}
	The fact that $\bar{\mathbf{h}}$ is defined as a random distribution modulo constants means that we can test it against any smooth \emph{compactly supported} function in $D$ with average zero. But we will want to consider, in an appropriate sense, the ``behaviour of the field on the boundary''. In principle, see for example \cite{adams}, there is not a unique way to extend a distribution modulo constants $\bar{h}$ to a continuous linear functional on the space $\tilde{H}^1(D)=\{f\in H^1(D): \int_D f(x) \, dx =0\}$ (say), which includes functions with \emph{non-zero} boundary conditions. However, we will see later that for any given function $f\in \tilde{H}^1(D)$ (in fact, for $f$ in an even larger class of objects) the sum \eqref{NeumannGFF2} converges almost surely when tested against $f$, and the law of the limit is independent of the choice of orthonormal basis. 
\end{rmk}
\end{comment}

\begin{proof}[Proof of \cref{TD:NeumannGFF2}] We will carry out the proof %of the claim in \cref{D:NeumannGFF2} 
	in two steps: first assuming that $D$ is the unit disc $\D$; and then extending to general $D$ by conformal invariance.
	
	\vspace{0.2cm}
	\textbf{Step 1 $(D=\mathbb{D}$).} Write $\overline{\Harm}(\D)$ \indN{{\bf Function spaces}! $\overline{\Harm}(D)$; harmonic functions on $D$ with finite Dirichlet energy, viewed modulo constants}  for the space of harmonic functions on $\D$ with finite Dirichlet energy, viewed modulo constants. By the same reasoning as in \cref{lem:markovprop_od}, we can decompose
	\begin{equation}\label{eqn:orthogH1}\bar H^1(\D)=H^1_0(\D)\oplus \overline{\Harm}(\D)\end{equation}
	as a direct orthogonal sum with respect to the Dirichlet inner product.
	
	We can therefore define $f_j^0$ and $\bar{f}_j^H$ to be the projections onto $H_0^1(\D)$ and $\overline{\Harm}(\D)$ respectively, of each $\bar{f}_j$ in our orthonormal basis of $\bar H^1(\D)$. Accordingly, we set $\mathbf{h}_n^0:= \sum_{j=1}^n X_j f_j^0$ and $\bbh_n^H=\sum_{j=1}^n X_j \bar{f}_j^H$, so that $\bbh_n=\mathbf{h}_n^0+\bbh_n^H$ for each $n$. %Note that since $(f_j^0,\bar{f}_j^H)_\nabla=0$, we have that $(f_j^0,f_j^0)_\nabla \le 1$, $(\bar{f}_j^H,\bar{f}_j^H)_\nabla\le 1$ for each $j\ge 1$. 
	
	First, we claim that $\mathbf{h}_n^0$ converges almost surely in the space $H_0^s(\D)$ (for any $s<0$) to a limit $\mathbf{h}$ with the law of a zero boundary condition GFF in $\D$. The proof is very similar to that of \cref{T:GFFseriesSob}, but we reproduce it here, since there are some technical differences. We let $(e_m)_{m\ge 1}$ be an orthonormal basis of $L^2(\D)$ made up of eigenfunctions of $-\Delta$, with corresponding eigenvalues $(\lambda_m)_{m\ge 0}$, and recall that $\lambda_m \asymp m$ as $m\to \infty$ by Weyl's law. Then we have, for $n\ge 1$, 
	\begin{equation}\label{eq:h0cauchy} \mathbb{E}\big((\mathbf{h}_n^0,\mathbf{h}_n^0)_{H_0^s}\big) = \sum_{j=1}^n (f_j^0,f_j^0)_{H_0^s} \end{equation} where, using the Gauss--Green formula, the fact that $e_m\in H_0^1(\D)$ for each $m$, and Fubini:
	\begin{align*}
		\sum_{j\ge 1} (f_j^0,f_j^0)_{H_0^s} &  = \sum_{j\ge 1} \sum_{m\ge 1} \lambda_m^s (f_j^0,e_m)^2_{L^2}  = \sum_{j\ge 1} \sum_{m\ge 1} \lambda_m^{-1+s}(f_j^0,\tfrac{e_m}{\sqrt{\lambda_m}})^2_{\nabla} 
			 = \sum_{j\ge 1} \sum_{m\ge 1} \lambda_m^{-1+s}(\bar{f}_j,\tfrac{e_m}{\sqrt{\lambda_m}})^2_{\nabla} 
				 \\ &  = \sum_{m\ge 1} \lambda_m^{-1+s} \sum_{j\ge 1}(\bar{f}_j,\tfrac{e_m}{\sqrt{\lambda_m}})^2_{\nabla}  = \sum_{m\ge 1} \lambda_m^{-1+s}  <\infty;
		\end{align*}
the finiteness holding as long as  $s<0$. We deduce that $\mathbf{h}_n^0$ converges almost surely to a limit $\mathbf{h}^0$ in $H_0^s$, by exactly the same reasoning as in the proof of \cref{T:GFFseriesSob}. 
Moreover, as a limit of centered and jointly Gaussian random variables, $((\mathbf{h}^0,e_j))_{j\ge 1}$ are centered and jointly Gaussian, with \begin{align}\label{eq:FSlawunique} \mathbb{E}((\mathbf{h}^0,e_j) (\mathbf{h}^0,e_k))=\lim_{N\to \infty} \sum_{n=1}^N(f_n^0,e_j)(f_n^0,e_k) & = \lim_{N\to \infty} \sum_{n=1}^N (\bar{f}_n,\lambda_j^{-1}e_j)_\nabla (\bar{f}_n, \lambda_k^{-1} e_k)_\nabla \nonumber \\ & = (\lambda_j^{-1}e_j, \lambda_k^{-1}e_k)_{\nabla } =\lambda_j^{-1} \mathbf{1}_{j=k}\end{align} 
for each $j,k\ge 1$. This implies that $\mathbf{h}^0$ is equal in law (as a random element of $H_0^s$ and therefore as a distribution)  to a zero boundary GFF in $\D$ (indeed, it is immediate from the definition as a Fourier series that the above holds for such a GFF).

Next, we will show that $\bbh^H_n = \sum_{j=1}^n X_j \bar{f}_j^H$ converges almost surely in $\bar\cD_0'(\D)$, to a random element $\bbh^H$ of $\overline{\Harm}(\D)$. 
For this we will make use of a specific orthonormal basis of $\overline{\Harm}(\D)$, which will play a similar role to the basis $(e_m)_m$ of eigenfunctions of the Laplacian used above. This basis is given by 
\[ \bar{u}_j(z) = \frac{1}{\sqrt{\pi j}}\Re(z^j) \text{ and } \bar{v}_j(z)=\frac{1}{\sqrt{\pi j}}\Re(iz^j) \text{ for }  j\ge 1\]
(viewed modulo constants), which are easily checked to be orthonormal with respect to $(\cdot, \cdot)_\nabla$. They also span the space $\overline{\Harm}(\D)$, because any harmonic function on $\D$ is the real part of an analytic function, and therefore admits a Taylor series expansion of the form $a+\sum_{j=1}^\infty b_j \Re(z^j) + \sum_{j=1}^\infty c_j \Re(i z^j)$ with $a,\{b_j,c_j\}_j \in \R$.

For each $j\ge 1$, let us denote by $f_j^H$ the representative of $\bar{f}_j^H$ with $f_j^H(0)=0$. Similarly, we set $u_m=(\pi m)^{-1/2}\Re(z^m), v_m=(\pi m)^{-1/2} \Re (iz^m)$ for $m\ge 1$. Another simple calculation verifies that $((u_m,v_m))_{m\ge 1}$ are orthogonal with respect to the $L^2$ inner product on $\D$, and $(u_m,,u_m)_{L^2}=(v_m,v_m)_{L^2}=1/(2m(m+1))$ for each $m\ge 1$. We write 
	\begin{equation}\label{eq:hHcauchy} \mathbb{E}\big( \|\sum_{j=1}^n X_j f_j^H \|_{L^2(\D)^2}\big) = \sum_{j=1}^n (f_j^H, f_j^H)_{L^2(\D)} \end{equation} and by Parseval's identity, can express
\begin{equation}\label{eq:hHcauchy2}
	\sum_{j\ge 1} (f_j^H,f_j^H)_{L^2(\D)}  = \sum_{j\ge 1} \sum_{m\ge 1} 2m(m+1)\big((f_j^H,u_m)^2_{L^2}+(f_j^H,v_m)^2_{L^2}\big).\end{equation}
Now, for each $j\ge 1$, since $\bar{f}_j^H\in \overline{\Harm}(\H)$, it has an expansion $\bar{f}_j^\H=\sum_{m\ge 1} (\bar{f}^H_j,\bar{u}_m)_\nabla \bar u_m + (\bar{f}_j^H, \bar{v}_m)_\nabla \bar v_m$ which converges in $\bar{H}^1(\D)$, and therefore (since the $(u_m,v_m)$ are also orthogonal for the $L^2$ inner product and since $(\bar{f}_j^0,\bar{u}_m)_\nabla=0$) $$(f_j^H,u_m)_{L^2}=(\bar{f}_j^H,\bar{u}_m)_\nabla(u_m,u_m)_{L^2}=\frac{1}{2m(m+1)}(\bar{f}_j^H,\bar{u}_m)_{\nabla}=\frac{1}{2m(m+1)}(\bar{f}_j,\bar{u}_m)_{\nabla}$$ for each $j,m\ge 1$. The analogous equation holds for $v_m$. Thus the right hand side of \eqref{eq:hHcauchy2} becomes
\begin{align*}\sum_{j\ge 1} \sum_{m\ge 1} \frac{1}{2m(m+1)}\big((\bar f_j,\bar u_m)^2_{\nabla}+(\bar f_j,\bar v_m)^2_{\nabla}\big) & =\sum_{m\ge 1} \frac{1}{2m(m+1)}\sum_{j\ge 1} \big((\bar f_j,\bar u_m)^2_{\nabla}+(\bar f_j,\bar v_m)^2_{\nabla}\big) \\
	& =\sum_{m\ge 1} \frac{1}{2m(m+1)} (\|\bar u_m\|_\nabla^2+\|\bar v_m\|_\nabla^2) <\infty\end{align*}
where we applied Fubini, and the fact that $(\bar f_j)_{j\ge 1}$ are an orthonormal basis of $\bar{H}^1(\D)$. By the same argument used in the case of $(\mathbf{h}_n^0)_n$, this implies that $\sum_j X_j f_j^H$ converges almost surely in $L^2(\D)$, and in particular, $\bbh_n^H= \sum_{j=1}^n X_j \bar{f}_j^H$ converges almost surely in $\bar\cD_0'(\D)$. 

By analogous reasoning to \eqref{eq:FSlawunique}, it holds that the almost sure $L^2(\D)$ limit of $\sum_j X_j\bar{f}^H_j $ has to be  independent of the choice of $\{\bar f_j\}_j$. It remains to justify that the almost sure limit of
\begin{equation}\label{eq:harmonic_lim}\sum_1^n \sqrt{\frac{1}{\pi j}}\Re((\alpha_j + i\beta_j) z^j) \; ; \; \alpha_j, \beta_j \overset{\text{i.i.d}}{\sim} \mathcal{N}(0,1)\end{equation}
is harmonic. This simply follows from the fact limits in $\cD'(\mathbb{D})$ of harmonic functions are harmonic (by definition, distributional limits of weakly harmonic functions are weakly harmonic, and then true harmonicity follows by elliptic regularity). %that $\sum_j \sqrt{1/\pi j}(\alpha_j + i\beta_j)z^j$ almost surely defines a Taylor series with radius of convergence one (since, for example, $\limsup_n (\log n)^{-1}\max_{1\le j \le n}|\alpha_j+i\beta_j|\le 1$ almost surely). Thus for any $r<1$, the limit of \eqref{eq:harmonic_lim} is the real part of an analytic function in $B(0,r)$, and so the limit is harmonic in $\cup_{r\in (0,1)} B(0,r)=\D$.

Finally, we claim that $\bbh^H=\lim_{n\to \infty} \bbh^H_n$ and $\mathbf{h}^0=\lim_{n\to \infty} \mathbf{h}_n^0$ are independent. In other words, that for any $\rho,\eta\in \tilde\cD_0(\D)$, $(\mathbf{h}^0,\rho)$ and $(\bbh^H,\eta)$ are independent. Since $(\mathbf{h}^0,\rho),(\bbh^H,\eta)$ are the almost sure limits of $\sum_{j=1}^n (X_j f_j^0,\rho)$ and $\sum_{j=1}^n (X_j \bar{f}_j^H, \eta)$ respectively, they are centered and jointly Gaussian. Hence it suffices to check that \[ \lim_{n\to \infty} \mathbb{E}[(\sum_{j=1}^n X_j f_j^0,\rho)(\sum_{k=1}^n X_k \bar{f}_k^H, \eta)]=\lim_{n\to \infty} \sum_{j=1}^n (f_j^0,\rho)(\bar{f}_j^H,\eta) = 0.\] 
For this, recall the definitions of $(e_m,\bar{u}_m,\bar{v}_m)_{m\ge 1}$ appearing previously in the proof. Then for each $j$ we can write 
$$
f_j^0=\sum_{m\ge 1} (f_j^0,(\lambda_m)^{-1/2}e_m)_{\nabla}(\lambda_m)^{-1/2}e_m, \text{ and } \bar{f}_j^H=\sum_{m\ge 1}( (\bar{f}_j^H,\bar{u}_m)_\nabla \bar{u}_m+(\bar{f}_j^H, \bar{v}_m)_\nabla
\bar{v}_m,
$$
  with both sums converging in $\bar{H}^1(\D)$. By Parseval, we can therefore write
\begin{align*}  \sum_{j\ge 1} (f_j^0,\rho)(\bar{f}_j^H,\eta) & = \sum_{j\ge 1} \sum_{m,n\ge 1} \frac{(f_j^0,e_m)_\nabla}{\sqrt{\lambda_m}} (\bar{f}_j^H,\bar{u}_n)_\nabla (e_m,\rho)(\bar{u}_n,\eta) \\
& \quad \quad + \frac{(f_j^0,e_m)_\nabla}{\sqrt{\lambda_m}}(\bar{f}_j^H,\bar{v}_n)_\nabla (f_j^0,\rho)(\bar{v}_m,\eta) \\
	& = \sum_{m,n\ge 1}(e_m,\rho)(\bar{u}_n,\eta)\sum_{j\ge 1} \frac{(f_j^0,e_m)_\nabla}{\sqrt{\lambda_m}} (\bar{f}_j^H,\bar{u}_n)_\nabla  \\
	& \quad \quad  + \sum_{m,n\ge 1}(e_m,\rho)(\bar{v}_n,\eta)\sum_{j\ge 1} \frac{(f_j^0,e_m)_\nabla}{\sqrt{\lambda_m}} (\bar{f}_j^H,\bar{v}_n)_\nabla  
	\end{align*}
where in the second line we also applied Fubini. Note that this is justified since (restricting to the first of the two double sums by symmetry)
\begin{multline*} \sum_{m,n\ge 1}\sum_{j\ge 1} |(e_m,\rho)(\bar{u}_n,\eta)\frac{(f_j^0,e_m)_\nabla}{\sqrt{\lambda_m}} (\bar{f}_j^H,\bar{u}_n)_\nabla| \\
\le \sum_{m,n\ge 1} |(e_m,\rho)(\bar{u}_n,\eta)| \sum_{j\ge 1} \frac{(f_j^0,e_m)^2_\nabla}{{\lambda_m}} \sum_{j\ge 1} (\bar{f}_j^H,\bar{u}_n)_\nabla^2 
\end{multline*}
where the two sums over $j$ are bounded by $(e_m,e_m)_\nabla/\lambda_m = 1$ and $(\bar{u}_n,\bar{u}_n)_\nabla=1$ (using Parseval) and $\sum_{m,n\ge 1} |(e_m,\rho)(\bar{u}_n,\eta)|<\infty$ since $\rho,\eta\in \cD_0(\D)$.

Noticing that 
$$(\bar{f}_j,\bar{u}_n)_\nabla=(\bar{f}_j^H,\bar{u}_n)_\nabla, (\bar{f}_j,\bar{v}_n)=(\bar{f}_j^H,\bar{v}_n) \text{ and }(\bar{f}_j, \lambda_m^{-1/2}e_m)_\nabla=(f^0_j, \lambda_m^{-1/2}e_m)_\nabla$$ 
for each $j,m,n$ by orthogonality of $H_0^1(\D)$ and $\overline{\mathrm{Harm}}(\D)$, we conclude that
	\begin{align*} &\sum_{j\ge 1} (f_j^0,\rho)(\bar{f}_j^H,\eta)  \\
		 & = \sum_{m,n\ge 1}(e_m,\rho)(\bar{u}_n,\eta)\sum_{j\ge 1} (\bar{f}_j,\frac{e_m}{\sqrt{\lambda_m}})_\nabla (\bar{f}_j,\bar{u}_n)_\nabla  + \sum_{m,n\ge 1}(e_m,\rho)(\bar{v}_n,\eta)\sum_{j\ge 1} (\bar f_j,\frac{e_m}{{\sqrt{\lambda_m}}})_\nabla (\bar{f}_j,\bar{v}_n)_\nabla  \\
	 & = \sum_{m,n\ge 1}(e_m,\rho)(\bar{u}_n,\eta)(e_m,\bar{u}_n)_\nabla  + \sum_{m,n\ge 1}(e_m,\rho)(\bar{v}_n,\eta)(e_m,\bar{v}_n)_\nabla  = 0
\end{align*}
as required.

	\medskip
	\textbf{Step 2 (general $D$).} Suppose that $D\subsetneq \C$ is simply connected, and let $\{\bar f_j\}_j$ be an orthonormal basis for $\bar H^1(D)$. We would like to show that $\bar{\bf h}_n=\sum_1^n X_j f_j$ converges almost surely in $\bar{\cD}_0'(D)$ (when the $X_j$ are i.i.d.\ $\cN(0,1)$).
	
	For this, we are going to use Step 1 and conformal invariance. Let $T:\D\to D$ be a conformal isomorphism (which exists by the Riemann mapping theorem). Then by conformal invariance of the Dirichlet inner product, $\{\bar f_j \circ T\}_j$ forms an orthonormal basis of $\bar H^1(\D)$. We therefore know, by Step 1, that $\bar{\bf h}_n \circ T := \sum_1^n X_j  (\bar f_j \circ T)$ converges almost surely in $\bar\cD'_0(\D)$. That is, with probability one, there exists $\bar{\bf h}^\D\in \bar \cD_0'(\D)$ such that $(\bar{\bf h}_n \circ T, g)\to (\bar{\bf h}^\D, g)$ as $n\to \infty$ for all $g\in \tilde \cD_0(\D)$.
	
	Since for $f\in \tilde{\cD}_0(D)$ the function $g(z)=|T'(z)|^2 (f\circ T)(z)$ is in $\tilde \cD_0(\D)$, this tells us -- in particular -- that with probability one:
	$$ (\bar{\bf h}_n \circ T, |T'|^2(f\circ T))\to (\bar{\bf h}^\D, |T'|^2 (f\circ T)) \text{ as } n\to \infty, \; \forall f\in \tilde{\cD}_0(D).$$ Defining $\bar{\bf h}\in \bar{\cD}_0'(D)$ by $(\bar{\bf h}, f)=(\bar{\bf h}^\D, |T'|^2 (f\circ T))$ for all $f\in \bar \cD_0(D)$, this is exactly saying that with probability one,
	$(\bar{\bf h}_n, f)\to (\bar{\bf h}, f) \text{ as } n\to \infty$ for all $f\in \tilde \cD_0(D)$. That is, $\bar{\bf h}_n \to \bar{\bf h}$ in $\bar\cD_0'(D)$, almost surely as $n\to \infty$.
	
	Finally, by the same argument, if $T:\D\to D$ is conformal then the law of $\bar{\bf h}$ must be given by the law of $\bar{\bf h}^\D \circ T^{-1}$, where $\bar{\bf h}^\D$ is the (unique in law) limit of \eqref{D:NeumannGFF2} when $D=\D$. Note that this does not depend on the choice of $T$, since the law of $\bar{\bf h}^\D$ is conformally invariant (we can see this by applying the reasoning of the previous sentence with $D=\D$, together with the uniqueness in Step 1). Thus, the law of $\bar{\bf h}$ is unique for general $D$.
	
	Using the description of this law when $D=\D$ from Step 1, plus conformal invariance of the Dirichlet GFF (\cref{thm:CIDGFF}) and the fact that conformal isomorphisms preserve harmonicity, we see that in general the law of $\bar{\bf h}$ satisfies the description in \cref{D:NeumannGFF2}.
\end{proof}

By conformal invariance of the Dirichlet inner product, we obtain the following (the details are spelled out in the proof above):

\begin{cor}\label{C:NGFFci}
	Let $\bar{\bf h}^D$ be the Neumann GFF (viewed modulo constants) in $D$, as in \cref{D:NeumannGFF2}. Then the law of $\bar{\bf h}^D$ is conformally invariant. That is, if $T:D\to D'$ is a conformal isomorphism between simply connected domains, then
	$$\bar{\bf h}^{D'} \overset{(d)}{=} \bar{\bf h}^{D}\circ T^{-1}$$
	where $(\bar{\bf h}^{D}\circ T^{-1},f):=(\bar{\bf h}^D, |T'|^2 (f\circ T))$ for all $f\in \tilde{\cD}_0(D')$.
\end{cor} 

Straight from the definition, we also know that if $\bar{\bf h}$ is the Neumann GFF (viewed as a distribution modulo constants) in $D$, then $\bar{\bf h}$ can be written as the sum ${\bf h}_0+u$, where ${\bf h}_0$ has the law of a zero boundary GFF in $D$, and $u$ is an independent harmonic function modulo constants in $D$. By applying the Markov property of the Dirichlet GFF (\cref{T:mp}) to $\bar{\bf h}$ we get an analogous decomposition for the Neumann GFF.

\ind{Markov property!Neumann GFF}

\begin{theorem}[Markov property]\label{T:Nmp} Fix $U \subset D$, open. Let $\bar{\bf h}$ be a Neumann GFF viewed as a distribution modulo constants in $D$, as in \cref{D:NeumannGFF2}. Then we may write
	$$
	\bar{\bf h} = {\bf h}_0 + \varphi
	$$
	where:
	\begin{enumerate}
		\item ${\bf h}_0$ is a zero boundary condition GFF in $U$, and is zero outside of $U$;
		
		\item $\varphi$ is a harmonic function viewed modulo constants in $U$;
		
		\item ${\bf h}_0$ and $\varphi$ are independent.
	\end{enumerate}
\end{theorem}

For a more explicit Markov decomposition in the case $D=\H$ and $U$ a semidisc centered on the real line, see \cref{P:MP2}.

\medskip
Recall that we defined a distribution modulo constants to be a continuous linear functional on the space $\tilde{\cD}_0(D)$ of test functions with average $0$. Equivalently, we could define it to be an equivalence class of distributions (elements of $\cD_0'(D)$), under the equivalence relation identifying distributions $\phi_1$ and $\phi_2$ whenever $\phi_1-\phi_2\equiv C$ for $C\in \R$.

\begin{rmk}[Fixing the additive constant, see also \cref{D:NGFFnorm}]
	\label{R:NGFF_fixed}
	With the latter perspective, it is quite natural (and will sometimes be useful) to fix the additive constant for the GFF in some way (that is, to pick an equivalence class representative). For example, we could define the Neumann GFF $\mathbf{h}$ with average zero when tested against some fixed test function $\rho_0\in \cD_0(D)$, by setting
	$$({\bf h},\rho)=(\bar{\bf h}, \rho -\frac{\int_{\bar D} \rho(\dd x)}{\int_{\bar D}\rho_0(\dd x)}\rho_0) \quad \text{ for } \rho\in \mathcal{D}_0(D),$$
	where $\bar{\bf h}$ is as in \cref{D:NeumannGFF2}.
Since $\bar{\bf h}$ is almost surely a random distribution modulo constants, the above can be defined simultaneously for all $\rho\in \cD_0(D)$, and almost surely defines an element of $\cD_0'(D)$, that is, a distribution on $D$. In fact, by \cref{C:H1loc}, it almost surely defines an element of $H^{-1}_{\mathrm{loc}}(D)$: the \textbf{local Sobolev space} of distributions whose restriction to any $U\Subset D$ \indN{{\bf Miscellaneous}! $A\Subset B$; closure of $A$ is a subset of $B$}(that is such that $\bar{U}$ is a compact subset of $D$) is an element of $H^{-1}_0(U)$.
	
	\noindent Note that the choice of constant, or equivalence class representative, changes the resulting element of $\cD_0'(D)$, but not how it acts when tested against functions (with average zero) in  $\tilde \cD_0(D)$.
\end{rmk}

\begin{rmk}
	Although it is sometimes helpful to fix the additive constant for the Neumann GFF, one should take care with the conformal invariance and Markovian properties discussed above. In particular:
	\begin{itemize}
		\item if $\bar{\bf h}$ is a Neumann GFF in $D$ with additive constant fixed in some way, then it is \emph{no longer} conformally invariant;
		\item in this case one can still write $\bar{\bf h}={\bf h}_0$$+u$ with ${\bf h}_0$ a Dirichlet GFF in $D$ and $u$ a harmonic function, but ${\bf h}$ and $u$ \emph{need not be independent};
		\item on the other hand, if one starts with a Neumann GFF modulo constants, decomposes it as a Dirichlet GFF plus a harmonic function modulo constants, and then fixes the constant for the GFF in a way that only depends on the harmonic function (for example, by specifying the value of the harmonic function at a point), then the two summands \emph{will be} independent.
	\end{itemize}	
\end{rmk}

\subsection{Covariance formula: the Neumann Green function}

Recalling the definition of the Dirichlet GFF in a domain $D$, it is quite natural to guess that the Neumann GFF will have ``covariance'' given by a version of the Green function with Neumann boundary conditions in $D$. This is indeed the case, but with the caveat that the Green function with Neumann boundary conditions is not \emph{uniquely} defined (see discussion below).

We say that a function $G$ on $D\times D$ is \textbf{a covariance function for the Neumann GFF} in $D$ if 
\begin{equation}\label{eq:GFN} \E((\bar{\bf h},\rho_1)(\bar{\bf h},\rho_2)) =\int_{D\times D} \rho_1(x) G(x,y) \rho_2(y) \, \dd x \dd y\end{equation}
for every $\rho_1,\rho_2\in \tilde \cD_0(D)$, where $\bar{\bf h}$ is a Neumann GFF (viewed as a distribution modulo constants in $D$) as in \cref{D:NeumannGFF2}.

Let us immediately make a couple of remarks. 

\begin{itemize}
	\item This need not uniquely define $G$, because the equality is only required to hold for test functions with average $0$. For example, adding any nice enough functions $v(x)$ and $w(y)$ to $G$ will not affect the value of $\int_{D\times D} \rho_1(x) G(x,y) \rho_2(y) \dd x \dd y$. This ill definition is also an inherent property of the Neumann Green function (see below).
	\item As a consequence of \cref{C:NGFFci}, if $G(x,y)$ is a covariance function for the Neumann GFF in $D$, and $T:D'\to D$ is conformal, then $G'(x,y)=G(T(x),T(y))$ is a covariance function for the Neumann GFF in $D'$.
\end{itemize}

We will now show that any choice of 
\textbf{Neumann Green function} in $D$ (if it exists), will be a valid covariance function for the GFF in $D$. 
To explain this, we first need to introduce the \textbf{Neumann problem}\ind{Neumann problem} in $D$.
This is the problem, given $\{\psi, v\}$, of
\begin{equation}\label{E:NProb}
	\text{finding $f$ such that:}
	\begin{cases}
		\Delta f = \psi & \text{ in } D\\
		\displaystyle\frac{\partial f}{\partial n} = v & \text{ on } \partial D,
	\end{cases}
\end{equation}
subject to suitable regularity conditions on $D$. A requirement for the existence of a (weak) solution is that $\psi,v$ satisfy the Stokes condition:
\begin{equation}\label{Stokes}
	\int_D  \psi = \int_{\partial D} v.
\end{equation}
This condition comes from the divergence theorem; the integral of $v$ along the boundary measures the total flux of $\nabla f$ across the boundary, while the integral of  $ \Delta f= \text{div} ( \nabla f) $ inside the domain measures the total divergence of $\nabla f$. 
This solution is then (subject to appropriate conditions on the regularity of $D$) unique, up to a global additive constant. That is, this solution is unique in $\bar H^1(D)$. Existence of a solution is also known, for example, when $D$ is smooth and bounded and $v=0,\psi\in L^2(D)$, \cite[\S 6]{evans} (but we will not use any of these facts).

\ind{Green function!Neumann}
\begin{definition}[Neumann Green function] We say that $G$ is a \textbf{(choice of) Neumann Green function} in $D$, if for every $\rho\in \tilde \cD_0(D)$:
	\begin{equation}\label{Neumannsol}
		f(x) := \int_{D} G(x,y) \rho(y) \dd y
	\end{equation}
	is a solution of the Neumann problem \eqref{E:NProb} in $D$, with $\psi=- \rho$ and $v=0$.
\end{definition}

\begin{prop}\label{P:NGFcov}
	Suppose that $D\subset \C$ is simply connected and has $C^1$ smooth boundary. Then if $G$ is a choice of Neumann Green function, it is a valid choice of covariance for the Neumann GFF $\bar{\bf h}$ in $D$. That is for every $\rho\in {\tilde \cD_0(D)}$ \begin{equation}\label{E:vareq} \E((\bar{\bf h},\rho)^2) =\int_{D\times D} \rho(x) G(x,y) \rho(y) \, \dd x \dd y.\end{equation}
\end{prop}

\begin{rmk}
	Note that adding an arbitrary function of $x$ to $G$ will not affect whether $f$ defined in \eqref{Neumannsol} is a solution to the Neumann problem. In other words, we have the same lack of uniqueness for $G$ as for the covariance of the Neumann GFF.
\end{rmk}

\begin{proof}[Proof of \cref{P:NGFcov}]
	We need to check that if $\rho\in {\tilde  \cD_0(D)}$ and $\bar{\bf h}$ is a Neumann GFF in $D$, then
	\begin{equation}\label{E:goalcov}
		\E((\bar{\bf h},\rho)^2)=\int_{D\times D} \rho(x) G(x,y) \rho(y) \, \dd x \dd y.
	\end{equation}
	Defining $
	f(x) := \int_{D} G(x,y) \rho(y) \dd  y
	$, we will show that both sides are equal to $\|f\|_\nabla^2$.
	
	Note that by assumption the right hand side of \eqref{E:goalcov} is equal to
	$$ \int_D -{\Delta f}(x) f(x)\, \dd x,$$
	which by applying the Gauss--Green formula and the Neumann boundary condition for $f$ is equal to $$\int_D \nabla f(x) \cdot \nabla f(x) \, \dd x = \|f\|_\nabla^2.$$
	
	For the left hand side we use the construction of $\bar{\bf h}$ as the limit as $n\to \infty$ of $\sum_1^n X_j \bar{f}_j$ where the $X_j$s are i.i.d.\ $\mathcal{N}(0,1)$ and the $\bar f_j$s are an orthonormal basis of $\bar{H}^1(D)$. Since this is an almost sure limit in the space of distributions modulo constants, we have that $$(\bar{\bf h},\rho)=\lim_{n\to \infty} \sum_{j=1}^{n} X_j (\bar{f}_j,\rho) \text{ almost surely.}$$ Furthermore, by the Gauss--Green formula again, we have that $(\bar{f}_j, \rho)=(\bar{f}_j,f)_\nabla$ for each $j$, and so
	$$ \E((\sum_{j=1}^{n} X_j (\bar{f}_j,\rho))^2)=\sum_{j=1}^n (\bar{f}_j,f)_\nabla^2. $$
	Note that this is bounded above by $\|f\|_\nabla^2$ for every $n$. Hence, $\sum_{j=1}^n X_j (\bar{f}_j,\rho)$ defines a martingale that is bounded in $L^2$, and so
	$$ \E((\bar{\bf h}, \rho)^2)=\E(\lim_{n\to \infty}(\sum_{j=1}^{n} X_j (\bar{f}_j,\rho))^2)=\lim_{n\to \infty}\E((\sum_{j=1}^{n} X_j (\bar{f}_j,\rho))^2)=\|f\|_\nabla^2,$$
	as desired.
\end{proof}

\ind{Green function!Neumann}

\begin{example}\label{ex:NGFnu}
	We can define a choice of Neumann Green function in the unit disc $\D$ by $$G_N^\D(x,y)=-(2\pi)^{-1}\log|(x-y)(1-x\bar{y})|;\quad x\ne y \in \D.$$
	\indN{{\bf Function spaces}! $G_N^D(x,y)$; choice of Neumann Green function on $
		D$} 
	Indeed, a {tedious but straightforward calculation can be used} to verify that if {$g_y(x):=G_N^\D(x,y)$ for fixed $y\in \D$}, then (in the sense of distributions on $\D$)
	%\nb{NOTE: there seems to be confusion in $x,y$.}
	$$
	\begin{cases}
		\Delta g_y&=- \delta_y \\
		\tfrac{\partial g_y}{\partial n} &= -1/(2\pi) \text{ on } \partial \D.
	\end{cases}
	$$
	This implies that if $\rho\in {\tilde{D}_0(\D)}$ then $f(x) = \int_{D} G_N^\D(x,y) \rho(y) \dd y$ as in \eqref{Neumannsol} is a solution of the Neumann problem with $\psi=-\rho$ and $v=0$. Indeed,
	\begin{itemize}
		\item $\Delta f(x) = \int_\D \Delta g_y(x) \rho (y) \, \dd y = - \int_\D \delta_y(x) \rho(y) \, \dd y = - \int_\D \delta_x(y) \rho(y) \dd y =-\rho (x)$;
		
		\item and for $x\in \partial \D$, $(\partial f /\partial n) (x) = \int_{\D} (\partial g_y /\partial n)(x) \rho(y) \, \dd y = - (2\pi)^{-1} \int_{ \D} \rho (y) \, \dd y = 0$.
	\end{itemize}
	Hence $G^\D_N$ is a choice of Neumann Green function in $\D$, and so also a valid choice of covariance for the Neumann GFF in $\D$.
\end{example}

\begin{example}\label{ex:NGFnu2}
	Define  \begin{equation} \label{eqn:GH} G^\H_N(x,y)=-\frac{1}{2\pi}\log|(x-y)|-\frac{1}{2\pi}\log|(x-\bar{y})|;\quad x\ne y \in \H.\end{equation} 
		 In this case, defining the conformal isomorphism $T: \H\to \D$ by $T(z)=(i-z)(i+z)^{-1}$, we have that if $g_y(x):=G^\H(T^{-1}(x),T^{-1}(y))$, then $\Delta g_y=- \delta_y $ and $\partial g_y/\partial n = - \delta_{-1}$ on $\partial \D$. Similarly to in the previous example, this  implies that $G^\H(T^{-1}(\cdot), T^{-1}(\cdot))$ is a valid choice of Neumann Green function on $\D$. Hence by \cref{P:NGFcov}, it defines a valid choice of covariance function for the Neumann GFF on $\D$. Finally, by conformal invariance of the Neumann GFF (\cref{C:NGFFci}), we see that $G^\H$ is a valid choice of covariance function for the Neumann GFF in $\H$. 
	
	\textbf{Note}: it may seem that we have taken a rather long winded approach in this example. Indeed, one can easily verify that if $g_y^\H(x)=G^\H(x,y)$ then $\Delta g_y^\H = - \delta_y$ on $\H$ and $\partial g_y^\H/\partial n=0$ on $\R$. It is tempting to say that $G^\H$ therefore defines a choice of Green function on $\H$ and so by \cref{P:NGFcov}, is a valid covariance for the Neumann GFF on $\H$. However, one needs to take care that there is an extra ``point at $\infty$'' on the boundary of $\H$ (where, as one can see from the calculations in \cref{ex:NGFnu2}, we actually have a Dirac mass for $\partial g_y^\H/\partial n$). To make this example rigorous it is therefore necessary to map to the unit disc and appeal to conformal invariance -- as carried out above.
\end{example}

\begin{rmk}
	Recall that the Green's function $G_0^D$ for a GFF with zero boundary conditions on a domain $D$ could be defined in terms of  the expected occupation time of ($\sqrt{2}$ times a) Brownian motion killed when leaving $D$: $$G_0^D(x,y)= \int_0^\infty p_t^D(x,y)\, \dd t \; \;\; (x\ne y).$$
	There is a similar relationship between the Neumann Green's function and Brownian motion reflected on the boundary of $D$. The fact that the Neumann Green's function is not uniquely defined is related to the fact that reflected Brownian motion is recurrent. This means if $\tilde{p}_t^D(x,y)$ is the transition density for this reflected Brownian motion, then $\int_0^\infty \tilde{p}_t^D(x,y) \, \dd t$ does not actually converge, so one needs to normalize in some way to obtain a finite quantity. There are many possible ways to do this -- hence the non-uniqueness.
	
	Let us describe this more precisely in the case where $D=\H$. Denoting by $p_t(x,y)$ the transition density of Brownian motion in $\C$, it is easy to see that $\tilde{p}_t^\H(x,y)=p_t(x,y)+p_t(x,\bar{y})=(4\pi t)^{-1}\left(\exp(-|x-y|^2/4t)+\exp(-|x-y|^2/4t)\right)$, which does not have finite integral over $t\in [0,\infty)$.  However, if we look at $\tilde{p}_t^\H(x,y)-\tilde{p}_t^\H(x_0,y)$ for some fixed $x_0\in \H$ (for instance) then the corresponding integral does converge: to $G^\H(x,y)$ as defined in \eqref{eqn:GH} plus the function $\log|x_0-y|$. It is straightforward to check that this integral, for any choice of $x_0$, does define a valid choice of Neumann Green function on $\H$.
\end{rmk}

\begin{rmk}[A choice of covariance for general $D$]\label{rmk:def_NGF}
	Let us remark again that by \cref{C:NGFFci},  if $G^D_N$ is a valid choice of covariance function for the Neumann GFF on some domain $D$, and $T:D'\to D$ is conformal, then $G^D_N(T(\cdot),T(\cdot))$ is a valid choice of covariance function for the Neumann GFF on $D'$.
	
	From this observation and the above examples, we obtain a recipe to define a valid covariance function for the Neumann GFF in any simply connected domain $D$. This works even when the boundary of $D$ is too rough to make sense of the Neumann problem.
	
	We emphasise that any valid choice gives the same value for $\E((\bar{\bf h},\rho_1)(\bar{\bf h}, \rho_2))$ when $\bar{\bf h}$ is a Neumann GFF (viewed as a distribution modulo constants) in $D$ and $\rho_1,\rho_2\in \tilde \cD_0(D)$.
\end{rmk}

\subsection{Neumann GFF as a stochastic process}\label{NGFF:sp}

\ind{GFF!Neumann, stochastic process}

In this section, we will define the Neumann GFF as a stochastic process, similarly to the definition of the Dirichlet GFF in \cref{DGFF:sp}. As with the Dirichlet GFF, this will allow us to ``test'' the Neumann GFF against a wider range of functions: in particular, they need not be smooth or compactly supported inside the domain $D$ and can be non-zero near the boundary. However they can of course not be too singular near the boundary either. We formulate below a condition which, although not optimal, is easy to check in many examples and hence very practical. 

For $x \in \bar \D$ and $A \subset \partial \D$ a Borel set let $q_x(A)$ denote the harmonic measure viewed from $x$, that is, the law of Brownian motion at the first time it leaves $\D$ (if $x \in \partial \D$ then $q_x$ is just a Dirac mass at $x$). Given a Radon measure $m$ on $\bar \D$, let  $\nu_m(A) = \int_{x\in \bar \D} q_x(A) m(\dd x)$. (If $m$ is a probability measure, then $\nu_m$ is simply the law of a Brownian motion at its first exit of $\D$, starting from a point distributed according to $m$; if $m$ is supported on $\partial \D$ then note that $\nu_m = m$).

\begin{definition} \label{defMN}
%Let $\rho = \rho^+ - \rho^- \in \mathfrak{M}_0$ be a signed Radon measure on $D$, such that  $\rho_+(D),\rho_-(D)<\infty$.  Let $T$ be a conformal isomorphism between $D$ and $\D$. Let $T_*\rho$ denote the pushforward of $\rho$ by $T$, that is, $T_*\rho(A)=\rho(T^{-1}(A))$ for $A\subset \partial \D$. Suppose that $T_*\rho$ extends to a signed Radon measure $m$ on $\bar \D$ equivalently, $m = m^+ - m^-$ and 
Let $m$ be a non-negative Radon measure on $\bar \D$ (that is, a finite non-negative measure on $\bar \D$). We say that $m \in \mathfrak{M}_N^+(\D)$ 	\indN{{\bf Function spaces}! $\mathfrak{M}_N^+(\D)$; non-negative Radon measures on $\bar{\D}$ whose restriction to $\D$ is an element of $\mathfrak{M}_0^{\D}$ and such that integral of $m$ against the Poisson kernel on $\partial \D$ is an element $H^{-1/2}(\partial \D)$} if:
\begin{itemize}
\item 
 $m|_\D \in \mathfrak{M}_0^{\D}$, that is, $\iint_{\D^2} m (\dd x) m (\dd y) G_0^{\D} (x,y) < \infty$;   

\item and 
$\nu_{m} \in H^{-1/2}(\partial \D).$
\end{itemize}
If $m = m^+ - m^-$ is a signed Radon measure on $\bar \D$, let us say that $m \in \mathfrak{M}_N(\D)$ \indN{{\bf Function spaces}! $\mathfrak{M}_N(\D)$; difference of two elements in $\mathfrak{M}_N^+(\D)$} if $m^\pm \in \mathfrak{M}_N^+(\D)$. 

Finally, let $D$ be a simply connected domain with a locally connected boundary $\partial D$. Fix $T$ a conformal isomorphism $T:\D\to D$. By definition we say that $\rho \in \mathfrak{M}_N (D)$ \indN{{\bf Function spaces}! $\mathfrak{M}_N(D)$; pushforward of a signed measure in $\mathfrak{M}_N(\D)$ under a conformal isomorphism from $\D$ to $D$} if $\rho = T_*m$ for some $m \in \mathfrak{M}_N(\D)$, where $T_*m$ \indN{{\bf Miscellaneous}! $T_*\mu$; pushforward of a measure $\mu$ by a map $T$} denotes the pushforward of $m$ by $T$, that is, $T_*m( T(A))=m(A)$ for $A\subset  \bar \D$.
\end{definition}

For convenience, we note that the condition $\nu_m \in H^{-1/2} ( \partial \D)$ is implied by the more concrete condition $ \nu_m \in L^2 (\partial \D)$. 

This definition calls for a few comments. First of all, when $D$ is simply connected with a locally connected boundary, a conformal isomorphism from $\D$ to $D$ extends to a continuous map from $\bar \D$ to $\bar D$ (\cite{Pommerenke}). In fact, in terms of the so called ``prime ends of $D$'' (equivalently, the Martin boundary of $D$, \cite{SLEnotes}) the extended map is a homeomorphism. The pushforward $T_*m$ should therefore be viewed as a measure on $D$ and its boundary in this sense. Secondly, given such a measure $\rho$, to check if $\rho \in \mathfrak{M}_N(D)$ we therefore need to check if $m = T^{-1}_* \rho \in \mathfrak{M}_N(\D)$. 
Notice that this does not depend on the choice of the conformal isomorphism $T$: indeed, any two such conformal isomorphisms differ by a conformal automorphism of $\D$ which is a M\"obius map and therefore extends analytically to a neighbourhood of $\D$.

\begin{example}\label{Ex:Neumann_Arc}
As an example, any $\rho \in \mathfrak{M}_0$ compactly supported in $D$ is clearly in $\mathfrak{M}_N$. 
As another example, suppose $D = \D$ and $m$ is the uniform measure on a circular arc of $\partial \D$ of positive length. Then $m \in \mathfrak{M}_N(\D)$. Indeed, in this case clearly $\nu_m = m \in L^2 (\partial \D) \subset H^{-1/2} ( \partial \D)$. As a final example, the measure $\rho$ on the intersection $\gamma \cap D$ of a smooth curve $\gamma$ and a Jordan domain $D$ satisfies $\rho \in \mathfrak{M}_N(D)$, even if $\gamma$ is not fully in $D$. 
\end{example}

\medskip We can now define the index set of the Neumann GFF, which we denote by $\tilde{\mathfrak{M}}_N(D)\subset \mathfrak{M}_N(D)$. By definition this consists of those measures $\rho=\rho^+-\rho^-$ with $\rho^+(\bar D)=\rho^-(\bar D)$ (which corresponds to requiring that the total mass of $\rho$ is zero). More precisely, $\rho \in \tilde{\mathfrak{M}}_N(D)$ if $\rho^\pm = T_* m^\pm$ with $T$ a conformal isomorphism from $\D$ to $D$, $m^\pm \in \mathfrak{M}_N(\D)$, and $m^+ ( \bar \D) = m^-(\bar \D)$.  
%In particular, it allows us to define circle and semi  circle averages of the field.

%Let us first introduce the index set. Analagously to \cref{T:ito} for the Dirichlet GFF, a natural choice for this will be the dual space of $\bar{H}^1(D)$. That is, the space $(\bar{H}^1(D))'$\index[spaces]{$(\bar{H}^1(D))'$; the dual space of $(\bar{H}^1(D),(\cdot,\cdot)_\nabla)$} of continuous linear functionals on the Hilbert space $(\bar{H}^1(D),(\cdot, \cdot)_\nabla)$. Indeed, the Riesz representation theorem immediately implies the following:

\begin{theorem}\label{T:NGFFseries}
	Let $D$ be simply connected with a locally connected boundary. Let $\rho\in \tilde{\mathfrak{M}}_N(D)$.% with action $f\mapsto (f,\rho)$ on $\bar{H}^1(D)$. 
	Then if $\bar{\mathbf{h}}_n=\sum_{j=1}^n X_j \bar{f}_j$ is as in \eqref{NeumannGFF2}, \[ \lim_{n\to \infty}(\bar{\mathbf{h}}_n, \rho)=:(\bar{\mathbf{h}},\rho)\]
exists almost surely and in $L^2(\mathbb{P})$.
%Moreover, if $G$ is any valid choice of covariance for the Neumann GFF, as in \eqref{eq:GFN}, we have 
%\begin{equation}\label{eq:GNdef} \var((\bbh,\rho))=\int_{D\times D} G(x,y) \rho (\dd x) \rho( \dd y)\end{equation}. %Moreover,  $\var(\bbh,\rho)$ is the $(\bar{H}^1(D))'$ norm of $\rho$.
\end{theorem}

%\begin{rmk}
%	\label{rmk:H1dual}

\begin{proof}[Proof of Theorem \ref{T:NGFFseries}]By conformal invariance of $H^1(D)$ and the definition of $\mathfrak{M}_N(D)$, we assume without loss of generality that $D = \D$. 
	The first potential issue to address in the above theorem, is whether $(\bbh_n,\rho)$ makes sense for each \emph{fixed} $n$. We will check that $\rho(\bar{g})$ makes sense for general $\bar{g}\in \bar{H}^1(D)$. (In fact, we will check that our definition of $\mathfrak{M}_N(\D)$ implies that a measure $\rho \in \mathfrak{M}_N(\D)$ defines a continuous linear functional on $\bar H^1(\D)$ and thus is an element in the dual space of $\bar H^1(\D)$; this will imply the result.)
	%By applying a conformal isomorphism $T: D\to \D$, we have that $\bar{g}=\bar f\circ T^{-1}\in \bar{H}^1(\D)$. 
	
	By \eqref{eqn:orthogH1}, $\bar{g}$ can be decomposed as $\bar{g}=g^0+\bar{g}^H$ with $g^0\in H_0^1(\D)$ and $\bar{g}^H\in \overline{\Harm}(\D)$. The assumptions on $\rho$ in Definition \ref{defMN} mean that $\rho|_\D \in \mathfrak{M}_0\subset H^{-1}_0(\D)=(H_0^1(\D))'$ and therefore $\rho(g^0) = \rho|_\D ( g^0)$ is well defined, see for example Remark \ref{R:DH0}. In fact, 
	\begin{equation}\label{eq:norms_NGFF}
	|\rho(g^0)| \le \|g^0\|_{\nabla} \|  \rho\|_{H^{-1}_0}.  
	\end{equation}
	
	Let $g^H$ denote the representative of $\bar g^H$ which has mean zero over $\D$ (since $\rho^+ ( \bar \D) = \rho^-(\bar D)$, the choice of representative does not affect the value of $\rho (g^H)$). Then by the trace theorem, see for example \cite[Theorem 5.36]{adams}, ${g}^H$ is the harmonic extension of a function $ g_\partial $ on the boundary with $g_\partial \in H^{1/2}( \partial \D)$. Moreover, 
	\begin{equation}\label{eq:trace}
	\| g_\partial \|_{H^{1/2} (\partial \D)} \le C \| {g}^H\|_\nabla  
	\end{equation}
	 Since $g^H$ is harmonic and $g_\partial \in H^{1/2} ( \partial \D)$, we have for a \emph{fixed} $x \in \D$, $g^H(x)=\mathbb{E}_x(g_\partial(B_{\tau_\D}))$ (with $B$ a Brownian motion started from $x$ under $\mathbb{E}_x$ and $\tau_D$ its hitting time of $\partial \D$). (Here the regularity of $g_\partial$ on $\partial \D$ is not important, it would suffice that $g_\partial$ is for example an $L^1$ function on $\partial \D$.)
	 
	 By Fubini's theorem and the Cauchy--Schwarz inequality, we have  
	 \begin{align}\label{eq:trace2} 
	 \big|\rho(\bar{g}^H)\big|  = \big|\rho(g^H)\big| 
	& =\big|\int_{\bar \D} \mathbb{E}_x(g_\partial(B_{\tau_\D})) \rho^+(\dd x)-\int_{\bar \D} \mathbb{E}_x(g_\partial(B_{\tau_\D})) \rho^-(\dd x)\big| \nonumber 
	\\ 
	& =  \big|\nu_{\rho^+}(g_\partial)-\nu_{\rho^-}(g_\partial)\big| \nonumber 
	\\ 
	&\le  \big( \|\nu_{\rho^+}\|_{H^{-1/2}(\partial \D)}+\|\nu_{\rho^-}\|_{H^{-1/2}(\partial \D)}\big)\|g_\partial\|_{H^{1/2}(\partial \D)}  \nonumber \\
	&\le C_\rho \| {g}^H\|_\nabla    = C_\rho\| \bar g^H\|_\nabla,
	%\big( \|\nu_{T_*(\rho^+)}\|_{L^2(\partial \D)}+\|\nu_{T_*\rho^-}\|_{L^2(\partial \D)}\big)
		\end{align} 
	for some $C_\rho<\infty$. In the last line we also used \eqref{eq:trace} and the assumption that $ \nu_{\rho^\pm} \in H^{-1/2} (\partial \D)$ in Definition \ref{defMN}. 
	Hence $\rho(\bar{g})$ is well defined, and combining with \eqref{eq:norms_NGFF}, $\rho $ defines a continuous linear functional on $\bar H^1(\D)$. 
	%This in turn implies that $\rho(\bar{f})$ is well defined for any $\bar{f}\in\bar{H}^1(D)$.  
%	\end{rmk}

%As in the above remark, we first notice that if $\rho\in \tilde{\mathfrak{M}}_N$, then applying a conformal isomorphism $T$ from $D$ to the unit disk $\D$, the pushforward measure $T_* \rho$ will be an element of $\tilde{\mathfrak{M}}_N^\D$, and we can write $(\bbh_N,\rho)=(\sum_{n=1}^N X_n (\bar{f}_n\circ T^{-1}),  T_*\rho)$ for each $N$. Since $(\bar{f}_n\circ T^{-1})_n$ is an orthonormal basis for $\bar{H}^1(\D)$ (by conformal invariance of the Dirichlet inner product) it thus suffices to consider the case $D=\D$. 

%In this case, \eqref{eq:trace} and the argument just above it implies that $\rho$ defines a continuous linear functional on $(\bar{H}^1(\D), (\cdot, \cdot)_\nabla)$. 
By the Riesz representation theorem, there exists $\bar g_\rho\in \bar{H}^1(\D)$ with $\rho(\bar{f})=(\bar g_\rho,\bar{f})_\nabla $ for all $f\in \bar{H}^1(\D)$.  This means that \[(\bbh_n,\rho)=\sum_{j=1}^n X_j (\bar{f}_j,\bar g_\rho)_\nabla\]
 is a martingale with mean zero, and uniformly bounded variance:
 \[ \var(\bbh_n,\rho)=\sum_{j=1}^n (\bar{f}_j,\bar g_\rho)_\nabla^2 \le (\bar g_\rho,\bar g_\rho)_\nabla^2 \quad \forall n\ge 1.\]
The martingale convergence theorem yields the result.
\begin{comment}
To see the claim about the variance we observe that, as in the proof of Theorem \ref{GF}, we can find a sequence $(\rho_k)_{k\ge 1}$ with $\rho_k\in \tilde{\cD}_0(D)$ for each $k$, and satisfying $\Gamma_0(\rho-\rho_k,\rho-\rho_k)\to 0$ as $k\to \infty$, as well as $\rho^\pm(D)=\lim_{k\to \infty} \int 
\rho^\pm_k(x)dx$ where $\rho^\pm_k$ are the positive and negative parts of $\rho_k$ for each $k$. By definition of $G$ being a covariance as in \eqref{eq:GFN}, $\var((\bbh, \rho_k))=\int_{D\times D} G(x,y)\rho_k(\dd x)\rho_k(\dd y)$ for each $k$. Moreover, the argument of the previous paragraph implies that $\var(\bbh, \rho_k-\rho)\to 0$ as $k\to \infty$. 
\end{comment}
%	Fix $\rho\in (\bar{H}^1(D))'$. By the Riesz representation theorem, there exists
% an element $g_\rho\in \bar{H}^1(D)$ with $(f,\rho)=(f,g_\rho)_{\nabla}$ for each $f\in \bar{H}^1(D)$. This means that \[(\bbh_n,\rho)=\sum_{j=1}^n X_j (\bar{f}_j,g_\rho)_\nabla\]
% is a martingale with mean zero, and uniformly bounded variance:
% \[ \var(\bbh_n,\rho)=\sum_{j=1}^n (\bar{f}_j,g_\rho)_\nabla \le (g_\rho,g_\rho)_\nabla \quad \forall n\ge 1.\]
%The martingale convergence theorem yields the result.
 \end{proof}

For $\rho_1,\rho_2\in \tilde{\mathfrak M}_N(D)$, we denote
\begin{equation}\label{D:Gip} 
\Gamma_N(\rho_1,\rho_2) = \Gamma_N^D(\rho_1,\rho_2):=\cov((\bbh,\rho_1),(\bbh, \rho_2))
\end{equation} \indN{{\bf Green functions}!$\Gamma_N$; bilinear form, covariance of the Neumann GFF} 
where $(\bbh,\rho_1),(\bbh, \rho_2)$ are the almost sure limits from Theorem \ref{T:NGFFseries}. 
This brings us to the following definition.

\begin{definition}[Neumann GFF modulo constants as a stochastic process]\label{D:NeumannGFF1}
Let $D$ be a simply connected domain with locally connected boundary.	
There exists a unique stochastic process $$(\bar{\bf h}_\rho)_{\rho\in \tilde{{\mathfrak M}}_N}=((\bbh,\rho))_{\rho\in \tilde{{\mathfrak M}}_N},$$ indexed by $\tilde{\mathfrak M}_N$, such that for every choice of $\rho_1,\cdots, \rho_n\in \tilde{\mathfrak M}_N$, $(\bar{\bf h}_{\rho_1},\cdots, \bar{\bf h}_{\rho_n})$ is a centered Gaussian vector with covariance 
	$$\cov(\bbh_{\rho_i},\bbh_{\rho_j}) = \Gamma^D_N(\rho_1,\rho_2).$$
\end{definition}

%Due to \cref{T:NGFFseries}, \cref{D:NeumannGFF1} agrees with the definition of the Neumann GFF as a random distribution modulo constants, in the following sense.

\begin{comment}\begin{lemma}\label{L:Ndefsagree}
	Let $\bar{\bf h}_n = \sum_{j=1}^n X_j \bar{f}_j$ be as in \eqref{NeumannGFF2}. Then for any $\rho\in \tilde\cM$, $(\bar{\bf h}_n, \rho)$ has a limit in $L^2(\P)$ and hence in probability as $n\to \infty$. This limit has the law of $(\bar{\bf h}, \rho)$, where $\bar{\bf h}$ is as in \cref{D:NeumannGFF1}.
\end{lemma}
In particular:
\end{comment} 

%\begin{cor}
By construction, if we restrict the process in \cref{D:NeumannGFF1} to $$(\bar{\bf h},\rho)_{\rho\in {\tilde \cD_0(D)}},$$ then there exists a version of this process defining a random distribution modulo constants, with the same law as the Neumann GFF in \cref{D:NeumannGFF2}.

\begin{rmk}[Conformal invariance]\label{rmk:ciNGFF}
We can also talk about conformal invariance of the Neumann GFF viewed as a stochastic process. Indeed, suppose that $T:D\to \D$ is conformal. Then, as discussed in the proof of \cref{T:NGFFseries}, $\rho \in \tilde{\mathfrak M}_N(D)$ if and only if the pushforward measure $T_*\rho\in \tilde{\mathfrak M}_N(\D)$, and by conformal invariance of the Dirichlet inner product, we have $\Gamma_N^D(\rho_1,\rho_2)=\Gamma_N^\D(T_*\rho_1,T_*\rho_2)$ for all $\rho_1,\rho_2 \in \tilde{\mathfrak{M}}_N(D)$. It follows that if $\bbh^D$ and $\bbh^\D$ are the stochastic processes from \cref{D:NeumannGFF1}, corresponding to the domains $D$ and $\D$, then
$$ ((\bbh^D,\rho))_{\rho\in \tilde{{\mathfrak M}}_N^D}\overset{(\mathrm{law})}{=}  ((\bbh^\D,T_*\rho))_{\rho\in \tilde{{\mathfrak M}}_N^D}.$$
\end{rmk}

With this definition of the Neumann GFF as a stochastic process, it still makes sense to speak of \emph{fixing the additive constant} for the field. In fact, let us now make this notion more precise.

\begin{definition}[Neumann GFF with fixed additive constant]\label{D:NGFFnorm}
Let $D$ be simply connected with locally connected boundary.
Suppose that $\rho_0 \in {\mathfrak M}_N(D)\setminus \tilde{\mathfrak M}_N(D)$. The Neumann GFF ${\bf h}$ with additive constant fixed so that $({{\bf h}},\rho_0)=0$ is the stochastic process defined from $\bar{\bf h}$ in  \cref{D:NeumannGFF1} by setting
	$$({\bf h},\rho)=(\bar{\bf h}, \rho -\frac{\int_{\bar D} \rho(\dd x)}{\int_{\bar D}\rho_0(\dd x)}\rho_0)$$
	for each $\rho \in {\mathfrak M}_N(D)$ where, with an abuse of notation we write $\int_{\bar D} \rho(\dd x)$ for $\int_{\bar \D} T_*\rho (\dd x)$, and $T:D\to \D$ is a conformal isomorphism.
\end{definition}

\begin{rmk}
	For any $\rho_0\in \mathfrak{M}_N(D)\setminus \tilde{\mathfrak{M}}_N(D)$, the Neumann GFF with additive constant fixed so that $(\bf{h}, \rho_0)=0$ has a version which almost surely defined a random distribution, that is, an element of $\mathcal{D}_0'(D)$. Indeed, suppose without loss of generality that $I_{\rho_0}:=\int \rho_0 =1$, and fix an arbitrary $\rho'\in \mathcal{D}_0(D)$ with $I_{\rho'}=1$. Then, by \cref{D:NeumannGFF2} and \cref{T:NGFFseries}, there exists a probability space and a version of $\bar{\bf h}$ defined on this probability space such that $\bar{\bf h}$ defines a distribution modulo constants and $(\bar{\bf h},\rho'-\rho_0)$ is also defined. Then
	\[ (\bf{h}, \rho) := (\bar{\bf h}, \rho-I_\rho \rho') + I_\rho (\bar{\bf h},\rho'-\rho_0)\] 
	is defined simultaneously for all $\rho\in \cD_0(D)$, and defines a version of the Neumann GFF with fixed additive constant from \cref{D:NGFFnorm}. Moreover, $\rho\mapsto (\bf{h},\rho)$ is clearly linear in $\rho$ and $(\mathbf{h},\rho_n)\to 0$ for any sequence $\rho_n$ converging to $0$ in $\mathcal{D}_0(D)$. Thus $\bf h$ defines a random element of $\cD_0'(D)$.
\end{rmk}

In the following, whenever we talk of a Neumann GFF with \emph{arbitrary} fixed additive constant, we mean a Neumann GFF with additive constant fixed -- as defined above -- for some \emph{arbitrary, deterministic} $\rho_0\in {\mathfrak M}_N(D) {\setminus \tilde {\mathfrak M}_N(D)}$.

\begin{example}	[Semicircle averages]\label{ex:semicircle}
	Suppose that $D=\H$ and for $x\in \R$ and $\eps>0$, let $\rho_{x, \eps}$ be the uniform probability distribution on  $\partial B(x, \eps) \cap \H$ of radius $\eps$ about $x$. 
	Then it is straightforward to check that $\rho_{x,\eps}\in {\mathfrak M}_N(\H)$.
	Therefore if $\bf h$ is a Neumann GFF with a fixed additive constant, we can define the $\eps$-semicircle average $({{\bf h}},\rho_{x, \eps})$ of ${\bf h}$ about $x$. 
\end{example}

\begin{rmk}\label{rmk:fixedmodcorr}
	Notice that if $\rho_1,\rho_2\in {\mathfrak M}_N(D)$ with $\int_{\bar D} \rho_1 = \int_{\bar D} \rho_2$, then $\rho_1-\rho_2\in \tilde{{\mathfrak M}}_N(D)$. Hence we can define $(\bar{\bf h}, \rho_1-\rho_2)$ when $\bar{\bf h}$ is a Neumann GFF modulo constants. We can also define $({\bf h},\rho_1-\rho_2)$ whenever ${\bf h}$ is a Neumann GFF with fixed additive constant, and its law will not depend on the choice of additive constant: it will be exactly that of $(\bar{\bf h},\rho_1-\rho_2)$.
\end{rmk}

\subsection{Other boundary conditions}\label{sec:interlude}

\subsubsection{Whole plane GFF}

\ind{GFF!Whole plane}

In this section, we will discuss the \textbf{whole plane GFF}, which we will define as:
\begin{itemize}
\item 
 a distribution modulo constants on the whole complex plane $\C$ whose \emph{odd} and \emph{even} parts are given by reflecting the Dirichlet GFF and Neumann GFF respectively in the $x$ axis. 
 
 \end{itemize}
 
 \noindent Equivalently, we will see that the whole plane GFF coincides with:
 \begin{itemize}
 
 \item a stochastic process with covariance $-(2\pi)^{-1} \log |x-y|$ in a suitable sense, 
 
 \item a local limit of the Dirichlet GFF on large disks, 
 
 \item the spherical GFF constructed in Chapter \ref{S:LCFT} (when the latter is viewed modulo constants and the sphere is identified with the extended complex plane $\hat \C = \C \cup \{\infty\}$). 
\end{itemize}

Just as before, but now with $D=\C$, we define the space of distributions modulo constants on $\C$, $\bar\cD_0'(\C)$, to be the space of continuous linear functionals on $\tilde{\cD}_0(\C)=\{f\in C^\infty(\C) \text{ with compact support and } \int_\C f = 0\}$, equipped with the weak-* topology. \indN{{\bf Function spaces}!$\tilde{\cD}_0(\C)$; smooth functions with compact support and zero average on $\C$} \indN{{\bf Function spaces}!$\bar{\cD}_0'(\C)$; distributions modulo constants on $\C$, that is, the dual space of $\tilde{\cD}_0(\C)$} 

\begin{definition}\label{D:WPGFF} The whole plane GFF, $\bbh^\infty$, is the random distribution modulo constants on $\C$ defined by 
	\begin{equation}\label{E:WPGFF}
		\bbh^\infty=\bbh^\infty_{\mathrm{even}}+\bbh^\infty_{\mathrm{odd}}		
	\end{equation}
where for every $f\in \tilde\cD_0(\C)$ with conjugate $f^*:z\mapsto f(\bar{z})$, \indN{{\bf Miscellaneous}!$f^*$; the conjugate function $z\mapsto f(\bar{z})$} 
\begin{equation*}
	(\bbh^\infty_{\mathrm{even}},f)=\frac{(\bbh^\H,f|_\H+f^*|_\H)}{\sqrt{2}} \quad ; \quad  (\bbh^\infty_{\mathrm{odd}},f)=\frac{(\mathbf{h}_0^\H,f|_\H-f^*|_\H)}{\sqrt{2}}.
\end{equation*}
Here $\bbh^\H$, $\mathbf{h}_0^\H$ are independent Neumann (modulo constants) and Dirichlet GFFs in $\H$ respectively.
\end{definition}

The definition \eqref{E:WPGFF} is natural and should be compared with the fact that any function on $\C$ can be written as the sum of an even and an odd function respectively (where even and odd refer to reflection with respect to the real axis).

\medskip Recalling that the Dirichlet Green function in $\H$ is given by $G^\H_0(x,y)=\frac{1}{2\pi}(-\log|x-y|+\log|x-\bar{y}|)$ and a valid covariance for the Neumann GFF in $\H$ is given by $G^\H_N(x,y)=\frac{1}{2\pi}(-\log|x-y|-\log|x-\bar{y}|)$, a simple calculation (that we leave as an exercise) gives that 
\begin{equation}\label{E:WPGFFcov}
\cov((\bbh^\infty,f_1)(\bbh^\infty,f_2))=\frac{1}{2\pi} \iint_{\C\times \C} \log(\tfrac{1}{|x-y|}) f_1(x) f_2(y) \dd x \dd y
\end{equation}
for each $f_1,f_2\in \tilde{\cD}_0(\C)$. In other words, the covariance of the whole plane GFF (modulo constants) is equal to 
$-\frac{1}{2\pi}\log(|x-y|).$

Recall also from \cref{L:sphereWholePlaneGFF} that the zero average GFF with respect to a Riemannian metric $g$ on the sphere $\hat \C$, had covariance function 
 \begin{equation*}
	G^{\hat \C ,g} (x, y) = \frac1{2\pi} \Big[ - \log (|x-y|)
	+ \bar v_g \big( \log (|x-\cdot|)\big) +  \bar v_g \big( \log (|y-\cdot|)\big)  - \theta_g\Big],
\end{equation*}
where $\bar{v}_g$, $\theta_g$ were defined in that lemma. In particular, for any $f_1,f_2$ such that $\int_{\C} f_1 (x) \dd x = \int_\C f_2(x) \dd x = 0$, we will have 
\begin{equation}
	\iint_{\C\times \C} G^{\hat \C,g}(x,y)f_1(x)f_2(y) \dd x \dd y = \iint_{\C\times \C} \frac{1}{2\pi}\log(\tfrac{1}{|x-y|}) f_1(x) f_2(y) \dd x \dd y.
\end{equation}

This implies the following:
\begin{lemma}
	Let $g$ be a Riemannian metric on the sphere and $\bbh^{\hat{\C},g}$ be $\mathbf{h}^{\hat{\C},g}$ viewed as a distribution modulo constants. Then 
	\[ \bbh^\infty \overset{(\mathrm{law})}{=} \bbh^{\hat{\C},g}\]
\end{lemma}
In fact, just as with the Neumann GFF, we can define the whole plane GFF as a distribution on $\hat{\C}$ (not modulo constants) by fixing the additive constant in some way. For example, we can take the equivalence class representative of $\bbh^{\infty}$ which has average $0$ with respect to $g(z) \, \dd z$. We leave it as an exercise to check that this has precisely the same law as $\mathbf{h}^{\hat{\C},g}$.\\

As mentioned at the start of this subsection, there is another natural way to describe the whole plane GFF, and that is as the local limit of Dirichlet GFFs in large domains. This limit actually exists in a strong sense, and for this we need to recall the definition and some basic properties of the \textbf{total variation distance.}
For two random variables $X,Y$ taking values in the same measurable space $(E,\mathcal{E})$, with respective laws $\mu$ and $\nu$, we define the total variation distance \indN{{\bf Miscellaneous}! $\mathrm{d}_{\tv}$; total variation distance between two measures} between them by \ind{Total variation distance}
\[ \mathrm{d}_{\tv} ( \mu, \nu) = 
\sup_{A\in \mathcal{E}}|\mu( A)-\nu( A)|.\]
With an abuse of notation we also write $\mathrm{d}_{\tv}(X,Y)$ for $ \mathrm{d} _{\tv} ( \mu, \nu)$. Suppose that $\nu$ is absolutely continuous with respect to $\mu$, with Radon--Nikodym derivative $Z = \dd \nu/ \dd \mu$. Then 
for any $A \in \cE$, 
\begin{align*}
| \nu(A) - \mu (A)|  & = \E_\mu[ 1_A ( Z - 1) ]  \le \E_\mu [ | Z - 1|].
\end{align*}
Since $A \in \cE$ is arbitrary, we deduce that
\begin{equation}
\label{eq:tvRN}
\mathrm{d}_{\tv} ( \mu, \nu) \le  \E_\mu [ | Z - 1|].
\end{equation}
Now suppose that $E$ is a metric space and $\cE$ is the associated Borel $\sigma$-algebra. Then it is well known (and easy to check) that given two laws $\mu$ and $\nu$ on $(E, \cE)$ and a coupling $(X,Y)$ (measurable with respect to the product Borel $\sigma$-algebra) of these two laws, then $\mathrm{d} _{\tv} ( \mu, \nu) \le \P ( X \neq Y)$. 
Conversely, if $(E,\mathcal{E})$ is a separable metric measure space, then there necessarily exists a \emph{maximal} coupling of $\mu$ and $\nu$, that is, a coupling $(X,Y)$ measurable with respect to the product $\sigma$-algebra, such that $\P(X\ne Y)=\mathrm{d}_{\tv}(X,Y)$. See, for example, \cite[\S 5.1]{TVcoupling}.

It is straightforward to check that the set of measures on $(E, \cE)$, equipped with the total variation distance, is a complete metric space. The main point is to verify that if $\mu_n$ forms a Cauchy sequence with respect to the total variation distance then $\mu_n(A)$ converges to a limit $\ell(A)$ for any fixed set $A \in \cE$ (and in fact the convergence is uniform). As a consequence the limits $\ell(A)$ necessarily satisfy $\sigma$ additivity (with respect to $A$) and thus define a probability measure on $(E, \cE)$. 

Finally, if $(E, \cE)$ is a metric measure space then convergence of a sequence of measures $\mu_n$ on $(E, \cE)$ to $\mu$ in the sense of total variation distance implies weak convergence: indeed, by the portmanteau theorem, the latter is equivalent to convergence of $\mu_n(A)$ to $\mu(A)$ for every $\mu$-continuity set $A \in \cE$, whereas convergence in the total variance is equivalent to the uniform (in $A \in \cE$) convergence of $\mu_n(A)$. 

We can now state the result.

\begin{theorem}\label{lem:wplim}
	Fix $a>0$ and let $R>a$. Let ${\bf h}^R_0$ be a Dirichlet (zero) boundary condition GFF on $R\D$. Then as $R\to \infty$, 
	\[ \mathrm{d}_{\tv}({\bf h}_0^{R}|_{a\D},\bbh^\infty|_{a\D})\to 0,\]
	when ${\bf h}^{R}_0|_{a\D}$ and $\bbh^\infty|_{a\D}$ are considered as distributions \emph{modulo constants} in $a\D$. In fact the same statement holds when both of these are considered as elements of $H^{-1}_{\mathrm{loc}} ( a\D)$ modulo constants. 
%	In particular, for $\eps>0$ arbitrary, any equivalence class representative of $\bbh^{\infty}|_{\D}$ is almost surely an element of $H^{-\eps}(\D)$.
\end{theorem}

\begin{rmk}
%	It will be useful later on to be able to say that two GFF-like objects that are close in total variation distance can be coupled so that they are close with high probability. This could appear potentially problematic, since the space of (Schwartz) distributions is not a separable metric space.
%	On the other hand, taking the specific setting of the lemma above to illustrate this point, we know that ${\bf h}^{R}|_\D$ almost surely.\ defines an element of the metric space $H^{-\eps}(\D)$ for every $R$ and any $\eps>0$. Since $H^{-\eps}(\D)$ is a separable metric space, this implies that with arbitrarily high probability as $R\to \infty$, we can couple  ${\bf h}^{R}|_\D$ and $\bbh^\infty|_\D$ so that, as distributions, they differ by exactly a constant when restricted to $\D$.
The fact that $\bar {\bf h}^\infty|_{a\D}$ may be viewed as an element of $H^{-1}_{\mathrm{loc}}(a\D)$ modulo constants follows from the definition of the whole plane GFF in \eqref{E:WPGFF} and Remark \ref{R:NGFF_fixed}. 
\end{rmk}

\begin{proof}[Proof of \cref{lem:wplim}] 
	We will first show that 
	\begin{equation}\label{E:TVcauchy} \sup_{R_1,R_2\ge R} \mathrm{d}_{\tv}({\bf h}_0^{R_1}|_{a\D},{\bf h}_0^{R_2}|_{a\D})\to 0,\end{equation}
	when ${\bf h}^{R_1}_0|_{a\D}$ and ${\bf h}^{R_2}_0|_{a\D}$ are considered as distributions \emph{modulo constants} in $a\D$. Without loss of generality, suppose that $R_2\ge R_1$. Then by the Markov property of the Dirichlet GFF (\cref{T:mp}), we can write ${\bf h}^{R_2}=\tilde{{\bf h}}^{R_1}+\varphi$, where $\tilde{{\bf h}}^{R_1}$ has the law of ${\bf h}^{R_1}$, and $\varphi$ is independent of $\tilde{{\bf h}}^{R_1}$ and almost surely harmonic in $R_1 \D$. The proof of this lemma will essentially follow from the fact that, when viewed modulo constants and restricted to $a\D$, $\varphi$ is very small. 
	
Indeed, if we define $\varphi_0=\varphi-\varphi(0)$ then by independence and harmonicity, $\var({\bf h}_1^{R_2}(w)-{\bf h}_1^{R_2}(0))=\var({\bf h}_1^{R_1}(w)-{\bf h}_1^{R_1}(0))+\var(\varphi_0(w))$ for any $w\in \partial (8a\D)$ (say). Since we have the explicit expressions
\begin{equation}\label{eq:Greenlargedisc}
2\pi G^{R_i \D}_0(x,y)=\log R_i +\log|1-(\bar{x}y/R_i^2)|-\log (|x-y|) 
\end{equation}
%$$ 2\pi G^{R_2 \D}_0(x,y)=\log R_2+\log|1-(\bar{x}y/R_2^2)|-\log (|x-y|) $$
for $x\ne y \in R_i\D$  ($i=1,2$), it follows easily that 
	\begin{equation}\label{eq:varphi_0} \sup_{R_1,R_2\ge R}\sup_{w\in \partial(8a\D)} \var(\varphi_0(w))\to 0 \text{ as } R\to \infty. \end{equation}

	Now, note that ${\bf h}^{R_2}$ and $\tilde{{\bf h}}^{R_1}+\varphi_0$ differ by exactly a constant in $R_1\D$. So we would be done with the proof of \eqref{E:TVcauchy} if we could show that the laws of $$\tilde{{\bf h}}^{R_1}+\varphi_0 \text{ and } \tilde{{\bf h}}^{R_1}$$ are close in total variation distance when restricted to $a\D$ (uniformly in $R_2\ge R_1\ge R$ as $R\to \infty$). The idea for this is to use the explicit expression for the Radon--Nikodym derivative between a zero boundary GFF and a zero boundary GFF plus an $H_0^1$ function; see \cref{lem:CMGFF}.
	
	The first obstacle here is that $\varphi_0$ is not actually $H_0^1(R_1\D)$. To get around this, we introduce $\tilde \varphi(z)=\psi(|z|)\varphi_0(z)$ for $z\in R_1\D$, where $\psi:[0,R_1]\to [0,1]$ is smooth, equal to $1$ on $[0,a]$, and equal to $0$ on $[2a,R_1]$. Note that $\tilde{\varphi}\in H_0^1(R_1\D)$ and that $\tilde{\varphi}=\varphi_0$ in $a\D$. Moreover, \emph{ conditionally }on $\tilde{\varphi}$, the Radon--Nikodym derivative between the laws of $\tilde{{\bf h}}^{R_1}$ and $\tilde{{\bf h}}^{R_1}+\tilde{\varphi}$ is given by
	\begin{equation}\label{eq:RN1} Z:=\frac{\exp((\tilde{{\bf h}}^{R_1},\tilde{\varphi})_\nabla)}{\exp((\tilde{\varphi},\tilde{\varphi})_\nabla)},\end{equation}
	see \cref{lem:CMGFF}.
To complete the proof of \eqref{E:TVcauchy} it suffices (by the definition of total variation distance) to show that
	\eqref{eq:RN1} tends to $1$ in $L^1(\P)$, uniformly over $R_2\ge R_1\ge R$ as $R\to \infty$.
	
	To show this, we will first prove that
	\begin{equation}\label{eqn:eminus1} \sup_{R_1,R_2\ge R}\E(\e^{(\tilde{\varphi},\tilde{\varphi})_\nabla}-1)\to 0\end{equation}
	as $R \to \infty$. 
	To see this, note that $\nabla \tilde \varphi = 0$ outside $2a\D$ and for $x \in 2a\D$, $|\nabla \tilde \varphi| \le c_1(  |\nabla \varphi_0| + \sup_{x \in 2a\D} |\varphi_0| ) $, where the constant $c_1$ depends only on $\psi$. 
	
	We now make use of the fact that $\varphi_0$ is harmonic in $2a\D$ and of two well known inequalities for harmonic functions:
	
	\begin{lemma}
Let $u$ be a harmonic function in $4a\D$. Then there exists a universal constant $C>0$ such that
$$
\sup_{x \in 2a\D} | \nabla u | \le C \sup_{x\in \partial(4a\D)} |u|
$$	
	\end{lemma}

This follows for example from Theorem 7 in \cite[\S 2.2]{evans} and the maximum principle for harmonic functions. 
The second inequality we use is a consequence of Harnack's inequality:

\begin{lemma}
Let $u$ be a harmonic function in $8a\D$. Then there exists a universal constant $C>0$ such that for any $x \in 4a \D$, 
$$
|u(x) | \le C | u (0)|.
$$
\end{lemma}

See, for example, Theorem 11 in \cite[\S 2.2]{evans} for a proof when $u$ is assumed to be non-negative; the general case follows by considering $u - \inf_{8 a\D} u$.  
	
Combining these two estimates, we deduce 
$$
\sup_{x\in 2a \D} | \nabla \tilde \varphi(x) | \le c_2 |\varphi_0(0) | = c_2 | \int_{\partial (8a \D) } \varphi_0(x) \rho(\dd x)|, 
$$
 where $\rho$ is the uniform measure on the circle $\partial (8a\D)$. Therefore applying Cauchy--Schwarz, 
  $$
  \E(e^{(\tilde{\varphi},\tilde{\varphi})_\nabla}-1) \leq \E(e^{c_2\int_{\partial(8a\D)} |\varphi_0(w)|^2 \, \rho(\dd w)}-1)
  $$
   which by Jensen's inequality is less than 
   $$
   \E(\int_{\partial (8a\D)}(e^{c_2|\varphi_0(w)|^2}-1)\, \rho(\dd w))\le \int_{\partial (8a\D)}\E(e^{c_2|\varphi_0(w)|^2}-1)\, \rho(\dd w). 
   $$ 
   Note that since $c_2$ is a fixed constant and $\varphi_0(w)$ is a centred Gaussian random variable with arbitrarily small variance (uniformly over $\partial (8a \D)$) as $R\to \infty$, these expectations will all be finite for $R_2\ge R_1\ge R$ large enough. Moreover, the right hand side of the above expression will go to $0$ uniformly in $R_2\ge R_1\ge R$ as $R\to \infty$.
	To conclude, we simply observe that conditionally on $\tilde{\varphi}$, the random variable $Z$ from \eqref{eq:RN1} is log normal with parameters $(-(\tilde{\varphi},\tilde{\varphi})_\nabla^2/2 , (\tilde{\varphi},\tilde{\varphi})_\nabla)$. Hence 
	\begin{equation*} 
	\E(|Z-1|^2)=\E\Big(\E\Big(|Z-1|^2 \, \Big| \, \tilde{\varphi}\Big) \Big)=\E(e^{(\tilde{\varphi},\tilde{\varphi})_\nabla}-1).
	\end{equation*}
	By \eqref{eqn:eminus1}, this shows that $Z$ converges to 1 in $L^2(\P)$ (hence in $L^1(\P)$ and thus completes the proof of \eqref{E:TVcauchy}. 
	
With \eqref{E:TVcauchy} in hand, we know that ${\bf h}^R_0|_{a\D}$ (viewed as an element of $H^{-1}_{\mathrm{loc}}( a\D)$, modulo constants) is a Cauchy sequence with respect to total variation distance, and so its law has a limit (say $\mu$) in total variation distance as $R\to \infty$.  It remains to identify $\mu$ with the law of $\bbh^\infty|_{a\D}$.

Fix a test function $\ph \in \tilde \cD_0(a\D)$. Then 
$$
(h_0^R, \ph) \sim \cN (0, \sigma^2_\ph) \text{ where } \sigma^2_\ph = \iint_{(R\D)^2} G_{0}^{R\D} (x,y) \ph(x) \ph(y) \dd x \dd y.
$$
Using the expression in \eqref{eq:Greenlargedisc} for $G_0^{R\D}$ and the fact that $\int_{a\D} \ph(x) \dd x = 0$, we see that 
\begin{equation}\label{var_WP}
\sigma^2_\ph \to \frac{-1}{2\pi} \iint_{\C^2} \log | x- y| \ph (x) \ph(y) \dd x \dd y = \var ( (\bbh^\infty|_{a\D}, \ph) ). 
\end{equation}
Now let $\ph_1, \ldots, \ph_k$ be arbitrary test functions in $\tilde \cD_0(a\D)$ and fix $x_1, \ldots, x_k \in \R$. Consider the event 
$A = \{ h \in H^{-1}_{\mathrm{loc}}(a \D) : ( \h, \ph_1)< x_1, \ldots ( h, \ph_k) < x_k\}$ and let $\cA$ denote the set of events of this form. Since the law of $\bf h_0^R$ converges to $\mu $ in total variation, we immediately deduce that for all $A \in \cA$,
$$
\mu(A) = \lim_{R\to \infty} \P ( {\bf h}_0^R \in A), 
$$
but this also agrees with $\P ( \bbh^\infty|_{a\D} \in A) $ by \eqref{var_WP} and properties of Gaussian random variables. Thus $\mu$ agrees with the law of $\bbh^\infty|_{a\D}$ on $\cA$. However the latter is a $\pi$-system which clearly generates the Borel $\sigma$-field on $H^{-1}_{\mathrm{loc}}(a \D)$ modulo constants, hence we conclude by Dynkin's lemma. 
%
%
%More precisely, to show that for any $f_1,\ldots, f_n\in \tilde{\cD}_0(\D)$ the joint limit $\{X_1,\ldots, X_n\}$ of $\{({\bf h}^R_0,f_1), \ldots, ({\bf h}^R_0,f_n)\}$ is Gaussian and centered with covariance $$ \cov(X_i,X_j)=\frac{1}{2\pi}\int_{\C\times \C} f_i(x) f_j(y) \log(\tfrac{1}{|x-y|}) \dd x \dd y$$ for every $1\le i,j\le n$. The fact that the limit is centered and Gaussian is immediate, since $\{(\bf h^R_0,f)\}_{f\in \tilde\cD_0(\D)}$ is centered and Gaussian for every $R$. The limiting covariance is also easily identified from the explicit expression for $G_0^{R\D}$ above.
\end{proof}

As a corollary, we deduce that the whole plane GFF restricted to $\D$ inherits from the Dirichlet GFF the same Markov property: 

\begin{cor}[Markov property for the whole plane GFF]\label{cor:wp}
	$${\bf h}^\infty|_{\D}={\bf h}^\D+\varphi,$$ where ${\bf h}^\D$ has the law of a Dirichlet boundary condition GFF in $\D$, and $\varphi$ is a harmonic function modulo constants that is independent of ${\bf h}^\D$.
\end{cor}

\subsubsection{Dirichlet--Neumann GFF}

\ind{GFF!Dirichlet--Neumann}

Another variant of the GFF that is important, because it appears in a natural Markov property for the Neumann GFF, is the GFF with ``mixed'' boundary conditions. Here we will discuss one specific version, which is a distribution defined in $\D_+=\D\cap \H$ and (heuristically speaking) has free/Neumann boundary conditions on $[0,1]$ and zero/Dirichlet boundary conditions on $\partial \D\cap \H$.

\begin{definition}[Dirichlet--Neumann GFF]
	Suppose that ${\bf h}_0^\D$ is a Dirichlet GFF in $\D$. Then the Dirichlet--Neumann GFF, ${\bf h}^{\mathrm{DN}}$, is defined to be $\sqrt{2}$ times its even part
	\[{\bf h}^{\mathrm{DN}}:=\sqrt{2}({\bf h}^\D_0)_{\mathrm{even}}, \text{ where }  (({\bf h}^\D_0)_{\mathrm{even}},\rho):=\frac{({\bf h}_0^\D,\rho)+({\bf h}_0^\D,\rho^*)}{2} \text{ for }  \rho\in \mathcal{D}_0(\D_+) \] which is a random distribution on $\D_+$.
\end{definition}

Putting this together with \cref{lem:wplim} and \cref{D:WPGFF} we obtain a useful boundary Markov property for the Neumann GFF. Indeed recall that by definition of the whole plane GFF,
$$
\bbh^{\H}|_{\D_+} = \sqrt{2} \bbh^\infty_{\mathrm{even}}|_{\D_+}
$$
where we recall that $\bbh^{\H}$ is a Neumann boundary condition GFF in $\H$, modulo constants. On the other hand, by the Markov property of the whole plane GFF (\cref{cor:wp}) we also know that 
$$
\sqrt{2} \bbh^\infty_{\mathrm{even}} |_{\D_+} = \sqrt{2} ({\bf h}_0^\D)_{\mathrm{even}} |_{\D_+} + \sqrt{2} \ph_{\mathrm{even}}|_{\D_+}, 
$$
where $\ph_{\mathrm{even}}$ is the even part of the harmonic function $\ph$ appearing in \cref{cor:wp} and is thus also harmonic over all of $\D$, and $({\bf h}^0_\D)_{\mathrm{even}} $ is the even part of a Dirichlet GFF in $\D$. 
%let us write $\varphi_{\mathrm{even}} $ for the even part of a distribution $\ph$. 
By definition, the first term on the right hand side is the Dirichlet--Neumann GFF on $\D_+$. We thus obtain the following (since $\sqrt{2} \ph_{\mathrm{even}}$ is also a harmonic function in $\D$, and changing notations slightly for later convenience). 

% we see that on the one hand, the even part of the whole plane GFF $\bbh^\infty$ restricted to $\D_+$ is equal to $1/\sqrt{2}$ times the Neumann GFF $\bbh^\H$ restricted to $\D_+$. On the other hand, $\bbh^\infty$ restricted to $\D$ can be written as $\mathbf{h}^\D_0+\varphi$  where $\varphi$ is a harmonic function modulo constants and $\mathbf{h}^\D_0$ is a Dirichlet GFF in $\D$. Thus it has even part equal to the even part of $\mathbf{h}^\D_0$, plus a harmonic function $\varphi_{\mathrm{even}}$ modulo constants in $\D_+$ with Neumann boundary conditions on $(-1,1)$. (Note that the even part of a harmonic function in $\D$ defines a harmonic function in $\D_+$ with Neumann boundary conditions on $(-1,1)$). From the definition of the Dirichlet--Neumann GFF we therefore get the following.

\ind{Markov property!Neumann GFF: boundary}
\begin{prop}
	[Boundary Markov property]\label{P:MP2}
	Let ${\bf h}^\H$ be a Neumann GFF on $\H$ (considered modulo constants). Then we can write
	\[{\bf h}^\H|_{\D_+}={\bf h}^{\mathrm{DN}}+\varphi_{\mathrm{even}}\]
	where the two summands are independent, ${\bf h}^{\DN}$ has the law of a Dirichlet--Neumann GFF in $\D_+$, and $\varphi_{\text{even}}$ is a harmonic function modulo constants in $\D_+$, smooth up to and including $(-1,1)$ and satisfying  Neumann boundary conditions along $(-1,1)$.
\end{prop}

We conclude this section with one further comment, that will be useful at the end of this chapter and in \cref{S:zipper}. It can be used to say, roughly speaking, that any (nice enough) way of fixing the additive constant for a Neumann GFF in $\H$ will produce a field with the same behaviour when looking very close to the origin. Moreover, this will still be true if we condition on the realisation of the field far away from the origin.

\ind{Total variation distance}

\begin{lemma}\label{lem:hDNcircgood}
	Suppose that ${\bf h}$ is a Neumann GFF in $\H$, with additive constant fixed so that it has average $0$ on the upper unit semicircle (this makes sense by Example \ref{ex:semicircle}). Let ${\bf h}^{\DN}$ be an independent Dirichlet--Neumann GFF in $\D_+$. 
	Then, for any $K>1$, the total variation distance between
	\begin{itemize}\setlength\itemsep{0em}
		\item the joint law of $({\bf h}|_{K\D_+\setminus \D_+}, {\bf h}|_{\delta\D_+})$  and
		\item the (independent product) law $({\bf h}|_{K\D_+\setminus \D_+},{\bf h}^{\DN}|_{\delta\D_+})$,
	\end{itemize} 
	tends to $0$ as $\delta\to 0$. Note that the fields can be viewed as distributions here, rather than just distributions modulo constants.
	
	The same thing is true if $\mathbf{h}$ is replaced by $\mathbf{h}+\mathfrak{h}$ where $\mathfrak{h}$ is harmonic in a neighbourhood of $0$ and has $\mathfrak{h}(0)=0$. 
\end{lemma}

\begin{rmk}\label{rmk:hDNcircgood}Note that the Lemma in particular implies that $\mathrm{d}_{\tv}(\mathbf{h}|_{\delta\D_+},\mathbf{h}^{\mathrm{DN}}|_{\delta \D_+})\to 0$ as $\delta\to 0$. So the Lemma also holds if we replace the second pair $({\bf h}|_{K\D_+\setminus \D_+},{\bf h}^{\DN}|_{\delta\D_+})$ by $({\bf h}|_{K\D_+\setminus \D_+},\tilde{{\bf h}}|_{\delta\D_+})$ where $\tilde{{\bf h}}$ is independent of $\bf h$ but with the same marginal law. \end{rmk}
\begin{proof}
%After scaling by $R=1/\delta$ our goal is to
% compare the joint laws of:
%	\begin{itemize}\setlength\itemsep{0em}
%		\item the joint law of $({\bf h}|_{RK\D_+\setminus R\D_+}, {\bf h}|_{\D_+})$  and
%		\item the (independent product) law $({\bf h}|_{ RK\D_+\setminus R\D_+},{\bf h}^{\DN}|_{\D_+})$,
%	\end{itemize}
%	and show that their total variation distance tends to zero as $R \to \infty$. 
The proof essentially follows along the same lines as \cref{lem:wplim} and taking even parts.
		More precisely, let $R\gg K$ be large, and write 
	$$
	\tilde{{\bf h}}^{R\D}={\bf h}^{R\D}-{\bf h}_1^{R\D} (0),
	$$ 
	for ${\bf h}^{R\D}$ a Dirichlet GFF in $R\D$ and ${\bf h}_1^{R \D}(0)$ its unit circle average around $0$. By 
	\cref{P:MP2} and considering even parts (and multiplying by a factor $\sqrt{2}$), it suffices to prove that as $R \to \infty$, and for $  {\bf h}^\D $ a Dirichlet GFF in $\D$ that is independent of  $\tilde{{\bf h}}^{R \D}$
	\begin{equation}\label{eq:asymp_indep}
\lim_{\delta\to 0} \lim_{R\to \infty}	\mathrm{d} _{\tv}\left((\tilde{{\bf h}}^{R\D}|_{K\D\setminus \D},\tilde{{\bf h}}^{R\D}|_{ \delta\D}), (\tilde{{\bf h}}^{R\D}|_{K\D\setminus \D},{\bf h}^{\D}|_{ \delta\D})\right)\to 0. 
	\end{equation} 
	 This follows from the same argument as \cref{lem:wplim}. 
	 That is, we use the Markov property of the GFF to write  $\tilde{\mathbf{h}}^{R\D}|_{\D}=\mathbf{h}+\varphi$ where $\varphi$ is measurable with respect to $\tilde{\mathbf{h}}^{R\D}|_{R\D\setminus \D}$ and harmonic in $\D$ with $\varphi(0)=0$, while $\mathbf{h}$ is independent of $\tilde{\mathbf{h}}^{R\D}|_{R\D\setminus \D}$ with the law of $\mathbf{h}^\D$. 
	To show \eqref{eq:asymp_indep},
	 %the total variation distance in \eqref{eq:asymp_indep} converges to $0$ as $R\to \infty$ and then $\delta\to 0$, 
	  as in the proof of \cref{lem:wplim}, it suffices to prove that 
	  \begin{equation}\label{eq:abitlikeThm627}
	  \limsup_{R\to \infty} \sup_{z \in \partial (8 \delta \D)} \var( \ph(z)) \to 0
	  \end{equation}
	  %the second moment of $\varphi$ converges to $0$ uniformly on $\partial(8\delta \D)$
	   as  $\delta\to 0$ (analogous to \eqref{eq:Greenlargedisc}). Indeed, once \eqref{eq:abitlikeThm627} is shown, we can use the same considerations as in \cref{lem:wplim} on harmonic functions to control the Radon--Nikodym deriatve, so that  \eqref{eq:asymp_indep} follows. On the other hand, \eqref{eq:abitlikeThm627} follows from a straightforward Green's function calculation as in the proof of \cref{lem:wplim}.

	 If a deterministic function $\mathfrak{h}$ that is harmonic in a neighbourhood of the origin with $\mathfrak{h}(0)=0$ is added to the field, one must replace $\varphi$ with $\varphi+\mathfrak{h}$ in the above argument. However, since $\mathfrak{h}|_{\partial(8\delta \D)}\to 0$ as $\delta\to 0$, the same argument goes through.
	 	 %(after rescaling) that the Radon--Nikodym derivative, between the law of ${\bf h}^{\D}|_{\delta \D}$ and the conditional law of $\tilde{{\bf h}}^{R\D}|_{\delta \D}$ given $\tilde{{\bf h}}^{R\D}|_{K\D\setminus \D}$, tends to $0$ in $L^1(\P)$ as $\delta\to 0$. This gives the desired claim.
\end{proof}

\begin{comment}\begin{rmk}\label{rmk:hDNcircgood}
	Note that the proof (and therefore the Lemma) will still hold if we replace ${\bf h}$ by ${\bf h}+\mathfrak{h}-\mathfrak{h}(0)$ where $\mathfrak{h}$ is a deterministic harmonic function in $\D^+$ with Neumann boundary conditions on $[-1,1]$. Moreover the convergence will be uniform over $\{\mathfrak{h}: \sup_{z\in \D^+}|\mathfrak{h}(z)|\le C\}$ for any $C$.
\end{rmk}
\end{comment}

\subsection{Semicircle averages and boundary Liouville measure}

Let \begin{equation}\label{hNLQG} \bar h=\sqrt{2\pi}\bar{\bf h}\end{equation} where $\bar{ \bf h}$ is a Neumann GFF on $\H$ modulo constants (recall that we use a bar in order to distinguish statements concerning the Neumann GFF modulo constants and Neumann GFFs with fixed additive constants). We will refer to both $\bar{h}$ and $\bar{\bf h}$ as ``a Neumann GFF'' in what follows: the use of bold font distinguishing between the different multiples as in the Dirichlet GFF setting.   An immediate consequence of our previous considerations is the following fact. Recall our notation from \cref{ex:semicircle} that for $x \in \R$ and $\eps>0$, $\rho_{x,\eps}$ denotes the uniform distribution on the upper semicircle of radius $\eps$ around $x$ (and recall also that $\rho_{x,\eps} \in \mathfrak{M}_N(\H)$).
%For $x \in \R$, let $h_\eps(x) $ denote the ``average of $h$'' on the upper semicircle $\partial B(x, \eps) \cap \H$ of radius $\eps$ about $x$; that is if $\rho^+_{x, \eps}$ is the uniform probability distribution on this set, then $h_\eps(x)=(h, \rho^+_{x, \eps})$. \change{Note that $\int |\rho_{x, \eps}^+(dw)| |\rho_{x,\eps}^+(dz)||G(z,w)|<\infty$, so that this is well defined.}

\begin{theorem}
	\label{T:BM}
	For any $x \in \R$, the finite dimensional distributions of the process $$(X_t)_{t\in \R}:=((\bar h, \rho_{x,e^{-t}}-\rho_{x,1}))_{t\in \R}$$ are those of a two-sided Brownian motion with variance 2 (so $\var (X_t) = 2|t|$).
\end{theorem}
\ind{Circle average}
\noindent Note that the statement of the theorem makes sense, since for any $\eps>0$, $\rho_{x, \eps} - \rho_{x, 1}\in \tilde{\mathfrak M}_N^\H$.
By  \cref{rmk:fixedmodcorr}, this also means that if $h$ is a Neumann GFF in $\H$ with additive constant fixed \emph{in any way}, and $h_{\eps}(x):=(h,\rho_{x, \eps})$, %for $x\in \R$ and $\eps>0$,
then $$(h_{e^{-t}}(x) - h_1(x))_{t\in \R}$$ is a two-sided Brownian motion with variance $2$.

\begin{proof}[Proof of \cref{T:BM}]
	Without loss of generality we may take $x = 0$. Then by conformal invariance (actually just scale invariance) of $\bar{h}$, it follows that $X$ has stationary increments. Moreover, by applying the Markov property (a scaled version of \cref{P:MP2}) in the semidisc of radius $e^{-t}$ about 0 for any $t$, we see that $(X_r)_{r\le t}$ and $(X_s-X_t)_{s\ge t}$ are independent. Hence, $X$ has stationary and independent increments. % and fix the normalisation for $h$  so that $h_1(0) = 0$.
	%Then by scaling properties of $h$ (with this normalisation) it is clear that $X_t := h_{e^{-t}}(0)$ has stationary increments. Moreover by an easily verified Markov property (see exercises), these increments are independent.
	
	Since the increments are also Gaussian with mean zero and finite variance, it must be that $X_t = B_{\kappa t }$ for some $\kappa >0$, where $B$ is a standard Brownian motion. It remains to check that $\kappa = 2$, but this follows from the fact that a choice of Neumann GFF covariance in the upper half plane is given by $G^\H(0,y) = (2\pi)^{-1} \times 2 \log (1/\eps)$ if $|y| = \eps$: see \eqref{eqn:GH}. %The result follows.
\end{proof}

Having identified the ``boundary behaviour'' of the Neumann GFF, we can now construct a random measure supported on the boundary of $\H$. As it turns out, the measure of interest to us is again given by an ``exponential of the Neumann GFF'', but the multiplicative factor in the exponential is $\gamma/2$ rather than $\gamma$. It might initially seem that the factor $(1/2)$ appearing in this definition comes from the fact that we are measuring lengths rather than areas. We want to emphasise that this is however \emph{not} the real reason: instead, it is more related to the fact that the variance of the Brownian motion describing circle averages on the boundary has variance two (see \cref{T:BM}). This will guarantee that the boundary measure enjoys the same change of coordinate formula as the bulk measure, as should be the case. Alternatively, this can be seen as a consequence of the fact that the so called ``quantum length of SLE'' can be measured via this boundary length, and an application of the KPZ formula shows that the corresponding quantum scaling exponent is, as it turns out, always $\Delta = 1/2$.

%The reason for this choice is rather deep, and has to do with the fact that we plan to use it to measure the ``quantum length of an SLE''. (Another justification comes from the fact that it satisfies the same KPZ equation as in the bulk case).

\begin{theorem}[Boundary Liouville measure for the Neumann GFF on $\H$]
	\label{T:boundary}
	Let $h$ a Neumann GFF in $\H$ as in \eqref{hNLQG} but with additive constant fixed in an arbitrary way. Define a measure $\mathcal{V}_\eps$ on $\R$ by setting $\mathcal{V}_{\eps}(\dd x) = \eps^{\gamma^2/4} e^{(\gamma/2) h_\eps(x)} \dd x$.\indN{{\bf Gaussian multiplicative chaos}! $\mathcal{V}_\eps$; approximation of $\mathcal{V}$ at spatial scale $\eps$}  Then for $\gamma <2$, the measure $\mathcal{V}_\eps$ converges almost surely along the dyadic subsequence $\eps = 2^{-k}$ to a non-trivial, non-atomic measure $\mathcal{V}$ \indN{{\bf Gaussian multiplicative chaos}! $\mathcal{V}$; boundary Gaussian multiplicative chaos on the boundary for a field on a domain} called the boundary Liouville measure.
\end{theorem}\ind{Liouville measure!Boundary}\ind{Boundary Liouville measure|see{Liouville measure}}

\begin{proof}
	This can be proved as in %follows from \cref{T:conv} with $h$ the (normalised) Neumann GFF, $\sigma$ equal to Lebesgue measure on $\R$ and $\theta$ the semi-circle uniform measure. Equivalently, one can prove the convergence as in
	\cref{S:Liouvillemeasure}, proof of \cref{T:liouville}, using the Markov property and \cref{T:BM}. We leave the details as Exercise \ref{Ex:GMCboundary}.
\end{proof}

Note the scaling in $\mathcal{V}_{\eps}$, which is by $\eps^{\gamma^2/4}$. This is because, as proved in  \cref{T:BM} (also see the discussion below), when $x \in \R$ and $h=\sqrt{2\pi}{\bf h}$ for $\bf h$ a Neumann GFF with arbitary fixed additive constant, we have $\var h_\eps(x) = 2 \log (1/\eps) + O(1)$. %the use of the Markov property is justified with appropriate modifications -- though in fact, as remarked before, this property is an unnecessary assumption, see for example \cite{BerestyckiGMC}).
\begin{rmk}\label{rmk:measuremodconstants}
	The law of $\mathcal{V}$ above \emph{does} depend on the choice of additive constant for $h$. If one starts with a Neumann GFF modulo constants, then the boundary Liouville measure can be defined as a measure \emph{up to a multiplicative constant}.
\end{rmk}

\dis{As with $\cM$, we will sometimes also use the notation $\cV_h$ or $\cV_h^\gamma$ to indicate the dependence of $\cV$ on the underlying field $h$ or the field $h$ and the parameter $\gamma$.}
\indN{{\bf Gaussian multiplicative chaos}! $\cV_h$; Gaussian multiplicative chaos on the boundary associated with a field $h$ on a domain}
\indN{{\bf Gaussian multiplicative chaos}! $\cV_h^\gamma$; Gaussian multiplicative chaos on the boundary associated with a field $h$ on a domain and parameter $\gamma$}

For general $D$, $h=\sqrt{2\pi}{\bf h}$ a Neumann GFF in $D$ with arbitrary fixed additive constant, and $z,\eps$ such that $B(z,\eps)\subset D$, we can also define the circle average $(h,\rho_{z,\eps})=:h_\eps(z)$. Although we use the same notation $h_\eps(\cdot)$ for circle averages and semicircle averages it should always be clear which one we refer to, depending whether the argument lies, respectively, in the bulk or on the boundary of $D$.

\begin{definition}[Bulk Liouville measure for the Neumann GFF]
	When $h=\sqrt{2\pi}{\bf h}$ is a Neumann GFF with some arbitrary fixed additive constant and $\gamma<2$, we can also define the bulk Liouville measure $$\mathcal{M}(\dd z):=\lim_{\eps\to 0}\eps^{\gamma^2/2}\e^{\gamma h_\eps(z)}\, \dd z,$$ exactly as for the Dirichlet GFF.
\end{definition}

The existence of this limit follows from the construction of GMC measures for general log-correlated Gaussian processes in \cref{S:GMC}. The analogue of \cref{rmk:measuremodconstants} also applies in this case.

\begin{rmk}
	\label{rmk:finite_mass}
	Adapting the results of \cref{S:GMC}, it is not hard to see that for any fixed compact set of $\R$ (respectively $D$) the boundary (respectively bulk) Liouville measure of a Neumann GFF on $\H$ (respectively  $D$) will assign finite and strictly positive mass to that set with probability one. On the other hand, note that for the Neumann GFF on $D$, since there is no bounded function $g$ such that \eqref{cov} holds uniformly on $\bar{D}$, it is not obvious that the bulk Liouville measure of $D$ is finite. This will be addressed in Section \ref{SS:finitenessdisc}, more specifically in Theorem \ref{T:GMCneumann_finite}.
\end{rmk}

The conformal covariance properties of the boundary and bulk Liouville measures are not quite as straightforward as for the Dirichlet GFF. The first problem is that conformal invariance of the Neumann GFF only holds when we view it as a distribution modulo constants. %  This forces us to view $\mathcal{V}$ as measure modulo global multiplicative constant.
The second is that we have only defined the boundary measure on the domain $\H$, where semicircles centred on the boundary can be defined. We could extend this definition to linear boundary segments of other domains, but it is unclear what to do when the boundary of the domain is very wild.

Let us start with the bulk measure, where we only need to deal with the first problem. In this case, the statement
\[ \mathcal{M}_h \circ T^{-1}  = \mathcal{M}_{h \circ T^{-1} + Q \log |(T^{-1})'|}\]  of \cref{T:ci} still holds (by absolute continuity with respect to the Dirichlet GFF) when $T:D\to D'$ is a deterministic, conformal isomorphism and we replace the Dirichlet GFF with $\sqrt{2\pi}$ times a Neumann GFF $h$ in $D$ with some arbitrary fixed additive constant. However, now $h\circ T^{-1}$ is  a Neumann GFF in $D'$ with a \emph{different} additive constant. The exact analogue of \cref{T:ci} only holds if we consider Neumann GFFs modulo additive constants, and their associated bulk Liouville measures modulo multiplicative constants (see exercises).

Now for the boundary measure,  suppose that $h$ is  a Neumann GFF on $\H$ with some arbitrary fixed additive constant, and $T:\H\to D$ is a conformal isomorphism. Then $h':=h\circ T^{-1}$ is  a  Neumann GFF on $D$ with another additive constant. Moreover, if $\partial D$ contains a linear boundary segment $L\subset \partial D\cap \R$, the measure $ \mathcal{V}_{h'}(\dd x)=\lim_{\eps\to 0} \e^{(\gamma/2)h'_\eps(x)} \eps^{\gamma^2/4} \, \dd x $ is well defined and
\begin{equation}
	\label{eq:cocoB}
	\mathcal{V}_h \circ T^{-1} = \mathcal{V}_{h \circ T^{-1} + Q \log |(T^{-1})'|}
	= e^{\gamma Q  \log |(T^{-1})' |} \mathcal{V}_{h'}.
\end{equation}
on $L$ with probability one.
In fact, by \cite[Theorem 4.3]{SheffieldWang}, the measure is well defined and the above formula holds with probability one for all conformal $T:\H \to D$ with $\partial D\cap \R \ne \emptyset $ simultaneously.

We will use this formula to \emph{define} the boundary Liouville measure for GFF-like fields on the conformal boundary\footnote{The conformal boundary of a simply connected domain $D$, equivalent to the Martin boundary (see \cite[\S 1.3]{SLEnotes}), is the set of limit points of $D$ with respect to the metric $d(x,y)=d(\phi(x),\phi(y))$ for $\phi:D\to \D$ a conformal isomorphism.} of an arbitrary simply connected domain.

\begin{definition}[Boundary Liouville measure for the GFF on $D$]\label{def:BLM}
	Suppose that $h$ is a random variable in $\cD'_0(D)$, and that for some conformal isomorphism $T:\H\to
	D$ the field $h\circ T+Q\log|T'|$ has the law of  a Neumann GFF (with some fixed additive constant) plus an almost surely continuous function on some neighbourhood in $\H$ of $L\subset \R$. Then the measure $\mathcal{V}_{h \circ T+Q\log|T'|}$ is almost surely well defined on $L$, and we may define
	\begin{equation}
		\label{eqn:bmeasure_coc}
		\mathcal{V}_{h}:= \mathcal{V}_{h \circ T+Q\log|T'|}\circ T^{-1}
	\end{equation}
	to be the Liouville measure for $h$, on the part of the conformal boundary of $D$ corresponding to the image of $L$ under $T$. With probability one, this defines the same measure simultaneously for all choices of $T$.
\end{definition}

Note that the behaviour of conformal isomorphisms near the boundary of a domain can be very wild. For instance if $D$ is a domain whose boundary is only H\"older with a certain exponent, then the boundary Liouville measure defined as above may not be easy to construct directly by approximation.  %for an appropriate field $h$ on $D$ corresponds to some kind of limit as $\eps\to 0$ of
%measures $\e^{\frac{\gamma}{2} (h,\sigma_{x,\eps})}\eps^{\gamma^2/4} dx$, but now $\sigma_{x,\eps}$ will be (roughly) a semicircle of radius $\eps^a$ around $x$, for some $a$ not necessarily equal to $1$.

% In a way this is what happens in the case of SLE. (On the other hand, if the domain is smooth enough that the conformal map extends across the boundary say, then the conformal covariance of the boundary Liouville measure, viewed as a measure up to global multiplicative constant, does hold).

\subsection{Finiteness of the GMC on a disc with Neumann boundary conditions}
\label{SS:finitenessdisc}

When studying quantum surfaces and the mating of trees approach to Liouville quantum gravity in the following chapters, the following question comes up naturally. Let $\mathbf{h}$ be a Gaussian free field with Neumann boundary conditions on the unit disc $\D$, normalised to have zero unit circle average (say) and let $h = \sqrt{2\pi} \mathbf{h}$. For $0 \le \gamma < 2$ let $\cM (\dd x)  = \lim_{\eps \to 0} \eps^{\gamma^2/2} e^{\gamma h_\eps(x) } \dd x$ be the bulk Gaussian multiplicative chaos measure (with circle average normalisation) associated with $h$. This limit exists in probability as $\eps \to 0$ for the topology of vague (as opposed to weak) convergence of measures, as a consequence of the standard GMC theory. Indeed, when we restrict $h$ to an open subset $A$ such that $\bar A \subset \D$, $h|_A$ is a logarithmically correlated Gaussian field in the sense of Section \ref{setup_gmc}. Theorem \ref{T:conv} thus shows that the above limit exists on any such subset $A$ of $\D$, from which vague convergence in probability easily follows. This defines unambiguously a GMC measure (which we denote by $\cM$) on $\D$.

From the above, it is clear that $\E[ \cM(A) ] < \infty$ for any compact subset $A \subset \D$. It is however not obvious if $\cM(\D) < \infty$ almost surely. An easy to way to be convinced that there is something nontrivial to show is to observe that $\E[ \cM(\D) ] = \infty$ as soon as $\gamma \ge \sqrt{2}$. 
The almost sure finiteness of this measure is the content of the next theorem (shown by Huang, Rhodes and Vargas in \cite{HRVdisk}). 

\begin{theorem}\label{T:GMCneumann_finite}
Let $0\le \gamma <2$ and let $\cM$ be the above GMC measure on the disc. Then $\cM(\D) < \infty$ almost surely, and {$\mathbb{E}[\cM(\D)^\alpha]<\infty$ for $\alpha=\min(1,1/\gamma)$.}
\end{theorem}

\begin{proof}
The proof below is partly inspired by \cite{HRVdisk} (but we use a different strategy for the key estimate \eqref{keyestimateNeumannsupercritical} below: indeed, we appeal to \cref{T:GMCsuper_poly}, based on a rough scaling argument, rather than through a comparison to branching random walk and the delicate work of Madaule \cite{Madaule}). We will show that there is some $0<\alpha<1$ such that $\E[ \cM(\D)^\alpha] < \infty$, which clearly implies the result. For $n\ge 1$, let $A_n$ be the dyadic annulus $A_n = \{ z \in \D:  2^{-n-1} < 1- |z| < 2^{-n} \}$. Since the boundary of $A_n$ has zero Lebesgue measure, it is easy to see that $\cM(\D) = \sum_{n\ge 1} \cM(A_n)$. Moreover, since $0\le \alpha<1$, the function $x>0 \mapsto x^\alpha$ is subadditive, hence
$$
\E[ \cM(\D)^\alpha] \le \sum_{n=0}^\infty \E [ \cM(A_n)^\alpha ]. 
$$
To estimate $ \E [ \cM(A_n)^\alpha ]$ we will make a comparison via Kahane's inequality (Theorem \ref{T:Kahane}) to a field defined on the unit circle and which is scale-invariant. However, in order to do this, it is important to first rewrite $\cM$ in terms of the variance normalisation of $h$, rather than in terms of the global factor $\eps^{\gamma^2/2}$ used in the statement of the theorem. For this, we observe from \cref{ex:NGFnu} that the covariance of $h$ at $x\ne y$ is equal to 
$$
2\pi G_N^{\D} (x,y) = - \log (|x-y ||1- \bar x y | ).
$$
Thus if $x \in \D$ and $B(x, \eps) \subset \D$ we see that
$$
\var h_\eps(x) =  \log (1/\eps) - \int_{\partial B (x, \eps)} \log (| 1- \bar x y |)\rho_{x, \eps} (\dd y ). 
$$
The above integral can be computed explicitly, since the function $z \mapsto \log |1- \bar x z| $ is harmonic in $\D$, and we find that
$$
\var h_\eps(x) = \log (1/\eps) - \log (1- |x|^2). 
$$
Consequently, for any open subset $A$ such that $\bar A \subset \D$, 
\begin{equation}\label{Mnormalisedvar}
\cM(A) = \lim_{\eps \to 0} \int_A e^{\gamma h_\eps(x) - \tfrac{\gamma^2}{2} \var (h_\eps(x) )} \frac1{(1- |x|^2)^{\gamma^2/2}}\dd x.
\end{equation}
Now, we fix $n\ge 1$ and consider $\cM(A_n)$. 
We let $(Y(x))_{x \in  \mathbb{U}}$ denote a rotationally invariant centered Gaussian field on the unit circle $\mathbb{U}$, with covariance $K_Y$ satisfying 
$$
K_Y(x,x') = - \log | x - x'| + O(1),
$$
where $|x-x'|$ denotes the distance between $x $ and $x'$ in $\R^2$ or equivalently (possibly up to a change of the implied constant $O(1)$ on the right hand side) their distance viewed as points on the unit circle. (An example of such a field is for instance provided by an appropriate multiple of the Gaussian free field with Neumann boundary conditions on $\D$, normalised to have zero unit circle average, and restricted to $\mathbb{U}$, see \cref{Ex:Neumann_Arc}).
Then an approximation to $Y$ at scale $2^{-n}$ can be compared to the field $h$ on $A_n$ in the following sense.

\begin{lemma}\label{lem:ngff0boundary}
There exists a constant $C>0$ independent of $n$ and $x,y \in A_n$ such that for $x = r e^{i \alpha}, y = s e^{i \beta} \in A_n$, if $Y_{2^{-n}}(z)$ denotes the average of $Y$ over an arc of circle of length $2^{-n}$ whose midpoint is $z \in \mathbb{U}$, 
\begin{equation}\label{eq:comparisoncircle}
%\E[ h(x) h(y) ] 
2\pi G_N^{\D} (x,y)= - \log (|x-y ||1- \bar x y | ) \ge 2 \E[ Y_{2^{-n}}(e^{i \alpha}) Y_{2^{-n}}(e^{i \beta})] - C.
\end{equation}
%for all $\alpha'\in(\alpha-2^{-n},\alpha+2^{-n})$, $\beta'\in(\beta-2^{-n},\beta+2^{-n})$.
\end{lemma}
\begin{proof}
Note that $\E[ Y_{2^{-n}}(e^{i \alpha}) Y_{2^{-n}}(e^{i \beta})] = - \log (|x-y|\vee 2^{-n}) + O(1) \le - \log |x-y| + O(1) $ so to show \eqref{eq:comparisoncircle} it suffices to show the inequality 
$$
- \log | 1- \bar x y | \ge -\log (|e^{i \alpha}-e^{i \beta} |\vee 2^{-n})  - C
$$
for some $C>0$ large enough, uniformly over $x,y \in A_n$, and over $n\ge 0$. 

{For this, we observe that 
	\begin{align*}
	|1-\bar{x}y|=|1-rse^{i(\beta-\alpha)}|& \le |1-e^{i(\beta-\alpha)}|+|e^{i(\beta-\alpha)}-rse^{i(\beta-\alpha)}| \\
	& = |e^{i\alpha}-e^{i\beta}|+|(1-r)(1-s)+r(1-s)+s(1-r)| \\
	& \le  |e^{i\alpha}-e^{i\beta}| + 5.2^{-n} \\
%	& \le  |e^{i\alpha'}-e^{i\beta'}| + 7.2^{-n} \\
	& \le 6( |e^{i\alpha}-e^{i\beta}|\vee 2^{-n}),
	\end{align*}
from which we conclude by taking logarithms.}
%A simple proof of this proceeds by contradiction. Supposing not, we would get $C_n \to \infty$ tending to infinity at least along a subsequence, and points $x_n = r_n e^{i \alpha_n}, y_n= s_n e^{i \beta_n} \in A_n$ such that 
%$$
%- \log | 1- \bar x_n y_n | \le -\log ( |e^{i \alpha_n}-e^{i \beta_n} |\vee 2^{-n})  - C_n
%$$
%By rotational invariance we can assume that $x_n = r_ne^{i \alpha_n}$ and $y_n =s_n  $ (i.e., $\beta_n = 0$), with $r_n, s_n \in [1- 2^{-n}, 1-2^{-n-1}]$ and $\alpha_n \in [-\pi, \pi]$. Then $|\alpha_n|\asymp |e^{i \alpha_n } - 1|$, so we would obtain that
%$$
%- \log | 1- r_n s_n e^{i \alpha_n} | \le - \log (| \alpha_n | \vee 2^{-n}) - C_n.
%$$
%We conclude by Taylor expansion after considering separately the case where $|\alpha_n |2^n \to 0$  (resp. $+\infty$, resp. is bounded), after extracting a further subsequence if necessary.
 \end{proof}
 From \eqref{eq:comparisoncircle} we deduce the following. For $x \in A_n$ and $\eps < 2^{- n-1}$, let $ \rho_{ \eps, x}$ be the uniform distribution on the circle of radius $\eps$ around $x \in A_n$. Define the field $\tilde Y_{\eps,2^{-n}}$ on $A_n$ by setting for $x \in A_n$,
 $\tilde Y_{\eps, 2^{-n}}(x) = \int Y_{2^{-n}} (e^{i \arg (x')}) \rho_{\eps, x} ( \dd x')$, and note that if $n$ is fixed but $\eps \to 0$, $\tilde Y_{\eps, 2^{-n} } (x) \to Y_{2^{-n}} (e^{ i \arg (x)} )$ almost surely uniformly over $A_n$, and this limit is a field which is purely angular (i.e., depends only on the argument). 
 Furthermore, by \cref{lem:ngff0boundary} there exists $C>0$ such that for all $n\ge 0$ and $\eps<2^{-n-1}$, and for all $x, y \in A_n$,
 $$
 \E [h_\eps(x) h_\eps(y) ]  \ge 2 \E[ \tilde Y_{\eps, 2^{-n}} (x) \tilde Y_{\eps, 2^{-n} } (y )] -C.
 $$
 Thus, applying Kahane's inequality (Theorem \ref{T:Kahane}) to the field $h_\eps(x) + \mathcal{N}$ (where $\mathcal{N}$ is an independent normal random variable) and the field $\sqrt{2}\tilde Y_{\eps, 2^{-n}} (x)$, we deduce that there exists $C>0$ such that for all $0< \alpha<1$ and for all $n\ge 1$, and for all $\eps< 2^{-n-1}$,
 \begin{align*}
 \E [ \left( \int_{A_n}  e^{\gamma h_\eps(x) - \tfrac{\gamma^2}{2} \E[ h_\eps(x)^2] } \dd x \right)^\alpha] & \le C \E [ \left(  \int_{A_n} e^{\sqrt{2} \gamma \tilde Y_{\eps, 2^{-n} } (x)  - \gamma^2 \E[\tilde Y_{\eps, 2^{-n}} (x )^2] } \dd x \right)^\alpha ].
% & \le C 2^{-n\alpha} \E [\left( \int_{\mathbb{U} } e^{\sqrt{2} \gamma \tilde Y_{\eps, 2^{-n} } (z)  - \gamma^2 \E[ \tilde Y_{\eps, 2^{-n}  } (z)^2 ] }\dd z  \right)^\alpha ] .
 \end{align*}
Letting $\eps \to 0$ in this inequality, and recalling (see \eqref{Mnormalisedvar}) that $\sigma (\dd x) = \dd x / ( 1- |x|^2)^{\gamma^2/2}  \le C 2^{n \gamma^2/2} \dd x$ on $A_n$,  we deduce that 
$$
 \E [\cM(A_n)^\alpha] \le C  2^{n \alpha \gamma^2/2}  \E \left[ \left( \int_{A_n} e^{  \sqrt{2} \gamma Y_{2^{-n} }( e^{i \arg (x) } )   - \gamma^2 \E[  Y_{ 2^{-n}  } (e^{i \arg (x) } )^2 ] } \dd x \right)^\alpha \right]  
$$
The integral inside the expectation on the right hand side is the mass on $A_n$ of a multiplicative chaos of a field which is purely angular. As a result, and by Fubini's theorem, we obtain 
$$
\E [\cM(A_n)^\alpha] \le C 2^{n \alpha \gamma^2/2}  2^{-n \alpha} \E \left[ \left( \int_{\mathbb{U}} e^{  \sqrt{2} \gamma Y_{2^{-n} }( z )   - \gamma^2 \E[  Y_{ 2^{-n}  } (z )^2 ] } \dd z \right)^\alpha \right]  
$$
where $\dd z$ denotes (with a small abuse of notation) Lebesgue measure on $\mathbb{U}$. 
%We deduce that 
% $$
% \E [ \cM(A_n)^\alpha] \le C 2^{ n\alpha  \tfrac{\gamma^2}{2} - n \alpha }  \E [\left( \int_{\mathbb{U} } e^{ \sqrt{2} \gamma Y_{2^{-n} } (z)  - \gamma^2 \E[ Y_{2^{-n}  } (z)^2 ] }\dd z  \right)^\alpha ] .
% $$
 When $\gamma < \sqrt{2}$ then we may simply take $\alpha = 1$ in the above; then the expectation in the right hand side is constant and $\sum_{n \ge 1} \E [ \cM(A_n) ] < \infty$ so $\cM(\D) < \infty$ almost surely. This no longer applies when $\gamma \in [\sqrt{2}, 2)$. Instead, we use the fact that when $\gamma >1$, $\sqrt{2}\gamma$ is greater than the critical value for the multiplicative chaos associated to the one dimensional field $Y$, which is $\sqrt{2d} = \sqrt{2}$. As a consequence\footnote{We note that this theorem was stated for Gaussian multiplicative chaos in $\R^d$ with respect to Lebesgue measure. However the reader can verify readily that this applies without any changes to Gaussian multiplicative chaos on the unit circle of $\R^2$ with respect to uniform measure, and more generally to smooth manifolds.} of \cref{T:GMCsuper_poly} (with $d=1$ and $\gamma$ replaced by $\sqrt{2}\gamma$), taking $$
 	\alpha=\sqrt{2d}/(\sqrt{2}\gamma)=1/\gamma<1,$$ we have that %$\E [\left( \int_{\mathbb{U} } e^{ \sqrt{2} \gamma Y_{2^{-n} } (z)  - \gamma^2 \E[ Y_{2^{-n}  } (z)^2 ] } \dd z \right)^\alpha ]  $ %decays like a certain power of $2^{-n}$. %Unfortunately the estimate in \cref{T:GMCsuper_poly} is not optimal and this is not sufficient to show that $\sum_n \E[ \cM(A_n)^\alpha] < \infty$. Instead, we appeal to the estimate mentioned in \cref{R:Freezing}, coming from the work of Madaule, Rhodes and Vargas \cite{MadauleRhodesVargas}, 
% which implies in particular that for $\alpha < \gamma^{-1}$ and $\gamma >1$, 
 \begin{equation}\label{keyestimateNeumannsupercritical}
  \E [\left( \int_{\mathbb{U} } e^{ \sqrt{2} \gamma Y_{2^{-n} } (z)  - \gamma^2 \E[ Y_{2^{-n}  } (z)^2 ] } \dd z \right)^\alpha ] \le C 2^{- n \alpha ( \gamma - 1)^2 }.
 \end{equation}
for $C$ not depending on $n$.
 Thus %if $\gamma >1$ and $\alpha < \gamma^{-1} < 1$, 
 we obtain that
 $$
 \E [ \cM(A_n)^\alpha] \le C 2^{n\alpha ( \tfrac{\gamma^2}{2} - 1 - (\gamma-1)^2 ) }.
 $$
 Since $ \tfrac{\gamma^2}{2} - 1 - (\gamma-1)^2 <0$ for all $0\le \gamma < 2$, this shows that  $\E[\cM(\D)^\alpha]\le \sum_n \E[ \cM(A_n)^\alpha] < \infty$ and hence $\cM(D) < \infty$ almost surely. This concludes the proof of \cref{T:GMCneumann_finite}. 
\end{proof}

\subsection{Exercises}

\begin{enumerate}[label=\thesection.\arabic*]

\item Let $D= (0,1)^2$ be the unit square. Find an orthonormal basis of $L^2(D)$ consisting of eigenfunctions of $-\Delta$ in $D$, with \emph{Neumann} boundary conditions, and write down their eigenvalues. 
Now consider the setting of \cref{P:GFFsquare}, and set $V_N=D\cap(\mathbb{Z}^2/N)$ for $N\ge 1$. Come up with a definition of the discrete Neumann GFF in $V_N$ with Neumann boundary conditions, and prove that it converges as $N\to \infty$ to a continuum GFF in $D$ with Neumann boundary conditions in a suitable sense.

\item Consider the Hilbert space completion $(\bar{H}_{\C},(\cdot, \cdot)_\nabla)$ of the set of smooth functions modulo constants in $\C$ with finite Dirichlet norm.
Let $$
\bar{H}_{\text{even}} = \{ h \in \bar{H}_{\C}: h(z)-h(0) = h( \bar z)-h(0), z \in \C\}
$$
and likewise let
$$
\bar{H}_{\text{odd}} = \{ h\in \bar{H}_{\C} : h(z)-h(0)=-(h(\bar{z})-h(0)), z \in \C\}
$$
(note that $h(z)-h(0)$ is well defined for a function modulo constants).
Show that $\bar{H}_{\C} = \bar{H}_{\text{even}} \oplus \bar{H}_{\text{odd}}$. (Hint: orthogonality follows from the change of variables $z\mapsto \bar z$).

Show that the series	$$\sum_n X_n \bar{f}_n,$$ where $X_n$ are i.i.d. standard normal random variables, and $\bar{f}_n$ is an orthonormal basis of $\bar{H}_{\C}$, converges almost surely in the space $\bar\cD_0'(\C)$ and the limiting distribution modulo constants agrees in law with the whole plane GFF, $\bbh^\infty$.

\item Prove \eqref{E:WPGFFcov} using \cref{D:WPGFF} of the whole plane GFF, and the explicit expressions for the Neumann and Dirichlet Green functions in $\H$.

\item Give a rigorous definition (that is, as a random distribution) of the whole plane GFF with additive constant fixed to have average $0$ with respect to a Riemannian volume form $g(z) \dd z$ on the Riemann sphere $\C\cup \{\infty\}$, in a manner analogous to \cref{D:NGFFnorm}. Show that this is equal in law to the spherical GFF $\mathbf{h}^{\Sp,g}$.

\item Write down a definition of the Dirichlet--Neumann GFF in the upper unit semidisc $\D_+$ as a stochastic process, giving an explicit expression for its covariance function in the upper unit semidisc.

\item \label{Ex:GMCboundary} Give a complete proof of \cref{T:boundary}, using the same strategy as in \cref{S:Liouvillemeasure}. Explain briefly why \cref{T:conv} does not apply directly to this setting.

\item Prove \eqref{eq:cocoB} (see the proof of  \cref{T:ci}).
Check that the boundary Liouville measure $\nu$ satisfies the same KPZ relation as the bulk Liouville measure.

\item Let $\mathcal{V}^\sharp$ be the the boundary Liouville measure for a Neumann GFF on $\H$ with some fixed choice of additive constant, restricted to (0,1), and renormalised so that it is a probability distribution. Sample $x$ from $\mathcal{V}^\sharp$. Is the point $x$ thick for the field (in terms of semi-circle averages)? If so, how thick?

\end{enumerate}

\newpage

\section{Quantum wedges and scale-invariant random surfaces}\label{S:SIsurfaces}

% !TEX root = master.tex

\subsection{Convergence of random surfaces}
\label{sec:rs_conv}

%The Neumann GFF is conformally invariant (and in particular scale invariant) when viewed as a distribution modulo constants. However when we fix a normalisation this is no longer true, and so the abstract surface described by $(\H, h)$ is not scale invariant (in the sense that applying the conformal map $z\mapsto rz$ and applying the conformal change of coordinates formula does not yield the same surface {(do you mean field -- by definition they're giving the same surface?)} in distribution). {I'm not sure this is the right notion of scale invariance - should it rather be that you want to apply the conformal map $z\mapsto rz$, apply the change of coordinates formula \emph{and} add a constant so that the LQG mass of some fixed set remains the same, and then get a field with the same law?}  In order to construct a surface $(\H, h)$ which is invariant under scaling, Sheffield introduced the notion of \emph{quantum wedge}, which will play an important role in \change{our study of the ``quantum gravity zipper'' in next section}. Roughly speaking, a quantum wedge is the surface that one obtains by ``zooming in'' close to a point at the boundary. Since this surface is obtained as a scaling limit it is then automatic that it will be invariant under scalings.

\textbf{Note}: from this point onwards, we will almost exclusively work with the multiplicative normalisation $h=\sqrt{2\pi}\mathbf{h}$ as in \eqref{hNLQG} for the Neumann GFF and its variants.\\

Recall that we defined a \emph{random surface} to be an equivalence class of pairs $(D,h)$ where $D$ is a simply connected domain and $h$ is a distribution on $D$, under the relation identifying $(D_1,h_1)$ and $(D_2,h_2)$ if for some $f: D_1\to D_2$ conformal, $$h_2=h_1\circ f^{-1}+Q\log |(f^{-1})'|.$$ The reason for this was that if $h_1$ is a Dirichlet Gaussian free field in $D_1$, then all members of the equivalence class of $(D_1,h_1)$ describe the same Liouville measure up to taking conformal images.

Now, we have seen that the same thing is true when $h_1$ is  a Neumann GFF with an arbitrary fixed additive constant. And indeed if we want to view the Neumann GFF as a quantum surface then we have to  fix the additive constant, since the definition of quantum surface involves distributions and not distributions modulo constants. But the Neumann GFF is only really \emph{uniquely} defined as a distribution modulo constants.  This  manifests itself in the following problem: different ways of fixing the additive constant do not yield the same quantum surface in law (see Example \ref{ex:zooming} below). So if we want to view the Neumann GFF as a quantum surface, which way of fixing the additive constant should  we pick? The lack of a canonical answer to this suggests that, at least when working with quantum surfaces, it is perhaps more natural to look at a slightly different object.

Another point of view is the following: if we consider a Neumann GFF $h$ with some arbitrary fixed additive constant, and also the field $h+C$ for some $C$, then the Liouville measure for $h+C$ is just $e^{\gamma C}$ times the Liouville measure for $h$. So we can think that the quantum surface described by $h+C$ represents ``zooming in'' on the quantum surface defined by $h$. (Note that this is distinct from rescaling space by a fixed factor and applying the change of coordinate formula, since by definition this does not change the quantum surface).
For some purposes, it will be natural to work with quantum surfaces that are invariant (in law) under such a zooming operation. Such a property can be thought of as a type of \textbf{scale invariance for quantum surfaces}.

In order to construct a surface $(\H, h)$ which does have this invariance property (once again, by Example \ref{ex:zooming} below this is not true when $h$ is a Neumann GFF with some arbitrary fixed additive constant), Sheffield \cite{zipper} introduced the notion of \textbf{quantum wedge}. This will play an important role in our study of the \emph{quantum gravity zipper} in next chapter. Roughly speaking, a quantum wedge is the limiting surface that one obtains by ``zooming in'' to a Neumann GFF close to a point on the boundary. Since this surface is obtained as a scaling limit, it automatically satisfies the desired scale invariance. Later on we will also study scale invariant quantum surfaces without boundaries (quantum cones) and finite volume versions of both wedges and cones (namely, so called quantum discs and spheres).

In order to make proper sense of the above discussion, we first need to provide a notion of convergence for random surfaces -- and more precisely, for surfaces with marked points.% on their boundaries. %({(?)} which will be 0 and $\infty$ in the case of the upper half plane). %Recall first our definition of an abstract surface -- a pair ($D,h)$, where $D$ is a domain and $h$ a distribution in $D$, modulo the equivalence relation $(D,h) \sim (D',h')$ iff there exists a conformal map $\phi:D\to D'$ such that $h' = h \circ \phi^{-1} + Q \log |( \phi^{-1})'|$, where $Q = \gamma/2 + 2/\gamma$. Note here that $h$ really is a distribution, and not a distribution modulo constants, so if we want to view the Neumann GFF as such an object we must first pick a particular normalisation.

\begin{definition}[Quantum surface with $k$ marked points]\label{def:markedqs}
A quantum surface with $k$ marked boundary points is an equivalence class of tuples $(D,h,x_1,\cdots, x_k)$ where $D\subset \C$ is a domain, $h\in \mathcal{D}_0'(D)$, and $x_1,\cdots, x_k$ are points on the (conformal) boundary or in the interior of $D$, under the equivalence relation
 $(D,h,x_1,\cdots, x_k) \sim (D',h',x_1',\cdots, x_k')$ if and only if for some $T: D\to D'$ conformal with $T(x_i)=x_i'$ for $1\le i \le k$ (note that $T$ extends to a to map between conformal boundaries by definition):
\begin{equation}
\label{eq:concoc}
h'=h\circ T^{-1}+Q\log |(T^{-1})'|.
\end{equation}
\end{definition}
\ind{Quantum surface}

We recall that $Q=Q_\gamma=2/\gamma+\gamma/2$ depends on the LQG parameter $\gamma$, and therefore so does the notion of \emph{quantum surface}, but we  drop this from the notation for simplicity. Note that since $h$ is assumed to be in the space of distributions $\mathcal{D}_0'(D)$, this definition may be applied to a Neumann GFF with an arbitrary fixed additive constant.

\vspace{.2cm}
In order to define a quantum surface $\mathcal{S}$ with $k$ marked points, we need only specify a single equivalence class representative $(D,h,x_1,\cdots, x_k)$. We will call such a representative an \textbf{embedding} or \textbf{parametrisation} of the quantum surface.

This means that our usual topology on the space of distributions induces a topology on the space of quantum surfaces (with $k$ marked points).
\begin{definition}[Quantum surface convergence] \label{def:qsconv} A sequence of quantum surfaces $\mathcal{S}^{n}$ converges to a quantum surface $\mathcal{S}$ as $n\to \infty$ if there exist representatives $(D,h^n,x_1,\cdots, x_k)$ of $\mathcal{S}^n$ and $(D,h,x_1,\cdots, x_k)$ of $\mathcal{S}$, such that $h_n \to h$ in the space of distributions as $n\to \infty$.
\end{definition}

(We note that this notion of convergence is somewhat different from the notions used in \cite{zipper} or \cite{DuplantierMillerSheffield}, but this definition has the advantage that it makes sense for all deterministic distributions viewed as quantum surfaces rather than a special class of random ones. It is also, in any case, the one that actually used to verify convergence statements for quantum surfaces, as will be discussed below.) \\
\index{Random surface!Convergence|BH}
%\ind{Random surface!Half-plane like}
%\ind{Half-plane like! see {Random surface}}

%Such an abstract surface may be embedded in the domain $D = \H$ with the two marked points being $z_0 = 0$, $z_1 = \infty$. Such an embedding is not unique: indeed, the set of conformal maps mapping $\H$ to itself and mapping 0 and $\infty$ to themselves consists precisely of the dilations $z \mapsto rz$.

Now, when we are actually working with quantum surfaces, it will often be very useful to specify a surface by describing a particular \emph{canonically chosen} embedding. Of particular interest are \emph{random} surfaces (like the Neumann GFF or the quantum wedges defined below), and this allows for certain special choices of embedding (we will see several in the rest of this chapter and the next).

\begin{example}\label{ex:cd}
	Suppose that $h$ is equal to a continuous function plus a Neumann GFF (with some  fixed additive constant) in $D$ simply connected, and $z_0,z_1\in \partial D$ are such that the \emph{bulk} Liouville measure $\cM_h$ for $h$ assigns finite mass to any finite neighbourhood of $z_0$, and infinite mass to any neighbourhood of $z_1$.\footnote{If $h$ is just a Neumann GFF with arbitrary fixed additive constant in an unbounded domain $D$, then this will be the case by \cref{T:GMCneumann_finite} whenever $z_1=\infty$ and $z_0$ is another $(\ne \infty)$ boundary point where the boundary is smooth (say).} Then the doubly marked quantum surface $(D, h, z_0, z_1)$ has a unique representative  $(\H, \tilde h, 0 , \infty)$ such that  $\cM_{\tilde h} (\D\cap \H) = 1$. The distribution $\tilde h$ is called the \textbf{canonical description} of the quantum surface in \cite{zipper}.\footnote{but bear in mind that it is only well defined when $h$ is in a particular class of distributions for which the Liouville measure makes sense.} In fact, in practice this is a difficult embedding to work with and we usually prefer others; this will be discussed further in the following section.
\end{example}
\ind{Canonical description|see {Random surface}}
\ind{Random surface!Canonical description}
%In other words, for a half-plane like surface, we fix its embedding in $\H$ by choosing the dilation parameter $r$ so that $\D$ has Liouville mass 1 (such a choice is always possible by our definition of half-plane like surface).

\ind{Random surface!Zooming in}
\begin{example}[Zooming in -- important!]\label{ex:zooming} Let $h$ be a Neumann GFF in $\H$, for concreteness, normalised to have average zero in $\D\cap \H$. Then the canonical descriptions of $h$ and of $h+100$ (say), viewed as quantum surfaces in $\H$ with marked points at $0$ and $\infty$, are \emph{very different}. This can be confusing at first, since $h$ is in some sense defined ``up to a constant'', but the point is that ``equivalence as quantum surfaces'' and ``equivalence as distributions modulo constants'' \emph{are not the same}.

Indeed to find the canonical description of $h$ we just need to find the (random) $r$ such that $\mathcal{M}_h( B(0,r)\cap \H)= 1$, and apply the conformal isomorphism $z\mapsto z/r$; the resulting field $$
\tilde h(z) = h(rz) + Q \log (r)$$ defines the canonical description $\tilde h$ of the surface $(\H, h,0,\infty)$. On the other hand, in order to find the canonical description of $h+100$, we need to find $s>0$ such that $\mathcal{M}_{h+100} (B(0,s)\cap \H) = 1$. That is, we need to find $s>0$ such that $\mathcal{M}_h(B(0,s)\cap \H) = e^{-100 \gamma}$. The resulting field
$$
h^*(z) = h (sz) + Q \log (s) +100
$$
defines the canonical description of $(\H, h+100,0,\infty)$.

Note that in this example, the ball of radius $s$ is much smaller than the ball of radius $r$. Yet in $\tilde h$, the ball of radius $r$ has been scaled to become the unit disc, while in $h^*$ it is the ball of radius $s$ which has been scaled to become the unit disc. In other words, and since $s$ is much smaller than $r$, the surface $(\H,h+100,0,\infty)$ is obtained by taking the surface $(\H, h,0,\infty)$ and \textbf{zooming in} at $0$. %Hence, even when we consider a particular distribution modulo constants $\bar h \in \bar \cD(D)$ and $h \in \cD(D)$ is a specific normalisation of $\bar h$, then the canonical descriptions of $h$ and $h+100$ are not the same, even though both define the same element of $\bar \cD(D)$.
\end{example}
%\medskip Having fixed a canonical description (or canonical embedding), we can formulate a notion of convergence of surfaces.

%\begin{definition}\label{D:conv}
%	We say that a sequence of half-plane like surfaces $(D_n, h_n, z_0^n, z_1^n)$ (where $h_n \in \cD(D_n)$) converges to a half-plane like surface $(D, h, z_0, z_1)$ if in the canonical description of these surfaces $(\H, h_n, 0, \infty)$ and $(\H, h, 0, \infty)$ we have that the measures $\mu_{h_n}$ converge weakly towards $\mu_h$.
%\end{definition}

\subsection{Thick quantum wedges}

\ind{Quantum wedges}
As we will see very soon, a (thick) quantum wedge is the abstract random surface  that arises as the $C\to \infty$ limit %(in the sense discussed above)
of the doubly marked surface $(h + C,\H,0,\infty)$, %when $C\to \infty$,
when $h$ is a Neumann GFF in $\H$ with some fixed additive constant plus certain logarithmic singularity at the origin.
% (with any normalisation) and the marked points are $0$ and $\infty$.
Thus, as explained in the example above, it corresponds to zooming in near the origin of $(h,\H,0,\infty)$.

In practice however, we prefer to work with a concrete definition of the quantum wedge and then prove that it can indeed be seen as a scaling limit. It turns out to be most convenient to define it in the infinite strip $S = \R \times (0,\pi)$ rather than the upper half plane, with the two marked boundary points being $+\infty (=:\infty)$ and $-\infty$ respectively. A conformal isomorphism transforming $(S, \infty, - \infty)$ into $(\H, 0, \infty)$ is given by $z \mapsto - e^{-z}$, and under this conformal isomorphism, vertical line segments are mapped to semicircles. To be precise, the segment $\{z: \Re(z) = s\}$ is mapped to $\partial B(0, e^{-s}) \cap \bar{\H}$ for every $s\in \R$.

The following lemma will be used repeatedly in the rest of this chapter and the next.

\begin{lemma}[Radial decomposition]\label{L:radialdecomposition}
	Let $S$ be the infinite strip $S = \{ z = x + i y\in \C: y \in (0, \pi)\}$.
	Let $\bar \cH_{\rr}$ \indN{{\bf Function spaces}! $\bar{\cH}_{\rr}$; closure of smooth functions modulo constants on the infinite strip (resp. cylinder) $S$ (resp. $\cC$) which are on vertical segments} be the subspace of $\bar H^1(S)$ obtained as the closure of smooth functions which are constant on each vertical segment, viewed modulo constants. Let $\cH_{\cc}$  \indN{{\bf Function spaces}! $\bar{\cH}_{\cc}$; closure of smooth functions  on the infinite strip (resp. cylinder) $S$ (resp. $\cC$) which have mean zero  on vertical segments}be the subspace obtained as the closure of smooth functions which have mean zero on all vertical segments. Then
	$$\bar H^1(S) = \bar \cH_{\rr} \oplus \cH_{\cc}.$$
\end{lemma}

\begin{proof}
	Suppose that $g_1$ is a smooth function modulo constants in $S$, that is constant on vertical lines, and that $g_2$ is a smooth function in $S$ that has mean zero on every vertical line. Then it is straightforward to check that $\iint_S \nabla g_1 \cdot \nabla  g_2 = 0$. Indeed $\nabla g_1=(\partial_{x} g_1, 0)$ and $\nabla g_2 = (\partial_{x} g_2, \partial_{y} g_2)$ where the partial derivative $\partial_{x} g_1$ is constant on vertical lines and $\partial_{x} g_2$ has average $0$ on vertical lines. This means that $\nabla g_1 \cdot \nabla g_2$ has average $0$ on every vertical line, and consequently has average 0 over $S$. By definition of  $\bar{\cH}_{\rr}$ and $\cH_{\cc}$ (as closures with respect to $(\cdot, \cdot)_\nabla$) the two spaces are therefore orthogonal with respect to $(\cdot, \cdot)_{\nabla}$.
	
To check that they span $\bar{H}^1(S)$, note that if we consider smooth $f\in \bar\cD(S)$ and we set $f_{\rr}(z)$ to be the average of $f$ on the line $\Re z + i[0,2\pi]$, then $f_{\rr}\in \bar{\cH}_{\rr}$. Moreover, defining $f_{\cc}=f-f_{\rr}$, it is clear that $f_{\cc}\in \cH_{\cc}$. From this it follows that if $f\in \bar{H}^1(S)$ then we can write $f=\lim_n f_n $ for a sequence $(f_n)_n\in \bar\cD(S)$, and by decomposing each $f_n$ we have  $\lim_n f_n= \lim_n ((f_n)_{\rr}+(f_n)_{\cc})$. By orthogonality, the sequences $(f_n)_{\rr}$ and $(f_n)_{\cc}$ are each Cauchy and have individual limits $f_{\rr}\in \bar{\cH}_{\rr}$ and $f_{\cc}\in \cH_{\cc}$. Hence $f=f_{\rr}+f_{\cc}$, and the two spaces do indeed span $\bar{H}^1(S)$.
\end{proof}

Similarly to the domain Markov property for the Neumann GFF (that we saw arises from the orthogonal decomposition $\bar{H}^1(D)=H_0^1(D)\oplus \overline{\Harm}(D)$), this results in another representation of the Neumann GFF on $S$ modulo constants. Namely, as a stochastic process indexed by $\tilde{\mathfrak M}_N^S$, it %normalised so that its average on $(0, i\pi)$ is zero,
can be written as $\bar{h} = \bar{h}_{\rr}^{S} + h_{\cc}^{S}$ where:
\begin{itemize}
	\item $\bar{h}_{\rr}^{S}, h_{\cc}^{S}$ are independent;
	\item  $\bar{h}_{\rr}^{S}(z) = \bar{B}_{2\Re(z)}$, where $\bar{B}$ is a standard Brownian motion modulo constants (by \cref{T:BM} and conformal invariance);
	\item $h_{\cc}^{S}(z)$ has mean zero on each vertical segment.
\end{itemize}

To justify the above, notice that given $\bar{h}$, we can define $\bar{h}^{S}_{\rr}$ to be constant on each vertical segment with value equal to the average of $\bar{h}$ on that segment. Then we know by \cref{T:BM} and conformal invariance that $\bar{h}^{S}_{\rr}$ has the law described in the second bullet point. Thus it remains to justify is that $(\bar{h}-\bar{h}^{S}_{\rr},\rho)$ and $(\bar{h}^{S}_{\rr},\rho)$ are independent for any $\rho\in \tilde{{\mathfrak M}_N^S}$. For this, observe that if $(\bar{h}_n)_n$ are as in \cref{T:NGFFseries} (but multiplied by $\sqrt{2\pi}$) then $(\bar{h}_n,\rho)$ converges in $L^2(\P)$ and in probability, to a random variable with the law of $(\bar{h},\rho)$. Moreover $((\bar{h}_n)_{\rr},\rho)$ and $(\bar{h}_n-(\bar{h}_n)_{\rr},\rho)$ are independent for every $n$, with $\var((\bar{h}_n)_{\rr},\rho)\le \var((\bar{h}^{S}_{\rr},\rho))$ and $\var(\bar{h}_n-(\bar{h}_n)_{\rr},\rho)\le \var(\bar{h}-\bar{h}^{S}_{\rr},\rho)$. This implies that $(\bar{h}-\bar{h}^{S}_{\rr},\rho)$ and $(\bar{h}^{S}_{\rr},\rho)$ are uncorrelated and hence, by Gaussianity, independent.

Note that the $\bar{h}_{\rr}^{S}$ part is defined modulo constants, while the $h_{\cc}^{S}$ part really has additive constant fixed. As such, we can actually define $h_{\cc}^{S}$ to be a stochastic process indexed by ${\mathfrak M}_{N}^S$ rather than just $\tilde{\mathfrak M}_{N}^S$.\footnote{Concretely, we can set $h_{\cc}^{S}$ to be $h-h_{\rr}$ where $h$ is $\bar{h}$ with additive constant fixed so that it's average on $(0,i\pi)$ is zero, as in \cref{D:NGFFnorm}, and $h_{\rr}$ is constant on each vertical segment with value equal to the average of $h$ on that segment.}

Also observe that all the roughness of $h$ is contained in the $h_{\cc}^{S}$ part, as $\bar{h}_{\rr}^{S}$ is a nice continuous function modulo constants. %(which is even constant on vertical segments).
On the upper half plane, this would correspond to a decomposition of $\bar{h}$ into a part which is a radially symmetric continuous function (modulo constants), and one which has zero average on every semicircle (hence the notation). \ind{Radial decomposition}% We will often fix the normalisation for the Neumann GFF on $S$ by declaring that $B_0=0$.

\begin{rmk}[Translation invariance of $h^{S}_{\cc}$]\label{rmk:tiNGFF}
Note that the Neumann GFF $h$ on $S$ is invariant under horizontal translations (modulo constants), as it is conformally invariant (modulo constants). %As a result, both radial and circle parts are translation invariant too (modulo constants).
Since the radial part is simply a two-sided Brownian motion, the translation invariance of this part modulo constants is also clear. Thus, we may deduce that the circular part $h_{\cc}^{S}$ is translation invariant as well. (Note that the additive constant here is specified).
%Since it is only $h^{\textrm{GFF}}_{\rr}$ that is defined modulo constants, we have a unique covariance function for $h^{\textrm{GFF}}_{\cc}$. This can be calculated from any choice of covariance function for $h$, and using $G^S(z, w)=-\log (|e^{-z}-e^{-w}||e^{-z}-\overline{e^{-w}}|)$ (see Example \ref{ex:NGFnu2}) we obtain that $h_{\cc}$ has covariance kernel
%$$
%G^{S, \cc}(z,w)=-\log(|e^{-\Im z}+e^{-\Im w}||e^{-\Im z}+e^{\Im w}|).
%$$ In particular, the law of $h^{\textrm{GFF}}_{\cc}$ is stationary under translations of $S$. As before, $h^{\textrm{GFF}}_{\cc}$ defines a stochastic process indexed by $\cM=\cM^S$.
\end{rmk}

Let $ 0\le \alpha \le Q=2/\gamma+\gamma/2$. We will define an \textbf{$\alpha$-(thick) quantum wedge} to be a quantum surface $(S,h,+\infty,-\infty)$, where the law of the representative field $h$ on $S$ will be defined by specifying, separately, its averages on vertical line segments, and % on $\bar \cH_{\rr}$ and
what is left when we subtract these. The second of these components will be an element of  $\cH_{\cc}$, %We will actually specify the particular normalisation for the $\bar \cH_{\rr}$ part where the mean on $(0, i \pi)$ is required to be zero, corresponding to a set of functions which we call $\cH_{\rr}$.
having exactly the same law as the corresponding projection $h_{\cc}^{S}$ of the standard Neumann GFF. It is only the ``radially symmetric part'' which is different.

The bound $\alpha\le Q$ corresponds to the fact that we are defining a so called ``thick'' quantum wedge. When $\alpha>Q$ it is possible to define something called a ``thin'' quantum wedge, as introduced in \cite{DuplantierMillerSheffield}, but we will discuss this separately later on.

\begin{definition}
	\label{D:wedge}
	Let
	\begin{equation}\label{E:wedge_radial}
	h^{\mathrm{wedge}}_{\rr}(z) =
	\begin{cases}
	B_{2s} + (\alpha - Q) s & \text{ if } \Re(z) = s \text{ and } s \ge 0\\
	\widehat B_{- 2s} + (\alpha - Q) s & \text{ if } \Re(z)=s \text{ and } s < 0
	\end{cases}
	\end{equation}
	where $B=(B_t)_{t\ge 0}$ is a standard Brownian motion, and $\widehat B=(\widehat B_t)_{t\ge 0}$ is an independent Brownian motion conditioned so that $\widehat B_{2t} + (Q - \alpha) t >0$ for all $t>0$ (see below for what this means precisely).
	
	Let $h^{\mathrm{wedge}}_{\cc}$ be a stochastic process indexed by $\mathfrak{M}_N^S$, that is independent of $h^{\mathrm{wedge}}_{\rr}$ and has the same law as $h^{S}_{\cc}$.
Finally, set $h^{\mathrm{wedge}} = h^{\mathrm{wedge}}_{\rr} + h^{\mathrm{wedge}}_{\cc}$ (which, since $h^{\mathrm{wedge}}_{\rr}$ is just a continuous function, can again be defined as a stochastic process indexed by $\mathfrak{M}_N^S$). We call $h^{\mathrm{wedge}}=h^{\mathrm{wedge}}_{\rr}+h^{\mathrm{wedge}}_{\cc}$ the \textbf{$\alpha$-quantum wedge} field in $(S,+\infty,-\infty)$. 

The $\alpha$-quantum wedge itself is defined to be the doubly marked quantum surface represented by $(S,h^{\mathrm{wedge}},+\infty,-\infty)$ (see \cref{R:wedge_RV} below).
\end{definition}

The conditioning defining the process $\hat{B}$ above bears on an event of probability zero and thus requires a careful definition. For instance, when $\alpha=Q$, this is the limit (as $C \to \infty$ and then $\eps \to 0$) of a (speed two) Brownian motion starting from $\eps$ condition to hit $C$ before zero. This turns out to be identical to a (speed two) Bessel process of dimension $3$ (this classical fact can be proved using standard techniques in stochastic calculus and in particular Doob's $h$-transform; for a proof see, e.g., Lemma 8.2 in \cite{SLEnotes}).
On the other hand, if $\alpha < Q$, then the process $\hat B$ can be defined as the limit, as $\eps\to 0$, of a (speed two) Brownian motion with drift $(Q-\alpha)$, started from $\eps>0$ and conditioned to stay positive for all time. As observed in particular by Williams \cite{Williams}, under this conditioning, $(\hat B_t + (Q - \alpha) t)_{t\ge 0} $ (equivalently, $h_{\rr} (z)$ for $\Re(z) = - t$, viewed as a function of $t\ge 0$), is a strong Markov process with generator $ \frac{\dd^{{\hspace{.05cm}2}}}{\dd x^2} + b \coth( b x) \frac{\dd}{\dd x}$, where $b = Q - \alpha>0$. This again follows from Doob's $h$-transform after observing that $x \mapsto \phi(x) = e^{-2bx}$ is harmonic for Brownian motion with drift $-b$. 

However, a simpler definition of $h^{\mathrm{wedge}}_{\rr} (z)$ for $\Re(z) = t$ which works for \emph{all} $t\in \R$ whether positive or negative, (and works for all $\alpha \le Q$, including $\alpha = Q$), is the following. Fix $C>0$, and let $Y_t$ be a Brownian motion with speed two, starting from $C$ at time 0, and with drift $- b = - (Q- \alpha)\le 0$. That is, $Y_t = B_{2t} - bt$. Let $\tau = \inf \{ t \ge 0: Y_t = 0\}$. Then we set (for $ \Re(z) = t$) 
\begin{equation}\label{hradwedge}
h^{\mathrm{wedge}}_{\rr} (z): = \tilde Y_t ; \text{ where } \tilde Y_t = \lim_{C \to \infty} Y_{\tau + t}.
\end{equation} 
As this procedure provides a consistent definition of $(Y_{\tau+t })_{t\in \R}$ as $C$ increases (thanks to the strong Markov property of $Y$), it is obvious that this limit exists. Note that this provides a definition which works for all $t \in \R$ (i.e., there is no need to distinguish between $t\ge 0$ and $t \le 0$). We will refer to $\tilde Y$ as Brownian motion with speed two and drift $-b = - (Q-\alpha)$, from $+ \infty$ to $- \infty$ (although this terminology is a slight abuse of language when $\alpha = Q$). It turns out to give the same process as the one described above by the conditioning using Doob's $h$-transform; however this equivalence is not obvious. Instead this is a consequence of William's path reversal theorem \cite[Theorem 2.5]{Williams}, which implies that the time-reversal of a Brownian motion with drift $- b$, started from a value $C>0$, considered until its first hitting time of 0, has the same law as the process with generator $\frac12 \frac{\dd^{\hspace{.05cm}2}}{\dd x^2} + b \coth( b x) \frac{\dd}{\dd x}$, considered until its last hitting time of $C$ (the extra factor 1/2 in the generator comes from the speed of Brownian motion). 
In later arguments we will see that using the definition \eqref{hradwedge} leads to simple proofs, and these do not rely on the above mentioned equivalence, so we bypass in this manner the use of results requiring delicate stochastic calculus techniques. This is therefore the definition we adopt.

We warn the reader that the process $\tilde Y$ defined in \eqref{hradwedge} is somewhat different (even modulo horizontal translation) from a two-sided Brownian motion with drift $-b$ (i.e., $(B_{2t} - bt)_{t\in \R}$, where $(B_t)_{t\in \R}$ is a two-sided Brownian motion. Indeed, the latter has a local time at zero given by the sum of \emph{two} independent exponential random variables (corresponding to the local time accumulated in the negative and positive times respectively), whereas $\tilde Y$ has a local time at the origin which is a single exponential random variable. Essentially, the difference comes from the fact that the latter can be viewed as the process $\tilde Y$ biased by its local time at the origin. 
 
 \medskip
To emphasise once more, our definition of quantum wedge fields is such that they come with a specific way of fixing the additive constant; in other words, they are stochastic processes indexed by $\mathfrak{M}_N^S$ rather than just $\tilde{\mathfrak{M}}_N^S$. We will \emph{not} want to consider these wedge fields modulo constants.

\begin{figure}\begin{center}
		\includegraphics[scale=.8]{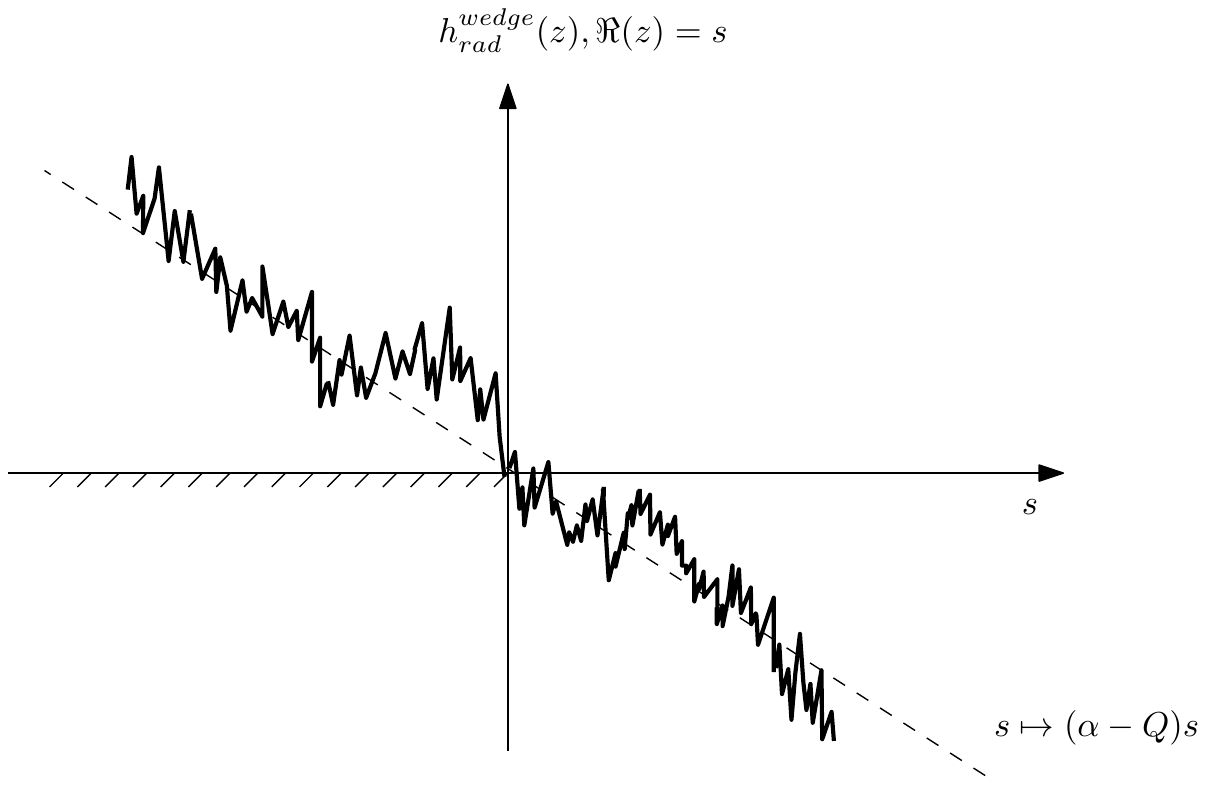}
	\end{center}
	\caption{Schematic representation of the radially symmetric part of a quantum wedge in a strip. When $s<0$, the function is conditioned to be positive.}
\end{figure}

%Note that we have defined a quantum wedge to be a particular distribution,

\begin{rmk}\label{R:wedge_RV}
	Observe that by the corresponding property of the Neumann GFF, if we restrict the index set of $h$ defined above to $\cD_0(S)$, it gives rise to a stochastic process having a version that almost surely defines a distribution in $S$, that is, an element of $\cD'_0(S)$. In fact, by \cref{R:NGFF_fixed}, it is almost surely an element of $H^{-1}_{\mathrm{loc}}(S)$.
	
	We can then define the $\alpha$-quantum wedge as a doubly marked random surface, by letting it be the equivalence class of $(S,h^{\mathrm{wedge}},+\infty,-\infty)$, in the sense of \cref{def:markedqs}.\footnote{Note that this is the same definition as in \cite{Sheffield, DuplantierMillerSheffield} for the thick quantum wedge as a doubly marked quantum surface, but in \cite{Sheffield, DuplantierMillerSheffield} it is represented by $(S,\tilde{h},-\infty,+\infty)$ instead, where $\tilde{h}(\cdot)=h(-\cdot)$.} Using the change of coordinate formula  we could thus also view it as being parametrised by the upper half plane, and we would obtain a distribution $\hat h$ defined on $\H$. %If we instead started with the distribution modulo constants $\bar h$, this would similarly incude a distribution modulo constants on $\H$ using the same change of coordinate formula.
	However the expression for $\hat h$ is not particularly nice, and makes the following proofs more difficult to follow, which is why we usually take the strip $S$ as our domain of reference.
\end{rmk}

Note also that, when we parametrise this quantum surface by the infinite strip $S=
\{x+iy: x\in \R, y\in (0,\pi)\}$ with the two marked points at $+\infty$ and $-\infty$, that is, restrict to equivalence class representatives of the form $(S,h^{\mathrm{wedge}},+\infty,-\infty)$ where $h^{\mathrm{wedge}}$ is a field on $S$, then there is still one degree of freedom in the choice of the field $h^{\mathrm{wedge}}$, given by translations. Namely, because the term $Q\log |(T^{-1})'|$ in the change of coordinates formula \eqref{eq:concoc} disappears when $T$ is a translation,
$$ (S,h^{\mathrm{wedge}},+\infty,-\infty) \text{ and } (S,h^{\mathrm{wedge}}(\cdot+a),+\infty,-\infty)$$ 
are equivalent as doubly marked quantum surfaces,  for any $a\in \R$. In other words, we should not distinguish between $h^{\mathrm{wedge}}$ and $h^{\mathrm{wedge}}(+\cdot a)$. In \cref{D:wedge} of the thick quantum wedge, we chose a specific representative field $h^{\mathrm{wedge}}$ by fixing the horizontal translation so that the radial part of the field hit $0$ for the first time at $0$. But, in light of the discussion above, we could alternatively think of the wedge as being the doubly marked quantum surface, which when parametrised by the strip $S$ with marked points at $+\infty,-\infty$, is represented by $h^{\mathrm{wedge}}_{\cc} + \tilde Y_{\Re (z)}$ \emph{modulo translation},
where $h^{\mathrm{wedge}}_{\cc}$ has the law of $h_{\cc}^{\GFF}$, and $\tilde Y $ is a  Brownian motion with speed two and drift $-b = - (Q-\alpha)$, from $+ \infty$ to $- \infty$, i.e., has the law \eqref{hradwedge}, as per \cref{D:wedge}.

In fact, later on it will be slightly more convenient to (equivalently) parametrise the wedge by the strip $S$ with marked points $-\infty,+\infty$ (that is, switched). In this case, the thick quantum wedge is represented by the field
\begin{equation} h_{\cc}+(B_{2\Re(\cdot)}+(Q-\alpha)\Re(\cdot)) \text{ modulo horizontal translation,} \label{E:wedgemodtranslation} \end{equation} 
where the radial part $(B_{2\Re(\cdot)}+(Q-\alpha)\Re(\cdot)) $ is a Brownian motion with speed two and drift  $b =  (Q-\alpha)$ from $- \infty$ to $+ \infty$, (i.e., with drift of the opposite sign), using the terminology below \eqref{hradwedge}.

We note finally that there is an embedding of a thick quantum wedge in the upper half plane for which the associated field has a nice description in $\D_+$:

\begin{rmk}\label{rmk:wedge_neumann} Note that when $s\ge 0$, $h^{\mathrm{wedge}}_{\rr}(s)$ is  simply a Brownian motion with a drift of coefficient $\alpha-Q \le  0$. %as $s\to \infty$.
This means that, embedding in the upper half plane using $z\mapsto - e^{-z}$ and taking into account the conformal change of variables formula, the obtained representative
$(\H,\hat h^{\mathrm{wedge}},0,\infty)$ of the quantum wedge has a logarithmic singularity of coefficient $\alpha$ near zero. %(that is, the mean in any normalisation is asymptotic to $\alpha \log (1/z)$ near zero). Note also that %if $h$ is a quantum wedge embedded in $\H$ as above, with the specific normalisation of \cref{D:wedge}, then
In fact,
$$\hat h^{\mathrm{wedge}}(z)\big|_{\D_+} \overset{(law)}{=} (h+\alpha \log1/|z|)\big|_{\D_+},$$ where $h$ has the law of a Neumann GFF in $\H$, normalised so that it has zero average on the semicircle of radius 1. %Moreover, the unit circle is special in this embedding for the following reason: it is the circle of greatest radius for which the circle average is zero (since $h_{\rr}$ is conditioned to be positive for $s<0$). As a result this embedding is commonly known as the \textbf{unit circle embedding} for the quantum wedge. We emphasise that this embedding is in no way the \emph{canonical description} of this wedge (see Example \ref{ex:cd}). However it is far more convenient to work with.
\end{rmk}

\begin{rmk}[Unit circle embedding]
	Suppose that $h$ is a distribution on $S$ of the form $h_{\cc}^{S}+h_{\rr}$ where $h_{\rr}$ is constant on each vertical segment $\{\Re(z)=s\}$, and these constant values define a continuous function $h_{\rr}(s)$ that is positive for all $s$ less than some $ s_0\in \R$ but touches zero at some $s_1 \ge s_0$. Consider the unique translation of the strip so that the image of $h_{\rr}$ under this translation hits $0$ for the first time at $s=0$, and let $\tilde{h}$ be the image of $h$ after applying this translation, mapping to $\H$ using the map $z\mapsto -\e^{-z}$ and applying the change of coordinates formula.
	
	If a quantum surface has a representative of the form $(S,h,\infty,-\infty)$ with $h$ as above, then we call $(\H,\tilde{h},0,\infty)$ the \textbf{unit circle embedding} of this quantum surface.\ind{Random surface! Unit circle embedding}
	
	The unit (semi)circle clearly plays a special role in this embedding since it is the image of the vertical segment with $\Re(z)=0$ on the strip. Note that if $\hat{h}^{\mathrm{wedge}}$ is defined as in \cref{rmk:wedge_neumann}, then $(\H,\hat{h}^{\mathrm{wedge}}, 0,\infty)$ is the unit circle embedding of the $\alpha$-quantum wedge.
	\label{R:unitcircleembedding}
\end{rmk}%Again note that such an embedding is not available for general deterministic quantum surfaces: it must be the case that circle averages are defined for the underlying field and there is a largest circle of average zero.
We can now state the result about the quantum wedge being the scaling limit of a Neumann GFF with a logarithmic singularity near the origin.

\begin{theorem}
	\label{T:wedge} Fix $0 \le \alpha < Q$. Then the following hold:
	
	\vspace{.1cm}
	\noindent \emph{(i)} Let $\tilde h$ be a Neumann GFF in $\H$ with  additive constant fixed so that its upper unit semicircle average is equal to $0$, %so that $(\tilde h, \rho_0)$ is equal to $0$ for some $\rho_0\in {\mathfrak M}_N^\H\setminus \tilde{\mathfrak M}_N^\H$ that is compactly supported away from the origin, 
	and set $h(z) = \tilde h(z) + \alpha \log (1/|z|)$. 
	Let $h^C$ be such that $(\H,h^C,0,\infty)$ is the unit circle embedding of $(\H,h+C,0,\infty)$, and let $(\H, \hat{h}^{\mathrm{wedge}},0,\infty)$ be the unit circle embedding of a quantum wedge. Then for any $R>0$, $h^C|_{R\D_+}$ converges in total variation distance to $\hat h^{\mathrm{wedge}}|_{R\D_+}$ as $C\to \infty$.
	
	\vspace{.1cm}
	\noindent \emph{(ii)} If $(\H, h^{\mathrm{wedge}}, 0 , \infty)$ is an $\alpha$-quantum wedge, then $(\H,  h^{\mathrm{wedge}}, 0,\infty)$ and $(\H,  h^{\mathrm{wedge}} +C,0, \infty)$ have the same law as quantum surfaces.	
	
%	\vspace{.1cm}
	%In fact, we will prove a stronger form of (i): namely that for any $R>0$, if $(\H,h^C,0,\infty)$ is the unit circle embedding of the surface $(\H,h+C/\gamma,0,\infty)$, and $(\H, h^{\mathrm{wedge}},0,\infty)$ is the unit circle embedding of an $\alpha$-quantum wedge, then the total variation distance between $h^C|_{B(0,R)}$ and $h^{\mathrm{wedge}}|_{B(0,R)}$ as elements of $\cD_0'(B(0,R))$ goes to $0$ as $C\to \infty$.
%	For a more concrete formulation of (ii): if $(\H, h^{\mathrm{wedge}}, 0,\infty)$ is the unit circle embedding of an $\alpha$-quantum wedge and $(\H,h^{\mathrm{wedge}}_C,0,\infty)$ is the unit circle embedding of $(\H, h^{\mathrm{wedge}}+C/\gamma, 0, \infty)$, then $h^{\mathrm{wedge}}_C = h$ in distribution.
\end{theorem}

To summarise in the language introduced at the start of this chapter: (ii) says that a quantum wedge is invariant under rescaling, while (i) says that a quantum wedge is the limit, zooming in near zero, of the surface described by $\tilde h(z) + \alpha \log 1/|z|$. The fact that the convergence holds in the strong sense of total variation is very useful (as we shall see in the next chapter). Note that for this theorem we have to restrict to the case $\alpha<Q$.

\begin{proof} Note that point (ii) of the theorem follows from point (i), since the limit in law of $(\H, h+ C, 0, \infty)$ and $(\H, h + C + C', 0 , \infty)$ as $C \to \infty$, must coincide (in other words, scaling limits are by definition invariant under scaling). It thus suffices to prove point (i). 

%We start with (i) in the case that $\rho_0$ is the uniform measure on the unit semicircle centred at the origin. That is, $h=\tilde{h}-\alpha \log|z|$ where $\tilde{h}$ is the Neumann GFF with additive constant fixed so that $(\tilde{h},\rho_0)=0$. 
For this, we embed the field $h$ into the strip $S$ using the conformal isomorphism $z\in S \mapsto \phi(z) = - e^{-z} \in \H$, and apply the change of coordinates formula \eqref{eq:concoc}. Then the law of the resulting field can be written as %as %the sum of a Gaussian random variable (according to the choice of  additive constant) plus 
%	a standard two-sided Brownian motion equal to $0$ at $0$, with drift $\alpha-Q$ ($+\alpha$ coming from the logarithmic singularity of $h$, and $-Q$ from the change of coordinates). %In general the Gaussian random variable and the Brownian motion are not independent. Nevertheless,
%the radial part will hit $-C/2$ at some finite time almost surely. We can therefore further apply a horizontal translation to the strip so that the new resulting field $h^{C/2}$ has radial part hitting $-C/2$ for the first time at time $0$.
	
%	By \cref{L:diffbcs_agree} (and scaling) it follows that $h^{C/2}+C/2$ restricted to $S^+=S\cap \{\Re(z)>0\}$ converges in total variation distance, as $C\to \infty$, to the field  
\begin{equation}
	\label{almost_wedge}
h = 	h^{S}_{\cc} + h_{\rr},
	\end{equation}
	 where $h^S_{\cc}$ has the law specified by the radial decomposition of the Neumann GFF (see after Lemma \ref{L:radialdecomposition}); $h_{\rr}$ is independent of $h^{S}_{\cc}$; and finally $h_{\rr}$ is constant equal to $B_{2s} + (\alpha - Q) s$ on each vertical segment $\{\Re(z)=s\}$, with $B$ a standard two-sided Brownian motion equal to $0$ at time $0$. Here the $+\alpha$ comes from the logarithmic singularity of $h$, and the $-Q$ from the change of coordinates formula. %If the field $\tilde{h}$ in the statement of the theorem is normalised so that $(\tilde{h},\rho)=0$ for some $\rho\in \cD_0(\H)$, this follows from the Markov property of the Neumann GFF, applied in a small semi-disc centered at the origin that does not intersect the support of $\rho$. For arbitrary $\rho\in \cM$ {*...*}

	Observe that, by definition, the unit circle embedded field of $h + C $  is the image (after mapping back to $\H$ using the change of coordinates formula with $\phi$) of	$h^{S}_{\cc}(\cdot + s_C) + h_{\rr}(\cdot + s_C)+C$, where  
$$
s_C = \inf \{ s \in \R: B_{2s} + (\alpha - Q) s = - C\}
$$
 is the first hitting time of $-C$ by the radial part of $h$, and is independent of $h_{\cc}^{S}$.
 %\footnote{We note for use in the general $\rho_0$ case, that for any $T$ (which we will want to take large) we could define $\hat{s}_C$ to be the first hitting time after time $T$ that $(B_{2s}+(\alpha-Q)s)$ hits $-C$, and would have $\mathrm{d}_{\tv}(s_C,\hat{s}_C)\to 0$ as $C\to \infty$.\label{scsc} }
 
Since $h^S_{\cc}$ is independent of $s_C$ and is translation invariant, and since $h^{\mathrm{wedge}}$ has an angular part also independent of its radial part and with the same law, it suffices to show that as $C \to \infty$, 
\begin{equation}\label{eq:TVconvergencewedge}
C + h_{\rr}(\cdot + s_C) \to h^{\mathrm{wedge}}_{\rr} ( \cdot),
\end{equation}
in total variation, uniformly over intervals of $\R$ of the form $[- \log R, \infty)$. However, from our definition of $h^{\mathrm{wedge}}_{\rr} ( \cdot)$ in \eqref{hradwedge}, the two objects can be exactly coupled on $[-\log R , \infty)$ provided that $s_C \ge \log R$. As this is an event of probability tending to 1 as $C \to \infty$ and $R$ is fixed, we deduce the claim. 
\end{proof}

In fact, part (i) of \cref{T:wedge} still holds if we take an arbitrary choice of additive constant for the Neumann GFF $\tilde{h}$. We state something slightly stronger in the next Lemma: not only does this still hold if we choose the additive constant of the field slightly differently, but also the the way the scaling needs to be done when we zoom in is asymptotically independent of the initial choice of this additive constant. 

\begin{lemma}\label{lem:Twedgenorm}
Let $\tilde h,h$ be as in Theorem \ref{T:wedge}. Also let $\tilde{k}$ be  a Neumann GFF in $\H$ with  additive constant fixed so that $(\tilde{k}, \rho_0)$ is equal to $0$ for some $\rho_0\in {\mathfrak M}_N^\H\setminus \tilde{\mathfrak M}_N^\H$ (i.e., with nonzero average)
that is compactly supported away from the origin, and let $\mathfrak{h}$ be a deterministic function, that is harmonic in a neighbourhood of the origin.
Set $k(z) = \tilde{k}(z) +  \mathfrak{h} (z) + \alpha \log (1/|z|)$. 
Let $\varphi_C$ (resp. $\psi_C$), be the scaling maps such that $h^C:=h\circ \varphi_C^{-1}+Q\log|(\varphi_C^{-1})'|+C$ (resp. $k^C:=k\circ \psi_C^{-1}+Q\log|(\psi_C^{-1})'|+C$) is the unit circle embedding of $h+C$ (resp. $k+C$). Then for any bounded set $K \subset \H$, 
\[
\mathrm{d}_{\tv}\left( (\varphi_C ,h^C|_K ),(\psi_{C},k^C|_K) \right)\to 0
\]
as $C\to \infty$. 
\end{lemma}

\begin{proof}
 	We use the same notations as in Theorem \ref{T:wedge}. We first check that for any fixed $a\in \R$,
	\begin{equation} \label{eq:scalingfixed}
\mathrm{d}_{\tv}\left( (\varphi_C ,h^C|_K ),(\varphi_{C+a},h^{C+a}|_K) \right)\to 0
\end{equation}
as $C \to \infty$. % (i.e., the law of $(\varphi_C ,h^C|_K )$ is Cauchy with respect to the total variation distance). 
This will in fact follow from the observation that 
\begin{equation}
\label{eq:sCa} 
\mathrm{d}_{\tv} (s_{C}, s_{C+ a} ) \to 0
\end{equation}
as $C \to \infty$, for any fixed $a \in \R$. Indeed, once \eqref{eq:sCa} is known, let us fix $\eps>0$, and suppose that $K\subset B(0,R)$ for some $R>0$. \eqref{eq:sCa} means that for large enough $C$ we can choose a coupling of $s_{C/2}$ and $s_{C/2+a}$ which succeeds with probability at least $1- \eps$ (and we may further assume that $s_{C/2}=s_{C/2+a}>0$). Then, conditionally on this event and using the Markov property of Brownian motion, we can further couple the two processes $(h_{\mathrm{rad}} (s_{C/2} + \cdot ))$ and $ (h_{\mathrm{rad}} (s_{C/2+ a}  + \cdot ))$ so that they differ exactly by the constant $a$ on $\{\Re(\cdot)\ge 0\}$. In particular, on this event of probability $>1-\eps$, the radial parts hit $-C$ and $-(C+a)$ respectively at the same time, say $s$, and unless $|s_{C/2} - s| < \log R$ (which has probability tending to zero as $C \to \infty$) the radial parts agree up to the additive constant $a$ on $[s - \log R, \infty) $. \eqref{eq:scalingfixed} then follows by definition of $h^C$ and $h^{C+a}$ (since the scaling maps $\ph_C$ and $\ph_{C+a}$ depend only on $s_C$ and $s_{C+a}$ respectively).

Let us therefore check \eqref{eq:sCa}. First, we may instead work with $\tau_C$ and $\tau_{C+a}$ which are the first hitting times of $-C$ and $-(C+a)$ \emph{after time 0}, for a Brownian motion equal to $0$ at $0$ and with drift $- b:=- (Q-\alpha)<0$ (so $b>0$). (This suffices since the probability that $s_C$  or $s_{C+a}$ is negative converges to $0$ as $C\to \infty$). Then we observe, by an elementary calculation using the reflection principle and Girsanov's transform, that the density $f_{\tau_C}$ of $\tau_C$ with respect to Lebesgue measure on $[0,\infty)$ is given by  
\[ f_{\tau_C}(t) = \frac{C}{\sqrt{4\pi t^3}}e^{-\frac{(C-b t)^2}{4t}}.\]
It is easy to check from this expression that 
\[
\int_{C/b+C^{2/3}}^\infty f_{\tau_C}(t) \dd t \le \int_{C/b+C^{2/3}}^\infty f_{\tau_{C+a}}(t) \dd t \to 0
\]
as $C\to \infty$. 
Moreover, 
\[\int_0^{C/b - C^{2/3}} f_{\tau_{C+a}}(t) \dd t \le \int_0^{C/b-C^{2/3}} f_{\tau_C}(t) \dd t\to 0
\]
as $C\to \infty$, since the right hand side can be bounded above by the probability that the maximum of Brownian motion at time less than a constant times $C$ is bigger than a constant times $C^{2/3}$. Thus, to prove \eqref{eq:sCa} it is enough to show that 
\[
\sup_{t\in [C/b-C^{2/3}, C/b+C^{2/3}]} \left|\frac{f_{\tau_{C+a}}(t)}{f_{\tau_{C}}(t)}-1\right|\to 0
\]
as $C\to \infty$. For this we use that 
\[
\left|\frac{f_{\tau_{C+a}}(t)}{f_{\tau_{C}}(t)}-1\right|=\left| 1- (1+\frac{a}{C})e^{ a(-\frac{a}{2t}- \frac{C}{ t} +b)  }
\right|
\]
and that $|(-\frac{a}{2t}- \frac{C}{ t} +b)|=O(C^{-1/3})$ when $t\in [C/b-C^{2/3}, C/b+C^{2/3}]$.

This concludes the proof of \eqref{eq:scalingfixed}. It is easy to see that this implies that if $X$ is an independent finite random variable, 
	\begin{equation} \label{eq:scalingfixed2}
\mathrm{d}_{\tv}\left( (\varphi_C ,h^C|_K ),(\varphi_{C+X},h^{C+X}|_K) \right)\to 0
\end{equation}
as $C \to \infty$. % (i.e., the law of $(\varphi_C ,h^C|_K )$ is Cauchy with respect to the total variation distance). 
To see this, suppose that $X_n, Y_n$ are two random variables on a Polish space such that $\mathrm{d}_{\tv} (X_n+a, Y_n) \to 0$ as $n\to \infty$. Let $\mu_n^a$ denote the law of $X_n+a$ and let $\nu_n$ denote the law of $Y_n$. Let $\lambda(\dd a)$ denote the law of $X$. Then 
\begin{align*}
\sup_{A\in \cF} |\P( X_n + X \in A) - \P(Y_n \in A) | & = \sup_{A\in \cF} \left| \int \mu_n^a ( A) \lambda(\dd a)   - \nu_n(A) \right| \\
& \le \sup_{A\in \cF} \left | \int [ \mu_n^a (A) - \nu_n(A) ] \lambda(\dd a) \right| \\  
& \le \int \sup_{A\in \cF} | \mu_n^a( A) - \nu_n(A) | \lambda(\dd a) \\
& = \int  \mathrm{d}_{\tv} (X_n+a, Y_n) \lambda(\dd a). 
\end{align*}  
The total variation distance inside the integral tends to zero pointwise by assumption, and is bounded by one. The right hand side thus tends to zero by dominated convergence, as desired. This proves \eqref{eq:scalingfixed2}.

Next, we let $k,h$ be as defined in the statement of the Lemma, and we assume without loss of generality that $\int \rho_0 = 1$. This means that $\tilde{k}$ is equal in law to $\tilde{h}-(\tilde{h},\rho_0)$, where $\rho_0$ is supported away from the origin. By Lemma \ref{lem:hDNcircgood} (in fact, Remark \ref{rmk:hDNcircgood}), we then have
$$
\mathrm{d}_{\tv}(k|_{\delta\D_+},(h+X)|_{\delta \D_+})\to 0
$$
as $\delta\to 0$, where $X$ is random but independent of $h$, with the law of $\mathfrak{h}(0)+(\tilde h, \rho_0)$. Now let $\eps>0$ be arbitrary. First, for $\delta>0$ small enough, we can construct a coupling of $k$ and $h+X$ so that they agree on $\delta\D_+$ with probability at least $1-\eps$. By taking $C$ large enough, we may also assume that $(\varphi_{C+X}, h^{C+X}|_K)$ and $(\psi_C, k^C|_K)$ depend only on $h|_{\delta\D_+}$ and $k|_{\delta \D_+}$ on this event.
Since this gives a coupling under which $(\psi_C,k^C|_K)$ and $(\varphi_{C+X},h^{C+X}|_K)$ agree with probability at least $1-\eps$, this implies that
$$
\mathrm{d}_{\tv}\left( (\psi_C ,k^C|_K ),(\varphi_{C+X},h^{C+X}|_K) \right)\to 0
$$
as $C\to \infty$.
We conclude by \eqref{eq:scalingfixed} and the triangle inequality, since $X$ is almost surely finite.
%$$
%\mathrm{d}_{\tv}\left( (\varphi_C ,h^C|_K ),(\varphi_{C+X},h^{C+X}|_K) \right)\to 0
%$$
%as $C\to \infty$. 
%The triangle inequality then completes the proof.
\end{proof}

As an example of application of this result, we mention that the quantum wedge field $\hat{h}^{\mathrm{wedge}}$ with parameter $\alpha<Q$ (in the unit circle embedding) has a well defined Liouville bulk measure $\mathcal{M}_{\hat{h}^{\mathrm{wedge}}}$ and boundary measure $\mathcal{V}_{\hat{h}^{\mathrm{wedge}}}$, since it can be coupled with arbitrarily high probability to a Neumann GFF (with a given logarithmic singularity) plus a constant. Note that these measures are locally finite and atomless almost surely, by the results of \cref{S:GMC} (with base measure $\sigma$ incorporating the log singularity).	\ind{Liouville measure!Boundary (quantum wedge)} 
{It is not, however, obvious that $\hat{h}^{\mathrm{wedge}}$ assigns finite mass to any neighbourhood of $0$, and thus has a canonical embedding as discussed in \cref{ex:cd}. It is not too hard to see (and this will be made explicit below) that this would follow from \cref{T:GMCneumann_finite} if $\hat{h}$ were a Neumann GFF in $\H$ with average $0$ on the upper unit semicircle. But we must also consider the log-singularity at the origin; the next lemma shows that this does not cause the mass to explode.}

\begin{lemma}\label{L:wedgefinitemass}
{ Let $\hat{h}^{\mathrm{wedge}}$ be the representative field of a quantum wedge in the unit circle embedding. Then the bulk measure $\mathcal{M}$ associated with $\hat{h}^{\mathrm{wedge}}$ has $\cM(\D\cap \H)<\infty$ almost surely. }
\end{lemma}

\begin{proof}
It is more convenient (and equivalent) to prove that the bulk measure $\cM$ associated to $h^{\mathrm{wedge}}$ defined in $S$ as in \cref{D:wedge} satisfies $\cM([0,\infty)\times (0,\pi))<\infty$ almost surely. For $n\ge 0$ we let $S_n$ be the strip $[n,n+1]\times (0,\pi)$
	and set $0\le q<\min(1,1/\gamma)\in [0,1)$. 
	
	We first claim that for each fixed $n\ge 0$ we have 
	\begin{equation}\label{finitemomentstrip}
	\E[ \cM(S_n)^q] < \infty.
	\end{equation}
	 Equivalently, let us show that $\E [ \hat \cM(\hat S_n)^q] < \infty$, where $\hat \cM$ is the measure associated to $\hat h^{\mathrm{wedge}}$, and $\hat S_n$ is the image of $S_n$ by the map $z\mapsto - e^{-z}$, i.e., the upper annulus centered at $0$ of radii $e^{-n}$ and $e^{ - n -1}$ respectively. Since $n\ge 0$ is fixed, the logarithmic singularity at the origin of $\hat h^{\mathrm{wedge}}$ is of no relevance, and it suffices to show that $\E[ \cM_h ( \hat S_n)^q ] < \infty$, where $\cM_h$ is the GMC associated to a Neumann GFF $h$ in $\H$, normalised to have zero unit upper-circle average. We will prove this by making use of  \cref{T:GMCneumann_finite}. To do this we will map the upper half plane to the unit disc using the map $z\mapsto \psi(z) = (z-i)/ (z+i)$. The image of $h$ under this map is a Neumann GFF in $\D$, say $\tilde h$, normalised to have a certain average (along a chord extending from $i$ to $-i$) equals zero. The image of $\hat S_n$ under this map is a certain subdomain, say $\tilde S_n$, of the unit disc $\D$, made of four arcs of circle (two arcs of the unit circle and two hyperbolic geodesics respectively). Since $|\psi'|$ is bounded away from 0 and $+\infty$ on $S_n$, it suffices to show that $\E [ \cM_{\tilde h} (\tilde S_n)^q ] < \infty$. Although $\tilde h$ is normalised differently from the  Neumann GFF on the unit disc in \cref{T:GMCneumann_finite}, but this is of no relevance: the two fields differ by a fixed Gaussian random variable, so that we can conclude by H\"older's inequality and \cref{T:GMCneumann_finite}. This proves \eqref{finitemomentstrip}.

Now, by definition, since we restrict to the strip with real positive part,  $h^{\mathrm{wedge}}=h^{\mathrm{wedge}}_{\rr}+h^{\mathrm{wedge}}_{\cc}$ where the two summands are independent, and $h^{\mathrm{wedge}}_{\rr}(z)=B_{2\Re(z)}+(\alpha-Q)\Re(z)$ for $B$ a standard Brownian motion.  Setting $m_n=\sup\{h_{\mathrm{rad}}^{\mathrm{wedge}}(z): z\in S_n\}$ we therefore have, by subadditivity (since $0\le \gamma \le 1$), 
\[
\mathbb{E}[\cM([0,\infty)\times (0,\pi))^q|h^{\mathrm{wedge}}_{\rr}]\le \sum_{n\ge 0} e^{qm_n}\mathbb{E}[\cM_{\cc}(S_n)^q],
\]
where $\cM_{\cc}$ is the chaos measure associated to $h^{\mathrm{wedge}}_{\cc}$, i.e. $\lim_{\eps\to 0}\eps^{\gamma^2/2}e^{\gamma (h^{\mathrm{wedge}}_{\cc})_\eps(z)} dz$. By \eqref{finitemomentstrip}, we have that $\mathbb{E}[\cM_{\cc}(S_n)^q]< \infty$ for all $n \ge 0$. On the other hand, this quantity cannot depend on $n\ge 0 $ by translation invariance of $h^{\mathrm{wedge}}_{\cc})$. Thus  $\mathbb{E}[\cM_{\cc}(S_n)^q] = c< \infty$ for all $n\ge 0$. 
Moreover, since $m_n$ is the maximum of a Brownian motion with drift on $[n,n+1]$ for each $n$, we have (since $B_t / t \to 0$ almost surely as $t\to \infty$), 
\begin{equation*}\label{e:BMsumfinite}
\sum_{n\ge 0} e^{qm_n}<\infty
\end{equation*}
almost surely. Putting these together, we obtain that $\mathbb{E}[\cM([0,\infty)\times (0,\pi))^q|h^{\mathrm{wedge}}_{\rr}]<\infty$ almost surely, and hence $\cM([0,\infty)\times (0,\pi))<\infty$ almost surely as well.  
\end{proof}

 Hence, we obtain the following strengthening of \cref{T:wedge}. We emphasise that we are making use of the strong convergence (i.e., in total variation distance) here, which allows us to couple things so that they are actually equal (when restricted to compacts) with high probability. We are also using that for a quantum surface parametrised by $\H$ with marked points at $0$ and $\infty$, as in Example \ref{ex:cd}, the scaling map that determines the canonical parametrisation only depends on the field in a  neighbourhood of the origin with unit LQG area.

\begin{cor}\label{cor:wedge_cd_conv}
	\noindent (i) Suppose that $\tilde{h},h$ are as in  \cref{T:wedge}. If $(\H, h^C,0,\infty)$ is the canonical description of $(\H, h+C, 0,\infty)$ and $(\H,\hat{h}^{\mathrm{wedge}},0,\infty)$ is the canonical description of an $\alpha$-quantum wedge, then for any $R>0$, $h^C|_{R\D_+}\to \hat{h}^{\mathrm{wedge}}|_{R\D_+}$ in total variation distance as $C\to \infty$.

	\noindent (ii) Let $\mathcal{M}_{h^C},\mathcal{M}_{\hat{h}^{\mathrm{wedge}}}$ be the respective Liouville measures of $h^C,\hat{h}^{\mathrm{wedge}}$ as in (i). Then for any $R>0$ $\mathcal{M}_{h^C}|_{R\D_+} \to \mathcal{M}_{\hat{h}^{\mathrm{wedge}}}|_{R\D_+}$ in total variation distance as $C\to \infty$.
\end{cor}

We remark that the convergence in point (ii) of the above Corollary (with weak convergence rather than total variation) was actually used in some of the earlier work of Sheffield, for example in \cite{zipper}, as the definition of convergence in law for quantum surfaces.

\subsection{Quantum cones}\label{SS:cones}
\ind{Quantum cones} 

The quantum wedges discussed above are sometimes referred to as \emph{infinite volume quantum surfaces with boundary}, because their associated GMC measures have infinite mass, and because they are parametrised by simply connected domains with boundary (for example, the upper half plane $\H$ or the strip $S$). In this section we will discuss surfaces known as \textbf{quantum cones}: these are still infinite volume surfaces but now \emph{without boundary}, and are sometimes referred to as having ``the topology'' of the sphere rather than the disc. There also important examples of quantum surfaces with \emph{finite volume} (with or without boundary); these are \textbf{quantum discs} and \textbf{quantum spheres} and will be discussed later. 

In fact, the theory of quantum cones is entirely parallel to that of quantum wedges, the only difference being that they are defined on the whole plane or the infinite cylinder rather than on a simply connected domain.  These quantum cones are obtained in essentially the same way as the quantum wedges, but starting from a whole plane GFF rather than a Neumann GFF. 

Let $\cC$ be the infinite cylinder $\cC:=\{ z = x + i y\in \C: y \in [0,2\pi i]\}/\sim$, where $\sim$ identifies points $x$ with $x+2\pi i$ for $x\in \R$.	Let $\bar{H}^1(\cC)$ be the Hilbert space completion, with respect to the Dirichlet inner product $(\cdot, \cdot)_\nabla$, of the set of smooth functions modulo constants on $\cC$ with finite $(\cdot,\cdot)_\nabla$ norm. We first need the analogue of the radial decomposition for $\bar{H}^1(S)$.

\begin{lemma}\label{L:H1Cdecomp}
 Let $\bar \cH_{\rr}(\cC)$ be the subspace of $\bar H^1(\cC)$ obtained as the closure of smooth functions which are constant on each vertical segment $\{x+iy; y\in [0,2\pi i]\}$, viewed modulo constants. Let $\cH_{\cc}(\cC)$ be the subspace obtained as the closure of smooth functions which have mean zero on all such vertical segments. Then
	$$\bar H^1(\cC) = \bar \cH_{\rr}(\cC) \oplus \cH_{\cc}(\cC).$$
\end{lemma}
\begin{proof}
	This is similar to the proof of the radial decomposition for $\bar{H}^1(S)$: we leave it to the reader as part of Exercise \ref{Ex:cone}.
\end{proof}

Recall that the \textbf{whole plane GFF} $\bar{\mathbf{h}}^\infty$ is the random distribution modulo constants on $\C$ (that is, a continuous linear functional on 
$\tilde\cD_0(\C)$, the set of $f\in C^\infty(\C) $ with compact support and $\int_\C f =0$) with covariance kernel $G^\infty(x,y)=-\tfrac{1}{2\pi}\log(|x-y|)$. We denote $$h^\infty:=\sqrt{2\pi} \mathbf{h}^\infty$$ as usual. The (whole plane) GFF on the cylinder $\cC$, $\bar{h}^\cC$, is then defined by 
$$\bar{h}^{\cC}:=\bar{h}^\infty\circ \psi^{-1}$$
where $\psi: \cC\to \C$ is the map $z\mapsto -\log(1/z)$, and the meaning of the above is that $(\bar{h}^{\cC},f)=(\bar{h}^\infty, |\psi'|^2 f\circ \psi)$ for every $f\in C^\infty(\cC)$ with compact support and $\int_{\cC} f =0$. Due to the covariance structure, similarly to in the Neumann GFF case, we can extend the definitions of $\bar{h}^{\cC},\bar{h}^\infty$ respectively to be stochastic processes indexed by a larger index sets; namely, the sets $\tilde{\mathfrak{M}}_\infty^{\cC},\tilde {\mathfrak{M}}_\infty^\C$  of signed Radon measures on $\cC,\C$ respectively, whose positive and negative parts $\rho^{\pm}$ have equal mass and satisfy $\int \log|x-y| \rho^{\pm}(\dd x)\rho^{\pm}(\dd y)<\infty$.

Just as in the case of the Neumann GFF, \cref{L:H1Cdecomp} means that we can decompose 

\begin{equation}\label{E:hCdecomp} \bar{h}^{\cC}=\bar{h}^{\cC}_{\rr}+h^{\cC}_{\cc}\end{equation}
where:
\begin{itemize}
	\item $\bar{h}^{\cC}_{\rr}$ and $h^{\cC}_{\cc}$ are independent;
	\item $\bar{h}_{\rr}^{\cC}=\bar{B}_s$ if $\Re(z)=s$, where $\bar{B}$ is a standard Brownian motion modulo constants; 
	\item ${h}_{\cc}^\cC$ has mean zero on each vertical segment.
\end{itemize}

Again we leave the details as an exercise for the reader. Notice that the Brownian motion is run at the standard speed in this decomposition (rather than speed two in the case of the Neumann GFF) since, after mapping to the whole plane, this corresponds circle averages around an interior rather (rather than a boundary point). 

This leads us to the definition of an $\alpha$-quantum cone, again for 
$$ 0 \le \alpha \le Q=\frac{2}{\gamma}+\frac{\gamma}{2}.$$

\begin{definition}
	\label{D:cone}
	Let
	\begin{equation}\label{E:cone_radial}
		h^{\mathrm{cone}}_{\rr}(z) =
		\begin{cases}
			B_{s} + (\alpha - Q) s & \text{ if } \Re(z) = s \text{ and } s \ge 0\\
			\widehat B_{- s} + (\alpha - Q) s & \text{ if } \Re(z)=s \text{ and } s < 0
		\end{cases}
	\end{equation}
	where $B=(B_t)_{t\ge 0}$ is a standard Brownian motion, and $\widehat B=(\widehat B_t)_{t\ge 0}$ is an independent Brownian motion conditioned so that $\widehat B_{t} + (Q - \alpha) t >0$ for all $t>0$ (see \eqref{hradwedge} for the meaning of this conditioning).
	
	Let $h^{\mathrm{cone}}_{\cc}$ be a stochastic process indexed by $\mathfrak{M}_\infty^{\cC}$, that is independent of $h^{\mathrm{cone}}_{\rr}$ and has the same law as $h^{\cC}_{\cc}$.
	Then we set $h^{\mathrm{cone}} = h^{\mathrm{cone}}_{\rr} + h^{\mathrm{cone}}_{\cc}$ (which, since $h^{\mathrm{cone}}_{\rr}$ is just a continuous function, can again be defined as a stochastic process indexed by $\mathfrak{M}_\infty^{\cC}$). We call $h^{\mathrm{cone}}$ an \textbf{$\alpha$-quantum cone} field in $\cC$. 
\end{definition}

\begin{rmk}
	Notice that the speed of the Brownian motion in \cref{D:cone} above is one, rather than two in \cref{D:wedge}. This is, roughly speaking, because the Neumann GFF has double the variance of the whole plane GFF near the real line.
\end{rmk}

We will again want to view the above definition as being a specific equivalence class representative of a quantum surface with two marked points (that we will also refer to as an $\alpha$-quantum cone with an abuse of notation). That is, if $h$ is as in \cref{D:cone}, we will associate with it the quantum surface with two marked points $(\cC,h,\infty,-\infty)$. Another quadruple $(D,h',a,b)$ represents the same quantum surface if there is a conformal isomorphism $T:\cC\to D$ with $T(\infty)=a$, $T(-\infty)=b$ and $h'=h\circ T^{-1}+Q\log|(T^{-1})'|$ as in \cref{def:markedqs}. Similarly to the quantum wedge, any such representative will have finite associated Gaussian multiplicative chaos mass in any neighbourhood of $a$, and infinite mass in any neighbourhood of $b$. 

One particularly nice equivalence class representative of the $\alpha$-quantum cone is obtained by conformally mapping to $\C$ using the map $z\mapsto -\e^{-z}$ which sends $\infty$ to $0$ and $-\infty$ to $\infty$. Under this mapping, the vertical segment $\{t+iy: y\in [0,2\pi]\}\subset \cC$ mapped to the circle of radius $\e^{-t}$ around $0$ in $\C$, and the shift from the conformal change of coordinates formula is given by $-Q\Re(z)$. As in the wedge case, the obtained representative $(\C,\hat{h}^{\mathrm{cone}},0,\infty)$ of the $\alpha$-quantum cone is said to be in the \emph{unit circle embedding} and the field restricted to the unit disc $\D$ has the same law as $h^\infty+\alpha \log({1}/{|z|})$ restricted to $\D$, where $h^\infty$ is a whole plane GFF with additive constant fixed so that its average on $\partial \D$ is equal to $0$.

Finally, we can state the analogue of \cref{T:wedge}, which identifies the quantum cone as a local limit of a whole plane GFF with an additional log singularity of strength $\alpha$ at the origin.

\begin{theorem}
	\label{T:cone} Fix $0 \le \alpha < Q$. Then the following hold:
	
	\vspace{.1cm}
	\noindent \emph{(i)} Let $\tilde h^\infty$ be a whole plane GFF (in $\C$) with  additive constant fixed so that $(\tilde h, \rho_0)$ is equal to $0$ for some signed Radon measure $\rho_0$ with compact support away from the origin in $\C$,  $\int_{\C\times \C} \log|x-y||\rho_0|(\dd x)|\rho_0|(\dd y) <\infty$ and $\rho_0(\C)\ne 0$. Set $h(z) =\tilde h^\infty(z) + \alpha \log 1/|z|$. Let $h^C$ be such that $(\C,h^C,0,\infty)$ is the unit circle embedding of $(\C,h+C,0,\infty)$, and let $(\H,h^{\mathrm{cone}},0,\infty)$ be the unit circle embedding of an $\alpha$-quantum cone. Then for any $R>0$, $h^C|_{R\D_+}$ converges in total variation distance to $h^{\mathrm{cone}}|_{R\D_+}$ as $C\to \infty$.
	
	\vspace{.1cm}
	\noindent \emph{(ii)} If $(\C, h^{\mathrm{cone}}, 0 , \infty)$ is an $\alpha$-quantum cone, then $(\C, h^{\mathrm{cone}}, 0,\infty)$ and $C\in \R$, then $(\C, h^{\mathrm{cone}} +C,0, \infty)$ have the same law as quantum surfaces.	
\end{theorem}

\begin{proof}
	Exercise \ref{Ex:cone}.
\end{proof}

\subsection{Thin quantum wedges}
\ind{Quantum wedges} 

Recall that for $0\le \alpha < Q$ we defined an (equivalence class representative) of the $\alpha$-quantum wedge to be the random distribution on the infinite strip $S$ whose circular (or angular) part $h^{\mathrm{wedge}}_{\cc}$ is equal in law to the circular part $h_{\cc}^{\GFF}$ of a Neumann GFF on $S$, and whose radial part $h^{\mathrm{wedge}}_{\rr}$ (which is constant on vertical line segments) is independent of $h^{\mathrm{wedge}}_{\cc}$ and evolves as a speed two infinite Brownian motion $B_{2s}$, plus a negative drift of $(\alpha-Q)s$ from $+ \infty$ to $-\infty$, translated to hit $0$ for the first time at time $0$ (recall \ref{hradwedge}).
\textbf{Thin quantum wedges} are the surfaces obtained when the parameter $\alpha$ is instead taken in the range $(Q,Q+\tfrac{\gamma}{2}$). We will see that in this case, it is not possible to represent the surface by a single random field defined on $\H$ or $S$, but the correct definition is rather as a Poisson point process of quantum surfaces, or \textbf{beads} of the quantum wedge in the terminology of \cite{DuplantierMillerSheffield}.

Let us now make a useful connection between thick quantum wedges and Bessel processes, in order to motivate the definition of thin quantum wedges.

\begin{definition}[Bessel process] Let $\delta>0$. We define the Bessel process of dimension $\delta$ started from $x\ge 0$ to be $Z_t = Y_t^{1/2}$, where $Y$ solves the square Bessel stochastic differential equation (SDE), namely
\begin{equation}\label{eq:BESQ}
\dd Y_t = 2\sqrt{Y_t} \dd B_t + \delta \dd t; \quad Y_0 = x^2.
\end{equation}
(See \cite[Section 3, Chapter IX]{RevuzYor} for the existence and uniqueness of solutions to this SDE).\index{SDE} 
\end{definition}
\ind{Bessel!Process}

Applying It\^o's formula, we can see that on intervals of time in which $Z$ is not equal to $0$, $Z$ satisfies its own SDE:
\begin{equation} \label{eq:SDEBES} 
\dd Z_t=\dd B_t + \frac{\delta-1}{2Z_t} \, \dd t \; ; \; Z_0 = x.
\end{equation}
However, defining $Z$ directly from \eqref{eq:SDEBES} is far from straightforward because of the singularity of the drift term when $Z$ gets close to zero. When $\delta\ge 2$, it is easy to check $Z_t>0$ for all $t>0$ and if $\delta>2$, $Z_t\to \infty$ as $t\to \infty$ with probability one (i.e., $Z$ is transient), and thus \eqref{eq:SDEBES} can be used as the definition of a Bessel process of dimension $\delta$. When $\delta<2$, $Z$ returns to $0$ infinitely often with probability one (see \cite[Chapter 11]{RevuzYor}), but 0 is instantaneously reflecting if $\delta >0$: that is, the Lebesgue measure of the set of times where $Z_t = 0$ is a.s. zero (\cite[Proposition (1.5), Chapter XI]{RevuzYor}). When $1< \delta<2$, it is still possible to think of the Bessel process of dimension $\delta$ as solution of \eqref{eq:SDEBES}, because it can be checked that the integral 
$$
\frac{\delta - 1}{2}\int_0^t\frac{\dd u}{Z_u}
$$
converges a.s. for all $t\ge 0$, and is equal to $Z_t - B_t$, where $B$ is the Brownian motion from \eqref{eq:SDEBES}. When $\delta \le 1$, the integral no longer converges and the SDE \eqref{eq:SDEBES} does not make sense on intervals of time during which $Z$ hits zero; in fact, it can be checked that $Z$ is then \textbf{not} even a semimartingale when $\delta<1$.  Nevertheless, the law of the Bessel process $Z$ of dimension $\delta>0$ is uniquely specified by the fact that that it is a Markov process on $[0, \infty)$ whose infinitesimal generator coincides on $C^2 ((0, \infty))$ with
$$
\frac12 \frac{\dd^{\hspace{.05cm}2}}{\dd x^2} + \frac{\delta -1}{2 x} \frac{\dd}{\dd x} 
$$
(i.e., it satisfies \eqref{eq:SDEBES} away from $x = 0$) and with instantaneous reflection at $x = 0$; see \cite[(1.a)]{PitmanYorDecompositionBridges}.

\medskip In what follows, we define a function on $(-\infty,\infty)$ modulo translation to be an equivalence class of functions $(-\infty,\infty)\to \R$, where two functions $x(t),x'(t)$ are equivalent if $x(t)=x'(t+a)$ for all $t$ and some $a\in \R$. The next lemma shows that the logarithm of a Bessel process with dimension $\delta\ge 2$,  started from $0$, reparametrised by its quadratic variation and viewed modulo translation, has the same law as a Brownian motion with positive drift from $-\infty$ to $+\infty$, i.e., precisely what arises in the definition of a thick quantum wedge. This will give us an alternative point of view on thick quantum wedges (first pointed out by Duplantier, Miller and Sheffield in \cite{DuplantierMillerSheffield}) which lends itself well to the generalisation required to treat the thin case of quantum wedges. 

%is simply a two-sided Brownian motion with positive drift depending on the dimension of the Bessel process.%\footnote{Recall that our definition of the thick quantum wedge field is a representation of the quantum surface when it is embedded in the strip with marked points at $+\infty,-\infty$ rather than $-\infty,+\infty$ as in \cite{DuplantierMillerSheffield}. This means that our Brownian motion actually has \emph{negative} drift, and we accordingly consider minus the logarithm of the Bessel process in the lemma below.}

\begin{lemma}\label{L:BessBM}
	Let $0\le \alpha < Q$ and let $(Z_t)_{t\ge 0}$ be a Bessel process of dimension 
	\begin{equation} \delta = \delta_{\mathrm{wedge}}(\alpha):= 2+ \frac{2(Q-\alpha)}{\gamma} \label{E:deltaalpha} \end{equation} 
	with $Z_0=0$. Consider the process 
	\begin{equation} \label{E:X} X_t:=\tfrac{2}{\gamma}\log(Z_{q(t)})\end{equation} where $q:(-\infty,\infty)\to (0,\infty)$ is defined by the requirement that $q(0)=q_0$ for some (arbitrary) $q_0>0$ and that the quadratic variation of $X$ satisfies $\dd \,[ X ]_t = 2 \dd t$ on $(-\infty,\infty)$.  Then \emph{as functions on $(-\infty,\infty)$ modulo translation}, 
	$$ \left(X_t\right)_{t\in \R}\overset{(\mathrm{law})}{=} \left(B_{2t}+(Q-\alpha)t\right)_{t\in \R},$$ 
	where the right hand side denotes a Brownnian motion with speed two and drift $b = Q- \alpha >0$, from $-\infty$ to $+\infty$, similar to \eqref{hradwedge}.
\end{lemma}

\begin{proof}
Fix $t_0 \in \R$ and let us condition on $(X_s)_{s\le t_0}$ (set $x = X_{t_0}$). Due to the strong Markov property of Brownian motion, it suffices to prove that $$\dd X_t = \dd B_{2t} + (Q-\alpha) \dd t \quad ; \quad t\ge t_0.$$ %The law of $(B_{2t}+(\alpha-Q)t)_{t\in \R}$ as a function modulo translation is determined by the property that for any $C\in \R$, if $(Y_t)_{t\in \R}$ is the unique equivalence class representative that hits $C$ for the first time at time $0$, then conditionally on $(Y_t)_{t\le 0}$, $(Y_t)_{t\ge 0}$ has the law of $(C+\hat{B}_{2t}+(\alpha-Q)t)_{t\ge 0}$ where $(\hat{B}_s)_{s\ge 0}$ is a standard linear Brownian motion started from $0$ at time $0$. 
%We will show that $X$, as defined in \eqref{E:X}, satsifies this same property. Let $C$ be fixed and let $\tau_C$ be the first time that $Z$ hits $\exp({-\tfrac{\gamma}{2}C})$. Then by the Markov property of $Z$ 
Since $\delta\ge 2$, $Z$ is a solution of the SDE \eqref{eq:SDEBES} for all time. Thus by It\^{o}'s formula, if $Y_t:=\tfrac{2}{\gamma}\log(Z_t)$, then
$$\dd  Y_t= \frac{2}{\gamma Z_t} \dd B_t + \frac{\delta-2}{\gamma Z_t^2} \dd t = \dd M_t+\tfrac{1}{2}(Q-\alpha) \dd \, [M]_t \quad \forall t,$$
where $$M_t:= \frac{2}{\gamma Z_t} B_t$$ is a continuous local martingale. By definition of $X$, we therefore have 
$$ \dd X_t = \dd \tilde{M}_t + \tfrac{1}{2}(Q-\alpha) \dd \, [\tilde{M}]_t \quad \forall t$$ where $\tilde{M}$ is a reparametrisation of $M$ such that $\dd\,  [ \tilde{M}]_t = 2\dd t$ for all $t$. By L\'{e}vy's characterisation of Brownian motion, it must be that $\dd \tilde{M}_t= \dd B_{2t}$. Substituting this into the expression for $\dd X_t$ concludes the proof.
%\begin{align*}
%	X_{\tau_C+t} & = C+\int_{q(\tau_C)}^{q(\tau_C+t)} \frac{-2}{\gamma Z_s} \dd B_s + \int_{q(\tau_C)}^{q(\tau_C+t)} \frac{1-\delta}{\gamma Z_s^2} \dd s - \int_{q(\tau_C)}^{q(\tau_C+t)} \frac{1}{\gamma Z_t^2} \dd s \\ 
%	& = C + M_{t} + \frac{2-\delta}{\gamma}\frac{\gamma^2}{2} t
%\end{align*}
%for $t\ge 0$
%where $$M_{t}:=\int_{q(\tau_C)}^{q(\tau_C+s)} \frac{-2}{\gamma Z_s} \dd B_s$$ is a continuous local martingale, independent of $(X_{s})_{s\le \tau_C}$, and with $M_0=0$, $d\langle M\rangle_t = 2 \dd t$ by definition. By L\'{e}vy's characterisation, $(M_t)_{t\ge 0}$ is therefore equal in law to $(\hat{B}_{2t})_{t\ge 0}$, for $\hat{B}$ a standard linear Brownian motion started from $0$ at time $0$, independent of $(X_{s})_{s\le \tau_C}$.
%Since $C$ was arbitrary, this completes the proof.
\end{proof}

\begin{rmk}\label{R:BessBM}
A similar argument applies when we do not assume $\delta \ge 2$ (i.e., if $\alpha > Q$). However in this case, to avoid conditioning on $q_0$ being less than the hitting time $\zeta$ of 0 by $Z$, which would affect the law of the process, we need to assume that $Z_0 = x>0$. %Suppose that $\alpha\in (Q,Q+\tfrac{\gamma}{2})$, and $Z$ is a Bessel process of dimension $\delta_{\mathrm{wedge}}(\alpha)\in (1,2)$ started from any $x>0$. Then exactly the same argument shows that 
The conclusion is that the process $X_t=\tfrac{2}{\gamma}\log(Z_t)$, reparametrised to have quadratic variation $[X]_t=2t$ for all time $t\ge 0$, is equal in law to $(B_{2t}+(Q-\alpha)t)_{t\ge 0}$ (with $B$ a standard Brownian motion started from $\tfrac{2}{\gamma}\log(x)$).
\end{rmk}

As a consequence of \cref{L:BessBM}, for $0\le \alpha < Q$ we can equivalently define the (thick) $\alpha$-quantum wedge to be the doubly marked quantum surface $(S,h,-\infty,+\infty)$, where
\begin{equation} h^{\mathrm{wedge}}=h^{\mathrm{wedge}}_{\cc}+X_{\Re(\cdot)} \text{ considered modulo horizontal translation,} \label{E:wedgemodtranslation2} \end{equation}
$X$ is as defined in \eqref{E:X}, and $h^{\mathrm{wedge}}_{\cc}$ is independent of $X$ with the law of $h_{\cc}^{S}$. The reason for rewriting the definition in this way is because, defining $X$ in terms of the $\delta_{\mathrm{wedge}}(\alpha)$ dimensional Bessel process $Z$, there will be a clear extension to the case $\alpha\in (Q,Q+{\gamma})$, corresponding to $\delta_{\mathrm{wedge}}(\alpha)\in (0,2)$. 

This extension will rely on the notion of excursion for the Bessel process, which we now introduce. As already mentioned, even when $\delta \in (0,2)$ the Bessel process of dimension $\delta$ is a strong Markov process for which $a=0$ is a recurrent point. It was already shown in the seminal work of It\=o \cite{Ito} how to attach to such a Markov process a collection of excursions which forms a Poisson point process. To state this properly requires a notion of \textbf{local time} for $Z$ at $a=0$ which, roughly speaking measures the amount of time spent by $Z$ near $a=0$. Traditionally (\cite{RevuzYor, kallenberg2}), local time is constructed for semimartingales and we have already mentioned that the semimartingale property for a Bessel process of dimension $\delta>0$ fails if $0<\delta<1$. Nevertheless, It\=o's theory does apply to the whole range of dimensions $\delta \in (0,2)$, and is based on a notion of local time which is called the Blumenthal--Getoor local time of $Z$ at $a$.  (An alternative would be to use the excursion theory for the squared Bessel process $Y_t = Z_t^2$, since that is both a recurrent process and a semimartingale for all $\delta \in (0,2)$, see \eqref{eq:BESQ}.) 

\ind{Local time}

The upshot is the following, which is both a definition and It\=o's result \cite{Ito} in this case:

\begin{definition}\label{D:bess}
Let $\delta\in (0,2)$, and $Z$ be a $\delta$ dimensional Bessel process. Then $Z$ has an associated It\={o} excursion measure $\nu_{\delta}^{\mathrm{BES}}$ \indN{{\bf Miscellaneous}! $\nu_\delta^{\mathrm{BES}}$; It\^o excursion measure for the $\delta$ dimensional Bessel process} on the space $\mathcal{E}$ of continuous paths from $0$ to $0$, equipped with the topology of uniform convergence, and a local time $l$ at $0$. It satisfies the classical It\={o} excursion decomposition with excursion measure $\nu_\delta^{\mathrm{BES}}$. That is, if $(e_i)_{i\ge 1}$ is any enumeration of the countable set of excursions that $Z$ makes from $0$, and for $i\ge 1$, $t_i$ is the common value of $l$ on the time interval associated to $e_i$, then 
$$
\sum_{i\ge 1} \delta_{(t_i,e_i)} 
$$ 
has the distribution of a Poisson point process on $\R_+\times \mathcal{E}$ with intensity measure $\dd u \otimes \nu_{\delta}^{\mathrm{BES}}$.
\end{definition}

\ind{Bessel!Excursion measure}

%Notice that with this viewpoint, the invariance of the quantum wedge under the addition of a constant  boils down to the Brownian scaling property of the Bessel process. 

%Note that because the $\delta(\alpha)> 2$ dimensional Bessel process is positive and transient, the radial part of the thick quantum wedge field tends to $-\infty$ at the boundary point $-\infty$, and tends to $+\infty$ at the boundary point $+\infty$ of $S$. 

\medskip 
Now suppose that $\alpha\in (Q,Q+\gamma)$, so that $\delta_{\mathrm{wedge}} (\alpha) \in (0,2)$. The above excursion decomposition means that we can extend the Definition \eqref{E:wedgemodtranslation2} of a quantum wedge to this range of $\alpha$. However, rather than a single surface we will actually get a certain \textbf{Poisson point process of quantum surfaces}. We are now ready to define the thin quantum wedge.

%We emphasise here that, although Bessel processes of dimension $\delta$ may be defined even when $\delta \le 1$ (and in fact even when $\delta \le 0$, at least until hitting zero for the first time), the classical excursion theory only gives the existence of $\nu_{\delta}^{\mathrm{BES}}$ when $\delta > 1$. This is because this theory applies to continuous Feller processes which are also \emph{semimartingales} (and satisfying some additional further assumptions), see \cite[Chapter 29]{kallenberg2}. However, as already mentioned, Bessel processes fail to have the semimartingale property for $\delta <1$. Nonetheless, we will later be able to \emph{extend} the definition of the excursion measure $\nu_{\delta}^{\mathrm{BES}}$ to all values of $\delta \in (-\infty,2)$, see \cref{D:bess_ext}. 

\begin{definition}\label{D:thinwedge}
	Let $\alpha\in (Q,Q+\gamma)$ and let $Z$ be a $\delta$ dimensional Bessel process with $\delta=\delta_{\mathrm{wedge}}(\alpha)=2+2(Q-\alpha)/\gamma$. Let $\sum_{i\ge 1} \delta_{(t_i,e_i)} $ be the Poisson point process of excursions of $Z$, as in \cref{D:bess}, and for each $i$ define $X^i$ to be the function on $(-\infty,\infty)$ modulo translation given by $(2/\gamma)\log e_i$ parametrised to have infinitesimal quadratic variation $2 \dd t$ (as in \cref{L:BessBM}). For each $i$, let $\cS_i$ be the doubly marked quantum surface $(S,h^i,+\infty,-\infty)$ where
	\begin{equation*} h^i=h^i_{\cc}+X^i_{\Re(\cdot)} \text{ considered modulo horizontal translation,}  \end{equation*}
and $\{h^{i}_{\cc} ; i\ge 1\}$ is a collection of independent copies of $h_{\cc}^{S}$, independent of $\{X^i; i\ge 1\}$. 

We define the (thin) $\alpha$-quantum wedge to be the Poisson point process
\[ \cW  = \sum_{i\ge 1} \delta_{(t_i, \cS_i)}.\]
\end{definition}

\begin{rmk}
In \cite{DuplantierMillerSheffield} the definition of thin quantum wedges is only given in the case $\alpha \in (Q, Q + \tfrac{\gamma}{2})$, corresponding to the case $\delta = \delta_{\mathrm{wedge}} (\alpha) \in (1,2)$. Indeed, for this range of $\alpha$ one gets an additional property for the law of the total mass (quantum area) of each surface in the Poisson point process above. This additional property is important in the mating of trees (see Proposition 4.4.4 in \cite{DuplantierMillerSheffield}, and see also the end of the proof of 
\cref{T:MOT} at the very end of \cref{S:MOT}); fortunately, in that case we will see that the relevant value of $\delta$ will be $\kp/4$ with $\kp>4$, so that indeed $\delta>1$.
\end{rmk}

To view a thin quantum wedge $\cW$ as a random variable in a nice (Polish) space, it is better to view each point in the Poisson point process (or ``bead'') $\cS_i$, as being embedded in the strip $S$ as above but with some fixed choice of translation for the field (for example, so that the radial part of the field has its maximum value at time $0$). In this case $\cW$ can be identified with a random variable in the space of (atomic) measures on $\R_+ \times C(\R, \R) \times H^{-1}_{\mathrm{loc}}$, where the last component describes the circular part of the field. \medskip

\begin{rmk}In the thick case, the potential relevance of the notion of quantum wedges is made clear by Theorem \ref{T:wedge}, which shows that the (thick) quantum wedge is the scaling limit obtained by zooming near zero in a Neumann GFF with appropriate logarithmic singularity at zero. There is however no immediate analogue of this result in the thin case, so that the reader may wonder at this point what is the relevance of this notion (or, put it another way, how to tell that the generalisation of thick quantum wedges using the Bessel process point of view given in Lemma \ref{L:BessBM} is the ``correct'' way to do this generalisation). Such a justification will be provided at least in \ref{S:MOT} where thin quantum wedges are an essential aspect of the mating of trees approach to LQG.
\end{rmk}

We have now described thick quantum wedges in terms of Bessel processes with dimensions $\delta>2$ and thin quantum wedges in terms of Bessel processes with dimensions $\delta\in (0,2)$. One nice consequence of this is a duality between thick and thin wedges corresponding to a duality between Bessel process of dimension $\delta$ and dimension $4-\delta$; see \cref{L:bessel_dual} below. 

Note that if $\alpha\in (Q-\gamma, Q)$ so that $\delta_{\mathrm{wedge}}(\alpha)\in (2,4)$, then 
$$4-\delta_{\mathrm{wedge}}(\alpha)=\delta_{\mathrm{wedge}}(2Q-\alpha)$$ In other words, the duality will be between $\alpha$-  and $\hat \alpha = (2Q-\alpha)$-quantum wedges. Note that  $\hat \alpha \in (Q, Q + \gamma)$ so that $\delta_{\mathrm{wedge}}( \hat \alpha) \in (0,2)$.
\ind{Duality! Bessel processes}
\ind{Duality! Quantum wedges}

Let us now describe this more precisely. 

\begin{lemma}\label{L:bessel_dual}
	For $\delta\in (0,2)$, decomposing a Bessel excursion according to its maximum value, we can write
	\begin{equation} \nu_{\delta}^{\mathrm{BES}} = c_\delta \int_{0}^\infty \nu_{\delta}^x  x^{\delta-3} \dd x \label{eq:bessel_dual}
	\end{equation}
	where:
\begin{itemize}
	\item $\dd x$ is Lebesgue measure on $\R_+$;
\item $c_\delta \in (0,\infty)$ depends only on $\delta$; and 
\item for each $x>0$, $\nu_\delta^x$ is a probability measure on excursions from $0$ to $0$ in $\R_+$ with maximum value $x$. A sample from $\nu_\delta^x$ corresponds to a Bessel process of dimension $4-\delta$ run until it first hits $x$, then concatenated with $x$ minus the time reversal of an independent copy of the same process.
\end{itemize}
\end{lemma}
\begin{proof}
See \cite[Theorem 1]{PitmanYor}. Note that the description of $\nu_\delta^x$ for given $x$ follows from \cref{L:BessBM} and \cref{R:BessBM}, plus the fact that conditioned on its maximum value, a Brownian motion with negative drift $-a$ can be written as a Brownian with drift $a$ until it hits this maximum value, and then concatenated with an independent Brownian motion with drift $-a$, conditioned to stay negative. See also \cref{L:BM_drift_rev} and \cite{Williams} for closely related statements.
\end{proof}

\begin{rmk}\label{R:thick_thin} As a consequence, we see that if $\alpha\in (Q,Q+\gamma)$ then, informally speaking, each of the quantum surfaces making up a (thin) $\alpha$-quantum wedge (when parametrised by the $S$) looks locally near $\pm \infty$ like a $(2Q-\alpha)$-quantum wedge in $S$ does near $-\infty$ (the marked point with neighbourhoods of finite quantum area). 
\end{rmk}

%\begin{definition}\label{D:bess_ext}
%We \emph{extend} the definition of  $\nu_{\delta}^{\mathrm{BES}}$ to all values of $\delta\in (-\infty, 2)$ by setting it equal to the right hand side of  \eqref{eq:bessel_dual} (noting that $\nu_\delta^x$ is well defined since $4-\delta >2$). 
%\end{definition}

\ind{Bessel!Excursion measure}

\subsection{Quantum discs }
\ind{Quantum discs}

 Having defined the thin quantum wedge for $\alpha\in (Q,Q+\gamma)$ as a Poisson point process of quantum surfaces, it is natural to ask about the ``law'' of each of these surfaces (although of course this actually corresponds to an infinite measure). This will lead us to the notion of quantum discs below.
  
 Recall that given an excursion $e_i$ of a $\delta$ dimensional Bessel process, we defined $X^i=:X^{e_i}$, a function on $(-\infty,\infty)$ modulo translation, to be given by $(2/\gamma)\log(e_i)$ parametrised to have infinitesimal quadratic variation $2 \dd t$. Since the excursion $e_i$ of the Bessel process starts at $0$, ends at $0$ and has finite maximum value, a natural way of fixing the horizontal translation of $X^{e_i}$ is to require that the maximum is reached at time $0$. Let us write $Y^{e_i}$ for this function on $(-\infty,\infty)$ (associated with the excursion $e_i$).
 
 \begin{definition}\label{D:qdisc}
 	Let $\alpha\in (Q,Q+\gamma)$, $\nu_{\delta}^{\mathrm{BES}}$ be as described in the previous subsection with $\delta=\delta_{\mathrm{wedge}}(\alpha)=2+ 2(Q-\alpha)/\gamma$, and $\P^{S}_{\cc}$ be the law of $h_{\cc}^{S}$ (obtained from a Neumann GFF in $S$ by subtracting its average value on each vertical line segment). We define the infinite $\alpha$-quantum disc measure $m^{\mathrm{disc}}_{\alpha}$ to be the measure on $H^{-1}_{\mathrm{loc}}(S)$ obtained by pushing forward the measure $\nu_{\delta}^{\mathrm{BES}}\otimes \P^{S}_{\cc}$ to $H^{-1}_{\mathrm{loc}}(S)$, via the map taking $(e,h_{\cc})$ to the field $h_{\cc}+Y^e_{\Re(\cdot)}$.
 \end{definition}

\begin{rmk}\label{R:DMSextendeddisc}
	In \cite{DuplantierMillerSheffield}, the definition of $\nu_{\delta(\alpha)}^{\mathrm{BES}}$ is extended to a larger range of $\alpha$, and thus so is the $\alpha$-quantum disc measure. However, we will stick to the case where the measure can be defined using classical  It\=o excursion theory; see \cref{D:bess}.
\end{rmk}

%We have defined $m_\alpha^{\textrm{disc}}$ to be a measure on doubly marked quantum surfaces, but we could also define it as a measure on $H^{-1}(S)/\sim_T$ where $\sim$ is the equivalence relation identifying two generalised functions $h,h'$ in $H^{-1}(S)$ if $h(\cdot)=h'(\cdot+a)$ for some $a\in \R$.
 
 \begin{rmk}\label{rmk:thinwedgepppdiscs} With this definition, if $\alpha \in (Q, Q + \gamma)$, the $\alpha$-quantum disc measure corresponds to the measure ``describing'' the individual quantum surfaces appearing in a (thin) $\alpha$-quantum wedge. Indeed, recalling Definitions \ref{D:bess} and \ref{D:thinwedge}, we see that an equivalent definition of the $\alpha$-quantum wedge is as a Poisson point process
 	\[ \cW  = \sum_{i\ge 1} \delta_{(t_i, \cS_i)}\]
 	with intensity $du \otimes \hat{m}_\alpha^{\mathrm{disc}}$, where $\hat{m}_\alpha^{\mathrm{disc}}$ is the pushforward of $m_\alpha^{\mathrm{disc}}$ by the map taking $h\in H_{\mathrm{loc}}^{-1}(S)$ to the doubly marked quantum surface $(S,h,-\infty,\infty)$.
 \end{rmk}

By \cref{R:thick_thin}, near each of the marked points, a sample from (some suitably conditioned version of) the $\alpha$-quantum disc measure  looks locally like a (thick) $(2Q-\alpha)$-quantum wedge at its apex (that is, near the marked point which has neighbourhoods of finite quantum mass). Or in other words, it looks locally like a free boundary Gaussian free field plus a $(2Q-\alpha)$-$\log$ singularity. (Of course this statement is informal on many levels!)
 
 \medskip 
 Notice that, due to Brownian scaling, a Bessel excursion with maximum $x$ (that is, sampled from $\nu_{\delta}^x$ with the notation of \cref{L:bessel_dual}) is equal in law to $x$ times a Bessel process with maximum $1$ (that is, sampled from $\nu_\delta^1$) modulo time change. However, since under the map $e\mapsto X^e$ we reparametrise time anyway (so that the infinitesimal quadratic variation is exactly $2 \dd t$) we see that the law of $Y^e$ when $e$ is sampled from $\nu_\delta^x$ is equal to the law of $((2/\gamma)\log x+Y^e)$ when $e$ is sampled from $\nu_\delta^1$. Hence from the decomposition in \cref{L:bessel_dual}, it follows that for any non-negative measurable function $F$ on $H^{-1}_{\mathrm{loc}}(S)$:

\begin{equation}\label{E:disint_m} m_\alpha^{\mathrm{disc}}(F) = c_\delta \int_{0}^\infty \mathbb{P}^{S}_{\cc} \otimes \nu_\delta^1(F(h_{\cc}+\tfrac{2}{\gamma}\log(x) + Y^e_{\Re(\cdot)})) \, x^{\delta-3} \dd x, \end{equation}
remembering that $\delta=\delta_{\mathrm{wedge}}(\alpha)=2+(2/\gamma)(Q-\alpha)$. %Note that although $h+X^e_{\Re(\cdot)}+\tfrac{2}{\gamma}$ is only defined modulo horizontal translation of $S$, $(S,h+X^e_{\Re(\cdot)}+\tfrac{2}{\gamma}\log(x),+\infty,-\infty)$ still uniquely defines a quantum surface.

From the description of $\nu_{\delta}^1$ in \cref{L:bessel_dual} we see that if $e$ is sampled from $\nu_\delta^1$ and $Y^e$ is as described above, then $Y^e_0=0$, and $(Y^e_t)_{t\ge 0}, (Y^e_{-t})_{t\le 0}$ are independent, each having the law of $(2/\gamma)\log Z$ reparametrised to have quadratic variation $2t$ at time $t$, where $Z$ is a Bessel process of dimension $4-\delta_{\mathrm{wedge}}(\alpha)>2$. By \cref{L:BessBM} we can rephrase this as follows. 

\begin{lemma} \label{L:Ydescription}Let $\alpha\in (Q,Q+\gamma)$ and 
	define $(Y_t)_{t\in \R}$ by setting $Y_0:=0$ and 
	\begin{itemize}
		\item $Y_t=B_{2t}+(Q-\alpha)t$ for $t>0$
		\item $Y_{t}=\hat{B}_{-2t}+(Q-\alpha)(-t)$ for $t<0$
	\end{itemize} 
where $B,\hat{B}$ are independent standard linear Brownian motions defined for $t\in [0,\infty)$, started from $0$ and conditioned that $B_{2t}+(Q-\alpha)t$ (resp. $\hat{B}_{2t}+(Q-\alpha)t$) is negative for all $t>0$. Then if $e$ is sampled from $\nu_\delta^1$
$$(Y_t)_{t\in \R}\overset{(d)}{=} (Y^e_t)_{t\in \R}.$$
\end{lemma}
 
 As promised, let us now justify that quantum discs really are \emph{finite volume} quantum surfaces (with boundary). 
 %Note that given a quantum surface $\cS$ parameterised by a simply connected domain, one can define its  total quantum area $\cM(\cS)$ by taking any representative $(D,h,a,b)$ of $\cS$ and defining $\cM(\cS)=\cM_h(D)$. By the definition of quantum surfaces (that is, of the equivalence relation) this will be independent of the choice of representative $(D,h,a,b)$. Similarly one can define its total quantum boundary length $\nu(\cS)$ by considering a representative in a domain with nice boundary (say $(S,h,a,b)$) and defining it to be the total mass of the GMC associated to $h$ on $\partial S$. 
 
 \begin{lemma}\label{L:finitedisc}
For $\alpha\in (Q,Q+\gamma)$, $\cM_{h}^\gamma(S)<\infty$ and $\mathcal{V}_h^\gamma(\partial S)<\infty$ for $m_\alpha^{\mathrm{disc}}$-almost every $h$.
 \end{lemma}

\begin{proof}
	We will verify the statement about $\cM$; leaving the boundary case as Exercise \ref{Ex:finitedisc}. By \eqref{E:disint_m}, it suffices to show that if $h_{\cc}$ is sampled from $\mathbb{P}_{\cc}^{S}$ (that is, has the law of $h_{\cc}^{S})$ and $e$ is sampled independently from $\nu_\delta^1$, then $$\cM_{h_{\cc}+Y_{\Re(\cdot)}^e}(S)<\infty$$ almost surely. 

{This follows from the same proof as \cref{L:wedgefinitemass}, using that under $\nu_{\delta}^1$, $Y^e$ is a two-sided Brownian motion with negative drift and $Y_0^e=0$: see \cref{L:Ydescription}.}
 \end{proof}

\paragraph{Conditioned quantum discs.}

Recall \eqref{E:disint_m}, which provides the decomposition 
\begin{equation*}
	m_\alpha^{\mathrm{disc}}(F) = c_\delta \int_{0}^\infty \mathbb{P}^{S}_{\cc} \otimes \nu_\delta^1(F(h_{\cc}+\tfrac{2}{\gamma}\log(x) + Y^e_{\Re(\cdot)})) \, x^{\delta-3} \dd x
\end{equation*} 
(for $F$ non-negative and measurable on $H^{-1}_{\mathrm{loc}}(S)$), of the $\alpha$-quantum disc measure on $H^{-1}_{\mathrm{loc}}$. Now we know that $m_\alpha^{\mathrm{disc}}$ is supported on quantum surfaces with finite quantum mass and boundary length, we can use the above decomposition to describe the pushforward of $m_\alpha^{\mathrm{disc}}$ via the map $h\mapsto \cM^\gamma_h(S)$ or $h\mapsto \mathcal{V}_h^\gamma(\partial S)$ very precisely. Indeed, we have (for example, working with the boundary length $\mathcal{V}$):

\[ m_\alpha^{\mathrm{disc}}(\mathcal{V}_h^\gamma(\partial S)\in A) = c_\delta \int_0^\infty \mathbb{P}_{\cc}^{S}\otimes \nu_\delta^1 (\mathcal{V}^\gamma_{h_{\cc}+Y^e_{\Re(\cdot)}+(2/\gamma)\log(x)}(\partial S)\in A) x^{\tfrac{2(Q-\alpha)}{\gamma}-1} \dd x.\]
%Here, if $h$ is a field defined modulo horizontal translation of $S$, then $\nu_h(\partial S)$ still makes sense because any way of fixing the horizontal translation will result in a field that gives the same GMC mass to $\partial S$. 
Notice that $$\mathcal{V}^\gamma_{h_{\cc}+Y^e_{\Re(\cdot)}+(2/\gamma)\log(x)}(\partial S)=x\, \mathcal{V}^\gamma_{h_{\cc}+Y^e_{\Re(\cdot)}}(\partial S)$$  for each $x$, so if we make the change of variables $u=x \, \mathcal{V}^\gamma_{h_{\cc}+Y^e_{\Re(\cdot)}}(\partial S)$ the above becomes
\begin{align*} m_\alpha^{\mathrm{disc}}(\mathcal{V}_h^\gamma(\partial S)\in A) & = c_\delta  \int_0^\infty \mathbb{P}_{\cc}^{S}\otimes \nu_\delta^1 \left(\indic{u\in A} u^{\tfrac{2(Q-\alpha)}{\gamma}-1} (\mathcal{V}^\gamma_{h_{\cc}+Y^e_{\Re(\cdot)}}(\partial S))^{-\tfrac{2(Q-\alpha)}{\gamma}} \right)\dd u \\ 
& = c_\delta \,  \mathbb{P}_{\cc}^{S}\otimes \nu_\delta^1 \left({(\mathcal{V}_{h_{\cc}+Y^e_{\Re(\cdot)}}(\partial S))}^{-\tfrac{2(Q-\alpha)}{\gamma}} \right) \int_{u\in A} u^{\tfrac{2(Q-\alpha)}{\gamma}-1} \dd u.\end{align*}
This yields the following conclusion.
 \begin{lemma}\label{L:volumepushforwarddisc} Let $\alpha\in (Q,Q+\gamma)$. The pushforward of $m_\alpha^{\mathrm{disc}}$ under the map $h\mapsto \mathcal{V}^\gamma_h(\partial S)$ is a constant multiple of the measure $u^{-1+2\gamma^{-1}(Q-\alpha)} \, \dd u$ on $[0,\infty)$. Similarly, the pushforward of $m_\alpha^{\mathrm{disc}}$ under the map $h\mapsto \cM^\gamma_h(\cS)$ is a constant multiple of the measure $u^{-1+\gamma^{-1}(Q-\alpha)} \, \dd u$ on $[0,\infty)$. \end{lemma}

To conclude the section, we are going to decompose the $\alpha$-quantum disc measure according to quantum boundary length and quantum area. It turns out to have a remarkable property: conditioned on the quantum boundary length or quantum area, if we subtract the correct constant from the field so that the area or boundary length becomes one, the law of the resulting field \emph{does not} depend on the mass or area that we conditioned on.

In the case of boundary length, we have the following:  
\begin{prop}  For $\alpha\in (Q,Q+\gamma)$
\begin{equation}
	\label{E:disc_disintegrated}
	m_\alpha^{\mathrm{disc}}=c_\delta \int_0^\infty \mathbb{P}_{\alpha}^{\mathrm{disc},u} u^{\tfrac{2(Q-\alpha)}{\gamma}-1} \dd u
\end{equation} 
where $\mathbb{P}_\alpha^{\mathrm{disc},u}$ is a \emph{probability measure} on $H^{-1}_{\mathrm{loc}}(S)$, such that for $\mathbb{P}_\alpha^{\mathrm{disc},u}$ every $h$, the boundary Gaussian multiplicative chaos measure $\mathcal{V}_h^\gamma(\partial S)$ is well defined and satisfies
$$\mathcal{V}_h^\gamma(\partial S)=u.$$ Moreover, if we write $h_\alpha^{\mathrm{disc},1}$ for the field on $S$  with the law of 
$$ h_{\cc}+ Y^e_{\Re(\cdot)}-\tfrac{2}{\gamma}\log \mathcal{V}^\gamma_{h_{\cc}+Y^e_{\Re(\cdot)}}(\partial S)  \text{ weighted by } \big(\mathcal{V}^\gamma_{h_{\cc}+Y^e_{\Re(\cdot)}}(\partial S)\big)^{-\tfrac{2(Q-\alpha)}{\gamma}},$$
 where $h_{\cc}\overset{(\mathrm{law})}{=} h_{\cc}^{S}$, and $e$ is an independent Bessel excursion of dimension $\delta_{\mathrm{wedge}}(\alpha)$ conditioned on taking maximum value $1$ (that is, with law $\nu_\delta^1$, so that $Y^e$ is as described in \cref{L:Ydescription}), then  $$h_{\alpha}^{\mathrm{disc},u}:=h_{\alpha}^{\mathrm{disc},1}+\tfrac{2}{\gamma}\log(u)$$ is a sample from $\mathbb{P}_\alpha^{\mathrm{disc},u}$.
\label{P:discdisBL}\end{prop}

 \begin{definition}\label{D:unitBLqd} We call the quantum surface $(S,h_{\alpha}^{\mathrm{disc},1},-\infty,+\infty)$ (when $h_{\alpha}^{\mathrm{disc},1}$ has law $\mathbb{P}^{\mathrm{disc},1}_\alpha$) a unit boundary length $\alpha$-quantum disc, and $(S,h_{\alpha}^{\mathrm{disc},u},-\infty,+\infty)$ (when $h_{\alpha}^{\mathrm{disc},u}$ has law $\mathbb{P}^{\mathrm{disc},u}_\alpha$) an $\alpha$-quantum disc with boundary length $u$. 
 \end{definition}

 With an abuse of terminology, we will also sometimes refer to just the field $h_{\alpha}^{\mathrm{disc},1}$ or $h_{\alpha}^{\mathrm{disc},u}$ as a unit boundary length $\alpha$-quantum disc, or an $\alpha$-quantum disc with boundary length $u$.
Let us emphasise that an $\alpha$-quantum disc with boundary length $u$ has the same law as a unit boundary length $\alpha$-quantum disc with a constant $\tfrac{2}{\gamma}\log(u)$ added to the field.

\begin{proof}[Proof of \cref{P:discdisBL}]
Let $F$ be a non-negative measurable function on $H^{-1}_{\mathrm{loc}}(S)$. To prove the proposition, it suffices to show that 
$$
m_\alpha^{\mathrm{disc}}(F(h-\tfrac{2}{\gamma}\log\mathcal{V}^\gamma_h(\partial S))\indic{\mathcal{V}_h^\gamma(\partial S)\in A}) =  m_{\alpha}^{\mathrm{disc}}(\mathcal{V}_h^\gamma(\partial S)\in A) \mathbb{P}_\alpha^{\mathrm{disc},1}(F)
$$
where $\mathbb{P}_{\alpha}^{\mathrm{disc},1}$ is the law of $h_\alpha^{\mathrm{disc},1}$, as described in the statement of the proposition.	%of a doubly marked quantum surface, and for a general doubly marked quantum surface $\cS$ represented by $(S,h,+\infty,-\infty)$ with $\nu_h(\partial S)<\infty$, write $\cS^*$ for the quantum surface represented by $(D,h-\tfrac{2}{\gamma}\log(\nu_h(\partial S)),+\infty,-\infty)$. 
	Applying the change of variables $u=x \, \mathcal{V}^\gamma_{h_{\cc}+Y^e_{\Re(\cdot)}}(\partial S)$, we have (very similarly to before):
	\begin{align*} & m_\alpha^{\mathrm{disc}}(F(h-\tfrac{2}{\gamma}\log\mathcal{V}^\gamma_h(\partial S))\indic{\mathcal{V}_h^\gamma(\partial S)\in A}) \\ & =c_\delta \,  \mathbb{P}_{\cc}^{\GFF}\!\otimes\! \nu_\delta^1 \left( \frac{F\big(h_{\cc}+Y_{\Re(\cdot)}^e - \tfrac{2}{\gamma}\log\mathcal{V}^\gamma_{h_{\cc}+Y_{\Re(\cdot)}^e }(\partial S) \big)}{(\mathcal{V}^\gamma_{h_{\cc}+Y^e_{\Re(\cdot)}}(\partial S))^{{2\gamma^{-1}(Q-\alpha)}}} \right) \int_A u^{\tfrac{2(Q-\alpha)}{\gamma}-1} \dd u \\
		& =  \frac{\mathbb{P}_{\cc}^{S}\!\otimes\! \nu_\delta^1 \left( {F\big(h_{\cc}\!+\!Y_{\Re(\cdot)}^e \!-\! \tfrac{2}{\gamma}\log\mathcal{V}^\gamma_{h_{\cc}+Y_{\Re(\cdot)}^e }(\partial S)\big)}(\mathcal{V}^\gamma_{h_{\cc}+Y^e_{\Re(\cdot)}}(\partial S))^{-\tfrac{2(Q-\alpha)}{\gamma}} \right)}{\mathbb{P}_{\cc}^{S}\!\otimes\! \nu_\delta^1 \left( (\mathcal{V}^\gamma_{h_{\cc}+Y^e_{\Re(\cdot)}}(\partial S))^{-\tfrac{2(Q-\alpha)}{\gamma}} \right)} \\
	&\quad \quad \times m_{\alpha}^{\mathrm{disc}}(\mathcal{V}_h^\gamma(\partial S)\in A) .\end{align*}
	The result then follows from the definition of $h_{\alpha}^{\mathrm{disc},1}$. 
\end{proof}

Similarly, we can make sense of a \emph{unit area quantum disc}.

\begin{prop}\label{prop:unitareaQD}
  For $\alpha\in (Q,Q+\gamma)$
	\begin{equation}
		m_\alpha^{\mathrm{disc}}=\hat{c}_\delta \int_0^\infty \hat{\mathbb{P}}_{\alpha}^{\mathrm{disc},a} a^{\tfrac{(Q-\alpha)}{\gamma}-1} \dd a
	\end{equation} 
	where $\hat{\mathbb{P}}_\alpha^{\mathrm{disc},a}$ is a \emph{probability measure} on $H^{-1}_{\mathrm{loc}}(S)$, such that for $\hat{\mathbb{P}}_\alpha^{\mathrm{disc},u}$ every $h$, the bulk Gaussian multiplicative chaos measure $\mathcal{M}_h^\gamma(S)$ is well defined and satisfies
	$$\mathcal{M}_h^\gamma(S)=a.$$ Moreover, if we write $\hat{h}_\alpha^{\mathrm{disc},1}$ for the field on $S$  with the law of 
	$$ h_{\cc}+ Y^e_{\Re(\cdot)}-\gamma^{-1}\mathcal{M}^\gamma_{h_{\cc}+Y^e_{\Re(\cdot)}}(S)  \text{ weighted by } \big(\mathcal{M}^\gamma_{h_{\cc}+Y^e_{\Re(\cdot)}}( S)\big)^{-\tfrac{(Q-\alpha)}{\gamma}},$$
	where $h_{\cc}\overset{(\mathrm{law})}{=} h_{\cc}^{S}$, and $e$ is an independent Bessel excursion of dimension $\delta(\alpha)$ conditioned on taking maximum value $1$ (that is, with law $\nu_\delta^1$, so that $Y^e$ is as described in \cref{L:Ydescription}), then  $$\hat{h}_{\alpha}^{\mathrm{disc},u}:=\hat{h}_{\alpha}^{\mathrm{disc},1}+\tfrac{1}{\gamma}\log(u)$$ is a sample from $\hat{\mathbb{P}}_\alpha^{\mathrm{disc},u}$.
\end{prop} 

\begin{proof}
	The proof is very similar to that of \cref{P:discdisBL}, and we leave it as an exercise.
\end{proof}

\begin{definition}\label{D:unitareaqd} We call the quantum surface $(S,\hat{h}_{\alpha}^{\mathrm{disc},1},+\infty,-\infty)$ (when $\hat{h}_{\alpha}^{\mathrm{disc},1}$ has law $\hat{\mathbb{P}}^{\mathrm{disc},1}_\alpha$) a unit area $\alpha$-quantum disc, and $(S,\hat{h}_{\alpha}^{\mathrm{disc},u},+\infty,-\infty)$ (when $\hat{h}_{\alpha}^{\mathrm{disc},u}$ has law $\hat{\mathbb{P}}^{\mathrm{disc},u}_\alpha$) an $\alpha$-quantum disc with quantum area $u$. 
\end{definition}

%To summarise, there is a ``duality'' between $\alpha$ quantum wedges with $\alpha\in (Q,Q+\tfrac{\gamma}{2})$ and $(2Q-\alpha)$ quantum wedges: informally speaking, an $\alpha$ quantum disc is a doubly marked quantum surface that looks locally, at each of its marked points, like a $(2Q-\alpha)$ quantum wedge at its apex (the marked point with neighbourhoods of finite mass). Note that if $\alpha \in (Q, Q+ \tfrac{\gamma}{2})$, $2Q - \alpha \in  (Q- \tfrac{\gamma}{2}, Q)$, which does not cover the full range of thick wedges. 

%However, this duality extends to quantum discs and quantum wedges for $\alpha \in (Q, 2Q)$ and this \emph{does} cover the full range of thick wedges.

%Recall that the $\alpha$ quantum disc for such an $\alpha$ is in fact described by an infinite measure on $H^{-1}_{\mathrm{loc}}(S)$ built from a 
%$$
%\delta(\alpha)=2+\tfrac{2}{\gamma}(Q-\alpha)=2+\tfrac{2}{\gamma}((2Q-\alpha)-Q) 
%$$
%dimensional Bessel excursion measure, plus the law of the circular part of a Neumann GFF on $S$. Note that $\delta(\alpha)\in (1,2)$ for $\alpha\in (Q,Q+\tfrac{\gamma}{2})$.
 
%This is why we restrict to this range of values for $\alpha$, even though the $(2Q-\alpha)$ quantum wedge can be defined for the entire range $2Q-\alpha\in (0,Q)$, equivalently, $\alpha\in (Q,2Q)$.

\subsection{Quantum spheres}\label{SS:spheres}
\ind{Quantum spheres} 

The final quantum surface we will introduce in this chapter is the so called \textbf{quantum sphere}, which has finite area like the quantum disc, but does not have a boundary. It can therefore be thought of as the finite volume analogue of the quantum cone introduced in \cref{SS:cones}. %We will see in Section \cref{SS:idLCFT} that the quantum sphere with a particular choice of parameter can in fact be identified with the quantum sphere from \cref{S:LCFT}.
As usual we consider the parameter $\gamma\in (0,2)$ to be fixed from now on, and the definition of quantum surfaces (that is, the change of coordinates formula) is with respect to this value of $\gamma$.

Quantum spheres will be defined for a parameter $\alpha\in (Q,Q+\tfrac{\gamma}{2})$ (note the difference with the case of quantum discs). The $\alpha$-quantum sphere will be defined as a doubly marked quantum surface (now with \textbf{interior} rather than boundary marked points) which looks, locally near the marked points, like a $(2Q-\alpha)$-quantum cone near its apex, at least in a suitable range of $\alpha$. 

Recall that we defined the the $\alpha^*$-quantum cone, for $0<\alpha^*<Q$, to be the doubly marked quantum surface represented by $(\cC,h_{\cc}+h_{\rr}, +\infty,-\infty)$ where  $\cC:=\{ z = x + i y: y \in \R/(2\pi \mathbb{Z})\}$ is the infinite cylinder, $h_{\cc}$ has the law of the whole plane GFF on $\cC$ minus its average value on each vertical segment $\{ x+iy\, : \, y\in [0,2\pi] \}$, and $h_{\rr}$ is independent of $h_{\cc}$ with
	\begin{equation*}
		h_{\rr}(z) =
		\begin{cases}
			B_{s} + (\alpha^* - Q) s & \text{ if } \Re(z) = s \text{ and } s \ge 0\\
			\widehat B_{- s} + (\alpha^* - Q) s & \text{ if } \Re(z)=s \text{ and } s < 0,
		\end{cases}
	\end{equation*}
for $B=(B_t)_{t\ge 0}$ is a standard Brownian motion, and $\widehat B=(\widehat B_t)_{t\ge 0}$ is an independent Brownian motion conditioned so that $\widehat B_{t} + (Q - \alpha^*) t >0$ for all $t>0$.
	
As was the case when describing thin quantum wedges and quantum discs, we can switch perspective slightly, and (equivalently) define the $\alpha^*$-quantum cone to be the doubly marked quantum surface which, when parametrised by $\cC$ with marked points at $-\infty,+\infty$ (note the switch in order) is represented by a field with the law of $h_{\cc}+h_{\rr}$ \emph{viewed modulo horizontal translation}, where 
\begin{equation}
	\label{E:hradcone}
h_{\rr}(z)= B_{s}+(Q-\alpha^*)s 
 \quad \text{for} \quad \Re(z)=s
\end{equation}
where $B$ is a standard two-sided Brownian motion, independent of $h_{\cc}$, and viewed modulo translation (that is, $(B_t)_{t\in \R}$ is identified with $(B_{t_0+t})_{t\in \R}$ for any $t_0\in \R$).

For $\alpha\in (Q,Q+ \tfrac{\gamma}{2})$, if we let
$$
\delta_{\mathrm{cone}}(\alpha)=2+\frac{4}{\gamma}(Q-\alpha)
$$
(notice the factor two difference compared to the quantum disc case), then $\delta_{\mathrm{cone}}(\alpha)=:\delta\in (0,2)$. Moreover, by \cref{L:BessBM} and \cref{R:BessBM}, if $\nu_{\delta}^{1}$ is the $\delta$ dimensional Bessel excursion measure conditioned on reaching maximum value $1$, and $e$ is sampled from $\nu_\delta^1$, then $V^e$ defined by taking $\tfrac{2}{\gamma}\log(e)$, reparametrised to reach its maximum at time $0$ and to have infinitesimal quadratic variation $\dd t$ (as described in \cref{L:BessBM} but with $2\dd t$ replaced by $\dd t$), then we have that 
$V_0=0$ and 
\begin{align}\label{E:Vdesc}
 V_t & =B_{t}+(Q-\alpha)t  & \text{ for } t>0 \notag \\
V_{t} & =\hat{B}_{-t}+(Q-\alpha)(-t)& \text{ for } t<0
\end{align}
where $B,\hat{B}$ are independent standard linear Brownian motions defined for $t\in [0,\infty)$, started from $0$ and conditioned that $B_{t}+(Q-\alpha)t$ (resp. $\hat{B}_{t}+(Q-\alpha)t$) is negative for all $t>0$.

This leads us to the following definition.

\begin{definition}\label{D:sphere}
	Let $\alpha\in (Q,Q+ \tfrac{\gamma}{2})$ and $\nu_{\delta}^{\mathrm{BES}}$ be the Bessel excursion measure with dimension $\delta=\delta_{\mathrm{cone}}(\alpha)=2+\tfrac{4}{\gamma}(Q-\alpha)$. Let $\P^{\cC}_{\cc}$ be the law obtained from a whole plane GFF in $\cC$ by subtracting its average value on each vertical line segment, as described in \eqref{E:hCdecomp}.
	
	We define the infinite $\alpha$ sphere measure $m^{\mathrm{sphere}}_{\alpha}$ to be the measure on $H^{-1}_{\mathrm{loc}}(\cC)$ obtained by pushing forward the measure $\nu_{\delta}^{\mathrm{BES}}\otimes \P^{\cC}_{\cc}$ to $H^{-1}_{\mathrm{loc}}(\cC)$, via the map taking $(e,h_{\cc})$ to the field $h_{\cc}+V^e_{\Re(\cdot)}$, where $V$ is constructed from $e$ as described above \eqref{E:Vdesc}.
\end{definition}

An analogous remark to \cref{R:DMSextendeddisc} holds in this case. We also have the following analogue of \cref{P:discdisBL}
\begin{prop}  For $\alpha \in (Q,Q+ \tfrac{\gamma}{2})$
	\begin{equation}
		\label{E:sphere_disintegrated}
		m_{\alpha}^{\mathrm{sphere}}=c^*_\delta \int_0^\infty \mathbb{P}_{\alpha}^{\mathrm{sphere},u} u^{\tfrac{2(Q-\alpha)}{\gamma}-1} \dd u
	\end{equation} 
	where $c_{\delta}^*$ is a constant, and $\mathbb{P}_\alpha^{\mathrm{sphere},u}$ is a \emph{probability measure} on $H^{-1}_{\mathrm{loc}}(\cC)$, such that for $\mathbb{P}_\alpha^{\mathrm{sphere},u}$-almost every $h$, the bulk Gaussian multiplicative chaos measure $\mathcal{M}_h^\gamma(\cC)$ is well defined and satisfies
	$$
	\mathcal{M}_h^\gamma(\cC)=u.
	$$ 
	Moreover, if we write $h_{\alpha}^{\mathrm{sphere},1}$ for the field on $\cC$  with the law of 
	$$ 
	h_{\cc}+ V^e_{\Re(\cdot)}-\tfrac{1}{\gamma}\log \mathcal{M}^\gamma_{h_{\cc}+V^e_{\Re(\cdot)}}(\cC)  \text{ weighted by } \big(\mathcal{M}^\gamma_{h_{\cc}+V^e_{\Re(\cdot)}}(\cC)\big)^{-\tfrac{2(Q-\alpha)}{\gamma}}
	$$
	where $h_{\cc}$ is distributed according to $\mathbb{P}^{\cC}_{\cc}$ and $e$ is an independent Bessel excursion of dimension $\delta_{\mathrm{cone}}(\alpha)$ conditioned on taking maximum value $1$, then  $$
	h_{\alpha}^{\mathrm{sphere},u}:=h_{\alpha}^{\mathrm{sphere},1}+\tfrac{1}{\gamma}\log(u)
	$$ is a sample from $\mathbb{P}_{\alpha}^{\mathrm{sphere},u}$.
	\label{P:spheredisBL}\end{prop} 

\begin{proof}
	The proof closely mirrors that of \cref{P:discdisBL} and we leave it as an exercise.
\end{proof}

 \begin{definition}\label{D:unitareasphere} We call the doubly marked quantum surface $(\cC,h_{\alpha}^{\mathrm{sphere},1},-\infty,+\infty)$ (when $h_{\alpha}^{\mathrm{sphere},1}$ has law $\mathbb{P}^{\mathrm{sphere},1}_{\alpha}$) a unit area $\alpha$-quantum sphere, and $(\cC,h_{\alpha}^{\mathrm{sphere},u},-\infty,+\infty)$ (when $h_{\alpha}^{\mathrm{sphere},u}$ has law $\mathbb{P}^{\mathrm{sphere},u}_{\alpha}$) an $\alpha$-quantum sphere with area $u$.
\end{definition}

From \eqref{E:Vdesc} and \eqref{E:hradcone} with $\alpha^*=2Q-\alpha$, we see that, at least informally speaking, an $\alpha$-quantum sphere looks locally, near each of its marked points, like an $\alpha^*$-quantum cone at its marked point with neighbourhoods of finite quantum area. 

\subsection{Special cases}

\cref{T:wedge} and \cref{T:cone} tell us that the $\alpha$-quantum wedge and $\alpha$-quantum cone can be obtained as local limits of Neumann and whole plane GFFs respectively, at boundary (respectively bulk) points where a deterministic $\alpha$-log singularity is added to the field. On the other hand, we know by Girsanov's theorem, similarly to \cref{T:rooted}, that if we take a Neumann (respectively whole plane) GFF and sample a point from the boundary (respectively bulk) $\gamma$ Liouville measures, this is closely related to sampling a point from Lebesgue measure and then adding a $\gamma$-log singularity to the field at this point. As such the $\gamma$-quantum wedge and the $\gamma$-quantum cones are particularly important examples of quantum surfaces. Indeed, we will see them appear prominently in the key theorems of \cref{S:zipper} and \cref{S:MOT}. The corresponding special parameters in the case of discs and spheres are, by the duality discussed in the previous two sections, when $\alpha=4/\gamma$.

\begin{rmk}[Weights]
	\label{R:weights}
	In \cite{DuplantierMillerSheffield}, an alternative parametrisation of wedges, cones, discs and spheres is used, in terms of their so called \textbf{weight} $W$. The reason for this is that parameterising by weight behaves well (in fact additively) under operations of cutting and welding surfaces; we will see such operations in \cref{T:sliced,T:cuttingcone1,T:cuttingcone2,T:weldingcone}. The conversion from $\alpha$ to $W$ goes as follows:
	
	\begin{itemize}
		\item \textbf{Wedge}: $W=\gamma(\tfrac{2}{\gamma}+Q-\alpha)$;
		\item  \textbf{Cone}: $W=2\gamma(Q-\alpha)$;
			\item  \textbf{Disc}: $W=\gamma(\tfrac{2}{\gamma}+\alpha-Q)$; 
		\item  \textbf{Sphere}: $W=2\gamma(\alpha-Q)$.
	\end{itemize}
For the special cases mentioned above we thus have:
	\begin{itemize}
	\item \textbf{Wedge}: $\alpha=\gamma \Rightarrow W=2$;
	\item  \textbf{Cone}: $\alpha=\gamma \Rightarrow W=4-\gamma^2$;
	\item  \textbf{Disc}: $\alpha=4/\gamma \Rightarrow W=2$;
	\item  \textbf{Sphere}: $\alpha=4/\gamma \Rightarrow W=4-\gamma^2$.
\end{itemize}
\end{rmk}
\ind{Random surfaces!Weights}

\subsection{Equivalence of quantum and Liouville spheres}

In this section, we show that the notion of quantum sphere introduced in this chapter actually coincides with the unit volume Liouville sphere (coming from Liouville CFT) defined in \cref{S:LCFT}, in the following sense. 

\begin{theorem}[Equivalence of spheres]\label{T:spheresequiv}
	Fix $\alpha\in (Q,Q+\tfrac{\gamma}{2})$ and suppose that $h$ is sampled from $\mathbb{P}^{\mathrm{sphere},1}_{\alpha}$ (defined above \cref{D:unitareasphere}). Given $h$, let $w\in \mathcal{C}$ be sampled from $\cM_h$, normalised to be a probability measure. Let 
	\[
	h'=h\circ \psi_w^{-1} + Q\log|(\psi_w^{-1})'(\cdot)|
	\]
where $\psi:\mathcal{C}\to\hat{\C}$ sends $-\infty\mapsto 0$, $\infty\mapsto \infty$ and $w\mapsto 1$. Then $h'$ has the law of the unit volume Liouville sphere $h_{\beta,\mathbf{z}}^{L,1}$ from \cref{D:unitLiouville} with $\beta=(2Q-\alpha,2Q-\alpha,\gamma)$ and $\mathbf{z}=(0,\infty,1)$.
\end{theorem}

\begin{rmk}
The restriction of $\alpha$ to the range $\alpha \in (Q, Q + \tfrac{\gamma}{2})$ is not only to guarantee the existence of the unit area quantum sphere, but also to guarantee that the insertion parameter $\beta = (2Q- \alpha, 2Q - \alpha, \gamma)$ satisfies the Seiberg bounds, which is necessary in order to define $h_{\beta,\mathbf{z}}^{L,1}$. 	
In the special case mentioned above, when $\alpha=4/\gamma$ (this is possible if $\gamma\in ( {\sqrt{2}}, 2)$), this produces the unit volume Liouville sphere with all weights equal to $\gamma$. 
\end{rmk}

\cref{T:spheresequiv} was first shown in the case $\alpha=4/\gamma$ by Aru, Huang and Sun \cite{AHSspheres} (this is stated for $\gamma \in (0,2)$, but as already noted the restriction to $\gamma \in ({\sqrt{2}}, 2)$ is important {in our framework; an extension to smaller values $\gamma$ is possible using Remarks \ref{R:extendedSB} and \ref{R:DMSextendeddisc} but would require more work}). The more general statement above is implicit in the work of Ang, Holden and Sun in \cite{AHSintegrability}. The proof below is based on \cite{AHSintegrability}; compared to \cite{AHSspheres} a key idea is to not work directly with the law of the unit volume objects (which are hard to manipulate directly owing to the singularity of the conditioning) and instead work in the setting of infinite measures from which these laws arise. Since we already know the corresponding disintegration statement with respect to the ``law'' of the volume in both cases this will greatly simplify the analysis (although it lends itself less to the probabilistic intuition).

\begin{proof}[Proof of Theorem \ref{T:spheresequiv}]
	As mentioned above it will be more convenient to work in the setting of infinite measures, where both spheres are more canonically defined. To this end, our first step will be to define two natural infinite measures $M^S$ and $M^L$ (corresponding to the quantum surface and Liouville CFT perspectives respectively) and reduce the the proof to showing that these measures are identical (up to a deterministic multiplicative constant).
	
	We start with the measure $M^S$ on $\mathcal{C}\times H^{-1}(\hat \C)$, defined by 
\[
	M^S(F(h)f(w))  := \int_\R m_\alpha^\mathrm{sphere}\left( \int_{\mathcal{C}} f(w) F(h^t\circ \psi_w^{-1}+Q\log|(\psi_w^{-1})'(\cdot)|)\cM_{h^t}(\dd w) \right) \, \dd t \\
\]
	for non-negative Borel functions $F$ on $H^{-1}(\hat \C)$ and $f$ on $\R$, where for $t\in \R$, $h^t$ denotes the field $h(\cdot+t)$. That is, we ``sample'' the field $h$ on $\mathcal{C}$ according to the infinite quantum sphere measure $m_\alpha^{\mathrm{sphere}}$, and ``independently sample'' a real number $t$ from Lebesgue measure. Then we horizontally shift the field by $t$, choose $w\in \mathcal{C}$ according to the quantum area measure associated with the shifted field $h^t$, and finally take the image of $h^t$ after conformally mapping $\mathcal{C}$ to $\hat{\C}$, sending $-\infty\mapsto 0$, $\infty\mapsto \infty$ and $w\mapsto 1$. Note that this is equivalent to simply choosing $w$ according to the quantum area measure associated with $h$ and mapping $h$ to $\hat{\C}$ (as we will justify and use below) but we want the measure $M^S$ to be represented in the form above. The reason for this is because, under $m_\alpha^{\mathrm{sphere}}$, the horizontal translation is fixed so that the maximum of the field is obtained at $\Re(\cdot)=0$, while this is not the case for the Liouville field. 
	
	Next we define $M^L$ on $ H^{-1}(\hat \C) \times \mathcal{C}$ by 
\[ 
M^L(F(h)f(w))= m_{\beta,\mathbf{z}}^L(F) \int_{\mathcal{C}} f(w) \dd w
\]
where $\beta=(2Q-\alpha,2Q-\alpha,\gamma)$, $\mathbf{z}=(0,\infty,1)$ and $m_{\beta,\mathbf{z}}^L(F)$ is as defined in Remark \ref{R:mu0}. The measure  $m_{\beta,\mathbf{z}}^L$ is an infinite measure on $H^{-1}(\hat \C) $ which defines the random area Liouville CFT sphere with insertions of strength $2Q-\alpha$ at $0$ and $\infty$, and an insertion of strength $\gamma$ at $1$.

We will show that $M^S$ and $M^L$ are proportional to one another; let us first explain why this yields the proof of the theorem. Observe that by changing variables $x=w+t$ (with $u\in \cC, t \in \R$), we have that (for arbitrary non-negative measurable functions $F$ and $f$)
\begin{align*}
M^S(F(h)f(w))  & = \int_\R m_\alpha^\mathrm{sphere}\left( \int_{\mathcal{C}} f(w) F(h^t\circ \psi_w^{-1}+Q\log|(\psi_w^{-1})'(\cdot)|)\cM_{h^t}(\dd w) \right) \, \dd t\\
&=  m_\alpha^\mathrm{sphere}\left( \int_{\mathcal{C}} (\int_\R f(x-t) \, \dd t) F(h\circ \psi_x^{-1}+Q\log|(\psi_x^{-1})'(\cdot)|)\cM_{h}(\dd x) \right) ,
\end{align*}
since $h^t\circ \psi_w^{-1} = h\circ \psi_x^{-1}$ and $|(\psi_w^{-1})'(\cdot)| = |(\psi_x^{-1})'(\cdot)|$. 
So, for example, choosing $f_0(z)=p(\Re(z))$ with $\int_\R p(y) \dd y=1$ with $p$ non-negative and measurable, we get that 
\begin{align*}
M^S(F(h)f_0(w))  & =  m_\alpha^\mathrm{sphere}\left( \int_{\mathcal{C}}  F(h\circ \psi_w^{-1}+Q\log|(\psi_w^{-1})'(\cdot)|)\cM_{h}(\dd w) \right)  \\
& = c \int_{u>0} u^{\tfrac{2(Q-\alpha)}{\gamma}} \mathbb{P}_\alpha^{\mathrm{sphere},u}(\int_{\mathcal{C}} F(h\circ \psi_w^{-1}+Q\log|(\psi_w^{-1})'|) \frac{\cM_h(\dd w)}{u})
\, \dd u
\end{align*}
by \eqref{P:spheredisBL}, where $c=c_\delta^*$ is a deterministic constant. 

Now consider $M^L$. By \cref{R:mu0}, we have that 
\begin{equation}\label{E:mbetadis} 
	M^L(F(h)f_0(w))= c'  \int_{u>0} u^{\tfrac{2(Q-\alpha)}{\gamma}} \mathbb{E} (F(h_{\beta,\mathbf{z}}^{L,u})) \dd u.
\end{equation}
where $c'$ is another constant (depending only on $\gamma$ and $\alpha$) and $h_{\beta,\mathbf{z}}^{L,u}$ is the volume $u$-Liouville sphere from Chapter \ref{S:LCFT} (equivalently $h_{\beta,\mathbf{z}}^{L,1}+\gamma^{-1}\log(u)$ where $h_{\beta,\mathbf{z}}^{L,1}$ is the unit volume Liouville sphere).

Therefore, if $M^S$ and $M^L$ are proportional to one another, it follows that the law of $h_{\beta,\mathbf{z}}^{L,u}$ and that of  $h\circ \psi_w^{-1}+Q\log|(\psi_w^{-1})'(\cdot)|)$ when $(h,w)$ is sampled from $ \mathbb{P}_\alpha^{\mathrm{sphere},u} \cM_h(\dd w)/u$, 
are proportional to each other for Lebesgue almost all $u>0$. Noting that both are probability measures (indeed, the total $\cM_h$-mass of $\cC$ under $\mathbb{P}_\alpha^{\mathrm{sphere},u} $ is a.s. equal to $u$), and noting that the dependence on $u$ is continuous, it follows that these two laws are equal to one another for all $u>0$. We thus obtain the desired statement by taking $u =1$.

 It therefore remains to prove that 
$M^S$ and $M^L$ are equal as (infinite) measures on  $ H^{-1}(\hat \C) \times \cC$, up to a multiplicative constant. 
We now state three key claims, whose proofs we postpone to the end.

\begin{claim}[An identity for shifted Brownian motions with drift]
	\[
\int_\R m_{\alpha}^{\mathrm{sphere}}(F(h^t)) \, \dd t = b_1	\int_\R e^{(2Q-2\alpha)c}P^\mathcal{C}(F(h+c)) \dd c 
	\]
	for some deterministic constant $b_1$, where under the probability measure $P^{\mathcal{C}}$, $h=h_{\mathrm{circ}}+B_{\Re(\cdot)}+(Q-\alpha)|\Re(\cdot)|$, $B$ is a two-sided Brownian motion equal to $0$ at $0$, and $h_{\mathrm{circ}}$ has the law $\mathbb{P}^\mathcal{C}_{\mathrm{circ}}$, as in the definition (\cref{D:sphere}) of the unit volume quantum sphere.
\label{C:sphere1} \end{claim}

\begin{claim}[A version of Girsanov for infinite measures]
	\begin{multline} 	\int_\R e^{(2Q-2\alpha)c}P^\mathcal{C}(\int_{\mathcal{C}} f(w) F(h+c)\cM_{h+c}(\dd w)) \dd c \\ = \int_\mathcal{C} f(w)e^{\gamma(Q+\gamma/2-\alpha)|\Re(w)|} \int_\R e^{(2Q-2\alpha+\gamma)c} P^{\mathcal{C}}(F(h+\gamma G(\cdot,w)+c)) \dd c \dd w
		\end{multline}
	where $G$ is the covariance kernel of $h_{\mathrm{circ}}+B_{\Re(\cdot)}$.
\label{C:sphere2} \end{claim}

\begin{claim}[An application of the Weyl anomaly] There exists a constant $b_2\neq 0$ depending only on $\gamma$ and $\alpha$, such that for each fixed $w \in \cC$, and for all non-negative Borel functions $F$  on $H^{-1}_{\mathrm{loc}}(\mathcal{C})$
	\label{C:sphere3}
	\[ m^{L}_{\beta,\mathbf{z}}(F(h\circ \psi_w +Q\log|\psi_w'|)) =  b_2 e^{\gamma(Q+\gamma/2-\alpha)|\Re(w)|} \int_\R e^{(2Q-2\alpha+\gamma)c} P^{\mathcal{C}}(F(h+\gamma G(\cdot,w)+c)) \, \dd c \]
\end{claim}

(Note that since $h $ is viewed as an element of $H^{-1} (\hat \C)$ under $m^L_{\beta, \mathbf{z}}$, which is a subset of $H^{-1}_{\mathrm{loc}}( \C)$, it is indeed the case that $h \circ \psi_w + Q\log|\psi_w'| \in H^{-1}_{\mathrm{loc}} ( \cC)$).

\medskip Let us check how this claims imply the desired proportionality result between $M^L$ and $M^S$, which will be obtained in \eqref{eq:proportionalLiouville}. For non-negative, measurable functions $F$ and $f$, set $\tilde F(h) = F(h) \int_{\cC} f(w) \cM_h ( \dd w)$. 
Then
	\begin{align*}
\int_\R m_{\alpha}^{\mathrm{sphere}}(\int_{\mathcal{C}}f(w) F(h^t)\cM_{h^t}(\dd w) ) \dd t &= \int_\R m_{\alpha}^{\mathrm{sphere}} ( \tilde F(h^t)) \dd t\\
&= b_1 \int_{\R} e^{(2Q - 2\alpha) c} P^{\cC} ( \tilde F( h+c)) \dd c \quad \text{ (by \cref{C:sphere1})}\\
& =b_1 \int_{\R} e^{(2Q - 2\alpha) c} P^{\cC} ( F( h+c) \int_{\cC} f(w) \cM_{h+c} ( \dd w) ) \dd c\\
& =  
\frac{b_1}{b_2} \int_\mathcal{C} f(w) m_{\beta,\mathbf{z}}^L (F(h\circ \psi_w +Q\log|\psi_w'|) \dd w,
\end{align*}
by \cref{C:sphere2} and \cref{C:sphere3}. Since this holds for measurable non-negative functions $f$ and $F$, a monotone class argument shows that for all jointly measurable non-negative functions $G$ of $h$ and $w$, 
$$
\int_{\R} m_{\alpha}^{\mathrm{sphere}}(\int_{\mathcal{C}} G(h^t, w) \cM_{h^t} (\dd w) \dd t )= \frac{b_1}{b_2} \int_\mathcal{C} f(w) m_{\beta,\mathbf{z}}^L (G(h\circ \psi_w +Q\log|\psi_w'|, w)) \dd w.
$$
Applying this identity with $G(h,w)=F(h\circ \psi_w^{-1}+Q\log|(\psi_w^{-1})'|) f(w)$ shows that 
$$
\int_{\R} m_{\alpha}^{\mathrm{sphere}}(\int_{\mathcal{C}}F(h^t\circ \psi_w^{-1}+Q\log|(\psi_w^{-1})'|) f(w)) \cM_{h^t} (\dd w) \dd t =  \frac{b_1}{b_2} \int_\mathcal{C}  m_{\beta,\mathbf{z}}^L ( F(h) f(w) ) \dd w.
$$
Referring to the definition of these measures, this means
 \begin{equation}\label{eq:proportionalLiouville}
 M^S({F}(h)f(w))  = \frac{b_1}{b_2} M^L({F}(h)f(w))
 \end{equation} 
 for arbitrary non-negative and measurable functions $F,f$. As discussed, this completes the proof.
\end{proof}

\begin{proof}[Proof of \cref{C:sphere1}]
Recall that
 \begin{align*}
 m_\alpha^{\mathrm{sphere}}(F) & = \mathbb{P}^{\cC}_{\cc}\otimes \nu_{\delta}^{\mathrm{BES}}(F(h_{\cc}+V^e_{\Re(\cdot)})) ,
\end{align*}
where in the above, $h_{\cc}$ has law $\mathbb{P}^{\cC}_{\cc}$ and $V^e$ is defined from the excursion $e$ (sampled from the Bessel excursion measure $\nu^{\mathrm{BES}}_{\delta}$ with $\delta=\delta_{\mathrm{cone}}(\alpha)=2+(4/\gamma)(Q-\alpha) \in (0, 2)$ for our range of values of $\alpha\in (Q, Q + \tfrac{\gamma}{2})$) as $(2/\gamma) \log(e)$ with time parametrised to have infinitesimal quadratic variation $\dd t$ and horizontal translation fixed so its maximum value occurs at time $0$. %where the second line follows from the change of variables $c=(2/\gamma)\log(x)$, and under $\tilde{P}^{\mathcal{C}}$, $h_{\cc}$ has the same law as under $P^{\mathcal{C}}$, and $\tilde{B}$ is a two-sided Brownian motion, equal to zero at $t=0$, and conditioned that $\tilde{B}_{t}+(Q-\alpha)|t|\le 0$ for all $t\in \R$. 
	
	By invariance of the law of $h_{\cc}$ under horizontal translations, the proof of this claim therefore amounts to showing that the measure $\eta$ on the space of continuous functions from $\R$ to $\R$ defined by \[\eta(A)=\int_\R e^{2(Q-\alpha)c} \mathbb{P}(B_{\cdot}+(Q-\alpha)|\cdot|+c \in A) \dd c\] and the measure $\tilde{\eta}$ defined by \[\tilde{\eta}(A)=\int_\R \nu^{\mathrm{BES}}_\delta(V^e_{\cdot+t}\in A) \dd t \]
are equal up to a multiplicative constant. Above and for the rest of this proof we use the notation $\mathbb{P}$ and $\mathbb{E}$ for the law of two-sided Brownian motion $B$ started from $0$ at time $0$ (that is, a standard Brownian motion run forward from time $0$ joined with an independent standard Brownian motion run backwards from time $0$).

To show the equivalence of these measures, we first show that the ``marginal law'' of the function at time $0$ under $\eta$ and $\tilde{\eta}$ are the same, up to deterministic multiplicative constant.
For this, observe that for any $C\in \R$, if we let $\tau_C=\{\inf s\in \R: V^e_{s}=C\}$ then the $\nu^{\mathrm{BES}}_\delta$ law of $(V^{e}_{\tau_C+s})_{s\ge 0}$ conditioned on $\sup e \ge \exp(\gamma C/2)$ is that of $(2/\gamma) \log(\mathfrak{e}^{(C)})$ parametrised to have quadratic variation $s$ at time $s$, where $\mathfrak{e}^{(C)}$ is a $\delta$-dimensional Bessel process, started from $\exp(\gamma C/2)$ and killed upon hitting zero. This is straightforward to see in the case of It\^o's excursion theory (see also \cite[Lemma 3.4]{RhodesVargas_tail}). By \cref{L:BessBM}, the $\nu^{\mathrm{BES}}_\delta$ law of $(V^{e}_{\tau_C+s})_{s\ge 0}$ conditioned on $\sup e \ge \exp(\gamma C/2)$ is simply the law of $(C+B_s+(Q-\alpha)s)_{s\ge 0}$ where $B$ is a standard Brownian motion started from $0$ at time $0$. In particular, if $A_C$ is the set of functions which are $\ge C$ at time $0$, then by definition of $\tilde\eta$, 
\begin{align} \tilde{\eta}(A_C) & = \int_\R \nu_{\delta}^{\mathrm{BES}}(V_t^e \ge C) \dd t \notag \\
& = \int_\R \nu_{\delta}^{\mathrm{BES}}(V_t^e \ge C ,\sup e\ge \exp(\gamma C/2), \tau_C \le t) \dd t \notag\\
	& = \int_0^\infty \nu_{\delta}^{\mathrm{BES}}(V_{t+\tau_C}^e \ge C, \sup e\ge \exp(\gamma C/2)) \dd t 
	\end{align}
by Fubini and changing variables $t \to t+\tau_C$. Hence
\begin{align}
 \tilde{\eta}(A_C)
&	= \int_0^\infty \mathbb{P}(B_t+(Q-\alpha)t \ge 0) \nu_{\delta}^{\mathrm{BES}}(\sup e \ge \exp(\gamma C/2)) \dd t \notag \\
& = \nu_{\delta}^{\mathrm{BES}}(\sup e \ge \exp(\gamma C/2)) \int_0^\infty \mathbb{P}(B_t\ge (\alpha-Q)t)  \dd t \notag \\
& \propto  \int_{e^{\gamma C/2}}^\infty x^{\delta-3} \dd x \, \times \, \int_0^\infty \int_{(\alpha-Q)\sqrt{t}}^\infty \tfrac{1}{\sqrt{2\pi}} e^{-x^2/2} \dd x \dd t \quad \hfill \text{ (by \cref{L:bessel_dual})} \notag \\
& = \frac{e^{2(Q-\alpha) C}}{\tfrac{4}{\gamma}(\alpha-Q)} \frac{1}{2(\alpha-Q)^2}\quad \text{(by Fubini and using the value of $\delta$)} \notag \\
& \propto \eta(A_C) \label{E:lawstime0}
	\end{align}
where the implied constants of proportionality above do not depend on $C$. In other words if we push forward the measures $\eta$ and $\tilde{\eta}$ via the map $X\mapsto X_0$, then the resulting infinite measures on $\R$ are multiples of one another.

Furthermore, for any $C\in \R$, and any non-negative measurable functions $F$ and $G$, by the same argument as above,
\begin{multline*}
\tilde \eta ( F((X_{-s})_{s\ge 0} ) G((X_r-X_0)_{r\ge 0}) )=	\int_\R \nu_{\delta}^{\mathrm{BES}}(F((V_{t-s}^e)_{s\ge 0})G((V_{t+r}^e-V_t^e)_{r\ge 0})\mathbf{1}_{V^e_t\ge C} )\dd t  \\
%	 = & \int_\R \nu_{\delta}^{\mathrm{BES}}(F((V_{\tau_C+t-s}^e)_{s\ge 0})G((V_{\tau_C+t+r}^e-V_{\tau_C+t}^e)_{r\ge 0})\mathbf{1}_{V_{\tau_C+t}^e\ge C}) \dd t \notag \\
	=  \int_0^\infty \nu_{\delta}^{\mathrm{BES}}(F((V_{\tau_C+t-s}^e)_{s\ge 0})G((V_{\tau_C+t+r}^e-V_{\tau_C+t}^e)_{r\ge 0})\mathbf{1}_{V^e_{\tau_C+t}\ge C}) \dd t .
\end{multline*}
%where the first equality follows from a change of variables ($t\mapsto t+\tau_C$) and the second equality holds because $V_{\tau+C+t}^e\ge C$ can only occur if $t\ge \tau_C$. 
Now, as noted previously, the $\nu^{\mathrm{BES}}_\delta$ law of $(V^{e}_{\tau_C+r})_{r\ge 0}$ conditioned on $\sup e \ge \exp(\gamma C/2)$ is that of $(C+B_r+(Q-\alpha)r)_{r\ge 0}$ where $B$ is a standard Brownian motion started from $0$ at time $0$, independent of $(V^{e}_{\tau_C-s})_{s\ge 0}$. Using the Markov property at time $t$ of the above Brownian motion with drift, we can therefore rewrite the final expression above as 
\[\mathbb{E}((G(B_r+(Q-\alpha)r)_{r\ge 0})  \int_0^\infty \nu_{\delta}^{\mathrm{BES}}(F((V_{\tau_C+t-s}^e)_{s\ge 0})\mathbf{1}_{V^e_{\tau+t}\ge C}) \dd t.\]
Considering the special case where $G=1$, we deduce that
\begin{align*}
\tilde \eta ( F((X_{-s})_{s\ge 0} ) G((X_r-X_0)_{r\ge 0}) ) & = \tilde \eta ( F((X_{-s})_{s\ge 0} ) ) \mathbb{E}((G(B_r+(Q-\alpha)r)_{r\ge 0}) 
\\
& =   \tilde \eta ( F((X_{-s})_{s\ge 0} ) ) \eta ( G( X_r - X_0)_{r\ge 0}).
\end{align*}
In other words (and somewhat informally), conditionally on the value of the function at time $0$ under $\tilde{\eta}$, the future evolution is the same as under $\eta$, and it is independent of the value at time $0$ and the evolution before time $0$. Since both $\tilde{\eta}$ and $\eta$ are manifestly invariant under reversal of time, this shows that the conditional laws of the evolution under $\eta$ and under $\tilde{\eta}$, given the value at time $0$, are equal. Putting this together with \eqref{E:lawstime0} completes the proof of the claim. 
	\end{proof}

\begin{proof}[Proof of \cref{C:sphere2}]
	
	For each $c\in \R$, applying \cref{T:Girsanov} to the field $h-(Q-\alpha)|\Re(\cdot)|$ under $P^{\mathcal{C}}$ (whose distribution is that of $h_{\mathrm{circ}}+B_{\Re(\cdot)}$), we see that 
	\begin{align*} P^\mathcal{C}(\int_{\mathcal{C}} f(w) F(h+c)\cM_{h+c}(\dd w)) & = 
		e^{\gamma c} P^\mathcal{C}(\int_{\mathcal{C}} f(w) e^{\gamma(Q-\alpha)|\Re(\cdot)|} F(h+c)\cM_{h-(Q-\alpha)|\Re(\cdot)|}(\dd w)) \\
		& = e^{\gamma c}\int_\mathcal{C} f(w) F(h+c+\gamma G(\cdot,w)) \sigma_\gamma(\dd w)) 
	\end{align*}
where $G$ is (as in the statement of \cref{C:sphere2}) the covariance kernel of the field $h_{\cc}(\cdot) + B_{\Re(\cdot)}$,  and
$\sigma_\gamma(\dd w)$ is the measure $A\mapsto \mathbb{E}[\cM_{h_{\mathrm{circ}}+B_{\Re(\cdot)}} (A)]$ for $A$ a Borel subset of $\mathcal{C}$. We can compute the expected mass of this GMC as follows:
\[ \sigma_\gamma(\dd w) = \lim_{\eps\to 0} \eps^{\tfrac{\gamma^2}{2}}e^{\tfrac{\gamma^2}{2} \var(h_\eps(w)+B_{\Re(w)})} \dd w. %=\lim_{\eps\to 0} \eps^{\gamma^2/2} e^{\tfrac{\gamma^2}{2} \var(h_{\eps |\exp(w)|}^\fc(\exp(w)))} \dd w
\]
The covariance of the field $h^* (\cdot) := h_{\cc}(\cdot) + B_{\Re(\cdot)}$ is easy to compute from the one of $h^\fc$ on $\hat \C$ via the exponential map $w \in \cC \mapsto e^w \in \hat \C$; by Lemma \ref{L:Ginfty}, we get 
\begin{align*}
\var (h^*_\eps (w) ) &= \log(1/\eps) -\log |e^w| + 2 \log (|e^w| \vee 1) + o(1) \\
& = \log (1/\eps) + | \Re ( w)| + o(1).
\end{align*}
Thus, 
\[
\sigma_\gamma(\dd w) =e^{\tfrac{\gamma^2}{2}|\Re(w)|}.
\]
Substituting this into the above and integrating over $c$ completes the proof of the claim.
\end{proof}

\begin{proof}[Proof of \cref{C:sphere3}]
	Recall that $\mathbf{z}=(0,\infty,1)$. We first use M\"{o}bius invariance, \cref{T:LfieldMobius} and the specialisation to one marked point at $\infty$ \eqref{eq:mob_infinity}, to observe that if we set  $\phi_w: \hat{\C}\to\hat\C$ to be the M\"obius map $z\mapsto z/(\exp(w))$
		\begin{align*} m^{L}_{\beta,\mathbf{z}}(F(h\circ \phi_w)-Q\log|\phi_w'|))  &=
		 |\phi_w'(0)|^{-2 \Delta_{2Q-\alpha}} |\phi'_w(1)|^{-2\Delta_{\gamma}} |\phi'_w(\infty)|^{2\Delta_{2Q-\alpha}} m^{L}_{\beta,\exp(w) \mathbf{z}}(F) \\
& =		 |\exp(w)|^{2\Delta_{\gamma}} m^{L}_{\beta,\exp(w) \mathbf{z}}(F). 
		\end{align*}
		 Now by
		\cref{R:3ptexplicit} with $z=\exp(w)$ we have that $m^L_{\beta,\exp(w)\mathbf{z}}(F)$ is equal to a multiple $C$ (not depending on $w$) of
		\begin{multline}
			 (|\exp(w)|\vee 1)^{-4\Delta_\gamma+2\gamma(2Q-\alpha)} |\exp(w)|^{-\gamma(2Q-\alpha)} \times \\ \int_{u>0} \mathbb{E}(F(\tilde{h}^\fc+\gamma^{-1}(\log u - \log \cM_{\tilde h^\fc}(\C)))\cM_{\tilde{h}^\fc}(\C)^{-s}) u^{s-1} \dd u
		\end{multline} 
where $ s = \tfrac{\sum_i \alpha_i-2Q}{\gamma} = \tfrac{2Q - 2\alpha + \gamma}{\gamma}$,
	\begin{equation} \tilde{h}^\fc=h^\fc+(2Q-2\alpha)\log(|\cdot|\vee 1)-(2Q-\alpha) \log|\cdot| +2\pi \gamma G^\fc(\exp(w),\cdot)\label{eq:hfcsphere3}
	\end{equation}
	and $h^\fc$ is the whole plane GFF with zero average on the unit circle. 
	Combining the powers of $|\exp(w)|$ and $|\exp(w)|\vee 1$, noting that $|\exp(w)|=e^{\Re(w)}$, and applying the change of variables $u=\cM_{\tilde{h}^\fc}(\C) e^{\gamma c}$ in the integral, after multiple cancellations we reach the expression
	\[
		m^{L}_{\beta,\mathbf{z}}(F(h\circ \phi_w-Q\log|\phi_w'|)) = e^{(-2\Delta_\gamma + \gamma(2Q-\alpha))|\Re(w)|} \int_\R e^{(2Q-2\alpha+\gamma)c} \mathbb{E}(F(\tilde{h}^\fc + c)) \dd c
	\]
	with $\tilde{h}^\fc$ as above. 
	
	From here, notice that $\psi_w=\phi_w\circ \exp$, so that $\psi'_w(\cdot) =  \exp (\cdot) \phi'_w \circ \exp ( \cdot)$. Thus if we let 
	$$
	\tilde F(h) = F( h \circ \exp + Q \log | \exp(\cdot) | ) = F( h \circ \exp + Q \Re ( \cdot)|), 
	$$
	we can write the left hand side of the identity in \cref{C:sphere3} as 
	\begin{align*}
	m^{L}_{\beta,\mathbf{z}}&(F(h\circ \psi_w-Q\log|\psi_w'|))   = m^{L}_{\beta,\mathbf{z}}(\tilde F(h\circ \phi_w-Q\log|\phi_w'|)) \\
	& = e^{(-2\Delta_\gamma + \gamma(2Q-\alpha))|\Re(w)|} \int_\R e^{(2Q-2\alpha+\gamma)c} \mathbb{E}(\tilde F(\tilde{h}^\fc + c)) \dd c\\
	& = e^{(-2\Delta_\gamma + \gamma(2Q-\alpha))|\Re(w)|} \int_\R e^{(2Q-2\alpha+\gamma)c} \mathbb{E} ( F ( \tilde{h}^\fc \circ \exp (\cdot) + Q \Re (\cdot) + c ) ) \dd c.
	\end{align*}
	Note furthermore using \eqref{eq:hfcsphere3} that 
	\[
	  \tilde{h}^\fc\circ \exp(\cdot) + Q\Re(\cdot) = h_{\mathrm{circ}}+B_{\Re(\cdot)}+(Q-\alpha)|\Re(\cdot)| + 2\pi \gamma G^\fc( \exp(w), \exp(\cdot)) \]
	where $h_{\mathrm{circ}}$ and $B$ are as in the statement of the claim, and $2\pi G^\fc(\exp(\cdot),\exp(\cdot))$ is the covariance of $h^\fc \circ \exp = h_{\mathrm{circ}}+B_{\Re(\cdot)}$, and is therefore equal to $G$ by definition. Combining this with the fact that $-2\Delta_\gamma+\gamma(2Q-\alpha)=-\gamma Q +\gamma^2/2+2\gamma Q - \gamma \alpha=\gamma(Q+\gamma/2-\alpha)$, we obtain the statement of the claim. 
\end{proof}

%special resampling property

 \subsection{Exercises}

\begin{enumerate}[label=\thesection.\arabic*]

	\item\label{Ex:eucl_wedge} Let $D = \{z: \arg (z) \in [0, \theta]\}$ be the (Euclidean) wedge of angle $\theta$, and suppose that $\theta \in (0,2\pi)$. Let $h$ be a Neumann GFF in $D$. Show that by zooming in $(D,h)$ near the tip of the wedge, we obtain a thick quantum wedge with $\alpha = Q(\theta/\pi - 1) $ (which satisfies $\alpha<Q$ if $\theta< 2\pi$).

	\item \label{wedge_robust} Show that  \cref{T:wedge} (i) remains true in the sense of convergence in distribution with respect to the topology of uniform convergence on compacts (as opposed to total variation)
	if we replace $h$ by $h = \tilde h + \alpha \log (1/|\cdot|) + \varphi$, where $\tilde h$ is a Neumann GFF on $\H$ with some fixed  additive constant and $\varphi$ is a  function which is independent of $\tilde{h}$ and continuous at $0$. That is, show that if $h$ is as above then as $C\to \infty$, the surfaces $(\H, h+ C,0,\infty)$ converge to an $\alpha$-thick wedge in distribution.

	\item (\cite[Proposition 4.2.5]{DuplantierMillerSheffield}): show the following characterisation of quantum wedges. Fix $\alpha<Q$ and suppose that $h$ is a fixed representative of a quantum surface that is parametrised by $\H$. Suppose that the following hold:
%	is a canonical description of a quantum surface parametrised by $\H$ such that the following hold:
	
	(i) The law of $(\H, h, 0,\infty)$ (as a quantum surface with two marked points 0 and $\infty$) is invariant under the operation of multiplying its area by a constant.
That is, if we fix $C\in \R$, then  $(\H, h+C/\gamma, 0, \infty)$ has the same law as $(\H, h, 0, \infty)$.

	(ii) The total variation distance between the law of $h$ restricted to $B(0, r)$ and the law of an $\alpha$-quantum wedge field $h_{\text{wedge}}$ (in its unit circle embedding in $\H$) restricted to $B(0, r)$ tends to 0 as $r\to  0$.
	
	Then $(\H, h, 0, \infty)$ has the law of an $\alpha$-quantum wedge; more precisely $h$ has the law of $h_{\text{wedge}}$.
	
	\item\label{Ex:cone} (\textbf{Quantum cones.}) Verify \cref{L:H1Cdecomp} and give a proof of Theorem \ref{T:cone}.
	
	State and prove the analogue of Exercise \ref{Ex:eucl_wedge}.
	\ind{Quantum cones}
	
\item Show that a $\delta$ dimensional Bessel process starting from zero cannot satisfy the SDE \eqref{eq:SDEBES} on $[0,t]$ when $\delta =1$. 
	
	\item Prove that a Bessel process enjoys the Brownian scaling property: if $Z$ is a $\delta$ dimensional Bessel process with  $Z_0 = x>0$, then for all $\lambda>0$, $(Z_{\lambda t}/\sqrt{\lambda})_{t \ge 0}$ is a $\delta$ dimensional Bessel process started from $x/\sqrt{\lambda}$. 

Show that the converse is also true; if $Z$ solves the SDE $\dd Z_t = \sigma(Z_t) \dd B_t + \beta(Z_t) \dd t$ until the first hitting time of zero, with $\sigma, \beta$ locally Lipschitz on $(0, \infty)$, and if $Z$ satisfies the above Brownian scaling property, then $\sigma(x) \equiv \sigma$ is constant and $\beta(x) \propto 1/x$ for all $x>0$.

\item For an excursion $e$ (i.e., a continuous path from an interval $(0, \zeta)$ to $(0, \infty$), with $\lim_{t \to 0} e(t) = \lim_{t\to \zeta} e(t) = 0$), set
$$
I(e) = \int_0^\zeta \frac{\dd t}{ e(t)}. 
$$ 	
Verify (using the Brownian scaling property of a Bessel process, see exercise above) that 
$$
 \nu_{\delta}^{\mathrm{BES}}  ( I(e) \in \dd x) \propto x^{\delta -3 } \dd x. 
$$
Deduce that $\sum_{i: t_i \le 1} I(e_i) < \infty$ if and only if $\delta >1$ where $(t_i, e_i)$ is a Poisson point process of intensity $\dd t \otimes  \nu_{\delta}^{\mathrm{BES}}  $. Explain how this is related to the fact that the Bessel SDE \eqref{eq:SDEBES} can only be solved for $\delta >1$.

	\item \label{Ex:finitedisc} Prove \cref{L:finitedisc} for the total quantum length of $\partial \cS$, and \cref{prop:unitareaQD} for the unit area quantum disc.
	
		\item\label{Ex:sphere} (\textbf{Quantum spheres.}) Give a proof of \cref{P:spheredisBL}.

\end{enumerate}

\newpage

\section{SLE and the quantum zipper}\label{S:zipper}

% !TEX root = master.tex

 In this section we discuss some fundamental results due to Sheffield \cite{zipper}, which have the following flavour.
\begin{enumerate}

\item \cref{T:coupling}: An SLE$_\k$ curve has a `nice' coupling with $e^{\gamma h}$, when $h$ is a certain variant of the Neumann GFF. This coupling can be formulated as a Markov property analogous to the domain Markov properties inherent to random maps. It makes the conjectures about convergence of random maps toward Liouville quantum gravity plausible, and in particular  justifies that the ``correct'' relationship between $\kappa$ and $\gamma$ is $\kappa = \gamma^2$.

\item \cref{T:quantumlength}: An SLE$_\k$ curve can be endowed with a random measure which can roughly be interpreted as $e^{\gamma h} \dd \lambda$ for $\dd\lambda$ a natural length measure on the curve. In fact, the measure $\dd\lambda$ is in itself hard to define, and the exponent $\gamma$ needs to be changed slightly from $\sqrt{\k}$ to take into account the quantum scaling exponent of  the SLE curve -- see \cite{BerestyckiSheffieldSun} for a discussion -- so we will not actually take this route to define the measure. We will instead use the notion of quantum boundary length. This has the advantage that measures on \emph{either side} of the SLE$_\k$ curve can be defined without difficulty, but we will have to do a fair bit of work to show that they are the same.

\item \cref{T:sliced}: An SLE$_\k$ curve divides the upper half plane into two independent random surfaces, glued according to boundary length. Thus, SLE curves are solutions of natural random \emph{conformal welding problems}. In fact, the existence of such solutions from a complex analytic view point is a highly non-trivial problem.

\end{enumerate}

\dis{We collect some relevant background material on SLE in Appendix \ref{app:sle}. Readers unfamiliar with the theory may wish to refer to this now.}

\subsection{SLE and GFF coupling; domain Markov property} \ind{Coupling!SLE and GFF}

Here we describe one of the two couplings between the GFF and SLE. This was first stated in the context of Liouville quantum gravity (although presented slightly differently from here) in \cite{zipper}. However, ideas for a related coupling go back to two seminal papers by Schramm and Sheffield \cite{SchrammSheffield} on the one hand, and Dub\'edat \cite{Dubedat} on the other.\\

\noindent \textbf{Notational remarks:} 
\begin{itemize}\item In what follows we will use the multiplicative normalisation for our Neumann GFFs as on the left hand side of \eqref{hNLQG}. That is, such that its covariance in the bulk grows like $\log$ (rather than $(2\pi)^{-1}\log)$) near the diagonal. 

\item Unless stated otherwise, in what follows the use of bars (for example, $\bar{h}$) indicates a distribution that is considered modulo constants.
\end{itemize}

Let $\bar h$ be a Neumann GFF on $\H$ (viewed modulo constants). Let $\kappa  >0$ and let 
\begin{equation}\label{gammakappa}
\gamma = \min \left( \sqrt{\kappa}, \sqrt{\tfrac{16}{\kappa}}\right) = 
\begin{cases}
\sqrt{\kappa} ; \text{ if } \kappa \le 4\\
\sqrt{\tfrac{16}{\kappa}} ; \text{ if } \kappa \ge 4.
\end{cases}
\end{equation}

Set
\begin{equation}\label{h0}
\bar h_0 =\bar h + \varphi \quad \quad \text{ where }  \varphi(z) = \frac{2}{\sqrt{\kappa}} \log |z|;\quad  z \in \H.
\end{equation}
Hence $\bar h_0$ is a Neumann GFF from which we have \emph{subtracted} (rather than added) a logarithmic singularity at zero. (The reason for the choice of multiple $2/\sqrt{\kappa}$ will become clear only gradually.)

Let $\eta =(\eta_t)_{t\ge 0} $ be an independent chordal SLE$_\kappa$ curve in $\H$, going from 0 to $\infty$ and parametrised by half plane capacity, where  $\kappa = \gamma^2$. Let $g_t$ be the unique conformal isomorphism $g_t: \H \setminus\{\eta_s\}_{s\le t} \to \H$ such that $g_t(z) = z + 2t /z + o(1/z)$ as $z\to \infty$ (we will call $g_t$ the Loewner map). Then
$$
\frac{\dd g_t(z)}{\dd t} = \frac{2}{g_t (z) - \xi_t};  \quad \quad z \notin \{\eta_s\}_{s\le t}
$$
where $(\xi_t)_{t \ge 0}$ is the Loewner driving function of $\eta$, and has the law of $\sqrt{\kappa}$ times a standard one dimensional Brownian motion. Let $\tilde g_t(z) = g_t(z) - \xi_t$ be the \emph{centred} Loewner map.

\begin{theorem}
\label{T:coupling}
Let $T>0$ be deterministic, and set $$\bar h_T = \bar h_0 \circ \tilde g_T^{-1} + Q \log |( \tilde g_T^{-1})'|, \text{ where } Q = \frac{2}{\gamma} + \frac{\gamma}2.$$
Then $\bar h_T$ defines a distribution in $\H$ modulo constants which has the same law as $\bar h_0$.
%, and indeed for all bounded stopping time of the filtration generated by $\eta$.
\end{theorem}
(Recall the meaning of $\circ$ when dealing with generalised functions:  $$(\bar h_0\circ \tilde g_T^{-1},\rho)=(\bar h_0, |(\tilde g_T^{-1})|^2(\rho\circ \tilde{g}_T^{-1}))$$ for any test function $\rho$.)

\begin{figure}
\begin{center}
\includegraphics{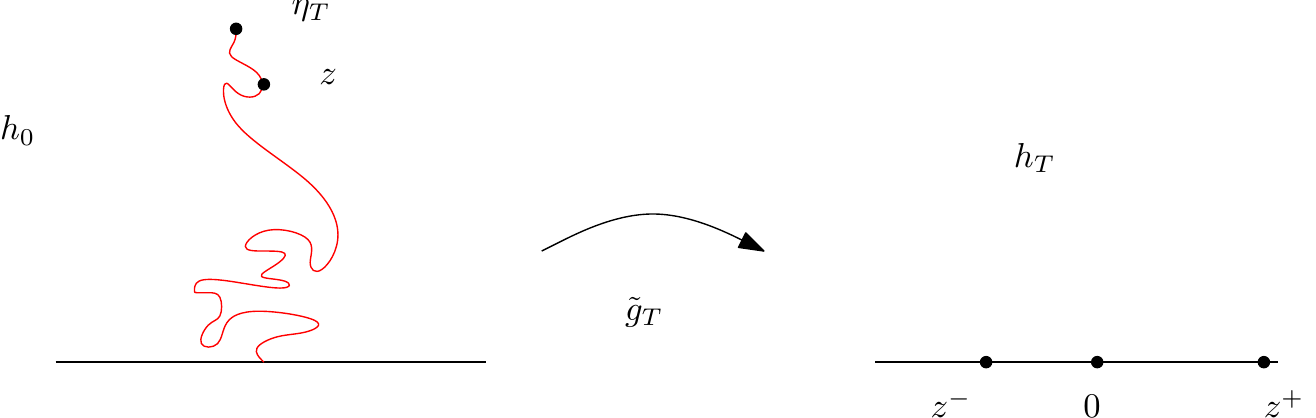}
\end{center}
\caption{Start with the field $\bar h_0$ and an independent SLE$_\kappa$ curve run up to some time $T$. After mapping $\bar h_0$, restricted to the complement of the curve $H_T$, by the Loewner map $\tilde{g}_T$ and applying the change of coordinate formula, we obtain a distribution modulo constants $\bar h_T$ in $\H$ which by the theorem has the same law as $\bar h_0$. This is a form of Markov property for random surfaces.}
\label{F:coupling}
\end{figure}

\begin{rmk}\label{R:welding}
Here we have started with a field $\bar h_0$ with a certain law (described in \eqref{h0}) and a curve $\eta$ which is \emph{independent} of $\bar h_0$. However, $\eta$ is \emph{not} independent of $\bar h_T$. In fact, we will see later on that $\bar h_T$ entirely \emph{determines} the curve $(\eta_s)_{0\le s \le T}$.  More precisely, we will see in  \cref{T:quantumlength} that when we apply the map $\tilde g_T$ to the curve $(\eta_s)_{0\le s \le T}$, the boundary lengths (measured with $\bar h_T$) of the two intervals to which $\eta$ is mapped by $\tilde g_T$ must agree: that is, on Figure \ref{F:coupling}, the $\gamma$ quantum lengths (with respect to $\bar h_T$) of $[z^-, 0]$ and $[0, z^+]$ are the same. (Note that these quantum lengths are only defined up to a multiplicative constant but their ratio is well defined, so this statement makes sense.) Then, in  \cref{T:uniqueness}, we will show that given $\bar h_T$, the curve $(\eta_s)_{0\le s \le T}$ is determined by the requirement that $\tilde g_T^{-1}$ maps intervals of equal quantum length to identical pieces of the curve $\eta$. This is the idea of \textbf{conformal welding} (we are welding $\H$ to itself by welding together pieces of the positive and negative real line that have the same quantum length).
\end{rmk} \ind{Conformal welding}

\begin{rmk}
	\label{R:welding2}
Suppose that instead of starting with $\bar{h}_0$, viewed modulo constants, we took $h_0$ to be an equivalence class representative of $\bar{h}_0$ with additive constant fixed in some arbitrary way (so that $(h_0,\rho_0)=0$ for some deterministic $\rho_0\in \cM_N$ with $\int \rho_0 =1$). Then $h_T:=h_0\circ \tilde g_T^{-1}+Q\log |( \tilde g_T^{-1})'|$ would be such that
\[h_T-(h_T,\rho_0)\overset{(\mathrm{law})}{=} h_0 \quad \text{as distributions.} \]
In other words, the laws of $h_T$ and $h_0$ would differ by a random constant.
\end{rmk}

\begin{rmk}
The proof of the theorem (and the statement which can be found in Sheffield's paper \cite[Theorem 1.2]{zipper}), involves the (centred) reverse Loewner flow $f_t$ rather than, for a fixed $t$, the map $\tilde g_t^{-1}$. In this context, the theorem is equivalent to saying that
$$
\bar h_T = \bar h_0 \circ f_T + Q \log |f_T'|, \text{ where } Q = \frac{2}{\gamma} + \frac{\gamma}2.
$$
Moreover, in this case the theorem is also true if $T$ is a bounded stopping time (for the underlying reverse Loewner flow). The current formulation of \cref{T:coupling} has been chosen because the usual forward Loewner flow is a simpler object and more natural in the context of the Markovian interpretation discussed below. On the other hand, the formulation in terms of the reverse flow will be the most useful when we actually come to prove things in this section.
\end{rmk}

\medskip \noindent \textbf{Discussion and interpretation.} Let $H_T = \H \setminus \{ \eta_s\}_{0 \le s\le T}$ and let $h_0$ be an equivalence class representative of $\bar{h}_0$ (defined by \eqref{h0}), with additive constant fixed in some arbitrary way. In the language of random surfaces,  \cref{T:coupling} (more precisely, \cref{R:welding2}) states that the random surface $(H_T, h_0|_{H_T},\eta(T),\infty)$ has the same distribution, \emph{up to multiplying areas by a random constant}, as $(\H, h_0, 0, \infty)$. This is because $h_T$ is precisely obtained from $h_0$ by mapping its restriction to $H_T$ through the centred Loewner map $\tilde g_T$ and applying the change of coordinates formula. %(the conformal covariance of Theorem \ref{T:ci}).
The meaning of ``up to multiplying areas by a random constant'' corresponds to the fact that the laws of $h_T$ and $h_0$ differ by a random constant: see  \cref{R:welding2}.

To rephrase the above, suppose we start with a surface described by $(\H, h_0,0,\infty)$. Then we explore a small portion of it using an independent SLE$_\kappa$, started where the logarithmic singularity of the field is located (here it is important to assume that $\gamma$ and $\kappa$ are related by \eqref{gammakappa}). In this exploration, what is the law of the surface that remains to be discovered after some time $T$? The theorem states that, after zooming in or out by a random amount\footnote{Recall from \cref{sec:rs_conv} that we can view the addition of a constant to the field describing a random surface, equivalently multiplying the area measure for the random surface by a constant, as ``zooming'' in or out of the surface.}, this law is the same as the original one.  Hence the theorem can be seen as a \textbf{Markov property for Liouville quantum gravity}.
\ind{Markov property (of Liouville quantum gravity)}

The fact that this invariance only holds up to additive constants for the field, or multiplicative constants for the area measure, is because the Neumann GFF is only really uniquely defined modulo constants. \emph{A more natural result comes if one replaces the Neumann GFF by a quantum wedge}, which is scale invariant by definition (meaning that if one adds a constant to the field, its law as a quantum surface does not change). In this context, we have a similar Markov property, but only if the exploration is stopped when the quantum boundary length of the curve reaches a given value: see  \cref{T:quantumlength} and \cref{T:QZ}. Of course at the moment, however, we do not even know that the quantum boundary length of SLE is well defined -- this will be addressed in \cref{sec:length_sle}.

\medskip \noindent \textbf{Connection with the discrete picture.}  This Markov property is to be expected from the discrete side of the story. To see this, consider for instance the uniform infinite half plane triangulation (UIHPT) \ind{Uniform infinite half plane triangulation} constructed by Angel and Curien \cite{Angel2003, Angel, AngelCurien}. This is obtained as the local limit of a uniform planar map with a large number of faces and a large boundary, rooted at a uniform edge along the boundary. One can further add a critical site percolation process on this map by colouring vertices black or white independently with probability 1/2 (as shown by Angel, this is indeed the critical value for percolation on such a map). We make an exception for vertices along the boundary, where those to the left of the root edge are coloured in black, and those to the right in white. This generates an interface and it is possible to use that interface to discover the map. Such a procedure is called \emph{peeling} and was used with great efficacy by Angel and Curien \cite{AngelCurien} to study critical percolation on the UIHPT. The important point for us is that conditionally on the map being discovered up to a certain point using this peeling procedure, it is straightforward to see that the rest of the surface that remains to be discovered also has the law of the UIHPT. An analogue also exists for FK models with $q\in (0,4)$ in place of critical percolation.
\ind{Peeling}

This suggests that a nice coupling between the GFF and SLE should exist, recalling the discussion of \cref{sec:RPM_LQG}.
However, identifying the exact analogue in the continuum requires a little thought.
First, observe that if one embeds the UIHPT into the upper half plane with the distinguished root edge sent to 0, there is a freedom in how the upper half plane is scaled. Roughly, it can be specified how many triangles should be mapped into the upper unit semidisc. The natural scaling limit to consider is then the one that arises by letting this number of triangles go to infinity, and rescaling the counting measure on faces appropriately. Note that such a scaling limit will be a ``scale invariant'' random surface by definition. Indeed, it is expected to be  %something called the Brownian half plane \ellen{(citation)}, which corresponds to
the $(\gamma=\sqrt{8/3})$ LQG measure associated with a certain quantum wedge.

In fact, it is known that in the abstract ``Gromov--Hausdorff--Prokhorov topology'', the UIHPQ\footnote{quadrangulation rather than triangulation here} equipped with its natural area measure converges under the rescaling described above to a metric measure space known as the Brownian half plane \cite{BMR19,GM17}. Furthermore, the aforementioned quantum wedge can be equipped with metric in such a way that it agrees in law with the Brownian  half plane as a metric measure space.
Conjectures also hold for other models of maps, and correspondingly, for wedges associated with different values of $\gamma$. This explains (arguably) why the most natural Markov property is actually the one that holds for quantum wedges.

%Of course, the above discussion only suggests that a nice coupling between the GFF and SLE must exist, when taking the scaling limit. But the strength of the logarithmic singularity at zero (the coefficient $2/\gamma = 2 / \sqrt{\kappa}$) cannot be guessed using the above arguments. Instead, this comes out of the calculations in the proof. In turn, this suggests that the uniform measure on the vertices of the UIHPT, embedded conformally (say using circle packing) into $\H$, converge in the scaling limit to $e^{\gamma \bar h_0}$ rather than $e^{\gamma \bar h}$ - in other words, this logarithmic singularity must be present in the scaling limit of discrete map models.

\ind{Reverse Loewner flow}
\begin{proof}[Proof of \cref{T:coupling}]
First, the idea is to use the \emph{reverse Loewner flow} rather than the ordinary Loewner flow $g_t(z)$ and its centred version $\tilde g_t(z) = g_t(z) - \xi_t$. Recall that while $\tilde g_t(z) : H_t \to \H$ satisfies the SDE:
$$
{\dd  \tilde g_t(z) }= \frac{2}{\tilde g_t(z)}\dd t  - \dd \xi_t,
$$
in contrast, the reverse Loewner flow is the map $f_t: \H \to H_t:= f_t(\H)$ defined by the SDE:
$$
{\dd  f_t(z)} = - \frac{2}{f_t(z)} \dd t - \dd  \xi_t.
$$
Note the change of signs in the $\dd t$ term, which corresponds to a change in the direction of time. This Loewner flow is building the curve from the ground up rather than from the tip. More precisely, in the ordinary (forward) Loewner flow, an unusual increment for $d\xi_t$ will be reflected in an unusual behaviour of the curve near its tip at time $t$. But in the reverse Loewner flow, this increment is reflected in an unusual behaviour near the origin. Furthermore, by using the fact that for any fixed time $T>0$, the process $(\xi_T- \xi_{T-t})_{ 0 \le t \le T}$ is a Brownian motion with variance $\kappa$ run for time $T>0$, the reader can check that $f_T = \tilde g_T^{-1}$ in distribution. Note that this is not necessarily true if $T$ is a stopping time: we will see an example of this later on.

\begin{lemma}\label{L:mart}
Suppose that $\gamma>0$ and $\kappa > 0$ are arbitrary. For $z\in \H$, let
$$
M_t = M_t(z): = \frac{2}{\sqrt{\kappa}} \log  |f_t(z)| + Q \log | f'_t(z)|; \qquad Q=\frac{2}{\gamma}+\frac{\gamma}{2}.
$$
Then for any fixed $z$, $(M_t(z); t\ge 0)$ is a continuous local martingale (with respect to the filtration generated by $\xi$) if and only if $ \gamma^2 = \kappa$ or $\gamma^2 = 16/\kappa$. (Thus if also $\gamma<2$, this holds if and only if $\gamma$ and $\kappa$ are related via \eqref{gammakappa}).
Furthermore, if $z,w \in \H$, then the quadratic cross variation between $M(z)$ and $M(w)$ satisfies
$$
\mathrm{d} [M(z), M(w)]_t = 4 \Re (\frac{1}{f_t(z)}) \Re(\frac1{f_t(w)}) \dd t.
$$
In particular, 
\begin{equation}\label{eq:Mzqv}
[M(z)]_t \le  \frac{4}{(\Im(z))^2}t
\end{equation}
for all $t$.
\end{lemma}

\begin{proof}
Set $Z_t = f_t(z)$. Then $\dd Z_t = -2/Z_t \dd t - \dd \xi_t$. Set $M^*_t =  \frac{2}{\sqrt{\kappa}} \log  f_t(z) + Q \log f'_t(z)$, so that $M_t = \Re(M^*_t)$. Applying It\^o's formula we see that
\begin{align*}
\dd \log Z_t  & = \frac{\dd Z_t}{Z_t}  - \frac{1}{2}\frac{\mathrm{d} [\xi]_t}{Z_t^2}\\
& = -\frac{\dd \xi_t}{Z_t} + \frac1{Z_t^2}( -2- \kappa/2) \dd t.
\end{align*}
To obtain $\dd  f_t'(z)$ we differentiate $\dd  f_t(z)$ with respect to $z$; the term $\dd \xi_t$ does not contribute to the derivative in $z$ (since it is the same driving function $\xi$ for different values of $z$). We find that
$$
\dd  f'_t(z) =  2 \frac{f'_t(z) }{Z_t^2} \dd t,
$$
and therefore
\begin{align*}
\dd  \log f_t'(z) & = \frac{\dd  f_t'(z)}{f'_t(z)}  =  \frac{2}{Z_t^2} \dd t.
\end{align*}
Putting the two pieces together we find that
\begin{equation}\label{E:c_mart}
\dd M^*_t = -\frac{2\dd \xi_t}{\sqrt{\kappa} Z_t} + \frac{2}{ Z_t^2} \left(\frac{1}{\sqrt{\kappa}} (-2 - \kappa /2) +  Q\right) \dd t.
\end{equation}
The $\dd t$ term vanishes if and only if $2/\sqrt{\kappa} + \sqrt{\kappa}/2 = Q$. Clearly this happens if and only if $\gamma = \sqrt{\kappa}$ or $\gamma = \sqrt{16/\kappa}$.

Furthermore, taking the real part in \eqref{E:c_mart}, if $z,w$ are two points in the upper half plane $\H$, then the quadratic cross variation between $M(z)$ and $M(w)$ is a process which can be identified as
$$
\dd \, [M(z), M(w)]_t = 4 \Re (\frac{1}{f_t(z)}) \Re(\frac1{f_t(w)}) \dd t,
$$
and so \cref{L:mart} follows.

To prove \eqref{eq:Mzqv}, note that
$$
\mathrm{d}[M(z)]_t=4(\Re(\frac{1}{f_t(z)})^2=4\frac{(\Re(f_t(z)))^2}{|f_t(z)|^4} \le \frac{4}{|f_t(z)|^2}\le \frac{4}{(\Im(f_t(z)))^2} \le \frac{4}{(\Im(z))^2}
 $$
 where we use the fact that for reverse SLE, $t\mapsto \Im (f_t(z))$ is increasing with time (as can be seen readily using the reverse Loewner equation for instance). 
\end{proof}

One elementary but tedious calculation shows that if $$G_t (z,w) = G^\H_N(f_t(z), f_t(w))=-\log(|f_t(z)-\overline{f_t(w)}|)-\log(|f_t(z)-f_t(w)|)$$ then $G_t (z,w)$ is a finite variation process and furthermore:

\begin{lemma}[Hadamard's formula]
\ind{Hadamard's formula}
\label{L:green_crossvar}
We have that
$$
\dd G_t(z,w) = - 4 \Re (\frac{1}{f_t(z)}) \Re(\frac1{f_t(w)}) \dd t.
$$
In particular, $\dd \, [M(z), M(w)]_t = - \dd G_t(z,w)$.
\end{lemma}

\begin{proof}This is proved in \cite[Section 4]{zipper}.
We encourage the reader to skip the proof here, which is included only for completeness. (However, the result itself will be quite important in what follows.)

Set $X_t = f_t(z)$ and $Y_t = f_t(w)$. From the definition of the Neumann Green function,
\begin{align*}
 \dd  G_t(x,y) & =  - \dd  \log (|X_t - \bar Y_t |) - \dd  \log( |X_t - Y_t|)\\
& = -  \Re ( \dd  \log (X_t - \bar Y_t)) - \Re ( \dd  \log (X_t - Y_t)).
\end{align*}
Now, $\dd X_t = (2/X_t) \dd t - \dd \xi_t $ and $\dd Y_t = (2/Y_t) - \dd \xi_t$ so taking the difference
$$
\mathrm{d} (X_t - Y_t) = \frac{2}{X_t} \dd t - \frac2{Y_t} \dd t = 2\frac{Y_t - X_t}{X_t Y_t} \dd t
$$
and so
$$
\dd  \log (X_t - Y_t) = -\frac2{X_t Y_t} \dd t; \ \ \
\dd  \log (X_t - \bar Y_t) = -\frac2{X_t \bar Y_t} \dd t.
$$
Thus we get
\begin{equation}\label{greendiff}
 \dd G_t(x,y) = -  2 \Re (\frac1{X_t Y_t} + \frac1{X_t \bar Y_t}) \dd t.
\end{equation}
Now, observe that for all $x,y \in \C$,
$$
\frac{1}{xy} + \frac1{x\bar y}  = \frac{\bar x \bar y + \bar x y }{|xy|^2} = \frac{\bar x ( \bar y + y )}{|xy|^2}
 = \frac{ 2 \Re(y)}{|xy|^2} \bar x.
$$
Therefore, plugging into \eqref{greendiff}:
$$
\dd G_t(x,y) = -4 \frac{\Re(X_t)\Re(Y_t)}{|X_tY_t|^2} = - 4 \Re( \frac1{X_t}) \Re( \frac1{Y_t})
$$
as desired.
\end{proof}

Equipped with the above two lemmas, we prove  \cref{T:coupling}. Set $\bar h_0 = \bar h + \varphi=\bar h + \frac{2}{\sqrt{\kappa}} \log |z|$, and let $(f_t; t\ge 0)$ be an independent reverse Loewner flow as above. Define
$$
\bar h_t = \bar h_0 \circ f_t + Q \log |f_t'|.
$$
Then, viewed as a distribution modulo constants, we claim that:
\begin{equation}\label{E:goalcoupling}
\text{$\bar h_t$ has the same distribution as $\bar h_0$. }
\end{equation}
Let $\rho$ be a test function with zero average, so $\rho \in \bar \cD(\H)$. To prove \eqref{E:goalcoupling}, it suffices to check that $(\bar h_t, \rho)$ is a Gaussian with mean $(\varphi,\rho)$ and variance as in \eqref{eq:GFN}, that is, $\sigma^2 = \int \rho(\dd z) \rho(\dd w) G(z,w)$ where $G(z,w)$ is a valid choice of covariance for the Neumann GFF in $\H$.

To do this, we take conditional expectations given $\cF_t = \sigma(\xi_s, s \le t )$ (note that $f_t$ is measurable with respect to $\cF_t$), and obtain
\begin{align}
\E[ e^{i (\bar h_t, \rho)}| \cF_t ] & = \E [ e^{i (\bar h_0 \circ f_t + Q \log |f'_t| , \rho)} | \cF_t] \nonumber \\
& = e^{i (\frac{2}{ \sqrt{\kappa}} \log |f_t| + Q \log |f'_t|, \rho)} \times \E(e^{i ( \bar h \circ f_t, \rho)} | \cF_t) \nonumber \\
& = e^{i M_t (\rho)} \E[ e^{i (\bar h \circ f_t, \rho)} | \cF_t], \label{eq:conditionalexpectationzipper}
%& = \E[ e^{i M_t(\rho)} \times e^{-\frac12 \int \rho(dz) \rho(dw) G_t(z,w)} ]
\end{align}
where 
\begin{equation}\label{eq:Mtrho}
M_t(\rho) = \int_{\H} M_t(z) \rho(z) \dd z.
\end{equation}
Note that by Fubini's theorem for conditional expectations, and \eqref{eq:Mzqv} (which implies that $\E (M_t(z)^2) \le Ct $ for some constat depending only on the support of $\rho$ in $\H$), $M_t(\rho)$ is in fact a martingale. 

Now we evaluate the term in the conditional expectation in \eqref{eq:conditionalexpectationzipper}.
By definition of $\bar h \circ f_t$, the term $(\bar h  \circ f_t, \rho)$ can be computed almost surely by changing variables, that is, is equal to $(\bar h, \rho_t)$, where the corresponding ``integration'' takes place on $f_t(\H) = H_t$ and
 where 
$$  \rho_t (z) =  |(f_t^{-1})' (z)|^2\rho \circ f_t^{-1}(z)$$
We may view $H_t$ as a subset of $\H$, and note that the test function $\rho_t$, which is defined a priori on $H_t$, can be extended to $\H$ by setting it to zero on $\H \setminus H_t = K_t$. Then this test function also has mean zero on $\H$ (by change of variable), and we deduce that $(\bar h \circ f_t,\rho) = (\bar h, \rho_t)$ is Gaussian with mean zero and variance
%By conformal invariance, $\bar h \circ f_t$ is a Neumann GFF (modulo constants) in $H_t$. The integral $(\bar h \circ f_t, \rho)$ is \emph{a priori} an integral over all of $\H$, but since the curve $\eta$ has a.s. zero Lebesgue measure, we can view it as an integral only over $H_t$ (and note that the restriction of $\rho$ to $H_t$ also has average value zero). Therefore, given $\cF_t$,
%By definition, $(\bar{h}_0\circ f_t, \rho)=(\bar{h}_0, |(f_t^{-1})'|^2\rho\circ f_t^{-1})$, which is a Gaussian random variable with mean zero and variance
\begin{align*}
\var (\bar h, \rho_t) &= \iint_{\H^2} \rho_t(z) \rho_t(w) G^{\H}_N (z,w) \dd z \dd w\\
& = 
\int_{H_t}\int_{H_t} G^{\H}_N(z,w) |(f_t^{-1})'(z)|^2|(f_t^{-1})'(w)|^2(\rho\circ f_t^{-1})(z) (\rho\circ f_t^{-1})(w) \dd z\dd w\\
&=\int_{\H}\int_\H G_t(z,w) \rho(z) \rho(w) \dd z \dd w.
\end{align*}
by change of variables, where we recall that $G_t(z,w):=G^\H_N(f_t(z),f_t(w))$.  
%Intuitively, the field $\tilde h \circ f_t$ is a Gaussian field with covariance $G(f_t(z),f_t(w))$, which allows us to compute the conditional expectation readily. This takes a bit of justification, as follows.
%Here $\tilde h \circ f_t$ is just the push-forward of the distribution (modulo constants) $\tilde h$ by the map $f_t$. By definition, if $F$ is a distribution and $T:D \to D'$ a conformal transformation, then the distribution $F\circ T$ is defined on $D'$ to be such that when integrating against a given test function $\rho \in \cD(D')$,
%$$
%( F \circ T , \rho) := (F, (\rho \circ  T^{-1})|(T^{-1})'|).
%$$
%This extends to distribution modulo constants readily (note that if $\rho$ has mean zero in $D'$, then so does %$(\rho \circ  T^{-1})|(T^{-1})'|$ in $D$).
%
%Therefore, $( \tilde h \circ f_t, \rho)$ is equal to $(\tilde h, \rho \circ f_t^{-1} |(f_t^{-1})'| )$. (Note that we apply the above to $\rho|_{H_t}$, which has mean zero because the Lebesgue measure of $\eta$ is zero a.s.) Consequently, given $\cF_t$, this is a Gaussian random variable with mean zero and variance
%$$\int_{H_t} G(x,y) (\rho \circ f_t^{-1} (x)) |(f_t^{-1})' (x)|  (\rho \circ f_t^{-1} (y)) |(f_t^{-1})' (y)| dy = \int_{\H} G(f_t(x), f_t(y)) \rho(x) \rho(y) dy $$
%after changing variables.
Hence, recalling our notation for $M_t(\rho)$ in \eqref{eq:Mtrho},  we deduce that
\begin{equation}\label{eqn:laplace}
\E( e^{i (\bar h_t, \rho)}| \cF_t ) = e^{i M_t(\rho)} \times e^{-\frac12 \iint \rho(z) \rho(w) G_t(z,w)\dd z\dd w } .
\end{equation}
%, and we have used that given $\cF_t$, $h\circ f_t$ is (by conformal invariance) a Neumann GFF in $H_t$, with covariance given precisely by $G_t$ (also by conformal invariance).
%Moreover, %by an easy application of the conditional Fubini theorem, it follows that $M_t(\rho)$ is  a local martingale and its quadratic variation is given by
%an application of Fubini's theorem gives that
Next, we claim that
\begin{equation}\label{eq:Mrhoqv}
\mathrm{d} [M(\rho)]_t = \int \rho(z) \rho(w) \mathrm{d} [M(z), M(w)]_t \dd z \dd w.
\end{equation}
To see this, we need to verify that 
\begin{align*}
Q_t & := M_t(\rho)^2 -   \int \rho(z) \rho(w)  [M(z), M(w)]_t \dd z \dd w\\
& = \int \rho(z) \rho(z) \left( M_t(z) M_t(w) - [M(z), M(w) ]_t \right) \dd z \dd w
\end{align*}
is a local (and in fact a true) martingale. Indeed, from \eqref{eq:Mzqv} and the Kunita--Watanabe inequality, we know that $M_t(z) M_t(w) - [ M(z), M(w)]_t$ is a true martingale, so
\begin{align*}
 \E \left[\left. Q_t \right|\cF_s \right] & = \int \rho(z) \rho (w) \E \left( \left. M_t(z) M_t(w) - [M(z), M(w) ]_t  \right | \cF_s\right) \dd z \dd w \\
& =  \int \rho(z) \rho (w) \left( M_s (z) M_s(w)   - [ M(z) , M(w)]_s \right) \dd z \dd w.
\end{align*}
Here in the first line  we have used Fubini's theorem for conditional expectations, whose use is once again justified by \eqref{eq:Mzqv} and the Kunita--Watanabe inequality. This proves \eqref{eq:Mrhoqv}.

Therefore, by \cref{L:green_crossvar},
$$
\int \rho(x) \rho(y) G_t(x,y) \dd x \dd y =  \int \rho(z) \rho(w) G_N^\H(z,w) \dd z \dd w - [M(\rho)]_t.
$$
Combining with \eqref{eqn:laplace} finally implies that
$$
\E( e^{i (\bar h_t, \rho)} ) = e^{-\frac12 \int \rho(z) \rho(w) G_N^\H(z,w) \dd z \dd w} \E( e^{i M_t (\rho) + \frac12 [M(\rho)]_t }).
$$
To conclude we observe that by It\^o's formula, $e^{i M_t (\rho) + \frac12 [M(\rho)]_t }$ is an exponential local martingale, and it is not hard to see that it is a true martingale ($[M(\rho)]_t\lesssim \int |\rho(z)| |\rho(w)| G^\H_N(z,w)$, which is finite for all $t$). We deduce that the expectation in the right hand side above is equal to $\E(e^{i(M_0,\rho)})=e^{i(\varphi,\rho)}$, and therefore
$$
\E( e^{i (h_t, \rho)} ) = e^{-\frac12 \int \rho(z) \rho(w) G^\H_N(z,w)\dd z\dd w} e^{i(\varphi,\rho)}.
$$
This proves \eqref{E:goalcoupling}. Arguing that $f_t$ and $\tilde g_t^{-1}$ have the same distribution finishes the proof of the theorem.
\end{proof}

%\begin{rmk}
%The reason why we expect $M_t$ to be a martingale for each $z$ is that it is the expected height of the field $h_t$ at $z$ given the filtration generated by the Loewner flow (or equivalently the driving Brownian motion) up to time $t$.
%\end{rmk}

\begin{rmk}
As mentioned earlier, since the proof relies on martingale computation and the optional stopping theorem, the theorem remains true if $T$ is a (bounded) stopping time for the \emph{reverse} Loewner flow.
\end{rmk}

\begin{rmk}
This martingale is obtained by taking the real part of a certain complex martingale. Taking its imaginary part (in the case of the forward flow) gives rise to the imaginary geometry developed by Miller and Sheffield in a striking series of papers \cite{MSIG1,IG2,IG3,IG4}.
\end{rmk}

\subsection{Quantum length of SLE }\label{sec:length_sle} \ind{Quantum length (of SLE)}
We start with one of the main theorems of this section, which allows us, given a chordal SLE$_\kappa$ curve and an independent Neumann GFF, to define a notion of quantum length of the curve unambiguously. The way this is done is by mapping the curve down to the real line with the centred Loewner map $\tilde g_t$, and using the quantum boundary measure $\mathcal{V}$ (associated with the image of the GFF via the change of coordinates formula) to define the length. % introduced in Chapter \ref{S:SIsurfaces}.
However, when we map away the curve using the map $\tilde{g}_t$, each point of the curve corresponds to two points on the real line (except for the tip of the curve which is sent to the origin since we consider the centred map). Hence, to define the length of the curve unambiguously, we first need to know that measuring the length on one side of $0$ almost surely gives the same answer as measuring the length on the other side of $0$.

This is basically the content of the next theorem. For ease of proof, the theorem is stated in the case where $h$ is not a Neumann GFF but rather the field of a certain wedge. However, we will see (\cref{C:length}) that this is no loss of generality.% since we can couple this wedge field with a Neumann GFF plus a log singularity and a large constant, so they are equal in any compact set with arbitrarily high probability (\cref{T:wedge}), this is no loss of generality.

\begin{theorem}
\label{T:quantumlength} Let $0 < \gamma<2$ and let $(\H, h, 0, \infty)$ be %any equivalence class representative of an
an $\alpha$-quantum wedge in the unit circle embedding (see \cref{R:unitcircleembedding}), with  $\alpha = \gamma - 2/\gamma$. Let $\zeta$ be an independent SLE$_\kappa$ with $\kappa = \gamma^2$. Let $\tilde g_t$ be the (half plane capacity parametrised) centred Loewner flow for $\zeta$, fix $t>0$, and consider the distribution $h_t = h\circ \tilde{g}_t^{-1} + Q \log |(\tilde{g}_t^{-1})'|$ as before. Let ${\mathcal{V}}_{h_t}$ be the boundary Liouville measure on $\R$ associated with the distribution $h_t$.
Finally, given a point $z \in \zeta([0,t])$, let $z^- < z^+$ be the two images of $z$ under $\tilde g_t$. Then
$$
{\mathcal{V}}_{h_t} ( [z^-, 0] ) = {\mathcal{V}}_{h_t}([0 , z^+]),
$$
almost surely for all $z \in \zeta([0,t])$.
\end{theorem}

\begin{rmk}\label{R:blmakessense}
	By \cref{R:welding2} and the fact that the slit domain formed by an SLE$_\kappa$ with $\kappa<4$ is almost surely H\"{o}lder continuous, we see that a Neumann GFF (with arbitrary normalisation) plus a $(\gamma-2/\gamma)$ log-singularity in such a slit domain does satisfy the conditions of \cref{def:BLM}. That is, the quantum boundary length on either side of the curve is well defined by mapping down to the real line. Since a $(\gamma-2/\gamma)$-quantum wedge in the unit circle embedding has the same law when restricted to $B(0,1)\cap \H$ as such a Neumann GFF (with normalisation fixed so that it has mean value 0 on the upper unit semicircle) this implies that the field $h$ of the above theorem also satisfies the conditions of \cref{def:BLM}, at least when restricted to $B(0,1)$. Scale invariance implies that this holds when the field is restricted to any large disc.  In other words, the boundary Liouville measure ${\mathcal{V}}_{h_t}$ for $h_t$ is well defined.
\end{rmk}

\begin{cor}\label{C:length}
\cref{T:quantumlength} is still true when $h$ is replaced by a Neumann GFF on $\H$, with arbitrary normalisation. Indeed, by the discussion in the previous remark, it is true until the curve exits the upper unit semidisc, when the normalisation for the GFF is such that it has average $0$ on the upper unit semicircle. This extends to arbitrary normalisations, since two Neumann GFFs with different normalisations (can be coupled so that they) differ by a random additive constant. Finally, scaling removes the need to restrict to the unit semidisc.   %whose normalisation is fixed so that $(h,\rho_0)=0$ for some $\rho_0\in \cM\setminus \tilde{\cM}$ supported away from the origin. Indeed,
\end{cor}

\begin{definition}
\label{D:quantumlength}
The quantity  $${\mathcal{V}}_{h_t} ( [\zeta(s)^-, 0] ) = {\mathcal{V}}_{h_t}([0 , \zeta(s)^+])$$ is called the \textbf{quantum length} of $\zeta([s,t])$ in the wedge $(\H, h, 0, \infty)$.
\end{definition}

\medskip \noindent \textbf{False proof} of  \cref{T:quantumlength}. The following argument does not work but helps explain the idea and why wedges are a useful notion. Let $\zeta$ be the infinite SLE$_\kappa$ curve parametrised by half plane capacity. Let $L(t)={\mathcal{V}}_{h_t}([\zeta(t)^{-},0])$ be the quantum length of left hand side of the curve $\zeta$ up to time $t$ (measured by computing the boundary quantum length on the left of zero after applying the map $\tilde g_t$) and likewise, let $R(t)$ be the quantum length of the right hand side of $\zeta$. Then it is \emph{tempting} (but wrong) to think that, because SLE is stationary via the domain Markov property, and the Neumann GFF is invariant by \cref{T:coupling},  $L(t)$ and $R(t)$ form processes with stationary increments. If that were the case, we would conclude from Birkhoff's ergodic theorem for stationary increments processes that $L(t)/ t $ converges almost surely to a possibly random constant, and $R(t)/ t$ converges also to a random constant. We would deduce that $L(t) / R(t)$ converges to a possibly random constant. Finally, we would argue that this constant cannot be random because of tail triviality of SLE (that is, of driving Brownian motion) and in fact must be one by left-right symmetry. On the other hand by scale invariance, the distribution of $L(t)/R(t)$ is constant. Hence  we would deduce that $L(t) = R(t)$.

This proof is wrong on at least two counts: first of all, it is not true that $L(t)$ and $R(t)$ have stationary increments.  This does not hold, for instance, because $h$ loses its stationarity (that is, the relation $h_T =h_0$ in distribution does not hold) as soon as a normalisation is fixed for the Neumann GFF. Likewise the scale invariance does not hold in this case. This explains the importance of the concept of wedges, for which scale invariance holds by definition, as well as a certain form of stationarity (see  \cref{T:QZ}). These properties allow us to make the above proof rigorous.

\subsection{Proof of Theorem \ref{T:quantumlength}}

Essential to the proof of \cref{T:quantumlength} is the definition of two stationary processes: the \emph{capacity zipper} and the \emph{quantum zipper}. As in the original paper of Sheffield \cite{zipper}, once the existence and stationarity of these processes is proven, \cref{T:quantumlength} follows relatively easily (in fact, using a similar argument to the ``false proof'' above).

In order to simplify notation in what follows, whenever $f$ is a conformal isomorphism and $h$ is a distribution or distribution modulo constants, we write
\begin{equation}
\label{eqn:coc_not}
f(h):=h\circ f^{-1} + Q \log|(f^{-1})'|.
\end{equation}

From now on we assume that $\gamma\in (0,2)$ is fixed, $Q=Q_\gamma$, and $\kappa=\gamma^2$. Recall that $\bar\cD_0'(\H)$ denotes the space of distributions modulo constants on $\H$, and we write $C([0,\infty),\H)$ for the space of continuous functions from $[0,\infty)$ to $\H$.

\begin{theorem}[Capacity zipper]\label{T:CZ}
There exists a two-sided \textbf{stationary} process $(\bar h^t,\eta^t)_{t\in \R}$, taking values in $\bar \cD_0'(\H) \times C([0, \infty);\H)$, such that:
\begin{itemize}
	\item \emph{\textbf{(Marginal law)}} $(\bar{h}^0,\eta^0)$ has the law of a Neumann GFF (modulo constants) plus the function $\varphi(z)=\frac{2}{\gamma}\log|z|$, together with an independent SLE$_\kappa$;
	\item \emph{\textbf{(Positive time)}} there exists a family of conformal isomorphisms $(f_t)_{t\ge 0}: \H\to \H \setminus \eta^t([0,t])$, whose (marginal) law is that of a reverse SLE$_\kappa$ Loewner flow parametrised by capacity, and such that $\bar h^t|_{\H\setminus \eta^t([0,t])}=f_t(\bar{h}^0)$ and $\eta^t([t,\infty))=f_t(\eta^0)$ for all $t \ge 0$;
	\item \emph{\textbf{(Negative time)}} for $t<0$, if $\tilde{g}_{-t}$ is the centred Loewner map corresponding to $\eta^0([0,-t])$ then $\eta^t=\tilde g_{-t}(\eta^0)$ and $\bar{h}^t=\tilde g_{-t}(\bar{h}^0)$.
\end{itemize}
\end{theorem}

Thus given a field $\bar h^0$ and an independent SLE$_\kappa$ infinite curve $ \eta^0$, we can either ``zip it up'' (weld it to itself) to obtain the configuration $(\bar h^t, \eta^t)$ for some $t>0$, or ``zip it down'' (cut it open along $ \eta^0$) to obtain the configuration $(\bar h^t, \eta^t)$ for some $t<0$. See \cref{F:CZ}. Beware that the relation between time $t$ and time $0$ is opposite to that of \cref{T:coupling} -- hence the change in notation from subscripts to superscripts for the time index.

 Also note that for $t>0$, $\bar{h}^t|_{\H\setminus \eta^t([0,t])}$ uniquely defines $\bar{h}^t$ as a distribution modulo constants on $\H$ (since $\eta^t([0,t])$ is independent of $\bar{h}^t$ and has Lebesgue measure zero). The term ``capacity'' in the definition refers to the fact that in any positive time $t$, we are zipping up a curve with $2t$ units of half plane capacity.
\begin{figure}
	\begin{center}
		\includegraphics{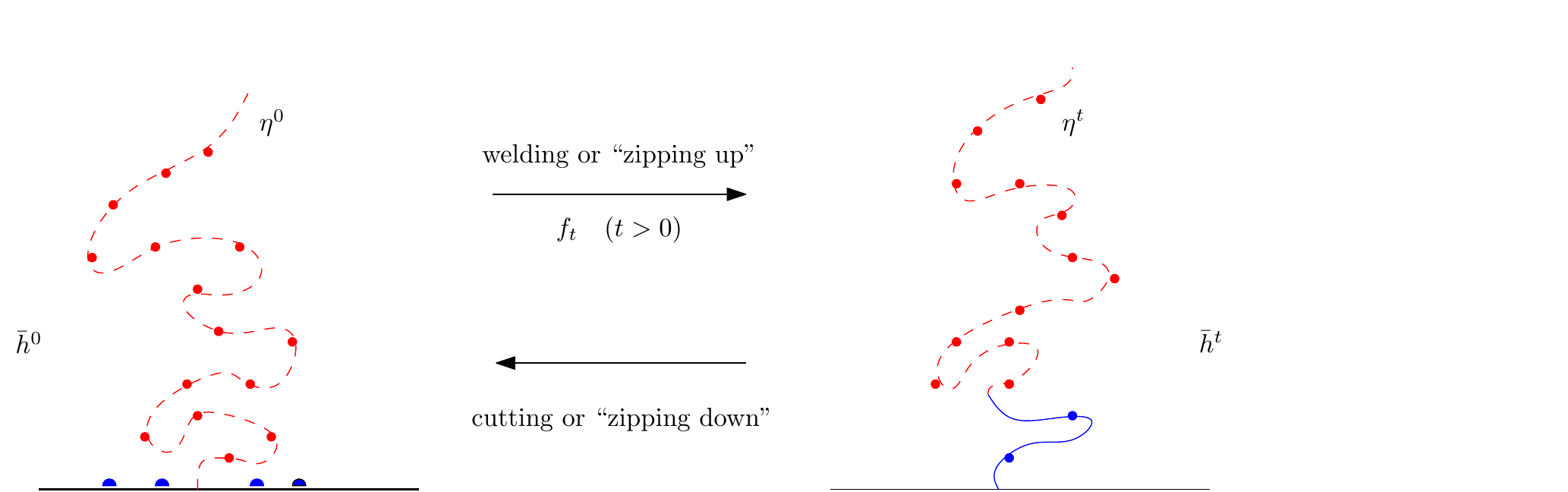}
	\end{center}
\caption{The capacity zipper}
	\label{F:CZ}
\end{figure}
\begin{rmk}\label{R:cznormed}
	Note that the capacity zipper of \cref{T:CZ} is defined to be a process taking values in $\bar \cD_0'(\H)\times C([0,\infty),\H)$. However, we can also define from $(\bar{h}^0,\eta^0,(f_t)_{t\ge 0})$ a version $(\tilde{h}^t,\eta^t)_{t\ge 0}$ of the capacity zipper indexed by positive times and taking values in $\cD_0'(\H) \times C([0,\infty),\H)$. That is, so that the field at any time is a distribution, not just a distribution modulo constants.
	
	To do this, we can just fix a normalisation of $\bar{h}^0$ to obtain $\tilde{h}^0\in \cD'_0(\H)$, and then for $t>0$ set $(\tilde{h}^t,\eta^t):=(f_t(\tilde{h}^0),\H\setminus f_t(\H\setminus \eta^0))$. Note that this process will be no longer stationary: for given $t$, $\tilde{h}^t$ will have the law of $\tilde{h}^0$ plus a random constant.
\end{rmk}

Now we move on to the definition of the \textbf{quantum zipper}. For this, we need the notion of doubly marked surface curve pair. This is just an extension of the definition of doubly marked surface, when the surface comes together with a chordal curve.
More precisely, suppose that for $i=1,2$, $D_i$ is a simply connected domain with marked boundary points $(a_i,b_i)$, $h_i$ is a distribution in $D_i$, and $\eta_i$ is a simple curve (considered up to time reparametrisation) from $a_i$ to $b_i$ in $D_i$. 
\begin{definition}\label{D:surfacecurve}
We say that $(D_1,h_1,a_1,b_1,\eta_1)$ and $(D_2,h_2,a_2,b_2,\eta_2)$ are equivalent if there exists a conformal isomorphism $f:D_1\to D_2$ such that $h_2=f(h_1)$, $a_2=f(a_1)$, $b_2=f(b_1)$ and $\eta_2=f(\eta_1)$. A doubly marked surface curve pair (from here on in just \emph{surface curve pair}) is an equivalence class of $(D,h,a,b,\eta)$ under this equivalence relation.
\end{definition}

\begin{theorem}[Quantum zipper]\label{T:QZ}
	There exists a two-sided process $$((\H,h^t,0,\infty,\zeta^t))_{t\in \R},$$ 
	(which we will sometimes write as $(h^t,\zeta^t)_{t\in \R} $ with a slight abuse of notation), which is \textbf{stationary} as a process of surface-curve pairs, and such that:
	\begin{itemize}
		\item $(\H,h^0,0,\infty)$ is a quantum wedge in the unit circle embedding;
		\item $(h^0,\zeta^0)$ has the law, as a surface curve pair, of a $(\gamma-2/\gamma)$-quantum wedge together with an independent SLE$_\kappa$;
		\item for any $t>0$, if  $\zeta^0$ is parametrised by half plane capacity, $$\sigma(t):=\inf\{s\ge 0: {\mathcal{V}}_{h^0}\left(\text{RHS of } \zeta^0([0,s])\right)\ge t \},$$  if $\tilde{g}_{\sigma(t)}$ is the centred Loewner map sending $\H\setminus \zeta^0{[0,\sigma(t)]}$ to $\H$, then we have that $h^{-t}=\tilde g_{\sigma(t)}(h^0)$ and $\zeta^{-t}=\tilde g_{\sigma(t)}(\zeta^0)$.
	\end{itemize}
Note that by stationarity, this defines the law of the process for all time (positive and negative).
\end{theorem}

%\begin{rmk}
%	Note that the boundary lengths with respect to $h^0$ of the left and right hand sides of $\eta^0$ are well defined.
%\end{rmk}

So this is a similar picture to that of the capacity zipper (moving backwards in time corresponds to ``cutting down'' and hence moving forward in time corresponds to ``zipping up'') but now a segment of $\zeta^0$ with \emph{right} $h^0$ LQG boundary length $t$ is cut out between times $0$ and $-t$. Hence the name ``quantum zipper'': the dynamic is parametrised by (right) quantum boundary length. Note that it makes sense to talk about the right boundary length of a segment of $\eta$, by conformally mapping to the upper half plane and applying the change of coordinate formula (see \cref{R:blmakessense}).  Also note the difference with the capacity zipper: here $h^0$ is a distribution (not a distribution modulo constants) while the stationarity is in the sense of quantum surface curve pairs.

Assuming for now that \cref{T:QZ} holds, we make the following claim.
\begin{claim} \label{claimLR} For any fixed $t$, the ${\mathcal{V}}_{h^0}$ boundary length of the \emph{left hand side} of $\zeta^0([0,\sigma(t)])$ is also equal to $t$.
	\end{claim}
This means that the parametrisation is really, unambiguously, by quantum boundary length. It also immediately implies \cref{T:quantumlength}.

\begin{proof}[Proof of  \cref{claimLR}, and hence \cref{T:quantumlength}, given
	\cref{T:QZ}]
	Denote by $L(t)$ the ${\mathcal{V}}_{h^0}$ boundary length of the left hand side of $\zeta^0[0,\sigma(t)]$, so our aim is to show that $L(t) \equiv t$. We begin by making the following observations.
	\begin{itemize}
		\item By stationarity of the quantum zipper, we have that $(L(s+t)-L(s))_{t\ge 0}$ is equal in distribution to $(L(t))_{t\ge 0}$ for any fixed $s\ge 0$.
		\item By scale invariance of SLE$_\kappa$ and the invariance property of quantum wedges (\cref{T:wedge}), 
		$$
		\frac{L(t)}{t} \; \overset{(d)}{=} L(1)
		$$
		for any $t>0$, and for any $A>0$ and $s<t$,
		$$
		\frac{L(At)}{At}-\frac{L(As)}{As}\; \overset{(d)}{=}\; \frac{L(t)}{t}-\frac{L(s)}{s}.
		$$
	\end{itemize}
	The first point means that we can apply Birkhoff's ergodic theorem \cite[Theorem 25.6]{kallenberg2}
\begin{equation}\label{eq:ergo}
\frac{L(n)}{n} \to X=\E(L(1)\,|\,\mathcal{I})  \text{ almost surely as } n\to \infty \text{ in } \N,
\end{equation} where $\mathcal{I}$ is the $\sigma$-field generated by invariant sets under the shift map $(L(1),L(2),\ldots)\mapsto $ $(L(2),L(3),\ldots)$.  Note that the theorem is often stated under the assumption that $\E(|L(1)|)$ is finite, but the conclusion is also true if we only know $L(1)\ge 0$ almost surely: in this case, conditional expectation can always be defined using monotone convergence, and the left hand side in \eqref{eq:ergo} converges to infinity on the event that the conditional expectation is infinite; see \cite[Theorem 25.6]{kallenberg2} for the proof.
	
	Note also that since $L(n)/n$ converges to $X$ in distribution, and since $L(n)/n$ is equal in distribution to $L(1)<\infty$ almost surely, we see that in fact $X<\infty$ almost surely.
	 We may then deduce that $$\frac{L(t)}{t}-\frac{L(s)}{s}=0 \text{ almost surely}$$ for any fixed $s,t\in \Q$ with $s\le t$. Indeed, the law $m$ of this difference is  equal to that of $L(At)/At-L(As)/As$ for any $A$, and by taking a sequence $A_k \uparrow \infty$ such that $A_k t\in \N, A_k s\in \N$ for all $k$, we obtain a sequence of random variables all having law $m$, which by \eqref{eq:ergo} tend to $0$ as $k\to \infty$.  Hence, with probability one we have that \begin{equation}
		\label{eq:ergo2}
		\frac{L(t)}{t}=X \;\; \forall t\in \Q\end{equation}
		(where $X$ is as in \eqref{eq:ergo}).  In particular, %$L(1)=X$ and $L(t)=tX$ for all $t\in \Q\cap (0,1)$, so that
		we have that
	  $$X=\lim_{t \downarrow 0, t\in \Q}\frac{L(t)}{t}.$$
	   Now by definition, the above limit (and therefore the random variable $X$) is measurable with respect to the $\sigma$-algebra
	  $$
	  \cT = \bigcap_{\eps>0} \sigma( \left.({h^0-h^0_\eps})\right|_{B(0,\eps)\cap \H}, \zeta^0|_{B(0,\eps)\cap \H})
	  $$
	  (here $h^0_\eps$ is the $\eps$-semicircle average of $h^0$ about the origin and can be subtracted since $L(t)/t$ is not affected by adding a constant to the field).
	  On the other hand, since the $h^0$ right/left quantum boundary lengths along $\zeta^0$ almost surely do not have atoms at $0^{\pm}$,
	  $X$ is also measurable with respect to
	  $$\sigma(\cA) \; ; \; \cA=\bigcup_{\eps>0} \sigma( \left.{(h^0-h^0_\eps)}\right|_{B(0,1)\setminus B(0,\eps)}, \zeta^0|_{B(0,1)\setminus B(0,\eps)}).$$
	  Hence the proof will be complete if we can show $\cT\cap \sigma(\cA)$ is trivial, because then $X$ must be almost surely constant, and by symmetry, this constant must be equal to $1$.
	
	  For this final step, since $\cA$ is a $\pi$-system, it suffices to show that for any  $\eps_0>0$, $A_0\in \cT$ and
	  $$A\in \sigma(h^0-h^0_{\eps_0}|_{B(0,1)\setminus B(0,\eps_0)}, \zeta^0|_{B(0,1)\setminus B(0,\eps_0)}),$$ we have  $\P(A\cap A_0)=\P(A)\P(A_0)$. However, this follows by independence of $h^0$ and $\zeta^0$, since the driving function of $\zeta^0$ is a Brownian motion, and by \cref{lem:hDNcircgood}.

\end{proof}

The rest of this section will be dedicated to proving \cref{T:CZ} and \cref{T:QZ}. In fact, \cref{T:CZ} is straightforward to obtain from \cref{T:coupling}. The idea to then deduce \cref{T:QZ} is to reparametrise time according to right quantum boundary length and appropriately ``zoom in'' at the whole capacity zipper picture at the origin. This step, however, is somewhat technical.

\subsubsection{The capacity zipper}

In this section we prove \cref{T:CZ}. That is, we construct the stationary two-sided capacity zipper, using the coupling theorem, \cref{T:coupling}.

Let $\bar h_0$ be as in the original  \cref{T:coupling} (that is, $\bar h_0$ has the distribution \eqref{h0}), and let $\eta = \eta_0$ be an independent infinite SLE$_\kappa$ curve from 0 to $\infty$. As in the coupling theorem, set $\bar h_t = \bar h_0 \circ \tilde g_{t}^{-1} + Q \log | ( \tilde g_t^{-1})'|$, where $\tilde{g}_t$ is the centred Loewner map corresponding to $\eta_0([0,t])$ for each $t$, and let $\eta_{t}$ be the image by $\tilde g_t$ of the initial infinite curve $\eta = \eta_0$. Then  \cref{T:coupling} says that $\bar{h}_t =\bar{h}_0$ in distribution, and in fact we can also see that the joint distribution $(\bar{h}_t, {\eta}_t)$ is identical to that of $(\bar{h}_0,\eta_0)$.

For $0 \le t \le T$, let $\bar{h}^t = \bar{h}_{T-t}$, and let $\eta^t = \eta_{T-t}$. Then it is an easy consequence of  \cref{T:coupling} that the following lemma holds:

\begin{lemma}
	\label{L:consistency}
	The laws of the process $(\bar{h}^t,\eta^t)_{0 \le t \le T}$ (with values in $\bar \cD'(\H) \times C([0, \infty)$) are consistent as $T$ increases.
\end{lemma}

By \cref{L:consistency}, and applying Kolmogorov's extension theorem, it is obvious that there is a well defined process $(\bar h^t, \eta^t)_{0 \le t  <\infty}$ whose restriction to $[0,T]$ agrees with the process described above. Hence for $t>0$, starting from $\bar h^0$ and an infinite curve $\eta^0$, there is a well defined, possibly random, procedure giving rise to $(\bar{h}^t,\eta^t)$, that we want to view as ``welding'' together parts of the positive and negative real lines, or ``zipping up''. The dynamic on the field is obtained by applying the change of coordinates formula to $\bar h^0$, with respect to a flow  $(f_s)_{s\le t}$ that has the \emph{marginal} law of a reverse Loewner flow, but we stress that here the reverse Loewner flow is not independent of $\bar h^0$ (rather, it will end up being uniquely determined by $\bar h^0$, while $(f_s)_{s\le t}$ will be independent of $\bar h^t$).

But we could also go in the other direction, cutting $\H$ along $\eta^0$, as in \cref{T:coupling}. Indeed we could define, for $t<0$ this time, a field $\bar h^t$ by considering the centred Loewner flow $(\tilde{g}_{|t|})_{t<0}$ associated to the infinite curve $\eta^0$, and setting
$$
\bar{h}^t = \bar h^0 \circ \tilde g^{-1}_{|t|}  + Q \log |(\tilde g_{|t|}^{-1})'| \quad (t<0).
$$
We can also, of course, get a new curve $\eta^t$ for $t<0$ by pushing $\eta^0$ through the map $\tilde g_{|t|}$. This gives rise to the two-sided stationary process $(\bar h^t, \eta^t)_{t \in \R}$ of \cref{T:CZ}.

\ind{Zipper! Capacity}
\begin{rmk}
	An equivalent way to define this process would be as follows. Start from the setup of  \cref{T:coupling}: thus $\bar{h}_0$ is a field distributed as in \eqref{h0}, and $\eta_0$ an independent infinite SLE$_\kappa$ curve. Set $\bar h_t = \bar h_0\circ \tilde g_t^{-1}  + Q \log |(\tilde g_t^{-1})'|$ as before, and $\eta_t = g_t(\eta^0 \setminus \eta_0[0,t])$. Then \cref{T:coupling} tells us that $(\bar h_t, \eta_t)_{t\ge 0}$ is a stationary process, so we can consider the limit as $t_0 \to \infty$ of $(\bar h_{t_0 + t}, \eta_{t_0+t})_{t \ge - t_0}$, which defines a two-sided process. The capacity zipper process $(\bar{h}^t,\eta^t)_{t \in \R}$ can then be defined as the image of this process under the time change $t\mapsto -t$.
\end{rmk}

\subsubsection{The quantum zipper}
\ind{Zipper! Quantum} 
\dis{We recall the notation \begin{equation}
		f(h):=h\circ f^{-1} + Q \log|(f^{-1})'|.
	\end{equation}
that will be used repeatedly in what follows.}

In this section we show the existence and stationarity of the quantum zipper: \cref{T:QZ}. In what follows, we will usually take our quantum wedges to be in the  \textbf{unit circle embedding} $(\H,h,0,\infty)$  (recall that the law of $h-\alpha \log(1/|z|)$ restricted to the upper unit semidisc is then just that of a Neumann GFF with additive constant fixed so that its average on the upper unit semicircle is equal to zero).

The key to the proof of \cref{T:QZ} is the following:
\begin{prop}\label{P:qwkey}
	Let $(h,\zeta)=((h,\H,0,\infty),\zeta)$ be a $(\gamma-2/\gamma)$-quantum wedge in the unit circle embedding, together with an independent SLE$_\kappa$. If $\zeta$ is parametrised by half plane capacity, let $\sigma$ be the smallest time such that the ${\mathcal{V}}_h$ boundary length of the right hand side of $\zeta([0,\sigma])]$ exceeds $1$\footnote{Recall that to measure this boundary length, we map the right hand side of the curve down to an interval $[0,x]$ of the positive real line using the centred Loewner map. Then we take the quantum boundary length of $[0,x]$ with respect to the field defined by applying the change of coordinates formula to $h$ with respect to this map.}. Let $g_\sigma$ be the centred Loewner map from $\H\setminus \zeta([0,\sigma])\to \H$. Then $(g_\sigma(h),g_\sigma(\zeta))$ is equal in law to $(h,\zeta)$ as a surface curve pair. That is, if $\psi$ is the unique conformal isomorphism such that $(\psi\circ g_\sigma)(h)$ is in the unit circle embedding, then $$(\psi \circ g_\sigma(h),\psi \circ g_\sigma(\zeta))\overset{(d)}{=} (h,\eta).$$
\end{prop}

\textbf{In words:} if we start with a $(\gamma-2/\gamma)$-quantum wedge and an independent SLE$_\kappa$, and ``zip'' down by one unit of right quantum boundary length, the law of the resulting quantum surface curve pair does not change.
\begin{proof}[Proof of \cref{T:QZ} given \cref{P:qwkey}]
	Note that there is nothing special about the choice to zip down by quantum boundary length one in \cref{P:qwkey}. Indeed we could replace one by any other $t>0$ and would obtain the result. Then the existence and stationarity of the quantum zipper follows in the same way that \cref{T:CZ} followed from \cref{T:coupling} (see the previous section). \end{proof}

The proof of \cref{P:qwkey} is quite tricky, and consists of several steps.

\mn \textbf{Step 1: Reweighting}	
We write $\P$ for the law of $(\tilde{h}^t,\eta^t)_{t\ge 0}$, the capacity zipper as in \cref{R:cznormed}, where the constant for $\tilde{h}^0$ has been fixed so that its unit semicircle average around the point $10$ is equal to 0 (this is fairly arbitrary, apart from the fact that the measure is supported a good distance away from the origin).
 We can extend this to define a law $\mathbf{P}$ on $\cD'(\H)\times C([0,\infty),\H)\times [1,2]$,  by setting
 $\mathbf{P}:=\P\times \mathrm{Leb}_{[1,2]}$ (so a sample from $\mathbf{P}$ consists of a capacity zipper $(\tilde{h}^t,\eta^t)_{t\ge 0}$ as just described, plus a point $Z$ chosen independently from Lebesgue measure on $[1,2]$). Define  $$c(z):=\E_{\mathbf{P}}(e^{\frac{\gamma}{2} \tilde h_\delta^0(z)}\delta^{\gamma^2/4}) \text{ for } z\in [1,2], $$
 which by \cref{T:BM} does not depend on $\delta>0$ and is a smooth function on $z\in[1,2]$.

We want to study the joint law of the capacity zipper plus a quantum boundary length typical point (in $[1,2]$). In fact, this is much easier to do if we reweight the law of the field $\tilde{h}^0$. To this end, we define a family of laws $(\mathbf{Q}_\eps)_{\eps>0}$ by setting
\begin{equation}
\label{eqn::defQe}
\frac{\dd \mathbf{Q}_\eps}{\dd \mathbf{P}}= \frac{\e^{\frac{\gamma}{2} \tilde h^0_\eps(Z)}\eps^{\frac{\gamma^2}{4}} }{\int_{[1,2]} c(z) \, dz}=: \frac{\e^{\frac{\gamma}{2} \tilde h^0_\eps(Z)}\eps^{\frac{\gamma^2}{4}} }{c([1,2])}
\end{equation}
for each $\eps$.

Under $\mathbf{Q}_\eps$, the marginal law of $(\tilde{h}^t,\eta^t)_{t\ge 0}$ is its $\P$ law weighted by $$\frac{ {\mathcal{V}}_{\tilde h^0_\eps}([1,2])}{c([1,2])}.$$ Moreover, given $(\tilde h^t,\eta^t)_{t\ge 0}$ the point $Z$ is sampled from the $\eps$-approximate measure ${\mathcal{V}}_{\tilde h^0_\eps}$ (restricted to $[1,2]$ and normalised to be a probability measure). Therefore, since ${\mathcal{V}}_{\tilde h^0_\eps}([1,2]) \to {\mathcal{V}}_{\tilde h^0}([1,2])$ in $\mathcal{L}^1$ as $\eps\to 0$ and the measure ${\mathcal{V}}_{\tilde h^0_\eps}$ converges weakly in probability to ${\mathcal{V}}_{\tilde h^0}$, we can deduce that
\[ \mathbf{Q}_\eps \Rightarrow \mathbf{Q}\] as $\eps \to 0$
where $\mathbf{Q}$ is the measure described by (a) and (b) of \cref{lem::q_alt} below.

This reweighting is analogous to the argument used to describe the GFF viewed from a Liouville typical point -- see \cref{T:rooted}. As in this proof, we can reverse the order in which $(\tilde{h}^0,\eta^0)$ and $Z$ are sampled, and this leads to the alternative description given by points (c) to (e) in following lemma.
\begin{lemma}\label{lem::q_alt}
	Under $\mathbf{Q}$, the following is true:
	\begin{enumerate}
		\item[(a)] the marginal law of $(\tilde h^t,\eta^t)_{t\in \R}$ is given by ${\mathcal{V}}_{\tilde h^0}([1,2])/c([1,2]) \dd \P$ (and is therefore absolutely continuous with respect to $\P$);
		\item[(b)] conditionally on $(\tilde{h}^t,\eta^t)_{t\in \R}$, $Z$ is chosen uniformly from ${\mathcal{V}}_{\tilde h^0}$ on $[1,2]$;
		\item[(c)] the marginal law of $Z$ on $[1,2]$ has density $c(z)/c([1,2])$ with respect to Lebesgue measure;
		\item[(d)] conditionally on $Z$, for every $0\le t\le\tau_Z$ (where $\tau_Z$ is the first time that $f_t(Z)=0$) the law of $\eta^t([0,t])$ is that of a reverse SLE$_\kappa(\kappa,-\kappa)$ curve with force points $(Z,10)$, run up to time $t$;
		\item[(e)] for any $t\ge 0$, conditionally on $\{Z,(f_s)_{0 \le s\le t}\}$, we have that the conditional law of $\tilde{h}^t$ as a distribution  modulo constants, \emph{on the event $t\le \tau_Z$}, is that of \[\bar h+\frac{2}{\gamma}\log(|\cdot|) +\frac{\gamma}{2}G^\H_N(\cdot, f_t(Z))-\frac{\gamma}{2}\int G^\H(\cdot, f_t(y)) \rho_{10,1}(dy)\] where $\bar h$ has the law of a Neumann GFF (modulo constants) that is independent of $(f_s)_{0\le s\le t}$. Here for $x\in \R, \delta>0$, $\rho_{x,\delta}$ denotes uniform measure on the upper semicircle of radius $\delta$ around $x$.
	\end{enumerate}
\end{lemma}

\begin{rmk}
	The force point at 10 in (d) and the final term in the expression for $\tilde{h}^t$ in (e) make these descriptions look rather complicated. However, we will really be interested in taking $t=\tau_Z$ and looking at $(\tilde{h}^t,\eta^t$) in small neighbourhoods of the origin. In such a setting, as we will soon see, these terms will have asymptotically negligible contribution to the behaviour. The only features in the descriptions (d) and (e) that are genuinely important, are the force point of weight $\kappa$ at $Z$, and the function $(2/\gamma)\log(|\cdot|)+(\gamma/2)G(\cdot, f_t(Z))$.
\end{rmk}
%Note from Ellen 07/19: I wasn't sure about what was written previously - that this reweighting does not depend on the choice of normalisation for the field. I think that is true when you consider the reweighted marginal law of the zipper, but not for the law of $Z$ or for the conditional law of the zipper given $Z$. This should be checked though!)

\begin{proof} (a) and (b) define the measure $\mathbf{Q}$ (see discussion above the lemma) and (c) follows since this is true under $\mathbf{Q}_\eps$ for every $\eps>0$.
	
	For (d), we first claim that for any $t\ge 0$ and for any measurable function $F$ of $(f_s; s\le t)$ we have
	\begin{equation}\label{eq:ford} \mathbf{Q}_\eps(F(f_s;s\le t)\mathbf{1}_{\{t\le \tau_{Z-\eps}\}})=\frac{\mathbf{P}(F(f_s;s\le t)\mathbf{1}_{\{t\le \tau_{Z-\eps}\}}\e^{\frac{\gamma}{2}(M_{t}(Z)- M_{t}(10))- \frac{\gamma^2}{8}[M(Z)-M(10)]_{t}})}{\mathbf{P}(\e^{\frac{\gamma}{2}(M_{t}(Z)- M_{t}(10))- \frac{\gamma^2}{8}[M(Z)-M(10)]_{t}})}.\end{equation}
	To see this, we note that by definition $\tilde h^0=\tilde h^{t}\circ f_{t}+Q\log|(f_{t})'|$ and, due to the normalisation we chose for $\tilde h^0$, $\tilde h^0_\eps = (\tilde h^0,\bar{\rho}_Z^\eps):=(\tilde h^0,\rho_{Z,\eps}-\rho_{10,1})$. Therefore
	\[ \mathbf{P}(e^{\frac{\gamma}{2}\tilde{h}^\eps_0(Z)}\,|\, Z, (f_s; s\le t))=e^{\frac{\gamma}{2}(M_t,\rho_{Z,\eps}-\rho_{10,1})}\mathbf{P}(e^{\frac{\gamma}{2}(\tilde{h}^t\circ f_t-(2/\gamma)\log(|\cdot|)\circ f_t,\rho_{Z,\eps}-\rho_{10,1})}\, |\,Z, (f_s;s\leq t)) \]
where, because the average value of $\rho_{Z,\eps}-\rho_{10,1}$ is equal to 0, $(\tilde{h}^t\circ f_t-(2/\gamma)\log(|\cdot|)\circ f_t,\rho_{Z,\eps}-\rho_{10,1})$ depends only on the equivalence class modulo constants of $\tilde{h}^t\circ f_t-(2/\gamma)\log(|\cdot|)$. Moreover, by stationarity of the capacity zipper, this law is that of a Neumann GFF $\bar{h}$ (modulo constants) that is independent of $(f_s;s\le t)$. Thus we are reduced to doing a simple Gaussian computation. This is very similar what was carried out in the proof of \cref{T:coupling} and yields that
	\[\mathbf{P}(e^{\frac{\gamma}{2}(\tilde{h}^t\circ f_t-(2/\gamma)\log(|\cdot|)\circ f_t,\rho_{Z,\eps}-\rho_{10,1})}\, |\,Z, (f_s;s\leq t))=e^{-\frac{\gamma^2}{8}[(M,\rho_{Z,\eps}-\rho_{10,1})]_t}.\] We may also note that when $t\le \tau_{Z-\eps}$, $M_t$ can be extended by Schwarz reflection to a harmonic function on a domain containing $B(Z,\eps)$ and $B(10,1)$, and so by the mean value theorem $(M_{t},\bar{\rho}_{Z,\eps}-\rho_{10,1})=M_t(Z)-M_t(10)$. Similarly, on the event that $t\le \tau_{Z-\eps}$,  $[(M,\rho_{Z,\eps}-\rho_{10,1})]_t=[M(Z)-M(10)]_t$. \eqref{eq:ford} then follows by definition of $\mathbf{Q}_\eps$ and conditioning.

	Next, recall from the proof of \cref{L:mart} that $\dd M^*_r=-(2/(\gamma f_r(z)))\dd W_r$ where $W$ is the driving function of $(f_r)_r$ (and is a Brownian motion run at speed $\gamma^2$). Hence, by \eqref{eq:ford} and the Cameron--Martin--Girsanov theorem we have that (under $\mathbf{Q}_\eps$, conditionally on $Z$ and up to time $\tau_{Z-\eps}$), $W_t-\frac{\gamma}{2}[ W, M (Z)- M(10)]_t$ is a (speed $\gamma^2$) Brownian motion, or equivalently
	\begin{equation*}
	\dd W_t=\gamma \dd B_t-\gamma^2 \Re(\frac{1}{f_t(Z)}) \dd t+\gamma^2 \Re(\frac{1}{f_t(10)}) \dd t.
	\end{equation*}
	Since this does not depend on $\eps$, the same must hold under $\mathbf{Q}^Z=\mathbf{Q}(\cdot | Z)$, at least up to time $\tau_{Z-\eps}$. However, as $\eps>0$ was arbitrary, it in fact holds until time $\tau_Z$. Since this is exactly the equation satisfied by the driving function of an SLE$_{\kappa}(\kappa,-\kappa)$ process with force points at $(Z,10)$, we conclude the proof of (d).

	Finally, we deal with (e). For this, we use the same rewriting of $\tilde{h}_0^\eps$ as above to see that
	\begin{equation}\label{eq:condRN}\mathbf{Q}^\eps(F(\tilde{h}^t)\, |\, (f_s)_{0\le s\le t},Z) =\frac{\mathbf{P}(F(\tilde{h}^t) e^{\frac{\gamma}{2}(\tilde h^t, (\rho_{Z,\eps}-\rho_{10,1})\circ f_t^{-1})}\,|\,(f_s)_{0\le s\le t},Z)}{\mathbf{P}( e^{\frac{\gamma}{2}(\tilde h^t, (\rho_{Z,\eps}-\rho_{10,1})\circ f_t^{-1})}\,|\,(f_s)_{0\le s\le t},Z)}\end{equation}
\begin{comment}$$\mathbf{Q}_\eps(F(\tilde{h}^t))=\frac{\eps^{\gamma^2/4}}{c([1,2])}\mathbf{P}(F(\tilde{h}^t)e^{\frac{\gamma}{2}(\tilde h^t\circ f_t + Q\log|f_t'|, \rho_{Z,\eps}-\rho_{10,1})})=\frac{\mathbf{P}(F(\tilde{h}^t)e^{\frac{\gamma}{2}(\tilde h^t, (\rho_{Z,\eps}-\rho_{10,1})\circ f_t^{-1})})}{\mathbf{P}(e^{\frac{\gamma}{2}(\tilde h^t, (\rho_{Z,\eps}-\rho_{10,1})\circ f_t^{-1})})}$$
\end{comment}
	for any bounded measurable function $F$ of $\tilde{h}^t$ modulo constants. %\footnote{If $d\mathbb{Q}/d\mathbb{P}=Z$, then for $Y$ a bounded random variable and $\cG$ a sub-sigma algebra $\mathbb{Q}(Y|\cG)=\frac{\P(YZ|\mathcal{G})}{\P(Z|\cG)}$ because for any event $A\in \cG$, $\mathbb{Q}(Y1_A)=\P(YZ1_A)=\P(Z\frac{\P(YZ|\cG)}{P(Z|\cG)}1_A)=\mathbb{Q}(\frac{\P(YZ|\mathcal{G})}{\P(Z|\cG)}1_A)$.}
	On the other hand, recall that under $\mathbf{P}$, $\tilde h^{t}$ viewed modulo constants is independent of $(f_s)_{s\le t}$, and is distributed like a Neumann GFF plus the function $(2/\gamma) \log |\cdot|$ (modulo constants).
	\begin{comment}*********we use again the fact that $\tilde h^0=\tilde h^{t}\cdot f_{t}+Q\log|f'_{t}|$ to write
	\[\frac{d\mathbf{Q}_\eps^Z}{d\P}= \frac{\e^{\frac{\gamma}{2}(\tilde h^t\circ f_t + Q\log|f_t'|, \bar{\rho}_Z^\eps)}\eps^{\frac{\gamma^2}{4}}}{c(Z)}. \]
	Since $\tilde h^{t}-(2/\gamma) \log |\cdot|$ viewed modulo constants is distributed like a Neumann GFF modulo constants, and independent of $(f_s)_{s\le t}$ under $\P$, we can compute that
	\[\left. \frac{d\mathbf{Q}_\eps^Z(\cdot \, | \, (f_s)_{s\le t})}{d\P} \right|_{\sigma(\tilde h^t)}= \frac{\e^{\frac{\gamma}{2}(\tilde h^t,\, \bar{\rho}_Z^\eps\circ f_t^{-1})}}{\E_\P(\e^{\frac{\gamma}{2}(\tilde h^t,\, \bar{\rho}_Z^\eps\circ f_t^{-1})})}.\]
**********
\end{comment}
		Thus, by the Cameron--Martin--Girsanov theorem applied conditionally on $(Z,(f_s)_{s\le t})$, the law of $\tilde h^t$ under $\mathbf{Q}_\eps$ and conditionally on $(Z,(f_s)_{s\le t})$, considered modulo constants, is that of a Neumann GFF (modulo constants) \emph{plus} the function $(2/\gamma) \log|\cdot|$, \emph{plus} the function 
		$w\mapsto \int G^\H_N(w, y) (\bar \rho_Z^\eps\circ f_{t}^{-1})(\dd y)$. %``drift distribution'' that sends smooth $g$ to $\E_{\P}((\tilde h^t,g)(\tilde h^t,\bar{\rho}_Z^\eps\circ f_t^{-1}))$.
	Now, for any $t\le \tau_{Z-\eps}$ and any $w\in \H\setminus f_t(B(Z,\eps)\cap \H)$ we have
	$\int G^\H_N(w, y) (\bar \rho_Z^\eps\circ f_{t}^{-1})(\dd y) = G^\H_N(w, f_{t}(Z))$, and so on the set $\H\setminus f_t(B(Z,\eps)\cap \H)$ we can write (as distributions modulo constants)
	\begin{equation} \tilde h^{t} \overset{(d)}{=} \bar{h} + \frac{2}{\gamma} \log|\cdot|+\frac{\gamma}{2} G^\H_N(\cdot, f_{t}(Z)) - \frac{\gamma}{2}\int G^\H_N(\cdot, f_t(y)) \rho_{10}^1 (\dd y),\end{equation}
	where the equality in distribution holds under $\mathbf{Q}_\eps$ conditionally on $Z$ and $(f_s; s\le t)$, and where $\bar{h}$ is as described in the statement of (e). Taking a limit as $\eps\to 0$ we obtain the result.
\end{proof}

\begin{cor}\label{rmk::htau} Taking $t\nearrow \tau_Z$ in the previous lemma, we see that under $\mathbf{Q}^Z=\mathbf{Q}(\cdot\, | \, Z)$, $(\tilde{h}^{\tau_Z},{\eta}^{\tau_Z})$ can be described as follows:
	\begin{itemize}
		\item $\eta^{\tau_Z}([0,\tau_Z])$ has the law of a reverse SLE$_{\kappa}(\kappa,-\kappa)$, with force points at $(Z,10)$, and run until the point $Z$ reaches $0$;
		\item given $\eta^{\tau_Z}([0,\tau_Z])$, one has $\tilde h^{\tau_Z}\overset{(d)}{=} \bar h+(\gamma-2/\gamma)\log(1/|\cdot|) -\frac{\gamma}{2}\int G^\H_N(\cdot, f_{\tau_Z}(y)) \rho_{10,1}(\dd y)$, where the equality is an equality of distributions modulo constants, and $\bar h$ has the law of a Neumann GFF that is independent of $(f_s)_{0\le s\le \tau_Z}$.		
\end{itemize} \end{cor}

\begin{rmk}
Taking $t \nearrow \tau_Z$ rigorously in \cref{lem::q_alt} requires some justification, since the statement of (e) is actually for deterministic $t\ge 0$. 

To do this, we first consider $\tau:=\tau_{\{Z-\delta\}}$ for arbitrary $\delta>0$. Then from \cref{lem::q_alt}(e) we have that for any deterministic $k,n\ge 0$, the conditional law of $\tilde{h}^{k/n}$ given $Z$ and $(f_s)_{0\le s \le k/n}$, \emph{on the event} that $\tau\in (\tfrac{k-1}{n}, \tfrac{k}{n}]$, is that of
\begin{equation}\label{eq:trandomrigorous}\bar h +F_{f_{k/n},Z}(\cdot); \quad  F_{f_{k/n},Z}(\cdot):=\frac{2}{\gamma}\log(|\cdot|)+\frac{\gamma}{2}G^\H_N(\cdot, f_t(Z))-\frac{\gamma}{2}\int G^\H_N(\cdot, f_t(y)) \rho_{10,1}(dy).\end{equation}
	Now write $\tau^{(n)}:=k/n$ for the unique $k\in \mathbb{N}$ such that $\tau\in (\tfrac{k-1}{n}, \tfrac{k}{n}]$. Then for arbitrary continuous functionals $H_1, H_2, H_3$ (defined on appropriate spaces) taking values in $[0,1]$ we have 
\begin{align*}
	\mathbf{Q}[H_1(\tilde{h}^\tau)H_2((f_s)_{0\le s \le \tau}) H_3(Z)] & = \lim_{n\to \infty} \mathbf{Q}[H_1(\tilde{h}^{\tau^{(n)}})H_2((f_s)_{0\le s \le \tau^{(n)}}) H_3(Z)] \\
	& = \lim_{n\to \infty} \sum_{k} \mathbf{Q}[H_1(\tilde{h}^{k/n})H_2((f_s)_{0\le s \le k/n}) H_3(Z)\mathbf{1}_{\tau^{(n)}=k/n}] \\
	& = \lim_{n\to \infty} \sum_k \mathbf{Q}[H_1(\bar h +F_{f_{k/n},Z}(\cdot))H_2((f_s)_{0\le s \le k/n}) H_3(Z)\mathbf{1}_{\tau^{(n)}=k/n}]  \\
	& = \lim_{n\to \infty} \mathbf{Q}[H_1(\bar h +F_{f_{\tau^{(n)}},Z}(\cdot))H_2((f_s)_{0\le s \le \tau^{(n)}}) H_3(Z)] \\
	& = \mathbf{Q}[H_1(\bar h +F_{f_{\tau},Z}(\cdot))H_2((f_s)_{0\le s \le \tau}) H_3(Z)]
\end{align*}
where the middle equality follows from \eqref{eq:trandomrigorous}, the first and final by continuity, and the second and fourth by definition of $\tau^{(n)}$.
In other words, \cref{lem::q_alt}(e) holds for the random time $t=\tau=\tau_{\{Z-\delta\}}$ for any $\delta>0$. We can similarly take $\delta\to 0$ to obtain the statement for $t=\tau_Z$. 
\end{rmk}
We will use this to show that when we zoom in at this weighted capacity zipper at time $\tau_Z$, we obtain a field and curve whose joint law is that in the statement of \cref{P:qwkey}.

\mn \textbf{Step 2: Zooming in to get a wedge and an independent SLE}  Suppose that $\eta$ is a simple curve from $0$ to $\infty$ in $\H$, considered up to time reparametrisation, and that $K\subset \H$ is compact. In what follows, by $\eta$ restricted to $K$, we mean the trace of $\eta$ run up to the first time that it exits the set $K$ (which does not depend on the choice of time parametrisation). If $h\in \cD_0'(\H)$, by $h$ restricted to $K$, we mean the restriction in the standard sense of restriction of distributions.
\begin{lemma}\label{lem::zip_to_wedge}
	Let $((\tilde{h}^t,\eta^t)_{0\le t\le \tau_Z},Z)$ be sampled from $\mathbf{Q}$. Let $\varphi_C$ be the unique conformal isomorphism $\H\to \H$ such that $(\H,
	\varphi_C(\tilde h^{\tau_Z}+C),0,\infty)$
	is the unit circle embedding of $(\H,\tilde h^{\tau_Z}+C,0,\infty)$. Then for any $K\subset \H$ compact, the law of $(\varphi_C(\tilde h^{\tau_Z}+C),\varphi_C(\eta^{\tau_Z}))$ restricted to $K$ converges in total variation distance to the law of $(h,\zeta)$ restricted to $K$, where $(h,\zeta)$ is as in \cref{P:qwkey}.
\end{lemma}

\noindent Note that $\{(\varphi_C(\tilde{h}^{\tau_Z}+C),\varphi_C(\eta^{\tau_Z})): C >0\}$  is completely determined by $(\tilde{h}^{\tau_Z},\eta^{\tau_Z},Z)$. 
\vspace{1em}

For the proof of \cref{lem::zip_to_wedge}, let us first explain why zooming in (i.e., applying $\varphi_C$ to $\eta^{\tau_Z}$) produces an SLE$_\kappa$ curve $\zeta$ as $C\to \infty$. 

\begin{lemma}\label{lem:zoomSLEkkk}
	Let $K\subset\H$ be compact, let $Z$ be sampled from its $\mathbf{Q}$ law, and let $(\eta^t)_{0\le t\le \tau_Z}$ be a reverse SLE$_\kappa(\kappa,-\kappa)$ curve with force points at $(Z,10)$, run until $f_t(Z)$ reaches $0$. Let $\phi_R(z)=Rz$ for $R>0$. Then $\phi_R(\eta^{\tau_Z})|_K$ converges in total variation distance to $\zeta|_K$ as $R\to \infty$, where $\zeta$ has the law of an SLE$_\kappa$ curve from $0$ to $\infty$ in $\H$. \end{lemma}

\begin{proof}
	It is equivalent to prove that $\eta^{\tau_Z}|_{\delta\mathbb{D}_+}$ converges in total variation distance to $\zeta|_{\delta\mathbb{D}_+}$ as $\delta\downarrow 0$, where we recall that $\mathbb{D}_+:=B(0,\delta)\cap \H$ for $\delta>0$. 
	
	To show this, since reverse Loewner evolutions grow ``from the base''  rather than the tip, we only need to concentrate on the evolution of the curve on a time interval just preceeding $\tau_Z$. So, if $(f_t)_t$ is the reverse flow associated with $(\eta^t)_t$, we first pick $a>0$ small and run $f_t$ until the first time $\tau_Z^a$ that $f_{\tau_Z^a}(Z)=a$. Note that by \cref{lem::reversesleprops}, $f_T(10)>1$. Moreover by \cref{lem::reversesleprops}, the total variation distance between the flow $(f_{\tau_Z^a+t})_{0<t \le \tau_Z-\tau_Z^a}$ and the flow $(\tilde{f}_t)_{0\le t \le \tau_a}$ of a reverse SLE$_\kappa(\kappa)$ $\tilde{\eta}$ with a single force point at $a$ run until $\tilde{f}_t(a)=0$, converges to $0$ as $a\to 0$. (To apply \cref{lem::reversesleprops} directly we need to rescale by $1/a$, sending $a\mapsto 1$ and $f_{\tau_Z^a}(10)\mapsto f_{\tau_Z^a}(10)/a>1/a$).

But for $(\tilde{f}_t)$ we can use the time reversal symmetry of SLE$_\kappa(\rho)$ -- \cref{C:slerevsym} -- which says that the curve generated the reverse Loewner flow $(\tilde{f}_t)_{0\le t \le \tau_a}$ is that of an ordinary forward SLE$_\kappa$ run until an almost surely positive time $\Lambda_a$. Roughly speaking, this concludes the proof since zooming in at the origin for this segment of forward SLE$_\kappa$ produces an infinite SLE$_\kappa$ as desired.

Let us now put all these pieces together rigorously. Given $\eps>0$, we first choose $a>0$ small enough that $(f_t)_{t\le \tau_Z}$ can be coupled with the hybrid flow -  that is, $(f_t)_{t\le \tau_Z^a}$ concatenated with $(\tilde{f}_t)_{t\le \tau_a}$ - so that they agree with probability $\ge 1-\eps$. Then for $\delta>0$ small enough, we also have a sub-event of probability $\ge 1-2\eps$ that $\tilde{f}_{\tau_a}^{-1}(B(0,\delta))$ does not intersect $\eta^{\tau_Z^a}$. In other words, on this subevent we have that  $\eta^{\tau_Z}|_{\delta \D_+}=\tilde{\eta}^{\tau_a}|_{\delta \D_+}$, where $\tilde{\eta}^{\tau_a}$ has the law of a reverse SLE$_\kappa(\kappa)$ with a force point at $a$, run until the reverse Loewner image of $a$ reaches the origin. Finally, we use that for fixed $a$, the total variation distance between $\tilde{\eta}^{\tau_a}|_{\delta \D_+}$ and $\zeta|_{\delta \D_+}$ converges to $0$ as $\delta\to 0$. So choosing $\delta$ small enough we can couple $\eta^{\tau_Z}|_{\delta \D_+}$ and $\zeta|_{\delta \D_+}$ so that they agree with probability $\ge 1-3\eps$. Since $\eps>0$ was arbitrary, this proves the desired convergence in total variation distance.
\end{proof}

Since we know that under $\mathbf{Q}$, $\tilde{h}^{\tau_Z}$ is essentially a Neumann GFF plus an appropriate log-singularity, independent of $\eta^{\tau_Z}$ we want to conclude the proof of \cref{lem::zip_to_wedge} using Theorem \ref{T:wedge}, which says that zooming in at the origin for such a field produces a quantum wedge. However, the law of $\tilde{h}^{\tau_Z}$ also includes a random additive constant and a harmonic function (see Corollary \ref{rmk::htau}), both of which in fact depend on the curve $\eta^{\tau_Z}$ itself. This dependence could be problematic when trying to describe what happens to the curve $\eta^{\tau_Z}$ when we zoom in. We thus appeal to a more refined version of Theorem \ref{T:wedge}, namely Lemma \ref{lem:Twedgenorm}.

\begin{proof}[Proof of \cref{lem::zip_to_wedge}] 
Corollary \ref{rmk::htau} and Lemma \ref{lem:Twedgenorm} imply that almost surely, conditionally on $(f_t)_{t\le \tau_Z}$, the conditional law of $(\varphi_C,\varphi_C(\tilde{h}^{\tau_Z}+C)|_K)$ becomes arbitrarily close, in total variation distance as $C\to \infty$, to the law of $(\psi_C,\psi_C(h'+C)|_K)$ where $h'$ is a Neumann GFF normalised to have average $0$ on $\partial \D_+$ plus a $(\gamma-2/\gamma)$ log singularity at the origin, and $\psi_C$ is the scaling map such that $\psi_C(h'+C)$ is the unit circle embedding of $h'+C$. It is not hard to see that this implies 
$$
\mathrm{d}_{\tv}((\eta^{\tau_Z},\varphi_C,\varphi_C(\tilde{h}^{\tau_Z}+C)|_K) ; (\eta^{\tau_Z},\psi_C,\psi_C(h'+C)|_K)\to 0
$$
as $C\to \infty$, where in the second triple above, $(\psi_C,\psi_C(h'+C))$ is independent of $\eta^{\tau_Z}$.  In turn, this implies that 
$$
\mathrm{d}_{\tv}((\varphi_C(\eta^{\tau_Z})|_K,\varphi_C(\tilde{h}^{\tau_Z}+C)|_K) ; (\psi_C(\eta^{\tau_Z})|_K,\psi_C(h'+C)|_K))\to 0
$$
as $C\to \infty$. Now by Theorem \ref{T:wedge}, as $C\to \infty$, the law of $\psi_C(h'+C)|_K$ converges (in total variation) to the law of the quantum wedge $h$ appearing in the statement of Lemma \ref{lem::zip_to_wedge}. Furthermore, since $\psi_C$ and $h'$ are independent of $\eta^{\tau_Z}$, Lemma \ref{lem:zoomSLEkkk} implies that 
$$
\mathrm{d}_{\tv}((\psi_C(\eta^{\tau_Z})|_K,\psi_C(h'+C)|_K);(\zeta|_K,h|_K))\to 0
$$
as $C\to \infty$, where $\zeta$ is an SLE$_\kappa$ independent of $h$. Applying the triangle inequality concludes the proof.
\end{proof}

\mn \textbf{Step 3: Stationarity} In Step 2 above, we have shown that if one zooms in at the capacity zipper with reweighted law $\mathbf{Q}$ at time $\tau_Z$, then one obtains a field/curve pair having the distribution of $(h,\zeta)$ as in \cref{P:qwkey}. In this step we will prove that the operation of ``zipping down right quantum boundary length one'' does not change this law, and hence prove \cref{P:qwkey}.

Given a sample $((\tilde{h}^t,\eta^t)_{0\le t\le \tau_Z},Z)$ from $\mathbf{Q}$, and $C>0$, let $Z_C\in [0,Z]$ be such that ${\mathcal{V}}_{\tilde{h}^0}([Z_C,Z])=\e^{-C\gamma/2}$. If this is not possible (ie. if ${\mathcal{V}}_{\tilde{h}^0}([0,Z])<\e^{-C\gamma/2}$), set $Z_C=0$. Set $\tau_C=\tau_{Z_C}$ and let $\phi_C$ be the unique conformal isomorphism such that $(\H,\phi_C(\tilde h^{\tau_C}+C),0,\infty)$ is the unit circle embedding of $(\H, \tilde{h}^{\tau_C}+C,0,\infty)$.

Recall the notation $g_\sigma, \psi$ from \cref{P:qwkey}.
\begin{lemma}\label{L:s1} For any $K\subset \H$ compact,
	$(\phi_C(\tilde{h}^{\tau_C}+C),\phi_C(\eta^{\tau_C}))$ restricted to $K$ converges in total variation distance to $(\psi\circ g_\sigma(h),\psi \circ g_\sigma(\zeta))$ restricted to $K$, as $C\to \infty$.
\end{lemma}
\begin{lemma}\label{L:s2}For any $K\subset \H$ compact,
		$(\phi_C(\tilde{h}^{\tau_C}+C),\phi_C(\eta^{\tau_C}))$ restricted to $K$ converges in total variation distance to $(h,\zeta)$ restricted to $K$, as $C\to \infty$.
\end{lemma}

\begin{proof}[Proof of \cref{P:qwkey}]  \cref{L:s1,L:s2} tell us that for any $K$ we can couple $(h,\zeta)$ and $(\psi\circ g_\sigma(h),\psi\circ g_\sigma(\zeta))$ together so that they agree when restricted to $K$ with as high probability as we like. Thus their laws, when restricted to $K$, must agree. Since $K$ was arbitrary, we can conclude.
\end{proof}

\begin{figure}
	\begin{center}
	\includegraphics[width=0.8\textwidth]{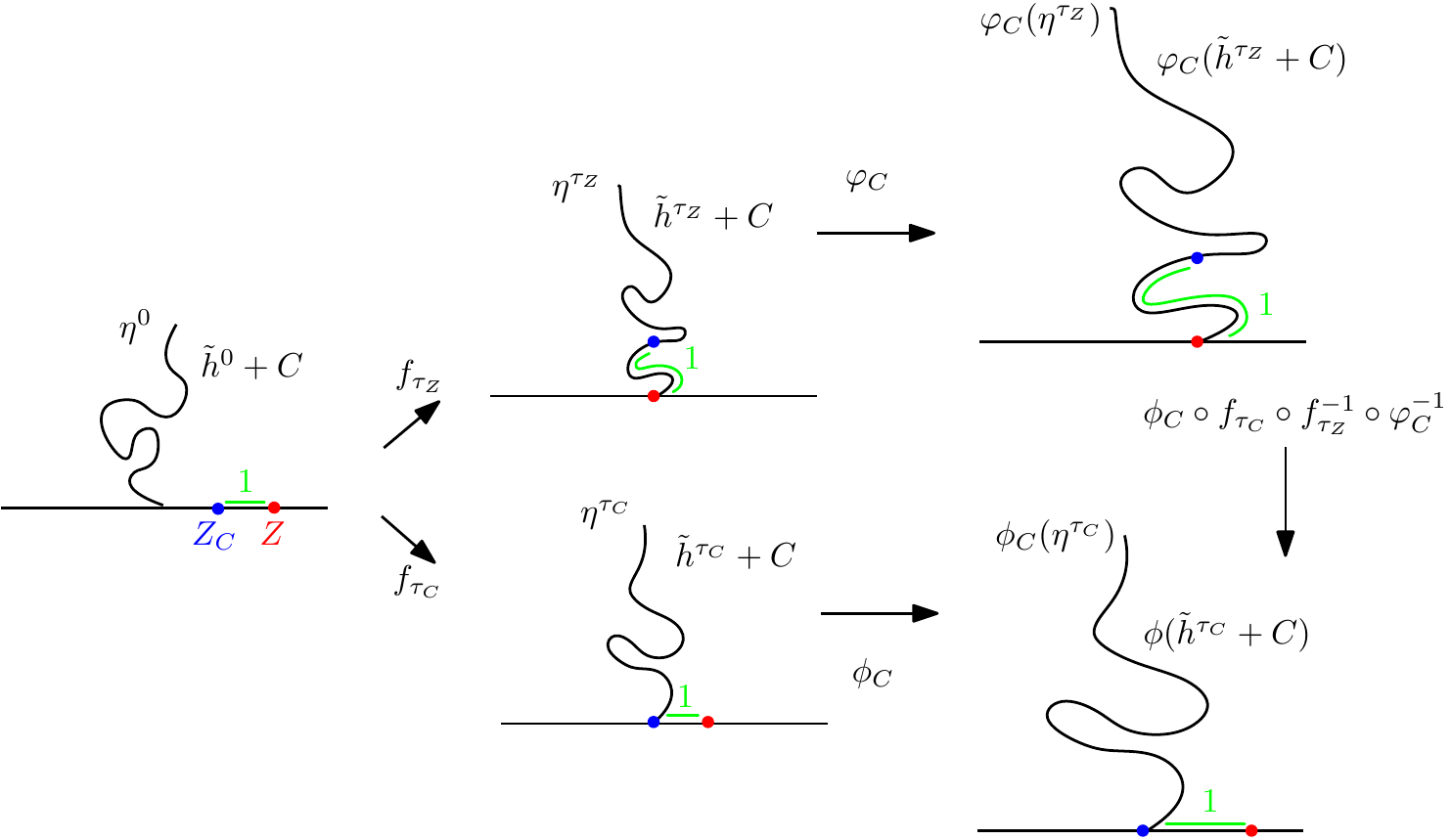}
	\end{center}
\caption{All the marked quantum boundary lengths (with respect to the field indicated on the relevant diagram) are equal to one. This is by definition of the conformal isomorphisms $f_{\tau_Z},f_{\tau_C},\phi_C$ and $\varphi_C$. Recall that $f(h)$ is obtained from $h$ by applying the conformal change of coordinates formula which preserves quantum boundary length.}
\end{figure}
\begin{proof}[Proof of \cref{L:s1}]
Let $\eps>0$ be arbitrary. First observe that we can choose $K_\eps\subset \H$ compact so that $K\subset g_\sigma(K_\eps)$ with probability greater than $1-\eps$.  By \cref{lem::zip_to_wedge}, for large enough $C$ we can also couple $(\varphi_C(\tilde h^{\tau_Z}+C),\varphi_C(\eta^{\tau_Z}))$ and $(h,\zeta)$ such that with probability $>1-\eps$ they are equal in $K_\eps$. We may also (by taking $C$ large enough) require that $Z_{C}\ne 0$ on this event. Then, since on this event we have that $$(\phi_C(\tilde h^{\tau_C}+C),\phi_C(\eta^{\tau_C}))=(\psi\circ g_\sigma(h),\psi\circ g_\sigma (\zeta))$$ (this is clear since these pairs are obtained from $(\varphi_C(\tilde h^{\tau_Z}+C),\varphi_C(\eta^{\tau_Z}))$ and $(h,\zeta)$ respectively by zipping down $1$ unit of right quantum boundary length and applying a conformal isomorphism so as to be in the unit circle parametrisation) the result follows.
\end{proof}

\begin{proof}[Proof of \cref{L:s2}]
For this, observe that if $\mu$ is the law of a uniform point in $[0,A]$ for $A>0$, and $\nu$ is the law of $U-\eps$ for $U\sim\mu$, then the total variation distance between $\nu$ and $\mu$ tends to $0$ as $\eps \to 0$. This means that we can couple the $\mathbf{Q}$ laws of $(Z,(\tilde h^t,\eta^t)_{t\ge 0})$ and $(Z_C,(\tilde h^t,{\eta}^t)_{t\ge 0})$ such that they are equal with probability tending to $1$ as $C\to \infty$ (by \cref{lem::q_alt} (b), definition of $Z_C$ and the fact that the $\tilde h^0$ boundary length of $[1,2]$ is finite almost surely). Hence we can couple the $\mathbf{Q}$ laws of $(\varphi_C(\tilde{h}^{\tau_Z}+C),\varphi_C({\eta}^{\tau_Z}))$ and $(\phi_C(\tilde{h}^{\tau_C}),\phi_C({\eta}^{\tau_C}))$ so they are equal with probability tending to $1$ as $C\to \infty$. Since the former law converges to that of $(h,\zeta)$ as $C\to \infty$ (\cref{lem::zip_to_wedge}), the same therefore holds for the latter.
\end{proof}

\subsection{Uniqueness of the welding}

\ind{Conformal welding}
Consider the \textbf{capacity zipper} $(\tilde h^t,\eta^t)_{t\in \R}$ of  \cref{R:cznormed} (where the additive constant for $\tilde{h}^0$ is fixed). The (reverse) Loewner flow associated to $ (\eta^t)_{t\ge 0}$ has the property that it zips  together intervals of $\R_+$ and $\R_-$ with the same ${\mathcal{V}}_{\tilde h^0}$ quantum length by  \cref{T:quantumlength}. It is natural to wonder if this actually determines the reverse flow. That is to ask:  could there be any other Loewner flow with the property that intervals of identical quantum length on either side of zero are being zipped together?

We will now show that the answer to this question is no, and hence the Loewner flow for $t\ge 0$ is entirely determined by $\tilde{h}^0$.

\begin{theorem}\label{T:uniqueness} Let $(\tilde h^t,\eta^t)_{t\in \R}$  be a capacity zipper as in \cref{R:cznormed}, with reverse Loewner flow $(f_t)_{t\ge 0}$. Then for $t>0$ the following holds almost surely. If $\hat f_t: \H \to \hat H_t := \hat f_t(\H)$ is a conformal isomorphism such that:
\begin{itemize}
\item  $\hat H_t$ is the complement of a simple curve $\hat\eta^t$,
\item $\hat f_t$ has the hydrodynamic normalisation $\lim_{z\to \infty} \hat f_t(z) - z = 0$;
\item $\hat f_t$ has the property that $\hat f_t(z^-) =\hat f_t(z^+)$ as soon as ${\mathcal{V}}_{\tilde h^0} ([z^- , 0 ] ) = {\mathcal{V}}_{\tilde h^0} ([0, z^+])$ and $f_t(z^-)\in \H\cup \{0\}$;
\end{itemize}
then $\hat f_t = f_t$ and $\hat \eta^t = \eta^t$. In particular, the reverse Loewner flow $(f_t)_{t\ge 0}$ is determined by $\tilde h^0$ only (and hence $((\tilde h^t, \eta^t))_{t\ge 0}$ is entirely determined by $(\tilde h^0,\eta^0)$).
\end{theorem}

\begin{proof}
Before we start the proof, we recall from the definition of the capacity zipper in \cref{T:CZ}, that we only have defined the reverse Loewner flow as being coupled to $\tilde h^0$ in a certain way specified by the application of Kolmogorov's theorem. Usually, proving that objects coupled to a GFF are determined by it can be quite complicated (for example, this is the case in the setup of imaginary geometry, or when making sense of level lines of the GFF).

Here the proof will turn out to be quite simple, given some  classical results from the literature. Indeed consider
$$
\phi = \hat f_t \circ f_t^{-1}.
$$
A priori, $\phi$ is a conformal isomorphism on $f_t(\H) = H_t$, and its image is $\phi(H_t) = \hat H_t$. However, because of our assumptions on $\hat f_t$ (and the properties of $ f_t$), the definition of $\phi$ can be extended unambiguously to all of $\H$. Moreover when we do so, the extended map is a homeomorphism of $\H$ onto $\H$, which is conformal off the curve $\eta^t([0,t])$. Thus the theorem will be proved if we can show that any such map must be the identity. In the terminology of complex analysis, this is equivalent to asking that the curve $\eta^t([0,t])$ is a \emph{removable} set. Now, by a result of Rohde and Schramm \cite{RohdeSchramm}, the complement $H_t$ of the curve is almost surely a H\"older domain for $\kappa<4$ (or $\gamma<2$), and by a result of Jones and Smirnov \cite{JonesSmirnov} it follows that $\eta^t([0,t])$ is a removable set. Hence the theorem follows.
\end{proof}

%\begin{rmk}
%In \cite{zipper}, the theorem states that almost surely for all $t>0$ there is a unique $f_t$ satisfying the above properties. However, a proof of that statement would require that if $f_t$ is a reverse SLE flow, then $f_t(\H)$ is almost surely a H\"older domain for all $t>0$ simultaneously. Such a statement does not appear to be known in the literature currently. We have thus decided to state the theorem for a fixed $t>0$ rather than for all $t>0$ simultaneously.
%\end{rmk}

\begin{rmk} By the same argument, it also holds that for the quantum zipper $(h^t,\zeta^t)_{t\in \R}$ of \cref{T:QZ} $(h^t,\zeta^t)$ is almost surely determined by $(h^0,\zeta^0)$ for any $t>0$.  In the language of conformal welding $(h^t,\zeta^t)$ is obtained from $(h^0,\zeta^0)$ by welding the interval on the left of $0$ with $h^0$ quantum length $t$ to the interval on the right of $0$ with $h^0$ quantum length $t$ (and pushing through $\zeta^0$ by the resulting conformal isomorphism). \end{rmk}

\subsection{Slicing a wedge with an SLE}

In this section we complement our previous discussion by the following remarkable theorem due to Sheffield \cite{zipper}. This result is fundamental to the theory developed in \cite{DuplantierMillerSheffield}, where the main technical tool is a generalisation of the result below.

\begin{figure}
	\begin{center}
		\includegraphics[scale=.6]{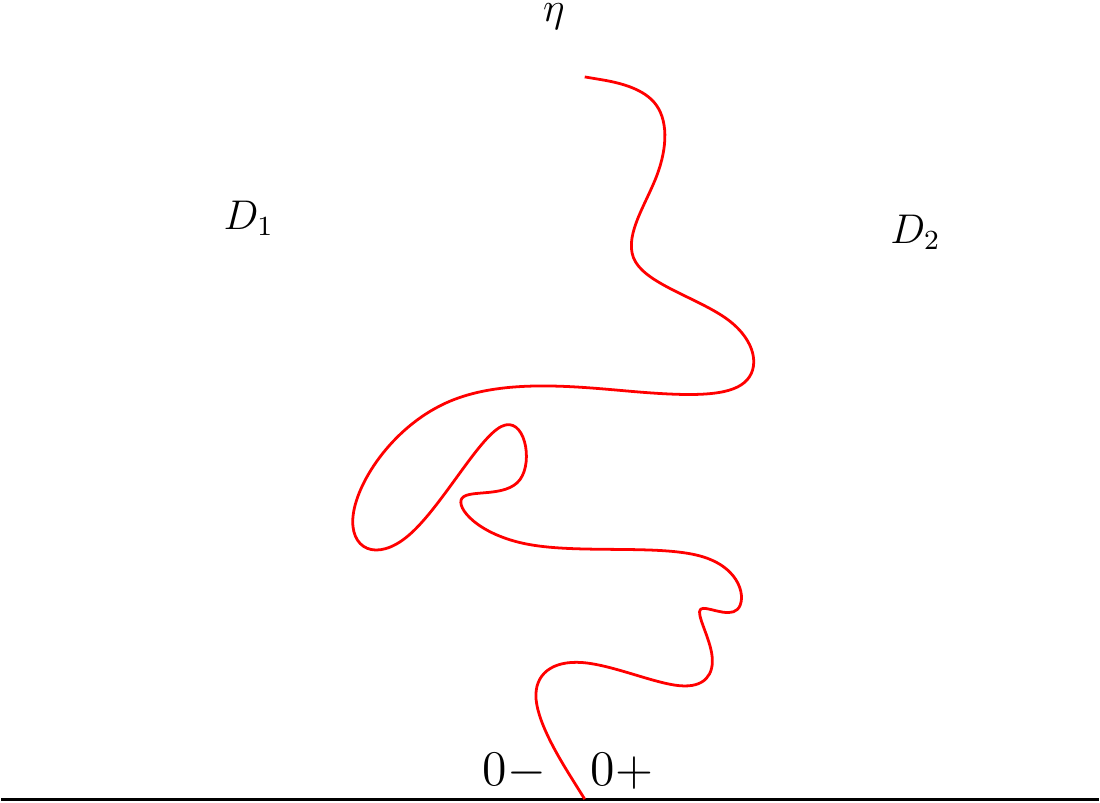}
	\end{center}
	\caption{An independent SLE slices an $(\gamma - 2/\gamma)$-thick wedge into two independent $\gamma$-thick wedges.}\label{fig:Wslice}
\end{figure}

Suppose we are given a $(\gamma - 2/\gamma)$-quantum wedge $(\H,h,0,\infty)$ in some embedding, and an independent SLE$_\kappa$ curve $\eta$ with $\kappa = \gamma^2 < 4$. Then the curve $\eta$ slices the wedge into two surfaces (see picture).
The result below says that \emph{as quantum surfaces} these are independent, and that they are both $\gamma$-thick wedges. See \cref{fig:Wslice}.

\begin{theorem}\label{T:sliced}
Suppose we are given an $(\gamma - 2/\gamma)$-quantum wedge $(\H,h,0,\infty)$ in the unit circle embedding, and an independent SLE$_\kappa$ curve $\eta$ with $\kappa = \gamma^2 < 4$. Let $D_1, D_2$ be the two connected components of $\H \setminus\eta$, whose boundaries contain the negative and positive real lines respectively. Let $h_1 = h|_{D_1}$ and $h_2 = h|_{D_2}$. Then the two surfaces $(D_1, h_1,0-,\infty)$ and $(D_2, h_2,0+,\infty)$ are independent $\gamma$-quantum wedges.
\end{theorem}

\begin{rmk}
	This does \emph{not} imply that the fields, or generalised functions, $h_1$ and $h_2$ are independent. It is a statement about the two doubly marked surfaces $(D_1, h_1,0-,\infty)$ and $(D_2, h_2,0+,\infty)$. So what it does say, for example, is that if $\tilde{h}_1$ and $\tilde{h}_2$ are the fields corresponding to the unit circle embeddings of these surfaces then $\tilde{h}_1$ and $\tilde{h}_2$ are independent.
\end{rmk}

\begin{rmk}
	By the same argument as in the previous subsection, the surfaces $$(D_1, h_1,0-,\infty) \text{ and } (D_2, h_2,0+,\infty)$$ determine $h$ and $\eta$ in the following sense. Suppose that $(\H,\tilde{h}_1,0,\infty)$ and $(\H,\tilde{h}_2,0,\infty)$ are the two unit circle embeddings of these surfaces, and that $(\hat f_1,\hat f_2,\hat{\eta})$ are such that:\begin{itemize} \item $\hat \eta$ is a simple curve from $0$ to $\infty$; \item $\hat f_1$ (resp. $\hat f_2$) is a conformal isomorphism from $\H$ to the left hand side (resp. right hand side) of $\hat\eta$;\item $\hat f_1,\hat f_2$ extend to $\R$ in such a way that for any $x^\pm\in \R_\pm$ with ${\mathcal{V}}_{\tilde{h}^1}([0,x^+])={\mathcal{V}}_{\tilde{h}^2}([x^-,0])$ we have $\hat f_1(x^+)=\hat f_2(x^-)$.\end{itemize} Then if $\hat{h}$ is defined by setting it equal to $\hat f_1(\tilde h_1)$ (resp. $\hat f_2(\tilde h_2))$ on the left hand side (resp. right hand side) of $\hat\eta$, we have that with probability one, $(\hat h,\hat \eta)=(\phi(h),\phi(\eta))$ for some simple scaling map $\phi:z\mapsto az$.
\end{rmk}

We also remark that the choice of embedding for the $(\gamma-2/\gamma)$-wedge in \cref{T:sliced} does not matter, which can be argued as follows. Suppose that $(\H,h,0,\infty)$ is some parametrisation of a $(\gamma-2/\gamma)$-quantum wedge and that $\eta$ is an SLE$_\kappa$ that is independent of $h$. Then there exists a scaling map $\varphi:\H\to \H$ such that $(\H, \varphi(h), 0,\infty)$ is the unit circle embedding of the quantum wedge. Since $\varphi$ is independent of $\eta$ and SLE is scale invariant, $\varphi(\eta)$ is an SLE$_\kappa$ that is independent of $\varphi(h)$. Thus, applying \cref{T:sliced}, we see that the two quantum surfaces obtained by slicing $\varphi(h)$ along $\varphi(\eta)$ are two independent $\gamma$ quantum wedges. On the other hand, these surfaces are by definition equivalent to the two surfaces obtained by slicing $h$ along $\eta$. This means that the latter pair also have the law (as doubly marked quantum surfaces) of two independent $\gamma$ quantum wedges.

\begin{proof}[Proof of \cref{T:sliced}]It is clear from the definition that $(D_1, h_1,0-,\infty)$, $(D_2, h_2,0+,\infty)$ almost surely have finite LQG areas in neighbourhoods of $0-$ and $0+$ respectively, and infinite LQG areas in neighbourhoods of $\infty$. Therefore, we can define unique conformal isomorphisms $\phi_1:D_1\to \H$ sending $0-\to 0$ and $\infty\to \infty$ and $\phi_2:D_2\to \H$  sending $0+\to 0$ and $\infty\to \infty$, so that $(\H,\phi_i(h_i),0,\infty)$ gives LQG area one to the upper unit semidisc $B(0,1)\cap \H$ for $i=1,2$.  Recall that we refer to $\phi_i(h_i)$ as the canonical description of the surface $(D_i,h_i,0\pm,\infty)$, and we continue to use the ``change of coordinate'' notation \eqref{eqn:coc_not} for conformal isomorphisms applied to fields. It clearly suffices to show that for any large semidisc $K\subset \H$, $(\phi_1(h_1)|_K, \phi_2(h_2)|_K)$ agrees in law with $(h^{\wg}_1|_K, h^{\wg}_2|_K)$ where $h^{\wg}_1$ and $h^{\wg}_2$ are independent, and each has the law of the canonical description of a $\gamma$-quantum wedge.  (The reason we choose to work with the canonical description rather than the unit circle embedding here is simply to avoid any ambiguity concerning the a priori existence of the maps $\phi_1$ and $\phi_2$.)  %In fact, by symmetry of the marginal laws, it suffices to show the same statement with $h_1^{\wg}$ replaced by any other field $\hat{h}_1$, that is independent of $h_2^{\wg}$.
	
To show this equality in law, we need to appeal to the results of the previous section: in particular \cref{lem::zip_to_wedge} and \cref{T:quantumlength}. Consider the process $((\tilde{h}^t, \eta^t)_{t\ge 0}, Z)$ under the law $\mathbf{Q}$ from \cref{lem::q_alt}, and in this set up, let $Y$ denote the point to the left of zero such that the $\tilde{h}^0$ boundary length of $[Y,0]$ is equal to that of $[0,Z]$. Write $h_Z^C$ for the canonical description of $ (H_Z,\tilde{h}^0+C, Z, \infty)$ and $h_Y^C$ for the canonical description of $ (H_Y,\tilde{h}^0+C, Z, \infty)$ where $H_Z$ and $H_Y$ are the connected components of $\H\setminus\tilde{\eta}^0$ containing $Z$ and $Y$ respectively.
Combining \cref{lem::zip_to_wedge} and \cref{T:quantumlength} gives that:
\begin{claim} We can couple pairs of fields with
\begin{itemize}\setlength{\itemsep}{0em}
	\item the joint law of $(h_Y^C, h_Z^C)$ under $\mathbf{Q}$, and
	\item the joint law of  $(\phi_1(h_1),\phi_2(h_2))$ described in the first paragraph,
\end{itemize}
so that they agree when restricted to $K$, with probability arbitrarily close to one as $C\to \infty$.
\end{claim} \begin{proof}[Proof of claim] (See \cref{fig:sliceproof}).
First we observe that one (slightly convoluted!) way to sample a pair with the law of $(h_Y^C,h_Z^C)$ under $\mathbf{Q}$ is to: \begin{enumerate} \item[(1)] consider the ``zipper'' $((\tilde{h}^t, \eta^t)_{t\ge 0}, Z)$ under $\mathbf{Q}$ and apply the conformal isomorphism $f_{\tau_Z}^C$ that zips up $Z$ to $0$ and then scales $\H$ so that $f_{\tau_Z}^C(\tilde{h}^0)$ is in the unit circle embedding;
		\item[(2)] then, restrict the field $f_{\tau_Z}^C(\tilde{h}^0+C)$ to the left and right of $f_{\tau_Z}^C(\eta^0)$, and apply conformal isomorphisms from these left and right hand sides to $\H$, such that the resulting fields (under the change of coordinates formula) are the canonical descriptions of these two surfaces.
	\end{enumerate}
Here we are using the fact, due to \cref{T:quantumlength}, that $Y$ is zipped up to $0$ at exactly the same time as $Z$.

On the other hand, \cref{lem::zip_to_wedge} says that we can couple $(f_{\tau_Z}^C(\tilde{h}^0+C),f_{\tau_Z}^C(\eta^0))$ with $(h,\eta)$ as in the statement of the present theorem, so that they agree in any large semidisc $K'$, with probability arbitrarily close to one as $C\to \infty$. Consequently, if we restrict the field $f_{\tau_Z}^C(\tilde{h}^0+C)$ to the left and right of $f_{\tau_Z}^C(\eta^0)$, and apply conformal isomorphisms as in the second step of the previous bullet point, then the resulting pair of fields can be coupled with $(\phi_1(h_1),\phi_2(h_2))$ so that they agree when restricted to $K$ with arbitrarily high probability.

Combining these two paragraphs yields the claim.
\end{proof}
\begin{figure}
	\begin{center}
		\includegraphics[scale=.8]{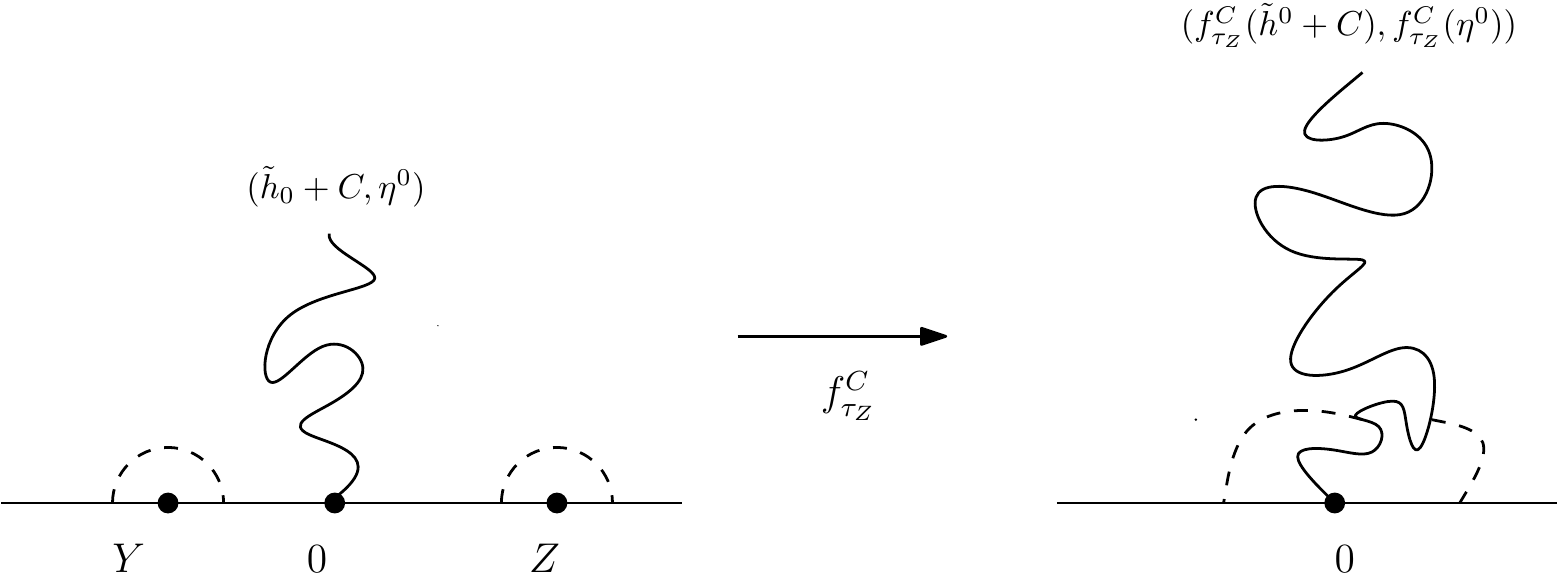}
	\end{center}
	\caption{The surfaces to the left and right of $\eta^0$ on the left hand picture (defined using the field $\tilde{h}^0+C$ and marked at $(Y,\infty)$ and $(Z,\infty)$) have canonical descriptions given by $(\H,h_Y^C,0,\infty)$ and $(\H,h^C_Z,0,\infty)$. So the same is true, by definition, for the surfaces to the left and right of the curve on the right hand picture (defined using the field $f_{\tau_Z}^C(\tilde{h}^0+C)$). But \cref{lem::zip_to_wedge} says that for $C$ large, the joint law of the field and curve on the right hand picture is very close to that of $(h,\eta)$ from the statement of \cref{T:sliced}. So, the law of the canonical descriptions of the surfaces to the left and right of the curve is very close to that of $(\H,\phi_1(h_1),0,\infty), (\H,\phi_2(h_2),0,\infty)$. Hence we can approximate the joint law of $(\phi_1(h_1),\phi_2(h_2))$ by that of $(h_Y^C,h_Z^C)$ for $C$ large. }
\label{fig:sliceproof}\end{figure}

\noindent So, with the claim in hand, it actually suffices to show that we can couple $(h_Y^C,h_Z^C)$ with $(h_1^{\wg},h_2^{\wg})$ (recall that the latter are an pair of independent $\gamma$-wedge fields in their canonical descriptions) %where
%\begin{itemize}\setlength{\itemsep}{0em}
%	\item the two components are independent,
%	\item $\hat{h}^1$ is some field that gives LQG mass one to $\D_+$, and
%	\item $h_2^{\wg}$ is the canonical description of a $\gamma$-quantum wedge
%\end{itemize} (as described before)
so that their restrictions to $K$ agree with probability arbitrarily close to $1$ as $C\to \infty$. The idea is that when $C$ is very large, the restrictions of $h^C_Y$ and $h^C_Z$ to $K$ will correspond to images -- under the conformal change of coordinates \eqref{eqn:coc_not} -- of $\tilde{h}^0+C$ restricted to very tiny neighbourhoods of $Z$ and $Y$. Roughly speaking, these restrictions become independent in the limit as the size of the neighbourhoods goes to $0$, and furthermore, the field near $Z$ (and by symmetry near $Y$) converges to a $\gamma$-quantum wedge field.

To be more precise, let us consider a sample $(\tilde{h}^0, Z)$ from $\mathbf{Q}$, together with a field $\tilde{h}=\hat h+(\gamma-2/\gamma)\log(|\cdot|^{-1})$, where $\hat h$ is a Neumann GFF normalised to have average $0$ on the upper unit semicircle that is \emph{independent} of $\tilde{h}^0$. Then we have the following: %and write $\nu^-_\delta=\nu_{\tilde h^0}([Z-\delta, Z])$, $\nu^+_\delta=\nu_{\tilde h^0}([Z,Z+\delta])$. We will show that:

\begin{lemma}\label{lem:zooming_twopoints}
%	For any $\delta>0$,
As above, let $(\tilde{h}^0, Z)$ have their $\mathbf{Q}$ joint law, and let $\tilde{h}=\hat h+(\gamma-2/\gamma)\log(|\cdot|^{-1})$, where $\hat h$ is a Neumann GFF normalised to have average $0$ on the upper unit semicircle, that is independent of $\tilde{h}^0$. Then the total variation distance between $$ (\tilde{h}^0|_{\overline{B(Z,\eps)\cap \H}}, \tilde{h}^0|_{\H\setminus B(Z,1)},{\mathcal{V}}_{\tilde{h^0}}[1,Z], {\mathcal{V}}_{\tilde{h^0}}[Z,2]) $$ and $$(\tilde{h}(\cdot-Z)|_{\overline{B(Z,\eps)
	\cap \H}}, \tilde{h}^0|_{\H\setminus B(Z,1)},{\mathcal{V}}_{\tilde{h^0}}[1,Z], {\mathcal{V}}_{\tilde{h^0}}[Z,2])$$ converges to $0$ as $\eps\to 0$. %(Here $(\tilde{h}^0, Z)$ are sampled from $\mathbf{Q}$, and $\tilde{h}$ is sampled independently with the correct law.)
\end{lemma}
 In words this says that \emph{conditionally} on $\tilde{h}^0$ outside of $B(Z,1)$ \emph{and} on the boundary lengths ${\mathcal{V}}_{\tilde{h^0}}[1,Z]$, ${\mathcal{V}}_{\tilde{h}^0}[Z,2]$, the law of $\tilde{h}^0$ restricted to $\overline{B(Z,\eps)\cap \H}$ is very close in total variation distance to the field $\tilde{h}$ recentred at $Z$ and restricted to $\overline{B(Z,\eps)\cap \H}$. %In particular, note that $\tilde{h}^0$ restricted to any small enough neighbourhood of $Y$ is measurable with respect to the information being conditioned on here.

 \medskip
 Before proving the lemma, let us first see how it allows us to conclude the proof of the theorem. From now on, we assume that $K\subset \H$ is large, fixed semidisc. Consider a pair $(h_Y^C,\tilde{h}^C)$ where $h_Y^C$ has its $\mathbf{Q}$ law, and $\tilde{h}^C$ is independent of $h_Y^C$ having the law of the canonical description of $(\H,\tilde{h}+C,0,\infty)$. The consequence of \cref{lem:zooming_twopoints} is that by taking $\eps$ very small and then $C$ sufficiently large, we can couple the joint law of $(h^C_Y,h_Z^C)$ with that  of the pair $(h_Y^C, \tilde{h}^C)$, so that the fields agree when restricted to $K$ with probability arbitrarily close to one. Since the law of $\tilde{h}^C|_K$ converges in total variation distance to $h^{\wg}_2|_K$ as $C\to \infty$, see \cref{cor:wedge_cd_conv}, this means that we can couple $(h_Y^C,h_Z^C)$ with $(h_Y^C,h_2^{\wg})$ (where the latter pair are independent) so that they agree when restricted to $K$ with arbitrarily high probability as $C\to \infty$.

To finish the proof, we observe that by symmetry, $h_Y^C$ has the same law as $h_Z^C$ for each $C$. Since the argument above clearly gives that $h_Z^C|_K\to h_2^{\wg}|_K$ in total variation distance as $C\to \infty$, it must therefore also be the case that $h_Y^C|_K$ converges in total variation distance to $h_1^{\wg}|_K$ as $C\to \infty$. Thus $(h_Y^C,h_2^{\wg})$ can be coupled with $(h_1^{\wg},h_2^{\wg})$ so that the fields agree when restricted to $K$ with arbitrarily high probability as $C\to \infty$. Putting this together with the previous paragraph, we obtain the desired result.
\end{proof}

\begin{proof}[Proof of \cref{lem:zooming_twopoints}]
We first claim that for any $\delta>0$, \begin{equation}\label{eq:wedgeconditionedout}d_{TV}\left((\tilde{h}^0|_{\overline{B(Z,\eps)\cap \H}}, \tilde{h}^0|_{\H\setminus B(Z,\delta)}) \, , \, (\tilde{h}(\cdot - Z)|_{\overline{B(Z,\eps)\cap \H}}, \tilde{h}^0|_{\H\setminus B(Z,\delta)})\right) \to 0\end{equation} as $\eps\to 0$.
 Indeed, by \cref{lem::q_alt}, the $\mathbf{Q}^Z$ (that is, $\mathbf{Q}(\cdot|Z)$) law of $\tilde{h}^0$ recentred at $Z$ is that of $h'+(\gamma - 2/\gamma)\log(|\cdot|^{-1})+\mathfrak{h}$, where $h'$ is a Neumann GFF normalised to have average $0$ on the upper unit semicircle centred at $10$ and $\mathfrak{h}$ is a harmonic function that independent of $h'$ and is deterministically bounded in $B(Z,1)$. Hence \eqref{eq:wedgeconditionedout} follows from \cref{lem:hDNcircgood} and \cref{rmk:hDNcircgood}.

We will now extend this in the following way. We are going to show that the law of ${\mathcal{V}}_{\tilde{h}^0}([1,Z])$ is basically the same (when $\eps$ is small enough) whether we condition on $\tilde{h}^0$ restricted to ${\H\setminus B(Z,1)}$ \emph{and} ${\overline{\H\cap B(Z,\eps)}}$, or just restricted to ${\H\setminus B(Z,1)}$: see \eqref{eq:lengthind0}. The basic idea for the proof is that, given the restriction of $\tilde{h}^0$ to ${\H\setminus B(Z,1)}$, the restriction of $\tilde{h}^0$ to $\overline{\H\cap B(Z,\eps)}$ has a very tiny influence on the boundary length of $[1,Z]$ when $\eps$ is small. On the other hand, there is quite a bit of variation in the boundary length coming from sources completely independent of $\tilde{h}^0|_{\overline{B(Z,\eps)\cap \H}}$. To argue this rigorously, we will use the Fourier decomposition of the free field, similarly to the argument \cite{zipper}. %to show that for any $\delta>0$
%\begin{equation} \label{eq:wedgeconditionedonL} d_{TV}\left((\tilde{h}^0|_{\overline{B(Z,\eps)\cap \H}}, \tilde{h}^0|_{\H\setminus B(Z,\delta)},\nu^-_\delta, \nu^+_{\delta}) \, , \,  (\tilde{h}(\cdot-Z)|_{\overline{B(Z,\eps)
	%	\cap \H}}, \tilde{h}^0|_{\H\setminus B(Z,\delta)},\nu^-_\delta, \nu^+_{\delta})\right) \to 0\end{equation}
%	as $\eps\to 0$, where
%$\nu^-_\delta=\nu_{\tilde h^0}([Z-2\delta, Z])$, $\nu^+_\delta=\nu_{\tilde h^0}([Z,Z+2\delta])$. Note that this immediately implies \cref{lem:zooming_twopoints}. For simplicity of presentation, let us show \eqref{eq:wedgeconditionedonL} only for the first three coordinates of each quadruple (i.e., we will omit $\nu_{\delta}^+$). This makes no difference whatsoever to the method of proof, apart from reducing the amount of notation.
%To prove \eqref{eq:wedgeconditionedonL} (with just $\nu_{\delta}^-$) the very basic idea is that, given $\tilde{h}^0|_{\H\setminus B(Z,\delta)}$, the restriction of $\tilde{h}^0$ to $\overline{B(Z,\eps)\cap \H}$ has a very tiny influence on $\nu_{\delta}^-$ when $\eps$ is small, while on the other hand, there is quite a bit of variation in $\nu_\delta^-$ coming from sources completely independent of $\tilde{h}^0|_{\overline{B(Z,\eps)\cap \H}}$. This idea can be made rigorous in the following way.

\begin{figure}
	\centering
	\includegraphics[width=0.7\textwidth]{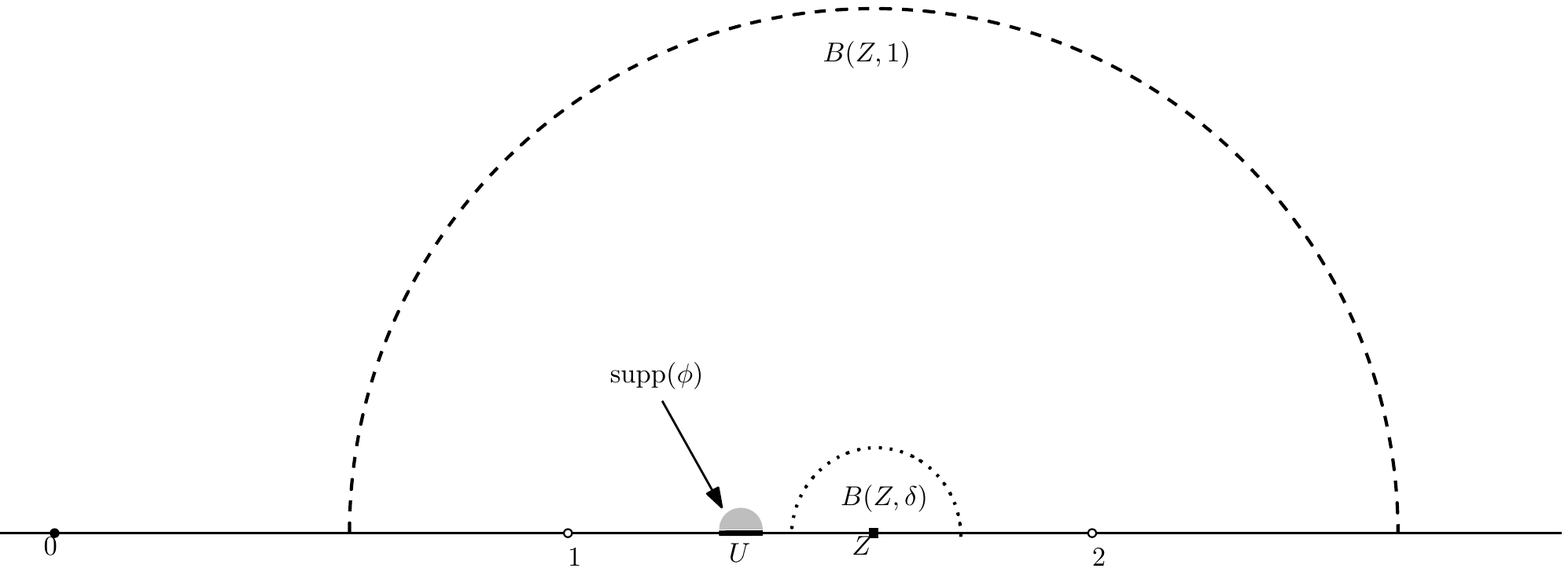}
\end{figure}
We take $\delta>0$ small, and fix a function $\phi$ that is smooth, positive and supported in the upper unit semidisc of radius $\delta/4$ centred at $Z-3\delta/2$, with $(\phi,\phi)_{\nabla}=1$. Let us write $U:=[Z-7\delta/4,Z-5\delta/4]$, $U^c=[1, Z]\setminus U$. This will be non-empty with arbitrarily high probability if $\delta$ is small enough, so let us assume from now on that $Z$ is such that this is the case. Then by
\cref{D:NeumannGFF2} and \cref{D:NGFFnorm}, we can decompose $\tilde{h}^0=X\phi+h$ where $X$ is Gaussian and $h$ is independent of $X$. %Notice that $X$ is independent of $\tilde{h}^0$ restricted to $\overline{B(Z,\eps)\cap \H}$ by definition.

Next, we observe that due to the decomposition of $\tilde{h}^0$, the conditional law of ${\mathcal{V}}_{\tilde h^0}(U)$ given $h|_{\H\setminus B(Z,\delta)}$ almost surely has smooth density $F^{h|_{\H\setminus B(Z,\delta)}}$ with respect to Lebesgue measure: indeed, given $h$ restricted to $\H\setminus B(Z,\delta)$, ${\mathcal{V}}_{\tilde h^0}(U)$ is almost surely smooth and increasing in $X$.
%Note that $F$ depends on $\hat{h}^0|_{\H\setminus B(Z,\delta/2)}$, but we don't include this in the notation.
In particular,  \begin{equation*}\label{eqn:leb_density}
\text{the conditional law of } {\mathcal{V}}_{\tilde{h}^0}([1,Z]) \text{ given } h\text{ has density } \propto F^{h|_{\H\setminus B(Z,\delta)}}(\cdot - {\mathcal{V}}_{h}(U^c))
\end{equation*} with respect to Lebesgue measure.
Using the fact that $F$ is smooth,  that ${\mathcal{V}}_{h}$ almost surely does not have an atom at $Z$, and  \eqref{eq:wedgeconditionedout} applied with $\delta'\ll \delta$, we may deduce from this that for any
$x\in \R$,  the quantity
$$\E( F^{h|_{\H\setminus B(Z,\delta)}}(x - {\mathcal{V}}_{h}(U^c)) \, \big| \, h_{\H\setminus B(Z,\delta)})-\E(F^{h|_{\H\setminus B(Z,\delta)}}(x - {\mathcal{V}}_{h}(U^c)) \, \big| \, h_{\H\setminus B(Z,\delta)},h|_{\overline{B(Z,\eps)\cap \H}})$$
tends to $0$ almost surely as $\eps\to 0$. This is important because it means that
\begin{equation}\label{eq:lengthind0}
\mathrm{d}_{\tv}\left( \cL( {\mathcal{V}}_{\tilde{h}^0}([1,Z])\, | \, h|_{\H\setminus B(Z,\delta)}), \cL({\mathcal{V}}_{\tilde{h}^0}([1,Z]) \, | \, h|_{\H\setminus B(Z,\delta)}, h|_{\overline{B(Z,\eps)\cap \H}})\right)\to 0
\end{equation}in probability $\eps\to 0$ (where $\cL(Y_1|Y_2)$ denotes the law of $Y_1$ conditioned on  $Y_2$). In fact, since  $$h|_{\overline{B(Z,\eps) 	\cap \H}}=\tilde{h}^0|_{\overline{B(Z,\eps) 	\cap \H}},$$ and by combining with \eqref{eq:wedgeconditionedout} this actually means that \begin{equation*}
d_{\tv}\left(({\mathcal{V}}_{\tilde{h}^0}([1,Z]),\tilde{h}^0|_{\overline{B(Z,\eps)\cap \H}}, h|_{\H\setminus B(Z,\delta)}) \, , \, ({\mathcal{V}}_{\tilde{h}^0}([1,Z]),\tilde{h}(\cdot - Z)|_{\overline{B(Z,\eps)\cap \H}}, h|_{\H\setminus B(Z,\delta)})\right) \to 0
\end{equation*} in probability $\eps\to 0$. This is extends with exactly the same argument (but a little more notation) to the same statement with ${\mathcal{V}}_{\tilde h^0}([1,Z]), {\mathcal{V}}_{\tilde h^0}([Z,2])$ in place of just ${\mathcal{V}}_{\tilde h^0}([1,Z])$. Putting this together with the fact that $h=\tilde{h}^0$ outside of $B(Z,1)$ completes the proof.
%The upshot of this is that the Radon--Nikodym derivative between the laws of
%$$ \left({\mathcal{V}}_{\delta}^- \text{ given } (\hat{h}^0|_{\H\setminus B(Z,\delta/2)}, \tilde{h}^0|_{\overline{B(Z,\eps)\cap \H}})\right) \;\; \text{ and } \;\; \left(
%{\mathcal{V}}_{\delta}^- \text{ given } \hat{h}^0|_{\H\setminus B(Z,\delta/2)}\right)$$ is equal to
%$$\frac{\E[F^{\hat{h}^0|_{\H\setminus B(Z,\delta/2)}}({\mathcal{V}}_\delta^- - {\mathcal{V}}_{\hat{h}^0}(U^c))\, | \, \tilde{h}^0|_{\overline{B(Z,\eps)\cap \H}},\hat{h}^0|_{\H\setminus B(Z,\delta/2)}]}{\E[F^{\hat{h}^0|_{\H\setminus B(Z,\delta/2)}}({\mathcal{V}}_{\delta}^- - {\mathcal{V}}_{\hat{h}^0}(U^c))\, | \, \hat{h}^0|_{\H\setminus B(Z,\delta/2)}]} $$

\end{proof}

%\subsection{Exercises}

%\begin{enumerate}

%\ind{SLE$_\kappa(\rho)$}
%\item Let $f_t$ be a reverse SLE$_\kappa(\rho)$ Loewner flow, let
%$$
%h_0 (\cdot) = \tilde h(\cdot) + \frac2{\sqrt{\kappa}} \log |\cdot | + \frac{\rho}{2\sqrt{\kappa}} G(x,\cdot)
%$$
% be independent of $(f_t)$, where $G(x,y) = -\log|x-y|$ is the whole plane Green function and $\tilde h$ is a Neumann GFF with some normalisation.

%For all $0\le t \le \tau$ set
%$$
%h^\rho_t = \tilde h \circ f_t +  \frac2{\sqrt{\kappa}} \log |f_t | + Q \log |f_t'| + \frac{\rho}{2\sqrt{\kappa}} G(f_t(x),f_t(z)).
%$$
%Show that
%$$
% h^\rho_{t\wedge \tau} = h_0
% $$
% in law, as distributions modulo additive constants. (This is Theorem 4.5 in \cite{zipper}). Give two proofs of this result: one based on stochastic calculus in the manner of the proof of Theorem \ref{T:coupling}, and another one using Lemma \ref{L:rho}.

%\ellen{*Add more exercises! Natural parametrisation?*}

%\end{enumerate} 

\newpage 

\section{Liouville quantum gravity as a mating of trees}\label{S:MOT}

% !TEX root = master.tex

\subsection{Orientation}

In this chapter we take forward the ideas developed in Chapter \ref{S:zipper} and obtain a beautiful and important description of the way that a certain quantum cone can be explored by an independent variant of SLE called \emph{space-filling SLE}. This description has many important implications, and we already emphasise the following points.

\begin{itemize}

\item On the one hand, this shows that a quantum cone, considered as a random surface decorated with a designated space-filling path, can rigorously be described as the ``mating'' (that is, gluing or welding) of two correlated (infinite) continuum random trees. This is the so called \textbf{peanosphere} description of a quantum cone. It also has analogues for other quantum surfaces; see Section \ref{SS:extension}.

\item On the other hand, this construction is the direct continuum analogue of the discrete bijection due to Sheffield in \cite{Sheffield_burger}  for random planar maps weighted by the self dual Fortuin--Kasteleyn percolation model, as was presented in Chapter \ref{S:maps}. Hence a particular consequence of this work (as developed in \cite{FKstory1, FKstory2} and \cite{FKstory3}) is that, at least in the so called peanosphere sense (which is a relatively weak notion of convergence), these random planar maps can be  proven to converge to  quantum cones. 

\end{itemize}

This very fruitful approach was developed in the seminal paper of Duplantier, Miller and Sheffield \cite{DuplantierMillerSheffield}. Since this paper is long and difficult, we will not aim to present complete proofs of their main theorems; rather, we will state precise results and hope to convey some of the key ideas that are used in their proofs.

To state the main theorems, we must first explain the construction of the aforementioned space-filling SLE. The details of the construction are not straightforward, and in fact rely on a whole other body of work (the so called \emph{imaginary geometry} of Miller and Sheffield, \cite{MSIG1, IG2, IG3} and especially \cite{IG4}), which falls outside of the scope of this book. 
\dis{Although we will give a complete and self contained introduction to whole plane SLE (Section \ref{SS:wpsle}), and space-filling SLE$_\kappa$ for $\kappa\ge 8$ (Section \ref{SS:spacefillingSLE}), the construction of space-filling SLE$_\kappa$ in the case $\kappa\in (4,8)$ (Section \ref{SS:sf48}) will be explained rather than fully justified.} In Section \ref{SS:cuttingwelding} we will state a cutting/welding theorem for (thick and thin) quantum wedges, analogous to but more complicated than the welding statement of Theorem \ref{T:sliced}, which is crucial to the ``mating of trees''  theorem of \cite{DuplantierMillerSheffield}, which we will state in Section \ref{SS:MOTstatement}. In Section \ref{SS:MOTdiscussion} we will discuss the implications of this theorem, in relation to the two bullet points above. In Section \ref{SS:proofMOT} we will give a proof of the main theorem in the case $\kp\in (4,8)$, admitting the welding theorems from Section \ref{SS:cuttingwelding} and a stationarity statement (analogous to but more complicated than the stationarity of the quantum zipper in Proposition \ref{P:qwkey}). This proof also partially covers the case $\kp\ge 8$, up to a certain step that we will explain properly in that section.
%We conclude with an important application (to so-called CRT mated maps) in Section \ref{SS:CRTmm}, followed by a short discussion of some extensions of Duplantier--Miller--Sheffield's results in Section \ref{SS:extension}.

\subsection{Whole plane \texorpdfstring{$\mathrm{SLE}_\kappa$}{TEXT} and \texorpdfstring{$\mathrm{SLE}_\kappa(\rho)$}{TEXT} }
\label{SS:wpsle}

\subsubsection{Whole plane \texorpdfstring{$\mathrm{SLE}_\kappa$}{TEXT}}

\begin{definition}[Whole plane SLE$_\kappa$] \label{D:WPSLE}
	For $\kappa> 0$, whole plane SLE$_\kappa$ in $\C$ from $0$ to $\infty$ is defined to be the 
	 collection of  maps $(g_t)_{t\in \R}$ that solve the \emph{whole plane Loewner equation} for each $z\in \C\setminus \{0\}$:
	\begin{equation}\label{E:WPLE}
		\partial_t g_t(z) = g_t(z) \frac{U_t+g_t(z)}{U_t-g_t(z)}; \quad \lim_{t\to -\infty} e^t g_t(z)=z; \quad t\in (-\infty, \zeta(z))
	\end{equation}
where $U_t = e^{i\sqrt{\kappa}B_t}$,  $B$ is a standard \emph{two-sided} Brownian motion, and for each $z$, $\zeta(z):=\inf\{t\in \R:U_t=g_t(z)\}$.
\end{definition}

We emphasise that the map $g_t$ in the definition above is defined for all $t\in \R$, not just for $t\ge 0$ as is the case for chordal or radial Loewner chains (see Appendix \ref{app:sle} and Appendix \ref{S:App_radial}).

For a given realisation of $B$, existence and uniqueness of $g_t(z)$ for each $z\in \C\setminus \{0\}$ and $t\in (-\infty,\zeta(z))$ follows from standard ODE theory. If $K_t:=\overline{\{z: \zeta(z)\le t\}}$ for $t\in \R$ are the \emph{whole plane Loewner hulls} generated by $B$, then it can be shown that $g_t(z)$ is indeed a conformal isomorphism from $\hat{\C}\setminus K_t\to \hat{\C}\setminus \bar{\D}$ and that $\mathrm{cap}(K_t):=\lim_{z\to \infty} z/g_t(z) = e^t$ for each $t$; see for example \cite{Lawlerbook}.

 In fact, many more properties can be deduced immediately from the following connection to \emph{radial} SLE$_\kappa$. In some sense, the following lemma suggests that whole plane SLE$_\kappa$ should be viewed as a bi-infinite time version of radial SLE$_\kappa$.

\begin{lemma}
Let $(g_t,K_t)_{t\in \R}$ be the conformal isomorphisms and Loewner hulls associated to a driving function $(U_t)_{t\in \R}$ as in \eqref{E:WPLE}. 
%a whole plane SLE$_\kappa$ process for $\kappa>0$. 
Then for any $t_0\in \R$, $\tilde{K}_t:=g_{t+t_0}(K_{t+t_0}\setminus K_t)$ has the law of a radial 
%SLE$_\kappa$ 
Loewner evolution in $\hat{\C}\setminus K_{t_0}$ from the point $ g_{t_0}^{-1}(U_{t_0})$ to $\infty$. More precisely, the
 hulls $1/g_{t+t_0}(K_{t+t_0}\setminus K_t)$ for $t\ge 0$ are described by a radial Loewner evolution in $\D$ whose driving function is given by $(\bar U_{t_0 + t})_{t\ge 0}$ (the complex conjugate function of $(U_{t_0+t})_{t\ge 0}$).

In particular, if $(K_t)_{t\in \R}$ are the whole plane hulls associated to an SLE$_\kappa$ and $t_0 \in \R$, then  $(K_{t_0 + t})_{t\ge 0}$ are the hulls of a radial SLE$_\kappa$ in $\C\setminus K_{t_0}$ from $g_{t_0}^{-1} (U_{t_0})$ to $\infty$.  
 %SLE$_\kappa$ in $\D$ from $1/z_{t_0}$ to $0$. 
 \label{L:WPRAD}
\end{lemma}

\begin{proof}
It suffices to check that if $\tilde{g}_s:=1/g_s(1/\cdot)$ for $s\in \R$ and $\hat{g}_{t}:=\tilde{g}_{t+t_0}\circ \tilde{g}_{t_0}^{-1}$ for $t>0$ (so that $\hat{g}_t$ is the unique conformal isomorphism from $\D \setminus \{1/g_{t+t_0}(K_{t+t_0}\setminus K_t)\}$ to $\D$ with $\hat{g}_t'(0)=e^t$ for each $t\ge 0$) then $(\hat{g}_t)_{t\ge 0}$ satisfies the radial Loewner equation \eqref{E:RLE} with driving function given by $\bar U_{t + t_0}$. This follows from a simple calculation using \eqref{E:WPLE}, which we leave to the reader (note that \eqref{E:RLE} and \eqref{E:WPLE} are identical, apart from the time domain and the ``initial'' conditions).
\end{proof}

It therefore follows from the corresponding results for radial SLE$_\kappa$ (see Appendix \ref{A:radialSLE}) that for each $\kappa>0$, and given $(B_t)_{t\in \R}$, there almost surely exists a continuous non self crossing curve $\gamma:(-\infty,\infty)\to \C$ such that the unique conformal isomorphism $g_t$ from the unbounded connected component of $\C\setminus \gamma((-\infty,t])$ to $\C\setminus \D$ with $g_t(\infty)=\infty$ and $g_t'(\infty)>0$, solves \eqref{E:WPLE} (and in fact has $g_t'(\infty)=e^{-t}$) as in \cref{D:WPGFF}. The curve starts at $0$ in the sense that $\lim_{t\to -\infty}\gamma(t)=0$ and is transient, that is, $\lim_{t\to \infty} \gamma(t)=\infty$. It also follows that whole plane SLE$_\kappa$ has the same distinct phases as radial (and chordal) SLE$_\kappa$: it is a simple curve for $\kappa\le 4$, is self intersecting but non self crossing and non space-filling for $\kappa\in (4,8)$, and is space filling for $\kappa\ge 8$.

The scaling property of Brownian motion also implies that if $\gamma$ is a whole plane SLE$_\kappa$ from $0$ to $\infty$ and $a\in \C\setminus \{0\}$, then $a\gamma$ (with time reparametrised appropriately) also has the law of a whole plane SLE$_\kappa$ from $0$ to $\infty$. This means that the following definition makes sense. 

\begin{definition}\label{D:WPSLE2}
Let $z_1,z_2\in \hat{\C}$. Whole plane SLE$_\kappa$ from $z_1$ to $z_2$ is defined to be the image of whole plane SLE$_\kappa$ from $0$ to $\infty$, under a M\"{o}bius transformation sending $0$ to $z_1$ and $\infty$ to $z_2$. The law of this process does not depend on the choice of M\"{o}bius transformation.
\end{definition}

With this definition, it is immediate that as a family indexed by $z_1,z_2\in \hat{\C}$, whole plane SLE$_\kappa$ from $z_1$ to $z_2$ is M\"{o}bius invariant (in law). For instance, the whole plane SLE$_\kappa$ from $\infty$ to 0 is obtained by applying the M\"obius inversion $\psi(z) = 1/z$ to the hulls $(K_t)_{t\in \R}$ of Definition \ref{D:WPSLE}. In doing so we obtain hulls $\tilde K_t = \psi(K_t), t \in \R$; note that the parametrisation of $\tilde K_t$ is then such that the capacity seen from 0 of $\tilde K_t$ is always equal to $e^{t}$. In other words, 
\begin{equation}\label{E:capCR}
\mathrm{CR}( 0; \C \setminus \tilde K_t) = e^{-t}
\end{equation}
where we recall that $\mathrm{CR}(x,D)$ stands for the conformal radius of $x$ in $D$.  In this sense, $(\tilde K_t)_{t\in \R}$ is just simply parametrised by \textbf{log conformal radius}. 

%For $\kappa\in (0,8]$, it also has a time reversal property: if $(\eta(t))_{t\in \R}$ is a whole plane SLE$_\kappa$ from $0$ to $\infty$ and $(\tilde{\eta}(t))_{t\in \R}$ is a whole plane SLE$_\kappa$ from $\infty$ to $0$, then 
%$(\tilde{\eta}(-t))_{t\in \R}$ has the same law as $(\eta(t))_{t\in \R}$. This was proven by Dapeng Zhan in \cite{ZhanWP} for $\kappa\in (0,4]$ and extended to $\kappa\in (4,8]$ in \cite{IG4}. For $\kappa>8$, it is known that this is \emph{not} the case (\cite{IG4}).

\medskip We caution that invariance of whole plane SLE$_\kappa$ from $0$ to $\infty$ under the inversion map $z \mapsto 1/z$ does not mean that it is \emph{reversible}. That is, it is \emph{a priori} not obvious that the image of whole plane SLE$_\kappa$ from 0 to $\infty$ under the inversion map coincides in law with its time-reversal. In fact, while chordal SLE$_\k$ is reversible for $\k\in (0, 8]$ (this was proved by Zhan in \cite{Zhan_chordal} for $\kappa\in (0,4]$ and extended to $\kappa\in (4,8]$ in \cite{IG3} using the tools of imaginary geometry), chordal SLE$_\k$ is known to be non reversible for $\kappa > 8$. (This is credited in \cite{IG4} to unpublished work of Rohde and Schramm.) Nevertheless, Viklund and Wang \cite{ViklundWang} conjectured reversibility of whole-plane SLE$_\k$ for $\k >8$ on the basis of an invariance of a certain large deviation functional in the limit $\k \to \infty$; this remarkable conjecture was proved in the no less remarkable work of Ang and Yu \cite{AngYu}.

\subsubsection{Whole plane \texorpdfstring{$\mathrm{SLE}_\kappa(\rho)$}{TEXT} }

\label{A:wpSLE}

In this section we will discuss the definition of SLE$_\kappa(\rho)$ for $\rho> -2$. We will only consider the case of one ``weight'' $\rho$, and the initial force point will (in some sense) be the same as the starting point.

To do this, we need the following lemma.
\begin{lemma}[\cite{IG4}, Proposition 2.1]
	Suppose that $\kappa>0, \rho> -2$ and that $(\tilde{U}_s,\tilde{V}_s)_{s\ge 0}$ solves \eqref{eqn:Rslekpdef} (with $m=1$)\footnote{That is, $(\tilde{U}_s)_s$ is the driving function for a radial SLE$_\kappa(\rho)$ from $\tilde{U}_0$ to $0$, with force point initially at $\tilde{V}_0$, and $(\tilde{V}_s)_s$ is the evolution of $\tilde{V}_0$ under the Loewner flow} and some choice of $\tilde{U}_0,\tilde{V}_0\in \partial \D$.
	There exists a unique time stationary law on continuous processes $(U_t,V_t)_{t\in \R}$, taking values on $\partial \D\times \partial \D$, for which $(U_t,V_t)_{t\ge t_0}$ is equal to the limit in law and in total variation distance of $(\tilde{U}_{t+T},\tilde{V}_{t+T})_{t\ge t_0}$ as $T\to -\infty$ for any $t_0\in \R$. This law does not depend on the choice of $\tilde{U}_0,\tilde{V}_0 \in \partial \D$. \label{L:stationaryUV}
\end{lemma}

\begin{proof}[Sketch of proof]The idea behind this lemma is simple. Let us write $U_t = e^{i \xi_t} $ and $V_t = e^{i \psi_t}$, where $(\xi_t)_{t \ge 0}$ and $(\psi_t)_{t\ge 0}$ are uniquely defined by continuity. Then it can be checked that, analogous to the Bessel equation \eqref{eqn:slekpdef}, the angle difference $\theta_t = \psi_t - \xi_t$ satisfies 
\begin{equation}
\label{eqn:cotangent}
\dd\theta_t = \frac{\rho +2}{2} \cot(\theta_t) \dd t + \sqrt{\kappa} \dd B_t 
\end{equation}
Recall that $\cot(\theta) \sim 1/\theta$ as $\theta \to 0^+$ so this diffusion looks like a Bessel diffusion of dimension $ 1 + 2( \rho +2) / \kappa >1 $ near zero, and the same is true as $\theta \to 2\pi^-$, with the drift now repelling $\theta_t$ away from $2\pi$. Even when the dimension of the Bessel process is such that these two boundary values are touched by the diffusion, the process $(\theta_t)_{t\ge 0}$ always takes values in $[0, 2\pi]$ and thus has a unique invariant distribution. This is the desired law.
\end{proof}

\begin{definition}\label{D:WPslekp}
	Let $(\bar U_t,\bar V_t)_{t\in \R}$ have the stationary law in \cref{L:stationaryUV} for some $\kappa>0, \rho> -2$. Whole plane SLE$_\kappa(\rho)$ from $0$ to $\infty$ is defined to be the family of whole plane Loewner hulls generated by $( U_t)_{t\in \R}$ via the whole plane Loewner equation, as described in \cref{D:WPSLE}.
\end{definition}

We use $\bar U_t$ and $\bar V_t$ in the definition instead of $U_t$ and $V_t$ even though this does not change the resulting law, but we do so in order to be consistent with Lemma \ref{L:WPRAD}. Indeed, is immediate from the definition and from \cref{L:WPRAD} that given $(U_t,V_t)_{t\in (-\infty,t_0]}$ for fixed $t_0\in \R$, the associated whole plane SLE$_\kappa(\rho)$ from $0$ to $\infty$ from time $t_0$ onwards has the law of a radial SLE$_\kappa(\rho)$ in $\hat{\C}\setminus K_t$, from $g_{t_0}^{-1}(U_{t_0})$ to $\infty$ and with marked point at $g_{t_0}^{-1}(V_{t_0})$. (Its driving function is exactly equal to $(\bar{U}_{t+t_0})_{t\ge 0}$). In particular, for every $\kappa>0$, there almost surely exists a continuous curve $(\gamma(t))_{t\in \R}$ such that $g_{t}^{-1}(\hat{\C}\setminus \bar{\D})=\hat{\C}\setminus K_t$ is the unbounded connected component of $\hat{\C}\setminus \gamma((-\infty,t])$ for each $t$, and as with ordinary whole plane SLE$_\kappa$, it satisfies $\lim_{t\to -\infty} \gamma(t)=0$ and $\lim_{t\to \infty} \gamma(t)=\infty$ (see for example \cite{LawlerWP}).

Whole plane SLE$_\kappa(\rho)$ is also scale invariant: if $a\in \C\setminus \{0\}$ and $\gamma$ is a whole plane SLE$_\kappa(\rho)$ from $0$ to $\infty$, then $a\gamma$ has the same law (modulo time parametrisation) as a whole plane SLE$_\kappa(\rho)$. This again allows us to define whole plane SLE$_\kappa(\rho)$ from $z_1$ to $z_2$ with $z_1\ne z_2\in \hat{\C}$ in a consistent way. 

\begin{definition}\label{D:WPSLErho2}
	Let $z_1,z_2\in \hat{\C}$ and $\rho>-2$. Whole plane SLE$_\kappa(\rho)$ from $z_1$ to $z_2$ is defined to be the image of whole plane SLE$_\kappa$ from $0$ to $\infty$, under a M\"{o}bius transformation sending $0$ to $z_1$ and $\infty$ to $z_2$. The law of this process does not depend on the choice of M\"{o}bius transformation.
\end{definition}

\subsubsection{Whole plane \texorpdfstring{$\mathrm{SLE}_\kappa(\kappa-6)$}{TEXT}}

We end this section on whole plane SLE$_\kappa(\rho)$ with a short discussion about some properties of the curve in the special case $\rho = \kappa -6$. These will be needed in the construction of space-filling SLE$_\kappa$.

\begin{lemma}[Target invariance of SLE$_\kappa(\kappa-6)$] Suppose that $\kappa>4$ and $b_1,b_2\in \C$. Then it is possible to couple a whole plane SLE$_\kappa(\kappa-6)$ curve from $0$ to $b_1$ in $\C$, and from $0$ to $b_2$ in $\C$ so that they coincide until $b_1,b_2$ are contained in separate components of the complement of the curve, and afterwards evolve independently.
\label{L:WPtarget}
\end{lemma}

\begin{proof}
This follows from the target invariance of radial SLE$_\kappa(\kappa - 6)$ (Lemma \ref{L:radialtarget}) and the relationship between whole plane and radial SLE (Lemma \ref{L:WPRAD}). This requires discovering a small (in terms of diameter, say) part of either whole plane curves and taking a limit as the diameter shrinks to zero; the details are left to the reader.
\end{proof}

Let $\k \ge 8$. The next statement shows that the whole plane SLE$_\k(\k - 6)$ from $0$ to $\infty$ does not fill the whole plane (even though, for example, chordal SLE$_\kappa$ in $\H$ does fill the whole of $\H$). Let $K$ denote the hull generated by $\eta$: this is the set of points for which solving the Loewner equation \emph{is not} possible for all times. Equivalently $K = \cup_{t \in \R} K_t$, with $K_t = \C\setminus D_t$, and $D_t$ the unique unbounded component of $\C \setminus \eta(-\infty, t]$. %The next lemma shows that $K \neq \C$ or equivalently $D =\cap_{t \in \R} D_t$ is non-empty. 

\begin{lemma}\label{L:WPnotSF}
Suppose $\kappa \ge 8$. The hull $K$ of a whole plane SLE$_\k(\k - 6)$ curve $\eta$ from 0 to $\infty$ is not all of $\C$. Moreover, $\eta$ is transient: $\eta(t) \to \infty$ as $t \to \infty$. 
\end{lemma}

% (which of course is, by definition, the image of the whole plane SLE$_\k ( \k - 6)$ from 0 to $\infty$ under the M\"obius inversion% $z \mapsto 1/z$). 
Note that this is in contrast with say, chordal SLE$_{\kappa}$ for $\kappa\ge 8$, which eventually swallows every point of the upper half plane. 

\begin{proof} 
To see that $D$ is not empty, suppose that we discover a small chunk $\eta(-\infty, t_0)$ of the whole plane SLE$_\kappa(\kappa - 6)$ from $0$ to $\infty$. The future of this curve is, by Lemma \ref{L:WPRAD} a radial SLE$_\kappa(\kappa-6)$ in the complement of the hull generated by $\eta(-\infty,t_0)$, started at $\eta(t_0)$ and targeted at $\infty$. Its force point is determined by  $V_{t_0}$; more precisely it is given by $z_0 = g_{t_0}^{-1}(V_{t_0})$, where $V_{t_0}$ is as in \cref{D:WPslekp}. By changing coordinates (that is, Lemma  \ref{L:coordinate_change}), we can also view it as a \emph{chordal} SLE$_\kp$ (with no force point) but targeted at $z_0 $ and run until it hits $\infty$ (which is just some interior point of the domain in which this chordal SLE$_\kp$ lives). In particular, the hull generated by the curve $\eta$ does not contain all of $\C$. 

Transience is shown in from \cite{LawlerWP}; alternatively, it follows from the above argument and elementary properties of chordal SLE$_\kp$.
\end{proof}

The following result is in some sense elementary but also very useful conceptually (and also technically, as we will see below). It states that whole plane SLE$_\k(\k-6)$ from $\infty$ to 0 can be viewed as the infinite volume  limit of standard \emph{chordal} SLE$_\k$ in a large domain between two arbitrary boundary points, and stopped when it reaches zero. (We will see later in the chapter that this has a useful implication for \emph{space-filling} SLE$_\kp$: namely, space-filling SLE$_\kp$ is the infinite volume limit of the same curve, \emph{without} stopping it when it reaches zero. See Theorem \ref{T:SFlimitchordal}). 

In order to state this result, we need to discuss the topology for which this convergence holds. This will be the topology of uniform convergence on intervals of the form $[t_0, \infty)$ for every $t_0 \in \R$ (we leave it to the reader to check this defines a metric space, in fact a complete separable metric space, although this is not needed here). In other words, $\eta_n$ converges to $\eta_t$ in this topology if for all $\eps>0$, for all $t_0 \in \R$, there exists $n_0$ such that $|\eta_n(t) - \eta(t)| \le \eps$ for all $n \ge n_0$ and $ t \ge t_0$. Since this is a metrizable topology (and in fact a Polish metrizable one, as mentioned above), it makes sense to talk about convergence in distribution with respect to this topology. 
   
 Let $(\eta(t))_{t\in \R}$ denote a whole plane SLE$_\k$ from $\infty$ to 0, and recall from \eqref{E:capCR} that $\eta$ is parametrised so that $\mathrm{CR}(0,\C\setminus \eta((-\infty,t]))=e^{-t}$ for all $t$.

 \begin{lemma}
 \label{L:WPlimitchordal}
 Suppose that $\k>4$ and let $D_n$ be a sequence of simply connected domains such that $D_n \subset D_{n+1}$ and $\cup_{n \ge 0} D_n = \C$. For each $n$, let $a_n, b_n$ be two prime ends of $D_n$, and let $\eta_n$ denote a chordal SLE$_\k$ in $D_n$ from $a_n$ to $b_n$. Then as $n \to \infty$, the law of $\eta_n$ converges to the law of $\eta$,  a whole plane SLE$_\k$ from $\infty$ to 0, in the sense described above. 
 
In fact, let $\eps>0$ and $t_0 \in \R$ be given. Then for all $n \ge n_0(\eps,t_0)$ large enough, there exists a coupling between $\eta_n$ and $\eta$, and an event $A_n$ of probability at least $1- \eps$ for this coupling, on which $\eta_n(t)=F_n(\eta(t))$ for every $t\ge t_0$, where $F_n$ is a conformal isomorphism defined in a neighbourhood $U$ of $\eta(t_0, \infty)$ satisfying $|F_n(w) - w | \le \eps$ for every $w\in U$.
 \end{lemma}

\begin{proof}
It suffices to prove the  second claim, since this clearly implies the first. 
Let $\eps>0$ and let $t_0 \in \R$. Let $R = C e^{-2t_0} / \eps$, where $C$ will be made precise later (it will in fact be allowed to depend on $\eps$) and define $t_1$ via $e^{-t_1} = R$. Observe that by the change of coordinate formula (Lemma \ref{L:coordinate_change}), until hitting zero $\eta_n$ has the same law as a radial SLE$_\k(\k-6)$ from $a_n$ to 0 with force point at $b_n$ (although since $\rho = \k -6$, by the target invariance property of Lemma \ref{L:radialtarget}, the precise location of this force point will is not relevant except to know that $\eta_n$ does not separate $0$ from $b_n$ until reaching 0). Let $g_n$ denote the conformal isomorphism from $D_n \setminus \eta_n((-\infty, t_1])$ to $\D$ with $g_n(0) = 0$ and $g_n'(0)>0$ and also let $g$ denote the conformal isomorphism from $\C\setminus \eta((-\infty, t_1])$ to $\D$ with $g(0) = 0$ and $g'(0)>0$. Let $\tilde \eta_n = g_n( \eta_n( [t_1, \infty)]$ and let $\tilde \eta = g( \eta( [t_1, \infty))$. Note that both $\tilde \eta_n$ and $\tilde \eta$ are radial Loewner evolutions in $\D$, whose driving functions we denote respectively by $(U^n_{t_1+t})_{t\ge 0}$ and $(U_{t_1 + t})_{t\ge 0}$. Note that $(U_t)_{t\ge t_1}$ has the equilibrium law of Lemma \ref{L:stationaryUV}. Note also that, using the convergence to equilibrium in Lemma \ref{L:stationaryUV}, we can choose $n_0$ large enough so that not only does $D_n$ contain the ball of radius $t_1$ for all $n \ge n_0$, but in fact, for all $n \ge n_0$, we can couple $\eta_n$ and $\eta$ so that $U^n_t = U_t$ for all $t \ge t_1$ on an event of probability at least $1-\eps/2$. Let $A'_n$ denote this event. 

Also choose a constant $k = k(\eps)>0$ large enough so that with probability at least $1- \eps/2$, $\eta(t_0, \infty)$ stays in a ball of radius $ke^{-t_0}$. Let $A_n$ denote the intersection of this event with $A'_n$ (and note that $A_n$ has probability at least $1- \eps$). It remains to show that on $A_n$, $\eta_n([t_0, \infty)$ and $\eta([t_0, \infty))$ are uniformly close.  This will follow from well known distortion estimates, for example from Proposition 3.26 in \cite{Lawlerbook}. Indeed, from this proposition we know that there exists a constant $C = C_{1/2}$ such that for any function $f$ defined on the unit disc which is analytic and one to one with $f(0) = 0$ and $f'(0) =1$, 
\begin{equation}\label{eq:distortion}
|f(z) -z| \le C_{1/2} |z^2|
\end{equation}
for $|z| \le 1/2$. Now consider the map 
$$F_n = g_n^{-1} \circ g:  \C\setminus \eta((-\infty, t_1]) \to D_n \setminus \eta_n ((-\infty, t_1]),$$
and observe that $\eta_n(t_1, \infty)$ is obtained from $\eta(t_1, \infty)$ by mapping it through $F_n$. So it suffices to prove that $F_n$ is close to the identity on the relevant region.
Let $R' = e^{-t_1}/4 = R/4$. By Koebe's quarter theorem, the domain where $F_n$ is defined contains at least $B(0, R')$. The map $z \mapsto F_n(zR')/ R'$ is therefore analytic and one to one on the unit disc, fixing zero and having unit derivative at zero. Hence by  \eqref{eq:distortion}, we deduce that for $r>0$ and $w \in B(0, r)$
$$
|F_n(w) - w | \le  C_{1/2} r^2/R'.
$$
Choosing $r = k e^{-t_0}$ and keeping in mind that $R' = e^{-t_1}/4$ and $e^{-t_1} = R  = C e^{-2t_0} / \eps $, with $C$ to be determined, 
 this means that
\begin{equation}\label{eq:devFn}
|F_n(w) - w | \le C_{1/2} k^2  \eps / C .
\end{equation}
for $w \in B(0,k e^{-t_0})$. We obtain the desired result by taking $C = C_{1/2} k^2$: indeed, on the event $A_n$, for $t \ge t_0$, $\eta_n(t) = F_n( \eta(t))$ and $\eta(t) \in B(0,k e^{-t_0})$ so the use of \eqref{eq:devFn} is justified. Consequently, $|\eta_n(t) - \eta(t) | \le \eps$ for all $t \ge t_0$ on the event $A_n$. 
\end{proof}

\begin{rmk}
\label{R:convLebesgue}
Note that Lemma \ref{L:WPlimitchordal} also holds if the curves are parametrised by Lebesgue area rather than log conformal radius (with respect to time zero). In this case the convergence holds uniformly on compact time intervals (rather than on sets of the form $[t_0, \infty)$). Indeed, the second claim of the lemma shows that the two curves are equal with high probability, up to a small uniform distortion. 

Likewise, the complement of $\eta_n(-\infty, \infty)$ in $D_n$ converges to the complement of $\eta(-\infty, \infty)$ in $\C$ in a very strong sense. For instance, the proof shows that given any neighbourhood $U$ of the origin, we can couple $\eta_n$ and $\eta$ so that with probability arbitrarily close to $1$ as $n\to \infty$, $D_n\setminus \eta_n(-\infty,\infty)$ is the image of $\C\setminus \eta(-\infty,\infty)$ under a conformal isomorphism $F_n$ (defined on a larger domain, including $U$) and is arbitrarily close to the identity on $U$.
\end{rmk}

It will also be useful to describe the boundary of $K=\cup_{t\in \R} K_t$, which is non-empty (by Lemma \ref{L:WPnotSF} in the case $\k\ge 8$ and by the corresponding property of radial SLE$_\k$ in the case $\k\in (4,8)$.)  Since the description of the boundary requires talking about both the value $\kappa$ and the dual parameter $16/\kappa$, we switch to the notation where $\kp > 4$ and $\kappa = 16/\kappa' < 4$. 
In fact, we will only give a description of the boundary in the case $\kp \ge 8$.

 \begin{lemma}
 \label{L:boundaryWP}
Let $\kappa' \ge 8$  and let $\kappa = 16/\kappa'<2$. Let $\eta$ denote a whole plane SLE$_\kp(\kp -6)$ from $\infty$ to 0.
%\footnote{This is, by definition, the image of the whole plane SLE$_\k ( \k - 6)$ from 0 to $\infty$ under the M\"obius inversion $z \mapsto 1/z$.}.
 Then the boundary of $K$ has the same law as $\eta_L(-\infty, \infty)  \cup \eta_R ( -\infty, \infty)$, where $\eta_L,\eta_R$ are defined as follows:
\begin{itemize}
\item $\eta_L$ is a whole plane SLE$_{\kappa}( 2-\kappa)$ from 0 to $\infty$ (note that $\eta_L$ is a simple curve by our assumption that $\kappa' \ge 8 > 6$ and Lemma \ref{L:rhohit})

\item Given $\eta_L$, $\eta_R$ is a chordal SLE$_{\kappa} (  - \kappa/2, - \kappa/2)$ from 0 to $\infty$ in $\C \setminus \eta_L ( -\infty, \infty)$ with force points on either side of the starting point 0 (note that since $\kappa'\ge 8$, $ - \kappa/2 \ge\kappa/2 - 2$ and so by Lemma \ref{L:rhohit}, $\eta_R$ does not hit  any part of $\eta_L$).

\end{itemize} \end{lemma}

\begin{proof}
We expect that such a statement might follow from known duality arguments for chordal SLE$_\kp$ via Lemma \ref{L:WPlimitchordal} (see, for example, \cite{Zhan_duality, Dub3}). However, we could not find such a result in the literature. Nonetheless, this description can be deduced from Theorems 1.4 and 1.6 in \cite{IG4}.
\end{proof}

\begin{rmk}
\label{R:leftrightboundary} Furthermore, given a whole plane SLE$_\kp(\kp -6)$ curve $\eta$ from $\infty$ to 0, with $\kp\ge 8$, it is possible to unambiguously associate to it two curves $\eta_L$ and $\eta_R$, whose union is the boundary of $\eta$, and such that $\eta_L$ lies to the \emph{left} of the curve as we traverse it from $\infty$ to zero, while $\eta_R$ lies to its right as we traverse it from $\infty$ to zero (this is a topological property of curves -- which are oriented by definition -- in two dimensions). The distribution of $(\eta_L, \eta_R)$ is as specified above. Interestingly however, the joint distribution of $(\eta_L, \eta_R)$ is the same as that of $(\eta_R, \eta_L)$. 
\end{rmk}

\subsection{Space-filling SLE in the case \texorpdfstring{$\kappa'\ge 8$}{TEXT}}
\label{SS:spacefillingSLE}

In order to state the mating of trees theorem, we first explain the definition and construction of a space-filling version of SLE in the whole plane (from $\infty$ to $\infty$). We first fix $\kappa' \ge 8$ and  stick with the convention that $\kappa'$ denotes a parameter greater than 4, that will take the value $16/\gamma^2$ when our curves are coupled with $\gamma$ Liouville quantum gravity. The notation $\kappa$ is reserved for the dual parameter $\kappa = 16/\kappa' \in (0, 4)$. In fact, when $\kappa'\ge 8$, the whole plane SLE$_{\kappa'}$, whose definition and properties we have studied in the sections above, already fills the entire hull that it generates. As a result, the construction is much simpler in this case than when $\kp\in (4,8)$. %A reader who is interested in understanding the spirit of the mating of trees theorem can keep this case in mind and skip many of the technical details of this section, which correspond to $\kp\in (4,8)$. 
We note that on the LQG side (when we eventually couple our space-filling curve with $\gamma$ LQG), choosing $\kappa'\ge 8$ amounts to restricting $\gamma$ to the interval $(0, \sqrt{2}]$ (which essentially corresponds to the $L^2$ phase of GMC). Unfortunately it is the interval $\gamma \in [\sqrt{2}, 2)$ which is believed to correspond to scaling limits of random planar maps weighted by the self dual FK percolation model described in Chapter \ref{S:maps}.

%%%%
% 1. \kappa' >8
% 2. The two points of views: ordering vs. exploration
% 3. description of ordering: branching tree, merging of boundaries. This gives an order which (thm) is.a space-filling curve.  
% 4. dynamical description of the path using space-filling loops and chords.   
% 
%

\subsubsection{Definition of space-filling \texorpdfstring{$\mathrm{SLE}_\kp$ ($\kp \ge 8$)}{TEXT} }

\ind{Space-filling SLE!Definition, $\kp\ge 8$}
Let $\kappa'\ge 8$. The whole plane SLE$_{\kappa'}$ defined in Section \ref{A:wpSLE} is a curve $(\tilde \eta_t)_{t\in (-\infty, \infty)}$ which ``starts'' at zero (meaning $\lim_{t\to -\infty} \tilde \eta_t = 0$) and is targeted at infinity. For the mating of trees theorem, however, we will need to define  a curve from $\infty$ to $\infty$, which visits zero at time 0; this will make it possible for the ``past'' and ``future'' of the curve with respect to 0 to play symmetric roles, which turns out to be an important feature of the theory.

We therefore cannot directly use $\tilde \eta$ as our space-filling SLE. Instead we proceed in two steps. Let $\eta^-$ denote a whole plane SLE$_\kp (\kp-6)$ from $\infty$ to 0, as defined in Definition \ref{D:WPSLErho2}. Let $K^- = \eta^-((-\infty,\infty))$, and let $D^- = \C \setminus K^-$. We will use $K^-$ as the ``past'' of time zero, and the closure of $D^-$ will be the future. The following property of $\eta^-$ will motivate the definition of the space-filling curve coming below.   
%It follows from the relationship with radial SLE$_\kappa(\kappa-6)$ (see below Definition \ref{D:WPSLErho2}) that both $D^-$ and $D^+$ are almost surely simply connected domains. In fact, let us note the following special property of $\eta^-$.

\begin{lemma}\label{L:WPstat}
	Let $\kp\ge 8$ and let $\tau$ be any almost surely finite stopping time for $\eta^-$ (with respect to the filtration generated by the curve $\eta^-$ itself, parametrised so that the conformal radius of 0 in the complement of the curve is $e^{-t}$)\footnote{Or equivalently, the filtration generated by the pair $(U_t, V_t)_{t \in \R}$ of \cref{D:WPslekp} after applying a M\"obius inversion.}. Then the complement, $D^-_\tau$, of $\eta^{-}(-\infty,\tau]$ in $\C$, is an unbounded simply connected set with probability one. Moreover, given $\eta^{-}(s), s\le \tau]$, the law of $\eta^{-}|_{[\tau,\infty)}$ is that of a {radial} SLE$_\kp$ in $D^-_\tau$ from $\eta^{-}(\tau)$ to $0$, parametrised by minus log conformal radius seen from $0$.
\end{lemma}

\begin{proof} 
Let $\eta(t)$ denote a whole plane SLE$_\kp (\kp -6)$ from 0 to $\infty$ (thus the laws of $\eta^-$ and $\eta$ are related to each other by M\"obius inversion).  We may assume without loss of generality that $\tau$ is a stopping time for $\eta$. Let $D_\tau$ denote the complement of $\eta( - \infty, \tau]$. To prove the lemma it suffices to show that: (1) $D_\tau$ is simply connected; (2) contains points arbitrarily close to zero; and (3) given $\eta(-\infty, \tau]$, the rest of the curve $\eta$ is distributed as a {radial} SLE$_\kp$ in $D_\tau$ from $\eta(\tau)$ to 0,  parametrised by logarithmic capacity (seen from infinity). By changing coordinates (that is, \cref{L:coordinate_change}), (3) is equivalent to saying that the conditional law of $\eta([\tau,\infty))$ is that of a {chordal} SLE$_\kp(\kp-6)$ in $D_\tau$ from $\eta(\tau)$ to $\infty$, with force point at $0$.
 
Let $(\bar U_t, \bar V_t)_{t\in \R}$ be the stationary (radial) process of Lemma \ref{L:stationaryUV} defining the whole plane curve $\eta$. For $t\in \R$, let $K_{t}$ be the hull generated by $\eta((-\infty,t])$, and let $g_{t}$ be the unique conformal isomorphism from $\hat{\C}\setminus K_{t}\to \hat{\C}\setminus \bar{\D}$ with $g_{t}(\infty)=\infty$ and $g_{t}'(\infty)=e^{-t}$ (that is, $g_{t}(z) = z e^{-{t}} + O(1)$ as $z \to \infty$). 
Then from the strong Markov property of $(U,V)$ and the relationship between whole plane and radial SLE$_\kp(\rho)$ (specifically the discussion just below Definition \ref{D:WPSLErho2} in the appendix) we learn given $\eta(-\infty, \sigma]$ for any $\eta$-stopping time $\sigma\in \R$, the remainder of $\eta$ is a radial SLE$_\kp ( \kp -6)$ in $\hat{\mathbb{C}}\setminus K_\sigma$, targeted at $\infty$, and with force point located at $g_\sigma^{-1}(V_\sigma)$. (Equivalently by Lemma \ref{L:coordinate_change}, it is a chordal SLE$_\kp$ targeted at $g_\sigma^{-1}(V_\sigma)$). Since $\kappa' \ge 8$, this means that the curve $\eta$ is in the ``space-filling phase''. In particular,  $D_\tau=\mathbb{C}\setminus K_\tau$ almost surely, and by properties of the Loewner evolution, $D_\tau$ is almost surely simply connected.

To see point (2) -- that $D_\tau$ contains points arbitrarily close to zero -- we simply observe that $D_\tau \supset D_\infty$ which itself satisfies this property by Lemma \ref{L:boundaryWP}. Indeed, the boundary of $\eta$ is given by an explicit pair of SLE curves and therefore the complement of $\eta$, that is $D_\infty$, contains points arbitrarily close to zero as desired. In fact this argument shows that $D_\infty$ is a Jordan domain for which all boundary points (including $0$) correspond to a unique prime end in the language of  \cite{Pommerenke}. As a consequence, note  that 0 corresponds to a unique prime 5 in $D_\tau$ as well, and not just in $D_\infty$, which is a consequence of transience and a zero-one argument left to the reader. This will be required below. 

Finally, we are left to show (3), which by the discussion above (with $\sigma=\tau)$), boils down to proving that $z_\tau = g_\tau^{-1}(V_\tau) = 0$.

To see this, fix a sequence $t_n \to - \infty $. As discussed above, given $(U_t,V_t)$ for $t\le t_n$, the conditional law of $\eta$ after time $t_n$ is that of a chordal SLE$_\kp$ in $D_{t_n}$, from $\eta(t_n)$ to $z_{t_n}$, reparametrised according to log capacity seen from $\infty$. In particular, it will not hit $z_{t_n}$ again after time $t_n$, which means that $z_t$ stays constant after time $t_n$.  Consequently, we have that $z_\tau = z_{t_n}$ almost surely. Since $n \ge 1$ was arbitrary, and $z_{t_n}$ 
lies on the boundary of $\eta((-\infty,t_n])$ (a set of diameter tending deterministically to 0), it follows that $z_{t_n} \to 0$ and thus $z_\tau = 0$ as desired.

Note also that this argument implies that $U_t$ and $V_t$ are determined by $\eta(-\infty, t)$, since $U_t$ is the driving of the Loewner evolution (explicitly, $U_t = g_t ( \eta(t))$) and $V_t = g_t(0)$. Therefore the filtrations generated by $(U, V)$ and by the curve are indeed equal, as claimed in the Lemma.
\end{proof}

In particular, this makes the following definition possible (and natural).

\begin{definition}\label{D:sfSLE}
Let $\kp\ge 8$. Given a whole plane SLE$_\kp ( \kp - 6)$ from $\infty$ to 0 which we denote by $\eta^-$, let $\eta^+$ denote a (conditionally independent) chordal SLE$_\kp$ in $\C \setminus \eta^- ( \R)$ from 0 to $\infty$. By definition, the whole plane \textbf{space-filling SLE$_\kp$} from $\infty$ to $\infty$, is the curve $\eta$ obtained by concatenating $\eta^-$ and $\eta^+$, and then reparametrising time so that $\eta(0) = 0$ and $\Leb (\eta([0,t]) ) = |t|$, that is,  so that $\eta$ is parametrised by its (Lebesgue) area. 
\end{definition}

Indeed, it can be checked that both $\eta^-$ and $\eta^+$ cover an area of positive Lebesgue measure in any finite-time interval, and that this area is in fact a continuous function of the length of the interval (this is well known for $\eta^+$ by properties of chordal SLE$_\kp$, and for $\eta^-$ it can be deduced from Lemma \ref{L:WPstat}.) This means that such a continuous reparametrisation by Lebesgue area is indeed possible.

Given \cref{L:WPstat}, it is not surprising that the space-filling SLE curve we have just defined is stationary in a strong sense; however, there are some subtleties in justifying this because of the way the curve is parametrised (since we know at time 0 it must visit 0, and visit exactly an area of size $t$ in any interval of length $t$). A precise statement of this sort will be given (but not proved) in Lemma \ref{L:conestationary} a bit later on.

%\medskip It is useful to define a notion of stopping time which does not depend on the parametrisation, such as the first time the curve $\eta$ visits a given region. Let us say that a time $\tau$ is a \textbf{stopping rule} if it remains a stopping time after an arbitrary change of time-parametrisation curve: that is, for any continuous increasing homeomorphism $f$ from $\R$ to $\R$, $f( \tau)$ is a stopping time of $(\eta(f(t)))_{t\in \R}$ in its natural filtration. \ind{Stopping rule}

\ind{Space-filling SLE!Markov property}
For now, we formulate a useful Markov property. 
\begin{lemma}\label{L:SF_markov}
Let $\kp\ge 8$, let $\eta$ be a space-filling SLE$_\kp$ from $\infty$ to $\infty$, let $U$ be a non-empty bounded subset of $\C$, and let $\tau$ be the first time that $\eta$ enters $U$.  Then conditionally on $(\eta(t), t \le \tau)$, the rest of the curve $(\eta(\tau+ t))_{t\ge 0}$ is, up to a change of time parametrisation, a chordal SLE$_\kp$ in $\hat\C \setminus \eta((-\infty, \tau])$ from $\eta(\tau) $ to $\infty$. 
\end{lemma}

\begin{proof}
	Let $g:\hat\C\setminus \eta((-\infty,\tau])\to \hat\H$ be the unique conformal isomorphism with $g(\eta(\tau))=0$ and $g(z)/z\to 1$ as $z\to \infty$. Let $\tilde{\eta}$ be the image of $(\eta(\tau+t))_{t\ge 0}$ under $g$, reparametrised by half plane capacity (that is, so that the infinite connected component of $\H\setminus \tilde{\eta}([0,t])$ has half plane capacity $2t$ for $t\ge 0$). The lemma is equivalent to the fact that, conditionally on $(\eta(t), t\le \tau)$ $\tilde{\eta}$ has the law of a chordal SLE$_\kp$ in $\H$ from $0$ to $\infty$.
	
	There are two events to consider. On the event that $0\in \eta((-\infty,\tau))$, let $\tau_0:=\inf\{t: 0\in \eta((-\infty,t])$, so that $\tau_0\le \tau$ and $U\subset \hat\C\setminus \eta((-\infty,\tau_0])$. Then $(\eta(t))_{t\ge \tau_0}$ is by definition a chordal SLE$_\kp$ in $\hat\C\setminus \eta((-\infty,\tau_0])$ from $\eta(\tau_0)$ to $\infty$, reparametrised by Lebesgue area. $\tau$ is simply the first time that this curve enters $U$, and the Markov property of chordal SLE$_\kp$ implies the desired statement in this case.
	
	On the event that $0\notin \eta((-\infty,\tau])$, write $\eta'$ for $\eta$ but in its usual whole plane Loewner evolution parametrisation, and $\tau'$ for the firs time it enters $U$. Notice that the sigma-fields generated by $(\eta'(t),t\le \tau')$ and $(\eta(t), t\le \tau)$ are the same. Therefore, by Lemma \ref{L:WPstat}, and after reparameterising by half plane capacity, $\tilde{\eta}$ has the law of a chordal SLE$_\kp$ targeted at $\infty$ up until the first time it hits $g(0)$. After this time, by definition, it has the law of a chordal SLE$_\kappa'$ in the remaining domain, targeted at $\infty$. But this two step description gives exactly the law of a chordal SLE$_\kappa$ from $0$ to $\infty$ in $\H$. Thus $\tilde{\eta}$ has this law, as required.
\end{proof}

One consequence is that any fixed point $z\in \C$ is almost surely not a double point of the space-filling SLE$_\kp$ curve $\eta$, since this is true of chordal SLE$_\kp$ (by the Markov property of chordal SLE$_\kp$ and properties of Bessel processes). 

\subsubsection{Space-filling SLE as an infinite volume limit of chordal SLE \texorpdfstring{$(\kp\ge 8)$}{TEXT} }

The following description of space-filling SLE is extremely useful for the intuition: it says that we can view space-filling SLE$_\kp$ as the infinite volume limit of standard, chordal SLE$_\kp$ in a domain $D_n$ tending to infinity between two arbitrary prime ends of $D_n$. This point of view is sometimes taken as a definition of space-filling SLE$_\kp$, although we were not able to find a reference for the existence of such a limit in the literature.

 \begin{theorem}
 \label{T:SFlimitchordal}
 Let $D_n$ be a sequence of simply connected domains such that $0 \in D_n \subset D_{n+1}$ and $\cup_{n \ge 0} D_n = \C$. For each $n$, let $a_n, b_n$ be two prime ends of $D$ ,and let $\eta_n$ denote a chordal SLE$_\kp$ in $D_n$ from $a_n$ to $b_n$, parametrised by Lebesgue area with $\eta_n(0) = 0$. Then as $n \to \infty$, the law of $\eta_n$ converges to the law of $\eta$, a space-filling SLE$_\kp$ from $\infty$ to $\infty$, for the topology of uniform convergence on compact intervals of time. 
 \end{theorem}

\begin{proof}
The proof of the theorem follows almost directly from Lemma \ref{L:WPlimitchordal} (see also Remark \ref{R:convLebesgue}). Indeed let $\eta_n, \eta$ be as in the theorem, and let $K_n^- = \eta_n(( - \infty, 0])$ (resp $K^- = \eta((-\infty, 0 ])$. Lemma \ref{L:WPlimitchordal} shows that $\eta_n|_{(-\infty ,0]}$ converges weakly to $\eta|_{(-\infty, 0]}$, uniformly on compact time intervals. Furthermore, given $K_n^-$, $\eta_n(0, \infty)$ is a chordal SLE$_\kp$ in $D_n \setminus K_n^-$, while given $K^-$, $\eta(0,\infty)$ is a chordal SLE$_\kp$ in $\C\setminus K^-$. 
Moreover, by Remark \ref{R:convLebesgue}, we can couple $D_n\setminus K_n^-$ and $\C\setminus K^-$ so that for any fixed neighbourhood $U$ of the origin, with probability arbitrarily close to $1$ as $n\to \infty$,  $D_n \setminus K_n^-$ is the image of $\C \setminus K^-$ under a conformal isomorphism $F_n$ (defined on a larger domain, including $U$), that is arbitrarily close to the identity on $U$. %from the statement of Lemma \ref{L:WPlimitchordal} to a chordal SLE$_\kp$ in $\C\setminus K^-$. 
This immediately implies the desired convergence.    
\end{proof}

\begin{rmk}[Reversibility of (whole plane) space-filling SLE$_\kp$] \ind{Space-filling SLE!Reversibility}
Although we will not need it can be checked that this theorem implies the reversibility of whole plane, space-filling SLE$_\kp$ (this is the only kind of space-filling SLE discussed in this book). This is not entirely straightforward because chordal SLE$_\kp$ is \emph{not} exactly reversible when $\kp \ge 8$. Instead, let us sketch the argument here (we emphasise the rest of the arguments in this chapter do not depend on this reversibility). The time reversal of an SLE$_\kp$ from $a_n$ to $b_n$ is a chordal SLE$_\kp (\rho, \rho)$ with $\rho = \kp/2 - 4$ and the two force points located on either side of $b_n$ (\cite{IG4}). However, as $n \to \infty$, the effect of these force points vanishes when we concentrate on a bounded window around zero. Indeed, even though the location of the force points changes whenever the chord swallows a force point, these remain constantly on the boundary of $D_n$ and thus uniformly far away from the bounded window.
\end{rmk}

\subsubsection{Alternative construction from a branching SLE \texorpdfstring{$(\kp\ge 8)$}{TEXT} )}

\label{SSS:branching}

Let $\cQ = \{ z_i\}_{ i \ge 1}$ denote a countable dense set in $\C$.  It is not hard to see that space-filling path $\eta$ that we have just defined almost surely induces an order on $\cQ$: indeed let us say that \begin{equation}\label{E:preceq_eta} z_i \preceq_\eta z_j\end{equation} if and only if $\eta$ visits $z_i$ before $z_j$. This is almost surely an order, since if $ z_i \preceq_\eta z_j$ and $z_j \preceq_\eta z_i$ then either $z_i = z_j$ or $z_i$ is a double point of $\eta$, where the latter event has probability zero (simultaneously for all $i$) by Lemma \ref{L:SF_markov}.
  
Let us suppose that $z_0 = 0$, and call the \textbf{past of 0} the set $K_\cQ^-(0) = \{ z_i :z_i \preceq 0\}$. Likewise let us call the \textbf{future of 0} the set $K^+_\cQ(0) = \{ z_i : 0 \preceq_\eta z_i\}$. Both these sets can be described directly using the whole plane SLE$_\kp ( \kp - 6)$ curve $\eta^-$ from $\infty$ to 0: namely, $K_\cQ^-(0) = K^- \cap \cQ$ with $K^-$ the hull of $\eta^-$, and $K_{\cQ}^+(0)=\cQ\setminus K_{\cQ}^-(0)$.
 
We will now give an equivalent description of the (law of the) ordering $\preceq_\eta $ on $\cQ$ defined in \eqref{E:preceq_eta}, in terms of what is known as \textbf{branching SLE}$_\kp ( \kp - 6)$. Conversely, this gives us an alternative (implicit) description of the  law of the space-filling curve $\eta$ in terms of such branching SLE$_\kp ( \kp - 6)$, which provides a useful alternative point of view. 

\paragraph{Branching SLE$_\kp(\kp-6)$.}
\ind{Branching SLE$_{\kp}(\kp -6)$}
We first give a definition, valid for every $\kappa' >4$, of the branching SLE$_\kp (\kp - 6)$ (branching SLE$_\kp(\rho)$ only makes sense in the case when the weight $\rho$ of the force point is equal to $\kp - 6$, since target invariance is a key part of the definition). Recall that by Lemma \ref{L:WPtarget}, given two points $z$ and $w $ in $\C$, it is possible to couple a whole plane SLE$_\kp (\kp -6)$ from $\infty$ to $z$ and $w$ respectively, in such a way that the two curves coincide (up to reparametrisation) up until $z$ and $w$ are separated from one another by the curve, after which the evolution of the two curves is independent. This coupling can immediately be extended to the dense countable set $\cQ$: that is, for each point $z_i \in \cQ$, we have a whole plane SLE$_\kp (\kp -6)$ curve $\eta_{z_i}$ from $\infty$ to $z_i$, and the joint law of $\eta_{z_i}$ and $\eta_{z_j}$ is as described above for all pairs $i,j$.

A concrete inductive construction when $\kp\ge 8$ goes as follows. Start with a whole plane SLE$_\kp ( \kp - 6)$ from $\infty$ to $z_1$, and call it 
$\eta_{z_1}$. 
Now consider $z_2 $. If $\eta_{z_1}$ visits $z_2$ (at time $\tau_{z_2}$, say) then we define $\eta_{z_2}$ to be $\eta_{z_1} ((-\infty , \tau_{z_2}])$, up to reparametrisation. Otherwise, we run an independent radial SLE$_\kp (\kp -6)$ in $\C \setminus \eta_{z_1} ( \R)$ from $z_1$ to $z_2$, with force point at $\infty$, and call the concatenation of $\eta_{z_1}$ and this additional curve. Now we proceed inductively as follows. Suppose that $\eta_{z_1}, \ldots, \eta_{z_n}$ have been constructed and that for each $1\le i \neq j \le n$, either $\eta_{z_i}$ is a subcurve of $\eta_{z_j}$ or the other way around; let $\eta_{z_m}$ denote the maximal curve. We construct $\eta_{z_{n+1}}$ as follows. If $z_{n+1}$ is visited by $\eta_{z_m}$, at time $\tau_{z_{n+1}}$ say, then $\eta_{z_{n+1}} = \eta_{z_m} (( -\infty, \tau_{z_{n+1}}])$ (up to reparametrisation). If not, then we append to $\eta_{z_m}$ an independent radial SLE$_\kp (\kp - 6)$ in $\C \setminus \eta_{z_m} ( \R)$ from $z_m $ to $z_{n+1}$ with force point at $\infty$. The validity of this construction is justified simply by the strong Markov property of whole plane SLE$_\kp (\rho)$ (see the discussion above \cref{D:WPSLErho2}) and the target invariance of \cref{L:WPtarget}.

\paragraph{Ordering from a branching SLE$_\kp(\kp-6)$ when $\kp\ge 8$.} We now return to the ordering of $\cQ$ associated with a branching SLE$_\kp ( \kp - 6)$, and assume that $\kp\ge 8$. Let $\{ \eta_{z_i}\}_{i \ge 1} $ be a branching SLE$_\kp ( \kp - 6)$ %(here we use a - superscript to indicate that these branches are only going to define the ``past'' of each $z_i$). 
and for each $i \ge 1$, let $K^-(z_i) $ denote the hull of $\eta_{z_i}$, and let $K^-_\cQ(z_i) = K^-(z_i) \cap \cQ$. We can use this to define an order on $\cQ$ almost surely: we say that $z_i \preceq_b z_j$ (the $b$ stands for branching) if $z_i \in K^- (z_j)$. It is not hard to see that this is indeed almost surely an order on $\cQ$: for instance, to check transitivity, one simply notes that since $\kp \ge 8$, $z_i$ becomes separated from $z_j$ if and only if $\eta_{z_j}$ hits $z_i$ on its way to $z_j$, almost surely. Transitivity follows immediately, as does antisymmetry. 

We can now verify that the two orders $\preceq_\eta$ and $\preceq_b$ on $\cQ$ coincide in law. 

\begin{lemma}
\label{L:2orders} Let $\kp\ge 8$. There is a coupling of a space-filling SLE$_\kp$, $\eta$, and a branching SLE$_\kp ( \kp -6)$, $(\eta_{z_i})_{z_i\in \cQ}$ such that $\eta_{z_i}=\eta((-\infty,\tau_{z_i}])$ for each $z_i\in \cQ$, where $\tau_{z_i}$ is the first time that $\eta$ visits $z_i$. In particular, in this coupling, $z_i \preceq_\eta z_j$ if and only if $z_i \preceq_b z_j$. 
\end{lemma}

\begin{proof}
Indeed if $\eta$ is a space-filling SLE$_\kp$, and $\tau_{z_i}$ is the first time that $\eta$ visits $z_i$, then the collection $\eta_{z_i}: = \eta((-\infty, \tau_{z_i}])$ ($z_i\in \cQ$) has the law of a branching SLE$_\kp (\kp - 6)$, up to reparametrisation of the curves. This follows from the inductive construction defining the branching SLE$_\kp(\kp - 6)$ on the one hand, and the Markov property of space-filling SLE$_\kp$ proved in Lemma \ref{L:SF_markov}. 
\end{proof}

\subsubsection{Imaginary geometry ordering; continuum trees \texorpdfstring{$(\kp\ge 8)$}{TEXT} }
\label{SSS:treeskp8}

There is another, perhaps slightly more geometric, description of the ordering defined by the branching SLE$_\kp ( \kp - 6)$ which can be phrased simply in terms of the left and right boundaries of each branch $\eta_{z_i}$, as defined in Lemma \ref{L:boundaryWP}. Recall from this lemma that if $z \in \C$ is fixed and $\eta_z$ is a whole plane SLE$_\kp(\kp - 6)$ from $\infty$ to $z$, then the boundary of $\eta$ has the law of the union of two curves $\eta^L_z$ and $\eta^R_z$, where $\eta^L_z$ has the law of an SLE$_\k ( 2- \k)$ and, given $\eta^L_z$,  $\eta^R_z$ has the law of a chordal SLE$_\k ( - \k/2, - \k/2)$ in the complement of $\eta_z^L$ from $0 $ to $\infty$ and with force points on either side of zero. (Recall that $\k = 16/\kp$ and that in the case of the whole plane curve $\eta^L_z$, we don't need to specify the location of the force point of weight $2-\kp$, which is in some sense the same as the starting point, that is, $\infty$ in this case). In fact, these two curves $\eta^L_z$ and $\eta^R_z$ are determined unambiguously by $\eta_z$, see \cref{R:leftrightboundary}.

\begin{rmk}\label{R:IG}
{It can be checked that the joint laws of the curves $\{\eta_z^L\}_{z\in \cQ}$ and $\{ \eta_z^R\}_{z \in \cQ}$ are identical to what Miller and Sheffield, \cite{IG4}, call the family of ``flow lines''  of a Gaussian free field $h$ in the whole plane with respective angles $-\pi/2$ and $\pi /2$. This will not be needed in the following but, together with the discussion below, it explains why our definition of space-filling SLE$_\kp$ coincides with that given in \cite{IG4} and \cite{DuplantierMillerSheffield}.}
\end{rmk}

\begin{figure}
	\centering
 \includegraphics[width=0.5\textwidth]{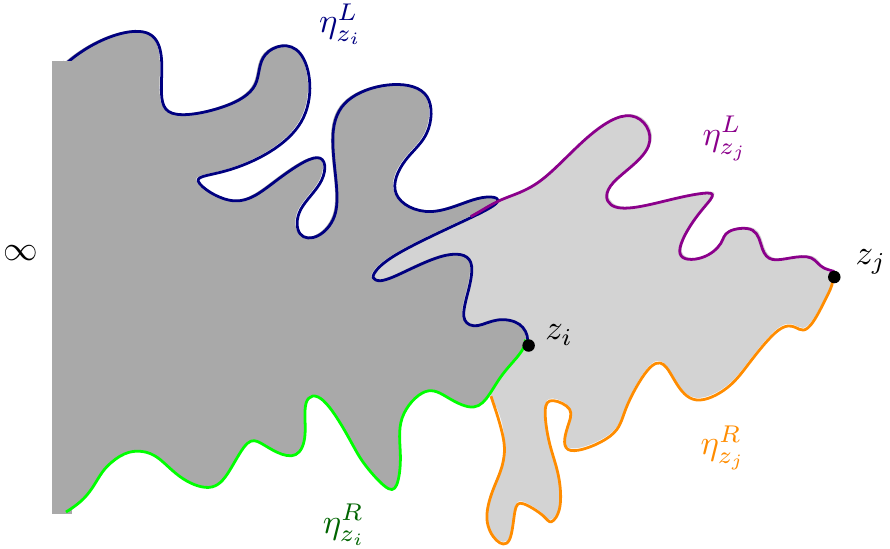}	
\caption{Illustration of a branching whole plane SLE$_\kp(\kp-6)$ $\{\eta_{z_i}^-\}_{z_i\in \cQ}$. The range of $\eta^-_{z_i}$ (the whole plane SLE$_\kp(\kp-6)$ branch from $\infty$ to $z_i$) is  shaded dark grey, and the range of $\eta_{z_j}$ (the whole plane SLE$_\kp(\kp-6)$ branch from $\infty$ to $z_j$) is shaded light grey. Their left and right outer boundaries are coloured in purple/blue and orange/green respectively. In this situation $z_i\preceq z_j$, since $\eta_{z_i}^L$ merges with $\eta_{z_j}^L$ from the left (equivalently, $\eta_{z_i}^R$ merges with $\eta_{z_j}^R$ from the right).\label{F:merging}}
\end{figure}

Continuing with this geometric definition, take two points $z_i, z_j \in \cQ$, and consider their associated left boundary paths (say). That is, the two curves $\eta^L_{z_i}$ and $\eta^L_{z_j}$, where we now view them as starting from $z_i$ and $z_j$ respectively and targeted at $\infty$. One can see from the inductive construction of the branching SLE$_\kp ( \kp - 6)$ that these two paths necessarily merge eventually. We will see that the way these paths merge actually determines the ordering between $z_i$ and $z_j$. 
Indeed, let us say that 
\begin{equation}\label{Eq:orderboundary}
z_i \preceq_{\mathrm{IG}} z_j \text{ iff  $\eta_{z_i}^L$ merges with $\eta^L_{z_j}$ from the \emph{left}.}  
\end{equation}
See Figure \ref{F:merging}. Equivalently we can  use the right boundaries; in this case $z_i \preceq_{\mathrm{IG}} z_j$ if and only if $\eta_{z_i}^R$ merges with $\eta_{z_j}^R$ from the right. We refer to this ordering as the \textbf{Imaginary Geometry ordering}. 
It is not hard to check that (for topological reasons) 
\begin{equation}\label{Eq:IGb}
z_i \preceq_{\mathrm{IG}} z_j \text{ iff } z_i \preceq_b z_j.
\end{equation}
\ind{Imaginary geometry}

By Lemma \ref{L:2orders}, the space-filling curve $\eta$ can be uniquely recovered from this ordering: it is the unique curve (up to reparametrisation) which traverses the points $z_i$ in an order compatible with \eqref{Eq:orderboundary}. This is the definition used in \cite{DuplantierMillerSheffield}, but we stress that it is not obvious at all why such a (continuous) curve should exist at all (in \cite{DuplantierMillerSheffield} the existence of the curve is imported from the theory of imaginary geometry and in particular \cite{IG4}).%, so the construction we have given here is slightly different although of course equivalent, see \cref{R:IG}). 
\medskip{}

Let us conclude this subsection by describing another, more heuristic, way to think of how the boundary curves $\{\eta^L_z\}_{z \in \cQ}$ (or $\{\eta^R_z\}_{z \in \cQ}$) determine the space-filling curve $\eta$.

\paragraph{Continuum trees.} Observe that the merging property of the curves $\{\eta^L_z\}_{z \in \cQ}$ described above, means that they define a topological tree $\cT^L$ embedded in the plane. This tree is simply the union of all the paths $\{\eta^L_z\}_{z \in \cQ}$ (by topological tree $\cT$, we simply mean that for every pair of points $z, w \in \cT$ there is a unique simple continuous path going from $z$ to $w$ up to reparametrisation). Likewise the right boundary curves $\{ \eta^R_{z}\}_{z \in \cQ}$ define a topological tree $\cT^R$ embedded in the plane. These two trees are dual to one another in the sense that, for instance, a curve in $\cT^L$ cannot cross another curve in $\cT^R$.

 These two trees $\cT^L$ and $\cT^R$ can be thought of as the continuum analogues, and indeed should be the scaling limits, of the two canonical trees arising from Sheffield's bijection described in Chapter \ref{S:maps} for random planar maps weighted by the self dual Fortuin--Kasteleyn percolation model. The space-filling SLE$_\kp$ defined above can then be thought of as the Peano curve ``snaking'' in between these two trees.

\subsubsection{Summary of the constructions for \texorpdfstring{$\kp\ge 8$}{TEXT} }
\label{SSS:summary}

We have now seen several equivalent viewpoints of (whole plane) space-filling SLE$_\kp$, which can be used as alternative equivalent definitions depending on the properties one cares about. As these points of views are quite different from one another, we summarise what we have just done with the following table.

\begin{center}
\begin{tabular}{| l |l|l|l|}
\hline
Direct construction                                   & branching ordering                                  & Imaginary Geometry ordering \\
\hline
\hline
Whole plane SLE$_\kp$,                        & $z_i\preceq_b z_j$ iff           & $z_i \preceq_{m} z_j$ iff   \\
followed by chordal SLE$_\kp$ &  $z_i \in \eta^-_{z_j}$             & $\eta_{z_i}^L$ merges with $\eta^L_{z_j}$ from \emph{left}. \\
\hline
Definition \ref{D:sfSLE}       & Lemma \ref{L:2orders}  & \eqref{Eq:IGb}\\
\hline
\end{tabular}
 \end{center}
 
 We also recall that we proved that space-filling SLE$_\kp$ can be obtained as the infinite volume limit of chordal SLE$_\kp$ in large domains from one arbitrary boundary point (prime end) to another, see Theorem \ref{T:SFlimitchordal}. This too could be used as a definition.

\medskip  It will be useful to contrast these definitions with the definitions we will give in the case where $\kp \in (4,8)$, which is much more delicate (and about which we will consequently prove less in this book).  
%n fact, the above lemma means that for given D,a,c, SLEκ(κ−6) can be simultaneously defined towards a countable dense set of target points in D, in such a way that the above description holds for any two given target points. The object created in this manner is referred to as an SLEκ(κ − 6) branching tree, or sometimes just a branching SLEκ.

\subsection{ Space-filling SLE for \texorpdfstring{$\kp\in(4,8)$}{TEXT} }
\label{SS:sf48}

As hinted above, the definition of space-filling SLE$_\kp$ is considerably easier when $\kp \ge 8$ than when $\kp \in (4,8)$: the reason for this is that chordal SLE$_\kp$ is already space-filling, so that the definition requires little modifications. By contrast when $\kp \in (4,8)$ a chordal SLE$_\kp$ is self touching but not space-filling, and the definition of the space-filling versions (chordal or whole plane) of SLE$_\kp$ requires sophisticated tools (note that the space-filling version of SLE does \emph{not} exist for $\kappa < 4$, and the case $\kappa =4$ is very delicate and will be partly discussed later on). In the case when $\kappa \in (4,8)$ let us say right from the start, to help the intuition, that the space-filling version of SLE$_\kp$ is believed to be the scaling limit of the discrete space-filling paths associated to decorated planar maps defined in Chapter \ref{S:maps}.

\ind{Imaginary geometry}
\medskip To this day the only tools which have been developed to define space-filling SLE$_\kp$ in the case where $\kp \in (4,8)$ are those coming from the theory of \emph{Imaginary Geometry} developed in \cite{MSIG1,IG2, IG3} and especially \cite{IG4}. In order to avoid going into such technical details, we have opted for a presentation of those aspects of the definition which do not rely on imaginary geometry, and can be understood without familiarity or knowledge of this theory. The downside of this approach, however, is that the proof of the theorem defining space-filling SLE$_\kp$ as a continuous curve will not be included in this book. 

\dis{\emph{Disclaimer:} we warn the reader that throughout Section \ref{SS:sf48} we intend to provide some explanations which we believe to be useful, but \textbf{these should \emph{not} be considered fully rigorous proofs}. }

\medskip Instead of providing a direct construction of space-filling SLE$_\kp$ for $\kp \in (4,8)$, we will define an ordering of a dense set of points in $\C$, in the manner of columns 2 and 3 in the table of Section \ref{SSS:summary}, and we will check that these two orderings are consistent with one another. However, we will not verify that this ordering is associated with a (unique, up to translation of time) continuous curve in the sense that points are traversed by the curve in the order specified above. We will, however, make a precise statement and give references for the proof. 

\begin{figure}
	\begin{center}
		\includegraphics[width=.4\textwidth]{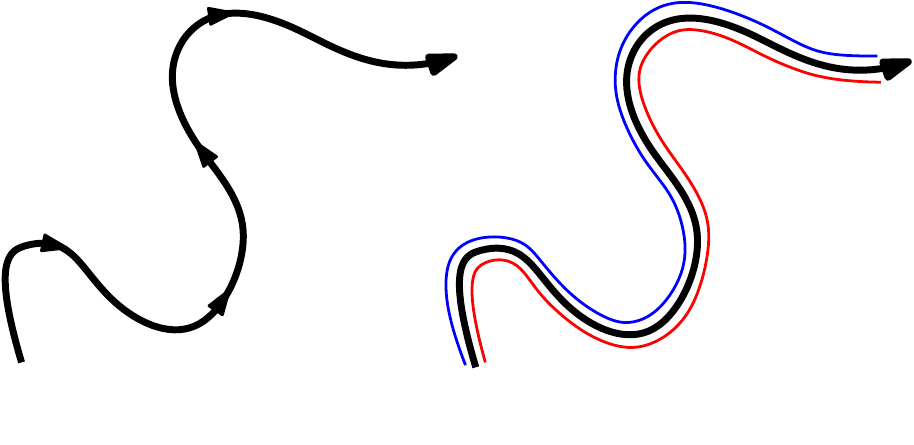} \quad \quad 
		\includegraphics[width=.4\textwidth]{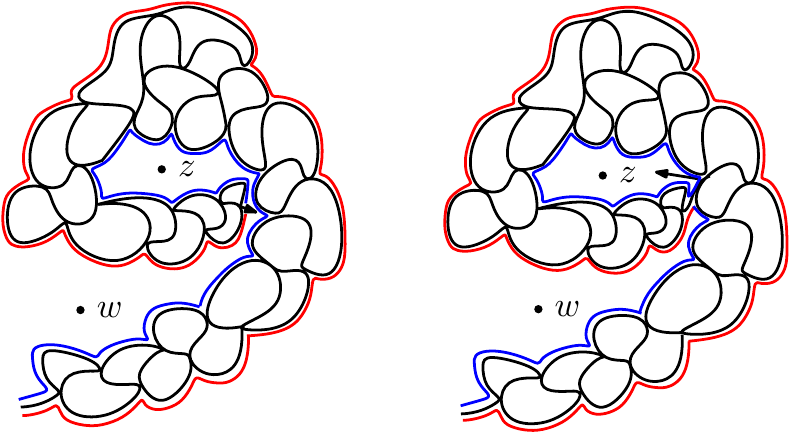}
		\caption{\emph{Left}: the colouring of the two sides of a planar curve. \emph{Right:} when $z$ and $w$ are first separated, $z$ is in a monochromatic component, while $w$ is in a bichromatic component (note that only part of the boundary of the component containing $w$ has been drawn, in fact this will be a finite, bounded component with probability one).}
		\label{F:monobi}
	\end{center}
\end{figure}

\subsubsection{Colouring} 
For both constructions it will be essential to have a notion of \textbf{colouring} of the two sides of an SLE$_\kp$ curve. Suppose $\eta$ is such a curve (or in fact, more generally, suppose $\eta$ is \emph{any} planar curve): we can colour the points immediately to the left of the curve with one colour (say, blue) and the points immediately to its right with another colour (say red). This choice of colours is made to match the conventions we adopted in Chapter \ref{S:maps} when discussing discrete space-filling paths on planar maps (recall, for example, Figure \ref{fig:bij}).

Formally, the colouring is a function defined on the \emph{prime ends} of $D_t = \C\setminus \eta((-\infty, t]) $ to $\{0,1\}$, for every $t\in \R$. We leave it to the reader to check that, when the curve is not space-filling, this assignment of colours is  \emph{consistent} as $t$ varies: that is, a prime end for $D_s$ also corresponds to a prime end for $D_t$ when $s<t$, and its colour at time $s$ also matches its colour at time $t$. (The assignment of colours can be defined precisely using Loewner theory, but we will leave the description at this informal level.) Because of the consistency of the colours, we can refer unambiguously to points on the left hand side of the curve (blue) and points on its right hand side (red). When the curve is space-filling, the set of prime ends of $\C\setminus \eta((-\infty,t])$ depends on the time $t$, but their colours, provided they exists, do not.

\medskip Now, suppose that $\eta$ takes values in $D$ (where $D$ could be a simply connected domain or could also be $\C$). Then for $t\in \R$, 
%is such that there exist two points $z,w$  in the same connected component of $D\setminus \eta((-\infty, s))$ for all $s<t$, but in separate connected components at time $t$: we call $t$ a \textbf{disconnection time}. If $C$ is one of the 
each connected component $C$ of $D\setminus \eta(-\infty,t)$ is either:

$\bullet$  \textbf{monochromatic}, if all the boundary points of $C$ (that is, all its prime ends) have the same colour; 

$\bullet$ or \textbf{bichromatic} (sometimes polychromatic), otherwise.  

See Figure \ref{F:monobi} for an illustration. If $t\in \R$ is such that there exist two points $z,w$  in the same connected component of $D\setminus \eta((-\infty, s])$ for all $s<t$, but in separate connected components of $D\setminus \eta((-\infty,t])$, we call $t$ a \textbf{disconnection time}. A useful observation is that at every such disconnection time, either the new components of $D\setminus \eta((-\infty,t])$ containing $z$ will be monochromatic and the one containing $w$ will be bichromatic, or vice versa. The behaviour of the space-filling SLE$_\kp$ path 
at this disconnection time will depend on which of these two situations occurs. Such a distinction is to be anticipated, bearing in mind the connection with the construction of the discrete path coming from Sheffield's bijection on decorated planar maps.\ind{Planar maps}

\subsubsection{Branching ordering, \texorpdfstring{$\kp\in(4,8)$}{TEXT} }
\ind{Space-filling SLE!Definition, $\kp \in (4,8)$}

 \begin{figure}
\begin{center}
\includegraphics[width=.4\textwidth]{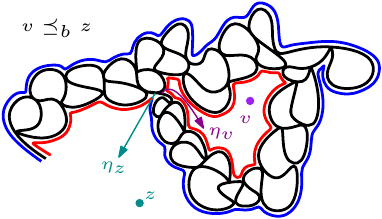} 
%\quad \quad 
\caption{The order defined by a branching SLE$_\kp$, $\kp \in (4,8)$.}
\label{F:brorder}
\end{center}
\end{figure}

We start with the branching ordering. We content ourselves with giving the definition in the whole plane (the chordal definition is completely analogous, we will outline how to adapt it to this case at the end). 
%except one needs to assign the boundary of the simply connected domain two distinct colours along both arcs separating the source and the target of the chordal SLE, the details are left to the reader). 
Let $\cQ$ denote a fixed dense countable set, and let $\{\eta_z\}_{z \in \cQ}$ denote a branching SLE$_\kp( \kp-6)$, which we defined in Section \ref{SSS:branching}.
\ind{Branching SLE$_{\kp}(\kp -6)$} 
\begin{definition}[Branching ordering] \label{D:branching48}
Fix $z \in \cQ$ and $w\in \cQ$ with $z \neq w$.  Let us say that $w \preceq_b z$ if, at the time that $w$ is disconnected from $z$ by $\eta_z$, then $w$ belongs to a \emph{monochromatic} component. 
\end{definition}

Equivalently, we could consider the branch $\eta_w$ targeted at $w$. Then $w \preceq_b z$ if and only if $z$ belongs to a \emph{polychromatic} component at the time when $\eta_w$ disconnects $z$ from $w$. By convention, we take $z \preceq_b z$ for any $z \in \cQ$. See Figure \ref{F:brorder} for an illustration. We leave it to the reader to check this does almost surely define a total order on $\cQ$.

\medskip For $z \in \cQ$, let $K^-_\cQ(z) = \{ w \in \cQ; w \preceq_b z\}$ be the past of $z$ (restricted to $\cQ$), and let $\cK^-_z$ denote its closure; this is the ``past'' of $z \in \cQ$. One can check that for a given $w\in \cQ$ and $z \in \cQ$ with $z\neq w$, the event $\{ w \preceq z\}$ is measurable. It is not hard to deduce that $\cK^-_z$ is also a random variable (on closed sets equipped with Hausdorff topology, say).

\subsubsection{Imaginary geometry ordering, \texorpdfstring{$\kp\in(4,8)$}{TEXT} }\label{SS:etaLR48}

We now wish to define the analogue of the curves $\eta^L_z$ and $\eta^R_z$ in the case $\kp \in (4,8)$, which was introduced in the case $\kp \ge 8$ in Section \ref{SSS:treeskp8} in order to introduce the alternative (``imaginary geometry'') description of the space-filling curve. 
Recall that when $\kp \ge 8$, the past $\cK^-_z$ coincides with the trace of the branch towards $z$, $\eta_z$, of the branching SLE$_\kp ( \kp-6)$. There was therefore no difficulty in talking about the boundary of the past, which is simply the boundary of $\eta_z$ and whose law is thus described by Lemma \ref{L:boundaryWP}. 

\medskip As everything else, the situation is more complicated when $\kp \in (4,8)$. Indeed, the past does not coincide with the trace of $\eta_z$ or even the filling in of the trace of $\eta_z$ (which one could define by filling in the components that are disconnected by $\eta_z$ and do not contain $z$). The issue is that some of these components, namely the bichromatic ones, will in fact be part of the future of the space-filling curve. 

\medskip We give the informal definition now. Consider the bichromatic components disconnected from $z$ by $\eta_z$, with the order they inherit from $\preceq_b$. That is, for $w \in \cQ$, let $C_z(w) $ be the component containing $w$ when $w$ is disconnected from $z$ by $\eta_z$, and consider the set $\cS = \{ C_z(w): w \in \cQ, z \preceq_b w\}$ which can be ordered as follows: $C_z(w_1) \preceq C_z(w_2) $ if $w_1 \preceq_b w_2$. Note that even though $w$ is not uniquely associated to its component $C_z(w)$, the above order is consistently defined.  

It can be checked that, almost surely, the boundary of each component $C = C_z(w)$ for some $w \in \cQ$ with $ z \preceq w$, consists of exactly two monochromatic arcs of opposite colours, which can be parametrised by curves. We call these respectively $\eta^L(C)$ and $\eta^R(C)$. The curves can be given a direction, which corresponds to the reverse chronological order with which these points were visited by the original curve $\eta_z$. Equivalently, $C$ lies to the right of $\eta^L(C)$ and to the left of $\eta^R(C)$. See Figure \ref{F:etaLR}. Then, by definition $\eta^L_z$ and $\eta^R_z$ is the result of the concatenation of these arcs $\eta^L(c)$ and $\eta^R(C)$ as $C$ varies across the set $\cS$, in the order defined above.  Note that it is not obvious (but true) that these give continuous simple curves, although this can be understood at a heuristic level by drawing enough pictures, and by considering the following situation.

\begin{figure}
	\centering
	\includegraphics[width=.5\textwidth]{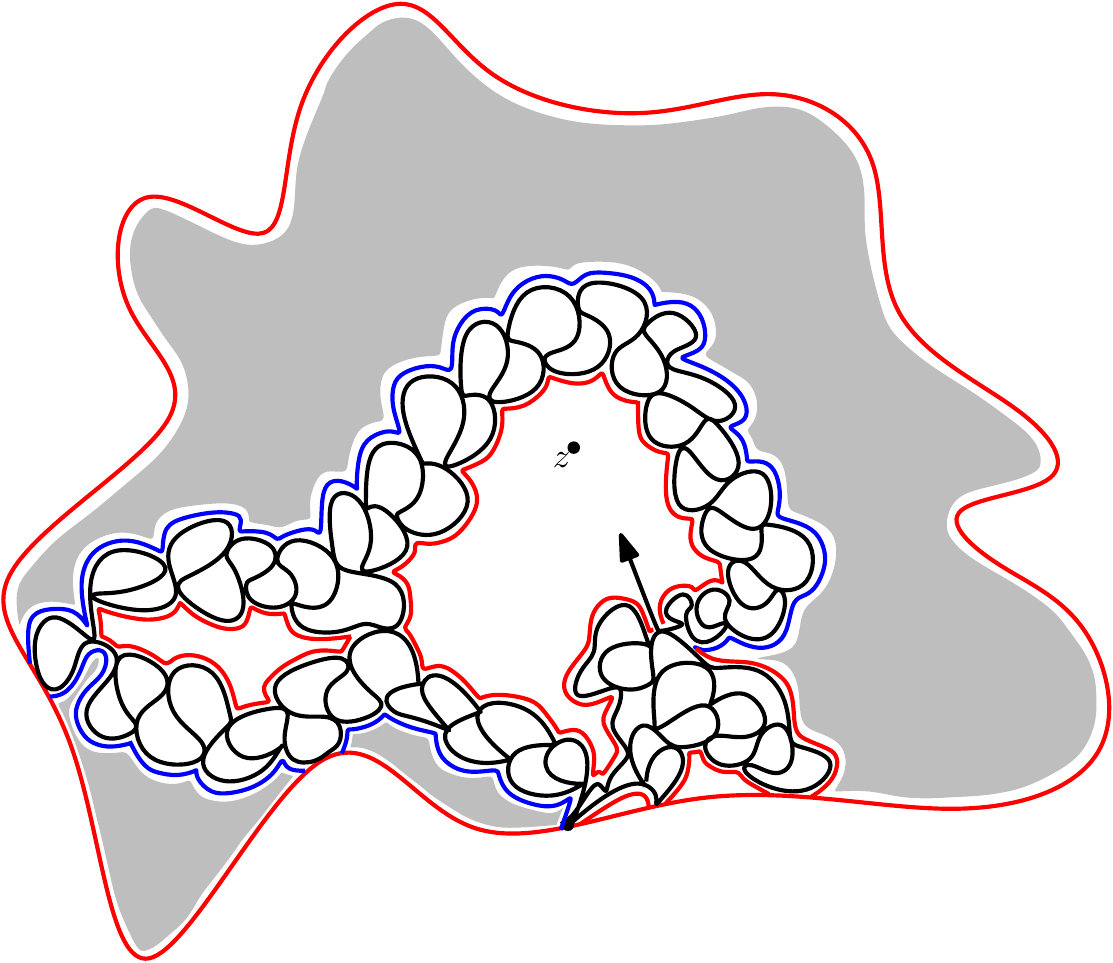}	
	\caption{The three bichromatic components cut off by $\eta_z$ between times $(s,t)$ are shaded in grey. \label{F:BCconnected}}
\end{figure}

Let $s\in \R$ be a time at which some $w \in \cQ$ is disconnected by $\eta_z$, and suppose that $z \preceq_b w$ so $w$ belongs to a bichromatic component, whereas $z$ belongs to the monochromatic one, call it $D$, and suppose without loss of generality that $D$ is coloured red. Consider the evolution of $\eta_z$ after time $s$, until the first time $t>s$ where the component containing $z$ (that is, the the connected component of the complement of $\eta_z((-\infty, t])$ containing $z$) does not share a positive proportion of the boundary with $D$. Between the times $s$ and $t$, the curve $\eta_z$ creates a number of bichromatic components which will be added to the set $\cS$. These components are simply created by the hits of $\eta_z$ to the boundary of $D$, but only those where the boundary is on the left curve when it hits it. These form a connected sequence of components, whose right boundary arcs will come from the boundary of $D$, and left boundary arc will come from the curve itself. See Figure \ref{F:BCconnected}.

% \begin{figure}
%\begin{center}
%\includegraphics[width=.6\textwidth]{figs/bicoloured} 
%\quad \quad 
%\caption{Bichromatic components separated from $z$ by $\eta_z$, and their boundary arcs. When concatenated, these forms the curves $\eta^L_z$ and $\eta^R_z$.}
%\label{F:bicoloured}
%\end{center}
%\end{figure}

A theorem from Miller and Sheffield \cite{IG4} shows that $\eta^L_z$ and $\eta^R_z$ are indeed curves and describes their joint law:

\begin{lemma}
 \label{L:boundarySFkp48}
Let $\kappa' \in (4,8)$  and let $\kappa = 16/\kappa' \in (2,4)$. Let $\eta^L = \eta^L_z, \eta^R = \eta^R_z$ be as above. Then, almost surely, $\eta^L_z, \eta^R_z$ are continuous curves. Furthermore, 
%\footnote{This is, by definition, the image of the whole plane SLE$_\k ( \k - 6)$ from 0 to $\infty$ under the M\"obius inversion $z \mapsto 1/z$.}.
\begin{itemize}
\item $\eta_L$ is a whole plane SLE$_{\kappa}( 2-\kappa)$ from $z$ to $\infty$  (note that $\eta_L$ is a simple curve when $\kappa' \ge  6$ but not when $\kp \in (4,6)$, cf. Lemma \ref{L:rhohit})

\item Given $\eta_L$, $\eta_R$ is a chordal SLE$_{\kappa} (  - \kappa/2, - \kappa/2)$ from $z$ to $\infty$ in $\C \setminus \eta_L ( -\infty, \infty)$ with force points on either side of the starting point $z$. (Note that since $\kappa'< 8$, $ - \kappa/2 <\kappa/2 - 2$ and so by Lemma \ref{L:rhohit}, $\eta_R$ hits both sides of $\eta_L$).

\end{itemize} 
\end{lemma}

Although the definitions of $\eta^L_z$ and $\eta^R_z$ are much more complicated in the case $\kp \in (4,8)$ than in the case $\kp \ge 8$, the above description is formally exactly the same in both cases, see Lemma \ref{L:boundaryWP}. As in that result, the joint law of $\eta^L$ and $\eta^R$ is actually symmetric: $(\eta^R, \eta^L)$ has the same law of $(\eta^L, \eta^R)$. Hence $\eta^R$ is also a whole plane SLE$_{\kappa}( 2-\kappa)$ from $z$ to $\infty$ and in particular is simple precisely when $\kp \in [6,8)$. Note that $\eta^L$ and $\eta^R$ hit themselves when $\kp \in (4,6)$ but not when $\kp \in [6,8)$, but they \emph{always} hit each other. 

\begin{rmk}\label{R:IG2}
It can be also be checked that the curves $\eta^L$ and $\eta^R$ have the same law as a pair of so called flow lines (\cite{IG4}) of a whole plane Gaussian free field with respective angles $-\pi/2,\pi/2$. As in the case $\kp \ge 8$, this will not be needed in any proof in the following. 
\end{rmk}

 \begin{figure}
\begin{center}
\includegraphics[width=.7\textwidth]{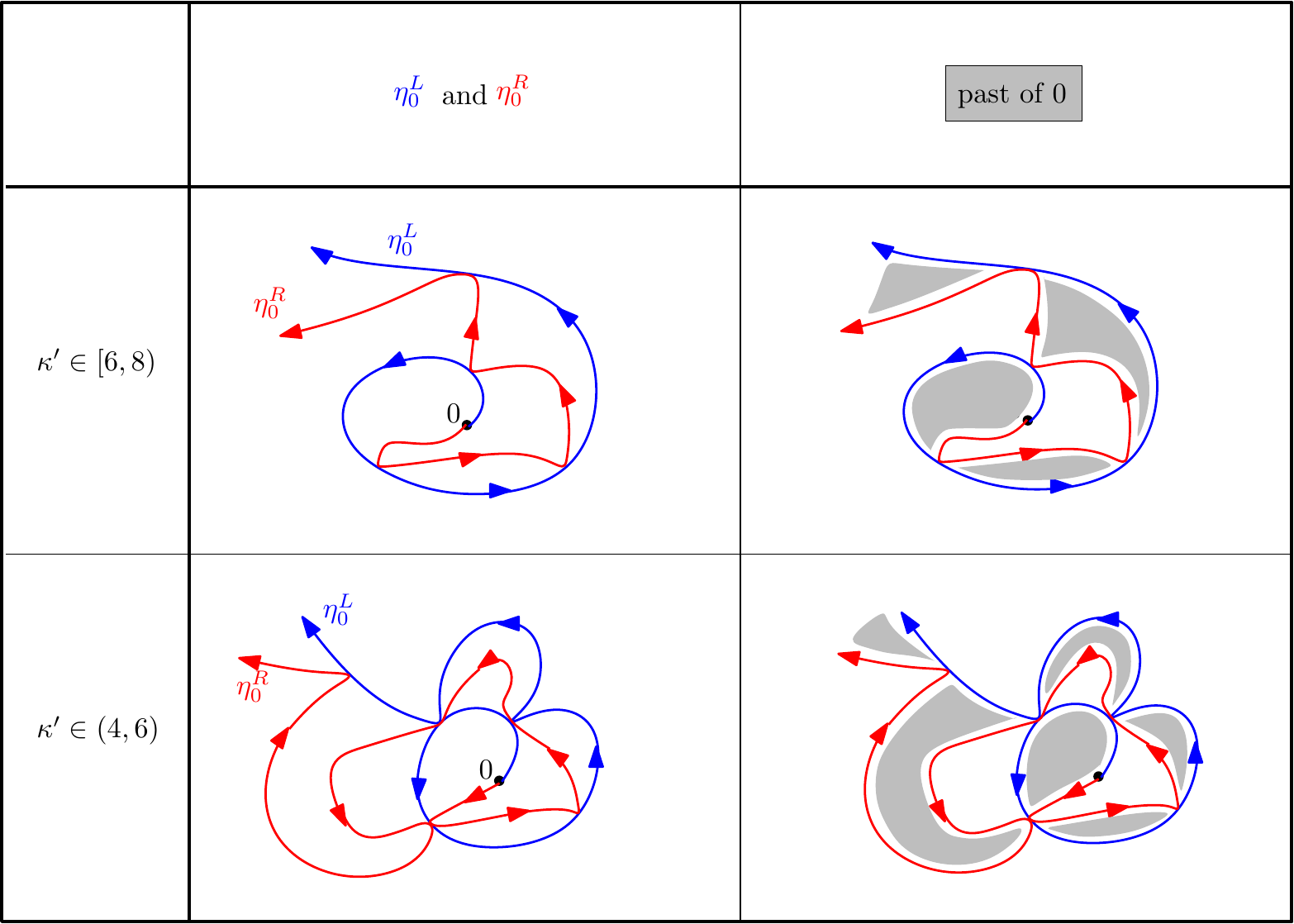} 
%\quad \quad 
\caption{The left and right boundary of the space-filling curve when $\kp \in(4,8)$.}
\label{F:etaLR}
\end{center}
\end{figure}

Having defined the curves $\eta^L_z$ and $\eta^R_z$ for $z \in \cQ$, we can finally give the description of the ``imaginary geometry'' ordering induced by the branching SLE$_\kp(\kp -6)$, $\kp \in (4,8)$. 
Take two points $z,w \in \cQ$, and consider their associated left-boundary paths (say), $\eta^L_z$ and $\eta^L_w$. Since the two branches going to $z$ and $w$ coincide for sufficiently negative times (up to parametrisation), it is straightforward to check that $\eta^L_z$ and $\eta^L_w$ must also merge eventually. 

Let us say
\begin{equation}\label{Eq:orderboundarykp48}
w \preceq_{\mathrm{IG}} z \text{ iff  $\eta_{w}^L$ merges with $\eta^L_{z}$ from the \emph{left}.}  
\end{equation}
Then we claim this order is identical to the branching order: that is, 
\begin{equation}\label{Eq:IGb48}
w \preceq_{\mathrm{IG}} z \text{ iff } w \preceq_{b} z.
\end{equation}
This is in fact easier to check in the case $\kp \in (4,8)$ than in the case $\kp \ge 8$, as here one can use the fact that the component of $w$ disconnected by $\eta_z$ is bounded, forcing the paths $\eta_z^L$ and $\eta_w^L$ to coincide once they both leave this component. 
%Equivalently we can  use the right-boundaries; in this case $z_i \preceq z_j$ if and only if $\eta_{z_i}^R$ merges with $\eta_{z_j}^R$ from the right. We refer to this ordering as the \textbf{Imaginary Geometry ordering}. \ind{ImaginaryGeometry}

Thus, the branching and Imaginary Geometry orders that are associated to the branching SLE$_\kp(\kp -6)$ coincide. What is left to say is that both of these orders define a unique continuous, space-filling curve $(\eta(t))_{t \in \R}$. This is the content of the following theorem, which can be derived from results in \cite{IG4} (which, roughly speaking, applies because of Remark \ref{R:IG2}).

\begin{theorem}
\label{T:sfkp48}
Let $\kp \in (4,8)$. There almost surely exists a unique curve $(\eta(t))_{t\in \R}$ which is space-filling and is continuous, such that $\eta(0) =0$ and $\mathrm{Leb} \eta(s,t) = |t-s|$ for any $s, t \in \R$, and for every $z,w \in \cQ$, if $t_z$ (resp. $t_w$) is the first time that $\eta$ visits $z$ (resp. $w$), then 
$$
z \preceq w \text{ iff } t_z \le t_w.
$$
Furthermore, the curve $\eta$ does not depend on the choice of $\cQ$. $\eta$ is called the (whole plane) space-filling SLE$_\kp$ from $\infty$ to $\infty$. 
\end{theorem} 

The space-filling SLE$_\kp$ is invariant under translation, rotation and scaling (up to time-reparametrisation). It is also, in fact, invariant under M\"obius inversion (again, up to reparametrisation) and thus under all M\"obius transformations of the Riemann sphere.

\subsection{Cutting and welding theorems}\label{SS:cuttingwelding}

We are now ready to describe the framework of \cite{DuplantierMillerSheffield} and to state some of the main theorems. The ultimate goal is to describe the exploration of a $\gamma$-quantum cone with an independent space-filling SLE$_\kp$ path, where $\gamma$ and $\kp$ are related by 
$$
\kp = 16/ \gamma^2,
$$
so $\k = 16/\kp = \gamma^2$ as was already the case in Chapter \ref{S:zipper}. We will present the results covering both cases $\kp \in (4,8)$ and $\kp \ge 8$ here; we thus only assume in what follows that $\kp >4$, and that we are given the existence and continuity of the space-filling SLE$_\kp$ from $\infty$ to $\infty$ in $\C$.

\dis{The proofs of the results in this section fall outside of the scope of this book, although we provide references to the proofs in the literature. They can be proved using similar tools to the proof of Theorem \ref{T:sliced} in Chapter \ref{S:zipper}, but the proofs are more involved. In Section \ref{SS:proofMOT} we will explain how they lead to the main result of \cite{DuplantierMillerSheffield}.}

Let $h$ denote the field of a $\gamma$-quantum cone, as defined in \cref{D:cone}, but embedded in $(\C, 0 , \infty)$ via the map $w \in \cC \mapsto z = -e^{- w} \in \C$ (so a neighbourhood of zero has finite $\gamma$ LQG mass, but a neighbourhood of $\infty$ has infinite $\gamma$ LQG mass). Recall that we say that $h$ is a \emph{unit circle embedding} of a $\gamma$-quantum cone.

Let $\eta$ be an independent space-filling SLE$_\kp$ curve from $\infty$ to $\infty$. A priori, $\eta$ comes parametrised so that $\eta(0) = 0$ and $\Leb ( \eta (s, t) ) = t- s$ for any $s<t$. However, it is crucial in the theorem below to reparametrise $\eta$ by its quantum area: that is, we define a reparametrisation $\eta'$ of $\eta$ such that 
\begin{equation}\label{Eq:qparam}
\mu_h ( \eta' (s,t)) = t - s
\end{equation}
for any $s< t$, where $\mu_h$ is the $\gamma$ Liouville measure (or area measure) associated to $h$, that is, $\mu_h(\dd x) = \lim_{\eps \to 0} \eps^{ \gamma^2/ } e^{\gamma h_\eps(x)} \dd x$ where the limit is  in probability (or almost surely along an appropriate subsequence) and $h_\eps (x)$ is some regularisation of $h$ at scale $\eps$, as described in Chapter \ref{S:GMC}. %for some specific choice of log-correlated Gaussian field and reference measure.

It is not obvious, but true, that one can reparametrise $\eta$ so that \eqref{Eq:qparam} holds. This follows from the fact that almost surely, for any $s< t$, $\eta (s,t)$ contains an open ball, and any open ball has positive mass (since it contains a ball centred at a point with rational coordinates and with rational radius). Likewise, the reparametrisation $\eta'$ is such that there are no intervals of constancy, which follows from the fact that $\mu_h$ has no atoms almost surely (itself a consequence of, for example, Exercise \ref{Ex:atom} in Chapter \ref{S:GMC}). 

Let $\eta( - \infty, 0] = K^-_0$ denote the past of zero, and let $\eta^L_0$ and $\eta^R_0$ denote its left and right boundaries (which, we recall are given by a pair of SLE curves whose joint law is specified by Lemma \ref{L:boundaryWP} for $\kp \ge 8$ and Lemma \ref{L:boundarySFkp48} for $\kp \in (4,8)$). The first result below states that the boundaries $\eta^L_0 $ and $\eta^R_0$ divide the quantum cone into two regions (namely, the past and the future of zero), on which the restrictions of $h$ define two independent quantum wedges with parameter $\alpha = 3\gamma/2$ (in the terminology favoured by \cite{DuplantierMillerSheffield}, the wedges have ``weight'' $W = 2 - \gamma^2/2$).%, see Remark \ref{R:weights}). 

Note that when $\kp\ge 8$ the wedge parameter satisfies $\alpha \le Q$ and is therefore ``thick'' in the terminology of Chapter \ref{S:SIsurfaces}, whereas $\alpha >Q$ is ``thin'' for $\kp \in (4,8)$. Recall that such a thin wedge corresponds to an ordered collection of beads; these beads correspond precisely to the bichromatic components created by the branch $\eta_0$ (targeted at 0) of the branching SLE$_\kp(\kp -6)$.

The first theorem we present describes the surface that one obtains by ``cutting'' the quantum cone $(\C,h,0,\infty)$ with just one half of the boundary of $\eta_0$, say $\eta^L_0$. 

\begin{theorem}
\label{T:cuttingcone1} Suppose $\kp >4$, and let $h$ and $\eta$ be as described just above \eqref{Eq:qparam}.

%\begin{itemize}
%\item 
Let $D$ denote $\C \setminus \eta^L_0(\R)$. Then $\cW = ( D, h|_D, 0, \infty)$ has the law of a quantum wedge of parameter $\alpha  = 2\gamma - 2/\gamma$, which is thick when $\gamma^2 \le 8/3$, equivalently $\kp \ge 6$; and is thin  when $\kp \in (4,6)$.
%
%\item Let $D^-$ denote the interior of $K^-_0$ and let $D^+ = \C \setminus K^-_0$. Let $\cW^- = ( D^-, h|_{D^-}, 0 , \infty)$ and let $\cW^+ = (D^+, h|_{D^+}, 0, \infty)$. Then $\cW^-$ and $\cW^+$ are independent quantum wedges with parameter $\alpha = 3\gamma/2$. 
%\end{itemize}
\end{theorem}

\begin{rmk} This is \cite[Theorem 1.2.4]{DuplantierMillerSheffield}. \end{rmk}

When $\kp \in (4,6)$, $\eta_0^L$ touches itself almost surely, as already noted in Lemma \ref{L:boundarySFkp48}. Thus $D= \C\setminus \eta^L_0(\R)$ is not simply connected, but consists instead of a countable (ordered) collection of simply connected domains. In this case the theorem has to be understood as saying that the restriction of $h$ to these domains form a wedge of parameter $\alpha = 2 \gamma - 2/\gamma$. This corresponds to the fact that the wedge is thin in this case. 
When we cut along the second half of the boundary, the result is the following.

\begin{theorem}
\label{T:cuttingcone2} Suppose $\kp >4$, and let $h$ and $\eta$ be as described just above \eqref{Eq:qparam}.

%\begin{itemize}
%\item 
%Let $D$ denote $\C \setminus \eta^L_0(\R)$. Then $\cW = ( D, h|_D, 0, \infty)$ has the law of a quantum wedge of parameter $\alpha  = 2\gamma - 2/\gamma$ (which is thick when $\gamma^2 \le 8/3$, equivalently $\kp \ge 6$; and is thin  when $\kp \in (4,6)$)
%
 Let $D^-$ denote the interior of $K^-_0$ and let $D^+ = \C \setminus K^-_0$. Let $\cW^- = ( D^-, h|_{D^-}, 0 , \infty)$ and let $\cW^+ = (D^+, h|_{D^+}, 0, \infty)$. Then $\cW^-$ and $\cW^+$ are independent quantum wedges with parameter $\alpha = 3\gamma/2$. 
%\end{itemize}
\end{theorem}

\begin{rmk} This is \cite[Theorem 1.2.1]{DuplantierMillerSheffield}. \end{rmk}

Once again, the way to read this theorem properly depends on the value of $\kp$. Indeed when $\kp \in (4,8)$, neither $D^-$ nor $D^+$ are simply connected, instead they consist of ordered collections of simply connected domains. In that the case the theorem states that the restriction of $h$ to these domains form independent quantum wedges, each with parameter $\alpha = 3\gamma/2$ (this corresponds to the thin case precisely when $\kp \in (4,8)$). 

\begin{figure}
	\centering
	\includegraphics[width=\textwidth]{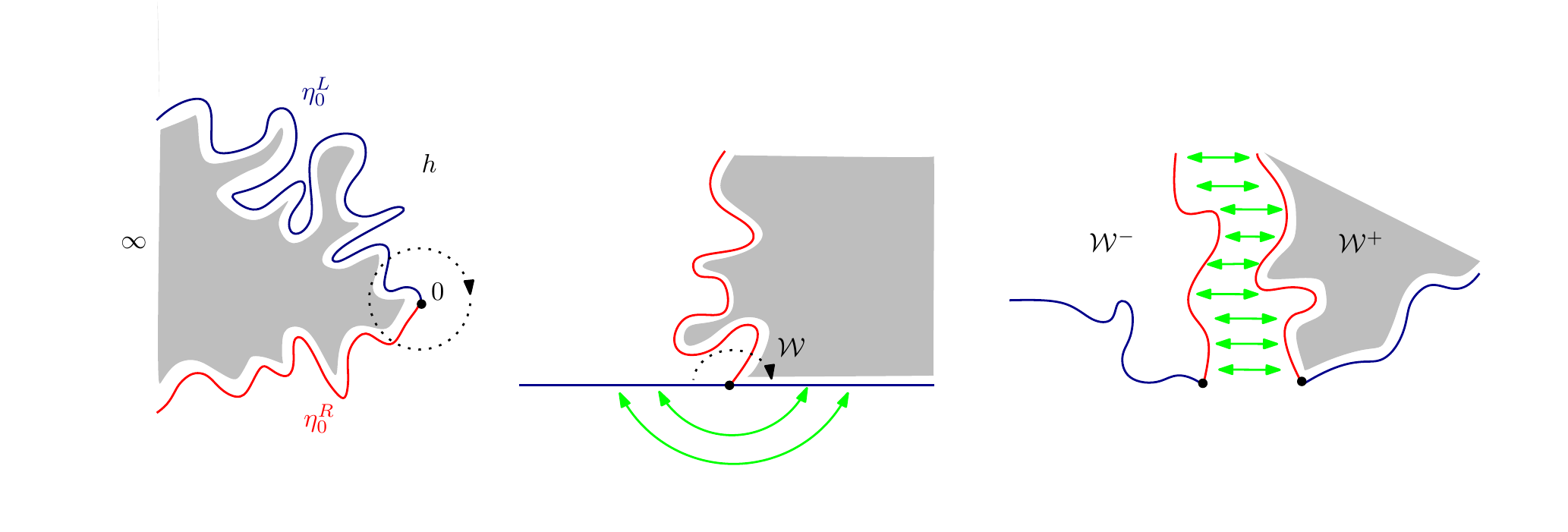}
	\includegraphics[width=.8\textwidth]{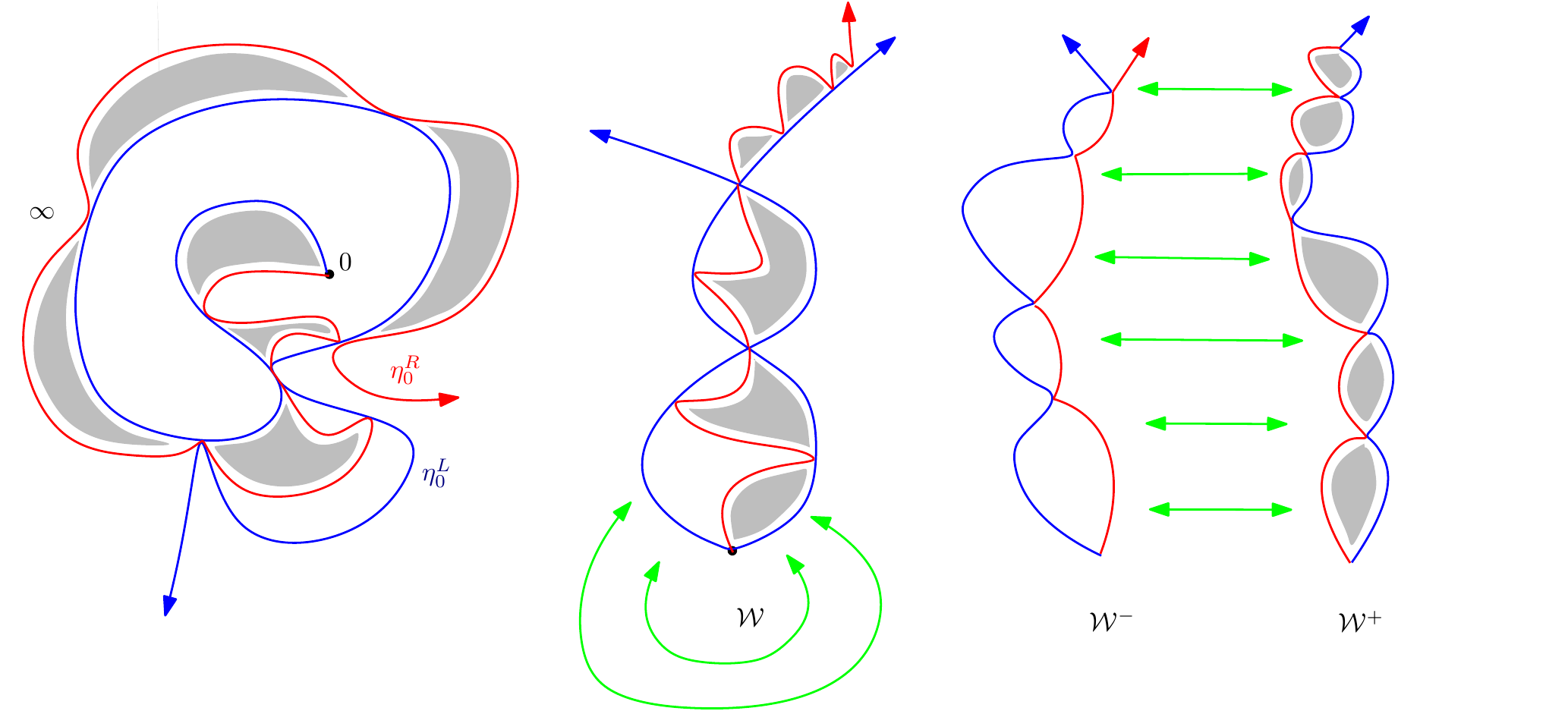}		
	\caption{Illustrations of the cutting and welding theorems when $\kp\ge 8$ (top) and when $\kp\in (4,8)$ (bottom). To get from the left picture to the central picture, one ``cuts'' along the curve $\eta_0^L$ -- more precisely, considers the quantum cone field $h$ in $D=\C\setminus \eta_0^L$, and views $\mathcal{W}:=(D,h,0,\infty)$ as a quantum surface. Theorem \ref{T:cuttingcone1} says that $\mathcal{W}$ has the law of a quantum wedge with parameter $\alpha=2\gamma-2/\gamma$ (we have conformally mapped $D$ to the upper half plane in the central figure, which illustrates a different embedding, or equivalence class representative in the sense of doubly marked quantum surfaces, of $\mathcal{W}$). To get from the central picture to the right picture, one cuts further along the image of $\eta_0^R$, and considers the restriction of the field to either side, to define a pair of quantum surfaces $\mathcal{W}^-,\mathcal{W}^+$. Theorem \ref{T:cuttingcone2} says that these are independent quantum wedges with parameter $\alpha=3\gamma/2$. The welding theorem, Theorem \ref{T:weldingcone} (which is stated only in the case $\kp\ge 8$) describes what happens when goes from right to left in the above pictures.The operation at each step is illustrated by the identification, or ``welding'', of points at equal quantum boundary length away from the black dot (this identification is depicted with green arrows).
		\label{F:welding}}
\end{figure}

\begin{rmk}\label{R:cuttingwedge}
A consequence of this statement and the description of the conditional law of $\eta^R_0$ given $\eta^L_0$ (provided in Lemma \ref{L:boundaryWP}) is that if one takes a quantum wedge $\cW = ( \H , \tilde h, 0, \infty)$ with parameter $\alpha=2\gamma-2/\gamma$ and cuts it with an independent chordal SLE$_\k ( - \k /2, - \k/2)$ curve from $0$ to $\infty$ with force points on either side of zero, then the restriction of $\tilde h$ to the complement of this curve defines two independent quantum wedges $\cW^-$ and $\cW^+$ with parameter $\alpha = 3 \gamma /2$, as in the theorem. In reality, the proof of the theorem goes in the converse direction: that is, one first establishes this fact and Theorem \ref{T:cuttingcone1} in order to deduce Theorem \ref{T:cuttingcone2}. However, it is Theorem \ref{T:cuttingcone2} which is the key input for the mating of trees theorem. 
\end{rmk}

\medskip The identification of $\cW$ and of $\cW^-, \cW^+$ as quantum wedges means we can talk about their boundary length measures. Theorems \ref{T:cuttingcone1} and \ref{T:cuttingcone2} can be complemented by a result showing that the boundary lengths naturally match with one another along $\eta^L_0$ and $\eta^R_0$. In other words, the quantum cone $(\C,  h, 0, \infty)$ can be viewed as a conformal welding of $\cW$ with itself (where points on either side of $0$ of equal boundary length are identified with one another), and can also be viewed as a conformal welding of $\cW^-$ with $\cW^+$, where points at equal (signed) distance from 0 along the boundary in $\cW^-$ and $\cW^+$ are identified with one another. 

\medskip Put it another way, using Remark \ref{R:cuttingwedge}, we can conformally weld $\cW^-$ and $\cW^+$ along just one half of their boundary (identifying points at equal distance from 0) to get $\cW$. Subsequently, we can conformally weld the two halves of the boundary of $\cW$ to obtain the quantum cone $(\C, h, 0, \infty)$. This can all be encapsulated in the following welding theorem, which for simplicity we only state in the case $\kp \ge 8$.

\begin{theorem}
\label{T:weldingcone}
In the same settings as Theorem \ref{T:cuttingcone1} and \ref{T:cuttingcone2}, let $\kp \ge 8$ and let $g^+$ (resp. $g^-$) be a conformal isomorphism from $D^+$ (resp. $D^-$) to $\H$ sending 0 to 0 and $\infty$ to $\infty$, with $\eta^L_0$ being mapped to $(-\infty, 0]$ and $\eta^R_0$ being mapped to $[0, \infty)$. Let $h^+$ (resp. $h^-$ denote the image of $h|_{D^+}$ under $g^+$ (resp. of $h|_{D^-}$ under $g^-$) under the change of coordinate formula \eqref{E:coordinate_change}, and let $\cV_{h^+}$ (resp. $\cV_{h^-}$) denote the boundary length measure of $h^+$ (resp. $h^-$) in $\H$. 

Almost surely the following statement holds for all points $z$ on either $\eta^L_0$ or $\eta^R_0$. Let $z^+$ (resp. $z^-$) denote the image of $z$ under $g^+$ (resp. $g^-$) in $\R$. Let $I^+$ (resp. $I^-$) denote the interval between 0 and $z^+$ (resp. $z^-$).Then 
$$
\cV_{h^+} (I^+) = \cV_{h^-} ( I^-).
$$
\end{theorem}

\begin{rmk} This follows from \cite[Theorems 1.2.1 and 1.2.4]{DuplantierMillerSheffield}. \end{rmk}

Note that in particular, Theorem \ref{T:weldingcone} allows us to unambiguously define the quantum length of any measurable portion of $\eta^L_0$ or $\eta^R_0$ with respect to the quantum cone $(\C, h, 0, \infty)$.

%Recall that the proofs of these results will be delayed to Section \ref{SS:proofMOT}.

\subsection{Statement of the mating of trees theorem }
\label{SS:MOTstatement}

As in the previous subsection, we let $(\C, h, 0, \infty)$ be a $\gamma$-quantum cone and $\eta$ be an independent space-filling SLE$_\kp$ ($\kp = 16/\gamma^2$), parametrised by Lebesgue area in such a way that $\eta(0)=0$. We let $\eta'$ be the reparametrisation of  its quantum area $\mu_h$ relative to time $0$ (which induces a dependence between $h$ and $\eta'$). 

We have explained above how Theorem \ref{T:weldingcone} allows us to unambiguously define the quantum length of any measurable portion of $\eta^L_0$ or $\eta^R_0$ with respect to the quantum cone $(\C, h, 0, \infty)$. In order to state the mating of trees theorem we will also need to measure the quantum lengths (with respect to $(\mathbb{C},h,0,\infty)$) of the curves $\eta^L_z$ and $\eta^R_z$, when $z$ is of the form $z = \eta'(t)$ for $t \in \R$ fixed. (Recall that when $\kp\ge 8$, $z=\eta'(t)$, $\eta^L_z,\eta^R_z$ are the left and right boundaries of $\eta'[0,t]$, and when $\kp\in (4,8)$ they are slightly more complicated to define, see Section \ref{SS:etaLR48}, but still correspond to the left and right boundaries of $\eta'[0,t]$ in an appropriate sense). The following key lemma shows that the quantum cone decorated with the space-filling path $\eta'$, viewed from $\eta'(t)$, is in fact stationary, and this (in particular) implies that the quantum lengths described above are well defined.

To state the lemma, recall the notion of a curve decorated random surface $[(D, h, a, b); \eta]$ from  \cref{D:surfacecurve}. The stationarity will be in the sense of such objects.

\begin{lemma}
\label{L:conestationary} Let $t \in \R$, and let $z= \eta'(t)$. 
Then the law of the curve decorated surface $  [ (\C, h , z , \infty ); \eta' ( t + \cdot) ]$ is the same as that of $ [( \C, h, 0 , \infty ); \eta'(\cdot))]$. 
\end{lemma}

\begin{rmk}
	For this statement it is crucial that $\eta'$ is parametrised by its quantum area.
\end{rmk}

To spell out what the statement really says, we warn the reader that it would be incorrect to say that the joint law of $( h ( z+ \cdot), \eta'( t + \cdot))$ is the same as that of $( h, \eta')$. Indeed, the laws of $h$ and $h ( z + \cdot)$ are \emph{not} the same \emph{as fields}; we would have to applying a random rescaling to $h(z+\cdot)$ for this to be the case. Nonetheless, the objects in the lemma have the same law \emph{as curve decorated random surfaces.} 

\dis{One can prove Lemma \ref{L:conestationary} in a similar manner to the proof of Proposition \ref{P:qwkey} in Chapter \ref{S:zipper}, but we will not provide the details in this book (we direct the interested reader to \cite[Proof of Lemma 8.1.3]{DuplantierMillerSheffield}). We will instead focus, in Section \ref{SS:proofMOT}, on how this leads to the proof of the main theorem of \cite{DuplantierMillerSheffield} (Theorem \ref{T:MOT} below).}

\medskip We now turn to the statement of one of the main theorems of this chapter. Let $h$ and $\eta'$ be as above. Let us define a process $(L_t, R_t)_{t\in \R}$ as follows. Informally, $L_t$ tracks the change in the length of the left (outer) boundary of $\eta'(t)$, relative to time zero, whereas $R_t$ tracks the same change but for the right (outer) boundary.  To define it formally, fix $s<t$, and let $w = \eta'(s), z = \eta'(t)$. Then we define the increment
\begin{equation}\label{E:boundarylengthLR}
L_t - L_s : = \cV_{h} ( \eta_z^L \setminus \eta_w^L) - \cV_h ( \eta_w^L \setminus \eta_z^L),
\end{equation}
and make the same definition for $R_t$ except that $\eta^L$ is replaced with $\eta^R$ in all occurrences. If we also set $L_0 = R_0 = 0$, then \eqref{E:boundarylengthLR} specifies a unique two-sided process $(L_t, R_t)_{ t \in \R}$. Note that the meaning of the random variables in \eqref{E:boundarylengthLR} measuring the lengths of various boundary curves is provided by Theorem \ref{T:weldingcone} (see the discussion immediately below that theorem) and the stationarity of Lemma \ref{L:conestationary}. With these definitions we can finally state the main theorem below.

\begin{theorem}
\label{T:MOT} Let $h, \eta'$ be as above and let $(L_t , R_t)_{t\in \R}$ denote the boundary length process \eqref{E:boundarylengthLR}. There exists $a>0$ depending solely on $\gamma\in (0,2)$ such that $(L_{at}, R_{at})_{t \in \R}$ is a two-sided correlated Brownian motion in $\R^2$, with 
$$
\var (L_{at}) = \var (R_{at}) = |t| ;  \quad \quad \cov ( L_{at}, R_{at} ) = - \cos \left( \frac{4\pi}{\kp} \right)|t|.
$$
 \end{theorem}

Observe that the Brownian motions are negatively correlated for $\kp \ge 8$, positively correlated when $\kp \in (4,8)$, and independent when $\kp =8$ (which corresponds to the case of the uniform spanning tree).

\begin{rmk}
The value of the constant $a$ appearing in the statement of that theorem was unknown for some time, until a recent work of of Ang, R\'emy and Sun \cite{AngRemySun} who computed it using tools coming from Liouville conformal field theory. 
\end{rmk}

\subsection{Discussion and uniqueness}\label{SS:MOTdiscussion}

Theorem \ref{T:MOT} should be compared with Theorem \ref{T:scalinglimitFK}. In that theorem, we showed that the scaling of the left and right relative boundary lengths of a space-filling path (these are precisely the hamburger and cheeseburger counts) exploring the infinite volume random planar map weighted by the self dual critical Fortuin--Kasteleyn percolation model, is also given by a pair of correlated Brownian motions. Identifying the limiting covariance of Theorem \ref{T:scalinglimitFK} with that in Theorem \ref{T:MOT} gives a relation between $q$ and $\gamma$ (or equivalently $\kp$) which is the same as the one announced in \cref{Table:values}
$$
q = 2 + 2 \cos ( \gamma^2 \pi/2)  = 2  \cos^2 ( 4\pi/\kp).
$$ 
This is consistent with the physics prediction discussed in Chapter \ref{S:maps}.

\subsubsection{A mating of trees?}
\label{SS:MOTdescription}

Before we explain some of the key steps going into the proof of Theorem \ref{T:MOT}, we spend some time explaining why this theorem is related to a ``mating of trees''. The word ``mating'' (that is, gluing) originates from the field of complex dynamics. It was coined by Douady and Hubbard \cite{Douady} who spoke of matings of polynomials to describe a way to glue together two 
(connected and locally connected) filled Julia sets along their boundaries. In a sense, Theorem \ref{T:MOT} describes a similar construction where the role of the Julia sets is played by two infinite continuum random trees (CRT). \ind{Continuous Random Tree (CRT)}

Let us first briefly explain the notion of Continuous Random Tree, originally due to Aldous \cite{Aldous}; we refer to the lecture notes by Le Gall \cite{LeGallsurvey} for a much more complete discussion and additional references, including in particular the history and applications of this important subject. Traditionally the theory is defined from a Brownian excursion $(e_t)_{0 \le t \le 1}$; the resulting continuum tree would then be a compact metric space. However, for our purposes it will be more natural to consider infinite volume analogues of this CRT, in which case the Brownian path defining the tree is simply a (real valued) two-sided Brownian motion $(B_t)_{t\in \R}$. The definition is simply the following and can be made path by path, that is almost surely given a fixed continuous function $f:\R\to \R$ (which will later be taken to have the law of the two-sided Brownian motion $B$). Given $s, t \in \R$, let us define an equivalence relation $\sim_f$ 
\begin{equation}\label{eq:equivrel}
s\sim_f t \text{ if } f(s) = f(t) = \inf_{u \in [s,t]} f(u).
\end{equation}
(Here one can have $s\le t$ or $t\le s$). It is easy to see that this defines an equivalence relation. By definition, the (infinite) continuous random tree associated to $B$ is simply equal to the quotient space $\cT_f = \R/\sim_f$. We can turn $\cT_f$ into a metric space by considering, for $s, t\in \R$ (say with $s \le t$ without loss of generality),
$$
d_{f} (s,t) = f(s) + f(t) - 2 \inf_{u \in [s,t]} f(u).
$$
(Note that $d_f(s,t)) = 0$ if and only if $s\sim_f t$, as required for a metric). This metric turns $\cT_f$ into what is known as a real or $\R$-tree: that is, any two simple curves $\sigma_1$ and $\sigma_2$ in $\cT$ (that is, injective continuous maps from $[0,1]$ to $\cT_f$) with same starting and endpoints must be reparametrisations of one another. 

A convenient way of visualising the tree $\cT_f$ associated to $f$ is to imagine painting the underside of the graph of $f$ with glue, and then squishing this graph horizontally (see below for a more precise description). Indeed the points that are glued with one another in this process correspond exactly to those that are identified via \eqref{eq:equivrel}. This suggests another way of describing $\cT_f$ which will here be more natural. Consider the portion of the $(t,x)$ plane lying below the graph of $f$: that is 
\begin{equation}\label{eq:underside}
\Gamma_f = \{ (t, x) \in \R^2: x \le f(t)\}. 
\end{equation}
Define an equivalence relation $\approx_f$ on $\Gamma_f$ as follows: for every $s,t\in \R$ such that $s\sim_f t$ put a horizontal segment between $(s, f(s))$ and $(t, f(t))$ (note this segment lies entirely in $\Gamma_f$ by definition of $\sim_f$) and identify all the points of $\Gamma_f$ lying on this segment; these identifications describe the equivalence classes of $\approx_f$. Now, $\Gamma_f$ inherits a topological structure from $\R^2$, thus turning the quotient space $\Gamma_f/\approx_f$ into a topological space. Furthermore, the equivalence classes of $\approx_f$ are clearly in bijection with those of $\sim_f$, hence $$(\Gamma_f /\approx_f) \quad  = \quad \cT_f,$$
in the sense, for example, that these two topological spaces are homeomorphic.
%
%It is furthermore intuitive that this process results in a tree which is nothing else but $\cT_f$. 

\begin{figure}
	\centering
	\includegraphics[width=.5\textwidth]{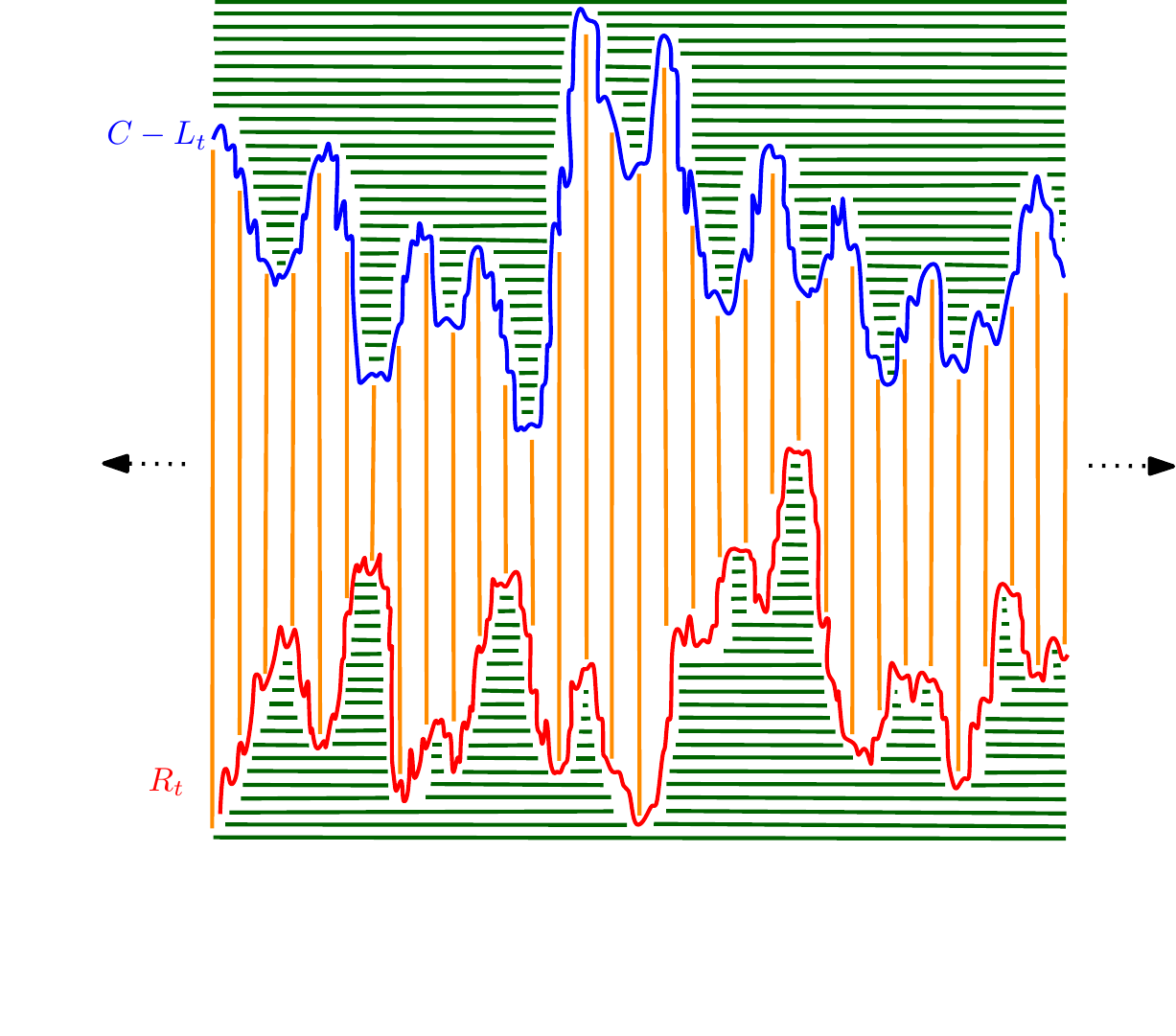}
\caption{The gluing of two CRT produces a topological surface (plane or sphere) equipped with a space-filling path, that is, a \textbf{peanosurface}.}
\label{F:peano}	
\end{figure}

\medskip Coming back to Theorem \ref{T:MOT}, let $(L_t)_{t\in \R}$ and $(R_t)_{t\in \R}$ denote the correlated two-sided Brownian motions describing the relative left and right boundary lengths associated with the quantum cone $h$ and the space-filling path $\eta'$ as in Theorem \ref{T:MOT}. As mentioned above, each of $L$ and $R$ separately encode an infinite CRT, which we denote by $\cT_L$ and $\cT_R$. In addition, the space-filling path $\eta'$ gives a natural way to identify (and hence glue) points on $\cT_L$ and $\cT_R$. More precisely, for $t\in \R$, let $\ell (t) \in \cT_L$ denote the point of $\cT_L$ corresponding to time $t$ (that is, the equivalence class of $t$ for  $\sim_L$). Similarly, let $r(t) \in \cT_R$ denote the point of $\cT_R$ corresponding to time $t$ (the equivalence class of $t$ for $\sim_R$). Since $\eta'(t)$ visits both $\ell(t)$ and $r(t)$, it is natural to identify $\ell(t)$ and $r(t)$. Somewhat miraculously this identification can be seen to give rise to a topological surface $M$ (in fact a topological plane) in the sense that $M$ is a topological space, almost surely homeomorphic to the plane. 

\medskip Let us explain this construction in more detail. In our infinite volume setting we will in fact first describe a finite volume approximation. To this end we fix $T>0$ and consider the restriction of $L$ and $R$ to $[-T, T]$. Pick a constant $C>0$ (depending on $T$ as well as on $L$ and $R$) such that $C> \sup_{|t| \le T } L_t+ \sup_{|t| \le T} R_t $. Consider the graphs of $(R_t)_{-T \le t\le T}$ and $(C- L_t)_{-T \le t \le T}$, drawn simultaneously as in Figure \ref{F:peano}. By our choice of $C$ these two graphs do not intersect, and in fact the graph $C-L$ sits entirely above the graph of $R$; beyond this the value of $C$ will not matter. Consider the closed rectangle $\cR$ of the $(t,x)$-plane containing the graphs of $R$ and $C-L$; that is, 
$$
\cR= \{ (t,x) \in \R^2: -T \le t \le T, \inf_{|u| \le T} R_u \le x \le \sup_{|u|\le T } C - L_u \}. 
$$
$\cR$ inherits a topological structure from $\R^2$. We will now consider an equivalence relation $\cong$ on $\cR$, defined as follows. On the underside of the graph $R$ (restricted to $[-T, T]$) we draw the horizontal segments $\approx_R$ as explained in \eqref{eq:underside}. On the upperside of the graph of $C-L$ we can draw the analogous horizontal segments (see Figure \ref{F:peano}). To these horizontal segments we add a \emph{vertical} segment joining $(t, R_t)$ to $(t, C- L_t)$ for each $-T \le t \le T$. Having drawn these horizontal segments, the equivalence relation $\cong$ is defined by identifying any two points lying on the same horizontal segment and any two points lying on the same vertical segments. 

While most equivalence classes in $\cong$ consist of just one segment (vertical), it is possible for an equivalence class to contain more than one segment. For instance, if $s\sim_R t$ then at least three segments are identified with one another: the horizontal segment containing the $(s, R_s)$ and $(t, R_t)$ but also the vertical segments containing these two points. In principle, doing this construction with arbitrary pairs of continuous functions we could have equivalence classes with arbitrary many segments. However it is possible to see that, when $L$ and $R$ are correlated Brownian motions, an equivalence class has at most five segments almost surely (of which two are then horizontal, corresponding to a branch point in the CRT, and three are vertical). Importantly, no equivalence class consists of four segments forming a rectangle (two vertical and two horizontal segments).

The equivalence relation $\cong$ on $\cR$ is furthermore \emph{topologically closed}: that is, if $x_n \cong y_n$ and $x_n \to x$, $y_n \to y$ then necessarily we have $x\cong y$. For such closed equivalence relations, there is a very nice criterion due to Moore \cite{Moore} (see Milnor \cite{Milnor} for a more modern formulation),  which can be used to check whether the quotient space retains the topology of $\cR$ (that is, a closed disc here). Namely, no equivalence class should disconnect $\cR$ into more than one connected component; indeed such an equivalence class would correspond to a pinch point in $\cR/\cong$, and would prevent the quotient from being homeomorphic to a closed disc. Here it can be checked this is almost surely the case, precisely because no equivalence class may consist of a rectangle. See \cite{DuplantierMillerSheffield} for details of these arguments.

The identification of $\cT_R$ with $\cT_L$ over $[-T, T]$ thus gives us a topological space $\cR/\cong$, which is homeomorphic to a closed disc. This closed disc also comes equipped with a natural space-filling path (call it $\tilde \eta(t), t\in [-T, T]$), which at time $t\in [-T, T]$ visits the equivalence class corresponding to the vertical segment joining $(t, R_t)$ with $(t, C- L_t)$. The pair $( \cR/\cong, \eta)$ is what we call a \textbf{peanosurface} (here a closed ``peanodisc''). Sending $T\to \infty$ in the natural way gives us an ``infinite volume'' version of this construction.

\medskip The upshot is that Theorem \ref{T:MOT} gives us access to a peanosurface, constructed as above from the gluing of the two trees $\cT_L$ and $\cT_R$ associated to the relative left  and right boundary length of the decorated quantum cone. At this point the parallel with the Theorem \ref{T:sheffieldbij} coming from Sheffield's bijection for FK-weighted random planar maps should be clear. Indeed these discrete planar maps could also be described as a gluing of two discrete trees whose scaling limit is given by two correlated CRTs (Theorem \ref{T:scalinglimitFK}).

\subsubsection{Uniqueness}

Theorem \ref{T:MOT} and the above discussion make it clear that we can associate to a quantum cone $h$, decorated by a space-filling SLE path $\eta'$, a peanosurface which is obtained from the gluing of two infinite correlated CRTs. The parallel with the discrete theory described in Chapter \ref{S:maps} raises the following question: do the processes $(L_t, R_t)_{t\in \R}$ characterise the pair $(h, \eta')$ uniquely? This question is natural because, in the discrete, there is a bijection between the trees and the decorated map.  Remarkably this turns out to be the case, as stated in Theorem 1.4.3 of \cite{DuplantierMillerSheffield}.

\begin{theorem}
\label{T:MOTuniqueness}
In the setting of Theorem \ref{T:MOT}, the pair $(L_t, R_t)_{t\in \R}$ almost surely determines $(h, \eta')$ uniquely up to a rotation of the plane. That is, suppose that $(h_1, \eta'_1)$ and $(h_2, \eta'_2)$ are two quantum cones (with a unit circle embedding) defined on the same probability space, and $\eta'_i$ is a space-filling SLE$_{\k'}$ independent of $h_i$ parametrised by the respective quantum area ($i=1,2$). Let $(L^i_t, R^i_t)_{t\in \R}$ ($i=1,2$) be their associated left and right relative boundary lengths, and suppose also that $L^1_t = L^2_t$ and $R^1_t = R^2_t$ for all $t \in \R$. Then $(h_2, \eta'_2)$ is obtained from $(h_1, \eta'_1)$ by applying a fixed rotation around the origin. 
\end{theorem}

When $(h_1,\eta'_1)$ and $( h_2, \eta'_2)$ are defined on the same probability space, the four dimensional process $(L^1, R^1, L^2, R^2)$, is a priori (as will follow from the proof described below) a generic four dimensional two-sided Brownian motion with some correlation matrix. The assumption that $L^1_t = L^2_t$ and $R^1_t = R^2_t$ for all $t \in \R$ corresponds to the assumption that this correlation matrix is block diagonal.

The proof of this result is highly technical and we will therefore not cover it here. Instead we refer the reader to Section 9 of \cite{DuplantierMillerSheffield}.

\subsection{Some elements of the proof of Theorem \ref{T:MOT}}
\label{SS:proofMOT}

%\nb{Proof of stationarity. }

\medskip We now have all the tools in hand to begin the proof of Theorem \ref{T:MOT} per se, given the stationarity of Lemma \ref{L:conestationary} and the cutting theorem of Theorem \ref{T:cuttingcone2}.

The proof consists of two fairly distinct steps. 
\begin{itemize}
	\item[\bf Step 1.] Show that $(L_t, R_t)_{t\in \R}$ has stationary and independent increments as well as the Brownian scaling property, and is therefore a two-sided Brownian motion with some covariance matrix. This actually works for all $\kp  >4$ and not just $\kp \in (4,8)$.
	\item[\bf Step 2.] Identify the covariance matrix using the notion of \emph{cone times}. This argument we will present comes from \cite{DuplantierMillerSheffield} and only works for $\kp \in (4,8)$; in the case $\kp \ge 8$ a related but more complicated argument was given separately by Holden, Gwynne, Miller and Sun in \cite{HoldenGwynneMillerSun}.
\end{itemize}

\begin{proof}[Proof of Step 1 for all $\kp>4$] Define a filtration $\cF_t$ by considering, for any $t\in \R$, the sigma-algebra generated by $\eta(s), s\le t$ and $h|_{\eta(-\infty, t)}$. Let $D_t^-$ denote the interior of $\eta(-\infty, t)$ and let $D_t^+$ denote the interior of $\eta(t, \infty)$. Observe that by the cutting theorem (Theorem \ref{T:cuttingcone2}), the doubly marked quantum surfaces
$$\cW_t^\pm = (D^\pm_t, h|_{D_t^\pm}, \eta(t), \infty)$$ 
are quantum wedges of parameter $\alpha = 3\gamma/2$ with $\cW^+_t $ independent of $\cW^-_t$. Indeed, for $t=0$, this follows from the second bullet point in that theorem, and for other values of $t\in \R$, the same can be deduced from the stationarity of the quantum cone viewed from $\eta(t)$ (Lemma \ref{L:conestationary}). Recall that this means that if $g^\pm_t$ is a map from $D_t^\pm$ to $\H$ sending $\eta(t)$ to 0 and fixing $\infty$ with some scaling chosen so that the resulting fields are in the unit circle embedding, then the fields $g^\pm_t(h)$ obtained from $h$ by applying the change of coordinates formula are independent fields in $\H$ (recall also that this does not require considering these fields as being defined modulo additive constant), with laws of a (thick) quantum wedge as specified in Chapter \ref{S:SIsurfaces}. 

Recall also that by Lemma \ref{L:WPstat}, given $\cF_t$ , the curve $(\eta(t+s))_{s\ge 0}$ is just a chordal SLE$_\kp$ in its domain $D_t^+$, from $\eta(t)$ to $\infty$. The key observation is that the increments of $(L,R)$ over $[t, \infty)$ can be described intrinsically in terms of the surface $\cW_t^+$. More precisely, since the boundary length can be computed by conformally changing the coordinates, we can compute the conditional law given $\cF_t$ of the increment $(L_{t+u}- L_t,R_{t+u} - R_t)$ for $u\ge 0$ as follows:

\begin{itemize}

\item[--]
  Take a quantum wedge of the appropriate parameter $\alpha = 3\gamma/2$ embedded in $\H$, and consider a chordal SLE$_\kp$ curve $\eta$ in $\H$ from 0 to $\infty$, reparametrised by Liouville area, and run it for $u$ units of time. 
  
 \item[--]Compute the relative boundary lengths of $\H \setminus \eta (0, u)$ to the left and right of $\eta(u)$, compared those of $\H$, left and right of zero (note that these could be both positive or negative!) 
 
 \end{itemize}

As the reader can see, this description is independent of $\cF_t$ and depends only on $u$. This immediately gives the desired independence and stationarity of the increments.

To conclude Step 1, it remains to check that $(L_t, R_t)_{t\in \R}$ obeys the Brownian scaling property. Namely, if $\lambda>0$, we want to show that 
$$
\frac1{\sqrt{\lambda}} ( L_{\lambda t}, R_{\lambda t})_{t\in \R}
$$
has the same law as $(L,R)$. Informally, this will follow from the fact that the volume (which parametrises $L$ and $R$) scales like $e^{\gamma h}$, while the length, which gives the values of $L$ and $R$, scales like $e^{(\gamma/2) h}$. More precisely, recall that if $C \in \R$, the quantum cone $( \C, h, 0, \infty)$ and $(\C, h+C, 0, \infty)$ have the same laws as quantum surfaces.  
 Let $(L^C, R^C)$ denote the process of left and right boundary lengths associated to to the field $h+C$ along the curve $\eta$ parametrised by $\mu_{h+C}$.  Since $h+C$ and $h$ define the same quantum surfaces in law, and because $\eta$ has the scale invariance property and is independent of $h$, 
 $$(L^C, R^C) \text{ has the same law as the original process $(L,R)$} . $$
 On the other hand, it is clear that $L^C$ can be obtained from $L$ simply by time changing and scaling: more precisely, 
$$
L^C_t = \frac1{\sqrt{\lambda}} L_{\lambda t}
$$
with $\lambda = e^{\gamma C}$, since quantum areas for $h+C$ are multiplied by $\lambda$, and quantum lengths of $h+C$ are multiplied by $\sqrt{\lambda}$. The same holds for $R$ as well, which concludes the proof of Brownian scaling and thus of Step 1. 

\end{proof}

By symmetry of $L$ and $R$, Step 1 implies that (for all $\kp>4$) we can write 

\begin{equation}\label{E:LRtransform}
	\left(\begin{matrix} L_t \\ R_t \end{matrix}\right) = a\left(\begin{matrix} \sin(\theta) & -\cos(\theta) \\ 0 & 1 \end{matrix}\right) \left(\begin{matrix} X_t \\ Y_t \end{matrix} \right)
\end{equation}
where $(X_t,Y_t)_{t\in \R}$ is a standard two-sided planar Brownian motion (started at the origin and with independent coordinates), $a>0$ is such that $\var(L_1)=\var(R_1)=a^2$ and $\theta\in [0,\pi]$ is such that $\cov(L_1,R_1)=-a^2\cos(\theta)$.

\begin{proof}[Proof of Step 2 for $\kp\in(4,8)$.]
This step consists of identifying $\theta$ in \eqref{E:LRtransform}, and the method we present will only work for the case $\kp\in (4,8)$, as will become clear shortly. This range of $\kp$ corresponds to $\theta\in (\pi/2,\pi)$, equivalently, $-\cos(\theta)=a^{-2}\cov(L_1,R_1)>0$. That is, the case where $L$ and $R$ are positively correlated. 

The argument we present will use the notion of \emph{cone times}. We say that $t$ is a local $\theta$-cone time for a process $(X_s,Y_s)_{s\in \R}$ if there exists $\eps>0$ such that $(X_s,Y_s)$ remains in the set $(X_t+Y_t)+C_\theta$ for all $s\in [t,t+\eps]$, where $C_\theta=\{z\in \C: \arg(z)\in [0,\theta]\}$. It is straightforward to see that if $(L,R)$ and $(X,Y)$ are related by \eqref{E:LRtransform}, then the set of local $\pi/2$-cone times for $(L,R)$ correspond precisely to the set of local $\theta$-cone times for $(X,Y)$. 

The key idea is to identify $\theta$ using a result of Evans, \cite{Evanscones}, which states that the almost sure Hausdorff dimension of the set of local $\theta$-cone times of $(X,Y)$ is equal to $0$ for $\theta\in [0,\pi/2]$, and equal to $1-\pi/2\theta$ for $\theta\in (\pi/2,\pi)$. On the other hand, as we will explain below, it is possible to compute the almost sure Hausdorff dimension of the set of local $\pi/2$-cone times for $(L,R)$, using the definition of $(L,R)$ in terms of space-filling SLE on a quantum cone. This will be equal to $0$ when $\kappa'\ge 8$, and $1-\kappa'/8$ when $\kappa'\in (4,8)$. Hence, we learn nothing if $\kappa'\ge 8$, but for $\kappa'\in (4,8)$ we see that necessarily: 
\[ 1-\frac\kp8=1-\frac{\pi}{2\theta}, \text{ equivalently } \theta=\frac{4\pi}{\kp},\] as required.

So, it remains to argue that when $\kp\in (4,8)$ the Hausdorff dimension of the set of local $\pi/2$-cone times of $(L,R)$ is almost surely equal to $1-\kp/8$. In fact, since $(L_s,R_s)_{s\in \R}$ and $(\hat{L}_s,\hat{R}_s)_{s\in \R}:=(L_{-s},R_{-s})_{s\in \R}$ are identical in law, it suffices to consider the set $\mathcal{A}$ of local $\pi/2$-cone times for $(\hat{L},\hat{R})$ and show that the Hausdorff dimension of $\mathcal{A}$ is almost surely equal to $1-\kp/8$.

We are going to identify the Hausdorff dimension of the set $\mathcal{A}$ with the Hausdorff dimension of a set of times determined by the geometry of $\eta'$. To set up for this, first notice that by the definition of $\pi/2$-cone times, if $s$ is a local $\pi/2$-cone time for $(\hat{L},\hat{R})$, then 
there exists some $t\in \mathbb{Q}$ with $(\hat{L}_r,\hat{R}_r)\in (\hat{L}_s, \hat{R}_s)+C_{\pi/2}$ for all $r\in (s,t)$. This implies that  $(L_{-t+u}-L_{-t},R_{-t+u}-R_{-t})_{u\ge 0}$ has a simultaneous running infimum at $u=t-s$. Conversely, for any $t\in \mathbb{Q}$, each simultaneous running infima of $(L_{-t+u}-L_{-t},R_{-t+u}-R_{-t})_{u\ge 0}$ corresponds to a local $\pi/2$ cone time for $(\hat{L},\hat{R})$. Thus, we can write $\mathcal{A}=\cup_{q\in \mathbb{Q}} \{-q-\mathcal{A}_q\}$, where $\mathcal{A}_q $ is the set of simultaneous running infima of $(L_{q+u}-L_q,R_{q+u}-R_q)_{u\ge 0}$ . Notice also that by the stationarity of Lemma \ref{L:conestationary},  the almost sure Hausdorff dimension of $\mathcal{A}_q$ does not depend on $q\in \mathbb{Q}$. In particular, this implies that the Hausdorff dimensions of $\mathcal{A}$ and $\mathcal{A}_0$, say, are equal almost surely.  We claim that
\begin{equation}\label{E:conetimestocurve}
	\mathcal{A}_0=\{s>0: \eta'(s)\in \eta_0^L\cap \eta_0^R\}
\end{equation}
with probability one.
To see this, recall from Section \ref{SS:sf48} that $\eta_0^L$ and $\eta_0^R$ are the concatenation left and right boundaries of an ordered collection of simply connected domains, ordered consistently with the ordering of $(\eta'(r))_{r<0}$. But we can also reverse this ordering, and view $\eta_0^L$ and $\eta_0^R$ as the boundaries of an ordered collection of simply connected domains forming $\C\setminus \eta'((-\infty,0])=\eta'((0,\infty))$, ordered according to when they are visited by $(\eta'(s))_{s>0}$. 
%The intervals of time on which $\eta'(r)\notin \eta_0^L\cap \eta_0^R$ correspond precisely to the intervals during which $\eta'|_{[0,\infty)}$ visits one of these domains. Note that the end point of such an interval $(r_1,r_2)$ does correspond to a simultaneous running infima of $(L_r,R_r)$ since $\eta_{\eta'(r_2)}^L\subset \eta_{\eta'(r)}^L$ for all $r\in (0,r_2)$ by definition, with the same holding when $L$ is replaced by $R$. Moreover, 
It is not hard to convince oneself that $L$ has a running infimum at time $r>0$ if and only if $\eta'(r)=z\in \eta_0^L$, and $r=\sup\{u: \eta'(u)=z\}$. The analogous statement holds when $L$ is replaced by $R$. Thus, $(L,R)$ have a simultaneous running infima at $s>0$ if and only if $\eta'(s)=z\in \eta_0^L\cap \eta_0^R$ and $s=\sup\{u: \eta'(u)=z\}$. It also follows from the definitions that points of $\eta_0^L \cap \eta_0^R$ are visited exactly once by $\eta'$, and so in this case $(L,R)$ have a simultaneous running infima at $s>0$ if and only if $\eta'(s)=z\in \eta_0^L\cap \eta_0^R$. 

Next, we recall the statement of Theorem \ref{T:cuttingcone2}. This says that if $(\C,h,0,\infty)$ is the $\gamma$-quantum cone used to define $(L,R)$, then viewed as a quantum surface, $(\eta'[0,\infty]), h, 0, \infty)$ has the law of a quantum wedge with parameter $\alpha=3\gamma/2$. This is a thin wedge for $\kp\in (4,8)$, meaning that it is an ordered collection of quantum surfaces, and this ordered collection of surfaces correspond precisely to the ordered collection of simply connected domains described in the above paragraph. 
The intervals of time on which $\eta'(r)\notin \eta_0^L\cap \eta_0^R$ correspond precisely to the intervals during which $\eta'|_{[0,\infty)}$ visits one of these domains. By definition of the parametrisation of $\eta'$, the lengths of these intervals are exactly the quantum areas of the quantum surfaces making up the thin wedge. It therefore follows from Lemma \ref{L:volumepushforwarddisc} (and an additional scaling property together with the finiteness of a certain moment, see Proposition 4.4.4 in \cite{DuplantierMillerSheffield}) %and Proposition \ref{prop:unitareaQD} 
that the ordered collection of lengths of these intervals are equal in law to the durations of excursions away from $0$, for a Bessel process of dimension $\delta=\kp/4$. It turns out that the finite moment assumption boils down to requiring that $\delta >$, which fortunately is the case when $\delta = \kp/4$ and $\kp >4$. 
By classical excursion theory, see for instance \cite[Chapter III]{BertoinLevy}, it follows that  $\mathcal{A}_0$ is the range of a $1-\delta/2$ stable subordinator, and the Hausdorff dimension of such a set is almost surely equal to $1-\delta/2=1-\kp/8$.
\end{proof}

The concludes the proof of Theorem \ref{T:MOT}.

%\subsection{Application: CRT-mated maps \ellen{*to be added}}\label{SS:CRTmm}

%\subsection{Mating of trees in other topologies \ellen{*to be added}}
\label{SS:extension}

%\newpage 
%
%\section{Critical Gaussian multiplicative chaos \ellen{*to be added}}\label{S:cGMC}

%\input{CriticalGMC.tex}

%\newpage 

%\section{Liouville Brownian motion \ellen{*to be added}}\label{S:LBM}

%\input{LBM.tex}

\newpage

\appendix

% !TEX root = master.tex

\section{Chordal Loewner chains and chordal SLE}\label{app:sle}

The aim of this appendix is to collect some relevant background material on Schramm--Loewner evolutions (SLE), primarily to accompany \cref{S:zipper,S:MOT}. For a much more detailed and pedagogical exposition, the reader is referred to \cite{SLEnotes, Kem17, Lawlerbook}. The presentation here most closely follows \cite{SLEnotes}. 

\subsection{Chordal Loewner chains} \ind{Loewner chain}

\paragraph{Complex analysis basics.} First, we fix some basic notation and terminology.
\begin{itemize}
	\item  $K\subset \H$ is said to be a \emph{compact $\H$-hull} if it is bounded and $H:=\H\setminus K$ is a simply connected domain.
	\item For any such hull, by the Riemann Mapping Theorem, one can choose a conformal isomorphism $g_K:H\to \H$ such that $g_K(z)-z\to 0$ as $z\to \infty$. In fact, one can prove that for this $g_K$, the expansion $ g_K(z)=z+\frac{a_K}{z} + O(|z|^{-2}) $ holds as $z\to \infty$ for some $a_K\geq 0$. We call $g_K$ the \emph{Loewner map} of $K$.
	\item  $a_K$ is known as the \emph{half plane capacity} of $K$ and denoted by $\text{hcap}(K)$.
	\item In some sense, the half plane capacity measures the size of the hull $K$, when ``viewed from infinity". In particular, the half plane capacity increases as a hull increases: if $K\subset K'$ are two complex $\H$ hulls, then $\text{hcap}(K)\leq \text{hcap}(K')$.
\end{itemize}

\paragraph{Loewner Chains.} A Loewner chain is a family $(K_t)_{t\geq 0}$ of \emph{increasing} ($K_s\subsetneq K_t$ for $s\leq t$) compact $\H$-hulls which satisfy a \emph{local growth property}: for any $T\geq 0$, $$\sup_{s,t\in [0,T], |s-t|\leq h} \textrm{rad}\left( g_{K_s}(K_t\setminus K_s)\right) \to 0 \;\; \text{as} \;\; h\to 0.$$ Here the radius of a hull means the radius of the smallest semicircle in which it can be inscribed. For such a chain one can show that the half plane capacity is a strictly increasing bijection from $[0,\infty) \to [0,\infty)$, so we can always assume (by convention) that time is parametrised so that $\text{hcap}(K_t)=2t$ for all $t$.
	
	\begin{theorem}[Loewner's theorem] Loewner discovered that such chains (parametrised by half plane capacity) are in bijection with continuous real valued functions via the following correspondence.
		\begin{itemize}
			\item Given $(K_t)_{t\ge 0}$ a Loewner chain, there is a unique point $\xi_t\in \overline{\cap_{h>0}g_{K_t}(K_{t+h}\setminus K_t)}$ for each $t\ge 0$. $(\xi_t)_{t\ge 0}$ is a continuous real valued function called the \emph{driving function of $(K_t)_{t\ge 0}$}.
			\item Given $(\xi_t)_{t\geq 0}$ a continuous real valued function, define, for each $z\in \H$, $g_t(z)$ to be the maximal solution to the \emph{Loewner equation}
			\begin{equation}\label{eq:forward_Loewner} \frac{\partial g_t(z)}{\partial t} = \frac{2}{g_t(z)-\xi_t}, \;\; g_0(z)=z \end{equation}
			which exists on some time interval $[0, \zeta(z)]$ by classical ODE theory. Let $K_t=\{z\in \H: \zeta(z)\leq t\}$.
			Then $(K_t)_{t\geq 0}$ is a Loewner chain with driving function $\xi_t$. Moreover, $g_t=g_{K_t}$ for all $t$.
		\end{itemize}
\end{theorem}

 \begin{figure}
 	\begin{center}	\includegraphics[scale=.75]{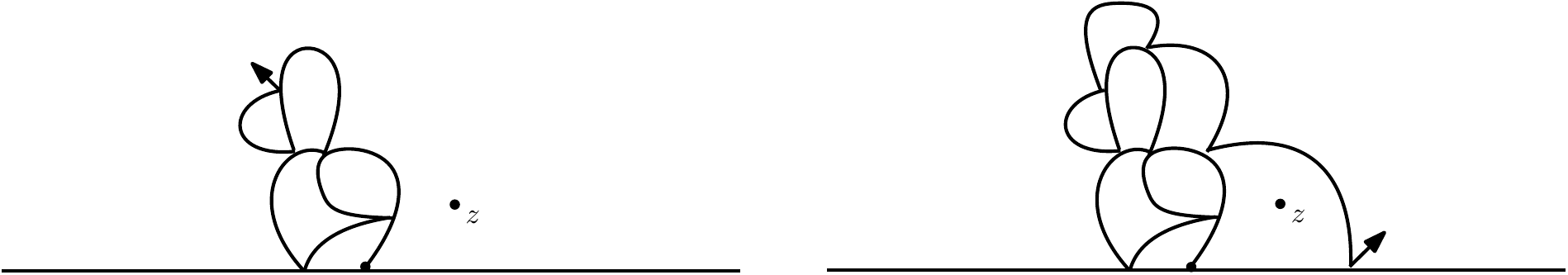} \end{center}
 	\caption{A Loewner chain drawn up to two times: on the left, a time before $\zeta(z)$, and on the right, just after $\zeta(z)$.}\label{fig:ls}
 \end{figure}

\noindent We call $(g_t)_{t\ge 0}$ the (forward) \emph{Loewner flow}. $\zeta(z)$ is the time that the growing hull $K_t$ ``swallows'' the point $z$. See \cref{fig:ls}.

\begin{rmk}
\label{rmk:LEcurve} Continuous curves $(\gamma(t))_{t\ge 0}=:(\gamma_t)_{t\geq 0}$ in $\H$ which do not cross themselves and have $|\gamma_t|\to \infty$ as $t\to \infty$ provide examples of Loewner chains. More precisely, one defines $H_t=\H\setminus K_t$ for each $t$ to be the connected component of $\H\setminus \gamma([0,t])$ containing $\infty$. Then $(K_t)_{t\ge 0}$ is an increasing family of compact $\H$-hulls satisfying the local growth property. The map $g_t$ sends the tip of the curve, $\gamma_t$, to the point $\xi_t$ (where $g_t$ is extended by ``continuity'', or more precisely this holds in the sense of prime ends).\end{rmk}

\subsection{Chordal \texorpdfstring{$\mathrm{SLE}_\kappa$}{TEXT}} \ind{SLE! Chordal} Chordal SLE$_\kappa$ processes, for $\kappa>0$\indN{{\bf Parameters}! $\kappa$; SLE parameter}, were introduced by Oded Schramm \cite{Sch00} as a family of potential scaling limits for interfaces in critical statistical physics models. As we will soon see, they satisfy two very natural properties that make them appropriate candidates for such limits: conformal invariance and a certain domain Markov property.

It turns out  (\cite{Sch00}) that these two properties actually \emph{characterise} SLE$_\kappa$ as a one parameter family, which means that there really can be no other candidates. On the other hand, proving convergence of discrete interface models to SLE is typically very challenging. To date it has been verified for just a few special values of $\kappa$; for example, critical percolation interfaces, \cite{Smirnov}, and the loop-erased random walk, \cite{LSW}.

\begin{definition}[Chordal SLE in $\H$ from $0\to \infty$]
	\label{defn::chordal_sle_H}
	For $\kappa>0$, SLE$_\kappa$ in $\H$ from $0$ to $\infty$ is defined to be the Loewner chain driven by $\xi_t=\sqrt{\kappa}B_t$ where $B_t$ is a standard Brownian motion. 	
\end{definition}

One of the first things to note about SLE is that, due to the scaling property of Brownian motion ($B_t$ has the same law as $\sqrt{t}B_1$ for any $t$), SLE  is itself scale invariant. That is, for any $r\geq 0$ if $(K_t)_{t\geq 0}$ is an SLE$_\kappa$ process, then the rescaled process $(r^{-1/2} K_{rt})_{t\geq 0}$ also has the law of an SLE$_\kappa$. This says that SLE is invariant under conformal isomorphisms of $\H$ that fix $0$ and $\infty$. This allows us to define SLE, by conformal invariance, in any simply connected domain and between any two marked boundary points.

\begin{definition}[Chordal SLE]
	\label{defn::chordal_sle}
	$SLE_\kappa$ is a collection $(\mu_{D,a,b})_{D,a,b}$ of laws on Loewner chains, indexed by triples $(D,a,b)$ where $D$ is a simply connected domain and $a$ and $b$ are two marked boundary points. The law $\mu_{\H,0,\infty}$ is that given by \cref{defn::chordal_sle_H}. For any other triple $(D,a,b)$, $\mu_{D,a,b}$ is defined to be the image of $\mu_{\H,0,\infty}$ under (any) conformal isomorphism sending $\H$ to $D$, $0$ to $a$ and $\infty$ to $b$.
\end{definition}

The choice of conformal isomorphism above does not impact the law of the curve up to time-change, since any two such conformal isomorphisms differ by scaling and SLE$_\kappa$ in the upper-half plane from 0 to $\infty$ is invariant under scaling. See \cite{SLEnotes} for a somewhat more rigorous measure-theoretic framework. 

\paragraph{Chordal SLE: properties.}
\begin{itemize}
	\item Chordal SLE$_\kappa$ is generated by a curve $\gamma$ (in the sense of \cref{rmk:LEcurve}) for every $\kappa>0$: due to \cite{RohdeSchramm} for $\kappa\ne 8$, and \cite{LSW} for $\kappa=8$.
	\item \emph{Conformal invariance:} if $\gamma$ is an SLE$_\kappa$ in $D$ from $a$ to $b$ and $\psi:D\to D'$ is a conformal isomorphism with $\psi(a)=a'$ and $\psi(b)=b'$, then $\psi(\gamma)$ (up to time reparametrisation) has the law of an SLE$_\kappa$ in $D'$ from $a'$ to $b'$.
	\item \emph{Domain Markov property:} if $\gamma$ is an SLE$_\kappa$ from $a$ to $b$ in $D$ and $T$ is a bounded stopping time that is measurable with respect to $\gamma$, then conditionally on $\gamma([0,T])$, writing $D_T$ for the connected component containing $b$ of $D\setminus \gamma([0,T])$, $\gamma([T,\infty))$ has the law of an SLE$_\kappa$ from $\gamma(T)$ to $b$ in $D_T$.
	\item It has three distinct \emph{phases}: for $\kappa\in [0,4]$ SLE$_\kappa$ is almost surely generated by a simple (non-self touching and non-boundary touching) curve; for $\kappa\in (4,8)$ it almost surely hits (but doesn't cross) itself and the boundary of the domain; and for $\kappa\ge 8$ it is almost surely space filling. See \cref{fig:sle}.
\end{itemize}

\begin{figure}
	\begin{center} \includegraphics[scale=0.4]{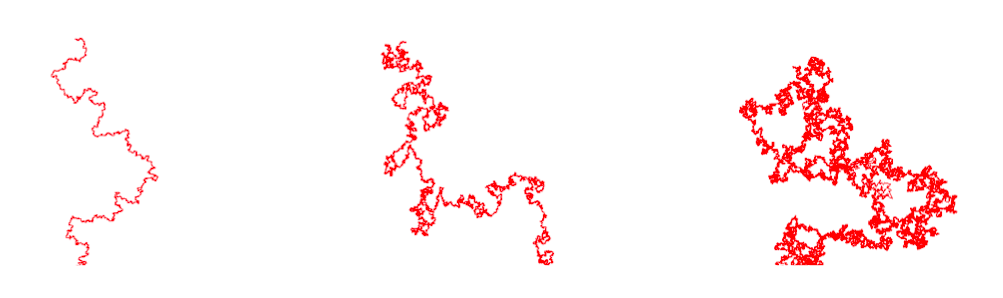} \end{center}
	\caption{From left to right: SLE$_2$, SLE$_4$, SLE$_6$. Simulations by Tom Kennedy.}\label{fig:sle}
\end{figure}

\subsection{Chordal \texorpdfstring{$\mathrm{SLE}_\kappa(\underline{\rho})$}{TEXT}} \ind{SLE! Chordal with force points}
It is best to view chordal SLE$_\kappa$ as a family of laws $\mu_{(D, a, b)}$ on random \emph{chords} in the domain $D$ connecting one boundary point $a$ to another boundary point $b$ (where boundary is understood in a conformal sense). Both the conformal invariance and domain Markov properties of chordal SLE are then easily formulated through this notation: for instance, the requirement of conformal invariance is that $\mu_{(\phi(D); \phi(a) , \phi(b))}$ is the push-forward of the measure $\mu_{D, a, b}$ by the conformal isomorphism $\phi$.

 It is also very natural to consider random curves in which domain Markov property and conformal invariance are satisfied only provided we specify additional information such as the location of a specified number of points in the domain or on its boundary. For concrete examples, consider the scaling limits of discrete interface models in which there is a change of boundary conditions at a specified number of points along the boundary: the law of the scaling limit will depend on the location of these special points.  
 %As we will see below (in the case of one boundary point, but the argument can probably be generalised to more marked points), such scaling limits are necessarily of the form SLE$_\kappa( \underline\rho)$, which we are about to define. 

As it turns out (see Remark \ref{rmk:bessel_sle} for a proof in the case of one marked point on the boundary), such curves are described by variants of SLE$_\kappa$, which have an additional attraction or repulsion from certain \emph{marked points} (also sometimes known as force points) in the domain or on its boundary. These are known as SLE$_\kappa(\underline\rho)$ and first appeared in \cite{conformalrestriction}; see also \cite{SW05,MSIG1}. The vector $\underline\rho$ encodes how strong this attraction or repulsion is, and in which direction.

Let $\kappa >0$. We will take again the upper half plane as a reference domain and $a = 0, b = \infty$ for the start and target points on the boundary. Let $v^1, \ldots, v^m$ be $m$ marked points on the boundary (we will discuss interior points below) and corresponding weights $\rho^1,..., \rho^m\in \R$ are such that \begin{equation} \label{eqn:rho_nice} \sum_{i\in S} \rho^i \ge -2 \text{ for every } S\subset\{1,..., M\}\end{equation} 
To define the law $\mu_{(\H, 0, \infty); (v^1, \ldots, v^m)}$  of SLE$_\kappa(\underline{\rho})=$ SLE$_\kappa(\rho^1, ..., \rho^m)$ we proceed as follows.  It will be a Loewner chain 
%that is almost surely generated by a curve and can be defined for all time. \footnote{In fact, this holds under a weaker assumption on $\underline{\rho}$ (see for example \cite{cleperc}) but some additional complications arise that we will not address in this brief overview.}
and hence
by Loewner's theorem, can be defined by specifying its driving function. As for ordinary SLE$_\kappa$, this driving function will be a random function closely related to Brownian motion.  However, the Brownian motion now comes with a \emph{drift}. This drift will 
depend on the position of the marked points $(V_t^1,...,V_t^m)_{t\ge 0}$ after applying the Loewner flow $(g_t)_{t \ge 0}$ as follows.

\begin{definition}[SLE$_\kappa(\rho^1,..., \rho^m)$ in $\H$ from $0$ to $\infty$]\label{D:chordalSLEkp1}
	Suppose that $v^1,..., v^m\in (\R\cup \{\infty\}) \setminus \{0\}$ are distinct and $\rho^1,..., \rho^m$ satisfy the condition \eqref{eqn:rho_nice}. SLE$_\kappa(\rho^1,..., \rho^m)$
	with marked points at $v^1, ..., v^m$ is  the Loewner chain with driving function $(\xi_t)_{t\ge 0}$ satisfying the following system of SDEs:
	\begin{eqnarray}\label{eqn:slekpdef}
	\xi_t & = & \sqrt{\kappa} B_t + \sum_i \int_0^t \frac{\rho^i}{\xi_s-V_s^{i}} \, ds  \nonumber \\
	V_t^{i} & = & v^i + \int_0^t \frac{2}{V_s^i - \xi_s} \, ds \text{ for } 1\le i \le M.
	\end{eqnarray}
\end{definition}

The second equation in \eqref{eqn:slekpdef} is simply Loewner's equation describing the evolution $V_t^i$  of the marked point $v^i$ under the Loewner flow, at least until $v^i$ is swallowed by the chain (in fact, the evolution can be extended beyond this point). The first equation in \eqref{eqn:slekpdef} describes the driving function of the Loewner flow as usual.

\medskip For any value of $\underline{\rho}$, the strong existence and uniqueness of solutions to \eqref{eqn:slekpdef} is clear until the first time that one of the $V^i$ collides with $\xi$, that is, $\sup\{t\ge 0: \xi_t\ne V_t^i \; 1\le i \le m\}$. This is also the first time that one the marked points $v^i$ is swallowed by the hull generated by the Loewner chain. 
%When $\underline{\rho}$ satisfies \eqref{eqn:rho_nice}, this evolution (meaning the existence of unique solutions) can be extended to show,
In fact, it can be shown (see \cite{MSIG1}) that when $\underline{\rho}$ satisfies \eqref{eqn:rho_nice}, there exists an almost surely continuous Markovian process $(\xi,V^1,\dots, V^m)$ that satisfies the integrated equation \eqref{eqn:slekpdef} for \emph{all} time, and for which the set of times $t$ with $\xi_t=V_t^i$ for some $i$ almost surely has Lebesgue measure $0$. It is also shown in \cite{MSIG1} that the law of this process is unique. Consequently, the corresponding chordal SLE$_\kappa(\underline{\rho})$ Loewner chain in \cref{D:chordalSLEkp1} is well defined (for all time).

\begin{rmk}\label{rmk:bessel_sle}
In the case of one marked point on the boundary ($m=1$), the process $(V^1_t-\xi_t)_{t\ge 0}$ describing the distance between the driving function and the evolution of the marked point, is $\sqrt{\kappa}$ times a Bessel process. When $\underline{\rho}=\rho^1=\rho$, the dimension of the Bessel process is \begin{equation}\label{dimBessel}
\delta = 1+\frac{2(\rho+2)}{\kappa}.
\end{equation}
 This formalises the notion that SLE$_\kappa(\underline{\rho})$ processes have an additional attraction/repulsion from the marked points.

In fact, for a Loewner chain to satisfy a conformal Markov property with an extra marked point (that is, the property that for any stopping time $\sigma$, the future evolution after applying the Loewner map at time $\sigma$ has the same law as the original process, with the marked point now located at the image of the original marked point) one finds that the difference between the driving function and the evolution of the marked point must be a continuous Markov process satisfying Brownian scaling. This implies that it actually has to be a Bessel process of some dimension. One can take this an explanation for the form of the SDEs \eqref{eqn:slekpdef}.
\end{rmk}

\begin{figure}
	\begin{center} \includegraphics[scale=0.75]{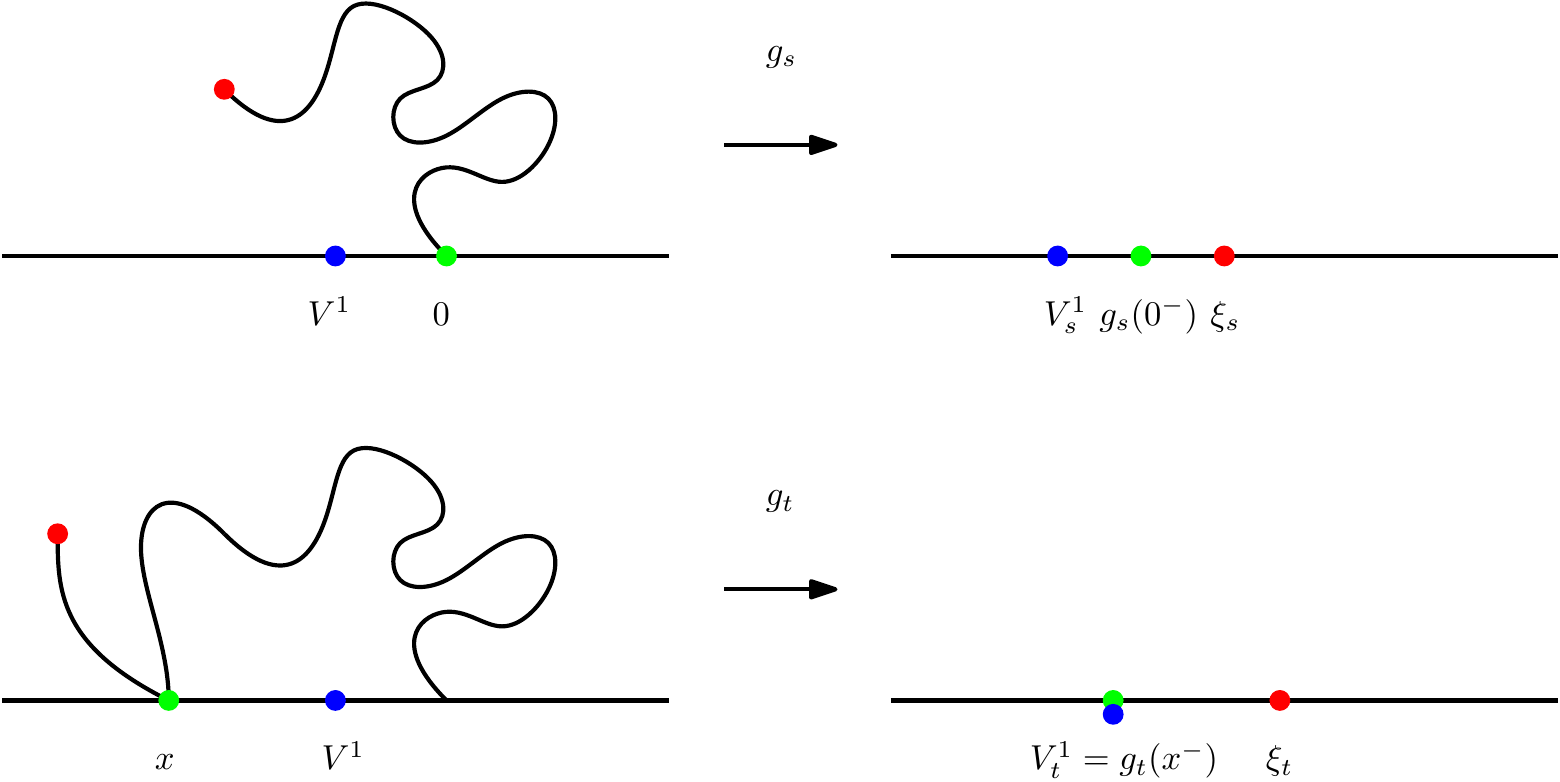} \end{center}
	\caption{A schematic picture of an SLE$_\kappa(\rho)$ with one marked point, drawn up to two times $s,t$ with $s<t$. At time $s$ the marked point has not been swallowed, but at time $t$ it has. After time $t$ the evolution of $V^t_1$ coincides (by definition) with the evolution under the Loewner flow of the point infinitesimally to the left of $x$.} \label{F:FP}
\end{figure}

\begin{rmk}
	The definition can also be extended to the case where there are marked points located infinitesimally to the left and/or right of $0$ (denoted $0^-$ and $0^+$). This is done by taking a limit in law (with respect to the Carath\'{e}odory topology on Loewner chains\footnote{this is the topology for which a sequence of chordal Loewner chains (from $0$ to $\infty$ in $\H$) with Loewner flow $(g^{n}_t)_{t\ge 0}$ converges to a Loewner chain with Loewner flow $(g_t)_{t\ge 0}$ as $n\to \infty$ if and only if $(g_t^{n})^{-1}(z)\to g_{t}^{-1}(z)$ uniformly on compact subsets of time and subsets of space that are compactly contained in $\H$} as one of the marked points approaches $0$ from the left and/or one of the marked points approaches $0$ from the right. Again this gives rise to unique laws on Loewner chains that are defined for all time.

When there is just one marked point, this boils down to starting a Bessel process of positive dimension from zero; in fact by \eqref{dimBessel} the dimension of this Bessel process is greater than 1 when $\rho > -2$. The reason why we assume the dimension to be greater than 1 (and so $\rho$ to be $> -2$) is to ensure that the integral in \eqref{eqn:slekpdef} is convergent. When $\rho\le-2$, assigning a meaning to this integral is less straightforward, though there are known procedures, including for example a principal value correction, see \cite{ExpTree}.
\end{rmk}

\begin{rmk}
	\cref{D:chordalSLEkp1} can also be extended to include \emph{interior force points}. That is, with some of the $v^i=V^i_0$ located in $\H$ rather than in $\R$. The definition is exactly the same, but in this case, existence and uniqueness of solutions to \eqref{eqn:slekpdef} is only guaranteed until the first time that $\xi_t=V_t^i$ for some $i$ such that $v^i\in \H$. As such, the chordal SLE$_\kappa(\underline{\rho})$ with interior force points is a well defined random Loewner chain, but only up to the first time that one of the interior force points is ``swallowed''.
\end{rmk}

Due to the scaling property of Brownian motion, it follows easily that SLE$_\kappa(\underline{\rho})$ from $0$ to $\infty$ in $\H$ also satisfies a form of scale invariance. More precisely, if $(K_t)_{t\ge 0}$ is an SLE$_\kappa(\underline{\rho})$ process with force points at $v^1,..., v^m$, then the rescaled process $(r^{-1/2}K_{rt})_{t\ge 0}$ has the law of an SLE$_\kappa(\underline{\rho})$ process with force points at $r^{-1/2}v^1,..., r^{-1/2}v^m$ for any $r>0$. This allows us to extend the definition of SLE$_\kappa(\underline{\rho})$ to arbitrary domains with finitely many marked boundary points.

\begin{definition}[SLE$_\kappa(\rho^1,..., \rho^m)$ in $D$ from $a$ to $b$]\label{D:chordalSLEkp2}
	Suppose that $\rho^1, ... , \rho^m$ are as in  \cref{D:chordalSLEkp1} and $(D,a,b, v^1,..., v^m)$ is a given domain with $(m+2)$ marked points. Let $\psi:\H\to D$ be a conformal isomorphism sending $a$ to $0$ and $b$ to $\infty$.
	
	SLE$_\kappa(\rho^1,..., \rho^m)$ from $a$ to $b$ in $D$
	with marked points at $v^1, ..., v^m$ is defined to be the image under $\psi$ of SLE$_\kappa(\rho^1,..., \rho^m)$ from $0$ to $\infty$ in $\H$, with marked points at $\psi^{-1}(v^1),... ,\psi^{-1}(v^m)$.
	
	This definition also extends to the case of interior force points, with both of the above SLE$_\kappa(\underline{\rho})$ curves being defined up to the first time that an interior force point is swallowed.
\end{definition}

\begin{rmk}[Properties]\label{R:slekp}
SLE$_\kappa(\underline{\rho})$ possesses many properties similar to those of SLE$_\kappa$, along with some additional features.
\begin{itemize}
	\item For any $\kappa>0$ and $\underline{\rho}$ satisfying \eqref{eqn:rho_nice}, SLE$_\kappa(\underline{\rho})$ is almost surely generated by a continuous curve $\gamma$, with $\gamma(0)=a$ and $\gamma(t)\to b$ as $t\to \infty$: see \cite{MSIG1}. 
	\item By definition, if $\psi:D\to D'$ is a conformal isomorphism sending $$(a,b,v^1,..., v^m) \text{ to }  (a',b',(v^1)',...,(v^m)'),$$ then the image of SLE$_\kappa(\underline{\rho})$ from $a$ to $b$ in $D$ with force points at $v^1,\cdots, v^m$ has the law of SLE$_\kappa(\underline{\rho})$ from $a'$ to $b'$ in $D'$ with force points at $(v^1)',...,(v^m)'$.
	\item Going back to the set up in the upper half plane, the processes $(V_t^i)_{t\ge 0}$ from \eqref{eqn:slekpdef} describe the evolution of the force points $v^i$ under the Loewner flow. More precisely, for each $i$ and until the first time $\tau^i$ that $v^i$ is ``swallowed'' by the curve, $V^i_t$ is equal to $g_t(v^i)$ (where $g_t$ is continuously extended to the boundary if necessary). After this time, if $v^i\in \R_+$ (respectively $\R_-$), $V^i_t$ will be equal to the image under $g_t$ of the furthest right (resp. furthest left) point on the real line that has been swallowed at time $\tau^i$. See \cref{F:FP}.
	\item As we have mentioned already, SLE$_\kappa(\underline{\rho})$ satisfies a domain Markov property, that now involves the marked points. To state this precisely, suppose that $\gamma$ is an SLE$_\kappa(\underline{\rho})$ from $a$ to $b$ in $D$ with force points  at $v^1,...,v^m$ and that $T$ is a bounded stopping time for $\gamma$.
	Write $D_T$ for the connected component of $D\setminus \gamma([0,T])$ containing $b$. Then conditionally on $\gamma([0,T])$, $\gamma([T,\infty))$ has the law of an SLE$_\kappa(\underline{\rho})$ from $\gamma(T)$ to $b$ in $D_T$, with force points at $V_T^1,...,V_T^m$.
	\item By inspecting \eqref{eqn:slekpdef}, it follows that for SLE$_\kappa(\rho)$ in $\H$, putting any weight $\rho$ at the boundary point $\infty$ does not affect the law of the curve. This observation will be useful when studying the relationship between chordal and \emph{radial} SLE. 
\end{itemize}
\end{rmk}

Recall that chordal SLE$_\kappa$, $(\gamma(t))_{t\ge 0}$ has three distinct phases. In terms of its interaction with the boundary $\partial \H=\R$, this can be described as follows: if $\kappa\in [0,4]$, $\gamma([0,\infty))\cap \R = \emptyset$ almost surely;  if $\kappa\in (4,8)$, $\gamma([0,\infty))\cap \R$ is almost surely non-empty, unbounded but has Lebesgue measure $0$; and if $\kappa\ge 8$, $\gamma([0,\infty))\cap\R=\R$ almost surely. In the case of SLE$_\kappa(\underline{\rho})$, where there is additional attraction or repulsion from force points on $\R$, this behaviour may be modified. 

Indeed, consider the case of SLE$_\kappa(\rho)$ with one force point at $v^1\in \R$ of weight $\rho$. Then we have already seen that the distance between the driving function and the evolution of $v^1$ under the Loewner flow, $(V_t^1-\xi_t)_{t\ge 0}$, is a Bessel process of dimension $1+2\kappa^{-1}(\rho+2)$. This means that if $\rho\ge \tfrac{\kappa}{2}-2$ then $(V_t^1-\xi_t)$ will almost surely be positive for all $t>0$, that is, the SLE$_\kappa(\rho)$ will almost surely not hit the half closed interval between $v_1$ and $\infty$ at any time $t>0$. If $\rho<\tfrac{\kappa}{2}-2$ then the Bessel process will hit $0$, which means that the SLE$_\kappa(\rho)$ \emph{will} hit this half closed interval. We have proved the following lemma:

\begin{lemma}\label{L:rhohit}
Let $\eta $ be a chordal SLE$_\kappa(\rho)$ with $\kappa>0$ and $\rho>-2$ for some boundary marked point $v$. Then $\eta$ hits $v$ (or more precisely $v$ is swallowed by the hull generated by $\eta$) if and only if $\rho < \kappa/2 - 2$.
\end{lemma}

For multiple force points, a description of the interaction of SLE$_\kappa(\underline{\rho})$ with the real line (depending on $\underline{\rho})$) can be found in \cite[Lemma 15]{Dub3}.

\section{Reverse Loewner flow and reverse SLE} \ind{Reverse Loewner flow}

\subsection{Definitions}

Until now this appendix has focused on standard Loewner evolutions, describing increasing families of compact hulls: in nice cases, growing curves. However, these should really be referred to as \emph{forward} Loewner evolutions, because they also have a counterpart: \emph{reverse Loewner evolutions}.  A reverse Loewner evolution is no longer a family of hulls that increases in time, but rather a family of hulls where in each infinitesimal increment of time, an infinitesimal new piece of hull is added ``at the root''. The whole of the previous hull is then conformally mapped to something slightly different (one might envisage the new piece of hull as ``pushing'' the existing one further into the domain). See \cref{fig:reverseLE}. Note that one cannot therefore speak of a ``single curve'' associated to a reverse Loewner evolution. \\

\begin{figure}
	\begin{center} \includegraphics[scale=0.75]{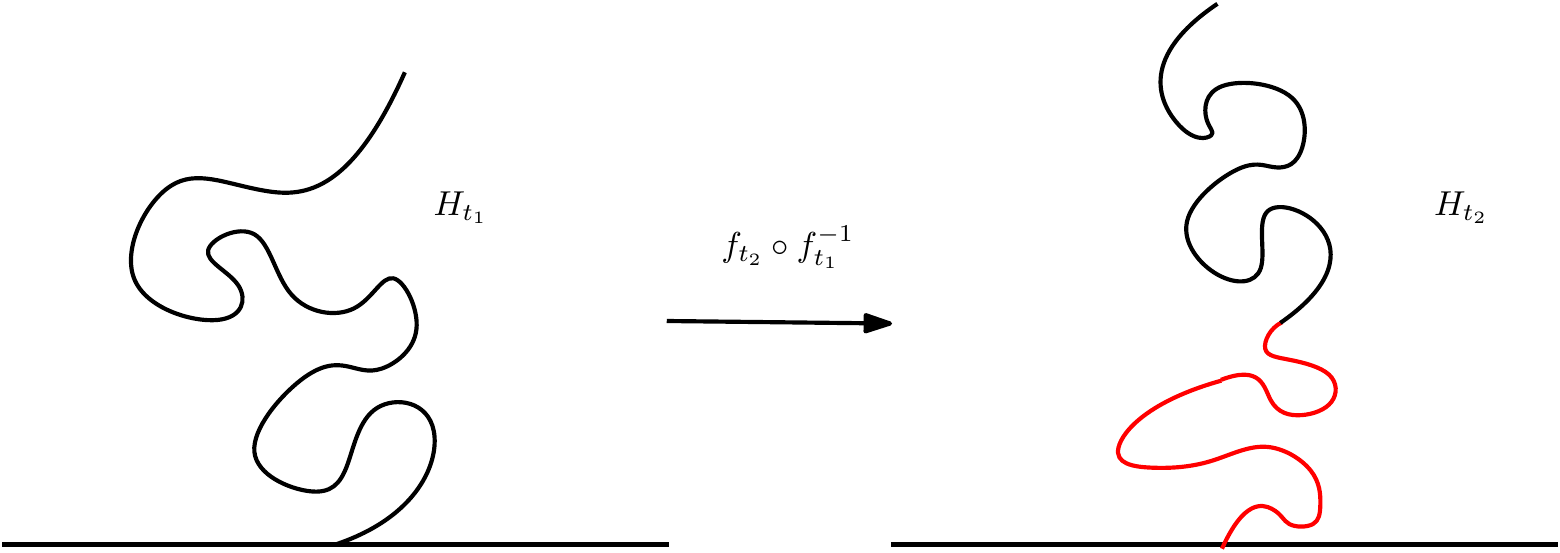} \end{center}
	\caption{A reverse Loewner evolution at two times $t_1,t_2$ with $t_1<t_2$. In both cases $H_{t_i}$ is the complement of the curve. One can see that between the two times, a new piece of curve (drawn in red) is added ``at the root'', and the existing curve (black) is conformally mapped into the domain formed by the complement of the red curve. In contrast, under the forward Loewner flow, new pieces of curve are always added ``at the tip'' of the existing curve. }\label{fig:reverseLE}
\end{figure}

In the following, we will only ever discuss \emph{centred} reverse Loewner evolutions. Informally, this means that new pieces of curve are always added at the origin.

\begin{definition}[Reverse Loewner evolution in $\H$]\label{def:RLE}	Let $(\xi_t)_{t\ge 0}$ be a continuous real valued function with $\xi_0=0$. The solution $(f_t(z))_{t\ge 0,z\in \H}$ to the family of equations
	\begin{equation} \label{eq:reverseflow}
	\frac{\partial (f_t(z)+\xi_t)}{\partial t} = \frac{-2}{f_t(z)}, \;\;\; f_0(z)= z \;;\; z\in \H
	\end{equation}
is called the reverse Loewner flow driven by $(\xi_t)_{t\ge 0}$. In contrast to the forward case, $f_t(z)$ is defined for all $t\ge 0$ and $z\in \H$. This means that $f_t$ defines a conformal isomorphism from $\H$ to some domain $H_t$ for all $t$ (and one can check that $f_t(z)\sim z$ as $z\to \infty$ for each $t$). $(\H\setminus H_t)_{t\ge 0}$ is called the reverse Loewner evolution driven by $(\xi_t)_{t\ge 0}$.
\end{definition}

We will now discuss the (deterministic) relation between forward and reverse Loewner evolutions. For this, it is helpful to consider the centred forward Loewner maps $\tilde g_t := g_t - \xi_t$ and associated with a given driving function $(\xi_t)_{t\ge 0}$.

\begin{lemma}[Forward/Reverse flow]\label{L:fr}
	Suppose that $(\tilde{g}_t)_{t\ge 0}$ is the centred forward Loewner flow with driving function $(\xi_t)_{t\ge 0}$. Fix $T>0$ and write $\hat{\xi}_t=\xi_{T-t}-\xi_T$ for $0\le t \le T$. Let $(\hat{f}_t)_{0\le t\le T}$ be the centred reverse Loewner flow with driving function $(\hat{\xi}_t)_{0\le t\le T}$. Then
	\[ \hat{f}_t(z):=\tilde{g}_{T-t}\circ \tilde{g}_T^{-1}(z) \; ; \; t\in [0,T]\, , \, z\in \H.\]
 In particular, $\tilde{f}_T\equiv \tilde{g}_T^{-1}$.
\end{lemma}

\begin{proof}
	Since $\tilde{g}_t=g_t-\xi_t$ by definition, the forward Loewner equation \eqref{eq:forward_Loewner} and then the substitution $t\mapsto T-t$ yields 
	\[ \mathrm{d}(\tilde{g}_t(z))=\frac{2}{\tilde{g}_t(z)} \dd t +\dd \xi_t\; ; \; \mathrm{d}(\tilde{g}_{T-t}(z))=-\frac{2}{\tilde{g}_{T-t}(z)}-\dd \hat{\xi}_t\]
for every $z$. Replacing $z$ with $\tilde{g}_T^{-1}(z)$, we may deduce that $\hat{f}_t(z)$ satisfies the reverse Loewner equation \eqref{eq:reverseflow} with driving function $\hat{\xi}$.
\end{proof}

\paragraph{Reverse SLE.} Now we have defined reverse Loewner evolutions, reverse SLE$_\kappa$ is simply defined in the analogous way to forward SLE$_\kappa$.

\begin{definition}[Reverse SLE$_\kappa$]\label{D:reversesle}
	Reverse SLE$_\kappa$ for $\kappa>0$ is the centred reverse Loewner evolution driven by a Brownian motion with diffusivity $\kappa$. That is, with driving function $(\xi_t)_{t\ge 0}=(\sqrt{\kappa}B_t)_{t\ge 0}$ where $B$ is a standard Brownian motion.
\end{definition}
\ind{SLE! Reverse}

\begin{figure}
	\begin{center} \includegraphics[width=\textwidth]{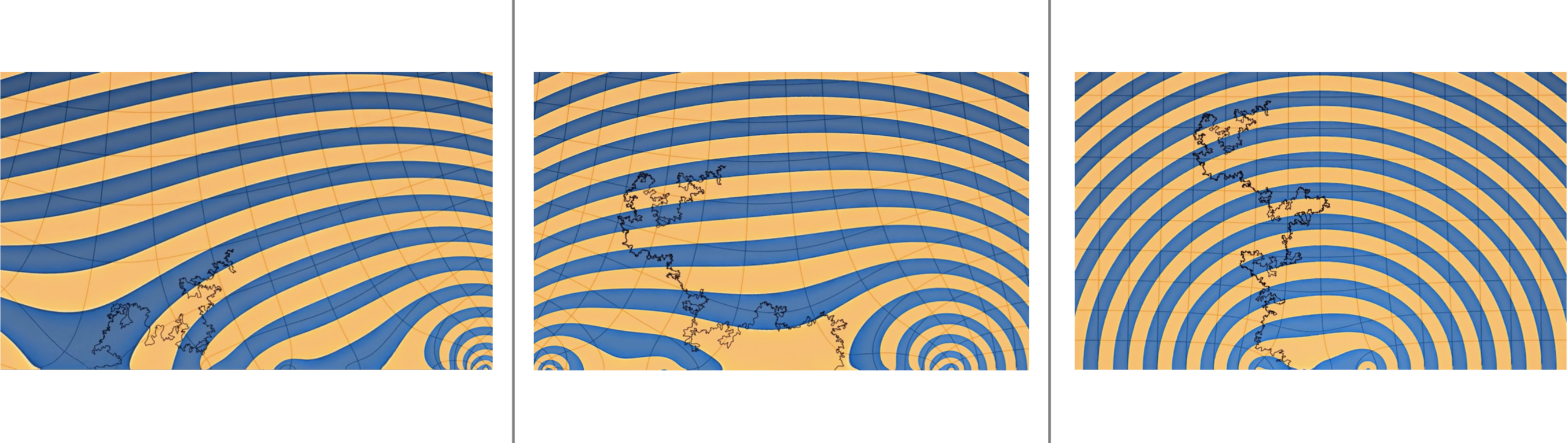} \end{center}
	\caption{A reverse SLE$_4$ at three increasing times. The background shows the deformation of the upper half plane under the reverse Loewner flow. Simulation by Henry Jackson.} \label{F:RSLE4}
\end{figure}

\begin{definition}[Reverse SLE$_\kappa(\underline{\rho})$ \cite{zipper}]\label{D:reverseslekp}
	 Suppose that $v^1,..., v^m\in \bar{\H}$ and $\rho^1,..., \rho^m$ are real numbers. Reverse SLE$_\kappa(\rho^1,..., \rho^m)$
	with force points at $v^1, ..., v^m$ is  the reverse (centred) Loewner evolution with driving function $(\xi_t)_{t\ge 0}$ satisfying:
	\begin{equation}\label{eqn:Revslekpdef}
	\xi_t  =  \sqrt{\kappa} B_t  -\sum_i \int_0^t \Re\left(\frac{\rho^i}{f_s(v^i) } \right)\, \dd s  \end{equation}
It is immediate that this has a unique solution in law, at least until the first time that $f_t(v^i)=0$ for some $i$. We will only consider the reverse SLE$_\kappa(\underline{\rho})$ up until this time.
\end{definition}

\begin{rmk}\label{rmk:bessel_rsle}
	In the case $m=1$ and $\rho^1=\rho$, a straightforward calculation shows that $f_t(v^1)$ is $\sqrt{\kappa}$ times a Bessel process of dimension
	\[ \delta = 1+ \frac{2(\rho-2)}{\kappa}.\]
	Note the difference with \cref{rmk:bessel_sle}. Roughly speaking, this is because the reverse SLE$_\kappa(\rho)$ generally pulls points towards the origin, while the forward version will be pushes them away (for intuition, consider the case $\rho=0$ and the way that the flow is defined).
\end{rmk}

The following properties of reverse SLE$_\kappa(\rho)$ will be needed for a technical discussion in \cref{S:zipper} of these notes. It says, roughly speaking, that specific SLE$_\kappa(\rho)$ curves are well behaved, in the sense that they do not create massive distortions, and that putting a force point very far away does not affect the law of the evolution at small times.
A reader simply wishing to learn about SLE would be safe to skip this.

\begin{lemma}\label{lem::reversesleprops}
	\begin{enumerate}
		\item[(1)] Let $(f_t)_{t\le \tau_1}$ be a reverse SLE$_\kappa(\kappa)$ flow, with a force point at $1\in \R$ and $\tau_1$ the first time that $f_t(1)=0$. Then as $R\to \infty$, the probability that $f_{\tau_1}(B(0,R))\supset B(0,1)$ tends to $1$.
		\item[(2)] Let $(\tilde{f}_t)_{t\le \tilde{\tau}_1}$ be a reverse SLE$_\kappa(\kappa,-\kappa)$ flow with force points at $(1,R)$ and $\tilde{\tau}_1$ the first time that $\tilde{f}_t(1)=0$. Then the total variation distance between $(f_t)_{t\le \tau_1}$ and $(\tilde{f}_t)_{t\le \tilde{\tau}_1}$ tends to $0$ as $R\to \infty$.
		\item[(3)] Let $(\tilde{f}_t)_{t}$ be a reverse SLE$_\kappa(\kappa,-\kappa)$ flow with force points at $(z,10)$, and $z\in [1,2]$. For $a\in (0,1]$ let $\tilde{\tau}_a$ be the first time that $\tilde{f}_t(z)=a$.
		Then with probability one,  $f_{\tilde{\tau}_a}(\{w\in \bar{\H}: |w-10|=1\})\subset \bar{\H}\setminus B(0, 1))$. %Then, uniformly in $z\in [1,2]$ $\mathbb{P}(f_{\tilde{\tau}_a}(\{w\in \bar{\H}: |w-10|=1\})\subset \bar{\H}\setminus B(0, \sqrt{a}))\to 1$ as $a\to 0$. \ellen{*probably want to change because proof implies something much stronger!*}
	\end{enumerate}
\end{lemma}

\begin{proof}
\begin{enumerate}
	\item[(1)] Note that $f_t(1)$ is $\sqrt{\kappa}$ times a Bessel process of dimension $3-(4/\kappa)<2$ started from $1$, and so the time $\tau_1$ 
	is almost surely finite.  Moreover, the driving function is continuous up to and including time $\tau_1$, because the integral $\int_0^{\tau_1} f_t(1)^{-1} \, dt$ converges almost surely. This implies that $f_{\tau_1}(z)\to \infty$ as $z\to \infty$ (see \cref{def:RLE}), which is the same thing as (1). 
	\item[(2)]For this we compute the Radon--Nikodym derivative between $(f_t)_{t\le \tau_1}$ and $(\tilde{f}_t)_{t\le \tilde{\tau}_1}$ using Girsanov's theorem. Let us write   $\mathbb{P}$ for the law of $(f_t)_{t\le \tau_1}$, and $\xi_t$ for the associated driving function, so for $z\in \H$ we have $\dd f_t(z)=-(2/f_t(z))\, \dd t $ $-\dd \xi_t$ and 
	\[ \dd \xi_t = \sqrt{\kappa}\dd B_t - \frac{\kappa}{f_t(1)}\, \dd t, \] where $B_t$ is a standard Brownian motion under $\mathbb{P}$.

We set 
	\[
	M_t:=R(R-1)^{-\kappa/2} (f_t(R)-f_t(1))^{\kappa/2}f_t(R)^{-1} f_t'(R)^{-1-\kappa/2}.
	\] 
	and first prove that $M_{\cdot \wedge \tau_1}$ is a local martingale under $\mathbb{P}$, (which is positive with unit mean by definition), with respect to the filtration $(\mathcal{F}_t)_{t\ge 0}$ generated by $(\xi_t)_{t\ge 0}$.
	
	This follows from It\^{o}'s formula. Indeed, for $A,B,C\in \R$, we have that for $t\le \tau_1$:
	\begin{align*}
		\mathrm{d}(f_t(R))^{A} & =Af_t(R)^{A-1}\dd f_t(R)+\frac{A(A-1)}{2} \mathrm{d}[f_t(R)]\\ & = \left(\frac{A(A-1)\kappa}{2}-2A\right) f_t(R)^{A-2}\dd t+\frac{\kappa A f_t(R)^{A}}{f_t(1)f_t(R)}\dd t -\sqrt{\kappa}A f_t(R)^{A-1} \dd B_t; \\
	\mathrm{d}(f_t'(R))^B & =Bf_t'(R)^{B-1}\dd f_t'(R) \\
	& = \frac{2Bf_t'(R)^B}{f_t(R)^2} \dd t; \quad \text{and} \\
\mathrm{d} (f_t(R)-f_t(1))^C & =C(f_t(R)-f_t(1))^{C-1}\left( -\frac{2}{f_t(R)}+\frac{2}{f_t(1)} \right) \dd t \\ & = \frac{2C(f_t(R)-f_t(1))^C}{f_t(1)f_t(R)} \dd t.
\end{align*}
	Thus if $M_t:=(f_t(R)^Af_t'(R)^B(f_t(R)-f_t(1))^C$, for $t\le \tau_1$ we have:
	\begin{align*}
		\mathrm{d}M_t & = (f_t(R)-f_t(1))^C \mathrm{d}(f_t(R)^Af_t'(R)^B)+f_t(R)^Af_t'(R)^B(f_t(R)-f_t(1))^C  \frac{2C}{f_t(1)f_t(R)} \dd t \\
		& = f_t(R)^A (f_t(R)-f_t(1))^C \mathrm{d}(f_t'(R)^B) + f_t'(R)^B (f_t(R)-f_t(1))^C \mathrm{d}(f_t(R)^A) \\
		& \quad \quad + f_t(R)^A f_t'(R)^B(f_t(R)-f_t(1))^C \frac{2C}{f_t(1)f_t(R)} \dd t \\
		& = f_t(R)^{A-2}f_t'(R)^B (f_t(R)-f_t(1))^C \left( -2A+2B+\frac{A(A-1)\kappa}{2}\right) \dd t \\ 
		& \quad \quad +  \frac{f_t(R)^Af_t'(R)^B(f_t(R)-f_t(1))^C}{f_t(1)f_t(R)}(\kappa A + 2C) \dd t - 
		\frac{\sqrt{\kappa} A M_t}{f_t(R)}\dd B_t.
	\end{align*}
Hence $M_{\cdot \wedge \tau_1}$ is a local martingale whenever $\kappa A + 2C=0$ and $-2A+2B+A(A-1)\kappa/2=0$ which is indeed the case when $A=-1$, $B=-1-\kappa/2$ and $C=\kappa/2$.

In this case we have 
	\[ 
	M_{t\wedge \tau_1}=1+ \sqrt{\kappa} \int_0^{t\wedge \tau_1} \frac{M_s}{f_s(R)} \dd B_s = \exp\left( Z_{t\wedge \tau_1}-\frac{1}{2}[Z]_{t\wedge \tau_1}\right)
	\]
	with \[ Z_{t\wedge \tau_1}= \sqrt{\kappa} \int_0^{t\wedge \tau_1} f_s(R)^{-1} \dd B_s,\]
where because
$
\frac{\dd}{\dd t} (f_t(R)-f_t(1))= \frac{2}{f_t(1)}-\frac2{f_t(R)}>0 $, we have
$f_t(R)>(R-1)$ for $t<\tau_1$.  This means that the quadratic variation is deterministically bounded above on any compact interval of time, so we can apply Novikov's condition to see that $M_{t\wedge \tau_1}$ is in fact a true martingale.

Furthermore, Girsanov's theorem tells us that if we define 
\[ 
\frac{\dd \tilde{\mathbb{P}}}{\dd \mathbb{P}}\Big|_{\mathcal{F}_{\tau_1}}:= M_{\tau_1},
\] 
then the process $\tilde{B}_t=B_t-\int_0^t (\sqrt{\kappa}/f_s(R)) \dd s$ will be a Brownian motion under $\tilde{\mathbb{P}}$.  Rewriting the expression for $\dd\xi_t$ in terms of $\tilde{B}$ we get 
\[\dd \xi_t = \sqrt{\kappa} \dd \tilde{B}_t -\frac{\kappa}{f_t(1)} 
	\dd t + \frac{\kappa}{f_t(R)} \dd t\]
	for $t\le \tau_1$, where $\tilde{B}$ is a standard Brownian motion under $\tilde{\mathbb{P}}$. Hence, under $\tilde{\mathbb{P}}$, the law of $\xi_{\cdot \wedge \tau_1}$ is that of the driving function of SLE$_{\kappa}(\kappa,-\kappa)$ with force points at $(1,R)$. Equivalently, the Radon--Nikodym derivative between the laws of $(\tilde{f}_t)_{t\le \tilde{\tau}_1}$ and $(f_t)_{t\le \tau_1}$ (that is, reverse SLE$_\kappa(\kappa,-\kappa)$ and SLE$_\kappa(\kappa)$) up to time $\tau_1$ is equal to $M_{\tau_1}$. 
	
To prove that the total variation distance between the laws of $(\tilde{f}_t)_{t\le \tilde{\tau}_1}$ and $(f_t)_{t\le \tau_1}$ converges to $0$ as $R\to \infty$, we first claim it suffices to prove that 
\begin{equation}\label{E:TVL2} \mathbb{E}(|M_{T\wedge \tau_1}-1|^2)=\mathbb{E}([M]_{T\wedge \tau_1})\to 0
	\end{equation}
as $R\to \infty$ for any fixed $T$. Indeed, suppose \eqref{E:TVL2} holds and let $\eps>0$ be given. We first choose $T$ large enough that $\mathbb{P}(\tau_1>T)<\varepsilon$. For this $T$, \eqref{E:TVL2} implies that the total variation distance between $(\tilde{f}_t)_{t\le \tilde{\tau}_1\wedge T}$ and $(f_t)_{t\le \tau_1\wedge T}$ converges to $0$ as $R\to \infty$. Therefore, for all $R$ large enough, we can couple the two processes so that they agree up to time $T\wedge \tau_1$ with probability greater than $1-\varepsilon$. But this implies that for all $R$ large enough, they can be coupled so that they agree up to time $\tau_1$ with probability $1-3\eps$. Since $\eps>0$ was arbitrary, this proves the desired convergence in total variation distance. 

To see \eqref{E:TVL2}, recall that $f_t(R)>R-1$ for $t\le \tau_1$. In addition, since 
$ \frac{\dd}{\dd t}f_t'(R)=\frac{2f_t'(R)}{f_t(R)^2}>0$ we have
$f_t'(R)>1$ for all $t\ge 0$.  Thus 
\[
\frac{M_t}{f_t(R)} \le R(R-1)^{-\kappa/2} f_t(R)^{\kappa/2-2} < \frac{R}{(R-1)^2}<\frac2R
\]
 for all $t\le \tau_1$ and hence, since $[M]_{T\wedge \tau_1}=\kappa \int_0^{T\wedge \tau_1} (M_t/f_t(R))^2 \dd t $,
		\[   [M ]_{T\wedge \tau_1} \le \frac{4 \kappa (T\wedge \tau_1) }{R^2} \]
		for any $T\ge 0$. Taking expectations and using the deterministic upper bound $4\kappa T/R^2$, which converges to $0$ as $R\to \infty$, proves \eqref{E:TVL2}.
	
	\item[(3)] For this, we claim that for $w\in \H$ with $\Re(w)>2$ and $z\in [1,2]$, the process $\Re(\tilde{f}_t(w))-\tilde{f}_t(z)$ is increasing for $t\le \tilde{\tau}_a$ (which clearly implies the result). To see the claim, observe that by definition of the reverse flow \[ \frac{\partial (\Re(\tilde{f}_t(w))-\tilde{f}_t(z))}{\partial t} = \frac{2}{\tilde{f}_t(z)} - \Re(\frac{2}{\tilde{f}_t(w)})= \frac{2}{\tilde{f}_t(z)} - \frac{2\Re(\tilde{f}_t(w))}{|\tilde{f}_t(w)|^2},\]  which is positive as long as $\Re(\tilde{f}_t(w))>\tilde{f}_t(z)>0$. Since this is true at time $0$ for $w$ with $\Re(w)>2$, it is therefore positive for all $t\le \tilde{\tau}_a$, and the process $\Re(\tilde{f}_t(w))-\tilde{f}_t(z)$ is increasing for this range of $t$. 
\end{enumerate}
\end{proof}

\subsection{Symmetries in law for forward/reverse \texorpdfstring{$\mathrm{SLE}_\kappa$ and $\mathrm{SLE}_\kappa(\rho)$}{TEXT}}
Now, because Brownian motion has time reversal symmetry, the relationship \cref{L:fr} between forward and reverse Loewner evolutions has particularly nice consequences for SLE.

More specifically, if $T>0$ is fixed and $(\xi_t)_{0\le t\le T}$ is $\sqrt{\kappa}$ times a Brownian motion, then $(\hat \xi_t)_{0\le t\le T}=(\xi_T-\xi_{T-t})_{0\le t\le T}$ also has the law of $\sqrt{\kappa}$ times a Brownian motion. Consequently:
	
	\begin{lemma}\label{L:forwardreverseSLE}
		For any fixed $T>0$ the curve generated by a reverse SLE$_\kappa$ run up to time $T$ and the curve generated by a forward SLE$_\kappa$ run up to time $T$ are equal in law.
	\end{lemma}
Mind that the \emph{processes} of the previous lemma, defined for all times $t\in [0,T]$, are \emph{not} the same in law. Indeed, we have seen that forward and reverse Loewner evolutions generate hulls via a completely different dynamic. Nonetheless, it is a very useful property that at any fixed time, the laws of the generated hulls are equal.\\

\begin{figure}
	\begin{center}
		\includegraphics[scale=.7]{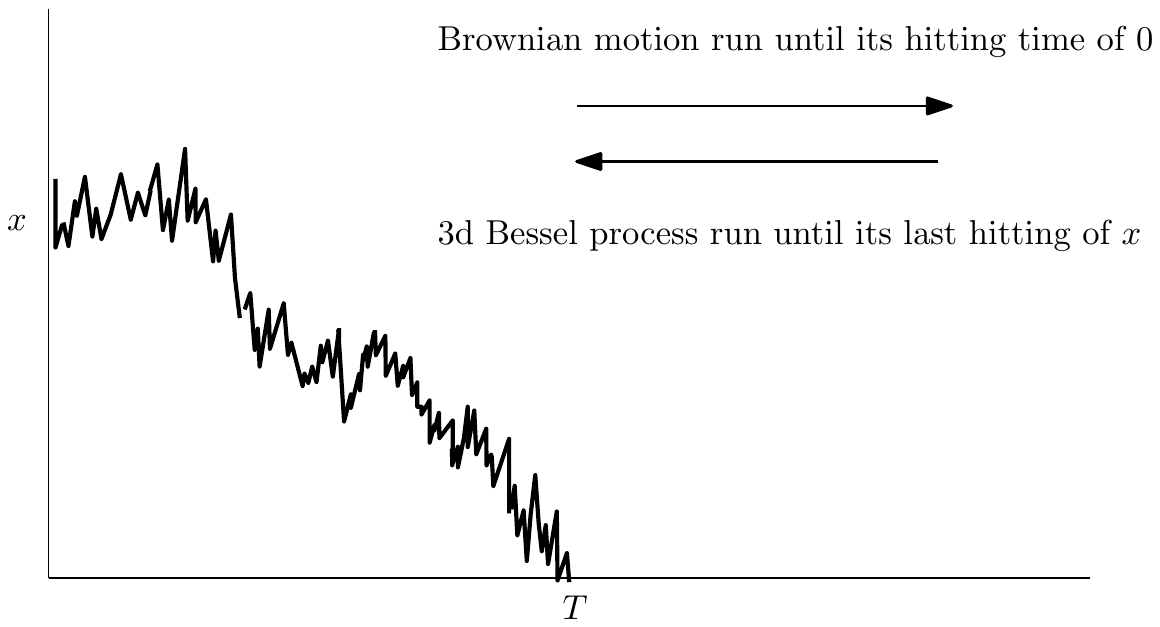}
		\caption{Illustration of Williams' path decomposition theorem. The classical result says that if $X$ is a Brownian motion started from $x>0$ and $T$ is its hitting time of zero, then its time-reversal $\hat X= (X_{T-t})_{0 \le t \le T}$ is distributed as three dimensional Bessel process, run until its last visit $\Lambda$ to $x$.}\label{F:williams}
	\end{center}
\end{figure}

There are similar consequences for SLE$_\kappa(\underline{\rho})$ processes, but the reversibility properties of solutions to \eqref{eqn:slekpdef} are somewhat more complicated. We will explain now what happens in the simplest case of one marked point. Due to remarks \cref{rmk:bessel_sle,rmk:bessel_rsle}, this requires understanding how Bessel processes behave under time reversal. 
%\ind{Bessel processes}

\begin{rmk}[Bessel process properties] Recall that the dimension $\delta$ of a Bessel process determines how often it returns to $0$: if $\delta\ge 2$, then the Bessel process will almost surely be strictly positive for all positive times; while if $\delta<2$ then from any starting point it will return to $0$ in finite time almost surely.
\end{rmk}

	\ind{Williams' path decomposition theorem} The following is an extension of a classical result about Brownian motion, due to Williams (see for example Corollary (4.6) in Chapter VII of \cite{RevuzYor} and  %which says that if $X$ is a Brownian motion started from $x>0$ and $T$ is its hitting time of zero, then its time-reversal $\hat X= (X_{T-t})_{0 \le t \le T}$ is distributed as three-dimensional Bessel process, run until its last visit $\Lambda$ to $x$
Figure \ref{F:williams}).

\begin{lemma}[Time reversal of Bessel processes]\label{lem:TR_Bessel}
	Suppose that $X$ is a Bessel process of dimension $\delta \in (0,2)$ started from $x>0$, run until its first hitting time $T$ of zero. Then its time reversal $\hat X = (X_{T-t}, 0 \le t \le T)$ is a Bessel process of dimension $\hat \delta =4 - \delta \in (2, 4)$, run until its last visit $\Lambda$ to $x$.	
	\ind{Bessel!Process}
	\ind{Bessel!Time reversal}
\end{lemma}

The proof of this will boil down to an analogous result for Brownian motion with drift, that we state and prove first.

\begin{lemma}\label{L:BM_drift_rev}
		Let $\mu> 0$. Then the time reversal of a Brownian motion with drift $\mu$, started from $0$ and stopped at its last hitting time of $y>0$, has the law of a Brownian motion with drift $-\mu$, started from $y$ and run up to its last hitting time of $0$.
\end{lemma}

\begin{proof}
	Let $(X_t)_{t\in \R}=(B_t+\mu t)_{t\in \R}$, where $B_t$ is a standard two-sided Brownian motion with $B_0=0$. Then $(\hat{X}_t)_{t\in \R}:=(X_{-t})_{t\in \R}$ is equal in law to $(B_t-\mu t)_{t\in \R}$.
	Define $\tau_0:=\{\inf: s\le 0: X_s = 0\}$ and $\tau_y = \sup\{s\ge 0: X_t=y\}$. Then by the strong Markov property at time $\tau_0$, $(X_{\tau_0+s})_{0\le s\le \tau_y-\tau_0 }$ has the law of a Brownian motion with drift $\mu$, started from $0$ and stopped at its last hitting time of $y$. So, we need to show that the time reversal $(X_{\tau_y-s})_{0\le s \le \tau_y-\tau_0}$ has the law of a Brownian motion with drift $-\mu$, started from $y$ and run up to its last hitting time of $0$.
	
	For this, we use the fact that, by definition of $\hat{X}$, $$(X_{\tau_y-s})_{0\le s \le \tau_y-\tau_0}=(\hat{X}_{s-\tau_y})_{0\le s \le \tau_y-\tau_0}=(\hat{X}_{s+\hat{\tau}_y})_{0\le s \le \hat{\tau}_0},$$ where $\hat{\tau_y}$ is the first time before $0$ that $\hat{X}$ hits $y$, and $\hat{\tau}_0$ is the last time that $(\hat{X}_{t+\hat{\tau}_y})_{t\ge 0}$ hits $0$. Since $\hat{X}$ is equal in law to a two-sided Brownian motion with drift $-\mu$, the law of the process on the right hand side above is (by the strong Markov property again, but this time for $\hat X$) indeed that of a Brownian motion with drift $-\mu$, started from $y$ and run up to its last hitting time of $0$. This concludes the proof.
\end{proof}

\begin{proof}[Proof of \cref{lem:TR_Bessel}]
	(\cite[Proposition 3.5]{DuplantierMillerSheffield})
We will make use of the following fact, which is just a rewriting of \cref{L:BessBM}/\cref{R:BessBM}:
\begin{itemize}
	\item Let $\tau(t)=\inf\{s>0: [\log(X)]_t>t\}$ and let $Z_t=\log(X_{\tau(t)})$ (recall that $[M]_t$ denotes the quadratic variation of the continuous semimartingale $M$). Note that because $\delta\in (0,2)$, $\tau(t)\uparrow T$ as $t \uparrow \infty$. Then \begin{equation}
	\label{eq:bessel_BM_corr}
(Z_t)_{t\ge 0} \overset{(\mathrm{law})}{=} (B_t + \frac{\delta-2}{2}t)_{t\ge 0} ,	\end{equation}
	where $B$ is a standard Brownian motion with $B_0=\log x$.
\end{itemize}
We now want to use this, along with the time reversal symmetry of Brownian motion,  to draw a conclusion similar to \eqref{eq:bessel_BM_corr} about the time reversal $\hat{X}$ of $X$, but with the opposite drift (corresponding to a dimension $\hat \delta = 4 - \delta$, as claimed). However, there is a slight technical complication that arises, since $\hat{X}_0=0$ and so $\log(\hat{X}_0)=-\infty$.

To get around this, we also define for any $\eps<x$, $T_\eps$ to be the last time before $T$ that $(X_t)_{t\ge 0}$ hits $\eps$. Then $(Z_t)_{t\in [0,[\log X]_{T_\eps}]}$ is a Brownian motion with drift as in \eqref{eq:bessel_BM_corr}, started from $\log x$ and stopped at its last hitting time of $\log(\eps)$. This implies (by \cref{L:BM_drift_rev}) that the time reversal of $Z$ with respect to this time interval is a Brownian motion with drift $-(\delta-2)/2=(\hat \delta -2)/2$, started from $\log \eps$ and run up to its last hitting time of $\log x$.

Reversing the argument for \eqref{eq:bessel_BM_corr} (that is, taking the exponential and reparametrising by quadratic variation), this implies that the time reversal of $(X_t)_{t\in [0,T_\eps]}$ is a Bessel process of dimension $(4-\delta)$, started from $\eps$ and run up to its last hitting time of $x$. Taking a limit as $\eps\to 0$ provides the result.
 \end{proof}
As a consequence of this and \cref{rmk:bessel_sle,rmk:bessel_rsle}, we obtain the following:

\begin{cor}[Symmetries for forward and reverse SLE$_\kappa(\rho)$]\label{C:slerevsym}
	Suppose that $(f_t)_{t\ge 0}$ is the reverse flow for a centred, reverse SLE$_\kappa(\rho)$ process with a single force point at $x>0$ of weight $\rho<\kappa/2 +2$. Consider the first time $\tau$ that $f_t(x)=0$. Then $H_\tau=f_\tau(\H)$ has the same law as $\H\setminus \eta([0,\sigma])$, where $\eta$ is a forward SLE$_\kappa(\kappa-\rho)$ curve with a force point at $0^+$, run until the \emph{last time} $\Lambda$ that the centered forward Loewner flow for $\eta$ sends $0^+$ to $x$. \end{cor}

\section{Radial Loewner chains and radial SLE}
\label{S:App_radial}
\subsection{Radial Loewner chains} While chordal Loewner chains describe ``locally growing'' sets started at one point on the boundary of a domain and targeted at another, \emph{radial} Loewner chains describe growing sets started on the boundary but targeted at a point in the \emph{interior} of the domain. The canonical configuration for chordal Loewner chains is the upper half plane $\H$, with starting point $0\in \partial \H$ and target point $\infty$. For radial Loewner chains, things turn out to be nicest if one works in the unit disc $\D\subset \C$ with starting point $1\in \partial \D$ and target point $0\in \D$.

\begin{definition}[Radial Loewner chain]
	Let $(U_t)_{t\ge 0}$ be a continuous process taking values in the unit circle $\partial \D$, with $U_0=1$. The \emph{radial Loewner chain driven by $U$} is the collection of maps $(g_t)_{t\ge 0}$ that solve the \emph{radial Loewner equation}:
	\begin{equation}\label{E:RLE}
		\frac{\partial g_t}{\partial t}(z) = g_t(z) \frac{U_t+g_t(z)}{U_t-g_t(z)}; \quad g_0(z)=z,
	\end{equation}
for each $z\in \D$ until time $\zeta(z):=\inf_{t>0} g_t(z)=U_t$. If one defines $D_t:=\{z\in \D: \zeta(z)>t\}$ for each $t\ge 0$, then $g_t$ is the unique conformal isomorphism  $$g_t:D_t\to \D \text{ with } g_t'(0)=e^t \text{ and } g_t(0)=0,$$
\cite{SW05}. The hulls generated by $U$, $K_t:=\D\setminus D_t$ for $t\ge 0$, are an increasing family of compact sets in $\D$. With a slight abuse of notation we will sometimes also refer to $(D_t)_{t\ge 0}$ or $(K_t)_{t\ge 0}$ as the Loewner chain driven by $U$.
\end{definition}

As in the chordal case, continuous non-crossing curves $(\gamma(t))_{t\ge 0}$ in $\bar{\D}$, with $\gamma(0)=1$, and parametrised so that $-\log \mathrm{CR}(\D\setminus \gamma([0,t]);0)=t$ for $t\ge 0$, provide examples of radial Loewner chains. That is, when $g_t$ is defined for each $t$ to be the unique conformal isomorphism from $\D\setminus \gamma([0,t])$ fixing $0$ and with positive real derivative at $0$.

\subsection{Radial \texorpdfstring{$\mathrm{SLE}_\kappa$ and $\mathrm{SLE}_\kappa(\rho)$}{TEXT}}\label{A:radialSLE}
\ind{SLE! Radial}
\begin{definition}[Radial SLE$_\kappa$]
	For $\kappa\ge 0$, radial SLE$_\kappa$ in $\D$ from $1$ to $0$ is defined to be the radial Loewner chain driven by $$(e^{i\sqrt{\kappa}B_t})_{t\ge 0}$$
	where $B$ is a standard one dimensional Brownian motion.
\end{definition}

For a general simply connected domain $D$ with marked boundary point $a\in \partial D$ and interior point $b\in D$, we define the radial SLE$_\kappa$ in $D$ from $a$ to $b$ to be the random process obtained by taking the image of a radial SLE$_\kappa$ in $\D$ from $1$ to $0$ under the unique conformal isomorphism from $\D$ to $D$ sending $1\mapsto a$ and $0\mapsto b$.

\paragraph{Radial SLE$_\kappa(\underline{\rho})$.} \ind{SLE! Radial with force points}

As with chordal SLE, we can generalise the definition of radial SLE$_\kappa$ by placing force points on the boundary or in the interior of the domain and keeping track of their evolution under the radial Loewner flow. The definition in the unit disc from $1$ to $0$ is as follows.

\begin{definition}[SLE$_\kappa(\rho^1,..., \rho^m)$ in $\D$ from $1$ to $0$]\label{D:radialSLEkp1}
	Suppose that $v^1,..., v^m\in \overline{\D}\setminus\{1\}$ are distinct and $\rho^1,..., \rho^m$ satisfy the condition \eqref{eqn:rho_nice}. Radial SLE$_\kappa(\rho^1,..., \rho^m)$ from $1$ to $0$ in $\D$ 
	with force points at $v^1, ..., v^m$ is  the radial Loewner chain, whose driving function $(U_t)_{t\ge 0}$ satisfies:
	\begin{eqnarray}\label{eqn:Rslekpdef}
		U_t & = & 1 + i\sqrt{\kappa} \int_0^t U_s \dd B_s -\int_0^t \frac{\kappa}{2} U_s \dd s + \sum_i \int_0^t \frac{\rho^i}{2} \hat{\Phi}(V_s^{i},U_s) \dd s \nonumber \\
		V_t^{i} & = & v^{i}+\int_0^t \Phi(U_s,V_s) \dd t \text{ for } 1\le i \le M.
	\end{eqnarray}
	Above we denote $\Phi(u,z)=z\tfrac{u+z}{u-z}$ and $\hat{\Phi}(u,z)=\tfrac{\Phi(u,z)+\Phi(1/\bar{u},z)}{2}$ for $z\in \D$ and $u\in \partial \D$.
\end{definition}

When $\underline{\rho}$ satisfies \eqref{eqn:rho_nice}, the existence and uniqueness of a continuous $(U,V^1,\dots, V^m)$ satisfying \eqref{eqn:slekpdef} up to the first time that $V^i_t=U_t$ for some $i$ with $v^i\in \D$ is proven in \cite{MSIG1}. In particular, there is a unique solution for all time when all of the force points $v^i$ are on the boundary $\partial \D$. %It is also proven there that such a radial SLE$_\kappa(\underline{\rho})$ is almost surely generated by a continuous curve.

As with ordinary radial SLE$_\kappa$, we define radial SLE$_\kappa(\underline{\rho})$ in a domain $D$ from $a\in \partial D$ to $b\in D$, with force points $v^1,\dots, v^m\in \overline{\D}$, to be the image of SLE$_\kappa(\underline{\rho})$ in $\D$ from $1$ to $0$ with force points at $\varphi(v^1), \dots, \varphi(v^m)$, where $\varphi$ is the unique conformal isomorphism from $D$ to $\D$ sending $a$ to $1$ and $b$ to $0$.

\begin{rmk}
Again we can extend the definition of SLE$_\kappa(\underline{\rho})$ to include force points located infinitesimally clockwise (respectively anticlockwise) from $1$ on $\partial \D$, by taking a limit (in the same way as for chordal SLE$_\kappa(\underline{\rho})$). 
\end{rmk}

Radial SLE$_\kappa(\underline{\rho})$ satisfies a very similar collection of properties to chordal SLE$_\kappa(\underline{\rho})$. Indeed, there is a simple connection between the radial and chordal variants, that can be verified using a careful stochastic calculus argument (omitted here).

\begin{lemma}\cite{SW05}
	Let $D$ be a simply connected domain, $a,b\in \partial D$ be boundary points, and $c\in D$ be an interior point. Let $\rho^1,\dots, \rho^m\in \R$ with $\sum_i \rho_i = \kappa-6$ satisfy \eqref{eqn:rho_nice}, and let $v^1,\dots, v^m\in \overline{D}$.
	
	 Suppose that $\eta$ is a radial $\mathrm{SLE}_\kappa(\underline{\rho})$ from $a\in \partial D$ to $b\in D$ (with force points at $v^1,\dots, v^m$) stopped at the infimum over $t$ for which $c$ and $b$ are in different connected components of $D\setminus \eta([0,t])$, or an interior force point is swallowed. Let  $\tilde{\eta}$ be a chordal $\mathrm{SLE}_\kappa(\underline{\rho})$ from $a\in \partial D$ to $c\in \partial D$ (with force points at $v^1,\dots, v^m$), stopped at the corresponding time. Then, as curves modulo reparametrisation of time, $\eta$ and $\tilde{\eta}$ agree in law. \label{L:coordinate_change}
\end{lemma}

\begin{rmk}
As already observed, adding a force point of any weight to the target point of a chordal or radial SLE$_\kappa(\underline{\rho})$ does not effect the law of the curve. So if we start with a given chordal or radial SLE$_\kappa(\underline{\rho})$, we can add such a force point so that the new weights add up to $\kappa-6$.
 \end{rmk}

For example, if we want to sample a radial SLE$_\kappa$ from $a\in \partial D$ to $b\in D$, then we can first run a chordal SLE$_\kappa(\kappa-6)$, $\eta_1$, in $D=:D_1$ from $a$ to some arbitrary $c_1\in \partial D$, with the force point at $b$, up until the first time $\tau_1$ that $c_1$ and $b$ are separated by $\eta_1$. Then, we can run an SLE$_\kappa(\kappa-6)$, $\eta_2$, in the connected component $D_2$ of $D\setminus \eta_1([0,\tau_1])$ containing $b$, from $\eta(\tau_1)$ to some other $c_2\in \partial D_2$ and with force point at $b$, and again stop it when $\eta_2$ first separates $c_2$ and $b$. Iterating this procedure, and reparametrisating the concatenated curve so that the conformal radius of $b$ in the to be explored domain is always $e^{-t}$, we obtain a curve with the law of radial SLE$_\kappa$ from $a$ to $b$. 

Similar procedures will work to generate radial SLE$_\kappa(\underline{\rho})$ with non-trivial $\underline{\rho}$.
In particular,  the following properties hold.

\paragraph{Radial SLE$_\kappa(\underline{\rho})$: properties.}
\begin{itemize}
	\item Suppose that $D$ is a simply connected domain, $\underline{\rho}$ satisfies \eqref{eqn:rho_nice}, and $v^1,\dots, v^m \in \partial D$. Then radial SLE$_\kappa(\underline{\rho})$ in $D$ from $a$ to $b$ is almost surely generated by a curve $\gamma$ (that is, there exists a curve $\gamma(t)$ defined for all time such that the connected component of $\D\setminus \gamma([0,t])$ containing $0$ is equal to $D_t=\{z\in \D: \tau_z>t\}$ for all $t$. We will also sometimes refer to the curve $\gamma$ as ``the radial SLE$_\kappa$''. When there are interior force points, the radial Loewner chain is generated by a continuous curve, until the first time that one of the interior force points is swallowed. Lawler proved in \cite{LawlerWP} that if $\gamma$ is a radial SLE$_\kappa$ with target point $b$, then $\lim_{t\to \infty} \gamma(t)=b$ almost surely. This was extended to the case of SLE$_\kappa(\underline{\rho})$ with boundary force points and $\underline{\rho}$ satisfying \eqref{eqn:rho_nice} in \cite{IG4}.
	\item \emph{Conformal invariance:} follows from the definition of radial SLE$_\kappa(\underline{\rho})$ in $D\ne \D$ from $a\in \partial D$ to $b\in D$ (see above).
	\item \emph{Domain Markov property:} suppose that $D$ is a simply connected domain, and $\gamma$ is a radial SLE$_\kappa(\underline{\rho})$ from $a\in \partial D$ to $b\in D$, with weights $\underline{\rho}=(\rho^1,\dots, \rho^m)$ satisfying \eqref{eqn:rho_nice} and force points $v
^1,\cdots, v^m\in \partial D$. Suppose that $T$ is a bounded stopping time that is measurable with respect to $\gamma$. Then, conditionally on $\gamma([0,T])$ and writing $D_T$ for the connected component containing $b$ of $D\setminus \gamma([0,T])$, $\gamma([T,\infty))$ has the law of an SLE$_\kappa(\underline{\rho})$ from $\gamma(T)$ to $b$ in $D_T$, with force points at $(V_T^1,\dots, V_T^m)$.
\end{itemize}

\paragraph{Target invariance.}\ind{SLE! Target invariance}
Finally, we consider the special case when $\eta$ is a radial SLE$_\kappa(\kappa-6)$ from $a\in \partial D$ to $b\in D$ for some domain $D$, and with force point $c\in \partial D$. Suppose that $\kappa\ge 4$ so that \eqref{eqn:rho_nice} holds with $m=1$ and $\rho^1=\kappa-6$. \cref{L:coordinate_change} then implies that $\eta$ (which is defined for all time) can be sampled as follows.
\begin{itemize}
	\item[-] Choose $x_1$ on $\partial D$ and run a chordal SLE$_\kappa(\kappa-6)$, (with force point at $c$) in $D$ from $a$ to $x_1$, stopped at the first time that $x_1$ and $b$ lie in separate connected components of complement of the curve. Reparametrise this curve so that the conformal radius of $b$ in its complement is equal to $e^{-t}$ for all $t$ up to the time that the curve is stopped. 
	\item[-] Repeat the first step with the new domain being the connected component of the complement of the first curve containing $b$, the new start point being the tip of the curve at the disconnection time, and the new force point being the image of $c$ under the radial Loewner flow generated by the first curve\footnote{This is well defined after conformally mapping $D$ to $\D$, $a$ to $1$ and $b$ to $0$.}.
	\item[-] Iterate the above procedure.
\end{itemize}

In particular, note that the only dependence on $b$ in the above is the choice of exploration domain at ``disconnection times''. This means that if $b'$ is another point in $D$, the above procedure (run until $b$ and $b'$ are first separated by the curve) also produces a sample of radial SLE$_\kappa(\kappa-6)$ (with the same force point) from $a\in \partial D$ to $b'$ (and stopped  when $b$ and $b'$ are first separated). More precisely, and also applying the Markov property of radial SLE after this separation time, we have the following.

\begin{lemma}[Target invariance of SLE$_\kappa(\kappa-6)$] Suppose that $\kappa>4$ and let $D$ be a simply connected domain with $a,c\in \partial D$. For $b_1,b_2\in D$, one can couple a radial SLE$_\kappa(\kappa-6)$ curve from $a$ to $b_1$ in $D$ (with force point at $c$), and from $a$ to $b_2$ in $D$ (with force point at $c$) so that they coincide until $b_1,b_2$ are contained in separate components of the complement of the curve, and afterwards evolve independently.\label{L:radialtarget}
\end{lemma}

In fact, the above lemma means that for given $D,a,c$, SLE$_\kappa(\kappa-6)$ can be simultaneously defined towards a countable dense set of target points in $D$, in such a way that the above description holds for any two given target points. The object created in this manner is referred to as an SLE$_\kappa(\kappa-6)$ branching tree, or sometimes just a branching SLE$_\kappa$.

\section{Convergence of random variables in the space of distributions}
\label{App:distributions}

In this appendix we prove Lemma \ref{L:measurability}, about the measurability of the convergence event for a sequence of random  variables in the space of distributions $\cD'_0(D)$, which we restate here for convenience:

\begin{lemma}
Let $D$ be a domain of $\R^d$. Let $\Conv$ denote the set of sequences in $\cD'_0(D)$ which are weak$-*$ convergent. Then $\Conv$ is a Borel set in $\cD'_0(D)^{\N}$ equipped with the product Borel $\sigma$-algebra. 
\end{lemma}

\begin{proof}
\ind{Schwartz space} \ind{Tempered distributions} The proof relies on some results in functional analysis, and in particular uses the Schwartz space $\cS^*$ of \textbf{tempered distributions}, whose definition is as follows. Let $\cS$ denote the space of rapidly decaying test functions
$$
\cS = \{ f \in C^\infty ( \R^d): \| f\|_j < \infty, \text{for all } j \ge 1\},
$$
where
$$
\|f\|_j: = \sup_{|\alpha | + |\beta| \le j} \| x^\alpha \partial^\beta f \|_\infty 
$$
and we use standard multi index notation for $\alpha = (\alpha_1, \ldots, \alpha_d)$ and $\beta = (\beta_1, \ldots, \beta_d)$. We equip $\cS$ with a topology defined by the requirement  $f_n \in \cS$ converges to $f\in \cS$ if and only if $\|f_n - f\|_j \to 0$ for every $j\ge 1$. 
Thus, the quantities $\|f\|_j$ define seminorms on $\cS$ (actually, norms) which together (by definition)  generate the topology on $\cS$. We define $\cS^*$ to be the space of continuous linear functionals on $\cS$, equipped once again with the weak-$^*$ topology of pointwise convergence: that is, a sequence $x_n^*\in \cS^*$ converges to $x^* \in\cS^*$ if and only if $x_n^* (x) \to x^*(x)$ for all $ x\in \cS$. 

The advantage of the space $\cS$ over $\cD_0(D)$ is that it is a (separable) \textbf{Fr\'echet space}, that is, a topological vector space which is locally convex, metrisable and complete. Equivalently, a Fr\'echet space is one for which there is a countable family of seminorms generating the topology (which is clearly the case for $\cS$), and which is complete and Hausdorff (which is also straightforward to check in the case of $\cS$). The separability of $\cS$ is also a standard fact.  We will first prove the lemma with $\cD'_0(D)$ replaced by $\cS^*$, and then explain how to go from $\cS^*$ to $\cD'_0(D)$. For $\cS^*$, the statement boils down to a general fact about Fr\'echet spaces, which we now introduce.

Fix $X$ a separable Fr\'echet space and fix a dense countable set $\cQ = \{x_i\}_{i\ge 1}$ in $X$.  Let $X^*$ be the set of continuous linear functionals on $X$. Let $\Conv(X^*)$ be the set of converging sequences in $X^*$: that is, 
$$
\Conv(X^*) = \{ (x_n^*)_{n\ge 1} \in (X^*)^{\N}: x_n^* \to x^* \text{ weak-$^*$ for some $x^* \in X^*$} \}.
$$
Let 
$$
V_j = \{ x \in X: \| x\|_k \le 1/j \text{ for all } 1\le k \le j\}
$$
where $(\|\cdot \|_j)_{j\ge 1}$ is a countable family of seminorms generating the topology of $X$. Then $V_j$ is what is called a countable basis of neighbourhoods of 0: that is, for any neighbourhood $V$ of 0 there exists $j\ge 1$ such that $V \supset V_j$.  

Let $K_j = {V}^\bullet_j$ denote the \textbf{polar set} of $V_j$, that is, 
$$
{V}^\bullet_j = \{ x^* \in X^*: | x^*(x) | \le 1 \text{ for all } x \in V_j\}. 
$$

\textbf{Claim.} 
\begin{equation}
\label{eq:Frechetconv}
\Conv(X^*) = \bigcup_{j\ge 1} \{ (x_n^*)_{n\ge 1} \in  K_j^{\N} \text{, such that } \ x_n^*(x) \text{ converges in $\R$ for all $x\in \cQ$}\}.
\end{equation}

\begin{proof}[Proof of \eqref{eq:Frechetconv}]
We have two inclusions to prove. We start by showing that the right hand side of \eqref{eq:Frechetconv} is contained in $\Conv(X^*)$. Let $(x_n^*)_{n\ge 1}$ denote a sequence in $X^*$ and suppose that there is some $j\ge 1$ such that $x_n^* \in K_j$ for all $n\ge 1$, and that $x_n^*(x)$ converges to a limit $\ell(x) \in \R$ for all $x \in \cQ$. By the Alaoglu theorem, see for example  \cite[Theorem 8.4.1]{TVS}, $K_j$ is compact and furthermore metrisable (as a bounded subset of the dual of a separable space), hence sequentially compact. Let $x^*$ be any weak-$^*$ subsequential limit. Then $x^* (x) = \ell(x)$ for all $x \in \cQ$. Since $\cQ$ is dense and $x^*$ is continuous, this identifies $x^*$ uniquely. Thus $x_n^*$ converges to $x^*$ in the weak-$^*$ sense. Thus $(x_n^*)_{n\ge 1} \in \Conv(X^*)$. 

Conversely, suppose $(x_n^*)_{n\ge1} \in \Conv(X^*)$. Clearly for $x \in \cQ$, $\langle x_n^*, x \rangle$ converges in $\R$ by definition of the weak-$^*$ topology. Therefore it suffices to show that there exists $j\ge 1$ such that $x_n^* \in K_j$ for all $n\ge 1$. We rely on the Banach--Steinhaus theorem in Fr\'echet spaces \cite[Theorem 11.9.1]{TVS}, which states that if $(x_n)^*_{n\ge 1}$ is pointwise bounded (that is, if $\sup_{n\ge 1} | x_n^*(x) | < \infty$ for any $x \in X$, which is the case since $x_n^*(x)$ converges in $\R$) then $x_n^*$ is ``bounded in the operator norm'' (more precisely, equicontinuous): that is, there is some neighbourhood $V$ of 0 such that $x_n^* \in {V}^\bullet$, the polar set of $V$. Since $(V_j)_{j\ge 1}$ is basis of neighbourhoods, we can find $j\ge 1$ such that $V_j \subset V$ and thus ${V}^\bullet \subset K_j$. This concludes the proof of \eqref{eq:Frechetconv}.
\end{proof}
An immediate consequence of \eqref{eq:Frechetconv} is that $\Conv(X^*)$ is a Borel set in $(X^*)^{\N}$ equipped with the power Borel $\sigma$-field. As already mentioned, this applies in particular to the case $X = \cS$, $X^* = \cS^*$. 

To conclude the proof of Lemma \ref{L:measurability}, it remains to reduce the convergence in the sense of distributions to convergence in the Schwartz space $\cS^*$ as follows. Let 
$$D_k = \{x \in D: \dist(x, \partial D) \ge 1/k \} \cap B(0, k).$$ 

Fix $(\varphi_k)_{k\ge 1}$ a sequence of test functions with compact support in $D$ such that:
\begin{itemize}
\item $\varphi_k \ge 0$, 

\item $\text{Supp} (\varphi_k) \subset D_{2k+1}$, 

\item $\varphi_k \equiv 1$ on $D_{2k}$. 
\end{itemize} 

For a distribution $T\in \cD'_0(D)$ and $k\ge 1$, let $T\varphi_k$ denote the distribution obtained by setting 
$$
(T\varphi_k, f) = ( T, f\varphi_k).
$$
Then note that, given a sequence $(T_n)_{n\ge 1} \in (\cD'_0(D))^{\N}$, we have 
\begin{equation}
\label{eq:conv1}
(T_n)_{n\ge 1} \in \Conv \iff (T_n \varphi_k)_{n\ge 1} \in \Conv, \text{ for all } k \ge 1.
\end{equation}
Indeed, given a test function $f$ with compact support in $D$, it is always possible to find a $k\ge 1$ such that $\text{Supp} (f) \subset  D_{2k}$. On the other hand, for a fixed $k\ge 1$,
\begin{equation}
\label{eq:conv2}
(T_n \varphi_k)_{n\ge 1} \in \Conv \iff (T_n \varphi_k)_{n\ge 1} \in \Conv( \cS^*).
\end{equation}
One implication in \eqref{eq:conv2} is trivial: if $(T_n \varphi_k)_{n\ge 1} \in \Conv( \cS^*)$ and $f$ is a test function with compact support in $D$, then clearly $f \in \cS$, so $(T_n \varphi_k , f)$ converges as desired. Conversely,  if $(T_n \varphi_k)_{n\ge 1} \in \Conv$ and $f \in \cS$, then 
$$
(T_n \varphi_k, f) = ( T_n\varphi_k  , f\varphi_{k+1}),
$$
since $D_{2k+1} \subset D_{2k+3}$. The test function on the right hand side now has compact support, so the left hand side converges as desired. Combining together \eqref{eq:conv1} and \eqref{eq:conv2} we deduce
$$
\Conv = \cap_{k\ge 1} \{ (T_n)_{n\ge 1}\in \cD_0'(\D)^\N: (T_n\varphi_k)_{n\ge 1} \in \Conv( \cS^*)\}
$$
 and thus $\Conv$ is a Borel set of the product $\sigma$-algebra by \eqref{eq:Frechetconv}, as desired.
\end{proof}

\newpage

\bibliographystyle{alpha}
\bibliography{RGM2}

\newpage

\printindex

\newpage
 \printindex[notation]

%\printnomenclature

\end{document}